\documentclass[a4paper,12pt]{article}
\usepackage[frenchb]{babel}
\usepackage[latin1]{inputenc}
\usepackage{pdfsync}
\usepackage[T1]{fontenc}
\usepackage{amsmath}
\usepackage{amssymb}
\newcommand{\cqfd}{\hfill $\Box$}
\newtheorem{defi}{Définition}[section]
      
      \newtheorem{souslem}[defi]{Sous-lemme}
      \newtheorem{prop}[defi]{Proposition}
      \newtheorem{prop-def}[defi]{Proposition-Définition}
      \newtheorem{def-prop}[defi]{Définition-Proposition}
      \newtheorem{thm}[defi]{Théorème}
      
      \newtheorem{lem}[defi]{Lemme}
      \newtheorem{cor}[defi]{Corollaire}

\def\del {\partial}

\def\ban{{\mathrm {ban}}}

\def\red{{\mathrm {red}}}
\def\geod{{\text {géod}}}
\def\top{{\mathrm {top}}}
\def\R{{\mathbb R}}
\def\Z{{\mathbb Z}}
\def\N{{\mathbb N}}
\def\C{{\mathbb C}}

\def\H{\mathcal{H}}

\def\P{\mathcal{P}}

\def\B{\mathcal{B}}

\def\K{\mathcal{K}}
\def\L{\mathcal{L}}

\def\Id{\mathrm {Id}}

\def\pp{\dots}

\def\tg{\text{-}\geod}
\def\de{\delta}

\def\vers{\rightarrow}

\newcommand{\s}[1]{\langle #1 \rangle}

\newcommand{\vag}[1]{\overset{#1}{\leftarrow}}
\newcommand{\vad}[1]{\overset{#1}{\rightarrow}}

\RequirePackage{graphics}

\begin{document}

\title{La conjecture de Baum-Connes à coefficients pour les groupes
  hyperboliques} 
\author{Vincent Lafforgue}
\maketitle

Le but de cet article est de montrer la conjecture de Baum-Connes à
coefficients  pour les groupes
hyperboliques au sens de Gromov~\cite{gromov-hyperb,delzanthyper,harpehyper}. 

Soit $G$ un groupe localement compact. La conjecture de Baum-Connes, formulée en 1982~\cite{preprintBaumConnes}  affirme que
$$\mu_{\red}^{G}: K_{*}^{\top}(G)\to K_{*}(C^{*}_{\red}(G))$$
est un isomorphisme de groupes abéliens. La conjecture de Baum-Connes à coefficients~\cite{baumconnes} affirme que pour toute $G$-$C^{*}$-algèbre $A$, 
$$\mu_{\red}^{G,A}: K_{*}^{\top}(G,A)\to K_{*}(C^{*}_{\red}(G,A))$$
est un isomorphisme de groupes abéliens.

Pour les groupes hyperboliques, 
l'injectivité de $\mu_{\red}^{G,A}$ a été montrée  par Kasparov et
  Skandalis~\cite{cras,ks}. 
La conjecture sans 
coefficients a été démontrée dans~\cite{kkban} et~\cite{mineyevyu}. La conjecture à coefficients commutatifs a été 
démontrée dans~\cite{kkbangpde} (mais la méthode de~\cite{kkbangpde} ne permet pas de montrer la conjecture de Baum-Connes à coefficients commutatifs pour un produit de deux groupes hyperboliques). 

Pour énoncer le théorème principal nous
avons besoin de quelques définitions. 

\begin{defi}\label{defi-hyperb-Hxyzt}
  Soit $\delta \geq 0$. Un espace métrique $(X,d)$ 
est dit $\delta$-hyperbolique
  si pour tout quadruplet $(x,y,z,t)$ de points de $X$ on a 
\begin{gather}\nonumber d(x,z)+d(y,t)\leq \max\big (d(x,y)+d(z,t),d(x,t)+d(y,z)\big)+\delta.\ \ \ \    (H_{\delta}(x,y,z,t))\end{gather}
\end{defi}
 %\label{prop-hyperb}
  \label{H{delta}(x,y,z,t)}
  
\begin{defi}
  Soit $\delta \geq 0$. Un espace métrique $(X,d)$ 
est dit faiblement $\delta$-géodésique si pour tous $x,y\in X$ et pour
tout $s\in [0,d(x,y)+\delta]$ il existe $z\in X$ tel que $d(x,z)\leq s$ et
$d(z,y)\leq d(x,y)-s+\delta $.
\end{defi}
% Dans~\cite{kkban} nous appelions faiblement 
% $\delta/2$-géodésique ce que nous appelons maintenant faiblement $\delta$-géodésique.
 
Un espace métrique $(X,d)$ est dit hyperbolique (resp. faiblement
géodésique) s'il existe $\delta \geq 0$ tel que $(X,d)$ soit 
$\delta$-hyperbolique (resp. faiblement $\delta$-géodésique). 

Lorsque $x\in X$ et $r\in \R_+$, on note 
$ B(x,r)=\{y\in X, d(x,y)\leq
r\}$.  \label{def-Bxr}

% (dans~\cite{kkban} cela désignait la boule ouverte mais ici nous n'aurons jamais besoin de considérer des boules ouvertes). 

\begin{defi}
Un espace métrique $(X,d)$ 
est dit uniformément localement fini si pour tout
  $r\in \R_+$ il existe $K\in \N$ tel que, pour tout $x\in X$, $
  B(x,r)$ contienne au plus $K$ points. 
\end{defi}

Le théorème principal est le suivant. 

\begin{thm}\label{enonce-ppal} Soit $G$ un groupe localement compact agissant de fa\c con
  isométrique, continue et propre sur un espace métrique  hyperbolique,
  faiblement géodésique et uniformément localement fini. 
Alors $G$ vérifie la conjecture de Baum-Connes à coefficients,
  c'est-à-dire que pour toute $G$-$C^*$-algèbre $A$, 
$\mu_{\red}^{G,A}:K_*^{\top}(G,A)\to K_*(C^*_{\red}(G,A))$ est une bijection. 
\end{thm}

%Le théorème découle des théorèmes~\ref{thhigson} et~\ref{thhyper}
On rappelle que l'injectivité de $\mu_{\red}^{G,A}$ est démontrée dans~\cite{cras,ks}. 

Tout groupe hyperbolique $\Gamma $ muni de la
métrique $d$ invariante à gauche associée à la longueur des mots
déterminée par  un système fini de générateurs est un espace métrique  hyperbolique, faiblement
géodésique et uniformément localement fini. Donc le théorème~\ref{enonce-ppal} implique la conjecture de Baum-Connes à coefficients pour les groupes hyperboliques. 

On aimerait remplacer dans le théorème~\ref{enonce-ppal} l'hypothèse ``uniformement localement fini''
 par l'hypothèse ``à géométrie grossière bornée'' qui est strictement plus faible (on renvoie à~\cite{ks} pour cette notion et on note que dans~\cite{ks} l'injectivité de  $\mu_{\red}^{G,A}$ est démontrée sous cette hypothèse plus faible). Cependant cela rendrait la démonstration encore  plus technique et nous y avons renoncé.

On aimerait aussi traiter le cas
général des groupoïdes hyperboliques, au moins à base compacte, mais cela serait très difficile, car il faudrait adapter la construction de
cet article (qui est assez combinatoire) aux techniques de Jean-Louis
Tu dans~\cite{tuhyper}, consistant à pondérer par des coefficients
tendant vers $0$ les éléments du groupoïde qui sont près de disparaître. 
En revanche, la conjecture sans coefficients pour ces groupoïdes est
beaucoup plus accessible par les méthodes de~\cite{kkbangpde}, elle
est d'ailleurs démontrée dans certains cas dans~\cite{kkbangpde}. 

D'après~\cite{contre-exemples} la conjecture de Baum-Connes à coefficients est fausse pour certains groupes 
aléatoires construits par Gromov dans~\cite{gromov}. 
Ces groupes sont des limites inductives de groupes hyperboliques, la limite étant indexée par $\N$, et les morphismes de transition étant surjectifs. Comme le membre de gauche de la conjecture de Baum-Connes commute aux limites inductives, on voit bien que ces contre-exemples sont ``dus'' au fait que $C^{*}_{\red}$ n'est pas fonctoriel en les morphismes de groupes (non nécessairement injectifs). Plus précisément  si $G\to H$ est un morphisme de groupes, et $A$ une $H$-$C^{*}$-algèbre, donc aussi une $G$-$C^{*}$-algèbre, le morphisme $C_{c}(G,A)\to C_{c}(H,A)$ ne se prolonge pas en général par continuité en un morphisme 
$C^{*}_{\red}(G,A)\to C^{*}_{\red}(H,A)$.

Voici maintenant quelques indications sur la démonstration du théorème~\ref{enonce-ppal}, qui occupe tout l'article. 
On commence par des rappels sur la méthode ``Dirac-dual Dirac'', inventée par Kasparov puis développée par Kasparov et Skandalis, et Higson et Kasparov. Cette méthode s'applique à une très large classe $\mathcal C$ de groupes localement compacts, dont la définition est rappelée dans l'introduction de~\cite{kkban}, et qui contient en particulier les groupes hyperboliques, les sous-groupes fermés des groupes réductifs sur un corps local, et les groupes ayant la propriété de Haagerup, c'est-à-dire possédant une action affine continue et propre sur un espace de Hilbert.  Pour tout groupe $G$ dans la classe $\mathcal C$, 
on possède un idempotent $\gamma\in KK_{G}(\C,\C)$ tel que 
 pour toute $G$-$C^{*}$-algèbre  $A$, $\mu_{\red}^{G,A}$ soit injectif et que l'image de l'action de $\gamma$ sur $K_*(C^*_{\red}(G,A))$ soit égale à l'image de $\mu_{\red}^{G,A}$. Donc si $G$ appartient à $\mathcal C$,  la conjecture de Baum-Connes à coefficients pour $G$ équivaut au fait que $\gamma$ agit par l'identité sur $K_*(C^*_{\red}(G,A))$. Si $G$ a la propriété de Haagerup, Higson et Kasparov ont montré dans~\cite{HigsonKasp} que $\gamma=1$ dans $KK_{G}(\C,\C)$.

Soit $G$ un groupe hyperbolique. 
 Comme certains groupes hyperboliques ont la propriété (T) de Kazhdan, on ne peut pas espérer montrer que $\gamma=1$ dans $KK_{G}(\C,\C)$. D'un autre côté
on sait d'après~\cite{kkban} que $\gamma$ est égal à $1$ dans $KK^{\ban}_{G,s\ell}(\C,\C)$ pour tout $s>0$ : cela permet de montrer la conjecture sans coefficients grâce à la propriété (RD) de Jolissaint~\cite{kkban}, et aussi à coefficients commutatifs grâce à un autre argument de stabilité par calcul fonctionnel holomorphe un peu plus subtil~\cite{kkbangpde}, mais  cela ne permet pas de montrer la conjecture à coefficients arbitraires.

L'idée  pour montrer la conjecture de Baum-Connes à coefficients pour les  groupes hyperboliques est que ceux-ci, même s'ils ont la propriété (T), ne vérifient pas la propriété (T) renforcée au sens de la définition 0.1 de~\cite{duke}. De fa\c con un peu imprécise
un groupe localement compact $G$ n'a pas la propriété (T) renforcée s'il existe une longueur $\ell$ sur $G$ (comme dans  la définition~\ref{deflong}) telle que pour tout $s>0$ il existe $C\in \R_{+}$ tel que la représentation triviale ne soit pas isolée parmi les représentations continues $\pi$ de $G$ dans des espaces de Hilbert vérifiant 
$\|\pi(g)\|\leq e^{C+s\ell(g)}$ pour tout $g\in G$. 
Le théorème 1.4 de~\cite{duke} affirme que si un groupe localement compact possède la propriété (T) renforcée  toute action continue et isométrique de ce groupe sur un espace métrique  
hyperbolique, faiblement géodésique et uniformément localement fini a des orbites bornées. La démonstration du théorème 1.4 de~\cite{duke} montre même que pour un groupe hyperbolique $\Gamma$ muni de la longueur $\ell$ associée à un système fini de générateurs,  il existe  un polynôme $P$ tel que la représentation triviale ne soit pas isolée parmi les représentations $\pi$ de $\Gamma$ dans des espaces de Hilbert vérifiant $\|\pi(g)\|\leq P(\ell(g))$  pour tout  $g\in \Gamma$.
 Notons que Ozawa~\cite{ozawa} a montré que les groupes  hyperboliques sont faiblement moyennables, ce qui amène à se demander si dans la phrase précédente on ne pourrait pas prendre pour $P$ un polynôme constant. Cependant cela n'apporterait rien pour la conjecture de Baum-Connes car la seule chose qui compte est que $P$ soit une fonction sous-exponentielle. 

Pour montrer le théorème~\ref{enonce-ppal} nous construisons une homotopie de $1$ à $\gamma$ en utilisant des représentations (continues) de $G$ dans des espaces de Hilbert qui ne sont pas unitaires mais à croissance exponentielle arbitrairement petite. Plus précisément  nous fixons une longueur $\ell$ sur $G$ et nous montrons que pour tout $s>0$ il existe $C\in \R_{+}$ tel que l'on puisse construire une homotopie de $1$ à $\gamma$ en utilisant des représentations $\pi$ de $G$ dans des espaces de Hilbert 
qui vérifient \begin{gather}\label{intro-croiss-rep}\|\pi(g)\|\leq e^{C+s\ell(g)} \text{ \  pour tout  } g\in G.\end{gather} 
Le théorème~\ref{thhyper} affirme l'existence pour tout $s>0$ d'une telle homotopie et 
le théorème~\ref{thhigson} (qui repose sur des idées de Nigel Higson) montre que cela  implique la conjecture de Baum-Connes à coefficients pour $G$.   La preuve du théorème~\ref{thhyper}  est ramenée à celle du théorème~\ref{hyperb-bon} où l'on suppose que $G$ agit proprement sur un espace hyperbolique $X$ vérifiant certaines propriétés supplémentaires (essentiellement que la métrique est associée à une structure de graphe). 
La preuve du théorème~\ref{hyperb-bon} repose sur l'acyclicité du complexe d'homologie simpliciale 
$$0\leftarrow \C^{(\Delta_0)}\vag{\partial }\C^{(\Delta_1)}\vag{\partial }\C^{(\Delta_2)}\dots \vag{\partial
  }\C^{(\Delta_{p_{\mathrm{max}}})}\leftarrow 0$$
du complexe de  Rips, où l'ensemble $\Delta_{p}$ des faces de dimension $p-1$ est formé des  parties de $X$ de cardinal $p$ et de diamètre $\leq N$, pour $N$ assez grand. On fixe un point base $x\in X$. La partie difficile est la construction de  $J_{x} : \C^{(\Delta_p)}\to \C^{(\Delta_{p+1})}$ tel que $\del J_{x}+J_{x}\del=1$
et de normes de Hilbert sur $\C^{(\Delta_p)}$ vérifiant \eqref{intro-croiss-rep}
 et telles que $\del$ et $J_{x}$ soient continus. L'homotopie de $1$ vers $\gamma$ se fait alors en conjuguant $\del+J_{x}$ par $e^{t\rho^{\flat}}$, où $\rho^{\flat}(a)$ est égal à $d(x,a)$ à une constante près et est   obtenu par un procédé de moyenne garantissant l'équivariance à compacts près 
  des opérateurs conjugués. Le lecteur qui voudrait se faire une idée rapide de la construction est invité à lire les paragraphes~\ref{para-structure-demo} et~\ref{para-cas-arbres}, l'introduction du paragraphe~\ref{construction-operateurs}, et les sous-paragraphes~\ref{sous-para-formule-norme} et~\ref{sous-para-enonce-resultats}.

La méthode que nous utilisons est semblable à celle utilisée par Pierre Julg pour montrer la conjecture de Baum-Connes à coefficients pour $Sp(n,1)$ (voir~\cite{julgcras}). En fait l'idée de  Julg d'utiliser des représentations non unitaires dans des espaces de Hilbert est très ancienne : dans~\cite{waterloo} Julg proposait  de construire une homotopie de $1$ à $\gamma$ en utilisant  des représentations uniformément bornées de $Sp(n,1)$.  
L'idée d'utiliser des représentations non pas uniformément bornées mais à petite croissance exponentielle a été dégagée lors de discussions avec Julg et Higson en 1999. 

Pour conclure voici un petit aper\c cu du statut actuel de la conjecture de Baum-Connes {\it à coefficients} $BC_{\mathrm {coeff}}$ pour des groupes $G$ de la classe $\mathcal C$ : 
\begin{itemize}
\item ``non T'' : si $G$ a la propriété de Haagerup, $G$ vérifie $BC_{\mathrm {coeff}}$ d'après~\cite{HigsonKasp}
\item ``T possible mais non T renforcé'' : $BC_{\mathrm {coeff}}$ est vrai si $G=Sp(n,1)$ d'après~\cite{julgcras}  ou si $G$ est un groupe hyperbolique par le présent article
\item ``T renforcé'' : dans ce cas, qui comprend probablement tous les groupes simples sur des corps locaux de rang déployé $\geq 2$ (et au moins ceux qui contiennent un $SL_{3}$ d'après~\cite{duke}), $BC_{\mathrm {coeff}}$ est totalement ouvert et ne pourra être résolu qu'avec des idées nouvelles comme le principe d'Oka (on renvoie à~\cite{anniversaire} pour plus de détails). 
\end{itemize}

Je remercie Georges Skandalis pour son aide et toutes les discussions que j'ai eues avec lui. Je remercie aussi Miguel Bermudez pour m'avoir indiqué le logiciel JPicEdit, avec lequel les dessins ont été réalisés, et Thomas Delzant pour m'avoir parlé du lemme d'approximation par les arbres, qui simplifie certaines démonstrations. Enfin je remercie vivement le rapporteur qui a tout lu en détail et indiqué de nombreuses corrections. 

\section{Structure de la démonstration}\label{para-structure-demo}

Le but de ce paragraphe est de ramener la démonstration du théorème~\ref{enonce-ppal} à celle du théorème~\ref{hyperb-bon}. Les paragraphes~\ref{construction-operateurs},~\ref{construction-normes} et~\ref{construction-fin} seront consacrés à la démonstration du théorème~\ref{hyperb-bon}.

Le théorème~\ref{thhigson} ci-dessous  affirme en gros que pour un groupe $G$ agissant proprement sur un espace hyperbolique,  l'existence d'homotopies de $1$ à $\gamma$, utilisant des représentations dans des espaces de Hilbert  dont la  croissance est contrôlée par une exponentielle arbitrairement petite,  implique la surjectivité de l'application de Baum-Connes à coefficients (l'injectivité est déjà connue
grâce à~\cite{ks}).  

\begin{defi}\label{deflong}
Soit $G$ un groupe localement compact. On appelle longueur sur $G$
une fonction continue $\ell:G\to \R_+$ vérifiant  $\ell(g^{-1})=\ell(g)$ et $\ell (g_1g_2)\leq
\ell (g_1)+\ell (g_2)$ pour tous $g,g_1,g_2\in G$. 
\end{defi}

Soit $G$ un groupe localement compact et $\ell$ une longueur sur $G$. 
Pour toutes $G$-$C^*$-algèbres $A$ et $B$  on
définit $E_{G,\ell}(A,B)$ comme l'ensemble des classes
d'isomorphisme de  $(E,\pi,T)$ où $E$ est un $(A,B)$-bimodule
hilbertien $\Z/2\Z$-gradué 
muni d'une action continue de $G$ vérifiant $\|\pi(g)\|\leq e^{\ell(g)}$ pour
tout $g\in G$, et d'un opérateur $T$ borné impair tel que pour tout $a\in A$
les opérateurs $[a,T]$ et $a(T^2-1)$ soient compacts et que
l'application $g\mapsto a(g(T)-T)$ soit une application normiquement continue de
$G$ dans $\K_{B}(E)$. On définit ensuite
$KK_{G,\ell}(A,B)$ comme l'ensemble des classes d'homotopie dans
$E_{G,\ell}(A,B)$ : deux éléments sont homotopes si ils sont les évaluations en $0$ et $1$ d'un élément de $E_{G,\ell}(A,B[0,1])$. 
On rappelle que  $B[0,1]=C([0,1],B)$ muni de la norme du sup. 
On peut montrer que la somme directe munit 
$KK_{G,\ell}(A,B)$  d'une structure de groupe abélien.

En particulier $E_{G,\ell}(\C,\C)$ est l'ensemble des classes
d'isomorphisme de  $(H,\pi,T)$ où $H$ est un espace de Hilbert $\Z/2\Z$-gradué 
muni d'une action continue de $G$ vérifiant $\|\pi(g)\|\leq e^{\ell(g)}$ pour
tout $g\in G$, et d'un opérateur $T$ borné impair tel que $(T^2-1)$ soit compact et que
l'application $g\mapsto g(T)-T$ soit une application normiquement continue de
$G$ dans  $\K(H)$.

\begin{thm}\label{thhigson}
 Soit $G$ un groupe localement compact agissant de fa\c con
  isométrique,  continue et propre  sur un espace métrique $(X,d)$ hyperbolique,
  faiblement géodésique et uniformément localement fini. Soit $\gamma \in KK_G(\C,\C)$ l'élément
  défini sous ces hypothèses par Kasparov et
  Skandalis~\cite{ks}.  Soit $x$ un
  point de $X$ et $\ell$ la longueur sur $G$ définie par
  $\ell(g)=d(x,gx)$. Supposons que  pour
  tout $s>0$ il existe $C\in \R_+$ tel que l'image de $1-\gamma$ dans
  $KK_{G,s\ell+C}(\C,\C)$ soit 
  nulle.
Alors $G$ vérifie la conjecture de Baum-Connes à coefficients,
  c'est-à-dire que pour toute $G$-$C^*$-algèbre $A$, 
$\mu_{\red}^{G,A}:K_*^{\top}(G,A)\to K_*(C^*_{\red}(G,A))$ est une bijection. 
\end{thm}

Ce  théorème est le  corollaire 2.12 de~\cite{anniversaire}, dont la preuve  repose sur des idées de Higson. En fait nous avons remplacé l'hypothèse ``à géométrie grossière bornée'' du  
 corollaire 2.12  de~\cite{anniversaire} par l'hypothèse ``uniformément localement fini'' qui est strictement plus forte, car nous appliquerons ce théorème à des espaces uniformément localement finis et qu'il n'est donc pas nécessaire de rappeler la notion de géométrie grossière bornée. 

Grâce au théorème~\ref{thhigson}, le théorème~\ref{enonce-ppal} est une conséquence du  théorème suivant. 

\begin{thm}\label{thhyper}
 Soit $G$ un groupe localement compact agissant de fa\c con
  isométrique, continue et propre sur un espace métrique $(X,d)$
  hyperbolique, faiblement géodésique et
  uniformément localement fini. Soit $\gamma\in KK_G(\C,\C)$ l'élément
  défini sous ces hypothèses par Kasparov et
  Skandalis~\cite{ks}. Soit $x\in  
  X$ et $\ell$ la longueur sur $G$ définie par
  $\ell(g)=d(x,gx)$.  Alors  pour
  tout $s>0$ il existe $C\in \R_+$ tel que l'image de $1-\gamma$ dans
  $KK_{G,s\ell+C}(\C,\C)$ est 
  nulle. 
\end{thm}
\noindent{\bf Remarque.} Le théorème~\ref{thhyper} est encore vrai sans l'hypothèse de propreté de l'action de $G$ sur $X$ (nous avons inclus cette hypothèse pour que $\gamma$ soit un ``élément $\gamma$'' pour $G$).  De toute fa\c con le théorème avec  l'hypothèse de propreté de l'action implique le théorème sans cette hypothèse car le groupe des automorphismes de $(X,d)$ est localement compact et agit proprement sur $X$.

Nous allons voir maintenant que le théorème~\ref{thhyper} résulte du théorème~\ref{hyperb-bon} ci-dessous qui utilise des hypothèses plus fortes sur l'espace métrique $(X,d)$.  

\begin{defi}
On dit qu'un espace métrique $(X,d)$ est un bon espace hyperbolique discret si
\begin{itemize}
\item $d$ prend ses valeurs dans $\N$ et est  géodésique, c'est-à-dire vérifie 
  $$\forall a,b\in X, \forall k\in \{0,...,d(a,b)\}, 
  \exists c\in X,
  d(a,c)=k, d(c,b)=d(a,b)-k, $$ autrement dit $d$ provient d'une structure de graphe connexe sur $X$, 
\item $(X,d)$ est uniformément localement fini (comme $d$ est  géodésique, cela équivaut à dire que le nombre  de points à distance $1$ d'un point, est borné indépendamment du point),
\item  $(X,d)$ est hyperbolique. 
\end{itemize}
\end{defi}
Remarquons que si $G$ est un groupe hyperbolique, et $d$ est la distance  invariante à gauche associée à la longueur des mots, pour un système fini de générateurs, $(G,d)$ est un bon espace hyperbolique discret, muni d'une action isométrique de $G$ par translations à gauche ($d$ provient de la structure de graphe de Cayley sur $G$ associée à ce système de générateurs).

Soit $\delta\in  \R_+^*$, et 
$(X,d)$  un espace métrique $\delta$-hyperbolique, faiblement $\delta$-géodésique, uniformément localement fini, et muni d'une action  isométrique d'un groupe  $G$. 
Munissons $X$ de la distance suivante : 
\begin{gather*}d'(a,b)=\min\{i\in \N, \exists\  a_0,...,a_i,\\ \text{\ \ tels  que \ }a_0=a,  a_i=b, \text{\ \ et \ }\forall
j\in \{0,...,i-1\}, d(a_j,a_{j+1})\leq \delta+1\}.\end{gather*}
Autrement dit $d'$ provient de la structure de graphe sur $X$ pour laquelle deux points distincts $x,y\in X$ sont voisins si $d(x,y)\leq \delta +1$. 
On a alors les
propriétés suivantes :
\begin{itemize}
\item l'action de $G$ sur $(X,d')$ est isométrique,
\item $d'$ est quasi-isométrique à $d$ : on a 
$d\leq (\de+1)d'$, et, en utilisant le fait que $(X,d)$ est faiblement $\delta$-géodésique, on montre facilement 
  que $d'\leq d+1$, 
  \item $(X,d')$ est un bon espace hyperbolique discret. 
\end{itemize}

L'hyperbolicité de $(X,d')$ résulte de la conservation de l'hyperbolicité
par quasi-isométrie pour des espaces faiblement géodésiques. Pour les espaces géodésiques la démonstration figure dans~\cite{harpehyper,delzanthyper} et pour les espaces faiblement géodésiques il n'y a pas grand chose à modifier. On peut aussi invoquer le théorème 3.18 de~\cite{vaisala} ainsi que  la remarque 3.19 de~\cite{vaisala} appliquée à $(X,d)$ et à l'espace total du graphe considéré précédemment (je remercie Yves Stalder qui m'a indiqué cette référence).  

Donc le théorème~\ref{thhyper} pour $(X,d)$  résulte du 
théorème~\ref{hyperb-bon} ci-dessous  appliqué à $(X,d')$.

\begin{thm}\label{hyperb-bon}
Soit $G$ un groupe localement compact agissant de fa\c con
  isométrique, continue et propre  sur un espace métrique $(X,d)$ qui est un bon espace hyperbolique discret.  Soit $\gamma\in KK_G(\C,\C)$ l'élément
  défini sous ces hypothèses par Kasparov et
  Skandalis~\cite{ks}.  Soit $x\in 
  X$ et $\ell$ la longueur sur $G$ définie par
  $\ell(g)=d(x,gx)$. Alors  pour
  tout $s>0$ il existe $C\in \R_+$ tel que l'image de $1-\gamma$ dans
  $KK_{G,s\ell+C}(\C,\C)$ est 
  nulle. 
\end{thm}

Pour montrer le théorème~\ref{enonce-ppal} nous sommes donc ramenés à montrer le théorème~\ref{hyperb-bon}. 
Les paragraphes~\ref{construction-operateurs},~\ref{construction-normes} et~\ref{construction-fin} sont consacrés à la démonstration du théorème~\ref{hyperb-bon}.

\section{Le cas des arbres}\label{para-cas-arbres}

Le but de  ce paragraphe est de  démontrer le théorème~\ref{hyperb-bon} dans le cas où $X$ est un arbre, afin d'introduire dans un cas
simple les idées de la démonstration du théorème~\ref{hyperb-bon}. 
Soit $X$ un arbre et $q$ un entier tel que chaque sommet de l'arbre
ait au plus $q+1$ voisins. Soit $G$ un groupe localement compact
agissant de fa\c con    isométrique, continue et propre sur $X$. 
Soit $\gamma\in KK_G(\C,\C)$ l'élément $\gamma$ de Julg et Valette, dont
la construction est rappelée ci-dessous. Soit $x \in X$ et
 $\ell$ la longueur  sur $G$ définie par
$\ell(g)=d(x,gx)$. 
La proposition suivante est le théorème~\ref{hyperb-bon}  dans le cas  où $X$ est un arbre.

\begin{prop} \label{cas-des-arbres}
Pour tout $s>0$ il existe $C\in \R_+$ tel que l'image de $1-\gamma$
dans $KK_{G,s\ell+C}(\C,\C)$ est nulle. 
\end{prop}

Bien sûr Julg et
Valette ont montré dans~\cite{julgvalette} que $\gamma=1$ dans
$KK_G(\C,\C)$, sans aucune hypothèse sur l'arbre (c'est-à-dire un résultat plus fort en partant d'une hypothèse plus faible).  

Rappelons la construction de l'élément $\gamma$ de Julg et Valette. 
Notons $\Delta _1=X$ l'ensemble des sommets de l'arbre et $\Delta _2$
l'ensemble des arêtes. Notons $\C^{(\Delta_1)}$ l'ensemble des
combinaisons finies d'éléments de $\Delta_1$. Pour simplifier les
notations, nous choisissons pour chaque arête une orientation, 
ce qui permet d'identifier $\C^{(\Delta_2)}$ au sous-espace vectoriel
de $\Lambda^2(\C^{(\Delta_1)})$ engendré par les $e_a\wedge e_b$ pour
$\{a,b\}\in \Delta_2$. De même $\ell^2(\Delta_2)$ s'identifie au sous-espace
vectoriel fermé de $\Lambda^2(\ell^2(\Delta_1))$ engendré par les
$e_a\wedge e_b$ pour $\{a,b\}\in \Delta_2$.

L'élément $\gamma\in KK_G(\C,\C)$ (associé au choix de $x$ comme origine) est  représenté par l'espace de Hilbert
$\Z/2\Z$-gradué $\ell^2(\Delta_1)\oplus \ell^2(\Delta_2)$ 
et par l'opérateur $T=\begin{pmatrix} 0 & u \\v & 0\end{pmatrix}$, où 
$u:\ell^2(\Delta_2)\to \ell^2(\Delta_1)$ et $v:\ell^2(\Delta_1)\to
\ell^2(\Delta_2)$ sont définis de la fa\c con suivante : si $\{a,b\}$ est
une arête telle que $b$ se trouve sur la géodésique entre $x$ et
$a$, alors $u(e_a\wedge e_b)=-e_a$ et $v(e_a)=-e_a\wedge e_b$ et enfin 
$v(e_{x})=0$. On voit que $1-v\circ u=0$ et que $1-u\circ v$ est le
projecteur orthogonal sur $e_{x}$. De plus, si $g\in G$, $g(u)-u$ et
$g(v)-v$ sont de rang fini.  Dans la suite nous écrirons toujours
$T=u+v$ au lieu de $T=\begin{pmatrix} 0 & u \\v & 0\end{pmatrix}$. 

\noindent{\bf Démonstration de la proposition~\ref{cas-des-arbres}.}
Nous allons montrer que pour tout $s>0$ il existe $C\in
\R_+$ tel que l'image de $1-\gamma$ 
dans $KK_{G,2s\ell+C}(\C,\C)$ soit nulle. 
Soit $s\in \R_+^*$. En fait nous allons montrer qu'il existe un
polynôme $P$ et une homotopie de $1$ à $\gamma$ faisant
intervenir des représentations $\pi$ de $G$ dans des espaces 
de Hilbert tels que $\|\pi(g)\|\leq P(\ell(g))e^{s\ell(g)}$. 

Notons $\del: \C^{(\Delta_2)}\to \C^{(\Delta_1)}$ l'opérateur bord de
l'homologie simpliciale. Pour toute arête $\{a,b\}$ on a
$\del(e_a\wedge e_b)=e_b-e_a$. Alors $\del$ est injectif et son
image  est de codimension $1$. Pour le montrer définissons 
$h:\C^{(\Delta_1)}\to \C^{(\Delta_2)}$ de la fa\c con suivante : 
pour tout point $a\in X$, notons $n=d(a,x)$ et $a_{0},\dots,a_{n}$ la suite de points reliant $a$ à $x$ (c'est-à-dire que $a_{0}=a$, $a_{n}=x$ et $a_{i}$ est voisin de $a_{i-1}$ pour tout $i\in \{1,\dots , n\}$) et posons 
$h(e_{a})=-(e_{a_0}\wedge e_{a_1}+...+e_{a_{n-1}}\wedge e_{a_n})$. En
particulier $h(e_{x})=0$. Alors $h\circ \del=\Id_{\C^{(\Delta_{2})}}$ et $\Id_{\C^{(\Delta_{1})}}-\del\circ h$
est de rang $1$. 

Plus généralement pour tout $t\in[0,1]$ on définit des opérateurs
$u_t: \C^{(\Delta_2)}\to \C^{(\Delta_1)}$ et 
$v_t:\C^{(\Delta_1)}\to \C^{(\Delta_2)}$
 de sorte que 
$u_1=\del$, $v_1=h$, et $u_0$ et $v_0$ soient les restrictions de $u$
 et $v$ à $\C^{(\Delta_2)}\subset \ell^2(\Delta_2)$ et $\C^{(\Delta_1)}\subset \ell^2(\Delta_1)$. 
La formule est la suivante : 
pour toute arête $\{a,b\}$ telle que $b$ se trouve sur la géodésique
entre $x$ et $a$ on pose $u_t(e_a\wedge e_b)=-(e_a-te_b)$. Pour tout point $a\in X$, on note $n=d(a,x)$ et $a_{0},\dots,a_{n}$ la suite de points reliant $a$ à $x$ (c'est-à-dire que $a_{0}=a$, $a_{n}=x$ et $a_{i}$ est voisin de $a_{i-1}$ pour tout $i\in \{1,\dots , n\}$) et  on pose 
$$v_t(e_{a})=-(e_{a_0}\wedge e_{a_1}+te_{a_1}\wedge
e_{a_2}+t^2e_{a_2}\wedge e_{a_3}+...+t^{n-1}e_{a_{n-1}}\wedge
e_{a_n}).$$ 

Nous allons compléter $\C^{(\Delta_2)}$ et $\C^{(\Delta_1)}$ pour certaines normes   de
 Hilbert  telles que $\|\pi(g)\|\leq P(\ell(g))e^{s\ell(g)}$ pour un
 certain polynôme $P$  et telles que les opérateurs
 $u_t$ et $v_t$ soient continus, uniformément en $t\in [0,1]$. 

Rappelons que dans~\cite{kkban} nous introduisons les espaces 
$\ell^1_{x,s}(\Delta_1)$ et $\ell^1_{x,s}(\Delta_2)$ comme les complétés de
$\C^{(\Delta_1)}$ et $\C^{(\Delta_2)}$ pour les normes $\ell^1$ 
pondérées suivantes : 
\begin{gather*}\Big\|\sum_{a\in \Delta_1} f(a)e_a\Big\|_{\ell^1_{x,s}(\Delta_1)}=\sum_{a\in \Delta_1}
|f(a)|e^{sd(x,a)}
\text{\ \ \ \ \  et }\\ \Big\|\sum_{\{a,b\}\in \Delta_2,\   b\in \geod(x,a)}
 f(a,b)e_a\wedge e_b\Big\|_{\ell^1_{x,s}(\Delta_2)}=\sum_{\{a,b\}\in \Delta_2,\   b\in \geod(x,a)}
|f(a,b)|e^{sd(x,a)}.\end{gather*}

On voit que pour tout $t\in [0,1]$ on a 
$\|u_t\|_{\L(\ell^1_{x,s}(\Delta_2),\ell^1_{x,s}(\Delta_1))}=1+te^{-s}$ et
$\|v_t\|_{\L(\ell^1_{x,s}(\Delta_1),\ell^1_{x,s}(\Delta_2))}=(1-te^{-s})^{-1}$. 

Rappelons la construction de~\cite{kkban} qui montre que l'image de $1-\gamma$ dans $KK^{\ban}_{G,s\ell}(\C,\C)$ est
nulle. On note $E_{x,s}(\Delta_1)$ la $\C$-paire $(c_{0,x,s}(\Delta_1),\ell^1_{x,s}(\Delta_1))$, où $c_{0,x,s}(\Delta_1)$ est le complété de 
$\C^{(\Delta_{1})}$ pour la norme $\|f\|=\sup _{a\in \Delta_1}
|f(a)|e^{-sd(x,a)}$ et  où le crochet entre 
$c_{0,x,s}(\Delta_1)$ et $\ell^1_{x,s}(\Delta_1)$ est donné par $\s{f,f'}=\sum_{a\in \Delta_{1}}f(a)f'(a)$.  On introduit de la même fa\c con la 
$\C$-paire $E_{x,s}(\Delta_2)$.  
Pour toute $\C$-paire $E=(E^{<},E^{>})$ on note 
$E[0,1]$ la $\C$-paire $(E^{<}[0,1],E^{>}[0,1])$ (où $E^{<}[0,1]=C([0,1],E^{<})$ muni de la norme du sup, et de même pour $E^{>}[0,1]$). 
 Alors 
  $(E_{x,s}(\Delta_1)[0,1]\oplus E_{x,s}(\Delta_2)[0,1], 
(u_{t}+v_{t})_{t\in [0,1]})$ définit un élément de $E^{\ban}_{G,s\ell}(\C,\C[0,1])$. 
D'autre part $(E_{x,s}(\Delta_1)\oplus E_{x,s}(\Delta_2), 
u_{1}+v_{1})$
est égal  à $1$ dans $KK^{\ban}_{G,s\ell}(\C,\C)$ grâce au lemme 1.4.2   de~\cite{kkban} et   $(E_{x,s}(\Delta_1)\oplus E_{x,s}(\Delta_2), 
u_{0}+v_{0})$  est égal  à 
$\gamma$ dans $KK^{\ban}_{G,s\ell}(\C,\C)$ (pour l'homotopie entre ces deux éléments, on garde les opérateurs, on complète $\C^{(\Delta_{1})}$ pour la norme $t\|.\|_{E_{x,s}(\Delta_1)}+(1-t) \|.\|_{\ell^2(\Delta_1)}$, on considère la $\C$-paire formée de cet espace de Banach et du complété de $\C^{(\Delta_{1})}$ pour la norme duale et on fait de même pour $\Delta_{2}$). 
Par conséquent on voit 
que l'image de $1-\gamma$ dans $KK^{\ban}_{G,s\ell}(\C,\C)$ est
nulle. 

Pour montrer que l'image de $1-\gamma$ dans $KK_{G,2s\ell+C}(\C,\C)$ est
nulle (avec $C$ une constante assez grande), la première idée qui vient est de remplacer  les normes $\ell^1$ par des normes $\ell^2$ pour avoir des espaces de
Hilbert. 

Nous définissons donc  
$\ell^2_{x,s}(\Delta_1)$ et $\ell^2_{x,s}(\Delta_2)$ comme les complétés de
$\C^{(\Delta_1)}$ et $\C^{(\Delta_2)}$ pour les normes $\ell^2$ 
pondérées suivantes : 
\begin{gather*}\Big\|\sum_{a\in \Delta_1} f(a) e_{a}\Big\|_{\ell^2_{x,s}(\Delta_1)}^2=\sum_{a\in \Delta_1}
|f(a)|^2e^{2sd(x,a)}
\text{\ \ \ \ \ \  et }\\ \Big\|\sum_{\{a,b\}\in \Delta_2,\   b\in \geod(x,a)} 
f(a,b)e_a\wedge e_b\Big\|_{\ell^2_{x,s}(\Delta_2)}^2=\sum_{\{a,b\}\in \Delta_2,\   b\in \geod(x,a)} 
|f(a,b)|^2e^{2sd(x,a)}.\end{gather*}

Il est clair que $u_t:\ell^2_{x,s}(\Delta_2)\to \ell^2_{x,s}(\Delta_1)$ est continu
(et que sa norme est majorée de fa\c con uniforme en $t$), mais ce n'est pas
le cas de $v_t$, si $s$ est  petit et $t$ proche de $1$, pour la raison suivante. 

Soit $\{z,z'\}$ une arête de $X$ (telle que $z'\in \geod(x,z)$) et  $$\xi_{z,z'} : 
\sum_{\{a,b\}\in \Delta_2,\   b\in \geod(x,a)} 
f(a,b)e_a\wedge e_b \mapsto f(z,z') $$ la forme linéaire sur 
$\ell^2_{x,s}(\Delta_2)$ qui donne le coefficient de cette arête. 
Cette forme linéaire est de norme $e^{-sd(x,z)}$. 
Or la forme linéaire ${}^{t}v_{t}(\xi_{z,z'})=\xi_{z,z'}\circ v_t$ sur $\ell^2_{x,s}(\Delta_1)$ est  $$\sum_a f(a)e_a\mapsto -\Big(\sum_{a\text{\   tel  que\ }z\in \text{géod}(x,a)}
t^{d(z,a)}f(a)\Big). $$
Si chaque sommet de l'arbre a exactement $q+1$ voisins, en notant
$k=d(x,z)$ on a 
$\|{}^{t}v_{t}(\xi_{z,z'})\|_{\ell^2_{x,s}(\Delta_1)^{*}}^2=\sum_{n\geq k} t^{2(n-k)}e^{-2sn}q^{n-k}$. 
On voit donc que cette forme linéaire ${}^{t}v_{t}(\xi_{z,z'})$ n'est bornée que si
$e^{-s}t\sqrt q<1$. 

Nous allons 
compléter 
 $\C^{(\Delta_1)}$ et $\C^{(\Delta_2)}$ pour  d'autres normes de Hilbert de telle sorte  que les opérateurs
$u_t$ et $v_t$ soient continus pour $t\in [0,1]$. 

Pour tout $n\in \N$,  on note $S^n_{x}=\{a\in X, d(x,a)=n\}$ la sphère de
rayon $n$ et de centre $x$. Pour $k,n\in \N$ avec $k\leq n$ et
pour tout $z\in S^k_{x}$ on note $I^{n,k,x}_z=\{a\in S^n_{x}, z\in \geod(x,a)\}$
de sorte que $S^n_{x}=\bigcup_{z\in S^k_{x}} I^{n,k,x}_z$ est une partition de la sphère $S^n_{x}$. Lorsque $k=n$ c'est la partition par les singletons et lorsque
$k=0$ c'est la partition grossière. 

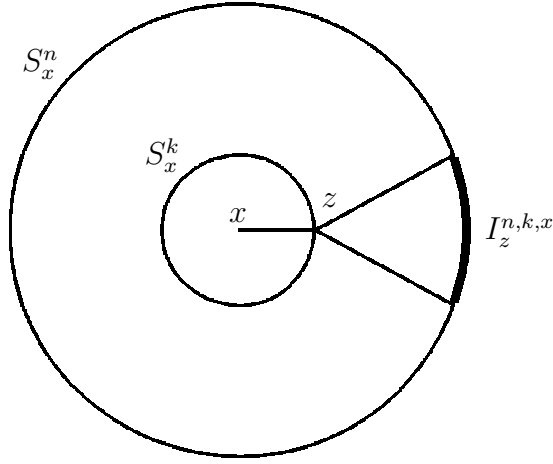
\ifx\JPicScale\undefined\def\JPicScale{1}\fi
\unitlength \JPicScale mm


On définit alors $H_{x,s}(\Delta_1)$ comme l'espace de Hilbert, complétion
de $\C^{(\Delta_1)}$ pour la norme 
$$\Big\|\sum_a f(a)e_{a}\Big\|^2_{H_{x,s}(\Delta_1)}=\sum_{n\in \N} e^{2sn}
\sum_{k=0}^n \sum_{z\in S^k_{x}} \Big|\sum_{a\in I^{n,k,x}_z} f(a)\Big|^2.$$

Alors pour tout $t\in [0,1]$, $v_t:H_{x,s}(\Delta_1)\to \ell^2_{x,s}(\Delta_2)$
est continu. En effet $v_t=\sum_{j\in \N}v_{t,j}$ où $v_{t,j}$ est défini de la manière suivante : pour 
tout point $a\in X$, en notant $n=d(a,x)$ et $a_{0},\dots,a_{n}$ la suite de points reliant $a$ à $x$ (c'est-à-dire que $a_{0}=a$, $a_{n}=x$ et $a_{i}$ est voisin de $a_{i-1}$ pour tout $i\in \{1,\dots , n\}$),  on pose 
$v_{t,j}(e_{a})=-t^je_{a_j}\wedge e_{a_{j+1}}$ si $j\leq n-1$ et $v_{t,j}(e_{a})=0$ si $j\geq n$. On vérifie facilement
que $v_{t,j}$ est continu de $H_{x,s}(\Delta_1)$ vers $\ell^2_{x,s}(\Delta_2)$ de norme inférieure ou égale à
$t^je^{-sj}$. La raison est que $\|.\|_{\ell^2_{x,s}(\Delta_2)}^{2}$ est une somme pondérée des carrés des formes linéaires $\xi_{z,z'}$, qu'en notant $k=d(x,z)$, ${}^{t}v_{t,j}(\xi_{z,z'})$ est le produit par $-t^{j}$ de la forme linéaire $f\mapsto \sum_{a\in I^{j+k,k,x}_z} f(a)$ et que $\|.\|_{H_{x,s}(\Delta_1)}^{2}$ est une somme pondérée des carrés de telles formes linéaires. 

Bien sûr $\del$ et plus généralement $u_t:\ell^2_{x,s}(\Delta_2)\to
H_{x,s}(\Delta_1)$ ne sont plus continus et nous devons remplacer
$\ell^2_{x,s}(\Delta_2)$ par un espace analogue à $H_{x,s}(\Delta_1)$. 

Pour $k\in \N, n\in \N^*$ avec $k\leq n$ et $z\in S^k_{x}$ on pose 
$$J^{n,k,x}_z=\{(a,b), \{a,b\}\in \Delta_2, a\in S^n_{x}, b\in S^{n-1}_{x},  z\in
\geod (x,a)\}.$$
On remarque que $\bigcup _{z\in S^k_{x}}J^{n,k,x}_z$ est une partition de
l'ensemble des arêtes à distance $n-1$ de $x$. Lorsque $k=n$ c'est la partition par les singletons et lorsque
$k=0$ c'est la partition grossière. 

On définit alors $H_{x,s}(\Delta_2)$ comme la complétion de
$\C^{(\Delta_2)}$ pour la norme 
$$\Big\|\sum_{\{a,b\}\in \Delta_2, \  b\in \geod(x,a)} 
f(a,b)e_a\wedge e_b\Big\|^2_{H_{x,s}(\Delta_2)}=\sum_{n\in \N^*} e^{2sn}
\sum_{k=0}^{n} \sum_{z\in S^k_{x}} \Big|\sum_{(a,b)\in J^{n,k,x}_z}
f(a,b)\Big|^2.$$ 
 
Il est alors très facile de voir que pour tout $t\in [0,1]$, les
opérateurs $u_t:H_{x,s}(\Delta_2) \to H_{x,s}(\Delta_1)$ et
$v_t:H_{x,s}(\Delta_1) \to H_{x,s}(\Delta_2)$ sont continus, de normes majorées
uniformément en $t$. Pour $v_{t}$ la raison est la suivante : l'image par ${}^{t}v_{t,j}$  de la forme linéaire $f\mapsto  \sum_{(a,b)\in J^{n,k,x}_z}
f(a,b)$ qui apparaît dans la formule ci-dessus est égale au produit par $-t^{j}$ de  la forme linéaire 
$f\mapsto  \sum_{a\in I^{n+j,k,x}_z} f(a)$ qui 
apparaît dans la formule pour 
$\|.\|_{H_{x,s}(\Delta_1)}^{2}$ et on en déduit facilement $\|v_{t,j}\|
_{\L(H_{x,s}(\Delta_1),H_{x,s}(\Delta_2))}\leq t^{j}e^{-sj}$.

Pour conclure il ne reste donc plus qu'à établir le résultat suivant : 
il existe un polynôme $P$ tel que pour tout $s>0$ les représentations
de $G$ sur $H_{x,s}(\Delta_1)$ et $H_{x,s}(\Delta_2)$ vérifient
$\|\pi(g)\|\leq P(\ell(g))e^{s\ell(g)}$ pour tout $g\in G$. 
 
Démontrons ce résultat pour $H_{x,s}(\Delta_1)$. 
Soit $g\in G$. On pose $x'=g(x)$, si bien  que
$d(x,x')=\ell(g)$.   
Il résulte de la définition de $H_{x,s}(\Delta_1)$ que l'on a une
isométrie $\Theta _{x,s}:H_{x,s}(\Delta_1)\to \ell^2(\{(n,k,z), 0\leq k\leq n,
z\in S^k_{x}\})$ définie par $$\Theta_{x,s}(\sum_a
f(a) e_{a})=\Big(e^{sn}\sum_{a\in I^{n,k,x}_z} f(a)\Big)_{(n,k,z)}.$$  
Pour tout $f\in \C^{(\Delta_1)}$ on a $\|\pi(g)f\|_{H_{x',s}(\Delta_1)}=
\|f\|_{H_{x,s}(\Delta_1)}$. D'autre part on a une  isométrie $\Theta _{x',s}$ de $H_{x',s}(\Delta_1)$ dans $\ell^2(\{(n',k',z'), 0\leq k'\leq n',
z'\in S^{k'}_{x'}\})$. 

 Le lemme suivant  permet de comparer les normes de $H_{x,s}$ et de $H_{x',s}$. 

\begin{lem}
Pour $(n,k,z)\in \{(n,k,z), 0\leq k\leq n,
z\in S^{k}_{x}\}$ on peut écrire $I^{n,k,x}_{z}$ comme une réunion disjointe 
finie d'au plus $(q-1)(d(x,x')+1)+2$ parties $I^{n',k',x'}_{z'}$ pour
$(n',k',z') \in \{(n',k',z'), 0\leq k'\leq n', z'\in S^{k'}_{x'}\}$,
de sorte que pour chaque $(n',k',z')$ la partie $I^{n',k',x'}_{z'}$
intervienne au plus $d(x,x')+2$  fois dans  ces
décompositions (c'est-à-dire que l'ensemble des $(n,k,z)$ tels que la partie $I^{n',k',x'}_{z'}$
intervienne dans la décomposition de $I^{n,k,x}_{z}$ est fini et de cardinal inférieur ou égal à $d(x,x')+2$)
\end{lem}
Ce lemme est le  lemme 1.5 de~\cite{duke} mais nous en rappelons la démonstration (avec une petite correction).   

\noindent{\bf Démonstration.}
 Si $z$ n'appartient pas à $\geod(x,x')$, on a $I^{n,k,x}_{z}=I^{n',k',x'}_{z'}$ avec $z'=z$, $k'=d(x',z)$,
$n'=n+k'-k$. Si $z$ appartient à $\geod(x,x')$, on peut écrire 
$I^{n,k,x}_{z}$ comme la réunion disjointe 
\begin{itemize}
\item 
des $I^{n',k',x'}_{z'}$, avec $z'$
n'appartenant pas à $\geod(x,x')$ mais  à distance $1$ d'un point
de $\geod(z,x')$,  $k'=d(x',z')$ et 
$n'=n-d(x,x')+2(k'-1)$, sous réserve que $n'\geq0$, 
\item et du singleton $I_{z'}^{k',k',x'}=\{z'\}$, où $z'\in \geod(z,x')$ vérifie $d(x,z')=n$, et avec $k'=d(x',z')$, si $n\leq d(x,x')$. 
\end{itemize} 

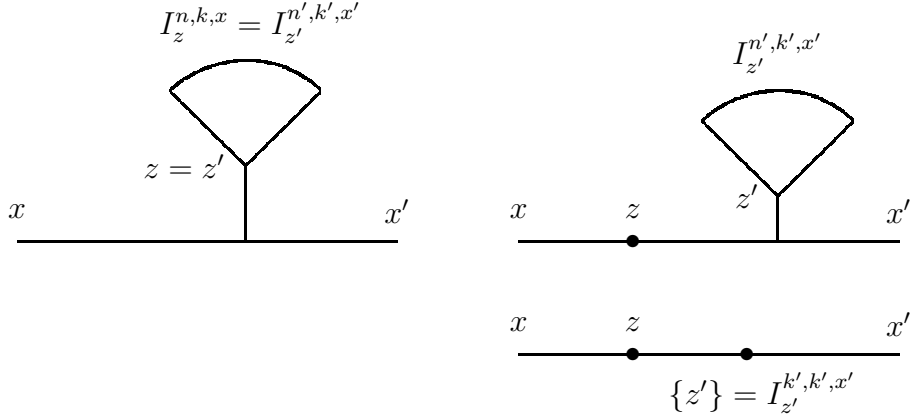
\ifx\JPicScale\undefined\def\JPicScale{1}\fi
\unitlength \JPicScale mm
\begin{picture}(126,70)(10,10)
\linethickness{0.3mm}
\put(10,40){\line(1,0){50}}
\linethickness{0.3mm}
\put(76,40){\line(1,0){50}}
\linethickness{0.3mm}
\multiput(49.63,60.34)(0.12,-0.11){3}{\line(1,0){0.12}}
\multiput(49.25,60.67)(0.13,-0.11){3}{\line(1,0){0.13}}
\multiput(48.86,60.98)(0.13,-0.1){3}{\line(1,0){0.13}}
\multiput(48.46,61.28)(0.2,-0.15){2}{\line(1,0){0.2}}
\multiput(48.04,61.56)(0.21,-0.14){2}{\line(1,0){0.21}}
\multiput(47.62,61.84)(0.21,-0.14){2}{\line(1,0){0.21}}
\multiput(47.19,62.09)(0.22,-0.13){2}{\line(1,0){0.22}}
\multiput(46.75,62.33)(0.22,-0.12){2}{\line(1,0){0.22}}
\multiput(46.3,62.56)(0.22,-0.11){2}{\line(1,0){0.22}}
\multiput(45.85,62.77)(0.23,-0.11){2}{\line(1,0){0.23}}
\multiput(45.39,62.96)(0.23,-0.1){2}{\line(1,0){0.23}}
\multiput(44.92,63.14)(0.47,-0.18){1}{\line(1,0){0.47}}
\multiput(44.44,63.3)(0.47,-0.16){1}{\line(1,0){0.47}}
\multiput(43.96,63.45)(0.48,-0.15){1}{\line(1,0){0.48}}
\multiput(43.48,63.58)(0.48,-0.13){1}{\line(1,0){0.48}}
\multiput(42.99,63.69)(0.49,-0.11){1}{\line(1,0){0.49}}
\multiput(42.5,63.78)(0.49,-0.09){1}{\line(1,0){0.49}}
\multiput(42,63.86)(0.5,-0.08){1}{\line(1,0){0.5}}
\multiput(41.5,63.92)(0.5,-0.06){1}{\line(1,0){0.5}}
\multiput(41,63.97)(0.5,-0.04){1}{\line(1,0){0.5}}
\multiput(40.5,63.99)(0.5,-0.03){1}{\line(1,0){0.5}}
\multiput(40,64)(0.5,-0.01){1}{\line(1,0){0.5}}
\multiput(39.5,63.99)(0.5,0.01){1}{\line(1,0){0.5}}
\multiput(39,63.97)(0.5,0.03){1}{\line(1,0){0.5}}
\multiput(38.5,63.92)(0.5,0.04){1}{\line(1,0){0.5}}
\multiput(38,63.86)(0.5,0.06){1}{\line(1,0){0.5}}
\multiput(37.5,63.78)(0.5,0.08){1}{\line(1,0){0.5}}
\multiput(37.01,63.69)(0.49,0.09){1}{\line(1,0){0.49}}
\multiput(36.52,63.58)(0.49,0.11){1}{\line(1,0){0.49}}
\multiput(36.04,63.45)(0.48,0.13){1}{\line(1,0){0.48}}
\multiput(35.56,63.3)(0.48,0.15){1}{\line(1,0){0.48}}
\multiput(35.08,63.14)(0.47,0.16){1}{\line(1,0){0.47}}
\multiput(34.61,62.96)(0.47,0.18){1}{\line(1,0){0.47}}
\multiput(34.15,62.77)(0.23,0.1){2}{\line(1,0){0.23}}
\multiput(33.7,62.56)(0.23,0.11){2}{\line(1,0){0.23}}
\multiput(33.25,62.33)(0.22,0.11){2}{\line(1,0){0.22}}
\multiput(32.81,62.09)(0.22,0.12){2}{\line(1,0){0.22}}
\multiput(32.38,61.84)(0.22,0.13){2}{\line(1,0){0.22}}
\multiput(31.96,61.56)(0.21,0.14){2}{\line(1,0){0.21}}
\multiput(31.54,61.28)(0.21,0.14){2}{\line(1,0){0.21}}
\multiput(31.14,60.98)(0.2,0.15){2}{\line(1,0){0.2}}
\multiput(30.75,60.67)(0.13,0.1){3}{\line(1,0){0.13}}
\multiput(30.37,60.34)(0.13,0.11){3}{\line(1,0){0.13}}
\multiput(30,60)(0.12,0.11){3}{\line(1,0){0.12}}
\linethickness{0.3mm}
\put(40,40){\line(0,1){10}}
\linethickness{0.3mm}
\multiput(40,50)(0.12,0.12){83}{\line(1,0){0.12}}
\linethickness{0.3mm}
\multiput(30,60)(0.12,-0.12){83}{\line(1,0){0.12}}
\linethickness{0.3mm}
\multiput(119.63,56.34)(0.12,-0.11){3}{\line(1,0){0.12}}
\multiput(119.25,56.67)(0.13,-0.11){3}{\line(1,0){0.13}}
\multiput(118.86,56.98)(0.13,-0.1){3}{\line(1,0){0.13}}
\multiput(118.46,57.28)(0.2,-0.15){2}{\line(1,0){0.2}}
\multiput(118.04,57.56)(0.21,-0.14){2}{\line(1,0){0.21}}
\multiput(117.62,57.84)(0.21,-0.14){2}{\line(1,0){0.21}}
\multiput(117.19,58.09)(0.22,-0.13){2}{\line(1,0){0.22}}
\multiput(116.75,58.33)(0.22,-0.12){2}{\line(1,0){0.22}}
\multiput(116.3,58.56)(0.22,-0.11){2}{\line(1,0){0.22}}
\multiput(115.85,58.77)(0.23,-0.11){2}{\line(1,0){0.23}}
\multiput(115.39,58.96)(0.23,-0.1){2}{\line(1,0){0.23}}
\multiput(114.92,59.14)(0.47,-0.18){1}{\line(1,0){0.47}}
\multiput(114.44,59.3)(0.47,-0.16){1}{\line(1,0){0.47}}
\multiput(113.96,59.45)(0.48,-0.15){1}{\line(1,0){0.48}}
\multiput(113.48,59.58)(0.48,-0.13){1}{\line(1,0){0.48}}
\multiput(112.99,59.69)(0.49,-0.11){1}{\line(1,0){0.49}}
\multiput(112.5,59.78)(0.49,-0.09){1}{\line(1,0){0.49}}
\multiput(112,59.86)(0.5,-0.08){1}{\line(1,0){0.5}}
\multiput(111.5,59.92)(0.5,-0.06){1}{\line(1,0){0.5}}
\multiput(111,59.97)(0.5,-0.04){1}{\line(1,0){0.5}}
\multiput(110.5,59.99)(0.5,-0.03){1}{\line(1,0){0.5}}
\multiput(110,60)(0.5,-0.01){1}{\line(1,0){0.5}}
\multiput(109.5,59.99)(0.5,0.01){1}{\line(1,0){0.5}}
\multiput(109,59.97)(0.5,0.03){1}{\line(1,0){0.5}}
\multiput(108.5,59.92)(0.5,0.04){1}{\line(1,0){0.5}}
\multiput(108,59.86)(0.5,0.06){1}{\line(1,0){0.5}}
\multiput(107.5,59.78)(0.5,0.08){1}{\line(1,0){0.5}}
\multiput(107.01,59.69)(0.49,0.09){1}{\line(1,0){0.49}}
\multiput(106.52,59.58)(0.49,0.11){1}{\line(1,0){0.49}}
\multiput(106.04,59.45)(0.48,0.13){1}{\line(1,0){0.48}}
\multiput(105.56,59.3)(0.48,0.15){1}{\line(1,0){0.48}}
\multiput(105.08,59.14)(0.47,0.16){1}{\line(1,0){0.47}}
\multiput(104.61,58.96)(0.47,0.18){1}{\line(1,0){0.47}}
\multiput(104.15,58.77)(0.23,0.1){2}{\line(1,0){0.23}}
\multiput(103.7,58.56)(0.23,0.11){2}{\line(1,0){0.23}}
\multiput(103.25,58.33)(0.22,0.11){2}{\line(1,0){0.22}}
\multiput(102.81,58.09)(0.22,0.12){2}{\line(1,0){0.22}}
\multiput(102.38,57.84)(0.22,0.13){2}{\line(1,0){0.22}}
\multiput(101.96,57.56)(0.21,0.14){2}{\line(1,0){0.21}}
\multiput(101.54,57.28)(0.21,0.14){2}{\line(1,0){0.21}}
\multiput(101.14,56.98)(0.2,0.15){2}{\line(1,0){0.2}}
\multiput(100.75,56.67)(0.13,0.1){3}{\line(1,0){0.13}}
\multiput(100.37,56.34)(0.13,0.11){3}{\line(1,0){0.13}}
\multiput(100,56)(0.12,0.11){3}{\line(1,0){0.12}}

\linethickness{0.3mm}
\put(110,40){\line(0,1){6}}
\linethickness{0.3mm}
\multiput(110,46)(0.12,0.12){83}{\line(1,0){0.12}}
\linethickness{0.3mm}
\multiput(100,56)(0.12,-0.12){83}{\line(1,0){0.12}}
\put(10,44){\makebox(0,0)[cc]{$x$}}

\put(60,44){\makebox(0,0)[cc]{$x'$}}

\put(76,44){\makebox(0,0)[cc]{$x$}}

\put(126,44){\makebox(0,0)[cc]{$x'$}}

\put(32,50){\makebox(0,0)[cc]{$z=z'$}}

\put(42,69){\makebox(0,0)[cc]{$I_z^{n,k,x}=I_{z'}^{n',k',x'}$}}

\put(110,65){\makebox(0,0)[cc]{$I_{z'}^{n',k',x'}$}}

\put(106,46){\makebox(0,0)[cc]{$z'$}}

\put(91,40){\makebox(0,0)[cc]{$\bullet$}}

\put(91,44){\makebox(0,0)[cc]{$z$}}

%\end{picture}
%
%\ifx\JPicScale\undefined\def\JPicScale{1}\fi
%\unitlength \JPicScale mm
%\begin{picture}(80,15)(-16,25)
\linethickness{0.3mm}
\put(76,25){\line(1,0){50}}
\put(76,29){\makebox(0,0)[cc]{$x$}}

\put(91,25){\makebox(0,0)[cc]{$\bullet$}}

\put(91,29){\makebox(0,0)[cc]{$z$}}

\put(106,25){\makebox(0,0)[cc]{$\bullet$}}

\put(108,20){\makebox(0,0)[cc]{$\{z'\}=I_{z'}^{k',k',x'}$}}

\put(126,29){\makebox(0,0)[cc]{$x'$}}

\end{picture}

\noindent 
Le dessin illustre les deux cas envisagés dans la démonstration (le dessin de gauche correspond au  cas où $z\not\in\geod(x,x')$ et les deux dessins de droite correspondent au cas où $z\in\geod(x,x')$). 
 \cqfd

\noindent{\bf Fin de la démonstration de la proposition~\ref{cas-des-arbres} en admettant le lemme~\ref{cauchy-schwarz}.}
Soit $A$ la  matrice de  $$\ell^2(\{(n',k',z'), 0\leq k'\leq n',
z'\in S^{k'}_{x'}\})\text{\ \  dans \ \ } \ell^2(\{(n,k,z), 0\leq k\leq n,
z\in S^k_{x}\})$$ dont le coefficient vaut $e^{s(n-n')}$ si $I^{n',k',x'}_{z'}$ intervient dans la décomposition de $I^{n,k,x}_{z}$ dans le lemme précédent, et $0$ sinon. On a alors 
 $A\circ \Theta _{x',s}=\Theta _{x,s}$. 
Dans chaque ligne de $A$ il y a au plus 
$(q-1)(d(x,x')+1)+2$ coefficients non nuls, dans chaque colonne au plus $d(x,x')+2$ et 
ces coefficients ont une valeur absolue  inférieure ou égale à $e^{sd(x,x')}$ (car pour $n$ et $n'$ comme dans le lemme avec $I^{n',k',x'}_{z'}$ intervenant dans la décomposition de $I^{n,k,x}_{z}$, on a $|n-n'|\leq d(x,x')$). 
D'après le lemme~\ref{cauchy-schwarz} ci-dessous, une matrice telle que dans chaque ligne il y ait au plus 
$C_1$ coefficients non nuls, dans chaque colonne au plus $C_2$ et que
ses coefficients aient une norme inférieure ou égale à $C_3$ a
 une norme d'opérateur inférieure ou égale à
$\sqrt{C_1C_2}C_3$. 
Donc 
$$\|A\|\leq \sqrt{(q-1)(d(x,x')+1)+2}\sqrt{d(x,x')+2}e^{sd(x,x')}. $$
Pour tout $f\in \C^{(\Delta_{1})}$, on a 
$\|f\|_{H_{x,s}}\leq \|A\| \|f\|_{H_{x',s}}$
et donc $$\|\pi(g)f\|_{H_{x,s}}\leq  \|A\| \|\pi(g)f\|_{H_{x',s}} =\|A\| \|f\|_{H_{x,s}}.$$ On voit que  $\pi(g)$ est continu de $H_{x,s}$ dans lui-même, de norme inférieure ou égale à $\sqrt{(q-1)(\ell(g)+1)+2}\sqrt{\ell(g)+2}e^{s\ell(g)}$.

On termine maintenant la démonstration de la proposition~\ref{cas-des-arbres}.
On décide que  $H_{x,s}(\Delta_{1})$ est pair et $H_{x,s}(\Delta_{2})$  impair, et alors 
$$\Big(\big(H_{x,s}(\Delta_{1}) \oplus H_{x,s}(\Delta_{2})\big)[0,1],
(u_{t}+v_{t})_{t\in [0,1]}\Big)$$
appartient à $KK_{G,2s\ell +C}(\C,\C[0,1])$ pour une constante $C$ assez grande, et réalise l'homotopie entre $1$ et $\gamma$ : en effet 
 $(H_{x,s}(\Delta_{1}) \oplus H_{x,s}(\Delta_{2}),
u_{1}+ v_{1})$ est égal à $1$ dans $KK_{G,2s\ell +C}(\C,\C)$ par la
proposition 1.4.2 de~\cite{kkban} et 
$(H_{x,s}(\Delta_{1}) \oplus H_{x,s}(\Delta_{2}),
u_{0}+v_{0})$ est homotope à 
$(\ell^{2}(\Delta_{1}) \oplus \ell^{2}(\Delta_{2}),
u_{0}+ v_{0})$ qui représente $\gamma$. Comme d'habitude cette homotopie est réalisée en introduisant les normes intermédiaires $\sqrt{t\|.\|^{2}_{H_{x,s}}+(1-t)\|.\|^{2}_{\ell^{2}}}$, et en conservant l'opérateur $u_{0}+v_{0}$ pendant l'homotopie. \cqfd

\begin{lem}\label{cauchy-schwarz}
Soit $B$ une $C^{*}$-algèbre, $(E_{i})_{i\in I}$ des $B$-modules hilbertiens et $E=\bigoplus E_{i}$. 
Pour $i,j\in I$ on se donne $a_{ij}\in \L_{B}(E_{j},E_{i})$. 
Alors si $$\sup_{i\in I}\big(\sum_{j\in I}\|a_{ij}\|_{\L_{B}(E_{j},E_{i})}\big)\text{\ \ \ et\ \ \ } \sup_{j\in I}\big(\sum_{i\in I}\|a_{ij}\|_{\L_{B}(E_{j},E_{i})}\big)$$ sont finies, $A=(a_{ij})_{i,j\in I}$  appartient à $\L_{B}(E)$ et on a 
\begin{gather}\label{ineg-l1-linfini}\|A\|_{\L_{B}(E)}^{2}\leq \Big( \sup_{i\in I}\big(\sum_{j\in I}\|a_{ij}\|_{\L_{B}(E_{j},E_{i})}\big) \Big)  \Big( \sup_{j\in I}\big(\sum_{i\in I}\|a_{ij}\|_{\L_{B}(E_{j},E_{i})}\big) \Big).\end{gather}
En particulier si $A$, considérée comme une matrice par blocs,  possède au plus $C_{1}$ blocs non nuls dans chaque ligne et $C_{2}$ blocs non nuls
dans chaque colonne et si ces blocs ont une norme inférieure ou égale à $C_3$, 
alors $\|A\|_{\L_{B}(E)}\leq \sqrt{C_1C_2}C_3$. 
\end{lem}
\noindent{\bf Démonstration.}
 Il suffit de montrer l'inégalité (\ref{ineg-l1-linfini}) pour $I$ fini. 
Comme $$A\mapsto  \sup_{i\in I}\big(\sum_{j\in I}\|a_{ij}\|_{\L_{B}(E_{j},E_{i})}\big)$$ est une norme d'algèbre sur $\L_{B}(E)$, on a alors pour tout $n\in \N^{*}$, en notant $(A^{*}A)^{n}=(b_{ij})_{i,j\in I}$, 
\begin{gather*}\|A\|_{\L_{B}(E)}^{2n}=\|(A^{*}A)^{n}\|_{\L_{B}(E)}
\leq |I|  \Big( \sup_{i\in I}\big(\sum_{j\in I}\|b_{ij}\|_{\L_{B}(E_{j},E_{i})}\big) \Big)\\ \leq 
|I| \Big( \sup_{i\in I}\big(\sum_{j\in I}\|a_{ij}\|_{\L_{B}(E_{j},E_{i})}\big) \Big)^{n}  \Big( \sup_{j\in I}\big(\sum_{i\in I}\|a_{ij}\|_{\L_{B}(E_{j},E_{i})}\big)\Big)^{n}, \end{gather*}
 d'où l'inégalité (\ref{ineg-l1-linfini}) en faisant tendre $n$ vers l'infini. \cqfd

%Lorsque $A=\C$ on peut le montrer  par interpolation mais on le démontre très simplement en général par Cauchy-Schwarz : pour  $(x_{i})_{i\in I}\in E$ on a 
%$$\sum_{i}\|\sum_{j}a_{i,j}(x_{j})\|^{2}\leq
%\sum_{i}\big( \sum_{j}\|a_{i,j}\|_{\L_{B}(E_{j},E_{i})} \|x_{j}\|_{E_{j}}\big)^{2}$$ $$\underset{C.S.}{\leq}
%\sum_{i}\big( \sum_{j}\|a_{i,j}\|_{\L_{B}(E_{j},E_{i})} \big)
%\big( \sum_{j}\|a_{i,j}\|_{\L_{B}(E_{j},E_{i})} \|x_{j}\|_{E_{j}}^{2} \big)
%$$ $$\leq \Big( \sup_{i\in I}\big(\sum_{j\in I}\|a_{i,j}\|_{\L_{B}(E_{j},E_{i})}\big) \Big)\sum_{i} \sum_{j}\|a_{i,j}\|_{\L_{B}(E_{j},E_{i})} \|x_{j}\|_{E_{j}}^{2}$$ $$ \leq\Big( \sup_{i\in I}\big(\sum_{j\in I}\|a_{i,j}\|_{\L_{B}(E_{j},E_{i})}\big) \Big) \Big( \sup_{j\in I}\big(\sum_{i\in I}\|a_{i,j}\|_{\L_{B}(E_{j},E_{i})}\big) \Big)\Big(\sum_{j}  \|x_{j}\|_{E_{j}}^{2}\Big).$$

\noindent {\bf Remarque. } Les espaces $H_{x,s}(\Delta_{1})$ ressemblent aux  espaces de Hilbert introduits par Julg et Valette dans~\cite{julgvalette} pour construire l'homotopie entre $\gamma$ et $1$. Plus précisément on rappelle que pour $\lambda\in ]0,+\infty[$ Julg et Valette notent $H_{\lambda}$ le complété de $\C^{(\Delta_{1})}$ pour le produit scalaire $\s{\xi_{a}^{\lambda}, \xi_{b}^{\lambda}}=e^{-\lambda d(a,b)}$, où $\xi_{a}^{\lambda}$ est la fonction qui vaut $1$ en $a$ et $0$ ailleurs. Nous écrirons ici $H^{JV}_{\lambda}$ au lieu de $H_{\lambda}$ pour éviter les confusions. Alors pour tout $\lambda>0$ et tout $n\in \N$, et $f\in \C^{(\Delta_{1})}$ supporté sur $S^{n}_{x}$, on a 
\begin{gather*}\|f\|^{2}_{H^{JV}_{\lambda}}=\sum_{a,b\in S^{n}_{x}} e^{-\lambda d(a,b)}
\overline{f(a)} f(b)\\ =
(1-e^{-2\lambda})e^{-2\lambda n} \sum_{1\leq k\leq n}e^{2\lambda k} 
\sum_{z\in S^{k}_{x}} \Big| \sum_{a\in I^{n,k,x}_{z}} f(a)\Big|^{2}
+e^{-2\lambda n}\Big| \sum_{a\in S^{n}_{x}} f(a)\Big|^{2}\end{gather*}
et comme $S^{n}_{x}=I^{n,0,x}_{x}$ on voit que cette formule ressemble à la formule pour $\|f\|^{2}_{H_{x,s}(\Delta_{1})}$, bien que la présence des facteurs $(1-e^{-2\lambda})$ et $e^{2\lambda k}$ empêche de les identifier, quelles que soient les valeurs de $s$ et de $\lambda$. 
Plus précisément les restrictions de $\|.\|^{2}_{H^{JV}_{\lambda}}  $ et $\|.\|^{2}_{H_{x,s}(\Delta_{1})}   $ à l'espace des fonctions supportées sur $S_{x}^{n}$ s'expriment comme des sommes des carrés des mêmes formes linéaires 
$f\mapsto  \sum_{a\in I^{n,k,x}_z} f(a)$, mais avec des pondérations un peu différentes. 
D'autre part les opérateurs $u_{t}+v_{t}$ que nous utilisons sont quasiment les mêmes que ceux introduits par Julg et Valette. 
On notera cependant que, dans~\cite{julgvalette},  $\C^{(\Delta_{2})}$ est complété pour la norme de $\ell^{2}(\Delta_{2})$ alors que nous le complétons pour la norme de $H_{x,s}(\Delta_{2})$.

\section{Construction des espaces et des opérateurs}\label{construction-operateurs}

Ce paragraphe et les deux suivants sont consacrés à la démonstration du thérorème~\ref{hyperb-bon}. 

Soit $G$ un groupe localement compact agissant de fa\c con
  isométrique,   continue et propre sur un espace métrique $(X,d)$ qui est un bon espace hyperbolique discret.  
 Soit  $\delta\in \N^{*}$ tel que 
$(X,d)$ soit $\delta$-hyperbolique. 
Dans toute la suite de cet article on utilisera les notations suivantes : 
si $A$ et $B$ sont des parties finies de $X$, on note $d(A,B)=\min_{a\in A,b\in B}d(a,b)$ et $d_{\max}(A,B)=\max_{a\in A,b\in B}d(a,b)$. Si $A$ est un singleton $\{x\}$ on note $d(x,B)$ et $d_{\max}(x,B)$ au lieu de 
$d(\{x\},B)$ et $d_{\max}(\{x\},B)$. 
\label{ddmax}

Soit $x\in 
  X$ et $\ell$ la longueur sur $G$ définie par
  $\ell(g)=d(x,gx)$. 
  Le point $x$ sert d'origine dans la construction. Cependant  on utilisera en d'autres points que $x$ les constructions faites pour $x$ (notamment pour vérifier les propriétés d'équivariance). On utilisera aussi $x$ comme variable dans certains lemmes.

  Soit $\gamma\in KK_G(\C,\C)$ l'élément
  défini sous ces hypothèses par Kasparov et
  Skandalis~\cite{ks}. Soit  $s\in ]0,1]$. 
  On va construire, pour une certaine constante $C$ assez grande, une homotopie de $1$ à $\gamma$ dans   $E_{G,2s\ell+C}(\C,\C[0,1])$. 
Cette homotopie sera composée d'une partie difficile (partant de $1$), et d'une partie facile (aboutissant à $\gamma$) dont la construction est reléguée au paragraphe~\ref{construction-fin}. 
Le but de ce paragraphe est de construire des espaces vectoriels et des  opérateurs pour la partie difficile de l'homotopie. Dans le paragraphe~\ref{construction-normes} nous construirons les normes sur ces espaces.  

On commence par fixer un entier $N$ tel que 
 $$(H_{N}) : N \text{  est assez grand en fonction de }\delta.$$ %Plus précisément on prend $N\geq 100\delta$.  
 \label{defHn}

Plus précisément nous utiliserons un nombre fini de fois l'inégalité $N\geq C\de$, avec $C$ un entier. 

On note $K$ un entier tel que pour tout point $a$ de $X$ le nombre de points de $X$ à distance $1$ de $a$ soit inférieur ou égal à $K$. 
Dans toute la suite on notera $C(\delta)$, 
$C(\delta,K)$, $C(\delta,K,N)$ ... des constantes qui ne dépendent que des variables indiquées, mais varient d'une formule à l'autre. 

\vskip2mm
On note  $\Delta$ l'ensemble des parties de $X$ dont le diamètre est inférieur ou égal à $N$. 
  On note $p_{\max}$ le cardinal maximal d'une partie de $X$ de diamètre
$\leq N$. On a  $p_{\max}\leq C(\delta,K,N)$. 
Pour $p\in \{1,...,p_{\max}\}$ on note $\Delta_p$ l'ensemble
des parties de $X$ à $p$ éléments de diamètre
$\leq N$. 
On note $\C^{(X)}$ le $\C$-espace vectoriel formé des fonctions à
support fini de $X$ dans $\C$, dont la base canonique est
notée $(e_a)_{a\in X}$. On choisit une orientation de chaque simplexe associé à un élément de $\Delta\setminus \{\emptyset\}$, ce
qui permet de noter, pour tout $p\in \{1,\dots
,p_{\max}\}$, $\C^{(\Delta_p)}$ le sous-espace vectoriel de $\Lambda^p (
\C^{(X)})$ engendré par les $e_{a_1}\wedge \dots \wedge e_{a_p}$ pour 
$\{a_1,...,a_p\}\in \Delta _p$. Pour $S=\{a_1,...,a_p\}\in \Delta_{p}$ on notera aussi $e_{S}=\pm e_{a_1}\wedge \dots \wedge e_{a_p}$ (suivant l'orientation choisie pour ce simplexe). Ces choix d'orientation ont simplement pour but d'alléger  les notations. 
On a $\C^{(\Delta_1)}=\C^{(X)}$. On note
encore $\Delta_0=\{\emptyset\}$, 
$\C^{(\Delta_0)}=\C=\Lambda^0 (
\C^{(X)})$ et $e_{\emptyset}=1$.   
 
On note  $\C^X$ l'espace des fonctions de $X$ dans $\C$, qui est le dual
algébrique de $\C^{(X)}$. Pour tout $p\in \{0,\dots ,p_{\mathrm{max}}-1\}$
on note  $\partial :\C^{(\Delta_{p+1})}\vers \C^{(\Delta_p)}$ la
contraction à gauche par $(\dots
,1,1,\dots )\in \C^X$, 
c'est-à-dire que pour $\{a_{0},...,a_{p}\}\in \Delta_{p+1}$, on a 
$$\del(e_{a_0}\wedge \dots \wedge e_{a_p})=\sum_{i=0}^{p}(-1)^{i}
e_{a_0}\wedge \dots \wedge e_{a_{i-1}}\wedge e_{a_{i+1}}\wedge \dots\wedge e_{a_p}.$$
\label{def-del-page}
 En d'autres termes $\partial : 
\C^{(\Delta_1)}\vers \C^{(\Delta_0)}=\C$ est défini par $\partial 
(\sum _{a\in X}f(a)e_{a})=\sum 
_{a\in X}f(a)$ et $(\C^{(\Delta_1)}\vag{\partial }\C^{(\Delta_2)}\dots
\vag{\partial 
  }\C^{(\Delta_{p_{\mathrm{max}}})})$ est le complexe
d'homologie  
simpliciale. 

Le complexe suivant est exact : 
$$0\leftarrow \C^{(\Delta_0)}\vag{\partial }\C^{(\Delta_1)}\vag{\partial }\C^{(\Delta_2)}\dots \vag{\partial
  }\C^{(\Delta_{p_{\mathrm{max}}})}\leftarrow 0.$$

En fait nous allons construire, pour tout point $x\in X$,  un premier parametrix $H_{x}$ pour ce complexe (le même que  dans~\cite{kkban}), puis un deuxième parametrix $u_{x}$, puis un troisième parametrix $J_{x}$ qui est un mélange de $H_{x}$ et $u_{x}$, et c'est celui-là qui nous servira.  Ici comme dans la suite on emploie le terme ``parametrix'' plutôt que ``homotopie''   pour éviter que ce mot ait un double sens. Pour construire l'homotopie de $1$ à $\gamma$ on commencera par représenter l'image de $1$ dans $KK_{G,2s\ell+C}(\C,\C)$ (avec $C$ assez grand) par l'opérateur $\partial +J_{x}$ agissant sur le complété de 
$\bigoplus_{p=1}^{p_{\max}}\C^{(\Delta_p)}$ pour certaines normes de Hilbert $\|.\|_{\H_{x,s}}$ très compliquées construites au paragraphe suivant, puis on déformera  l'opérateur $\partial +J_{x}$ en le conjuguant par $e^{\tau \theta^{\flat}_{x}}$ où $\theta^{\flat}_{x}:\bigoplus_{p=1}^{p_{\max}}\C^{(\Delta_p)}\to \bigoplus_{p=1}^{p_{\max}}\C^{(\Delta_p)} $ est   défini par $\theta^{\flat}_{x}(e_{S})=\rho^{\flat}_{x}(S)e_{S}$ et où 
$\rho^{\flat}_{x}$ est une variante moyennée de la distance à $x$ : c'est la partie difficile de l'homotopie, qui est l'objet de ce paragraphe et du suivant. 
Pour $\tau$ assez grand l'opérateur $\partial +J_{x}$ conjugué par $e^{\tau \theta^{\flat}_{x}}$ est continu sur $\bigoplus_{p=1}^{p_{\max}}\ell^{2}(\Delta_p)$, 
ce qui permet de remplacer les normes compliquées par les normes $\ell^{2}$ et d'arriver ensuite à $\gamma$  : c'est la partie facile de l'homotopie, qui  est traitée   au paragraphe~\ref{construction-fin}. 

L'opérateur $J_{x}$ vérifiera les deux conditions suivantes, qui ne sont pas formulées de manière précise et servent seulement d'heuristique pour la construction de $J_{x}$.  Il existe une constante $C$ telle que 

\noindent (C1) $J_{x}$ rapproche de l'origine, plus précisément si $S_{0},S_{1},...,S_{n}$ est une suite sans répétition d'éléments de $\Delta $  telle que $e_{S_{i+1}}$ apparaît avec un coefficient non nul dans $\del (e_{S_{i}})$ ou dans $J_{x}(e_{S_{i}})$, la suite $S_{0},S_{1},...,S_{n}$ se rapproche de $x$ en restant à distance $\leq C$ de la réunion des  géodésiques entre $x$ et les points de $S_{0}$, 

\noindent (C2) $J_{x}$ est une intégrale  sur un paramètre $\alpha$ d'opérateurs $J_{x,\alpha}$ tels qu'il existe des parties $Y_{x,\alpha,S}$  (pour $S\in \Delta $) vérifiant les propriétés suivantes : 
\begin{itemize}
\item $Y_{x,\alpha,S}$ est une partie finie de $X$, ne dépendant que de $x,\alpha,S$,  de cardinal $\leq C$, et contenant $x$ et $S$
\item tous les points de $Y_{x,\alpha,S}$ sont à distance $\leq C$  de la réunion des géodésiques entre $x$ et les points de $S$, 
\item les distances entre les points de $Y_{x,\alpha,S}$ sont déterminées à $ C$ près par $x,\alpha,S$ (ce qui fait que le nombre de possibilités pour l'ensemble des distances entre les points de $ Y_{x,\alpha,S}$ est borné par une constante) 
\item pour $T\in \Delta $, le coefficient de $e_{T}$ dans $J_{x,\alpha}(e_{S})$ est nul si $T$ n'est pas inclus dans $Y_{x,\alpha,S}$ et il ne dépend 
que de la connaissance des distances entre les points de $Y_{x,\alpha,S}$, c'est-à-dire, plus précisément :  si on se donne  \begin{gather*}S=\{a_{1},...,a_{p-1}\},  \ T=\{b_{1},...,b_{p}\}\subset Y_{x,\alpha,S}, \\
 x', \ S'=\{a'_{1},...,a'_{p-1}\}\text{ \  et \ } T'=\{b'_{1},...,b'_{p}\}\subset Y_{x',\alpha,S'}\end{gather*} tels qu'il existe une isométrie de $Y_{x,\alpha,S}$ dans $Y_{x',\alpha,S'}$ qui envoie $$x,a_{1},...,a_{p-1},b_{1},...,b_{p}\text{ \  sur \  }x',a'_{1},...,a'_{p-1},b'_{1},...,b'_{p},$$ alors le coefficient de $e_{b'_{1}}\wedge ... \wedge e_{b'_{p}}$ dans $J_{x',\alpha}(e_{a'_{1}}\wedge ... \wedge e_{a'_{p-1}})$  est égal au coefficient de $e_{b_{1}}\wedge ... \wedge e_{b_{p}}$ dans $J_{x,\alpha}(e_{a_{1}}\wedge ...\wedge e_{a_{p-1}})$. 
\end{itemize}

Dans la condition (C2) la donnée des points de $S$ sert à lever l'ambiguïté de signe  sur $e_{S}$ et de même pour $T,S',T'$.

Ces propriétés (C1) et (C2) serviront pour montrer la continuité de $J_{x}$.  D'après la formule (\ref{formule-norme}) ci-dessous, 
pour $f\in \C^{(\Delta_{p})}$, 
on aura \begin{gather}\label{explic-form-norme}\|f\|^{2}_{\H_{x,s}}=\sum_{Z}\kappa_{Z}|\xi_{Z}(f)|^2,\end{gather} où 
 la somme porte sur certaines classes d'équivalences $Z$ d'uplets 
 $(a_{1},...,a_{p},Y)$, avec $\{a_{1},...,a_{p}\}\in \Delta$ et $Y$ une partie finie de $X$ dont tous les points sont à distance $\leq C$ de la réunion des géodésiques entre $x$ et les points de $\{a_{1},...,a_{p}\}$
 (la somme sur $Z$ est infinie et le cardinal de $Y$ ne dépend que de $Z$ mais n'est pas borné indépendamment de  $Z$). La relation d'équivalence est telle que $Z$  détermine les distances entre les points de $\{x\}\cup 
 \{a_{1},...,a_{p}\}
 \cup Y$. Enfin  $\kappa_{Z}\in \R_{+}^{*}$ est une pondération et $$\xi_{Z}(f)=\sum_{(a_{1},...,a_{p},Y)\in Z}f(a_{1},...,a_{p})$$  
 où $f(a_{1},...,a_{p})$ désigne le coefficient de $e_{a_{1}}\wedge ...\wedge e_{a_{p}}$ dans  $f$. De plus si $Z$ est comme ci-dessus et $J_{x,\alpha}$ est comme dans (C2), ${}^{t}J_{x,\alpha}(\xi_{Z})$ sera une combinaison finie de formes linéaires $\xi_{Z'}$, avec $Z'$ une classe d'équivalence de $(a'_{1},...,a'_{p-1},Y')$ tels qu'il existe $(a_{1},...,a_{p},Y)\in Z$ vérifiant 
 \begin{itemize}
 \item $e_{a_{1}}\wedge ...\wedge e_{a_{p}}$ peut apparaître dans 
 $J_{x,\alpha}(e_{a'_{1}}\wedge ...\wedge e_{a'_{p-1}})$, en particulier 
 $a_{1},  ... , a_{p}$ sont à distance $\leq C$ de la réunion des géodésiques entre $x$ et les points de $\{a'_{1},...,a'_{p-1}\}$, 
 \item $Y'$ contient   $Y\cup  \{a_{1},...,a_{p}\}  \cup Y_{x,\alpha,\{a'_{1},...,a'_{p-1}\}}$. 
 \end{itemize}
   Donc $\|.\|^{2}_{\H_{x,s}}$ sera choisie de telle sorte que 
 si  $\xi_{Z}$ apparaît dans la formule pour $\|.\|^{2}_{\H_{x,s}}$, 
  ces nouvelles formes linéaires $\xi_{Z'}$ apparaissent aussi dans la formule pour $\|.\|^{2}_{\H_{x,s}}$. 
La condition (C2) fournit une constante $C'$ telle que   ${}^{t}J_{x,\alpha}(\xi_{Z})$ soit une combinaison d'au plus $C'$ formes linéaires $\xi_{Z'}$, ce qui permettra d'appliquer Cauchy-Schwarz
pour majorer $|\xi_{Z}(J_{x,\alpha}(f))|^{2}$ par une combinaison des 
$|\xi_{Z'}(f)|^{2}$.  Comme $J_{x}$ sera une intégrale  de tels opérateurs $J_{x,\alpha}$ on montrera, en utilisant encore Cauchy-Schwarz, que $\|J_{x} f\|_{\H_{x,s}}\leq C\|f\|_{\H_{x,s}}$ pour une certaine constante $C$.  
La formule  (\ref{formule-norme})   pour $\|.\|^{2}_{\H_{x,s}}$ que nous donnerons plus loin est plus ou moins déterminée par la condition que si une forme linéaire $\xi_{Z}$ figure dans la formule pour $\|.\|^{2}_{\H_{x,s}}$, les formes linéaires $\xi_{Z'}$ servant à décomposer ${}^{t}J_{x,\alpha}(\xi_{Z})$ doivent  apparaître dans $\|.\|^{2}_{\H_{x,s}}$
(ainsi que la même condition pour les autres opérateurs comme $e^{\tau \theta^{\flat}_{x}}\del e^{-\tau \theta^{\flat}_{x}}$ et $e^{\tau \theta^{\flat}_{x}}J_{x}e^{-\tau \theta^{\flat}_{x}}$). La condition  (C1) garantit que ce procédé ne diverge pas, en particulier que pour $f\in \C^{(\Delta_p)}$, $\sum_{Z}\kappa_{Z}|\xi_{Z}(f)|^2<\infty$. 

Nous allons voir que $H_{x}$ satisfait (C1) mais pas (C2). Inversement $u_{x}$ satisfera (C2) mais pas (C1). Heureusement $J_{x}$ vérifiera (C1) et (C2). 

\subsection{Construction du premier parametrix $H_{x}$}

Dans toute la suite   nous aurons besoin de la notation suivante. 
\begin{defi}
Si $(Y,d)$ est un espace métrique, $\epsilon \in \R_{+}$, et $x,y\in Y$, on note $\epsilon\text{-}\geod(x,y)$ l'ensemble des points $z$ de $Y$ tels que $d(x,z)+d(z,y)\leq d(x,y)+\epsilon$. 
On note $\geod(x,y)$ au lieu de $0\text{-}\geod(x,y)$.
\end{defi}
Donc $\geod(x,y)$ désigne simplement l'ensemble des points $z$ de $Y$ tels que $d(x,z)+d(z,y)=d(x,y)$. 

Les deux lemmes suivants sont valables pour tout espace métrique car ils n'utilisent que l'inégalité triangulaire. 

\begin{lem}\label{xx'yy'zz'}
Soient $(Y,d)$ un espace métrique, $x,x',y,y',z,z'\in Y$ et $\alpha\in \R_{+}$ tels que $z\in \alpha\tg(x,y)$. Alors $$z'\in \big(\alpha+2d(x,x')+2d(y,y')+2d(z,z')\big)\tg(x',y').$$
\end{lem}
\noindent{\bf Démonstration.}
Cela résulte des inégalités évidentes
$$d(x',z')\leq d(x,z)+d(x,x')+d(z,z'),\ \ \  d(z',y')\leq d(z,y)+d(y,y')+d(z,z')$$ et 
$d(x',y')\geq d(x,y)-d(x,x')-d(y,y')$. \cqfd

\begin{lem}\label{geod-comp-xabc}
Soient $(Y,d)$ un espace métrique, $x,a,b,c\in Y$ et $\alpha,\beta\in \R_{+}$ tels que $b\in \alpha\tg(x,a)$ et $c\in \beta\tg(x,b)$. Alors

\noindent a) $c\in (\alpha+\beta)\tg(x,a)$, 

\noindent b) $b\in (\alpha+\beta)\tg(a,c)$, 

\noindent c) si $d(b,c)\geq \alpha+\beta$, $d(a,c)\geq d(a,b)$. 
\end{lem}

\ifx\JPicScale\undefined\def\JPicScale{1}\fi
\unitlength \JPicScale mm
\begin{picture}(90,28)(20,35)
\linethickness{0.3mm}
\put(30,40){\line(1,0){80}}
\linethickness{0.3mm}
\multiput(92.5,52.5)(0.17,-0.12){104}{\line(1,0){0.17}}
\linethickness{0.3mm}
\multiput(30,40)(0.6,0.12){104}{\line(1,0){0.6}}
\linethickness{0.3mm}
\multiput(30,40)(0.27,0.12){148}{\line(1,0){0.27}}
\linethickness{0.3mm}
\multiput(70,57.5)(0.54,-0.12){42}{\line(1,0){0.54}}
\put(30,45){\makebox(0,0)[cc]{$x$}}

\put(110,45){\makebox(0,0)[cc]{$a$}}

\put(95,55){\makebox(0,0)[cc]{$b$}}

\put(70,61){\makebox(0,0)[cc]{$c$}}

\put(70,52.5){\makebox(0,0)[cc]{$\beta$}}

\put(89,45){\makebox(0,0)[cc]{$\alpha$}}

\end{picture}

\noindent{\bf Démonstration.}
 On a $$    d(a,b)+d(b,c)+d(c,x)\leq d(a,b)+d(b,x)+\beta\leq d(a,x)+\alpha+\beta,  $$ d'où l'on déduit a) et b) à l'aide des inégalités triangulaires pour les triangles $abc$ et $acx$ respectivement. Enfin  
 c) résulte immédiatement de b).   \cqfd

Voici maintenant un lemme qui utilise le fait que $(X,d)$ est $\de$-hyperbolique. 
\begin{lem}
Soient $x,a,b,c\in X$. 

\noindent a) Soit $\beta\in \N$ tel que $b\in \beta\tg(a,c)$. Alors 
$$d(x,b)\leq \max\big(d(x,a)-d(a,b),d(x,c)-d(c,b)\big)+\beta+\de. \ \ \ \   \ \ \ (H_{\delta}^{\beta}(x,a,b,c))$$

\noindent b) En particulier si  $b$ appartient à $ \geod(a,c)$ on a 
$$d(x,b)\leq \max\big(d(x,a)-d(a,b),d(x,c)-d(c,b)\big)+\de. \ \ \ \  \ \ \  \ \ \ (H_{\delta}^{0}(x,a,b,c))$$
\end{lem}

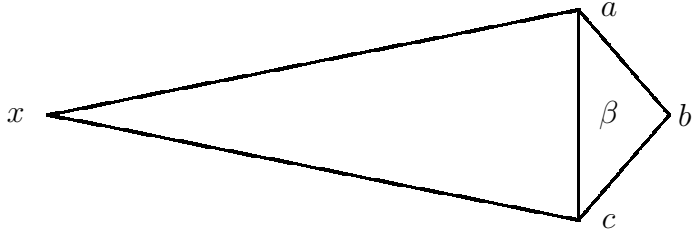
\ifx\JPicScale\undefined\def\JPicScale{1}\fi
\unitlength \JPicScale mm
\begin{picture}(100,35)(20,13)
\linethickness{0.3mm}
\multiput(40,30)(0.6,0.12){117}{\line(1,0){0.6}}
\linethickness{0.3mm}
\put(110,16){\line(0,1){28}}
\linethickness{0.3mm}
\multiput(40,30)(0.6,-0.12){117}{\line(1,0){0.6}}
\linethickness{0.3mm}
\multiput(110,44)(0.12,-0.14){100}{\line(0,-1){0.14}}
\linethickness{0.3mm}
\multiput(110,16)(0.12,0.14){100}{\line(0,1){0.14}}
\put(36,30){\makebox(0,0)[cc]{$x$}}

\put(114,44){\makebox(0,0)[cc]{$a$}}

\put(124,30){\makebox(0,0)[cc]{$b$}}

\put(114,16){\makebox(0,0)[cc]{$c$}}

\put(114,30){\makebox(0,0)[cc]{$\beta$}}

\end{picture}

\label{Hdelta-beta}
\noindent{\bf Démonstration.}
Le a)  est une conséquence directe de la propriété d'hyperbolicité $(H_{\delta}(x,a,b,c))$ de la définition~\ref{defi-hyperb-Hxyzt} et b) est le cas particulier de a) où $\beta=0$. \cqfd

Nous devons commencer par rappeler certains lemmes de~\cite{ks} et~\cite{kkban}. 
Comme nous avons affaire à des espaces hyperboliques et géodésiques, et non pas seulement boliques et faiblement géodésiques, les énoncés  et les démonstrations  des lemmes se simplifient, et nous repartirons donc de zéro. En plus nous éliminerons les paramètres $k$ et $t$ de~\cite{kkban}.

On rappelle que  $\Delta$ est l'ensemble des parties de $X$ dont le diamètre est inférieur ou égal à $N$. Pour $S\in \Delta\setminus \{\emptyset\}$ on note 
$$U_{S}=
\cap _{a\in S} B(a,N)=\{z\in X, \{z\}\cup S\in \Delta
\} .$$

\begin{lem}  \label{72}
(cf le LEMME 6.2 de~\cite{ks}, et le lemme 2.1.4 de~\cite{kkban})
Soient $x\in X$ et $S\in \Delta \setminus\{\emptyset\}$. 
Pour $z,z'\in U_{S}$ on a 
$$d(z,z')\leq (d(x,z)-d(x,U_{S}))+(d(x,z')-d(x,U_{S}))+4\delta.$$ En particulier le diamètre de $\{z\in U_S,
d(x,z)\leq d(x,U_S)+\delta \}$ est inférieur ou égal à $6\delta$.   
\end{lem}
Grâce à $(H_{N})$ on peut supposer $N\geq 6\de$, ce que l'on fait. Le lemme~\ref{72} implique alors que  le diamètre de $\{z\in U_S,
d(x,z)\leq d(x,U_S)+\delta \}$ est inférieur ou égal à $N$. 

\noindent{\bf Démonstration.} On commence par un résultat trivial. 
\begin{lem}\label{stab-boules}
Soit $x\in X$, $r\in \N$, $z_{1},z_{3}\in B(x,r)$ et $z_{2}\in \geod(z_{1},z_{3})$ avec $d(z_{1},z_{2})\geq \delta$ et $d(z_{2},z_{3})\geq \delta$. Alors $z_{2}\in B(x,r)$. 
\end{lem}
\noindent{\bf Démonstration. }
Par $(H_{\delta}^{0}(x,z_{1},z_{2},z_{3}))$ on a 
$$d(x,z_{2})\leq \max(d(x,z_{1})-d(z_{1},z_{2}),d(x,z_{3})-d(z_{2},z_{3}))+\de.$$ \cqfd

\noindent{\bf Fin de la démonstration du  lemme~\ref{72}.} 
Soient $z_{1},z_{3}\in U_{S}$ avec $$d(z_{1},z_{3})\geq (d(x,z_{1})-d(x,U_{S}))+(d(x,z_{3})-d(x,U_{S}))+4\delta.$$ Nous allons aboutir à une contradiction. 
Il existe $z_{2}\in \geod(z_{1},z_{3})$ tel que $$d(z_{1},z_{2})=(d(x,z_{1})-d(x,U_{S}))+2\delta. $$
Alors on a 
$d(z_{3},z_{2})\geq(d(x,z_{3})-d(x,U_{S}))+2\delta$. 
D'après le lemme~\ref{stab-boules}, on a $z_{2}\in  U_S$. Mais 
par $(H_{\delta}^{0}(x,z_{1},z_{2},z_{3}))$ on a 
$d(x,z_{2})\leq d(x,U_{S})-\delta$, d'où une contradiction. \cqfd

\begin{lem} \label{73}
(cf le lemme 6.3 de~\cite{ks} et le lemme 2.1.5 de~\cite{kkban})
Soient $x\in X$ et $S,T\in \Delta \setminus\{\emptyset\}$. Supposons que tout point $a$ dans la
différence symétrique de $S$ et $T$ vérifie $d(x,a)\leq d(x,U_S)+N-5\delta $. Alors 
 
 \noindent a)
  $d(x,U_T)=d(x,U_S)$,

\noindent  b)  pour tout $b$ dans la
  différence symétrique de $U_S$ et $U_T$, $$d(x,b)>d(x,U_S)+\delta.$$
 \end{lem}
\noindent{\bf Démonstration.}
Nous suivons la preuve du lemme 6.3 de~\cite{ks}. Par récurrence sur le cardinal de la différence symétrique de $S$ et $T$ on peut supposer que la différence symétrique de $S$ et $T$ est un singleton $\{a\}$. 

Supposons d'abord $a\in T$. Alors $U_{T}\subset U_{S}$. 
Comme $T=S\cup\{a\}$ appartient à $\Delta$, 
$a\in U_{S}$. Par le lemme~\ref{72}, $d(a,\{z\in U_{S},d(x,z)\leq d(x,U_{S})+\delta\})\leq N$. Donc
$$\{z\in U_{S},d(x,z)\leq d(x,U_{S})+\delta\} \subset
U_{T}\subset U_{S}$$
et les  assertions a) et b) du lemme en résultent. 

Supposons maintenant $a\in S$. Alors $U_{S}\subset U_{T}$ et comme $S=T\cup\{a\}$ appartient à $\Delta$, 
$a\in U_{T}$. 
Il suffit de montrer a) car en échangeant les rôles de $S$ et $T$, b) en résulte. 
Raisonnons par l'absurde et supposons $d(x,U_{T})<d(x,U_{S})$. Soit $b\in U_{T}$ tel que $d(x,b)=d(x,U_{T})$. On a alors $b\not\in U_{S}$ donc $d(a,b)>N$. On rappelle que   $N\geq 6\delta$. Soit $c\in \geod(a,b)$, $d(a,c)=N-3\delta$. Par le lemme~\ref{stab-boules},
comme $a,b\in U_{T}$, $d(a,c)\geq\de, d(b,c)\geq\de$ et $c\in \geod(a,b)$, 
 on a $c\in U_{T}$. Mais $d(a,c)\leq N$ donc $c\in U_{S}$. Or $(H_{\delta}^{0}(x,a,c,b))$ donne 
$$d(x,c)\leq\max(d(x,a)-d(a,c), d(x,b)-d(b,c))+\de<d(x,U_{S}),$$
ce qui est contradictoire. \cqfd

Pour $x\in X$, $S\in \Delta\setminus\{\emptyset\}$, on définit
$$A_{S,x}=\{z\in U_S, d(x,z)\leq d(x,U_S)+\delta \}$$ et
pour $r\in \N$, $$Y_{S,x,r}=
\{z\in U_{S},\exists y\in \delta\text{-}\geod(x,z),d(x,y)\leq r,d(y,z)=d(y,U_{S})\}.$$
\label{def-Asx}
\ifx\JPicScale\undefined\def\JPicScale{1}\fi
\unitlength \JPicScale mm
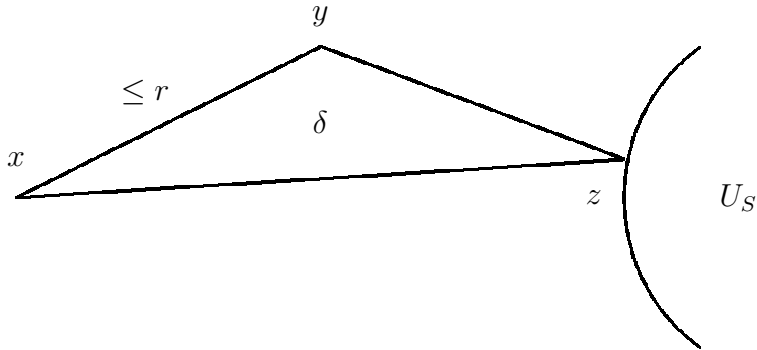
\begin{picture}(165,50)(40,18)
\linethickness{0.3mm}
\multiput(139.6,59.7)(0.13,0.1){3}{\line(1,0){0.13}}
\multiput(139.21,59.39)(0.13,0.1){3}{\line(1,0){0.13}}
\multiput(138.83,59.07)(0.13,0.11){3}{\line(1,0){0.13}}
\multiput(138.45,58.74)(0.13,0.11){3}{\line(1,0){0.13}}
\multiput(138.08,58.41)(0.12,0.11){3}{\line(1,0){0.12}}
\multiput(137.72,58.07)(0.12,0.11){3}{\line(1,0){0.12}}
\multiput(137.36,57.72)(0.12,0.12){3}{\line(1,0){0.12}}
\multiput(137.01,57.36)(0.12,0.12){3}{\line(0,1){0.12}}
\multiput(136.67,57)(0.11,0.12){3}{\line(0,1){0.12}}
\multiput(136.34,56.63)(0.11,0.12){3}{\line(0,1){0.12}}
\multiput(136.01,56.26)(0.11,0.13){3}{\line(0,1){0.13}}
\multiput(135.69,55.87)(0.11,0.13){3}{\line(0,1){0.13}}
\multiput(135.37,55.49)(0.1,0.13){3}{\line(0,1){0.13}}
\multiput(135.07,55.09)(0.1,0.13){3}{\line(0,1){0.13}}
\multiput(134.77,54.69)(0.15,0.2){2}{\line(0,1){0.2}}
\multiput(134.48,54.29)(0.14,0.2){2}{\line(0,1){0.2}}
\multiput(134.2,53.87)(0.14,0.21){2}{\line(0,1){0.21}}
\multiput(133.93,53.46)(0.14,0.21){2}{\line(0,1){0.21}}
\multiput(133.67,53.03)(0.13,0.21){2}{\line(0,1){0.21}}
\multiput(133.41,52.61)(0.13,0.21){2}{\line(0,1){0.21}}
\multiput(133.16,52.17)(0.12,0.22){2}{\line(0,1){0.22}}
\multiput(132.92,51.73)(0.12,0.22){2}{\line(0,1){0.22}}
\multiput(132.7,51.29)(0.11,0.22){2}{\line(0,1){0.22}}
\multiput(132.47,50.84)(0.11,0.22){2}{\line(0,1){0.22}}
\multiput(132.26,50.39)(0.11,0.23){2}{\line(0,1){0.23}}
\multiput(132.06,49.94)(0.1,0.23){2}{\line(0,1){0.23}}
\multiput(131.87,49.48)(0.1,0.23){2}{\line(0,1){0.23}}
\multiput(131.68,49.02)(0.09,0.23){2}{\line(0,1){0.23}}
\multiput(131.51,48.55)(0.18,0.47){1}{\line(0,1){0.47}}
\multiput(131.34,48.08)(0.17,0.47){1}{\line(0,1){0.47}}
\multiput(131.18,47.6)(0.16,0.47){1}{\line(0,1){0.47}}
\multiput(131.04,47.13)(0.15,0.48){1}{\line(0,1){0.48}}
\multiput(130.9,46.65)(0.14,0.48){1}{\line(0,1){0.48}}
\multiput(130.77,46.17)(0.13,0.48){1}{\line(0,1){0.48}}
\multiput(130.65,45.68)(0.12,0.48){1}{\line(0,1){0.48}}
\multiput(130.55,45.2)(0.11,0.49){1}{\line(0,1){0.49}}
\multiput(130.45,44.71)(0.1,0.49){1}{\line(0,1){0.49}}
\multiput(130.36,44.22)(0.09,0.49){1}{\line(0,1){0.49}}
\multiput(130.28,43.73)(0.08,0.49){1}{\line(0,1){0.49}}
\multiput(130.21,43.23)(0.07,0.49){1}{\line(0,1){0.49}}
\multiput(130.15,42.74)(0.06,0.49){1}{\line(0,1){0.49}}
\multiput(130.1,42.24)(0.05,0.5){1}{\line(0,1){0.5}}
\multiput(130.06,41.74)(0.04,0.5){1}{\line(0,1){0.5}}
\multiput(130.03,41.25)(0.03,0.5){1}{\line(0,1){0.5}}
\multiput(130.01,40.75)(0.02,0.5){1}{\line(0,1){0.5}}
\multiput(130,40.25)(0.01,0.5){1}{\line(0,1){0.5}}
\put(130,39.75){\line(0,1){0.5}}
\multiput(130,39.75)(0.01,-0.5){1}{\line(0,-1){0.5}}
\multiput(130.01,39.25)(0.02,-0.5){1}{\line(0,-1){0.5}}
\multiput(130.03,38.75)(0.03,-0.5){1}{\line(0,-1){0.5}}
\multiput(130.06,38.26)(0.04,-0.5){1}{\line(0,-1){0.5}}
\multiput(130.1,37.76)(0.05,-0.5){1}{\line(0,-1){0.5}}
\multiput(130.15,37.26)(0.06,-0.49){1}{\line(0,-1){0.49}}
\multiput(130.21,36.77)(0.07,-0.49){1}{\line(0,-1){0.49}}
\multiput(130.28,36.27)(0.08,-0.49){1}{\line(0,-1){0.49}}
\multiput(130.36,35.78)(0.09,-0.49){1}{\line(0,-1){0.49}}
\multiput(130.45,35.29)(0.1,-0.49){1}{\line(0,-1){0.49}}
\multiput(130.55,34.8)(0.11,-0.49){1}{\line(0,-1){0.49}}
\multiput(130.65,34.32)(0.12,-0.48){1}{\line(0,-1){0.48}}
\multiput(130.77,33.83)(0.13,-0.48){1}{\line(0,-1){0.48}}
\multiput(130.9,33.35)(0.14,-0.48){1}{\line(0,-1){0.48}}
\multiput(131.04,32.87)(0.15,-0.48){1}{\line(0,-1){0.48}}
\multiput(131.18,32.4)(0.16,-0.47){1}{\line(0,-1){0.47}}
\multiput(131.34,31.92)(0.17,-0.47){1}{\line(0,-1){0.47}}
\multiput(131.51,31.45)(0.18,-0.47){1}{\line(0,-1){0.47}}
\multiput(131.68,30.98)(0.09,-0.23){2}{\line(0,-1){0.23}}
\multiput(131.87,30.52)(0.1,-0.23){2}{\line(0,-1){0.23}}
\multiput(132.06,30.06)(0.1,-0.23){2}{\line(0,-1){0.23}}
\multiput(132.26,29.61)(0.11,-0.23){2}{\line(0,-1){0.23}}
\multiput(132.47,29.16)(0.11,-0.22){2}{\line(0,-1){0.22}}
\multiput(132.7,28.71)(0.11,-0.22){2}{\line(0,-1){0.22}}
\multiput(132.92,28.27)(0.12,-0.22){2}{\line(0,-1){0.22}}
\multiput(133.16,27.83)(0.12,-0.22){2}{\line(0,-1){0.22}}
\multiput(133.41,27.39)(0.13,-0.21){2}{\line(0,-1){0.21}}
\multiput(133.67,26.97)(0.13,-0.21){2}{\line(0,-1){0.21}}
\multiput(133.93,26.54)(0.14,-0.21){2}{\line(0,-1){0.21}}
\multiput(134.2,26.13)(0.14,-0.21){2}{\line(0,-1){0.21}}
\multiput(134.48,25.71)(0.14,-0.2){2}{\line(0,-1){0.2}}
\multiput(134.77,25.31)(0.15,-0.2){2}{\line(0,-1){0.2}}
\multiput(135.07,24.91)(0.1,-0.13){3}{\line(0,-1){0.13}}
\multiput(135.37,24.51)(0.1,-0.13){3}{\line(0,-1){0.13}}
\multiput(135.69,24.13)(0.11,-0.13){3}{\line(0,-1){0.13}}
\multiput(136.01,23.74)(0.11,-0.13){3}{\line(0,-1){0.13}}
\multiput(136.34,23.37)(0.11,-0.12){3}{\line(0,-1){0.12}}
\multiput(136.67,23)(0.11,-0.12){3}{\line(0,-1){0.12}}
\multiput(137.01,22.64)(0.12,-0.12){3}{\line(0,-1){0.12}}
\multiput(137.36,22.28)(0.12,-0.12){3}{\line(1,0){0.12}}
\multiput(137.72,21.93)(0.12,-0.11){3}{\line(1,0){0.12}}
\multiput(138.08,21.59)(0.12,-0.11){3}{\line(1,0){0.12}}
\multiput(138.45,21.26)(0.13,-0.11){3}{\line(1,0){0.13}}
\multiput(138.83,20.93)(0.13,-0.11){3}{\line(1,0){0.13}}
\multiput(139.21,20.61)(0.13,-0.1){3}{\line(1,0){0.13}}
\multiput(139.6,20.3)(0.13,-0.1){3}{\line(1,0){0.13}}

\linethickness{0.3mm}
\multiput(50,40)(0.24,0.12){167}{\line(1,0){0.24}}
\linethickness{0.3mm}
\multiput(90,60)(0.32,-0.12){125}{\line(1,0){0.32}}
\linethickness{0.3mm}
\multiput(50,40)(1.9,0.12){42}{\line(1,0){1.9}}
\put(90,64){\makebox(0,0)[cc]{$y$}}

\put(50,45){\makebox(0,0)[cc]{$x$}}

\put(67,54){\makebox(0,0)[cc]{$\leq r$}}

\put(90,50){\makebox(0,0)[cc]{$\delta$}}

\put(126,40){\makebox(0,0)[cc]{$z$}}

\put(145,40){\makebox(0,0)[cc]{$U_S$}}

\end{picture}

\noindent 
 Il est clair que $r\mapsto Y_{S,x,r}$ est une application croissante, c'est-à-dire que si $r\leq r'$ on a $Y_{S,x,r}\subset Y_{S,x,r'}$. De plus  $Y_{S,x,0}$ est non vide. 
Si $r\leq d(x,U_{S})$ on a 
\begin{gather}\label{def-eq-YSxr-22oct}Y_{S,x,r}=
\{z\in U_{S},\exists y\in \delta\text{-}\geod(x,z),d(x,y)= r,d(y,z)=d(y,U_{S})\}.\end{gather} 
En effet pour $z\in Y_{S,x,r}$ et $y\in \de\tg(x,z)$ vérifiant $d(x,y)\leq r$ et $d(y,z)=d(y,U_{S})$, $\{d(x,y'),y'\in \geod(y,z)\}$ est un intervalle contenant 
$d(x,y)$ et $d(x,z)$,  et comme 
$d(x,y)\leq r$ et $d(x,z)\geq d(x,U_{S})\geq r$,  il existe $y'\in \geod(y,z)$ tel que $d(x,y')=r$ et alors $d(y',z)=d(y',U_{S})$ et $y'\in \de\tg(x,z)$ par le a) du lemme~\ref{geod-comp-xabc}. 
%
%remplacer $y$ par un point de $\geod(y,z)$ à distance $r-d(x,y)$ de $y$. 

Plus tard on aura besoin des conventions : $A_{\emptyset,x}=A_{\{x\},x}=B(x,\de)$ et $Y_{\emptyset,x,r}=Y_{\{x\},x,r}$, qui est égal à $\{x\}$ pour $r=0$. 

Cet ensemble  $Y_{S,x,r}$ est très proche de celui défini par Kasparov et Skandalis dans~\cite{ks}, il sert pour l'astuce des ensembles emboîtés. 

Comme le diamètre de $U_{S}$ est inférieur ou égal à $2N$, il est clair que $U_{S}$ donc aussi $A_{S,x}$ et $Y_{S,x,r}$, ont des cardinaux bornés par $C(\delta,K,N)$. En fait le lemme~\ref{72} montre que le diamètre de $A_{S,x}$ est inférieur ou égal à $6\delta$, donc son cardinal est borné par $C(\delta,K)$. 

\begin{lem}\label{majoration-r}
Si $r\leq d(x,U_{S})-N$ ou si $r=0$, on a $Y_{S,x,r}\subset A_{S,x}$  et
$$Y_{S,x,r}=
\{z\in A_{S,x},\exists y\in \delta\text{-}\geod(x,z),d(x,y)=r,d(y,z)=d(y,A_{S,x})\}.$$ 
\end{lem}
\noindent{\bf Démonstration.}
Si $r=0$, c'est évident. Supposons donc $r\leq d(x,U_{S})-N$.
Soit $z\in Y_{S,x,r}$. Grâce à (\ref{def-eq-YSxr-22oct}) on a  $z\in U_{S}$  et il existe $y\in \delta\text{-}\geod(x,z)$, tel que $d(x,y)= r$ et 
$d(y,z)=d(y,U_{S})$. On a $d(x,z)\geq d(x,U_{S})\geq r+N$ donc $d(y,z)\geq N$. 

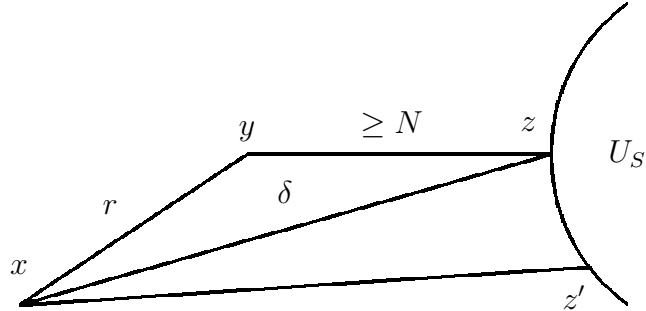
\ifx\JPicScale\undefined\def\JPicScale{1}\fi
\unitlength \JPicScale mm
\begin{picture}(140,55)(50,10)
\linethickness{0.3mm}
\multiput(149.6,59.7)(0.13,0.1){3}{\line(1,0){0.13}}
\multiput(149.21,59.39)(0.13,0.1){3}{\line(1,0){0.13}}
\multiput(148.83,59.07)(0.13,0.11){3}{\line(1,0){0.13}}
\multiput(148.45,58.74)(0.13,0.11){3}{\line(1,0){0.13}}
\multiput(148.08,58.41)(0.12,0.11){3}{\line(1,0){0.12}}
\multiput(147.72,58.07)(0.12,0.11){3}{\line(1,0){0.12}}
\multiput(147.36,57.72)(0.12,0.12){3}{\line(1,0){0.12}}
\multiput(147.01,57.36)(0.12,0.12){3}{\line(0,1){0.12}}
\multiput(146.67,57)(0.11,0.12){3}{\line(0,1){0.12}}
\multiput(146.34,56.63)(0.11,0.12){3}{\line(0,1){0.12}}
\multiput(146.01,56.26)(0.11,0.13){3}{\line(0,1){0.13}}
\multiput(145.69,55.87)(0.11,0.13){3}{\line(0,1){0.13}}
\multiput(145.37,55.49)(0.1,0.13){3}{\line(0,1){0.13}}
\multiput(145.07,55.09)(0.1,0.13){3}{\line(0,1){0.13}}
\multiput(144.77,54.69)(0.15,0.2){2}{\line(0,1){0.2}}
\multiput(144.48,54.29)(0.14,0.2){2}{\line(0,1){0.2}}
\multiput(144.2,53.87)(0.14,0.21){2}{\line(0,1){0.21}}
\multiput(143.93,53.46)(0.14,0.21){2}{\line(0,1){0.21}}
\multiput(143.67,53.03)(0.13,0.21){2}{\line(0,1){0.21}}
\multiput(143.41,52.61)(0.13,0.21){2}{\line(0,1){0.21}}
\multiput(143.16,52.17)(0.12,0.22){2}{\line(0,1){0.22}}
\multiput(142.92,51.73)(0.12,0.22){2}{\line(0,1){0.22}}
\multiput(142.7,51.29)(0.11,0.22){2}{\line(0,1){0.22}}
\multiput(142.47,50.84)(0.11,0.22){2}{\line(0,1){0.22}}
\multiput(142.26,50.39)(0.11,0.23){2}{\line(0,1){0.23}}
\multiput(142.06,49.94)(0.1,0.23){2}{\line(0,1){0.23}}
\multiput(141.87,49.48)(0.1,0.23){2}{\line(0,1){0.23}}
\multiput(141.68,49.02)(0.09,0.23){2}{\line(0,1){0.23}}
\multiput(141.51,48.55)(0.18,0.47){1}{\line(0,1){0.47}}
\multiput(141.34,48.08)(0.17,0.47){1}{\line(0,1){0.47}}
\multiput(141.18,47.6)(0.16,0.47){1}{\line(0,1){0.47}}
\multiput(141.04,47.13)(0.15,0.48){1}{\line(0,1){0.48}}
\multiput(140.9,46.65)(0.14,0.48){1}{\line(0,1){0.48}}
\multiput(140.77,46.17)(0.13,0.48){1}{\line(0,1){0.48}}
\multiput(140.65,45.68)(0.12,0.48){1}{\line(0,1){0.48}}
\multiput(140.55,45.2)(0.11,0.49){1}{\line(0,1){0.49}}
\multiput(140.45,44.71)(0.1,0.49){1}{\line(0,1){0.49}}
\multiput(140.36,44.22)(0.09,0.49){1}{\line(0,1){0.49}}
\multiput(140.28,43.73)(0.08,0.49){1}{\line(0,1){0.49}}
\multiput(140.21,43.23)(0.07,0.49){1}{\line(0,1){0.49}}
\multiput(140.15,42.74)(0.06,0.49){1}{\line(0,1){0.49}}
\multiput(140.1,42.24)(0.05,0.5){1}{\line(0,1){0.5}}
\multiput(140.06,41.74)(0.04,0.5){1}{\line(0,1){0.5}}
\multiput(140.03,41.25)(0.03,0.5){1}{\line(0,1){0.5}}
\multiput(140.01,40.75)(0.02,0.5){1}{\line(0,1){0.5}}
\multiput(140,40.25)(0.01,0.5){1}{\line(0,1){0.5}}
\put(140,39.75){\line(0,1){0.5}}
\multiput(140,39.75)(0.01,-0.5){1}{\line(0,-1){0.5}}
\multiput(140.01,39.25)(0.02,-0.5){1}{\line(0,-1){0.5}}
\multiput(140.03,38.75)(0.03,-0.5){1}{\line(0,-1){0.5}}
\multiput(140.06,38.26)(0.04,-0.5){1}{\line(0,-1){0.5}}
\multiput(140.1,37.76)(0.05,-0.5){1}{\line(0,-1){0.5}}
\multiput(140.15,37.26)(0.06,-0.49){1}{\line(0,-1){0.49}}
\multiput(140.21,36.77)(0.07,-0.49){1}{\line(0,-1){0.49}}
\multiput(140.28,36.27)(0.08,-0.49){1}{\line(0,-1){0.49}}
\multiput(140.36,35.78)(0.09,-0.49){1}{\line(0,-1){0.49}}
\multiput(140.45,35.29)(0.1,-0.49){1}{\line(0,-1){0.49}}
\multiput(140.55,34.8)(0.11,-0.49){1}{\line(0,-1){0.49}}
\multiput(140.65,34.32)(0.12,-0.48){1}{\line(0,-1){0.48}}
\multiput(140.77,33.83)(0.13,-0.48){1}{\line(0,-1){0.48}}
\multiput(140.9,33.35)(0.14,-0.48){1}{\line(0,-1){0.48}}
\multiput(141.04,32.87)(0.15,-0.48){1}{\line(0,-1){0.48}}
\multiput(141.18,32.4)(0.16,-0.47){1}{\line(0,-1){0.47}}
\multiput(141.34,31.92)(0.17,-0.47){1}{\line(0,-1){0.47}}
\multiput(141.51,31.45)(0.18,-0.47){1}{\line(0,-1){0.47}}
\multiput(141.68,30.98)(0.09,-0.23){2}{\line(0,-1){0.23}}
\multiput(141.87,30.52)(0.1,-0.23){2}{\line(0,-1){0.23}}
\multiput(142.06,30.06)(0.1,-0.23){2}{\line(0,-1){0.23}}
\multiput(142.26,29.61)(0.11,-0.23){2}{\line(0,-1){0.23}}
\multiput(142.47,29.16)(0.11,-0.22){2}{\line(0,-1){0.22}}
\multiput(142.7,28.71)(0.11,-0.22){2}{\line(0,-1){0.22}}
\multiput(142.92,28.27)(0.12,-0.22){2}{\line(0,-1){0.22}}
\multiput(143.16,27.83)(0.12,-0.22){2}{\line(0,-1){0.22}}
\multiput(143.41,27.39)(0.13,-0.21){2}{\line(0,-1){0.21}}
\multiput(143.67,26.97)(0.13,-0.21){2}{\line(0,-1){0.21}}
\multiput(143.93,26.54)(0.14,-0.21){2}{\line(0,-1){0.21}}
\multiput(144.2,26.13)(0.14,-0.21){2}{\line(0,-1){0.21}}
\multiput(144.48,25.71)(0.14,-0.2){2}{\line(0,-1){0.2}}
\multiput(144.77,25.31)(0.15,-0.2){2}{\line(0,-1){0.2}}
\multiput(145.07,24.91)(0.1,-0.13){3}{\line(0,-1){0.13}}
\multiput(145.37,24.51)(0.1,-0.13){3}{\line(0,-1){0.13}}
\multiput(145.69,24.13)(0.11,-0.13){3}{\line(0,-1){0.13}}
\multiput(146.01,23.74)(0.11,-0.13){3}{\line(0,-1){0.13}}
\multiput(146.34,23.37)(0.11,-0.12){3}{\line(0,-1){0.12}}
\multiput(146.67,23)(0.11,-0.12){3}{\line(0,-1){0.12}}
\multiput(147.01,22.64)(0.12,-0.12){3}{\line(0,-1){0.12}}
\multiput(147.36,22.28)(0.12,-0.12){3}{\line(1,0){0.12}}
\multiput(147.72,21.93)(0.12,-0.11){3}{\line(1,0){0.12}}
\multiput(148.08,21.59)(0.12,-0.11){3}{\line(1,0){0.12}}
\multiput(148.45,21.26)(0.13,-0.11){3}{\line(1,0){0.13}}
\multiput(148.83,20.93)(0.13,-0.11){3}{\line(1,0){0.13}}
\multiput(149.21,20.61)(0.13,-0.1){3}{\line(1,0){0.13}}
\multiput(149.6,20.3)(0.13,-0.1){3}{\line(1,0){0.13}}

\linethickness{0.3mm}
\put(100,40){\line(1,0){40}}
\linethickness{0.3mm}
\multiput(70,20)(0.18,0.12){167}{\line(1,0){0.18}}
\linethickness{0.3mm}
\multiput(70,20)(0.42,0.12){167}{\line(1,0){0.42}}
\linethickness{0.3mm}
\multiput(70,20)(1.79,0.12){42}{\line(1,0){1.79}}
\put(150,40){\makebox(0,0)[cc]{$U_S$}}

\put(137,44){\makebox(0,0)[cc]{$z$}}

\put(100,43){\makebox(0,0)[cc]{$y$}}

\put(70,25){\makebox(0,0)[cc]{$x$}}

\put(105,35){\makebox(0,0)[cc]{$\delta$}}

\put(143,21){\makebox(0,0)[cc]{$z'$}}

\put(82,33){\makebox(0,0)[cc]{$r$}}

\put(119,44){\makebox(0,0)[cc]{$\geq N$}}

\end{picture}
\noindent 
On va  montrer que $d(x,z)\leq d(x,U_{S})+\delta$. Soit $z'\in U_{S}$. On a $d(z,z')\leq 2N$ et $d(y,z')\geq d(y,z)\geq N$. 
 Par $(H_{\delta}(x,y,z,z'))$ on a 
$$d(x,z)\leq \max(d(x,z')+d(y,z)-d(y,z'),d(x,y)+d(z,z')-d(y,z'))+\delta.$$    Mais $d(x,z')+d(y,z)-d(y,z')\leq d(x,z')$ et $$d(x,y)+d(z,z')-d(y,z')\leq 
r+2N-N\leq d(x,U_{S}).$$ Donc $d(x,z)\leq d(x,z')+\delta$ pour tout  $z'\in U_{S}$. On a montré $Y_{S,x,r}\subset A_{S,x}$.

Soit maintenant $z\in A_{S,x}$ et $y\in \delta\text{-}\geod(x,z)$ tels que $d(x,y)=r$ et $d(y,z)=d(y,A_{S,x})$. On veut montrer $d(y,A_{S,x})=d(y,U_{S})$. Si ce n'est pas vrai, il existe $z'\in U_{S}\setminus A_{S,x}$ tels que $d(y,z')=d(y,U_{S})$. Comme $d(x,z')\geq d(x,U_{S})+\de \geq d(x,z)$ et $d(y,z')\leq d(y,z)$, 
on a $y\in \de\tg(x,z')$. Alors $z'\in Y_{S,x,r}\subset A_{S,x}$, ce qui amène une contradiction. 
 \cqfd

L'astuce des ensembles emboîtés repose sur le lemme suivant. 
\begin{lem}\label{Y-emboites}
Si $r\leq d(x,U_{S})-d(x,x')-\delta$, on a $Y_{S,x,r}\subset Y_{S,x',r+2d(x,x')+\delta}$.
\end{lem}
\noindent{\bf Démonstration du lemme~\ref{Y-emboites} en admettant le lemme~\ref{contract2}.}
Soit $z\in Y_{S,x,r}$. Par hypothèse il existe $y\in \delta\text{-}\geod(x,z)$, $d(x,y)\leq r$, $d(y,z)=d(y,U_{S})$.   
On a $$d(y,z)\geq d(x,z)-d(x,y)\geq d(x,U_{S})-r\geq d(x,x')+\de. $$
 Soit $y'\in \geod(y,z)$, $d(y,y')=d(x,x')+\delta$.

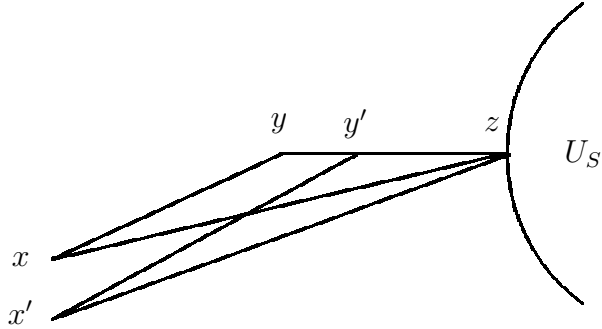
 \ifx\JPicScale\undefined\def\JPicScale{1}\fi
\unitlength \JPicScale mm
\begin{picture}(150,55)(20,11)
\linethickness{0.3mm}
\multiput(109.6,59.7)(0.13,0.1){3}{\line(1,0){0.13}}
\multiput(109.21,59.39)(0.13,0.1){3}{\line(1,0){0.13}}
\multiput(108.83,59.07)(0.13,0.11){3}{\line(1,0){0.13}}
\multiput(108.45,58.74)(0.13,0.11){3}{\line(1,0){0.13}}
\multiput(108.08,58.41)(0.12,0.11){3}{\line(1,0){0.12}}
\multiput(107.72,58.07)(0.12,0.11){3}{\line(1,0){0.12}}
\multiput(107.36,57.72)(0.12,0.12){3}{\line(1,0){0.12}}
\multiput(107.01,57.36)(0.12,0.12){3}{\line(0,1){0.12}}
\multiput(106.67,57)(0.11,0.12){3}{\line(0,1){0.12}}
\multiput(106.34,56.63)(0.11,0.12){3}{\line(0,1){0.12}}
\multiput(106.01,56.26)(0.11,0.13){3}{\line(0,1){0.13}}
\multiput(105.69,55.87)(0.11,0.13){3}{\line(0,1){0.13}}
\multiput(105.37,55.49)(0.1,0.13){3}{\line(0,1){0.13}}
\multiput(105.07,55.09)(0.1,0.13){3}{\line(0,1){0.13}}
\multiput(104.77,54.69)(0.15,0.2){2}{\line(0,1){0.2}}
\multiput(104.48,54.29)(0.14,0.2){2}{\line(0,1){0.2}}
\multiput(104.2,53.87)(0.14,0.21){2}{\line(0,1){0.21}}
\multiput(103.93,53.46)(0.14,0.21){2}{\line(0,1){0.21}}
\multiput(103.67,53.03)(0.13,0.21){2}{\line(0,1){0.21}}
\multiput(103.41,52.61)(0.13,0.21){2}{\line(0,1){0.21}}
\multiput(103.16,52.17)(0.12,0.22){2}{\line(0,1){0.22}}
\multiput(102.92,51.73)(0.12,0.22){2}{\line(0,1){0.22}}
\multiput(102.7,51.29)(0.11,0.22){2}{\line(0,1){0.22}}
\multiput(102.47,50.84)(0.11,0.22){2}{\line(0,1){0.22}}
\multiput(102.26,50.39)(0.11,0.23){2}{\line(0,1){0.23}}
\multiput(102.06,49.94)(0.1,0.23){2}{\line(0,1){0.23}}
\multiput(101.87,49.48)(0.1,0.23){2}{\line(0,1){0.23}}
\multiput(101.68,49.02)(0.09,0.23){2}{\line(0,1){0.23}}
\multiput(101.51,48.55)(0.18,0.47){1}{\line(0,1){0.47}}
\multiput(101.34,48.08)(0.17,0.47){1}{\line(0,1){0.47}}
\multiput(101.18,47.6)(0.16,0.47){1}{\line(0,1){0.47}}
\multiput(101.04,47.13)(0.15,0.48){1}{\line(0,1){0.48}}
\multiput(100.9,46.65)(0.14,0.48){1}{\line(0,1){0.48}}
\multiput(100.77,46.17)(0.13,0.48){1}{\line(0,1){0.48}}
\multiput(100.65,45.68)(0.12,0.48){1}{\line(0,1){0.48}}
\multiput(100.55,45.2)(0.11,0.49){1}{\line(0,1){0.49}}
\multiput(100.45,44.71)(0.1,0.49){1}{\line(0,1){0.49}}
\multiput(100.36,44.22)(0.09,0.49){1}{\line(0,1){0.49}}
\multiput(100.28,43.73)(0.08,0.49){1}{\line(0,1){0.49}}
\multiput(100.21,43.23)(0.07,0.49){1}{\line(0,1){0.49}}
\multiput(100.15,42.74)(0.06,0.49){1}{\line(0,1){0.49}}
\multiput(100.1,42.24)(0.05,0.5){1}{\line(0,1){0.5}}
\multiput(100.06,41.74)(0.04,0.5){1}{\line(0,1){0.5}}
\multiput(100.03,41.25)(0.03,0.5){1}{\line(0,1){0.5}}
\multiput(100.01,40.75)(0.02,0.5){1}{\line(0,1){0.5}}
\multiput(100,40.25)(0.01,0.5){1}{\line(0,1){0.5}}
\put(100,39.75){\line(0,1){0.5}}
\multiput(100,39.75)(0.01,-0.5){1}{\line(0,-1){0.5}}
\multiput(100.01,39.25)(0.02,-0.5){1}{\line(0,-1){0.5}}
\multiput(100.03,38.75)(0.03,-0.5){1}{\line(0,-1){0.5}}
\multiput(100.06,38.26)(0.04,-0.5){1}{\line(0,-1){0.5}}
\multiput(100.1,37.76)(0.05,-0.5){1}{\line(0,-1){0.5}}
\multiput(100.15,37.26)(0.06,-0.49){1}{\line(0,-1){0.49}}
\multiput(100.21,36.77)(0.07,-0.49){1}{\line(0,-1){0.49}}
\multiput(100.28,36.27)(0.08,-0.49){1}{\line(0,-1){0.49}}
\multiput(100.36,35.78)(0.09,-0.49){1}{\line(0,-1){0.49}}
\multiput(100.45,35.29)(0.1,-0.49){1}{\line(0,-1){0.49}}
\multiput(100.55,34.8)(0.11,-0.49){1}{\line(0,-1){0.49}}
\multiput(100.65,34.32)(0.12,-0.48){1}{\line(0,-1){0.48}}
\multiput(100.77,33.83)(0.13,-0.48){1}{\line(0,-1){0.48}}
\multiput(100.9,33.35)(0.14,-0.48){1}{\line(0,-1){0.48}}
\multiput(101.04,32.87)(0.15,-0.48){1}{\line(0,-1){0.48}}
\multiput(101.18,32.4)(0.16,-0.47){1}{\line(0,-1){0.47}}
\multiput(101.34,31.92)(0.17,-0.47){1}{\line(0,-1){0.47}}
\multiput(101.51,31.45)(0.18,-0.47){1}{\line(0,-1){0.47}}
\multiput(101.68,30.98)(0.09,-0.23){2}{\line(0,-1){0.23}}
\multiput(101.87,30.52)(0.1,-0.23){2}{\line(0,-1){0.23}}
\multiput(102.06,30.06)(0.1,-0.23){2}{\line(0,-1){0.23}}
\multiput(102.26,29.61)(0.11,-0.23){2}{\line(0,-1){0.23}}
\multiput(102.47,29.16)(0.11,-0.22){2}{\line(0,-1){0.22}}
\multiput(102.7,28.71)(0.11,-0.22){2}{\line(0,-1){0.22}}
\multiput(102.92,28.27)(0.12,-0.22){2}{\line(0,-1){0.22}}
\multiput(103.16,27.83)(0.12,-0.22){2}{\line(0,-1){0.22}}
\multiput(103.41,27.39)(0.13,-0.21){2}{\line(0,-1){0.21}}
\multiput(103.67,26.97)(0.13,-0.21){2}{\line(0,-1){0.21}}
\multiput(103.93,26.54)(0.14,-0.21){2}{\line(0,-1){0.21}}
\multiput(104.2,26.13)(0.14,-0.21){2}{\line(0,-1){0.21}}
\multiput(104.48,25.71)(0.14,-0.2){2}{\line(0,-1){0.2}}
\multiput(104.77,25.31)(0.15,-0.2){2}{\line(0,-1){0.2}}
\multiput(105.07,24.91)(0.1,-0.13){3}{\line(0,-1){0.13}}
\multiput(105.37,24.51)(0.1,-0.13){3}{\line(0,-1){0.13}}
\multiput(105.69,24.13)(0.11,-0.13){3}{\line(0,-1){0.13}}
\multiput(106.01,23.74)(0.11,-0.13){3}{\line(0,-1){0.13}}
\multiput(106.34,23.37)(0.11,-0.12){3}{\line(0,-1){0.12}}
\multiput(106.67,23)(0.11,-0.12){3}{\line(0,-1){0.12}}
\multiput(107.01,22.64)(0.12,-0.12){3}{\line(0,-1){0.12}}
\multiput(107.36,22.28)(0.12,-0.12){3}{\line(1,0){0.12}}
\multiput(107.72,21.93)(0.12,-0.11){3}{\line(1,0){0.12}}
\multiput(108.08,21.59)(0.12,-0.11){3}{\line(1,0){0.12}}
\multiput(108.45,21.26)(0.13,-0.11){3}{\line(1,0){0.13}}
\multiput(108.83,20.93)(0.13,-0.11){3}{\line(1,0){0.13}}
\multiput(109.21,20.61)(0.13,-0.1){3}{\line(1,0){0.13}}
\multiput(109.6,20.3)(0.13,-0.1){3}{\line(1,0){0.13}}

\linethickness{0.3mm}
\put(70,40){\line(1,0){30}}
\linethickness{0.3mm}
\multiput(40,26)(0.26,0.12){117}{\line(1,0){0.26}}
\linethickness{0.3mm}
\multiput(40,26)(0.51,0.12){117}{\line(1,0){0.51}}
\linethickness{0.3mm}
\multiput(40,18)(0.22,0.12){183}{\line(1,0){0.22}}
\linethickness{0.3mm}
\multiput(40,18)(0.33,0.12){183}{\line(1,0){0.33}}
\put(110,40){\makebox(0,0)[cc]{$U_S$}}

\put(80,44){\makebox(0,0)[cc]{$y'$}}

\put(98,44){\makebox(0,0)[cc]{$z$}}

\put(70,44){\makebox(0,0)[cc]{$y$}}

\put(36,26){\makebox(0,0)[cc]{$x$}}

\put(36,19){\makebox(0,0)[cc]{$x'$}}

\end{picture}
 
 \noindent 
 Alors $d(y',z)=d(y',U_{S})$, et par le lemme~\ref{contract2} (avec $\epsilon=\de$) on a $y'\in \delta\text{-}\geod(x',z)$. Enfin $$d(x',y')\leq d(x',x)+d(x,y)+d(y,y')\leq r+2d(x,x')+\delta. $$
\cqfd

\begin{lem}\label{contract1}
Pour tout $\epsilon >0$ et $x,z,y,y'\in X$, si $y\in
\epsilon\text{-}\geod(x,z)$, $y'\in \geod(y,z)$, $d(y,y')\geq
\epsilon/2$, alors $y'\in \delta\text{-}\geod(x,z)$.
\end{lem}

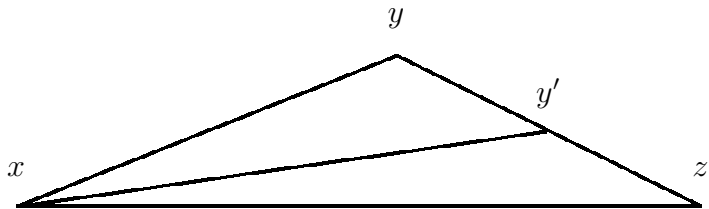
\ifx\JPicScale\undefined\def\JPicScale{1}\fi
\unitlength \JPicScale mm
\begin{picture}(120,45)(20,20)
\linethickness{0.3mm}
\put(30,30){\line(1,0){90}}
\linethickness{0.3mm}
\multiput(80,50)(0.24,-0.12){167}{\line(1,0){0.24}}
\linethickness{0.3mm}
\multiput(30,30)(0.3,0.12){167}{\line(1,0){0.3}}
\linethickness{0.3mm}
\multiput(30,30)(0.84,0.12){83}{\line(1,0){0.84}}
\put(80,55){\makebox(0,0)[cc]{$y$}}

\put(100,45){\makebox(0,0)[cc]{$y'$}}

\put(120,35){\makebox(0,0)[cc]{$z$}}

\put(30,35){\makebox(0,0)[cc]{$x$}}

\end{picture}

\noindent{\bf Démonstration.}
Par $(H_{\delta}^{0}(x,y,y',z))$ on a 
$$d(x,y')\leq
\max(d(x,y)-d(y,y'),d(x,z)-d(z,y'))+\delta. $$
D'où $d(x,y')+d(y',z)\leq
\max(d(x,y)-d(y,y')+d(y',z),d(x,z))+\delta$. 
Or $$d(x,y)-d(y,y')+d(y',z)=(d(x,y)+d(y,z))-2d(y,y')\leq
d(x,z)+\epsilon -2d(y,y').$$ 
Donc $d(x,y')+d(y',z)\leq
d(x,z)+\de$.
\cqfd

\begin{lem}\label{contract2}
Pour tous $\epsilon >0$, 
$x,x',z,y,y'\in X$, si $y\in 
\epsilon\text{-}\geod(x,z)$, $y'\in \geod(y,z)$, 
$d(y,y')\geq 
\epsilon/2+d(x,x')$, alors $y'\in \delta\text{-}\geod(x',z)$.
\end{lem}
\noindent{\bf Démonstration.} D'après le lemme~\ref{xx'yy'zz'},  $y\in (\epsilon+2d(x,x'))\text{-}\geod(x',z)$. On applique alors le lemme~\ref{contract1} à $(x',z,y,y')$ au lieu de $(x,z,y,y')$ et $\epsilon+2d(x,x')$ au lieu de $\epsilon$.   \cqfd

 Dans~\cite{ks} Kasparov et Skandalis définissent une mesure $\psi_{S,x}$ de masse $1$ en normalisant une moyenne paramétrée par $r$ des fonctions caractéristiques de $Y_{S,x,r}$. Nous allons faire de même, à ceci près que nous prendrons plutôt une moyenne sur $r$ de la fonction caractéristique de $Y_{S,x,r}$ normalisée. 
 Pour tout ensemble non vide $A$ on note $\nu_{A}=\frac{1}{\sharp A}\chi_{A}$, où $\chi_{A}$ est la fonction caractéristique de $A$. 
  
 On pose alors 
 $$\psi_{S,x}=\frac{1}{\max(1,d(x,U_{S})-N)}
 \sum_{r=0}^{\max(0,d(x,U_{S})-N-1)} 
 \nu_{Y_{S,x,r}}.$$
 \label{psiSx}
 
 \noindent D'après le lemme~\ref{majoration-r} le support de $\psi_{S,x}$ est inclus dans $A_{S,x}$. 
 
 Plus loin nous aurons besoin de la notation suivante :  on définit   $$\psi_{S,x,t}=\nu_{Y_{S,x,\max(0,E(t(d(x,U_{S})-N)))}}\text{\ \ \  pour }t\in [0,1]$$  de sorte que $\psi_{S,x}=\int_{0}^{1}\psi_{S,x,t}dt$. Dans tout cet article on note $E(.)$ la partie entière.

 \begin{lem}\label{support-psiSxt}
 Pour tout $t\in [0,1]$,  le support de $\psi_{S,x,t}$ est inclus dans $A_{S,x}$. 
 \end{lem}
 \noindent{\bf Démonstration.} 
Cela résulte immédiatement du lemme~\ref{majoration-r}. \cqfd
   
 Pour $r\in \{0,...,\max(0,d(x,U_{S})-N)\}$, on a, grâce au lemme~\ref{majoration-r},  
  $$Y_{S,x,r}=
\{z\in A_{S,x},\exists y\in \delta\text{-}\geod(x,z),d(x,y)= r,d(y,z)=d(y,A_{S,x})\}.$$ Donc  $Y_{S,x,r}$ ne dépend que de la connaissance des points de 
\begin{gather}\label{conn-YSxr-22oct09}S\cup A_{S,x}\cup \{y,\exists z\in A_{S,x},\ y\in \delta\text{-}\geod(x,z), d(x,y)=r\}\cup \{x\}\end{gather} et des distances mutuelles entre tous ces points.
De fa\c con plus précise, 
si $x'$ et $S'$ sont tels qu'il 
 existe une isométrie de \eqref{conn-YSxr-22oct09} vers l'ensemble correspondant pour $x'$ et $S'$, qui envoie $x$ et $S$ sur $x'$ et $S'$, alors $Y_{S',x',r}\subset A_{S',x'}$ est l'image de $Y_{S,x,r}\subset A_{S,x}$ par cette isométrie.

 En anticipant un peu justifions l'intérêt de cette propriété. 
D'après  la démonstration du lemme~\ref{hxt-eS-connaissance}, 
l'ensemble (\ref{conn-YSxr-22oct09}) est inclus dans l'ensemble figurant dans l'énoncé du
lemme~\ref{hxt-eS-connaissance} et on  déduit du lemme~\ref{cardinal-tranche-geod} 
 que le cardinal de (\ref{conn-YSxr-22oct09}) est borné par une constante $C(\delta,K,N)$, et de plus   la connaissance de $r$ et de $d(x,U_{S})$ détermine les distances  entre les points de (\ref{conn-YSxr-22oct09}) à une constante $C(\delta,N)$ près. Donc connaissant $r$ et  $d(x,U_{S})$, il n'y a qu'un nombre fini,  borné par $C(\delta,K,N)$, de possibilités pour la donnée de toutes les distances  entre les points de (\ref{conn-YSxr-22oct09}). Enfin $\psi_{S,x}$ est une moyenne entre $0$ et $\max(0,d(x,U_{S})-N-1)$ de $\nu_{Y_{S,x,r}}$. Au contraire avec la définition de Kasparov et Skandalis, où la normalisation est faite après la moyennne, le calcul de $\psi_{S,x}$ nécessiterait la connaissance simultanée de tous les points 
  $y$ tels qu'il existe $z\in A_{S,x}$ vérifiant  $y\in \delta\text{-}\geod(x,z)$.

  Le lemme  suivant, que nous venons d'utiliser, servira à de nombreuses reprises. 
   \begin{lem}\label{cardinal-tranche-geod}
Pour tout $r\in \N$ il existe  $C$ dépendant seulement de $\de,K$ et $r$ tel que pour $x,y\in X$, on ait, 
\begin{itemize} \item pour tout $l\in \{0,\pp,d(x,y)+r\}$, 
$\sharp\big\{z\in r\tg(x,y), d(x,z)=l\big\}\leq C$
\item et 
$\sharp\big(r\tg(x,y)\big)\leq C(d(x,y)+1)$.
\end{itemize}
\end{lem}
 \noindent{\bf Démonstration.}
Pour $x,y\in X$ on pose $d=d(x,y)$ et on choisit $x_{0}=x,x_{1},\pp,x_{d}=y$ des points de $\geod(x,y)$ tels que $d(x,x_{i})=i$. Soit $z\in r\tg(x,y)$ et $i=\min(d,d(x,z))$. 
Alors 
par $(H_{\delta}^{0}(z,x,x_{i},y))$ on a $$d(z,x_{i})\leq \max(d(z,x)-d(x,x_{i}),d(z,y)-d(y,x_{i}))+\de\leq r+\de$$
car $d(x,x_{i})=i$, $d(x_{i},y)=d-i$, $d(x,z)\in [i,i+r]$ et $d(z,y)\in [d-i,d-i+r]$. 
On  pose $C=1+K+K^{2}+\pp+K^{r+\de}$. Alors $\max_{x\in X}\sharp B(x,r+\de)\leq C$.   On voit que la première assertion est vraie et pour la deuxième assertion on remarque que 
 $r\tg(x,y)\subset \bigcup_{i=0}^{d}B(x_{i},r+\de)$. \cqfd

  Par l'astuce des ensembles emboîtés nous allons montrer que 
  pour  $x,x'\in X$, $\|\psi_{S,x}-\psi_{S,x'}\|_{1}$ tend vers $0$ en dehors des parties finies de $\Delta$ (où $\|.\|_1$ désigne la masse totale d'une mesure). Cela résulte du lemme suivant qui est plus fort. 
  
  \begin{lem}\label{lem-psix-psix'-8j1630}
   Il existe $C=C(\de,K,N)$ tel que la mesure de l'ensemble des $t\in [0,1]$ tels que $\psi_{S,x,t}\neq\psi_{S,x',t}$ soit $\leq \frac{Cd(x,x')}{1+d(x,S)}$.
 \end{lem}
  \noindent{\bf Démonstration.} 
  Le lemme est vrai si $x=x'$, donc on suppose $d(x,x')\geq1$. On prendra $C\geq 2N+2$ si bien que $\frac{Cd(x,x')}{1+d(x,S)}\geq 1$ si $d(x,S)\leq d(x,x')+2N$. Donc on  suppose 
 $$d(x,S)\geq d(x,x')+2N+1,$$ et il suffit de montrer le lemme sous cette hypothèse. 
  
  On a  alors  \begin{gather}\nonumber\min(d(x,U_{S})-N, d(x',U_{S})-N)\geq d(x,S)-2N-d(x,x')\geq 1.\end{gather} 
  
Comme le cardinal de $Y_{S,x,r}$ est non nul et borné par une constante $C_{0}=C(\delta,K,N)$, et que l'application $r\mapsto Y_{S,x,r}$ est croissante, l'intervalle $$[0,d(x,U_{S})-N-1]$$ se découpe en au plus $C_{0}$ intervalles où l'application $r\mapsto Y_{S,x,r}$ est  constante. Il résulte du  lemme~\ref{Y-emboites}
  que l'application $r\mapsto Y_{S,x',r}$ coïncide avec la précédente sur des intervalles (éventuellement vides) obtenus à partir des précédents en raccourcissant chaque extrémité de $2d(x,x')+\delta$. 
  
 Il en résulte que l'application croissante $$t\mapsto Y_{S,x,E(t(d(x,U_{S})-N))}$$ prend au plus $C_{0}$ valeurs et pour chaque valeur l'image inverse est un intervalle, et l'image inverse de cette même valeur par l'application 
$$t\mapsto Y_{S,x',E(t(d(x',U_{S})-N))}$$  est un autre intervalle dont les extrémités diffèrent au plus de $\frac{3d(x,x')+\delta+1}{
d(x,S)-2N}$ 
des extrémités du premier. 
En effet pour $t,t'\in [0,1]$, l'inégalité $$|E(t(d(x,U_{S})-N))-E(t'(d(x',U_{S})-N))|\leq 2d(x,x')+\delta$$ implique $|t-t'|\leq \frac{3d(x,x')+\delta+1}{
d(x,U_{S})-N}$
car $|d(x,U_{S})-d(x',U_{S})|\leq d(x,x')$ et on utilise le fait que $d(x,U_{S})-N\geq d(x,S)-2N$. 

 On en déduit que la mesure  de l'ensemble des  $t\in [0,1]$ tels que 
$\mu_{r,t}(x,a)\neq \mu_{r,t}(x',a)$ 
est inférieure ou égale à $\frac{2C_{0}(3d(x,x')+\delta+1)}{
d(x,S)-2N}$. Enfin on a  
  $$\frac{2C_{0}(3d(x,x')+\delta+1)}{
d(x,S)-2N}\leq \frac{2C_{0}(4+\de)(2N+2)d(x,x')}{1+d(x,S)}$$ car on a supposé  $d(x,x')\geq1$ et $d(x,S)\geq 2N+1$. On prend alors $C=2C_{0}(4+\de)(2N+2)$. 
   \cqfd

On définit un opérateur $h_x$ de degré $1$ sur 
$\oplus_{p=0}^{p_{\max}}\C^{(\Delta_{p})}$ par la formule suivante : $$h_x(e_S)=\psi_{S,x}\wedge e_{S}.$$ 
Plus loin nous aurons besoin de la notation suivante : pour $t\in [0,1]$,  $h_{x,t}$ est l'opérateur défini par  \begin{gather}\label{form-hxt-8j1629}h_{x,t}(e_S)=\psi_{S,x,t}\wedge e_{S},\end{gather} de sorte que $h_{x}=\int_{0}^{1}h_{x,t}dt$. 
On note que pour tout $t\in [0,1]$, $h_{x,t}(e_{\emptyset})=e_{x}$. 
\label{def-hx-hxs}
 \begin{lem}\label{lem-hxt-hx't-8j1630}
   Il existe $C=C(\de,K,N)$ tel que la mesure de l'ensemble des $t\in [0,1]$ tels que $(h_{x,t}-h_{x',t})(e_S)\neq 0$ soit $\leq \frac{Cd(x,x')}{1+d(x,S)}$.
 \end{lem}
  \noindent{\bf Démonstration.} 
  Cela résulte immédiatement de (\ref{form-hxt-8j1629}) et du lemme~\ref{lem-psix-psix'-8j1630}. \cqfd

Le lemme suivant servira tout à fait à la fin de l'article. 

\begin{lem}\label{lemmeh2=0}
On a $h_{x}^{2}=0$. 
\end{lem}
\noindent{\bf Démonstration.} D'après les lemmes~\ref{73} et~\ref{support-psiSxt}, si $e_{T}$ apparaît dans $h_{x}(e_{S})$ on a $A_{T,x}=A_{S,x}$. 
Le lemme~\ref{majoration-r} montre alors que $\psi_{T,x}=\psi_{S,x}$. \cqfd

Nous voulons montrer que l'opérateur $\del h_{x}+h_{x}\del$ est inversible de $\oplus_{p=0}^{p_{\max}}\C^{(\Delta_{p})}$ dans lui-même. En effet $H_x=h_x(\del h_x+h_x\del)^{-1}$ sera alors un parametrix pour $\del$, c'est-à-dire que l'on aura $\del H_x+H_x\del=1$. 
Pour cela nous introduisons, comme dans le paragraphe 2.5.1 de~\cite{kkban}, la ``distance moyenne tronquée de $S$ à $x$'' : 
$$\zeta_{x}(S)=\frac{1}{p_{\max}}\Big(\sum _{a\in S}
\max (d(x,U_{S}),d(x,a)-\delta)+(p_{\max}-\sharp S )d(x,U_{S})\Big),$$
si $S\neq\emptyset$ et $\zeta_{x}(\emptyset)=0$.
\label{zetaxS}
Cette distance est tronquée au sens où tous les points de $S$ dont la distance à $x$ est comprise entre $d(x,U_{S})$ et $d(x,U_{S})+\delta$ contribuent de la même fa\c con. 
On a clairement $d(x,U_{S}) \leq \zeta_{x}(S) \leq d(x,U_{S}) +N-\delta$. 

\begin{lem}\label{nilp}
Pour tout $S\in \Delta$ on a 
$$\del(e_{S})\in \bigoplus_{T\text{ tel que }
\zeta_{x}(T)\leq \zeta_{x}(S)}\C e_{T},
$$
$$h_x(e_{S})\in \bigoplus_{T\text{ tel que }
\zeta_{x}(T)\leq \zeta_{x}(S)}\C e_{T},
$$
$$(1-\del h_x-h_x\del)(e_{S})\in \bigoplus_{T\text{ tel que }
\zeta_{x}(T)< \zeta_{x}(S)-\frac{N-6\delta}{p_{\max}}}\C e_{T}.$$
\end{lem}
\noindent{\bf Démonstration.} 
La première assertion est évidente. La deuxième assertion est vraie car  le support de $\psi_{S,x}$ est inclus dans $A_{S,x}$, en vertu du lemme~\ref{support-psiSxt} et  pour $T=S\cup \{a\}$ avec $a\in A_{S,x}$ on a $d(x,U_{T})=d(x,U_{S})$ d'après le lemme~\ref{73}. Pour montrer la dernière assertion  on remarque que 
$$(1-\del h_x-h_x\del)(e_{S})=\sum_{a\in S}\pm (\psi_{S,x}-\psi_{S\setminus \{a\},x})\wedge e_{S\setminus \{a\}}.$$ Si $d(x,a)\leq d(x,U_{S})+N-5\delta$, on a $d(x,U_{S\setminus \{a\}})=d(x,U_{S})$ et 
$A_{S\setminus \{a\},x}=A_{S,x}$ par le lemme~\ref{73}, d'où 
$\psi_{S\setminus \{a\},x}=\psi_{S,x}$ par le lemme~\ref{majoration-r}. Si 
$d(x,a)> d(x,U_{S})+N-5\delta$, comme les supports de $\psi_{S,x}$ et $\psi_{S\setminus \{a\},x}$ sont inclus dans $A_{S,x}$ et $A_{S\setminus \{a\},x}$, $(\psi_{S,x}-\psi_{S\setminus \{a\},x})\wedge e_{S\setminus \{a\}}$ est une combinaison de $e_{S\setminus \{a\}\cup \{b\}}$ avec 
$$d(x,b)\leq \max( d(x,U_{S\setminus \{a\}})+\delta,d(x,U_{S})+\delta) = d(x,U_{S})+\delta,$$ ce qui fait que 
$\zeta_{x}(S\setminus \{a\}\cup \{b\})< \zeta_{x}(S)-\frac{N-6\delta 
}{p_{\max}}$.
%
%En effet si $a\in S$, on a soit $d(x,a)\leq d(x,U_{S})+N-5\delta$, soit 
%$d(x,a)> d(x,U_{S})+N-5\delta$. Dans le premier cas 
%$\psi_{S\setminus \{a\},x}=\psi_{S,x}$ et dans le deuxième cas 
%$h_x(e_{S\setminus \{a\}})$ est une combinaison de $e_{S\setminus \{a\}+\{b\}}$ avec 
%$d(x,b)\leq d(x,U_{S\setminus \{a\}})+\delta\leq d(x,U_{S})+\delta$, et donc 
%$\zeta_{x}(S\setminus \{a\}+\{b\})< \zeta_{x}(S)-\frac{N-6\delta 
%}{p_{\max}}$. Pour plus de détails voir les lemmes 2.5.3 et 2.5.4 de~\cite{kkban}.  
\cqfd

\begin{defi}
Pour  $p\in \{1,\dots,p_{\max}\}$ et $f=\sum_{S\in \Delta_{p}}f(S)e_{S}$, on note $$\mathrm{supp}(f)=\bigcup_{S\text{ tel que }f(S)\neq 0}S. $$
\end{defi}

On a $h_{x}=\int_{0}^{1}h_{x,t}dt$ et les lemmes suivants
permettent, pour $t\in [0,1]$ et $S\in \Delta$,  d'estimer le support de  $h_{x,t}(e_{S})$ et de savoir de quoi $h_{x,t}(e_{S})$ dépend.

\begin{lem}\label{ASx-geod-x}
Soit $x\in X$ et $S\in \Delta\setminus\{\emptyset\}$. 

\noindent a) Pour tout $t\in [0,1]$ le support de $h_{x,t}(e_{S})$ est inclus dans $S\cup A_{S,x}$. 

\noindent b)
%Pour $y\in A_{S,x}$ n'appartenant pas à $B(x,2\delta)$, il existe $a\in S$ n'appartenant pas à $B(x,2\delta)$ tel que $y\in 3\delta\text{-}\geod(x,a)$ et $d(a,y)> N-2\delta$. 
%
%Autrement dit 
%$$A_{S,x}\subset  \bigcup_{a\in S}\{y\in 3\delta\text{-}\geod(x,a), d(y,a)%\in ]N-2\delta,N]\} $$ si $d(x,S)> N+2\delta$ et 
%$$A_{S,x}\subset  B(x,2\delta)\cup \bigcup_{a\in S,a\not\in B(x,2\delta)}\{y\in 3\delta\text{-}\geod(x,a), d(y,a)\in ]N-2\delta,N]\} $$ en général. 
On a $$A_{S,x}\subset  \Big(B(x,2\delta)\cap U_{S}\Big)\cup \bigcup_{a\in S,a\not\in B(x,2\delta)}\{y\in 3\delta\text{-}\geod(x,a), d(y,a)\in ]N-2\delta,N]\} $$
En particulier on a toujours 
$A_{S,x}\subset \bigcup_{a\in S}4\delta\text{-}\geod(x,a)$.
\end{lem}
\noindent{\bf Démonstration.} 
Comme $h_{x,t}(e_S)=\psi_{S,x,t}\wedge e_{S}$, le a) résulte du lemme~\ref{support-psiSxt}. Pour montrer le b) on 
 rappelle que  $A_{S,x}=\{y\in U_S, d(x,y)\leq d(x,U_S)+\delta \}$.  Soit $y\in A_{S,x}$ n'appartenant pas à $B(x,2\delta)$. 
Soit $z\in \geod(x,y)$ tel que $d(y,z)=2\delta$. Comme $d(x,U_{S})\geq d(x,y)-\delta$ et $d(x,z)=d(x,y)-2\delta$, on a $z\not \in U_{S}$. 
Donc il existe $a\in S$ tel que $d(a,z)>N$.

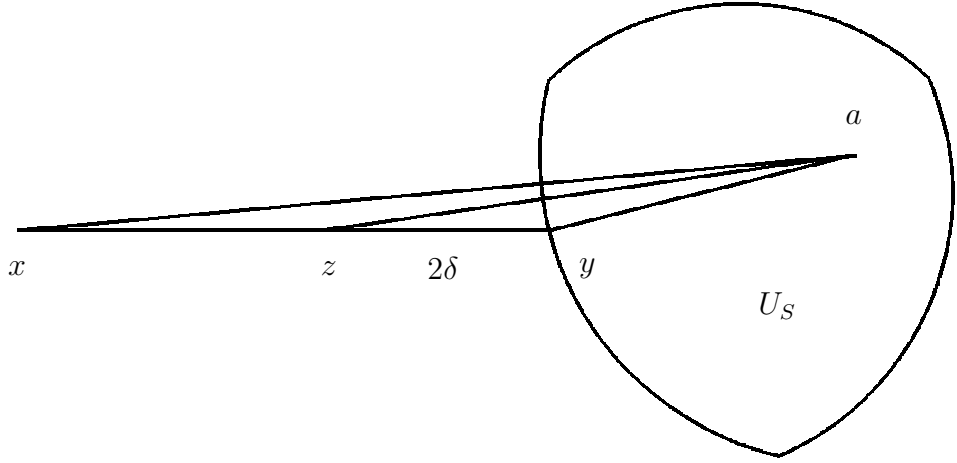
\ifx\JPicScale\undefined\def\JPicScale{1}\fi
\unitlength \JPicScale mm


\noindent  Comme $d(y,z)=2\delta$, 
on a $d(a,y)>N-2\delta$. Ensuite $(H_{\delta}^{0}(a,x,z,y))$ donne $$d(a,z)\leq \max(d(a,y)-2\delta, d(a,x)-d(x,y)+2\delta)+\delta.$$
Comme $d(a,z)>N$ et $d(a,y)\leq N$ on a 
  nécessairement 
  \begin{gather}\label{ineg-za-ax-xy-3de}
  d(a,z) \leq d(a,x)-d(x,y)+3\delta,
  \end{gather}
  et comme $d(a,y)\leq N<d(a,z)$ on en déduit $y\in 3\delta\text{-}\geod(x,a)$. 
  Par (\ref{ineg-za-ax-xy-3de}) on a  aussi $$d(x,a)\geq d(x,y)+N-3\delta \geq N-3\delta>2\delta$$ et $a\not\in B(x,2\delta)$. La seconde assertion de b) résulte immédiatement de la première : comme $S$ est non vide, $B(x,2\de)\subset \bigcup_{a\in S}4\delta\text{-}\geod(x,a)$.
  \cqfd

\begin{lem}\label{hxt-eS-connaissance}
Pour $t\in [0,1]$ et $S\in \Delta$, $h_{x,t}(e_{S})$ ne dépend que de la connaissance des points de 
\begin{gather}\nonumber B(x,2\delta)\cup S\cup \bigcup_{a\in S}\{y\in 3\delta\text{-}\geod(x,a), d(y,a)\in ]N-2\delta,N]\} \\ 
\label{ens-hxt-eS-connaissance}
\cup \bigcup_{a\in S}\{y\in 5\delta\text{-}\geod(x,a), 
d(x,y)\in [td(x,S)-2N-1,td(x,S)]\}\end{gather} et des distances  entre tous ces points. Plus précisément, si on note $S=\{a_{1},...,a_{p-1}\}$, si $T=\{b_{1},...,b_{p}\}$ est inclus dans \eqref{ens-hxt-eS-connaissance} (sans quoi le coefficient de $e_{T}$ dans $h_{x,t}(e_{S})$ est nul), 
et si $x'\in X$, $S'=\{a'_{1},...,a'_{p-1}\}$ et $T'=\{b'_{1},...,b'_{p}\}$ sont tels qu'il existe une isométrie de \eqref{ens-hxt-eS-connaissance}
dans l'ensemble correspondant pour $x'$ et $S'$, qui envoie 
$x,a_{1},...,a_{p-1},b_{1},...,b_{p}$ sur $x',a'_{1},...,a'_{p-1},b'_{1},...,b'_{p}$, alors le coefficient de $e_{b_{1}}\wedge ... \wedge e_{b_{p}}$ dans 
$h_{x,t}(e_{a_{1}}\wedge ...  \wedge e_{a_{p-1}})$ est égal au coefficient de 
$e_{b'_{1}}\wedge ... \wedge e_{b'_{p}}$ dans 
$h_{x',t}(e_{a'_{1}}\wedge ...  \wedge e_{a'_{p-1}})$. 

%Cet ensemble est lui-même inclus dans 
% $$ \bigcup_{a\in S}\Big\{y\in 5\delta\text{-}\geod(x,a), d(y,a)\in $$ $$\Big([0,N]\cup 
% [(1-t)d(x,a),(1-t)d(x,a)+3N+5\delta+1] \Big)\Big\}.$$
\end{lem}
\noindent{\bf Démonstration.} 
 On rappelle que $\psi_{S,x,t}=\nu_{Y_{S,x,r}}$ avec 
$$r=\max(0,E(t(d(x,U_{S})-N)))$$ et que, d'après le lemme~\ref{majoration-r},    $$Y_{S,x,r}=
\{z\in A_{S,x},\exists y\in \delta\text{-}\geod(x,z),d(x,y)= r,d(y,z)=d(y,A_{S,x})\}. $$
 Donc  $Y_{S,x,r}$ ne dépend que de la connaissance des points de 
$$S\cup A_{S,x}\cup \{y,\exists z\in A_{S,x},\ y\in \delta\text{-}\geod(x,z), d(x,y)=r\}\cup \{x\}$$ et des distances  entre tous ces points. D'après le    lemme~\ref{ASx-geod-x},  $$A_{S,x}\subset B(x,2\delta)\cup \bigcup_{a\in S}\{y\in 3\delta\text{-}\geod(x,a), d(y,a)\in ]N-2\delta,N]\}. $$
Soit maintenant $z\in A_{S,x}$ et $y\in \delta\text{-}\geod(x,z)$ tel que $d(x,y)=r$. D'après  le    lemme~\ref{ASx-geod-x}, il existe $a\in S$ tel que 
 $z\in 4\delta\text{-}\geod(x,a)$. 
Or    $y\in \delta\text{-}\geod(x,z)$ et $z\in 4\delta\text{-}\geod(x,a)$
impliquent 
%$d(x,y)+d(y,z)+d(z,a)\leq d(x,a)+5\de$, d'où 
$y\in 5\delta\text{-}\geod(x,a)$
par le a) du lemme~\ref{geod-comp-xabc}. Comme  $d(x,U_{S})\in [d(x,S)-N,d(x,S)]$ on a $r\in [td(x,S)-2N-1,td(x,S)]$. 
%Pour la deuxième partie on remarque que $d(x,y)\in [td(x,S)-2N-1,td(x,S)]\}$, 
%$y\in 5\delta\text{-}\geod(x,a)$ et l'inégalité évidente $d(x,S)\leq d(x,a)\leq d(x,S)+N$ impliquent 
%$$(1-t)d(x,a)\leq 
%d(x,a)-td(x,S)\leq 
%d(x,a)-d(x,y)\leq 
%d(y,a)\leq d(x,a)-d(x,y)$$ $$+5\delta\leq d(x,a)-td(x,S)+2N+5\delta+1
%\leq (1-t)d(x,a)+3N+5\delta+1
%.$$
\cqfd

Le lemme suivant est un lemme général sur les espaces hyperboliques, qui généralise le lemme~\ref{contract1}. 

\begin{lem}\label{iter2}
  Soient $\alpha,\beta\in \N$ et $x,a,b,c\in X$ tels que $b\in \alpha\text{-}\geod(a,x)$, $c\in \beta\text{-}\geod(b,x)$. Alors
   $$c\in \big(\max(\alpha+2\beta-2d(b,c),\beta)+\delta\big)\text{-}\geod(a,x), $$  en particulier $c\in (\beta+\delta)\text{-}\geod(a,x)$ si  $d(b,c)\geq \frac{\alpha+\beta}{2}$.
%
%\noindent b)   $b\in (\alpha+\beta)\text{-}\geod(a,c)$. 
%
%\noindent c) 
%Si $d(b,c)\geq \alpha+\beta$, $d(a,c)\geq d(a,b)$. 
%
%
  \end{lem}

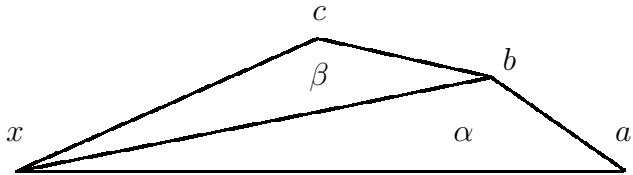
\ifx\JPicScale\undefined\def\JPicScale{1}\fi
\unitlength \JPicScale mm
\begin{picture}(90,28)(20,35)
\linethickness{0.3mm}
\put(30,40){\line(1,0){80}}
\linethickness{0.3mm}
\multiput(92.5,52.5)(0.17,-0.12){104}{\line(1,0){0.17}}
\linethickness{0.3mm}
\multiput(30,40)(0.6,0.12){104}{\line(1,0){0.6}}
\linethickness{0.3mm}
\multiput(30,40)(0.27,0.12){148}{\line(1,0){0.27}}
\linethickness{0.3mm}
\multiput(70,57.5)(0.54,-0.12){42}{\line(1,0){0.54}}
\put(30,45){\makebox(0,0)[cc]{$x$}}

\put(110,45){\makebox(0,0)[cc]{$a$}}

\put(95,55){\makebox(0,0)[cc]{$b$}}

\put(70,61){\makebox(0,0)[cc]{$c$}}

\put(70,52.5){\makebox(0,0)[cc]{$\beta$}}

\put(89,45){\makebox(0,0)[cc]{$\alpha$}}

\end{picture}

  \noindent{\bf Démonstration.} 
  En effet  $(H_{\de}^{\beta}(a,b,c,x))$ 
  s'écrit 
  $$d(a,c)\leq \max(d(a,b)-d(b,c), d(a,x)-d(x,c))+\beta+\delta,$$ d'où l'on déduit 
    $$d(a,c)+d(c,x)-d(a,x)\leq \max(d(a,b)-d(b,c)+d(c,x)-d(a,x),0)+\beta+\delta$$ et d'autre part  
  $$  d(a,b)+d(b,c)+d(c,x)\leq d(a,b)+d(b,x)+\beta\leq d(a,x)+\alpha+\beta.$$
% Comme $d(a,x)\leq d(a,c)+d(c,x)$, (\ref{triang-deux-fois}) implique b) et  enfin c) résulte immédiatement de b). 
  % s'écrit 
  %$$d(a,c)+d(b,x)\leq \max(d(a,b)+d(c,x), d(a,x)+d(b,c))+\delta,$$ d'où l'on déduit 
  %  $$d(a,c)+d(c,x)-d(a,x)$$ $$\leq \max(d(a,b)+2d(c,x)-d(b,x)-d(a,x),d(b,c)+d(c,x)-d(b,x))+\delta$$ $$ \leq \max(\alpha+2\beta-2d(b,c),\beta)+\delta.$$  
 %L'argument pour montrer b)   est encore plus simple et n'utilise même pas l'hyperbolicité : 
 %$d(a,b)+d(b,c)+d(c,x)\leq d(a,b)+d(b,x)+\beta\leq d(a,x)+\alpha+\beta$, donc par l'inégalité triangulaire à l'envers $d(a,b)+d(b,c)\leq d(a,c)+\alpha+\beta$. Enfin c) résulte immédiatement de b). 
  \cqfd

 Nous allons montrer qu'en appliquant $\partial$ et $h_{x}$ de fa\c con répétée à $e_{S}$ (pour $S\in \Delta \setminus  \{\emptyset\}$),  on reste à distance bornée de la réunion  des géodésiques reliant $x$ aux points de $S$. 

\begin{lem}\label{suite-S}
Soit $n\in \N$ et $S_{0},\dots,S_{n}$ une suite d'éléments de $\Delta$ telle que pour tout $i$, $e_{S_{i+1}}$ apparaît avec  un coefficient non nul 
dans $\partial(e_{S_{i}})$ ou dans $h_{x}(e_{S_{i}})$. Alors pour tout point  $y_{n}\in S_{n}$ n'appartenant ni à   $B(x,2\delta)$ ni à $S_{0}$, il existe $y_{0}\in S_{0}$ n'appartenant pas à $B(x,2\delta)$ tel que $y_{n}\in 4\delta\text{-}\geod(x,y_{0})$ et $d(y_{0},y_{n})> N-2\delta$. 
\end{lem}
\noindent{\bf Démonstration.} 
 D'après le lemme~\ref{ASx-geod-x} il existe $y_{n-1}\in S_{n-1}$,..., $y_{0}\in S_{0}$ n'appartenant pas à $B(x,2\delta)$ tels que $y_{i}=y_{i+1}$ ou bien $$y_{i+1}\in 3\delta\text{-}\geod(x,y_{i})\text{\  et\  }d(y_{i},y_{i+1})> N-2\delta. $$
 Par récurrence sur $i\in \{1,...,n\}$ on montre $y_{i}\in 4\de\tg(x,y_{0})$ en appliquant le lemme~\ref{iter2} à $(y_{0},y_{i-1},y_{i})$ 
 au lieu de $(a,b,c)$ et $(4\delta, 3\delta)$ au lieu de $(
 \alpha,\beta)$ et en utilisant le fait que $N-2\delta\geq 7\de/2$. 
  %
%On applique le lemme~\ref{iter2} avec $\alpha=4\delta$ et $\beta=3\delta$ : 
% comme on a supposé $N\geq 6\delta$,  on a $y_{i}\in 4\delta\text{-}\geod(x,a)$ pour tout $i\in \{0,\dots n\}$ (par récurrence ascendante sur $i$). 
 On suppose $N\geq 9\de$, ce qui est permis par  $(H_{N})$. 
 D'après le c) du    lemme~\ref{geod-comp-xabc} appliqué à $(y_{0},y_{i},y_{i+1})$ 
 au lieu de $(a,b,c)$ et à $(4\de,3\de)$ au lieu de $(\alpha,\beta)$, on a $d(y_{0},y_{i})\leq d(y_{0},y_{i+1})$ pour tout $i$, 
 et comme $d(y_{0},y_{i})>N-2\de$ si $i$ est le plus petit entier tel que $y_{i}\neq y_{0}$ on en déduit 
   $d(y_{0},y_{n})>N-2\delta$. 
\cqfd

Il résulte du lemme~\ref{hxt-eS-connaissance} que pour 
 $x\in X$ et $S\in \Delta\setminus \{\emptyset\}$, la connaissance de $h_{x}(e_{S})$ dépend  seulement de celle des points de  $\bigcup_{a\in S}5\delta\text{-}\geod(x,a)$. Grâce au lemme~\ref{suite-S}, on en déduit le corollaire suivant. 
 
 \begin{cor}\label{suite-S-connaissance}
Si $S_{0},\dots,S_{n}$ est une suite d'éléments de $\Delta$ telle que $S_{0}\neq \emptyset$ et que pour tout $i$, $e_{S_{i+1}}$ apparaît avec  un coefficient $c_{i}$ non nul 
dans $\partial(e_{S_{i}})$ ou dans $h_{x}(e_{S_{i}})$, alors $S_{0}\cup\dots\cup S_{n}\subset 
\bigcup_{a\in S_{0}}4\delta\text{-}\geod(x,a)$
et la connaissance de tous les $c_{i}$ 
ne dépend que de la connaissance des points de $\bigcup_{a\in S_{0}}9\delta\text{-}\geod(x,a)$ (et de leurs distances mutuelles).
\end{cor}
\noindent{\bf Démonstration.} Si $x,a,b,c$ sont tels que $b\in 4\delta\text{-}\geod(x,a)$ et $c\in 5\delta\text{-}\geod(x,b)$ alors $c\in 9\delta\text{-}\geod(x,a)$ par le a) du lemme~\ref{geod-comp-xabc}.  \cqfd

On rappelle que $\del h_{x}+h_{x}\del$ est inversible de $\oplus_{p=0}^{p_{\max}}\C^{(\Delta_{p})}$ dans lui-même et que l'on a posé 
$H_x=h_x(\del h_x+h_x\del)^{-1}$. 

 \begin{cor}\label{suite-S-connaissance-H-22oct09}
 Pour $p\in \{1,...,p_{\max}\}$ et $S\in \Delta_{p}$, le support de $H_{x}(e_{S})$ est inclus dans  
 $\bigcup_{a\in S}4\delta\text{-}\geod(x,a)$
et  $H_{x}(e_{S})$  
ne dépend que de la connaissance des points de $\bigcup_{a\in S}9\delta\text{-}\geod(x,a)$ (et de leurs distances mutuelles).
\end{cor}
\noindent{\bf Démonstration.} C'est une conséquence immédiate du corollaire~\ref{suite-S-connaissance}.  \cqfd

On note par ailleurs que $H_{x}(e_{\emptyset})=e_{x}$.

\subsection{Une conséquence de la construction de $H_{x}$}

Pour la construction du deuxième parametrix $u_{x}$, nous aurons besoin du lemme suivant qui est une variante du lemme 2.5.6 de~\cite{kkban}. 
On rappelle que pour  $f=\sum_{S\in \Delta_{p}}f(S)e_{S}$, on note $\mathrm{supp}(f)=\bigcup_{S\text{ tel que }f(S)\neq 0}S$.

\begin{lem}\label{remplissage}
Pour tout $p\in \{1,\dots ,p_{\mathrm{max}}\}$, il existe une application $G$-équivariante mais 
 non nécessairement linéaire $$\Phi _p: \{f\in \C^{(\Delta_p)}, 
\partial_{p-1}(f)=0\}\vers 
\C^{(\Delta_{p+1})}$$ 
\label{Phip-page}telle que  pour 
tout $f$ dans l'ensemble de départ de $\Phi_p$, 
\begin{itemize}
\item $\partial _p(\Phi_p(f))=f$, 
\item $\Phi_p(\lambda f)=\lambda\Phi_p(f)$ pour $\lambda\in \C$, 
\item  $\mathrm{supp}(\Phi_p(f))$
  est inclus dans 
  $\bigcup_{y,z\in \mathrm{supp}(f)}4\delta\text{-}\geod(y,z)$,
  %$\cup_{z\in \mathrm{supp}(f)}B(z,\max(N,\mathrm{diam}(\mathrm{supp}(f))))$,
\item $\Phi_{p}(f)$ ne dépend que de $f$ et de la connaissance des points de 
 \begin{gather}\label{cup-9degeod-suppf}\bigcup_{y,z\in \mathrm{supp}(f)}9\delta\text{-}\geod(y,z)\end{gather}
%$\cup_{z\in \mathrm{supp}(f)}B(z,\max(N,\mathrm{diam}(\mathrm{supp}(f))))$ 
et de leurs distances mutuelles, autrement dit si $f'$ est une autre fonction  une isométrie de \eqref{cup-9degeod-suppf} vers $\bigcup_{y,z\in \mathrm{supp}(f')}9\delta\text{-}\geod(y,z)$ envoyant $f$ sur $f'$ envoie 
$\Phi_{p}(f)$ sur $\Phi_{p}(f')$, 
\end{itemize}
et telle que   pour tout $R\in \R_{+}$ il existe $C\in \R_{+}$ ne dépendant que de $\de,K,N,R$,   tel que pour tout $f\in \C^{(\Delta_p)}$ avec 
$\partial_{p-1}(f)=0$ et $\mathrm{diam}(\mathrm{supp}(f))\leq R$ on ait  $\|\Phi_p(f)\|_{\ell^1(\Delta_{p+1})}\leq C\|f\|_{\ell^1(\Delta_{p})}$.  
  \end{lem}
\noindent{\bf Démonstration.}  On donne  une formule explicite pour $\Phi_{p}$, qui est la même que  dans~\cite{kkban}  : 
$$\Phi_{p}(f)=\frac{1}{\sharp(\mathrm{supp}(f))}
\sum _{z\in\mathrm{supp}(f) } H_{z}(f).$$ Pour montrer les propriétés de $\Phi_{p}$, on applique le corollaire~\ref{suite-S-connaissance-H-22oct09}. 
%Lipschitz en ponderant par $|f(x)|$? 
\cqfd

Pour la suite de l'article, on fixe de telles applications $\Phi_{p}$. 

\subsection{Construction d'un deuxième paramétrix $u_{x}$}

Le paramétrix $H_x$ ne nous convient pas car nous ne savons pas construire de normes sur $\bigoplus_{p=0}^{p_{\max}} \C^{(\Delta_p)}$ telles qu'il soit continu ainsi que $\del$ et que $\|\pi(g)\|\leq P(\ell(g))e^{s\ell(g)}$ pour un certain polynôme $P$. En fait nous avions déjà un problème dans~\cite{kkban}, puisqu'il avait fallu prendre l'opérateur $H_x$ associé à une valeur plus grande de $N$ et le modifier grâce au lemme 2.5.6. Ici le problème est encore aggravé 
%(vu les normes que nous devrons introduire ensuite, analogues à celles déjà construites dans le cas des arbres)
 :  pour $a\in X$ le coefficient dans $H_x(e_{a})$ d'une arête contenant $x$  dépend au moins de la connaissance de tous les points de $9\delta \text{-}\geod(x,a)$, et le nombre de possibilités pour les distances entre tous ces points est une  exponentielle en $d(x,a)$, comme le montre l'exemple suivant.
 
 \noindent{\bf Exemple.} Soit $\Gamma$ le produit libre de $\Z/3\Z$ avec $\Z$. On note $e_{1}$ et $e_{2}$ les générateurs de $\Z/3\Z$ et $\Z$. 
 Soit $\ell$ la longueur des mots sur $\Gamma$ et soit $X=\Gamma$ muni de la distance $d(a,b)=\ell(a^{-1}b)$. Alors pour tout $a\in X$, la réunion des 
 cycles de longueur $3$ passant par $a$ contient $ae_{1}$ et  $ae_{1}^{-1}$ mais ni $ae_{2}$ et  $ae_{2}^{-1}$. Donc la classe d'isométrie de 
 $9\delta \text{-}\geod(x,a)$ détermine l'emplacement de $e_{1}^{\mp 1}$ et
 $e_{2}^{\mp 1}$ dans l'écriture de $x^{-1}a$ comme mot réduit, et il y a $2^{d(x,a)}$ possibilités.

  Autrement dit $H_{x}$ 
 ne vérifie pas la condition (C2) (en revanche $H_{x}$ vérifie la condition (C1) grâce au lemme~\ref{nilp} et au corollaire~\ref{suite-S-connaissance-H-22oct09}). 
  %Heureusement cette nouvelle contrainte donne l'idée d'une construction différente pour un nouveau parametrix $u_{x}$. 

Nous allons maintenant construire un nouveau parametrix $u_{x}
$ vérifiant la condition (C2) (mais pas la condition (C1)). 
%
%Dans le sous-paragraphe suivant nous verrons que $u_{x}$ a aussi un inconvénient et le paramétrix définitif $J_x$ sera un  mélange de $H_x$ et $u_x$. 
%
%Construisons maintenant le nouveau parametrix $u_{x}$. 
Nous commen\c cons par définir, pour tous $x,a\in X$ et
$r\in\{1,...,d(x,a)-1\}$ une mesure $\mu_r(x,a)$ de masse $1$,
supportée par 
$\{y\in \delta\text{-}\geod(x,a), d(a,y)=r\}$.
% Si $y,y'$ appartiennent  à cet ensemble, $(H_{\delta}(y,a,y',x))$ montre que $d(y,y')\leq 2\delta$, donc cet ensemble est de diamètre $\leq 2\delta$ et son 
D'après le lemme~\ref{cardinal-tranche-geod}, le cardinal de cet ensemble   est borné
par une constante $C(\delta,K)$. Nous voulons que la propriété
suivante soit satisfaite : 
%\begin{gather}\label{fortement} \forall r>0,\forall \epsilon>0,
%\exists R>0, d(x,x')\leq r,  
%d(a,x)\geq R,\nonumber  \\ \text{ alors }
%\|\mu_r(x',a)-\mu_r(x,a)\|_1\leq 
%\epsilon .\end{gather} 
\begin{gather}\nonumber \forall \rho>0,\forall \epsilon>0,
\exists R>0, \text{ si }d(x,x')\leq \rho \text{ et }   
r\leq d(a,x)- R,  \\ \label{fortement}  \text{ alors }
\|\mu_r(x',a)-\mu_r(x,a)\|_1\leq 
\epsilon .\end{gather} 
 Dans la propriété ci-dessus, la condition que $d(a,x)-r $ est suffisamment grand est évidemment la plus faible possible, car $d(a,x)-r $ est essentiellement la distance entre $x$ et les points du support de $\mu_r(x,a)$ et cette distance doit nécessairement être grande pour que $\mu_r(x,a)$ varie peu en fonction de $x$. 

On rappelle que pour toute partie $A$ non vide de $X$, on note $\nu_A$ la mesure de
masse $1$ égale au produit par $(\sharp A)^{-1}$ de la fonction
caractéristique de $A$. 

Voici la formule : 
$$\mu_r(x,a)=\frac{1}{d(a,x)-r}\sum_{k=0}^{d(a,x)-r-1}\nu_{A_{x,a,r,k}}
\text{ \  où 
pour \ }k\leq d(a,x)-r \text{ \ on pose 
}$$
$$A_{x,a,r,k}=\{y, d(a,y)=r,\exists z\in \delta\text{-}\geod(x,a),
d(a,z)=r+k, y\in \geod(z,a)\}.$$ 
\label{murtxa}
\ifx\JPicScale\undefined\def\JPicScale{1}\fi
\unitlength \JPicScale mm
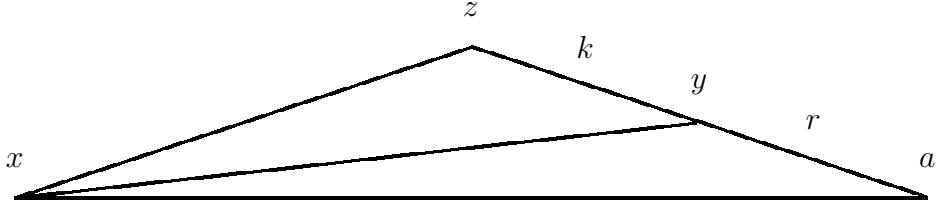
\begin{picture}(140,45)(20,20)
\linethickness{0.3mm}
\multiput(110,40)(0.36,-0.12){83}{\line(1,0){0.36}}
\linethickness{0.3mm}
\multiput(80,50)(0.36,-0.12){83}{\line(1,0){0.36}}
\linethickness{0.3mm}
\multiput(20,30)(0.36,0.12){167}{\line(1,0){0.36}}
\linethickness{0.3mm}
\put(20,30){\line(1,0){120}}
\linethickness{0.3mm}
\multiput(20,30)(1.08,0.12){83}{\line(1,0){1.08}}
\put(20,35){\makebox(0,0)[cc]{$x$}}

\put(80,55){\makebox(0,0)[cc]{$z$}}

\put(95,50){\makebox(0,0)[cc]{$k$}}

\put(110,45){\makebox(0,0)[cc]{$y$}}

\put(125,40){\makebox(0,0)[cc]{$r$}}

\put(140,35){\makebox(0,0)[cc]{$a$}}

\end{picture}

\noindent 
On a donc 
$$\mu_r(x,a)=
\int_{0}^{1}\mu_{r,t}(x,a)dt\text{\ \  où \ \ }
\mu_{r,t}(x,a)=
\nu_{A_{x,a,r,E(t(d(a,x)-r))}}.$$

D'après le a) du lemme~\ref{geod-comp-xabc},   $A_{x,a,r,k}\subset \{y\in \delta\text{-}\geod(x,a),
d(a,y)=r\}$ et d'après le lemme~\ref{cardinal-tranche-geod} le cardinal de ces ensembles est borné
par une constante de la forme $C(\delta,K)$. 
Pour tout $k\leq d(a,x)-r$, $A_{x,a,r,k}$ est non vide, puisque $X$ est géodésique.

\begin{lem}\label{comparaison-Axark-xx'}
a) Pour tout $k\in \{1,...,d(x,a)-r\}$ on a 
$A_{x,a,r,k}\subset A_{x,a,r,k-1}$.  

\noindent b) Pour $x'\in X$ et $k\in \N$ vérifiant $ d(x,x')+\delta\leq k \leq d(x,a)-r$ on a $A_{x,a,r,k}\subset A_{x',a,r,k-d(x,x')-\delta}$.
\end{lem}
\noindent{\bf Démonstration.} 
Montrons a). Soient $y,z$ tels que  \begin{gather}\label{cond-yz-22oct09}d(a,y)=r,\ \ z\in \delta\text{-}\geod(x,a),\ \ 
d(a,z)=r+k,\ \  y\in \geod(z,a).\end{gather} Il existe un  point $z'\in \geod(y,z)$ à distance $1$ de $z$. Alors $z'$ vérifie $$
d(a,z')=r+k-1, \ \ y\in \geod(z',a)$$ et $ z'\in \delta\text{-}\geod(x,a)$ par le a) du lemme~\ref{geod-comp-xabc}. 

\noindent Montrons b).  Soient  $y,z$ vérifiant (\ref{cond-yz-22oct09}).  Comme $d(y,z)=k\geq d(x,x')+\de$,  il existe  $z'\in \geod(y,z)$ à distance $d(x,x')+\delta$ de $z$. 

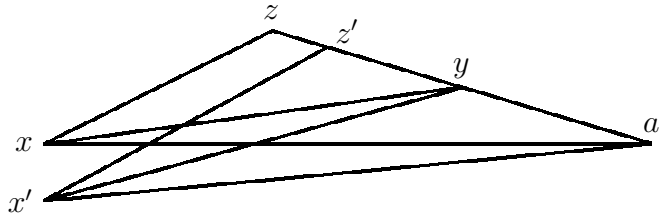
\ifx\JPicScale\undefined\def\JPicScale{1}\fi
\unitlength \JPicScale mm
\begin{picture}(110,35)(15,27)
\linethickness{0.3mm}
\put(30,40){\line(1,0){80}}
\linethickness{0.3mm}
\multiput(60,55)(0.4,-0.12){125}{\line(1,0){0.4}}
\linethickness{0.3mm}
\multiput(30,40)(0.87,0.12){63}{\line(1,0){0.87}}
\linethickness{0.3mm}
\multiput(30,40)(0.24,0.12){125}{\line(1,0){0.24}}
\linethickness{0.3mm}
\multiput(30,32.5)(0.22,0.12){169}{\line(1,0){0.22}}
\put(110,42.5){\makebox(0,0)[cc]{$a$}}

\put(85,50){\makebox(0,0)[cc]{$y$}}

\put(70,55){\makebox(0,0)[cc]{$z'$}}

\put(60,57.5){\makebox(0,0)[cc]{$z$}}

\put(27.5,40){\makebox(0,0)[cc]{$x$}}

\put(27,32.5){\makebox(0,0)[cc]{$x'$}}

\linethickness{0.3mm}
\multiput(30,32.5)(1.27,0.12){63}{\line(1,0){1.27}}
\linethickness{0.3mm}
\multiput(30,32.5)(0.44,0.12){125}{\line(1,0){0.44}}
\end{picture}

Alors 
 $z'$ vérifie $
d(a,z')=r+k-d(x,x')-\delta, y\in \geod(z',a)$ de fa\c con évidente et $z'\in \delta\text{-}\geod(x',a)$ 
 d'après le lemme~\ref{contract2}
 appliqué à $(x,x',a,z,z')$ au lieu de $(x,x',z,y,y')$ et $\de$ au lieu de $\epsilon$. \cqfd

La propriété (\ref{fortement}) résulte immédiatement du lemme suivant. 

\begin{lem}\label{lem-murtxx'-8j1848}
Il existe $C=C(\de,K)$ tel que pour $a,x,x'\in X$ et $r\in \N$ vérifiant 
$r< \min(d(a,x), d(a,x'))$ la mesure de l'ensemble des $t\in [0,1]$ tels que 
$\mu_{r,t}(x,a)\neq \mu_{r,t}(x',a)$ est $\leq \frac{Cd(x,x')}{d(a,x)-r}$.
\end{lem}
\noindent{\bf Démonstration.}
On utilise encore  l'astuce des ensembles
emboîtés.  Si $C_{0}=C(\delta,K)$ majore le cardinal des ensembles $A_{x,a,r,k}$ pour tous $a,x,r,k$, 
l'application décroissante $k\mapsto A_{x,a,r,k}$  de $\{0,...,d(x,a)-r\}$ dans l'ensemble des parties non vides de $X$ prend au plus $C_{0}$ valeurs, et pour chaque valeur, l'image inverse est un intervalle, et l'image inverse de cette même valeur par l'application 
$k\mapsto A_{x',a,r,k}$  de $\{0,...,d(x',a)-r\}$ dans l'ensemble des parties de $X$ est un autre intervalle dont les extrémités diffèrent au plus de $d(x,x')+\delta$ 
des extrémités du premier (il peut être vide si la longueur du  premier intervalle est inférieure ou égale à $2(d(x,x')+\delta)$).  
Comme $|d(a,x)-d(a,x')|\leq d(x,x')$, il en résulte que l'application décroissante $t\mapsto A_{x,a,r,E(t(d(a,x)-r))}$ prend au plus $C_{0}$ valeurs et pour chaque valeur l'image inverse est un intervalle, et l'image inverse de cette même valeur par l'application 
$t\mapsto A_{x',a,r,E(t(d(a,x')-r))}$  est un autre intervalle dont les extrémités diffèrent au plus de $\frac{2d(x,x')+\delta+1}{
d(a,x)-r}$ 
des extrémités du premier. 
 En effet pour $t,t'\in [0,1]$ vérifiant l'inégalité 
$|E(t(d(a,x)-r))-E(t'(d(a,x')-r))|\leq d(x,x')+\delta$ on a $|t-t'|\leq \frac{2d(x,x')+\delta+1}{
d(a,x)-r}$. 
On en déduit que la mesure  de l'ensemble des  $t\in [0,1]$ tels que 
$\mu_{r,t}(x,a)\neq \mu_{r,t}(x',a)$ 
est inférieure ou égale à $\frac{2C_{0}(2d(x,x')+\delta+1)}{
d(a,x)-r}$. 
Elle est donc inférieure ou égale à $\frac{2C_{0}(\de+3)d(x,x')}{
d(a,x)-r}$ car elle est nulle si $x=x'$. 
On peut prendre $C=2C_{0}(\de+3)$. 
\cqfd
 %$d(a,x)$ tend vers l'infini, alors que $r$ et $d(x,x')$ restent bornés. 

Cette formule pour $\mu_r(x,a)$ paraît artificiellement compliquée
mais son gros avantage pour nous est que $\mu_r(x,a)$ est une
certaine moyenne sur $k$ de $\nu_{A_{x,a,r,k}}$ et que
$A_{x,a,r,k}$ ne dépend que de la connaissance 
des points de $$\{a,x\}\cup\{y\in \delta\text{-}\geod(x,a), d(a,y)=r\}\cup\{z\in
\delta\text{-}\geod(x,a), d(a,z)=r+k\}$$ et de leurs distances mutuelles. 
%personnelle des
%points de $\{y\in \delta\text{-}\geod(x,a), d(a,y)=r\}$ et de la
%connaissance impersonnelle des points de 
%$\{z\in
%\delta\text{-}\geod(x,a), d(a,z)=r+k\}$ et des distances
%deux à deux entre les points de $\{y\in \delta\text{-}\geod(x,a), d(a,y)=r\}\cup\{z\in
%\delta\text{-}\geod(x,a), d(a,z)=r+k\}\cup\{a,x\}$. 
Or le nombre de ces points est borné par
une constante 
$C(\delta,K)$ et les distances entre ces points sont égales à $0$ ou
$d(x,a)$ ou  $d(x,a)-r$ ou $d(x,a)-r-k$ ou $r$ ou $k$ ou $r+k$ à
une constante $C(\delta)$ près. Donc lorsque $d(a,x),r$  et $k$ sont fixés,  le nombre de ces points et la donnée de leurs distances mutuelles n'admettent qu'un nombre fini  de possibilités, borné par $C(\delta,K)$. 

On pose enfin $\mu_{0,t}(x,a)=e_a$ et $\mu_{d(x,a),t}(x,a)=e_{x}$ pour tout $t\in [0,1]$. 

Nous mettons la remarque précédente sous forme d'un lemme, qui servira ensuite. 

\begin{lem}
Pour $x,a\in X$, $r\in \{0,\dots ,d(a,x)\}$ et $t\in [0,1]$, 
$\mu_{r,t}(x,a)$ est supporté par 
$$\{y\in \delta\text{-}\geod(x,a), d(a,y)=r\}$$ et 
dépend seulement de la connaissance 
des points de \begin{gather*}\{a,x\}\cup\{y\in \delta\text{-}\geod(x,a), d(a,y)=r\}\\ \cup\{z\in
\delta\text{-}\geod(x,a), d(a,z)=r+E(t(d(a,x)-r))\}\end{gather*} et de leurs distances mutuelles.
\end{lem}
\noindent{\bf Démonstration.} La preuve est incluse dans la remarque précédente. 
\cqfd

Nous allons maintenant construire un nouveau parametrix $u_{x}$ pour $\del$. 
Dans les trois pages qui suivent nous écrivons $u^{p}_x:\C^{(\Delta_{p-1})}\to
\C^{(\Delta_{p})}$ au lieu de $u_x$ pour rendre plus claire la construction, qui se fait par récurrence sur  $p$. Autrement dit on va construire des morphismes 
$$ \C^{(\Delta_0)}\vad{u_{x}^{1} }\C^{(\Delta_1)}\vad{u_{x}^{2} }\C^{(\Delta_2)}\dots \vad{u_{x}^{p_{\mathrm{max}}} }\C^{(\Delta_{p_{\mathrm{max}}})}$$ tels que 
$u^{p}_{x}\del+\del u^{p+1}_{x}=\Id_{\C^{(\Delta_{p})}}$ pour tout $p\in \{0,...,p_{\max}\}$. 

\label{def-u-v}
En fait on aura $$u_{x}^{p}=\int_{t\in [0,1]}u_{x,t}^{p}dt$$ et 
$u_{x,t}^{p}$ sera construit en même temps qu'une famille  indéxée par $r\in \N$ d'endomorphismes $(v_{x,r,t}^{p})_{p\geq 1}$ du complexe 
$$0\leftarrow \C^{(\Delta_0)}\vag{\partial }\C^{(\Delta_1)}\vag{\partial }\C^{(\Delta_2)}\dots \vag{\partial
  }\C^{(\Delta_{p_{\mathrm{max}}})}\leftarrow 0$$ qui consistent en gros à rapprocher de l'origine $x$ d'une longueur $r$, à l'aide des mesures $\mu_{r,t}(x,a)$. On aura
  $v_{x,0,t}^{p}=\Id_{\C^{(\Delta_{p})}}$ et quand  $r$ tend vers l'infini 
$v_{x,r,t}^{p}: \C^{(\Delta_{p})}\to \C^{(\Delta_{p})}$ tendra vers une limite $v_{x,\infty,t}^{p}$  égale à 
$\Id_{\C^{(\Delta_{0})}}$ si $p=0$, au morphisme de rang $1$ de $\C^{(\Delta_{1})}$ donné par $e_{a}\mapsto e_{x}$ si $p=1$ et à $0$ si $p>1$ (autrement dit $v_{x,\infty,t}^{p}$ est indépendant de $t$ et associé à la contraction de $X$ sur $x$). Pour $r\in \N^{*}$  et $p\geq 2$ on va construire 
$u_{x,r,t}^{p}:\C^{(\Delta_{p-1})}\to \C^{(\Delta_{p})}$ tel qu'en posant 
$u_{x,r,t}^{1}=0$ on ait 
$v_{x,r-1,t}^{p}-v_{x,r,t}^{p}=u_{x,r,t}^{p} \del+\del u_{x,r,t}^{p+1}$ pour tout $p$. 
On en déduira que \begin{gather}\label{parametrix-id-v-infini}\Id-v_{x,\infty,t}^{p}=\big(\sum_{r=1}^{\infty}u_{x,r,t}^{p}\big) \del+\del \big(\sum_{r=1}^{\infty} u_{x,r,t}^{p+1}\big) \text{ pour tout  } p.\end{gather}
On posera alors $u_{x,t}^{p}=\sum_{r=1}^{\infty}u_{x,r,t}^{p}$ pour $p\geq 2$ et 
$u_{x,t}^{1}(e_{\emptyset})=e_{x}$, de sorte que l'on aura 
$\del u^{p+1}_{x,t}+u_{x,t}^{p}\del=\Id_{\C^{(\Delta_{p})}}$ pour tout $p$
(grâce à \eqref{parametrix-id-v-infini} et au fait que $\del u_{x,t}^{1}=v_{x,\infty,t}^{0}$ et $ u_{x,t}^{1}\del =v_{x,\infty,t}^{1}$).

On passe maintenant à la construction des opérateurs $u_{x}^{p}$, 
$u_{x,t}^{p}$, $u_{x,r,t}^{p}$ et $v_{x,r,t}^{p}$.

On définit $u_{x}^{1}:\C^{(\Delta_0)}\to \C^{(\Delta_1)}$ par $u_{x}^{1}(e_{\emptyset})=e_{x}$ et on pose $u_{x,t}^{1}=u_{x}^{1}$. 

Voici la formule pour $u_x^{2}: \C^{(\Delta_1)}\vers \C^{(\Delta_2)}$: 
$$u_x^{2}(e_{a})=\int_{0}^{1}\sum
_{r=1}^{d(x,a)}u^{2}_{x,r,t}(e_{a})dt$$
où 
\begin{gather}\label{u2-Phi-mu-mu}u^{2}_{x,r,t}(e_{a})=
\Phi_1(\mu_{r-1,t}(x,a)-\mu_{r,t}(x,a)).\end{gather}
Il est évident que $\del u^{2}_x(e_a)=e_a-e_{x}=(1-u^{1}_{x}\del)(e_{a})$. 

Voici maintenant la formule pour $u_x^{3}:\C^{(\Delta_2)}\vers
\C^{(\Delta_3)}$. Par commodité nous étendons la fonction $r\mapsto
\mu_{r,t}(x,a)$ par $\mu_{r,t}(x,a)=e_{x}$ si
$r>d(x,a)$. La formule est, pour $\{a,b\}\in \Delta_{2}$ : 
\begin{gather}\nonumber u^{3}_{x}(e_a\wedge e_b)=\int_{0}^{1}\sum
_{r=1}^{\max(d(x,a),d(x,b))}u^{3}_{x,r,t}(e_a\wedge e_b)dt\text{ où}
\\ \nonumber 
u^{3}_{x,1,t}(e_a\wedge e_b)= \Phi_2\Big(e_a\wedge
e_b+\Phi_1\big(e_a-\mu_{1,t}(x,a)\big)-\Phi_1\big(e_b-\mu_{1,t}(x,b)\big) 
\\ \label{formule-u3ab-r=1} 
-\Phi_1\big(\mu_{1,t}(x,b)-\mu_{1,t}(x,a)\big)\Big),\text{ et pour }r>1,
 \\ \nonumber 
 u^{3}_{x,r,t}(e_a\wedge e_b)= \Phi_2\Big(
 \Phi_1\big(\mu_{r-1,t}(x,b)-\mu_{r-1,t}(x,a)\big)+
\Phi_1\big(\mu_{r-1,t}(x,a)-\mu_{r,t}(x,a)\big)
 \\ \label{formule-u3ab-r>1} 
  -
\Phi_1\big(\mu_{r-1,t}(x,b)-\mu_{r,t}(x,b)\big)-
 \Phi_1\big(\mu_{r,t}(x,b)-\mu_{r,t}(x,a)\big)
\Big).\end{gather}
On vérifie que $(\del u^{3}_{x}+u^{2}_{x}\del)(e_a\wedge e_b)=e_a\wedge e_b$. 

Afin de motiver la construction pour les plus grandes valeurs de $p$ on va réécrire les formules \eqref{u2-Phi-mu-mu}, \eqref{formule-u3ab-r=1}, \eqref{formule-u3ab-r>1} à l'aide des opérateurs $v_{x,r,t}^{1}:\C^{(\Delta_1)}\to \C^{(\Delta_1)}$ et $v_{x,r,t}^{2}:\C^{(\Delta_2)}\to \C^{(\Delta_2)}$ définis par 
\begin{gather*}
v_{x,0,t}^{1}(e_{a})=e_{a}, \ \ v_{x,r,t}^{1}(e_{a})=\mu_{r,t}(x,a), \\
v_{x,0,t}^{2}(e_{a}\wedge e_{b})=e_{a}\wedge e_{b}, \\
v_{x,r,t}^{2}(e_{a}\wedge e_{b})=\Phi_{1}(v_{x,r,t}^{1}(\del(e_{a}\wedge e_{b}))=\Phi_{1}(\mu_{r,t}(x,b)-\mu_{r,t}(x,a)).
\end{gather*}
On a alors 
\begin{gather*}u_{x,r,t}^{2}(e_{a})=\Phi_{1}(v_{x,r-1,t}^{1}(e_{a})-v_{x,r,t}^{1}(e_{a})), \\
u_{x,r,t}^{3}(e_{a}\wedge e_{b})=\Phi_{2}(v_{x,r-1,t}^{2}(e_{a}\wedge e_{b})-v_{x,r,t}^{2}(e_{a}\wedge e_{b})-u_{x,r,t}^{2}(\del(e_{a}\wedge e_{b}))).
\end{gather*}

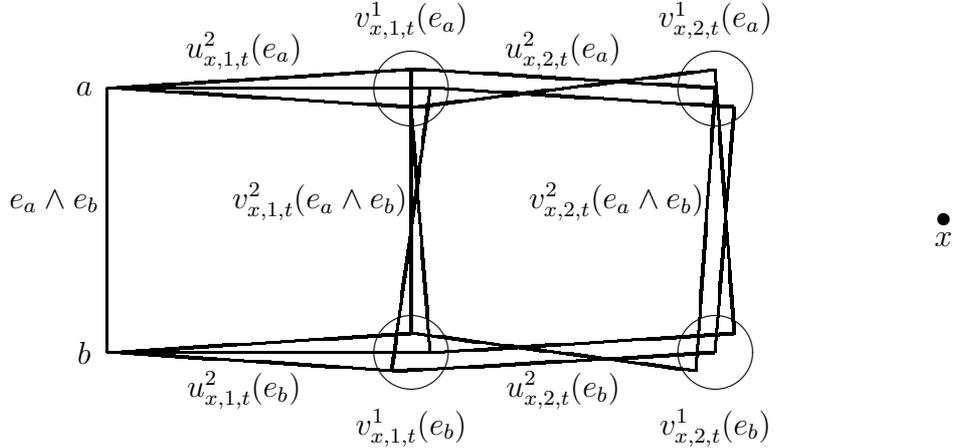
\ifx\JPicScale\undefined\def\JPicScale{1}\fi
\unitlength \JPicScale mm
\begin{picture}(135,80)(0,10)
\linethickness{0.3mm}
\put(15,30){\line(0,1){35}}
\linethickness{0.3mm}
\put(55,65){\circle{10}}

\linethickness{0.3mm}
\put(55,30){\circle{10}}

\linethickness{0.3mm}
\put(95,65){\circle{10}}

\linethickness{0.3mm}
\put(95,30){\circle{10}}

\linethickness{0.3mm}

\linethickness{0.3mm}

\linethickness{0.3mm}

\linethickness{0.3mm}

\multiput(54.99,67.5)(0.01,0){1}{\line(1,0){0.01}}
\multiput(54.96,67.5)(0.02,0){1}{\line(1,0){0.02}}
\multiput(54.92,67.5)(0.04,0){1}{\line(1,0){0.04}}
\multiput(54.87,67.49)(0.05,0){1}{\line(1,0){0.05}}
\multiput(54.81,67.49)(0.06,0){1}{\line(1,0){0.06}}
\multiput(54.74,67.48)(0.07,0){1}{\line(1,0){0.07}}
\multiput(54.65,67.48)(0.09,0.01){1}{\line(1,0){0.09}}
\multiput(54.55,67.47)(0.1,0.01){1}{\line(1,0){0.1}}
\multiput(54.44,67.47)(0.11,0.01){1}{\line(1,0){0.11}}
\multiput(54.32,67.46)(0.12,0.01){1}{\line(1,0){0.12}}
\multiput(54.18,67.45)(0.14,0.01){1}{\line(1,0){0.14}}
\multiput(54.04,67.44)(0.15,0.01){1}{\line(1,0){0.15}}
\multiput(53.88,67.43)(0.16,0.01){1}{\line(1,0){0.16}}
\multiput(53.71,67.42)(0.17,0.01){1}{\line(1,0){0.17}}
\multiput(53.52,67.41)(0.18,0.01){1}{\line(1,0){0.18}}
\multiput(53.33,67.4)(0.19,0.01){1}{\line(1,0){0.19}}
\multiput(53.13,67.38)(0.21,0.01){1}{\line(1,0){0.21}}
\multiput(52.91,67.37)(0.22,0.01){1}{\line(1,0){0.22}}
\multiput(52.68,67.36)(0.23,0.01){1}{\line(1,0){0.23}}
\multiput(52.44,67.34)(0.24,0.01){1}{\line(1,0){0.24}}
\multiput(52.19,67.32)(0.25,0.02){1}{\line(1,0){0.25}}
\multiput(51.93,67.31)(0.26,0.02){1}{\line(1,0){0.26}}
\multiput(51.66,67.29)(0.27,0.02){1}{\line(1,0){0.27}}
\multiput(51.38,67.27)(0.28,0.02){1}{\line(1,0){0.28}}
\multiput(51.09,67.26)(0.29,0.02){1}{\line(1,0){0.29}}
\multiput(50.79,67.24)(0.3,0.02){1}{\line(1,0){0.3}}
\multiput(50.48,67.22)(0.31,0.02){1}{\line(1,0){0.31}}
\multiput(50.16,67.2)(0.32,0.02){1}{\line(1,0){0.32}}
\multiput(49.83,67.18)(0.33,0.02){1}{\line(1,0){0.33}}
\multiput(49.49,67.16)(0.34,0.02){1}{\line(1,0){0.34}}
\multiput(49.14,67.13)(0.35,0.02){1}{\line(1,0){0.35}}
\multiput(48.79,67.11)(0.36,0.02){1}{\line(1,0){0.36}}
\multiput(48.42,67.09)(0.37,0.02){1}{\line(1,0){0.37}}
\multiput(48.05,67.07)(0.37,0.02){1}{\line(1,0){0.37}}
\multiput(47.66,67.04)(0.38,0.02){1}{\line(1,0){0.38}}
\multiput(47.27,67.02)(0.39,0.02){1}{\line(1,0){0.39}}
\multiput(46.88,66.99)(0.4,0.02){1}{\line(1,0){0.4}}
\multiput(46.47,66.97)(0.4,0.03){1}{\line(1,0){0.4}}
\multiput(46.06,66.94)(0.41,0.03){1}{\line(1,0){0.41}}
\multiput(45.64,66.92)(0.42,0.03){1}{\line(1,0){0.42}}
\multiput(45.22,66.89)(0.43,0.03){1}{\line(1,0){0.43}}
\multiput(44.78,66.86)(0.43,0.03){1}{\line(1,0){0.43}}
\multiput(44.35,66.83)(0.44,0.03){1}{\line(1,0){0.44}}
\multiput(43.9,66.81)(0.44,0.03){1}{\line(1,0){0.44}}
\multiput(43.45,66.78)(0.45,0.03){1}{\line(1,0){0.45}}
\multiput(43,66.75)(0.45,0.03){1}{\line(1,0){0.45}}
\multiput(42.54,66.72)(0.46,0.03){1}{\line(1,0){0.46}}
\multiput(42.07,66.69)(0.46,0.03){1}{\line(1,0){0.46}}
\multiput(41.61,66.66)(0.47,0.03){1}{\line(1,0){0.47}}
\multiput(41.13,66.63)(0.47,0.03){1}{\line(1,0){0.47}}
\multiput(40.66,66.6)(0.48,0.03){1}{\line(1,0){0.48}}
\multiput(40.18,66.57)(0.48,0.03){1}{\line(1,0){0.48}}
\multiput(39.69,66.54)(0.48,0.03){1}{\line(1,0){0.48}}
\multiput(39.21,66.51)(0.49,0.03){1}{\line(1,0){0.49}}
\multiput(38.72,66.48)(0.49,0.03){1}{\line(1,0){0.49}}
\multiput(38.23,66.45)(0.49,0.03){1}{\line(1,0){0.49}}
\multiput(37.73,66.42)(0.49,0.03){1}{\line(1,0){0.49}}
\multiput(37.24,66.39)(0.49,0.03){1}{\line(1,0){0.49}}
\multiput(36.74,66.36)(0.5,0.03){1}{\line(1,0){0.5}}
\multiput(36.25,66.33)(0.5,0.03){1}{\line(1,0){0.5}}
\multiput(35.75,66.3)(0.5,0.03){1}{\line(1,0){0.5}}
\multiput(35.25,66.27)(0.5,0.03){1}{\line(1,0){0.5}}
\multiput(34.75,66.23)(0.5,0.03){1}{\line(1,0){0.5}}
\multiput(34.25,66.2)(0.5,0.03){1}{\line(1,0){0.5}}
\multiput(33.75,66.17)(0.5,0.03){1}{\line(1,0){0.5}}
\multiput(33.26,66.14)(0.5,0.03){1}{\line(1,0){0.5}}
\multiput(32.76,66.11)(0.5,0.03){1}{\line(1,0){0.5}}
\multiput(32.27,66.08)(0.49,0.03){1}{\line(1,0){0.49}}
\multiput(31.77,66.05)(0.49,0.03){1}{\line(1,0){0.49}}
\multiput(31.28,66.02)(0.49,0.03){1}{\line(1,0){0.49}}
\multiput(30.79,65.99)(0.49,0.03){1}{\line(1,0){0.49}}
\multiput(30.31,65.96)(0.49,0.03){1}{\line(1,0){0.49}}
\multiput(29.82,65.93)(0.48,0.03){1}{\line(1,0){0.48}}
\multiput(29.34,65.9)(0.48,0.03){1}{\line(1,0){0.48}}
\multiput(28.87,65.87)(0.48,0.03){1}{\line(1,0){0.48}}
\multiput(28.39,65.84)(0.47,0.03){1}{\line(1,0){0.47}}
\multiput(27.93,65.81)(0.47,0.03){1}{\line(1,0){0.47}}
\multiput(27.46,65.78)(0.46,0.03){1}{\line(1,0){0.46}}
\multiput(27,65.75)(0.46,0.03){1}{\line(1,0){0.46}}
\multiput(26.55,65.72)(0.45,0.03){1}{\line(1,0){0.45}}
\multiput(26.1,65.69)(0.45,0.03){1}{\line(1,0){0.45}}
\multiput(25.65,65.67)(0.44,0.03){1}{\line(1,0){0.44}}
\multiput(25.22,65.64)(0.44,0.03){1}{\line(1,0){0.44}}
\multiput(24.78,65.61)(0.43,0.03){1}{\line(1,0){0.43}}
\multiput(24.36,65.58)(0.43,0.03){1}{\line(1,0){0.43}}
\multiput(23.94,65.56)(0.42,0.03){1}{\line(1,0){0.42}}
\multiput(23.53,65.53)(0.41,0.03){1}{\line(1,0){0.41}}
\multiput(23.12,65.51)(0.4,0.03){1}{\line(1,0){0.4}}
\multiput(22.73,65.48)(0.4,0.02){1}{\line(1,0){0.4}}
\multiput(22.34,65.46)(0.39,0.02){1}{\line(1,0){0.39}}
\multiput(21.95,65.43)(0.38,0.02){1}{\line(1,0){0.38}}
\multiput(21.58,65.41)(0.37,0.02){1}{\line(1,0){0.37}}
\multiput(21.21,65.39)(0.37,0.02){1}{\line(1,0){0.37}}
\multiput(20.86,65.37)(0.36,0.02){1}{\line(1,0){0.36}}
\multiput(20.51,65.34)(0.35,0.02){1}{\line(1,0){0.35}}
\multiput(20.17,65.32)(0.34,0.02){1}{\line(1,0){0.34}}
\multiput(19.84,65.3)(0.33,0.02){1}{\line(1,0){0.33}}
\multiput(19.52,65.28)(0.32,0.02){1}{\line(1,0){0.32}}
\multiput(19.21,65.26)(0.31,0.02){1}{\line(1,0){0.31}}
\multiput(18.91,65.24)(0.3,0.02){1}{\line(1,0){0.3}}
\multiput(18.62,65.23)(0.29,0.02){1}{\line(1,0){0.29}}
\multiput(18.34,65.21)(0.28,0.02){1}{\line(1,0){0.28}}
\multiput(18.07,65.19)(0.27,0.02){1}{\line(1,0){0.27}}
\multiput(17.81,65.18)(0.26,0.02){1}{\line(1,0){0.26}}
\multiput(17.56,65.16)(0.25,0.02){1}{\line(1,0){0.25}}
\multiput(17.32,65.14)(0.24,0.01){1}{\line(1,0){0.24}}
\multiput(17.09,65.13)(0.23,0.01){1}{\line(1,0){0.23}}
\multiput(16.87,65.12)(0.22,0.01){1}{\line(1,0){0.22}}
\multiput(16.67,65.1)(0.21,0.01){1}{\line(1,0){0.21}}
\multiput(16.48,65.09)(0.19,0.01){1}{\line(1,0){0.19}}
\multiput(16.29,65.08)(0.18,0.01){1}{\line(1,0){0.18}}
\multiput(16.12,65.07)(0.17,0.01){1}{\line(1,0){0.17}}
\multiput(15.96,65.06)(0.16,0.01){1}{\line(1,0){0.16}}
\multiput(15.82,65.05)(0.15,0.01){1}{\line(1,0){0.15}}
\multiput(15.68,65.04)(0.14,0.01){1}{\line(1,0){0.14}}
\multiput(15.56,65.03)(0.12,0.01){1}{\line(1,0){0.12}}
\multiput(15.45,65.03)(0.11,0.01){1}{\line(1,0){0.11}}
\multiput(15.35,65.02)(0.1,0.01){1}{\line(1,0){0.1}}
\multiput(15.26,65.02)(0.09,0.01){1}{\line(1,0){0.09}}
\multiput(15.19,65.01)(0.07,0){1}{\line(1,0){0.07}}
\multiput(15.13,65.01)(0.06,0){1}{\line(1,0){0.06}}
\multiput(15.08,65)(0.05,0){1}{\line(1,0){0.05}}
\multiput(15.04,65)(0.04,0){1}{\line(1,0){0.04}}
\multiput(15.01,65)(0.02,0){1}{\line(1,0){0.02}}
\multiput(15,65)(0.01,0){1}{\line(1,0){0.01}}

\multiput(15,65)(0.01,0){1}{\line(1,0){0.01}}
\multiput(15.01,65)(0.02,0){1}{\line(1,0){0.02}}
\multiput(15.04,65)(0.04,0){1}{\line(1,0){0.04}}
\multiput(15.08,65)(0.05,0){1}{\line(1,0){0.05}}
\multiput(15.13,65.01)(0.06,0){1}{\line(1,0){0.06}}
\multiput(15.19,65.01)(0.07,0){1}{\line(1,0){0.07}}
\multiput(15.26,65.02)(0.09,0.01){1}{\line(1,0){0.09}}
\multiput(15.35,65.02)(0.1,0.01){1}{\line(1,0){0.1}}
\multiput(15.45,65.03)(0.11,0.01){1}{\line(1,0){0.11}}
\multiput(15.56,65.03)(0.12,0.01){1}{\line(1,0){0.12}}
\multiput(15.68,65.04)(0.14,0.01){1}{\line(1,0){0.14}}
\multiput(15.82,65.05)(0.15,0.01){1}{\line(1,0){0.15}}
\multiput(15.96,65.06)(0.16,0.01){1}{\line(1,0){0.16}}
\multiput(16.12,65.07)(0.17,0.01){1}{\line(1,0){0.17}}
\multiput(16.29,65.08)(0.18,0.01){1}{\line(1,0){0.18}}
\multiput(16.48,65.09)(0.19,0.01){1}{\line(1,0){0.19}}
\multiput(16.67,65.1)(0.21,0.01){1}{\line(1,0){0.21}}
\multiput(16.87,65.12)(0.22,0.01){1}{\line(1,0){0.22}}
\multiput(17.09,65.13)(0.23,0.01){1}{\line(1,0){0.23}}
\multiput(17.32,65.14)(0.24,0.01){1}{\line(1,0){0.24}}
\multiput(17.56,65.16)(0.25,0.02){1}{\line(1,0){0.25}}
\multiput(17.81,65.18)(0.26,0.02){1}{\line(1,0){0.26}}
\multiput(18.07,65.19)(0.27,0.02){1}{\line(1,0){0.27}}
\multiput(18.34,65.21)(0.28,0.02){1}{\line(1,0){0.28}}
\multiput(18.62,65.23)(0.29,0.02){1}{\line(1,0){0.29}}
\multiput(18.91,65.24)(0.3,0.02){1}{\line(1,0){0.3}}
\multiput(19.21,65.26)(0.31,0.02){1}{\line(1,0){0.31}}
\multiput(19.52,65.28)(0.32,0.02){1}{\line(1,0){0.32}}
\multiput(19.84,65.3)(0.33,0.02){1}{\line(1,0){0.33}}
\multiput(20.17,65.32)(0.34,0.02){1}{\line(1,0){0.34}}
\multiput(20.51,65.34)(0.35,0.02){1}{\line(1,0){0.35}}
\multiput(20.86,65.37)(0.36,0.02){1}{\line(1,0){0.36}}
\multiput(21.21,65.39)(0.37,0.02){1}{\line(1,0){0.37}}
\multiput(21.58,65.41)(0.37,0.02){1}{\line(1,0){0.37}}
\multiput(21.95,65.43)(0.38,0.02){1}{\line(1,0){0.38}}
\multiput(22.34,65.46)(0.39,0.02){1}{\line(1,0){0.39}}
\multiput(22.73,65.48)(0.4,0.02){1}{\line(1,0){0.4}}
\multiput(23.12,65.51)(0.4,0.03){1}{\line(1,0){0.4}}
\multiput(23.53,65.53)(0.41,0.03){1}{\line(1,0){0.41}}
\multiput(23.94,65.56)(0.42,0.03){1}{\line(1,0){0.42}}
\multiput(24.36,65.58)(0.43,0.03){1}{\line(1,0){0.43}}
\multiput(24.78,65.61)(0.43,0.03){1}{\line(1,0){0.43}}
\multiput(25.22,65.64)(0.44,0.03){1}{\line(1,0){0.44}}
\multiput(25.65,65.67)(0.44,0.03){1}{\line(1,0){0.44}}
\multiput(26.1,65.69)(0.45,0.03){1}{\line(1,0){0.45}}
\multiput(26.55,65.72)(0.45,0.03){1}{\line(1,0){0.45}}
\multiput(27,65.75)(0.46,0.03){1}{\line(1,0){0.46}}
\multiput(27.46,65.78)(0.46,0.03){1}{\line(1,0){0.46}}
\multiput(27.93,65.81)(0.47,0.03){1}{\line(1,0){0.47}}
\multiput(28.39,65.84)(0.47,0.03){1}{\line(1,0){0.47}}
\multiput(28.87,65.87)(0.48,0.03){1}{\line(1,0){0.48}}
\multiput(29.34,65.9)(0.48,0.03){1}{\line(1,0){0.48}}
\multiput(29.82,65.93)(0.48,0.03){1}{\line(1,0){0.48}}
\multiput(30.31,65.96)(0.49,0.03){1}{\line(1,0){0.49}}
\multiput(30.79,65.99)(0.49,0.03){1}{\line(1,0){0.49}}
\multiput(31.28,66.02)(0.49,0.03){1}{\line(1,0){0.49}}
\multiput(31.77,66.05)(0.49,0.03){1}{\line(1,0){0.49}}
\multiput(32.27,66.08)(0.49,0.03){1}{\line(1,0){0.49}}
\multiput(32.76,66.11)(0.5,0.03){1}{\line(1,0){0.5}}
\multiput(33.26,66.14)(0.5,0.03){1}{\line(1,0){0.5}}
\multiput(33.75,66.17)(0.5,0.03){1}{\line(1,0){0.5}}
\multiput(34.25,66.2)(0.5,0.03){1}{\line(1,0){0.5}}
\multiput(34.75,66.23)(0.5,0.03){1}{\line(1,0){0.5}}
\multiput(35.25,66.27)(0.5,0.03){1}{\line(1,0){0.5}}
\multiput(35.75,66.3)(0.5,0.03){1}{\line(1,0){0.5}}
\multiput(36.25,66.33)(0.5,0.03){1}{\line(1,0){0.5}}
\multiput(36.74,66.36)(0.5,0.03){1}{\line(1,0){0.5}}
\multiput(37.24,66.39)(0.49,0.03){1}{\line(1,0){0.49}}
\multiput(37.73,66.42)(0.49,0.03){1}{\line(1,0){0.49}}
\multiput(38.23,66.45)(0.49,0.03){1}{\line(1,0){0.49}}
\multiput(38.72,66.48)(0.49,0.03){1}{\line(1,0){0.49}}
\multiput(39.21,66.51)(0.49,0.03){1}{\line(1,0){0.49}}
\multiput(39.69,66.54)(0.48,0.03){1}{\line(1,0){0.48}}
\multiput(40.18,66.57)(0.48,0.03){1}{\line(1,0){0.48}}
\multiput(40.66,66.6)(0.48,0.03){1}{\line(1,0){0.48}}
\multiput(41.13,66.63)(0.47,0.03){1}{\line(1,0){0.47}}
\multiput(41.61,66.66)(0.47,0.03){1}{\line(1,0){0.47}}
\multiput(42.07,66.69)(0.46,0.03){1}{\line(1,0){0.46}}
\multiput(42.54,66.72)(0.46,0.03){1}{\line(1,0){0.46}}
\multiput(43,66.75)(0.45,0.03){1}{\line(1,0){0.45}}
\multiput(43.45,66.78)(0.45,0.03){1}{\line(1,0){0.45}}
\multiput(43.9,66.81)(0.44,0.03){1}{\line(1,0){0.44}}
\multiput(44.35,66.83)(0.44,0.03){1}{\line(1,0){0.44}}
\multiput(44.78,66.86)(0.43,0.03){1}{\line(1,0){0.43}}
\multiput(45.22,66.89)(0.43,0.03){1}{\line(1,0){0.43}}
\multiput(45.64,66.92)(0.42,0.03){1}{\line(1,0){0.42}}
\multiput(46.06,66.94)(0.41,0.03){1}{\line(1,0){0.41}}
\multiput(46.47,66.97)(0.4,0.03){1}{\line(1,0){0.4}}
\multiput(46.88,66.99)(0.4,0.02){1}{\line(1,0){0.4}}
\multiput(47.27,67.02)(0.39,0.02){1}{\line(1,0){0.39}}
\multiput(47.66,67.04)(0.38,0.02){1}{\line(1,0){0.38}}
\multiput(48.05,67.07)(0.37,0.02){1}{\line(1,0){0.37}}
\multiput(48.42,67.09)(0.37,0.02){1}{\line(1,0){0.37}}
\multiput(48.79,67.11)(0.36,0.02){1}{\line(1,0){0.36}}
\multiput(49.14,67.13)(0.35,0.02){1}{\line(1,0){0.35}}
\multiput(49.49,67.16)(0.34,0.02){1}{\line(1,0){0.34}}
\multiput(49.83,67.18)(0.33,0.02){1}{\line(1,0){0.33}}
\multiput(50.16,67.2)(0.32,0.02){1}{\line(1,0){0.32}}
\multiput(50.48,67.22)(0.31,0.02){1}{\line(1,0){0.31}}
\multiput(50.79,67.24)(0.3,0.02){1}{\line(1,0){0.3}}
\multiput(51.09,67.26)(0.29,0.02){1}{\line(1,0){0.29}}
\multiput(51.38,67.27)(0.28,0.02){1}{\line(1,0){0.28}}
\multiput(51.66,67.29)(0.27,0.02){1}{\line(1,0){0.27}}
\multiput(51.93,67.31)(0.26,0.02){1}{\line(1,0){0.26}}
\multiput(52.19,67.32)(0.25,0.02){1}{\line(1,0){0.25}}
\multiput(52.44,67.34)(0.24,0.01){1}{\line(1,0){0.24}}
\multiput(52.68,67.36)(0.23,0.01){1}{\line(1,0){0.23}}
\multiput(52.91,67.37)(0.22,0.01){1}{\line(1,0){0.22}}
\multiput(53.13,67.38)(0.21,0.01){1}{\line(1,0){0.21}}
\multiput(53.33,67.4)(0.19,0.01){1}{\line(1,0){0.19}}
\multiput(53.52,67.41)(0.18,0.01){1}{\line(1,0){0.18}}
\multiput(53.71,67.42)(0.17,0.01){1}{\line(1,0){0.17}}
\multiput(53.88,67.43)(0.16,0.01){1}{\line(1,0){0.16}}
\multiput(54.04,67.44)(0.15,0.01){1}{\line(1,0){0.15}}
\multiput(54.18,67.45)(0.14,0.01){1}{\line(1,0){0.14}}
\multiput(54.32,67.46)(0.12,0.01){1}{\line(1,0){0.12}}
\multiput(54.44,67.47)(0.11,0.01){1}{\line(1,0){0.11}}
\multiput(54.55,67.47)(0.1,0.01){1}{\line(1,0){0.1}}
\multiput(54.65,67.48)(0.09,0.01){1}{\line(1,0){0.09}}
\multiput(54.74,67.48)(0.07,0){1}{\line(1,0){0.07}}
\multiput(54.81,67.49)(0.06,0){1}{\line(1,0){0.06}}
\multiput(54.87,67.49)(0.05,0){1}{\line(1,0){0.05}}
\multiput(54.92,67.5)(0.04,0){1}{\line(1,0){0.04}}
\multiput(54.96,67.5)(0.02,0){1}{\line(1,0){0.02}}
\multiput(54.99,67.5)(0.01,0){1}{\line(1,0){0.01}}

\linethickness{0.3mm}
\put(15,65){\line(1,0){42.5}}
\linethickness{0.3mm}
\multiput(15,65)(1.9,-0.12){21}{\line(1,0){1.9}}
\linethickness{0.3mm}
\multiput(15,30)(1.9,0.12){21}{\line(1,0){1.9}}
\linethickness{0.3mm}
\put(15,30){\line(1,0){42.5}}
\linethickness{0.3mm}
\multiput(15,30)(1.79,-0.12){21}{\line(1,0){1.79}}
\linethickness{0.3mm}
\multiput(52.5,27.5)(0.12,0.89){42}{\line(0,1){0.89}}
\linethickness{0.3mm}
\put(55,32.5){\line(0,1){35}}
\linethickness{0.3mm}
\multiput(55,62.5)(0.12,-1.55){21}{\line(0,-1){1.55}}
\linethickness{0.3mm}
\multiput(55,62.5)(0.95,0.12){42}{\line(1,0){0.95}}
\linethickness{0.3mm}
\multiput(57.5,65)(1.9,-0.12){21}{\line(1,0){1.9}}
\linethickness{0.3mm}
\multiput(55,67.5)(1.9,-0.12){21}{\line(1,0){1.9}}
\linethickness{0.3mm}
\multiput(95,67.5)(0.12,-1.67){21}{\line(0,-1){1.67}}
\linethickness{0.3mm}
\multiput(92.5,27.5)(0.12,1.79){21}{\line(0,1){1.79}}
\linethickness{0.3mm}
\multiput(95,30)(0.12,1.55){21}{\line(0,1){1.55}}
\linethickness{0.3mm}
\multiput(55,32.5)(0.89,-0.12){42}{\line(1,0){0.89}}
\linethickness{0.3mm}
\multiput(57.5,30)(1.9,0.12){21}{\line(1,0){1.9}}
\linethickness{0.3mm}
\multiput(52.5,27.5)(2.02,0.12){21}{\line(1,0){2.02}}
\put(12,65){\makebox(0,0)[cc]{$a$}}

\put(12,30){\makebox(0,0)[cc]{$b$}}

\put(8,50){\makebox(0,0)[cc]{$e_a\wedge e_b$}}

\put(43,50){\makebox(0,0)[cc]{$v_{x,1,t}^2(e_a\wedge e_b)$}}

\put(82,50){\makebox(0,0)[cc]{$v_{x,2,t}^2(e_a\wedge e_b)$}}

\put(33,25){\makebox(0,0)[cc]{$u_{x,1,t}^2(e_b)$}}

\put(75,25){\makebox(0,0)[cc]{$u_{x,2,t}^2(e_b)$}}

\put(55,20){\makebox(0,0)[cc]{$v_{x,1,t}^1(e_b)$}}

\put(95,20){\makebox(0,0)[cc]{$v_{x,2,t}^1(e_b)$}}

\put(95,74){\makebox(0,0)[cc]{$v_{x,2,t}^1(e_a)$}}

\put(55,74){\makebox(0,0)[cc]{$v_{x,1,t}^1(e_a)$}}

\put(33,70){\makebox(0,0)[cc]{$u_{x,1,t}^2(e_a)$}}

\put(75,70){\makebox(0,0)[cc]{$u_{x,2,t}^2(e_a)$}}

\put(125,47.5){\makebox(0,0)[cc]{$\bullet$}}

\put(125,45){\makebox(0,0)[cc]{$x$}}

\end{picture}

\noindent Le dessin illustre les calculs précédents. Les $4$ petits cercles contiennent les supports de $v_{x,r,t}^1(e_a)=\mu_{r,t}(x,a)$ et $v_{x,r,t}^1(e_b)=\mu_{r,t}(x,b)$ pour $r=1,2$ (en réalité ces supports ne seront pas disjoints mais on l'a supposé pour la clarté du dessin). Les six éléments $v_{x,r,t}^2(e_a\wedge e_b)$, $u_{x,r,t}^2(e_a)$, $u_{x,r,t}^2(e_b)$ (pour $r=1,2$) sont des combinaisons d'arêtes comme il est  indiqué sur le dessin. Enfin 
les deux éléments $u_{x,r,t}^3(e_a\wedge e_b)$ (pour $r=1,2$), qui ne sont pas représentés sur le dessin, sont des combinaisons de triangles remplissant les deux grands carrés.

En général on définit  $u^{p}_x:\C^{(\Delta_{p-1})}\to
\C^{(\Delta_p)}$ pour $p\in \{1,\dots,p_{\max}\}$ par récurrence ascendante sur $p$. La récurrence part de $p=1$ et fournit les formules ci-dessus pour $p=2,3$. 
On construit des applications linéaires  \begin{gather*}u_{x,r,t}^{p} : \C^{(\Delta_{p-1})}\vers
\C^{(\Delta_{p})}\text{ pour }r\in \N^{*}\text{ et }t\in [0,1] \\ \text{et }v^{p}_{x,r,t} : \C^{(\Delta_p)}\vers
\C^{(\Delta_{p})}\text{ pour }r\in \N\text{ et }t\in [0,1]\end{gather*} telles que 
$u_{x,r,t}^{p}(e_{a_{1}}\wedge ...\wedge e_{a_{p-1}})$ 
soit nul 
si $r>\max(d(x,a_{1}),...,d(x,a_{p-1}))$ et $p\geq 2$
 et 
$v_{x,r,t}^{p}(e_{a_{1}}\wedge ...\wedge e_{a_{p}})$
soit nul 
si $r\geq \max(d(x,a_{1}),...,d(x,a_{p}))$ et $p\geq 2$. 
Ces applications sont définies 
par $$v_{x,0,t}^{p}=\mathrm{Id}_{\C^{(\Delta_p)}} \text{ pour } p\geq 1, $$
 et pour $r\geq 1$, 
par récurrence ascendante sur $p$ à l'aide des formules suivantes : 
$$u^{1}_{x,r,t}=0, v^{1}_{x,r,t}(e_{a})=\mu_{r,t}(x,a) \text{ pour } r\geq 1     $$ pour initialiser et
\begin{gather*}v_{x,r,t}^{p}(e_{a_{1}}\wedge ...\wedge e_{a_{p}})=
\Phi_{p-1}(v_{x,r,t}^{p-1}(\del(e_{a_{1}}\wedge...\wedge e_{a_{p}}))) \text{ pour }r\geq 1 \text{ et } p\geq 2\\ 
u_{x,r,t}^{p}(e_{a_{1}}\wedge ...\wedge e_{a_{p-1}})=
\Phi_{p-1}\big(v_{x,r-1,t}^{p-1}(e_{a_{1}}\wedge ...\wedge e_{a_{p-1}})-v_{x,r,t}^{p-1}(e_{a_{1}}\wedge ...\wedge e_{a_{p-1}})
\\ 
-u_{x,r,t}^{p-1}(\del(e_{a_{1}}\wedge ...\wedge e_{a_{p-1}}))\big)\text{ pour }r\geq 1 \text{ et } p\geq 2.\end{gather*}

Notons que ces formules impliquent les deux égalités 
$$\del\circ v_{x,r,t}^{p}=v_{x,r,t}^{p-1}\circ \del \text{ pour } p\geq 2\text{ et } r\geq 0\text{ et }$$
$$\del\circ u_{x,r,t}^{p}+u_{x,r,t}^{p-1}\circ\del
=v_{x,r-1,t}^{p-1}-v_{x,r,t}^{p-1} \text{ pour } p\geq 2\text{ et } r\geq 1.$$

Pour $p\in \{2,\dots, p_{\max}\}$, on définit $u_{x}:\C^{(\Delta_{p-1})}\vers
\C^{(\Delta_{p})}$ par $$u_{x}=\int_{0}^{1}\sum_{r=1}^{\infty} u^{p}_{x,r,t}dt $$ (la somme sur $r$ est toujours finie, plus précisément  pour $u_{x}(e_{a_{1}}\wedge ...\wedge e_{a_{p-1}})$  la somme s'arrête à $\max(d(x,a_{1}),...,d(x,a_{p-1}))$). 
Alors $u_{x}$ et $\del$
sont des endomorphismes de degrés $1$ et $-1$ de 
$\bigoplus_{p=0}^{p_{\max}} \C^{(\Delta_p)}$ et on a $\del u_{x} + u_{x}
\del=1$.  

\noindent{\bf Remarque.} Comme les applications $\Phi_{i}$ ne sont pas linéaires, les formules précédentes pour $u_{x,r,t}$ en fonction de $\mu_{r,t}(x,a)$ ne donnent pas de formule pour $u_{x}$ en fonction de $\mu_{r}(x,a)$. Autrement dit $\mu_{r}(x,a)$ ne sert à rien. Plus loin pour montrer l'équivariance à compact près des opérateurs nous n'utiliserons pas  la propriété (\ref{fortement}) mais le lemme~\ref{lem-murtxx'-8j1848}. 

Le gros avantage du paramétrix $u_{x}$ sur $H_{x}$ apparaît dans la proposition suivante, où l'on voit  qu'il vérifie la condition (C2). Une propriété  semblable sera vraie aussi pour le paramétrix définitif $J_{x}$ et jouera un rôle crucial dans  la construction  de normes  telles que $J_{x}$ soit continu. 

\begin{prop}\label{supp-uxrt}
Pour $p\in \{2,\dots,p_{\max}\},t\in [0,1]$, $S=\{a_{1},\dots, a_{p-1}\}\in \Delta_{p-1}$,  
$$\mathrm{supp}(u_{x,r,t}(e_{a_{1}}\wedge ...\wedge e_{a_{p-1}}))\subset 
\bigcap_{a\in S ,
y\in \delta\text{-}\geod(x,a),d(y,a)=r} B(y,N+5p\delta)$$
et $u_{x,r,t}(e_{a_{1}}\wedge ...\wedge e_{a_{p-1}})$ ne dépend que de la connaissance des points de l'ensemble 
\begin{gather}\nonumber S\cup\{x\}\cup\bigcap_{a\in S ,
y\in \delta\text{-}\geod(x,a),d(y,a)=r} B(y,N+5\delta+5p\delta)\\ \nonumber \cup  \bigcup_{a\in S }\{y\in \delta\text{-}\geod(x,a),d(y,a)=r+E(t(d(x,a)-r))\} \\  \label{ens-3.28-2011}\cup  \bigcup_{a\in S }\{y\in \delta\text{-}\geod(x,a),d(y,a)=r-1+E(t(d(x,a)-(r-1)))\}\end{gather} 
et des distances  entre ces points. 
\end{prop}
Le nombre de  points dans l'ensemble ci-dessus  est  borné par une constante $C(\delta,K,N)$ et les distances entre ces points   sont elles-mêmes déterminées  par la connaissance de $r$, $t$ et $d(x,S)$ à une constante $C(\delta,K,N)$ près. Donc le nombre total de possibilités pour le nombre de ces points et leurs distances mutuelles (c'est-à-dire  le nombre de classes d'équivalence de \eqref{ens-3.28-2011} par les isométries qui préservent ses parties $S$ et $\{x\}$)  est borné par $C(\delta,K,N)$. 

Pour montrer la proposition~\ref{supp-uxrt} nous aurons  besoin de trois lemmes. 

\begin{lem}\label{abxyz-trapeze}
Soient $\alpha,\beta\in \N$ et $x,a,b,y,z\in X$ vérifiant $y\in \alpha\text{-}\geod(x,a) $ et $z\in \beta\text{-}\geod(x,b) $. Alors 
%\begin{itemize}\item a) 
 $$d(y,z)\leq d(a,b)+|d(a,y)-d(b,z)|+\alpha+\beta+2\delta .$$
%\item b) $d(y,z)\leq d(a,b)+|d(x,y)-d(x,z)|+\alpha+\beta+2\delta$. 
%\end{itemize}

\end{lem}

\ifx\JPicScale\undefined\def\JPicScale{1}\fi
\unitlength \JPicScale mm
\begin{picture}(125,65)(20,8)
\linethickness{0.3mm}
\multiput(30,40)(1.08,-0.12){83}{\line(1,0){1.08}}
\linethickness{0.3mm}
\multiput(30,40)(1.08,0.12){83}{\line(1,0){1.08}}
\linethickness{0.3mm}
\multiput(80,60)(0.48,-0.12){83}{\line(1,0){0.48}}
\linethickness{0.3mm}
\multiput(30,40)(0.3,0.12){167}{\line(1,0){0.3}}
\linethickness{0.3mm}
\multiput(30,40)(0.3,-0.12){167}{\line(1,0){0.3}}
\linethickness{0.3mm}
\multiput(80,20)(0.48,0.12){83}{\line(1,0){0.48}}
\put(25,40){\makebox(0,0)[cc]{$x$}}

\put(80,65){\makebox(0,0)[cc]{$y$}}

\put(125,50){\makebox(0,0)[cc]{$a$}}

\put(75,50){\makebox(0,0)[cc]{$\alpha$}}

\put(75,30){\makebox(0,0)[cc]{$\beta$}}

\put(125,30){\makebox(0,0)[cc]{$b$}}

\put(80,15){\makebox(0,0)[cc]{$z$}}

\end{picture}

\noindent{\bf Démonstration.} 
Par 
$(H_{\delta}^{\alpha}(b,a,y,x))$, on a 
$$d(b,y)\leq 
%\max(d(a,b)-d(a,x)+d(y,x),d(b,x)-d(a,x)+d(a,y))+\delta $$ $$
 \max(d(b,a)-d(a,y),d(b,x)-d(x,y))+\alpha+\delta$$ et par $(H_{\delta}^{\beta}(y,b,z,x)) $ on a 
$$d(y,z)
%\leq  \max(d(y,b)-d(b,x)+d(x,z),d(y,x)-d(b,x)+d(b,z))+\delta $$ $$
\leq \max(d(y,b)-d(b,z),d(y,x)-d(x,z))+\beta+\delta $$ $$
\leq \max\big(d(a,b)-d(a,y)-d(b,z)+\alpha+\beta+2\delta, d(b,x)-d(x,y)-d(b,z)+\alpha+\beta+2\delta, $$ $$  d(x,y)-d(x,z)+\beta+\delta\big)
\leq d(a,b)+
%\min(|d(a,y)-d(b,z)|,|d(x,y)-d(x,z)|)+
|d(a,y)-d(b,z)|+
\alpha+\beta+2\delta $$ 
car \begin{gather*}d(b,x)-d(x,y)-d(b,z)\leq d(a,b)+d(a,x)-d(x,y)-d(b,z)\\ \leq d(a,b)+d(a,y)-d(b,z) \leq 
d(a,b)+|d(a,y)-d(b,z)|,\end{gather*} 
%$$d(b,x)-d(b,z)-d(x,y)\leq d(b,x)-d(b,z)-d(x,z)+|d(x,y)-d(x,z)|$$ $$\leq |d(x,y)-d(x,z)|
% $$
et 
\begin{gather*}d(x,y)-d(x,z)\leq d(x,a)-d(a,y)+\alpha-d(x,b)+d(b,z)\\ \leq d(a,b)+|d(a,y)-d(b,z)|
+\alpha. \end{gather*}
\cqfd

\begin{lem}\label{boules-geod}
Soit $R,\alpha\in \R_{+}$ et $x,y,z,t\in X$ vérifiant 
$y,z\in B(x,R)$ et $t\in \alpha\text{-}\geod(y,z)$. Alors $t\in B(x,R+\alpha+\delta)$. 
\end{lem}
\noindent{\bf Démonstration.} 
Par $(H_{\delta}^{\alpha}(x,y,t,z))$, on a 
$$d(x,t)\leq \max(d(x,y)-d(y,t),d(x,z)-d(z,t))+\alpha+\delta\leq R+\alpha+\delta. $$  \cqfd

\begin{lem}\label{supp-remplissage}
Pour tout $p\in \{1,\dots ,p_{\mathrm{max}}\}$, $y\in X$, $R\in \R_{+}$, 
et $f\in \C^{(\Delta_p)}$ telle que $\mathrm{supp}(f)\subset B(y,R)$, on a 
$\mathrm{supp}(\Phi_p(f))\subset  B(y,R+5\delta)$ et  $\Phi_{p}(f)$ ne dépend que de $f$ et de la connaissance des points de 
  $B(y,R+10\delta)$ 
et de leurs distances mutuelles.
\end{lem}
\noindent{\bf Démonstration.} 
Cela résulte des lemmes~\ref{remplissage} et~\ref{boules-geod}. \cqfd

\noindent{\bf Démonstration de la  proposition~\ref{supp-uxrt}.}  
La réunion des supports des mesures $$\mu_{r,t}(x,a_{1}), ..., \mu_{r,t}(x,a_{p-1}), 
\mu_{r-1,t}(x,a_{1}),..., \mu_{r-1,t}(x,a_{p-1})$$ est incluse  dans 
\begin{gather*} \bigcup_{a\in S }\{y\in \delta\text{-}\geod(x,a),d(y,a)=r\} \\ 
\cup  \bigcup_{a\in S }\{y\in \delta\text{-}\geod(x,a),d(y,a)=r-1\}. \end{gather*} 
Le %a) du 
lemme~\ref{abxyz-trapeze} 
(avec $\alpha=\beta=\delta$, $a,b\in S$, $y\in \delta\text{-}\geod(x,a)$ à distance $r$ ou $r-1 $ de $a$, et $z\in \delta\text{-}\geod(x,b)$ à distance $r$ ou $r-1 $ de $b$) 
implique que cet ensemble est de diamètre $\leq N+4\delta+1\leq N+5\delta$ et qu'il est donc inclus dans $$\bigcap_{a\in S ,
y\in \delta\text{-}\geod(x,a),d(y,a)=r} B(y,N+5\delta).$$ Grâce au lemme~\ref{supp-remplissage}, on montre par récurrence sur $p$ l'assertion de la proposition~\ref{supp-uxrt} en même temps que l'assertion analogue pour $v_{x,r,t}$, à savoir que pour  
$p\in \{2,\dots,p_{\max}\},t\in [0,1]$, $S=\{a_{1},\dots a_{p}\}\in \Delta_{p}$,  
$$\mathrm{supp}(v_{x,r,t}(e_{a_{1}}\wedge ...\wedge e_{a_{p}}))\subset 
\bigcap_{a\in S,
y\in \delta\text{-}\geod(x,a),d(y,a)=r\text{\ ou\ }r+1} B(y,N+5p\delta)$$
et $v_{x,r,t}(e_{a_{1}}\wedge ...\wedge e_{a_{p}})$ ne dépend que de la connaissance des points de l'ensemble 
\begin{gather*}S\cup\{x\}\cup\bigcap_{a\in S,
y\in \delta\text{-}\geod(x,a),d(y,a)=r\text{\ ou\ }r+1} B(y,N+5\delta+5p\delta)\\ \cup  \bigcup_{a\in S}\{y\in \delta\text{-}\geod(x,a),d(y,a)=r+E(t(d(x,a)-r))\} \end{gather*} 
et des distances  entre tous ces points. \cqfd

\begin{lem}\label{28dec1029}
Il existe $C=C(\de,K,N)$ telle que pour tout $p\in \{1,...,p_{\max}\}$ et pour tout $S\in \Delta_{p-1}$, $r\in \N$ et $t\in [0,1]$ on ait $\|u_{x,r,t}(e_{S})\|_{\ell^{1}(\Delta_{p})}\leq C$.  
\end{lem}
 
 \noindent{\bf Démonstration.} Cela résulte de la dernière assertion du lemme~\ref{remplissage} et du fait que dans les formules pour $u_{x,r,t}$ on applique les $\Phi_{i}$  à des fonctions  dont le support est de diamètre $\leq 2(N+5\de p_{\max})$ d'après la démonstration de la proposition~\ref{supp-uxrt}. \cqfd 
 
%Le point essentiel de cette construction de $u_{x}$ est que 
%$u_{x,r,t}(e_{a_{1}}\wedge ...\wedge e_{a_{p-1}})$ ne dépend que de la connaissance de $\mu_{r,t}(x,a_{1})$,...,$\mu_{r,t}(x,a_{p-1})$, 
%$\mu_{r-1,t}(x,a_{1})$,..., $\mu_{r-1,t}(x,a_{p-1})$ et des points de la réunion des boules de rayon $N+20\delta$ et de centres les points de la réunion des supports de ces mesures, ainsi que des distances entre tous ces points.  En effet le diamètre de la réunion des supports de ces mesures est inférieur ou égal à $N+20\delta$ car le support de $\mu_{r,t}(x,a_{i})$ est inclus dans $\{y\in \delta\text{-}\geod(x,a_{i}),d(y,a_{i})=r\}$ et de même avec $r-1$ au lieu de $r$, et le diamètre de $S$ est inférieur ou égal à $N$. 
%Pour conclure $u_{x,r,t}(e_{a_{1}}\wedge ...\wedge e_{a_{p-1}})$ ne dépend que de la connaissance des points de l'ensemble 
%$$S\cup\{x\}\cup\bigcup_{i\in \{1,...,p-1\},
%y\in \delta\text{-}\geod(x,a_{i}),d(y,a_{i})=r} B(y,N+20\delta)$$ $$\cup  \bigcup_{i\in \{1,...,p-1\},y\in \delta\text{-}\geod(x,a_{i}),d(y,a_{i})=r+E(t(d(x,a_{i})-r))} B(y,N+20\delta)$$ 
%et des distances mutuelles entre tous ces points. 

\subsection{Construction du paramétrix définitif $J_{x}$}

Nous allons introduire un opérateur $\tilde H_x$, sorte de troncature de $H_x$, et $J_x$ sera donné par la formule 
$J_x=\tilde H_x+u_x(1-\del \tilde H_x -\tilde H_x \del)$. 

En effet le gros inconvénient de $u_{x}$ est que l'on sait seulement, par la proposition~\ref{supp-uxrt} et le lemme~\ref{xx'yy'zz'}, que, pour $S\in \Delta\setminus\{\emptyset\}$,  $u_{x}(e_{S})$ est supporté  par $\bigcap_{a\in S}\big( 2(N+5p_{\max}\delta)+\delta\big)\text{-}\geod(x,a)$. 
Autrement dit  $u_{x}$ ne vérifie pas la condition (C1) et la  seule chose que l'on puisse affirmer est que $u_{x}$ 
``n'éloigne pas de $x$ de plus que de $N+5p_{\max}\delta+\delta$''.
Cela compliquerait beaucoup la construction de normes pour lesquelles $\partial$ et $u_{x}$ soient continus, voire la rendrait impossible. 
Au contraire nous verrons dans la proposition~\ref{recap-supp-connaiss-H-uK} que l'opérateur $J_{x}$ vérifie les conditions (C1) et (C2). L'idée pour montrer que $J_{x}$ vérifie (C1) est que $(1-\del \tilde H_x -\tilde H_x \del)$ rapproche davantage  de $x$ que $u_{x}$ n'en éloigne. L'idée pour montrer que $J_{x}$ vérifie (C2) est que, contrairement  à $H_{x}$, $\tilde H_{x}$ vérifie (C2) grâce au fait que c'est une troncature.

On peut écrire formellement $H_x$ sous la forme $$H_x=\sum_{q=1}^{+\infty}h_x(1-(\del h_x+h_x\del))^{q-1}.$$ Cette somme est infinie (mais pour chaque $S\in \Delta$, $H_x(e_{S})$ est donné par une somme finie, grâce au lemme~\ref{nilp}). 

On 
fixe un  entier $Q$ tel que 
$$(H_{Q}) : Q \text{ soit assez grand en fonction de }\de,K,N.$$
Dans la suite nous utiliserons un nombre {\it fini} de fois l'inégalité $Q\geq C$ avec $C$ de la forme $C(\de,K,N)$. 
\label{hyp-HQ}

On 
pose $$\tilde H_x=\sum_{q=1}^{Q}h_x(1-(\del h_x+h_x\del))^{q-1}
\text{\ \ et\ \ } K_{x}=1-(\del \tilde H_x+\tilde H_x \del)
.$$ Il est clair que 
$K_{x}=
(1-(\del h_x+h_x\del))^{Q}$. 
Le lemme suivant résume les propriétés de $\tilde H_x$ qui nous serviront ensuite. La propriété principale est 
que $K_{x}$ rapproche strictement de l'origine. 
On pose $K_{x,0}=1$ et pour $q\in \{1,...,Q\}$ et $(t_{1},\dots ,t_{q})\in [0,1]^{q}$, on note 
\begin{gather*}K_{x,q,(t_{1},\dots ,t_{q})}=(1-(\del h_{x,t_{q}}+h_{x,t_{q}}\del))...(1-(\del h_{x,t_{1}}+h_{x,t_{1}}\del))\\ 
\text{et\ \ \ \ }\tilde H_{x,q,(t_{1},\dots ,t_{q})}=h_{x,t_{q}}K_{x,q-1,(t_{1},\dots ,t_{q-1})}\\ 
=h_{x,t_{q}}(1-(\del h_{x,t_{q-1}}+h_{x,t_{q-1}}\del))...(1-(\del h_{x,t_{1}}+h_{x,t_{1}}\del))\end{gather*} 
\label{Hxqttt} de sorte que 
\begin{gather*}\tilde H_{x}=\sum_{q=1}^{Q}
\int_{(t_{1},\dots ,t_{q})\in [0,1]^{q}}\tilde H_{x,q,(t_{1},\dots ,t_{q})} dt_{1}\dots dt_{q}\\ 
\text{et\ \ \ } K_{x}=
\int_{(t_{1},\dots ,t_{Q})\in [0,1]^{Q}}  K_{x,Q,(t_{1},\dots ,t_{Q})}  dt_{1}\dots dt_{Q}.\end{gather*}

\noindent {\bf Notation.} Dans toute la suite de l'article, si $A$ est une partie de $X$ et $r\in \N$, on notera 
\begin{gather}\label{def-B-30dec1647}B(A,r)=\{y\in X, d(y,A)\leq r\}=\bigcup_{a\in A}B(a,r).\end{gather}

\begin{lem}\label{support-connaissance-tildeH}
1) 
Pour $q\in \{1,...,Q\}$, 
$(t_{1},\dots ,t_{q})\in [0,1]^{q}$ et $S\in \Delta\setminus\{\emptyset\}$,

\noindent a)  $\tilde H_{x,q,(t_{1},\dots ,t_{q})}(e_{S})$
est une combinaison de $e_{T}$ où $T$ vérifie : 
$$T\subset  S\cup \Big(B(x,2\de)\cap B(S,qN)\Big)$$ $$
\cup\bigcup_{a\in S,a\not\in B(x,2\de)}\big\{y\in 4\delta\text{-}\geod(x,a), d(y,a)\in
]N-2\de
%[\max(N-2\de,(q-1)\frac{N-6\delta}{p_{\max}}-2N)
, q N]\big\},$$

\noindent b)   $\tilde H_{x,q,(t_{1},\dots ,t_{q})}(e_{S})$ ne dépend que de la connaissance des points de
\begin{gather}\nonumber
B(x,7\de)\cup\bigcup_{a\in S}\{y\in 4\delta\text{-}\geod(x,a), d(y,a)\leq qN\}\\ \label{conn-tildeH-form-22oct09}\cup 
\bigcup_{a\in S,i\in \{1,...,q\}}\{y\in 9\delta\text{-}\geod(x,a), |d(x,y)-t_{i}d(x,a)|\leq (q+2)N+1\}\end{gather}
et des distances entre ces points.

\noindent 2)   Pour $(t_{1},\dots ,t_{Q})\in [0,1]^{Q}$  et $S\in \Delta\setminus\{\emptyset\}$, 

\noindent a) 
$K_{x,Q,(t_{1},\dots ,t_{Q})}(e_{S})$ 
est une combinaison de $e_{T}$ où $T$ vérifie : 
$$T\subset \bigcup_{a\in S}\{y\in 4\delta\text{-}\geod(x,a), d(y,a)\in [Q\frac{N-6\delta}{p_{\max}}-2N,QN]\}.$$

\noindent b) 
$K_{x,Q,(t_{1},\dots ,t_{Q})}(e_{S})$ ne dépend que de la connaissance des points de
\begin{gather*}B(x,7\de)\cup\bigcup_{a\in S}\{y\in 4\delta\text{-}\geod(x,a), d(y,a)\leq QN\}\\ \cup 
\bigcup_{a\in S,i\in \{1,...,Q\}}\{y\in 9\delta\text{-}\geod(x,a), |d(x,y)-t_{i}d(x,a)|\leq (Q+2)N+1\}\end{gather*}
et des distances entre ces points.
\end{lem}
%Remarquons que les conditions 
%$$y\in 9\delta\text{-}\geod(x,a), |d(x,y)-t_{i}d(x,a)|\leq (q+1)N+1$$ qui apparaissent à la fin de 1)b) (et à la fin de 2)b) avec $q=Q$) impliquent 
%$|d(y,a)-(1-t_{i})d(x,a)|\leq (q+1)N+9\delta+1$. 

\noindent{\bf Démonstration.} 
Montrons 1)a). 
C'est essentiellement le lemme~\ref{suite-S} mais on doit en répéter les arguments parce qu'on a ici $h_{x,t_{i}}$ et non $h_{x}$. 
Soit $S_{0}=S$ et $S_{1},\dots , S_{q}$ tels que
 $e_{S_{i}}$  apparaisse avec un coefficient non nul  
dans $(1-(\del h_{x,t_{i}}+h_{x,t_{i}}\del))(e_{S_{i-1}})$ pour $i=1,\dots , q-1$ et que $e_{S_{q}}$  apparaisse avec un coefficient non nul  
dans $h_{x,t_{q}}(e_{S_{q-1}})$. Soit $y_{q}\in S_{q}$ n'appartenant pas à $B(x,2\delta)$. 
Par le lemme~\ref{ASx-geod-x}, il existe $y_{q-1}\in S_{q-1}$,..., $y_{0}\in S_{0}$  n'appartenant pas à $B(x,2\delta)$ tels que $y_{i}=y_{i+1}$ ou bien $y_{i+1}\in 3\delta\text{-}\geod(x,y_{i})$ et $d(y_{i},y_{i+1})\in ]N-2\delta,N] $
pour $i\in \{0,\dots ,q-1\}$. 
En répétant la preuve du lemme~\ref{suite-S} 
 on montre que $y_{q}\in 4\delta\text{-}\geod(x,y_{0})$ et $d(y_{0},y_{q})>N-2\de$ si $y_{0}\neq y_{q}$. 
Comme $d(y_{i},y_{i+1})\leq N$   pour $i\in \{0,\dots,q-1\}$ on a $d(y_{0},y_{q})\leq qN$. 

 Montrons b). Soient $q\in \{1,...,Q\}$, $S_{0}=S$ et $S_{1},\dots , S_{q-1}$ tels que
 $e_{S_{i}}$  apparaisse avec un coefficient non nul 
dans $(1-(\del h_{x,t_{i}}+h_{x,t_{i}}\del))(e_{S_{i-1}})$ pour $i=1,\dots , q-1$.
 Soit $i\in \{1,\dots,q\}$.  Grâce au lemme~\ref{hxt-eS-connaissance}, la connaissance de $(1-(\del h_{x,t_{i}}+h_{x,t_{i}}\del))(e_{S_{i-1}})$  si $i<q$ ou de $h_{x,t_{i}}(e_{S_{i-1}})$  si $i=q$ ne dépend que de la connaissance 
 des points de $B(x,2\de)$ et de la réunion pour $b\in S_{i-1}$ des ensembles 
\begin{gather} \nonumber \{b\}\cup \{y\in 3\delta\text{-}\geod(x,b), d(y,b)\in ]N-2\delta,N]\} \\ \label{con-hxti-form-22oct09}\cup \{y\in 5\delta\text{-}\geod(x,b), 
d(x,y)\in [t_{i}d(x,S_{i-1})-2N-1,t_{i}d(x,S_{i-1})+N]\} \end{gather}
(en effet si $T\in \Delta\setminus\{\emptyset\}$ est tel que $e_{T}$ apparaisse avec un coefficient non nul dans $\del(e_{S_{i-1}})$ on a $T\subset S_{i-1}$ donc $d(x,T)\in [d(x,S_{i-1}), d(x,S_{i-1})+N]$). 
 De plus la preuve de a) montre que  pour tout $b\in S_{i-1}$, 
 on a $b\in S_{0}\cup B(x,2\de)$ ou 
 il existe $a\in S_{0}$ tel que 
 \begin{gather}\label{cond-a-b-0-(i-1)-22oct09}a\not\in B(x,2\de), 
 b\in 4\de\tg(x,a)\text{\  et\ } d(a,b)\in ]N-2\de, (q-1)N].\end{gather} 
 Si $b\in S_{0}\cup B(x,2\de)$  on vérifie facilement que l'ensemble (\ref{con-hxti-form-22oct09}) est inclus dans 
 l'ensemble (\ref{conn-tildeH-form-22oct09})
 car $|d(x,S_{i-1})-d(x,a)|\leq iN\leq qN$. Si $a\in S_{0}$ vérifie (\ref{cond-a-b-0-(i-1)-22oct09}), l'ensemble (\ref{con-hxti-form-22oct09}) est également inclus dans 
 l'ensemble (\ref{conn-tildeH-form-22oct09}) car 
 \begin{itemize}
 \item si $y\in 3\de\tg(x,b)$ et $d(y,b)\in ]N-2\de,N]$, en appliquant le lemme~\ref{iter2} avec $\alpha=4\de$, $\beta=3\de$ on a $y\in 4\de\tg(x,a)$, et  de plus $d(a,y)\leq qN$,
 \item si $y\in 5\de\tg(x,b)$ et $d(x,y)\in [t_{i}d(x,S_{i-1})-2N-1,t_{i}d(x,S_{i-1})+N]$ le a) du lemme~\ref{geod-comp-xabc} montre $y\in 9\de\tg(x,a)$ et comme $|d(x,S_{i-1})-d(x,a)|\leq qN$ on a $|d(x,y)-t_{i}d(x,a)|\leq (q+2)N+1$. 
 \end{itemize}
 %fait intervenir les points $y\in X$ tels qu'il existe une partie $T$ de $S_{i-1}$ (en fait égale à $S_{i-1}$, ou à $S_{i-1}$ privé d'un point si $i<q$) et que 
%\begin{itemize}
%\item ou bien $y\in A_{T,x}$,
%\item ou bien il existe   $z\in A_{T,x}$ tel que $y\in \delta\text{-}\geod(x,z)$ et 
%$$d(x,y)=\max(0,
%E(t_{i}(d(x,U_{T})-N))).$$ 
%\end{itemize}
%Dans le premier cas, si $y\not\in B(x,2\de)$,  il existe $a\in S_{0}$ tel que  $y\in 4\delta\text{-}\geod(x,a)$ et $ d(y,a)\leq qN$,  par la preuve du a). Dans %le second cas, si $z\not\in B(x,2\de)$ on a 
%$z\in 4\delta\text{-}\geod(x,a)$ pour un certain $a\in S_{0}$ par les arguments de la preuve du a), et donc 
%$y\in 5\delta\text{-}\geod(x,a)$ par le a) du lemme~\ref{geod-comp-xabc}. On a aussi $$d(x,a)-qN\leq d(x,U_{T})\leq d(x,U_{S_{i-1}})\leq d(x,a)$$ et on en déduit 
%$$t_{i}d(x,a)-(q+1)N-1\leq d(x,y)\leq t_{i}d(x,a).$$ Enfin dans le second cas, si $z\in B(x,2\de)$ on a $y\in B(x,3\de)$. 
%% et comme $d(a,z)\leq qN$, on a 
%%$t_{q}d(x,a)-2qN\leq d(z,t)\leq t_{q}d(x,a)+qN$. 
%Ceci démontre 1)b)

%Pour montrer b) on applique, le corollaire~\ref{suite-S-connaissance}.  

La preuve de 2) est tout à fait similaire à celle de 1), sauf que 2)a) nécessite un argument supplémentaire. 
Si $e_{T}$ apparaît avec un coefficient non nul dans $K_{x,Q,(t_{1},\dots ,t_{Q})}(e_{S})$, le lemme~\ref{nilp} implique $\zeta_{x}(T)\leq \zeta_{x}(S)-Q\frac{N-6\delta}{p_{\max}}$. Le sous-lemme suivant, appliqué à $M=Q\frac{N-6\delta}{p_{\max}}$,  montre alors que pour tout $y\in T$ et $a\in S$ on a $d(a,y)\geq Q\frac{N-6\delta}{p_{\max}}-2N$. 
\begin{souslem}\label{zeta-estimees-dist}
Soient $M\in \R_{+}$ et $S,T\in \Delta$. Si $\zeta_{x}(T)\leq \zeta_{x}(S)-M$, pour $y\in T$ et $a\in S$, on a $d(x,y)\leq d(x,a)-M+2N$. 
\end{souslem}
\noindent{\bf Démonstration.} On a  $d(x,y)\leq d(x,U_{T})+N\leq \zeta_{x}(T)+N$ et  
$d(x,a)\geq d(x,U_{S})\geq \zeta_{x}(S)-N$. \cqfd

\noindent{\bf Fin de la démonstration du lemme~\ref{support-connaissance-tildeH}.}
On suppose $Q\frac{N-6\delta}{p_{\max}}-2N>N$, ce qui est permis par $(H_{Q})$, d'où $y\not\in S$. Enfin $S$ est non vide, donc $B(x,2\de)\subset \bigcup_{a\in S}4\de\tg(x,a)$. 
\cqfd 

\label{def-Jx}
Comme nous l'avons déjà dit nous posons 
$$J_x=\tilde H_x+u_xK_{x}.$$ Montrons que 
$$\del J_x+J_x\del =1.$$ D'abord comme $\del ^{2}=0$, $\del$ commute à $K_{x}=1-\del \tilde H_x -\tilde H_x \del$ et donc 
$$\del J_x+J_x\del =(\del \tilde H _x+\tilde H_x \del)
+(\del u_x+u_x\del)K_{x}=1$$ puisque $\del u_x+u_x\del=1$.
On a   \begin{gather*}J_x=\sum_{q=1}^{Q}
\int_{(t_{1},\dots ,t_{q})\in [0,1]^{q}}\tilde H_{x,q,(t_{1},\dots ,t_{q})} dt_{1}\dots dt_{q}
\\ +\sum_{r=1}^{+\infty} \int_{t,(t_{1},\dots ,t_{Q})\in [0,1]^{Q+1}} u_{x,r,t} K_{x,Q,(t_{1},\dots ,t_{Q})}  dtdt_{1}\dots dt_{Q}.\end{gather*}

%
%comme une somme sur $r$  et sur $q\in \{0,...,p_{\max}-1\}$ et une intégrale sur $t\in [0,1]$ et sur $(t_{i})\in [0,1]^{q+1}$ de 
%$J_{x,r,q,t,(t_{i})}=\tilde H_{x,q,(t_{i})}+u_{x,r,t}(1-\del \tilde H_{x,q,(t_{i})} -\tilde H_{x,q,(t_{i})} \del)$. 

Il résulte de tout ce qui précède que  $\tilde H_{x,q,(t_{1},\dots ,t_{q})}(e_{S})$ et 
$u_{x,r,t} K_{x,Q,(t_{1},\dots ,t_{Q})}  (e_{S})$
ne dépendent  que de la connaissance d'un nombre fini de points (borné par $C(\delta,K,N,Q)$)  et des distances entre ces points. De plus ces distances sont déterminées par $x,S,q,r,t$ et les $t_{i}$ à une constante $C(\delta,K,N,Q)$ près. 
Plus précisément le 1) du lemme~\ref{support-connaissance-tildeH} apporte toutes les informations relatives à $\tilde H_{x,q,(t_{1},\dots ,t_{q})}(e_{S})$, et le lemme suivant les donne pour $u_{x,r,t} K_{x,Q,(t_{1},\dots ,t_{Q})}  (e_{S})$. 

\begin{lem}\label{support-connaissance-uK}
Pour $t\in [0,1]$, $(t_{1},\dots ,t_{Q})\in [0,1]^{Q}$, $r\in \N^{*}$  et   $S\in \Delta\setminus\{\emptyset\}$, 

\noindent a) 
$u_{x,r,t} K_{x,Q,(t_{1},\dots ,t_{Q})}  (e_{S})$ 
est une combinaison de $e_{T}$ où $T$ vérifie : 
$$T\subset 
\bigcup_{a\in S, 
z\in 5\de\tg(x,a), d(z,a)\in [Q\frac{N-6\delta}{p_{\max}}-2N-5\de+r,QN+r]}
B(z,N+5p_{\max}\de),$$

\noindent b) 
$u_{x,r,t} K_{x,Q,(t_{1},\dots ,t_{Q})}  (e_{S})$ ne dépend que de la connaissance des points de
\begin{gather*}B(x,7\de)\cup \bigcup_{a\in S}\{y\in 4\delta\text{-}\geod(x,a), d(y,a)\leq QN\} \\ \cup 
\bigcup_{a\in S,i\in \{1,...,Q\}}\{y\in 9\delta\text{-}\geod(x,a), |d(x,y)-t_{i}d(x,a)|\leq (Q+2)N+1\} \\ \cup \bigcup_{a\in S, 
z\in 5\de\tg(x,a), d(z,a)\in [Q\frac{N-6\delta}{p_{\max}}-2N-5\de+r,QN+r]}
B(z,N+5p_{\max}\de+5\de) \\ \cup \bigcup_{a\in S}\big\{z\in 5\de\tg(x,a), d(x,z)\in \\ [(1-t)(d(x,a)-QN-r),(1-t)(d(x,a)-r)+2+\de]\big\}\end{gather*}
et des distances entre ces points.
\end{lem}
\noindent{\bf Démonstration.}  
Montrons a). D'après le 2)a) du lemme~\ref{support-connaissance-tildeH} et la proposition~\ref{supp-uxrt}, $$u_{x,r,t} K_{x,Q,(t_{1},\dots ,t_{Q})}  (e_{S})$$ est supporté par la réunion des $B(z,N+5p_{\max}\de)$ pour les $z$ tels 
que $z\in \de\tg(x,y)$ et $d(y,z)=r$, avec 
  $y\in 4\de\tg(x,a)$, $a\in S$ et  $ d(y,a)\in [Q\frac{N-6\delta}{p_{\max}}-2N,QN]$. 
Alors  $z\in 5\de\tg(x,a)$ par le a) du lemme~\ref{geod-comp-xabc}. On a aussi $y\in 5\de\tg(z,a)$ par le b) du lemme~\ref{geod-comp-xabc}, d'où  $$d(z,a)\in [Q\frac{N-6\delta}{p_{\max}}-2N-5\de+r,QN+r].$$ 

\noindent Montrons b). 
D'après le 2)b) du lemme~\ref{support-connaissance-tildeH}, $K_{x,Q,(t_{1},\dots ,t_{Q})}   (e_{S})$ dépend de la connaissance des points  de 
\begin{gather*}B(x,7\de)\cup \bigcup_{a\in S}\{y\in 4\delta\text{-}\geod(x,a), d(y,a)\leq QN\} \\ \cup 
\bigcup_{a\in S,i\in \{1,...,Q\}}\{y\in 9\delta\text{-}\geod(x,a), |d(x,y)-t_{i}d(x,a)|\leq (Q+2)N+1\}. \end{gather*}
D'après le 2)a) du lemme~\ref{support-connaissance-tildeH} et la proposition~\ref{supp-uxrt}, $$u_{x,r,t}K_{x,Q,(t_{1},\dots ,t_{Q})}   (e_{S})$$  dépend de la connaissance de $K_{x,Q,(t_{1},\dots ,t_{Q})}   (e_{S})$ et des points  de 
\begin{gather*}\{x\}\cup \bigcup_{a\in S}\{y\in 4\delta\text{-}\geod(x,a), d(y,a)\in [Q\frac{N-6\delta}{p_{\max}}-2N,QN]\} \\ 
%\cup \bigcup_{a\in S}\{y\in 4\delta\text{-}\geod(x,a), d(y,a)\leq QN\}$$ $$\cup 
%\bigcup_{a\in S,i\in \{1,...,Q\}}\{y\in 9\delta\text{-}\geod(x,a), |d(x,y)-t_{i}d(x,a)|\leq (Q+1)N\}$$
\cup \bigcup_{a\in S, 
y\in 4\de\tg(x,a), d(y,a)\in [Q\frac{N-6\delta}{p_{\max}}-2N,QN], z\in \de\tg(x,y),d(y,z)=r}
B(z,N+5p_{\max}\de+5\de)\\ 
%$$\cup \bigcup_{a\in S, 
%z\in 5\de\tg(x,a), d(z,a)\in [Q\frac{N-6\delta}{p_{\max}}-2N-5\de+r,QN+r]}
%B(z,N+5p_{\max}\de+5\de)$$
\cup \bigcup_{a\in S,y\in 4\de\tg(x,a), d(y,a)\in 
[Q\frac{N-6\delta}{p_{\max}}-2N,QN]}\big\{  z\in \de\tg(x,y),\\ d(y,z)=r+E(t(d(x,y)-r))\text{\  ou \  }d(y,z)=(r-1)+E(t(d(x,y)-(r-1)))\big\}\end{gather*}
%et des points $z$ tels qu'il existe  $a\in S$ et $y\in 4\de\tg(x,a)$ avec $d(y,a)\in 
%[Q\frac{N-6\delta}{p_{\max}}-2N,QN]$, $z\in \de\tg(x,y)$ et $d(y,z)=r+E(t(d(x,y)-r))$ ou  $d(y,z)=(r-1)+E(t(d(x,y)-(r-1)))$. 
Les deux derniers ensembles de la réunion ci-dessus sont inclus dans les deux derniers ensembles de la réunion figurant dans b) du lemme~\ref{support-connaissance-uK}, pour les raisons suivantes. 
La preuve de a) montre que pour $a\in S$ les conditions 
$$ 
y\in 4\de\tg(x,a), d(y,a)\in [Q\frac{N-6\delta}{p_{\max}}-2N,QN], z\in \de\tg(x,y),d(y,z)=r$$ impliquent 
$$z\in 5\de\tg(x,a), d(z,a)\in [Q\frac{N-6\delta}{p_{\max}}-2N-5\de+r,QN+r].$$
D'autre part  pour $a\in S$  les conditions 
\begin{gather*} y\in 4\de\tg(x,a), d(y,a)\in 
[Q\frac{N-6\delta}{p_{\max}}-2N,QN], z\in \de\tg(x,y),\\ d(y,z)=r+E(t(d(x,y)-r))\text{  ou   }d(y,z)=(r-1)+E(t(d(x,y)-(r-1)))\end{gather*}
impliquent $z\in 5\de\tg(x,a)$, 
\begin{gather*}r+t(d(x,y)-r)-2\leq d(y,z) \leq r+t(d(x,y)-r),\\ 
d(x,z)\geq d(x,y)-d(y,z)\geq (1-t)(d(x,y)-r), \\ 
d(x,z)\leq d(x,y)-d(y,z)+\de\leq (1-t)(d(x,y)-r)+2+\de, \\ 
d(x,y)\geq d(x,a)-d(y,a)\geq d(x,a)-QN\text{\ \ \ \ et} \\ 
d(x,y)\leq d(x,a)-d(y,a)+4\de\leq d(x,a)-(Q\frac{N-6\delta}{p_{\max}}-2N)+4\de\leq d(x,a)\end{gather*}
 car on suppose $
(Q\frac{N-6\delta}{p_{\max}}-2N)\geq 4\de$ (ce qui est permis par  $(H_{Q})$).
%Mais un tel point $z$ appartient au dernier ensemble figurant dans l'énoncé du lemme~\ref{support-connaissance-uK} car $d(x,z)\geq d(x,y)-d(y,z)\geq (1-t)(d(x,y)-r)$, $d(x,z)\leq d(x,y)-d(y,z)+\de\leq (1-t)(d(x,y)-r)+2+\de$, $d(x,y)\geq d(x,a)-d(y,a)\geq d(x,a)-QN$ et 
%$d(x,y)\leq d(x,a)-d(y,a)+4\de\leq d(x,a)-(Q\frac{N-6\delta}{p_{\max}}-2N)+4\de\leq d(x,a)$ puisque l'on suppose $Q$ assez grand pour que $
%(Q\frac{N-6\delta}{p_{\max}}-2N)\geq 4\de$. 
\cqfd
 
Le 1) et le 2) de 
la proposition suivante récapitulent la partie des lemmes~\ref{support-connaissance-tildeH}   et~\ref{support-connaissance-uK} qui nous sera utile ensuite sous une forme plus lisible. 
 On suppose \begin{gather}\label{def-F-hyp}
 N\geq 6\de+1, Q\geq 2\text{\  et\ }
 Q\frac{N-6\delta}{p_{\max}}\geq 2(3N+5\de p_{\max}+10\de),\end{gather} ce qui est permis par $(H_{N})$ et $(H_{Q})$. 
On pose \begin{gather}\label{def-F}F=15\de+2N+10\de p_{\max}.\end{gather} On a donc $F\leq C(\de,K,N)$. Il est important de souligner que $F$ est une simple notation permettant d'alléger les formules. Au contraire  les constantes $N,Q$ et d'autres qui seront introduites ensuite, sont des paramètres dans notre construction et  chacun de ces paramètres doit être choisi suffisamment grand par rapport à ceux introduits auparavant.  

\begin{prop}\label{recap-supp-connaiss-H-uK}
1) 
Pour $q\in \{1,...,Q\}$, 
$(t_{1},\dots ,t_{q})\in [0,1]^{q}$ et $S\in \Delta\setminus\{\emptyset\}$,

\noindent a)  $\tilde H_{x,q,(t_{1},\dots ,t_{q})}(e_{S})$
est une combinaison de $e_{T}$ où $T$ vérifie : 
\begin{gather*}T\subset  S\cup \Big(B(x,2\de)\cap B(S,qN)\Big)\\ \cup\bigcup_{a\in S}\big\{y\in 4\delta\text{-}\geod(x,a), d(y,a)\in
]N-2\de
%[\max(N-2\de,(q-1)\frac{N-6\delta}{p_{\max}}-2N)
, q N]\big\},\end{gather*}

\noindent b)   $\tilde H_{x,q,(t_{1},\dots ,t_{q})}(e_{S})$ ne dépend que de la connaissance des points de
\begin{gather*}B(x,7\de)\cup B(S,QN) \\ \cup 
\bigcup_{a\in S,i\in \{1,...,q\}}\{y\in F\text{-}\geod(x,a), |d(x,y)-t_{i}d(x,a)|\leq QF \}\end{gather*}
et des distances entre ces points.

2) Pour $t\in [0,1]$, $(t_{1},\dots ,t_{Q})\in [0,1]^{Q}$, $r\in \N^{*}$  et   $S\in \Delta\setminus\{\emptyset\}$, 

\noindent a) 
$u_{x,r,t} K_{x,Q,(t_{1},\dots ,t_{Q})}  (e_{S})$ 
est une combinaison de $e_{T}$ où $T$ vérifie : 
\begin{gather*}T\subset \bigcup_{a\in S}\{ 
z\in F\tg(x,a), d(z,a)\in [\frac{Q}{F}+r,QF+r]
\},\end{gather*}

\noindent b) 
$u_{x,r,t} K_{x,Q,(t_{1},\dots ,t_{Q})}  (e_{S})$ ne dépend que de la connaissance des points de
\begin{gather*}B(x,F)\cup B(S,QN) \cup \bigcup_{a\in S}\{y\in F\tg(x,a), d(y,a)\in [r,r+QF]\}\\ 
\cup \bigcup_{a\in S,i\in \{1,...,Q\}}\{y\in F\text{-}\geod(x,a), |d(x,y)-t_{i}d(x,a)|\leq QF \}\\ \cup
\bigcup_{a\in S}\{y\in F\tg (x,a), |d(x,y)-(1-t)(d(x,a)-r)|\leq QF \}\end{gather*}
et des distances entre ces points.

3) Il existe $C=C(\de,K,N,Q)$ tel que dans les notations de 1) et 2) on ait $$\|  \tilde H_{x,q,(t_{1},\dots ,t_{q})}(e_{S})\|_{\ell^{1}}\leq C\text{ \ et \ } \| u_{x,r,t} K_{x,Q,(t_{1},\dots ,t_{Q})} (e_{S})\|_{\ell^{1}}\leq C.$$ 
\end{prop}
On remarque que le nombre de points des ensembles apparaissant dans l'énoncé   est borné par $C(\de,K,N,Q)$ et que les distances entre ces points sont déterminées à $C(\de,K,N,Q)$ près par $d(x,S)$ et les divers paramètres  (à savoir $q, 
t_{1},\dots ,t_{q}$ dans 1) et  $r,t,t_{1},\dots,t_{Q}$ dans 2)).   Donc le nombre de combinaisons possibles pour le nombre de ces points et leurs distances mutuelles  est borné par $C(\de,K,N,Q)$. 

Par des arguments similaires à ceux de la preuve du lemme~\ref{suite-S} on peut montrer facilement, à l'aide de 1)a) et 2)a) du lemme précédent, que $J_{x}$ vérifie la condition (C1) (en fait ces arguments seront cachés dans l'étude des normes menée au paragraphe~\ref{construction-normes}). D'autre part 1)b) et 2)b) garantissent que $J_{x}$ vérifie (C2).  

\noindent 
{\bf Démonstration.} Le 1) résulte du 1) du lemme~\ref{support-connaissance-tildeH}   et le 2) résulte du  lemme~\ref{support-connaissance-uK}. En particulier pour montrer 2) on utilise  les inégalités suivantes, qui découlent de 
(\ref{def-F-hyp})  et (\ref{def-F}) : 
\begin{gather*}QN+N+5\de p_{\max}+5\de\leq QF, \ (Q+2)N+1\leq QF\\ 
\text{et\ \ \ } Q\frac{N-6\delta}{p_{\max}}-2N-5\de-N-5\de p_{\max}-5\de\geq  
Q\frac{N-6\delta}{2p_{\max}}\geq \frac{Q}{F}.\end{gather*}
On utilise aussi le fait que 
\begin{gather*} \bigcup_{a\in S, 
z\in 5\de\tg(x,a), d(z,a)\in [Q\frac{N-6\delta}{p_{\max}}-2N-5\de+r,QN+r]}
B(z,N+5p_{\max}\de+5\de)\\
\subset \bigcup_{a\in S} \{z\in (5\de+2(N+5p_{\max}\de+5\de))\tg(x,a), \\ 
d(z,a)\in [Q\frac{N-6\delta}{p_{\max}}-2N-5\de+r-(N+5p_{\max}\de+5\de),\\ QN+r+(N+5p_{\max}\de+5\de)]\}
\end{gather*} grâce au lemme~\ref{xx'yy'zz'}, ainsi que d'autres inclusions analogues ou plus faciles. 

Enfin 3) résulte du lemme~\ref{28dec1029} et 
du fait que $\|h_{x,t}(e_{S})\|_{\ell^{1}}\leq 1$ pour tout $S$ puisque $h_{x,t}(e_{S})=\psi_{S,x,t}\wedge  e_{S}$ et que $\psi_{S,x,t}$ est une mesure de probabilité. 
  \cqfd

\subsection {Construction d'une distance moyennée pour conjuguer les opérateurs}

On montrera dans le paragraphe~\ref{construction-normes} que l'image de $1$ dans $KK_{G,2s\ell+C}$ $(\C,\C)$
est représenté par l'opérateur impair $\del+J_{x}$ agissant sur le complété de $\bigoplus_{p=1}^{p_{\max}}\C^{(\Delta_{p})}$ pour une certaine norme de Hilbert $\|.\|_{\H_{x,s}}$. Pour construire la première partie de l'homotopie de $1$ à $\gamma$ (qui fait l'objet de ce paragraphe et du suivant), on  conjuguera $\del+J_{x}$ par $e^{\tau \theta^{\flat}_{x}}$ où  $\theta^{\flat}_{x}:\bigoplus_{p=1}^{p_{\max}}\C^{(\Delta_p)}\to \bigoplus_{p=1}^{p_{\max}}\C^{(\Delta_p)} $ est   défini par $\theta^{\flat}_{x}(e_{S})=\rho^{\flat}_{x}(S)e_{S}$,  où 
$\rho^{\flat}_{x}$ est une variante moyennée de la distance à $x$ et où $\tau$ varie de $0$ à $T$ (avec $T$ assez grand). 
 On doit  moyenner  la  distance   à $x$ pour qu'elle se comporte mieux quand on change l'origine $x$ en $g(x)$ pour $g\in G$.  Cela est nécessaire  pour que  l'opérateur $\del+J_{x}$ conjugué par $e^{\tau  \theta^{\flat}_{x}}$ soit $G$-équivariant à compact près. 

La suite de l'homotopie  (qui est reléguée au paragraphe~\ref{construction-fin}) sera facile : comme $T$ est assez grand, l'opérateur  $e^{T\theta^{\flat}_{x}}(\del+J_{x})e^{-T\theta^{\flat}_{x}}$ est continu pour la norme $\ell^{2}$, donc  par une homotopie et tout en gardant cet opérateur on peut remplacer la norme $\|.\|_{\H_{x,s}}$ par la norme $\ell^{2}$ sur $\bigoplus_{p=1}^{p_{\max}}\C^{(\Delta_{p})}$ et il est alors très simple de terminer l'homotopie en aboutissant à  $\gamma$ dans $KK_{G}(\C,\C)$. 

Le but de ce sous-paragraphe est la construction de $\rho^{\flat}_{x}$.

La proposition suivante renforce le théorème 17  de~\cite{mineyevyu}. 
On rappelle que $(X,d)$ est un bon espace discret $\delta$-hyperbolique. 
\begin{prop}\label{myu} 
%Soit $\delta>0$ et $(X,d)$ un espace $\delta$-hyperbolique et $\delta$-faiblement géodésique. 
Il existe une distance $G$-invariante $d''$ sur $X$ telle que $d-d''$ est
borné et que 
\begin{gather} \nonumber
\forall r\in \R_+, \forall \epsilon\in \R_+^*, 
\exists R\in \R_+, \forall x,x',y,y'\in X, d(x,x')\leq r,d(y,y')\leq
r,\\  \label{condition-d'} d(x,y)\geq R, \text{ on ait } 
|d''(x,y)-d''(x,y')-d''(x',y)+d''(x',y')|\leq \epsilon .\end{gather}
\end{prop}
La proposition~\ref{myu} est une conséquence du lemme suivant. 

\begin{lem}\label{myumesure} 
Il existe une famille de mesures positives $(\mu (x,y))_{(x,y)\in X\times X} $
de masse $1$ sur $X$ de sorte que 
\begin{itemize}
\item pour tout $g\in G$, $\mu(gx,gy)=g_{*}\mu(x,y)$, 
\item $\mu(x,y)=\mu(y,x)$, 
\item le support de $\mu(x,y)$ est inclus dans $7\delta\text{-}\geod(x,y)$, 
\item et \begin{gather}\nonumber \forall r\in \R_+, \forall \epsilon\in \R_+^*, 
\exists R\in \R_+, \\ \nonumber \forall x,x',y,y'\in X, \text{\ avec\ }d(x,x')\leq r,d(y,y')\leq r,d(x,y)\geq R,\\ \nonumber \text{on\ \ a\ \ }\|\mu(x,y)-\mu(x',y')\|_1 \leq \epsilon .\end{gather}
\end{itemize}
\end{lem}
La première condition est simplement la condition naturelle de $G$-équivariance. 

\noindent{\bf Démonstration de la proposition~\ref{myu} en admettant le lemme~\ref{myumesure}.} 
On pose \begin{gather}\label{ineg-28dec1056}d'(x,y)=\int_{X} (d(x,z)+d(z,y))d\mu(x,y)(z).\end{gather} On a 
$d(x,y)\leq d'(x,y) \leq d(x,y)+7\delta$ pour  $x,y\in X$ et 
$d'$ vérifie la condition (\ref{condition-d'}) de la proposition~\ref{myu}. 
En effet soit $r\in \R_{+}, \epsilon\in \R_{+}^{*}$. Soit $R$, et $x,y,x',y'$ comme dans la dernière assertion du lemme~\ref{myumesure}. 
Pour $z\in 7\delta\text{-}\geod(x,y)$, on a $z\in (7\delta+2r)\text{-}\geod(x,y')$, $z\in (7\delta+2r)\text{-}\geod(x',y)$, 
 et $z\in (7\delta+4r)\text{-}\geod(x',y')$ par le lemme~\ref{xx'yy'zz'}, donc 
\begin{gather*}|d'(x',y)-\int_{X} (d(x',z)+d(z,y))d\mu(x,y)(z)|\\ 
=\Big|\int_{X} (d(x',z)+d(z,y))(d\mu(x',y)-d\mu(x,y))(z)\Big| \\ 
= \Big|\int_{X} (d(x',z)+d(z,y)-d(x',y))(d\mu(x',y)-d\mu(x,y))(z)\Big|\leq \epsilon(7\delta+2r)
 \\ \text{et de même \ \ }|d'(x,y')-\int_{X} (d(x,z)+d(z,y'))d\mu(x,y)(z)|\leq \epsilon(7\delta+2r) \\ \text{et\ \ \ \ } |d'(x',y')-\int_{X} (d(x',z)+d(z,y'))d\mu(x,y)(z)|\leq \epsilon(7\delta+4r) .\end{gather*}
D'autre part \begin{gather*}(d(x,z)+d(z,y))-(d(x',z)+d(z,y))-(d(x,z)+d(z,y'))\\ +(d(x',z)+d(z,y'))=0\end{gather*} pour tout $z\in X$.
On en déduit 
$$|d'(x,y)-d'(x,y')-d'(x',y)+d'(x',y')|\leq \epsilon (21\delta+8r).$$ 

Cependant $d'$ n'est  pas nécessairement une distance. Comme  
$$d(x,y)\leq d'(x,y) \leq d(x,y)+7\delta\text{\   pour  \ }x,y\in X$$ on a 
$d'(u,w)\leq
d'(u,v)+d'(v,w)+7\delta$ pour  $u,v,w\in X$. On pose alors 
$d''(x,y)=0$ si $x=y$ et $d''(x,y)=d'(x,y)+7\delta$ si $x\neq y$. Alors  $d''$ est  une distance et on a  montré la proposition~\ref{myu} en admettant le lemme~\ref{myumesure}.  \cqfd
 
Avant de montrer le  lemme~\ref{myumesure} 
%Pour cela on a besoin des variantes assez triviales des lemmes~\ref{contract1} et~\ref{contract2}. 
%
%\begin{lem}\label{contract1'}
%Pour tout $\epsilon >0$ et $x,z,y,y'\in X$, si $y\in
%\epsilon\text{-}\geod(x,z)$, $y'\in \delta\text{-}\geod(y,z)$, $d(y,y')\geq
%(\epsilon+\delta)/2$, alors $y'\in \delta\text{-}\geod(x,z)$.
%\end{lem}
%
%En effet 
%$d(x,y')+d(y,z)\leq
%\max(d(z,y')+d(x,y),d(y,y')+d(x,z))+\delta$. 
%D'où $d(x,y')+d(y',z)\leq
%\max(d(z,y')+d(x,y)-d(y,y')+\delta,d(x,z))+\delta$. 
%Or $d(z,y')+d(x,y)-d(y,y')=(d(x,y)+d(y,z))-2d(y,y')\leq
%d(x,z)+\epsilon -2d(y,y')$. 
%
%On en déduit aussitôt le lemme suivant. 
%
%\begin{lem}\label{contract2'}
%Pour tous $\epsilon >0$, 
%$x,x',z,y,y'\in X$, si $y\in 
%\epsilon\text{-}\geod(x,z)$, $y'\in \delta\text{-}\geod(y,z)$, 
%$d(y,y')\geq 
%(\epsilon+\delta)/2+d(x,x')$, alors $y'\in \delta\text{-}\geod(x',z)$.
%\end{lem}
%
%En effet  $y\in (\epsilon+2d(x,x'))\text{-}\geod(x',z)$.  
 on commence par un lemme général très utile. 

\begin{lem}\label{lemmetubes}
Soient $x,y\in X$, $\alpha, \gamma\in \N$, $a\in \alpha\text{-}\geod(x,y)$, $c\in \gamma\text{-}\geod(x,y)$, 
et $b\in \geod(a,c)$, avec $d(a,b)\geq \frac{\alpha}{2}$ et $d(b,c)\geq \frac{\gamma}{2}$, alors $b\in 3\delta\text{-}\geod(x,y)$. 
\end{lem}

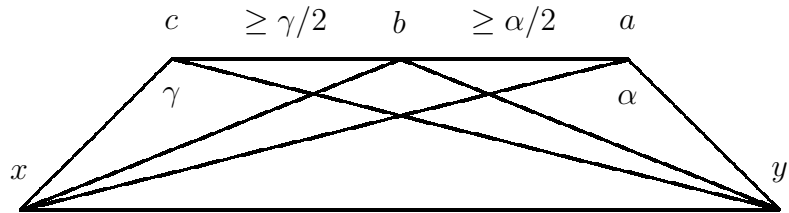
\ifx\JPicScale\undefined\def\JPicScale{1}\fi
\unitlength \JPicScale mm
\begin{picture}(120,35)(0,35)
\linethickness{0.3mm}
\put(20,40){\line(1,0){100}}
\linethickness{0.3mm}
\multiput(20,40)(0.12,0.12){167}{\line(1,0){0.12}}
\linethickness{0.3mm}
\multiput(40,60)(0.48,-0.12){167}{\line(1,0){0.48}}
\linethickness{0.3mm}
\multiput(100,60)(0.12,-0.12){167}{\line(1,0){0.12}}
\linethickness{0.3mm}
\multiput(20,40)(0.48,0.12){167}{\line(1,0){0.48}}
\linethickness{0.3mm}
\multiput(20,40)(0.3,0.12){167}{\line(1,0){0.3}}
\linethickness{0.3mm}
\multiput(70,60)(0.3,-0.12){167}{\line(1,0){0.3}}
\linethickness{0.3mm}
\put(40,60){\line(1,0){60}}
\put(20,45){\makebox(0,0)[cc]{$x$}}

\put(120,45){\makebox(0,0)[cc]{$y$}}

\put(40,65){\makebox(0,0)[cc]{$c$}}

\put(70,65){\makebox(0,0)[cc]{$b$}}

\put(100,65){\makebox(0,0)[cc]{$a$}}

\put(40,55){\makebox(0,0)[cc]{$\gamma$}}

\put(100,55){\makebox(0,0)[cc]{$\alpha$}}

\put(55,65){\makebox(0,0)[cc]{$\geq \gamma/2$}}

\put(85,65){\makebox(0,0)[cc]{$\geq \alpha/2$}}

\end{picture}

\noindent{\bf Démonstration.}
Par  $(H_{\delta}^{0}(x,a,b,c))$ et $(H_{\delta}^{0}(y,a,b,c))$ on a  
\begin{gather*}d(x,b)\leq \max(d(x,a)-d(a,b),d(x,c)-d(b,c))+\delta\\ \text{ et \ \ \ \ \ \ }d(y,b)\leq \max(d(y,a)-d(a,b),d(y,c)-d(b,c))+\delta.\end{gather*}
En additionnant ces deux inégalités on obtient 
\begin{gather*}d(x,b)+d(b,y)\leq 
\max(d(x,a)+d(a,y)-2d(a,b),d(x,c)+d(c,y)-2d(b,c),\\ d(x,a)+d(y,c)-d(a,c),d(x,c)+d(y,a)-d(a,c))+2\delta.\end{gather*}
Comme $d(a,b)\geq \frac{\alpha}{2}$ on a $d(x,a)+d(a,y)-2d(a,b)\leq d(x,y)$ et comme $d(b,c)\geq \frac{\gamma}{2}$ on a $d(x,c)+d(c,y)-2d(b,c)\leq d(x,y)$, donc 
\begin{gather*}d(x,b)+d(b,y)\leq \\ 
\max\big(d(x,y)+2\de, 
d(x,a)+d(y,c)-d(a,c)+2\delta,d(x,c)+d(y,a)-d(a,c)+2\delta\big).\end{gather*}
Si $d(x,a)+d(y,c)\leq d(a,c)+d(x,y)+\delta$ et $d(x,c)+d(y,a)\leq d(a,c)+d(x,y)+\delta$ on a fini. Dans le cas contraire, supposons par exemple $$d(x,a)+d(y,c)>d(a,c)+d(x,y)+\delta.$$ Grâce à $(H_{\delta}(x,y,a,c))$ qui s'écrit  
%Enfin, par la propriété d'hyperbolicité pour $x,y,a,c$, si l'un des nombres $d(x,a)+d(y,c)-d(a,c)$ et $d(x,c)+d(y,a)-d(a,c)$ est strictement supérieur à $d(x,y)+\delta$, disons par exemple 
%$d(x,a)+d(y,c)-d(a,c)>d(x,y)+\delta$, 
$$d(x,a)+d(y,c)\leq \max(d(x,y)+d(a,c),d(x,c)+d(y,a))+\de$$ on a alors  
$d(x,c)+d(y,a)\geq d(x,a)+d(y,c)-\delta>d(a,c)+d(x,y)$ d'où
$$\big(d(x,a)+d(y,c)\big)+\big(d(x,c)+d(y,a)\big)> 
2d(a,c)+2d(x,y)+\delta. $$
Or  $$\big(d(x,a)+d(y,c)\big)+\big(d(x,c)+d(y,a)\big)\leq 2d(x,y)+\alpha+\gamma\leq 2d(a,c)+2d(x,y). $$
Cette contradiction achève la démonstration du lemme~\ref{lemmetubes}. \cqfd

\noindent{\bf Démonstration du lemme~\ref{myumesure}.} Comme le lemme~\ref{myumesure} ne servira pas dans la suite nous abrégeons sa démonstration. 
On rappelle que pour toute partie $A$ non vide de $X$,  $\nu_A$ désigne la mesure de
masse $1$ égale au produit par $(\sharp A)^{-1}$ de la fonction
caractéristique de $A$. 
On pose alors $\phi(t)=E(t/8)$, 
 et $$\mu(x,y)=\frac{1}{(\phi(d(x,y))+1)^4}\sum_{k,\tilde k,l,\tilde l=0}^{\phi(d(x,y))}
\nu_{B_{x,y,k,\tilde k,l,\tilde l}}$$ où 
 pour $k,\tilde k,l,\tilde l\in \{0,\dots,\phi(d(x,y))\}$ on note 
\begin{gather*}B_{x,y,k,\tilde k,l,\tilde l}=\{z,\exists (\tilde x,\tilde y)\in 3\delta\text{-}\geod(x,y)^{2}, d(x,\tilde x)\leq k,  d(\tilde x,z)\geq \tilde k, \\ 
d(y,\tilde y)\leq l, d(z,\tilde y)\geq \tilde l,  z\in 
 \geod(\tilde x,\tilde y)\}.
  %, \\
%|d(x,z)-d(x,y)/2|\leq \phi(t)+k+1, |d(y,z)-d(x,y)/2|\leq
 % \phi(t)+l+1\}.
 \end{gather*}
 
 \ifx\JPicScale\undefined\def\JPicScale{1}\fi
\unitlength \JPicScale mm
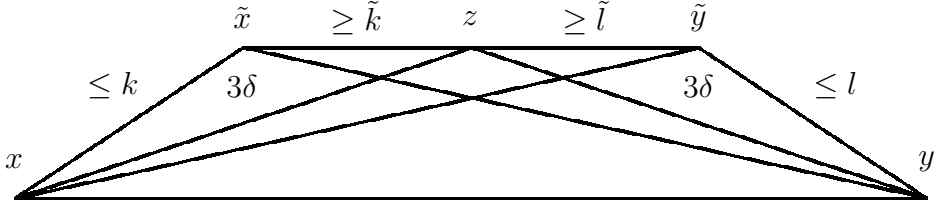
\begin{picture}(130,35)(30,35)
\linethickness{0.3mm}
\put(30,40){\line(1,0){120}}
\linethickness{0.3mm}
\multiput(120,60)(0.18,-0.12){167}{\line(1,0){0.18}}
\linethickness{0.3mm}
\put(60,60){\line(1,0){60}}
\linethickness{0.3mm}
\multiput(30,40)(0.18,0.12){167}{\line(1,0){0.18}}
\linethickness{0.3mm}
\multiput(60,60)(0.54,-0.12){167}{\line(1,0){0.54}}
\linethickness{0.3mm}
\multiput(30,40)(0.54,0.12){167}{\line(1,0){0.54}}
\put(30,45){\makebox(0,0)[cc]{$x$}}

\put(150,45){\makebox(0,0)[cc]{$y$}}

\put(60,64){\makebox(0,0)[cc]{$\tilde x$}}

\put(120,64){\makebox(0,0)[cc]{$\tilde y$}}

\linethickness{0.3mm}
\multiput(30,40)(0.36,0.12){167}{\line(1,0){0.36}}
\linethickness{0.3mm}
\multiput(90,60)(0.36,-0.12){167}{\line(1,0){0.36}}
\put(90,64){\makebox(0,0)[cc]{$z$}}

\put(75,64){\makebox(0,0)[cc]{$\geq \tilde k$}}

\put(105,64){\makebox(0,0)[cc]{$\geq \tilde l$}}

\put(43,55){\makebox(0,0)[cc]{$\leq k$}}

\put(138,55){\makebox(0,0)[cc]{$\leq l$}}

\put(60,55){\makebox(0,0)[cc]{$3\de$}}

\put(120,55){\makebox(0,0)[cc]{$3\de$}}

\end{picture}

 Il est clair que $\mu(x,y)=\mu(y,x)$. 
 Montrons que le support de $\mu(x,y)$ est inclus dans $7\delta\text{-}\geod(x,y)$. 
Soient $z,\tilde x$ et $\tilde y$ comme dans la définition de $B_{x,y,k,\tilde k,l,\tilde l}$. Si $d(\tilde x,z)\geq 2\delta $ et $d(\tilde y,z)\geq 2\delta $ le 
 lemme~\ref{lemmetubes} (avec $\alpha=\beta=3\delta$) implique 
 $z\in 3\delta\text{-}\geod(x,y)$. Si $d(\tilde x,z)\leq 2\delta $ ou $d(\tilde y,z)\leq 2\delta $ on a  $z\in 7\delta\text{-}\geod(x,y)$ par le lemme~\ref{xx'yy'zz'}. 
 Pour montrer la dernière assertion du lemme~\ref{myumesure} on 
 utilise l'astuce des ensembles emboîtés, qui apparaît
dans la démonstration de la proposition 6.9 de~\cite{ks}. Le
cardinal de $B_{x,y,k,\tilde k,l,\tilde l}$ est toujours compris entre $\frac{d(x,y)}{C_1}$ et
$C_2d(x,y)$ pour deux constantes $C_1$ et $C_{2}$ du type $C(\delta,K)$. De plus $(k,\tilde k,l,\tilde l)\mapsto B_{x,y,k,\tilde k,l,\tilde l}$ est une application croissante en $k$ et $l$ et décroissante en $\tilde k$ et $\tilde l$. 
%$B_{x,y,k,l}\subset
%B_{x,y,k+1,l}$, $B_{x,y,k,l}\subset
%B_{x,y,k,l+1}$ et 
%$B_{x',y',k,l}\subset
%B_{x,y,k+d(x,x')-\delta,l+d(y,y')-\delta}$ lorsque ces ensembles existent, par le lemme~\ref{contract2'}
Enfin étant donnés $x,y,x',y'\in X$, en posant $$d=d(x,x')+d(y,y')+2\delta$$ on a 
$$B_{x,y,k,\tilde k,l,\tilde l}\subset B_{x',y',k+2d,\tilde k-d,l+2d,\tilde l-d}$$  pour $\tilde k,\tilde l\geq d$.
En effet soit $z\in B_{x,y,k,\tilde k,l,\tilde l}$ et $\tilde x,\tilde y
\in 3\delta\text{-}\geod(x,y)^{2}$  tels que  $d(x,\tilde x)\leq k,  d(\tilde x,z)\geq \tilde k, 
d(y,\tilde y)\leq l, d(z,\tilde y)\geq \tilde l,  z\in 
 \geod(\tilde x,\tilde y)$.  
Par le lemme~\ref{xx'yy'zz'},  $\tilde x$ et $\tilde y$ appartiennent à $\big(3\delta+2d(x,x')+2d(y,y')\big)\text{-}\geod(x',y')$. 
  Soit 
 $\tilde x'$ un point de $\geod(\tilde x,z)$  à distance $d$ de $\tilde x$ et 
$\tilde y' $  un point  de $\geod(\tilde y,z)$  à distance $d$ de $\tilde y$. 
D'après le lemme~\ref{lemmetubes}, $\tilde x'$ et 
$\tilde y' $ appartiennent à $3\delta\text{-}\geod(x',y')$ et donc 
$z$ appartient à $B_{x',y',k+2d,\tilde k-d,l+2d,\tilde l-d}$. 
Pour conclure la démonstration du lemme~\ref{myumesure} on applique à 
$\Delta=d(x,y)$, $\Delta'=d(x',y')$, $n_{k,\tilde k,l,\tilde l}=\sharp B_{x,y,k,\tilde k,l,\tilde l}$ et 
$n'_{k,\tilde k,l,\tilde l}=\sharp B_{x',y',k,\tilde k,l,\tilde l}$ le sous-lemme suivant. 
\begin{souslem}\label{C1C2Delta}
Soient $d\in \N$, $C_{1},C_{2}\in \R_{+}^{*}$. Soit $\epsilon>0$. Il existe  $\Delta_{0}\in \N$ tel que pour  pour  $\Delta,\Delta'\in \N$ vérifiant $\Delta\geq \Delta_{0}$, $|\Delta'-\Delta|\leq d$ et pour $(n_{k,\tilde k,l,\tilde l})_{k,\tilde k,l,\tilde l\in \{0,\dots,\phi(\Delta)\}}$ des éléments de $\N\cap [\frac{\Delta}{C_{1}},C_{2}\Delta]$ et $(n'_{k,\tilde k,l,\tilde l})_{k,\tilde k,l,\tilde l\in \{0,\dots,\phi(\Delta')\}}$ des éléments de $\N\cap [\frac{\Delta'}{C_{1}},C_{2}\Delta']$, tels que 
\begin{itemize}
\item
$(k,\tilde k,l,\tilde l)\mapsto n_{k,\tilde k,l,\tilde l}$ soit croissant en $k$ et $l$ et décroissant en $\tilde k$ et $\tilde l$. 
\item $(k,\tilde k,l,\tilde l)\mapsto n'_{k,\tilde k,l,\tilde l}$ soit croissant en $k$ et $l$ et décroissant en $\tilde k$ et $\tilde l$. 
\item $n_{k,\tilde k,l,\tilde l}\leq n'_{k+2d,\tilde k-d,l+2d,\tilde l-d}$  
et $n'_{k,\tilde k,l,\tilde l}\leq n_{k+2d,\tilde k-d,l+2d,\tilde l-d}$  quand ces nombres ont un sens,
\end{itemize}
alors la proportion de $(k,\tilde k,l,\tilde l)$ dans $\{0,\dots,\phi(\Delta)\}^{4}$ tels que l'expression $\frac{n_{k,\tilde k,l,\tilde l}}{n'_{k+2d,\tilde k-d,l+2d,\tilde l-d}}$ ait un sens et appartienne à $[1-\epsilon,1]$ est supérieure ou égale à $1-\epsilon$.   \end{souslem}
\noindent{\bf Démonstration.} La démonstration est facile, et laissée au lecteur car ce sous-lemme ne servira pas dans la suite. 
\cqfd

\noindent{\bf Fin de la démonstration du lemme~\ref{myumesure}.} 
Si $A$ et $B$ sont deux parties finies de $X$ avec $A\subset B$ et $\sharp A/\sharp B \geq 1-\epsilon$ on a $\|\nu_{A}-\nu_{B}\|_{1}\leq 2\epsilon$. 
Le lemme~\ref{myumesure} en résulte facilement. \cqfd

%De nouveau nous supposons que $(X,d)$ est comme dans la proposition~\ref{prophyper}. 
 Nous pourrions poser $\rho^{\flat}_{x}(y)=d'(x,y)$ mais cela ne nous convient pas car $d'(x,y)$ fait intervenir une moyenne sur  $B_{x,y,k,\tilde k,l,\tilde l}$, donc nécessite de connaître le cardinal de cet ensemble. 
 Nous allons construire $d^{\flat}$ jouissant de  propriétés analogues à celles de $d'$ mais telle que  $(d^{\flat}-d)(x,y)$ soit une moyenne (sur un ensemble d'indices ne dépendant que de $
d(x,y)$) d'une fonction  à valeurs dans $[0,7\de]$,  qui ne dépend (c'est là le point important) que de la connaissance 
d'un nombre de points borné par  $C(\de,K)$ et des distances entre eux. Plus précisément ces points seront $x,y$, les $z\in 3\delta\tg(x,y)$ vérifiant $d(x,z)=r$ et les 
 $t\in 3\delta\tg(x,y)$ vérifiant $d(y,t)=s$,  pour $3(3\delta+1)$  valeurs différentes de $r$ et de $s$. 
Nous définirons alors $\rho^{\flat}_{x}$ en posant $\rho^{\flat}_{x}(y)=d^{\flat}(x,y)$. 
Le lecteur peut donc oublier la construction de $d'$ qui précède, c'est-à-dire la proposition~\ref{myu} et les lemmes~\ref{myumesure}  et~\ref{C1C2Delta} (en revanche le lemme~\ref{lemmetubes} sera réutilisé). L'idée de la construction de $d^{\flat}$ est de remplacer la moyenne par  $\mu(x,y)$ qui sert à construire $d'$  dans la formule~(\ref{ineg-28dec1056}) de  la preuve de la proposition~\ref{myu} par une moyenne sur un ``point virtuel'' défini comme la donnée de ses ``distances'' à certains autres points. Ces distances seront à valeurs dans $\Z$ et de plus on quotientera par une action naturelle de $\Z$ consistant à ajouter un entier relatif à certaines distances et à le retrancher aux autres. En particulier les distances d'un point virtuel à $x$ et $y$ ne seront pas bien définies, mais leur somme le sera. 
Les ensembles de points virtuels auront les deux propriétés suivantes : être  de cardinal  $\leq C(\de,K)$ et même temps varier  très peu souvent lorsqu'on bouge un peu $x$ ou $y$. Ces deux propriétés seraient contradictoires pour des points réels (dans la preuve du lemme~\ref{myumesure} c'est la deuxième propriété qui est en défaut car le cardinal du support de $\mu(x,y)$ tend vers l'infini quand $d(x,y)$ tend vers l'infini). 

%Pour tout entier $r\in \{0,...,E(\frac{d(x,y)}{2})\}$, on note $X_{x,y}^{r}$ l'ensemble  des points $z\in \delta\text{-}\geod(x,y)$ tels que $d(x,z)=r$. 
%On aurait bien aimé que le résultat suivant soit vrai : étant donnés trois entiers $r,r',r''\in \{0,...,E(\frac{d(x,y)}{2})\}$ avec $r\leq r' \leq r''$, et $z\in X_{x,y}^{r}$ et $z''\in X_{x,y}^{r''}$, 
%$\geod(z,z'')$ rencontre $X_{x,y}^{r'}$ (c'est-à-dire qu'il existe $z'\in X_{x,y}^{r'}$ tel que $d(z,z')+d(z',z'')=d(z,z'')$). 
%Ce n'est pas toujours vrai mais nous allons y remédier. Le lecteur peut oublier la notation $X_{x,y}^{r}$. 

La fonction $d^{\flat}:X\times X\to \R_{+}$ sera définie dans la 
formule (\ref{formule-d'}) ci-dessous. Jusqu'à (\ref{formule-d'}) on suppose $d(x,y)\geq 6\de$. 

Pour tout entier $r\in \{0,...,E(\frac{d(x,y)}{2})-3\delta\}$, on note $Y_{x,y}^{r}$ l'ensemble  des points $z\in 3\delta\text{-}\geod(x,y)$ tels que $$d(x,z)\in \{r,...,r+3\delta\}.$$

\begin{lem}\label{tubes} a) 
 Etant donnés trois entiers $r,r',r''\in \{0,...,E(\frac{d(x,y)}{2})-3\delta\}$ avec $r\leq r'\leq r''$, et $z\in Y_{x,y}^{r}$ et $z''\in Y_{x,y}^{r''}$, 
$\geod(z,z'')$ rencontre $Y_{x,y}^{r'}$. 

\noindent b)  Etant donnés deux entiers $r,r'\in \{0,...,E(\frac{d(x,y)}{2})-3\delta\}$ avec $r\leq r'$, et $z'\in Y_{x,y}^{r'}$, $\geod(x,z')$ rencontre $Y_{x,y}^{r}$. 
\end{lem} 
\noindent{\bf Démonstration.} 
 Pour a), on  doit montrer l'existence de $z'\in \geod(z,z'')$ appartenant à $Y_{x,y}^{r'}$. 
Si $z$ ou $z''$ appartient à $Y_{x,y}^{r'}$ on prend $z'=z$ ou $z'=z''$. Supposons donc  $d(x,z)<r'$ et $d(x,z'')>r'+3\delta$. Soit $z'\in \geod(z,z'')$ tel que $d(x,z')=r'+E(3\delta/2)$ (un tel $z'$ existe car lorsque $z'$ parcourt $
\geod(z,z'')$ les valeurs prises par  $d(x,z')$ forment un intervalle de $\N$).  
On a $d(z,z')\geq 3\delta/2$ et $d(z',z'')\geq 3\delta/2$ donc le lemme~\ref{lemmetubes} implique $z'\in 3\delta\text{-}\geod(x,y)$. Pour b) 
on  doit montrer l'existence de $z\in \geod(x,z')$ appartenant à $Y_{x,y}^{r}$.
Si $d(x,z')\leq r+3\de$ on prend $z=z'$. Sinon soit $z\in \geod(x,z')$ tel que $d(x,z)=r$. Comme $d(z,z')\geq 3\de$, le lemme~\ref{contract1} montre que 
$z$ appartient à $ \de\tg(x,y)$ et donc à $Y_{x,y}^{r}$.  
  \cqfd

\label{Yrxy}
Ensuite pour $r,s\in  \{0,...,E(\frac{d(x,y)}{2})-3\delta\}$ on note 
$\Lambda_{x,r}^{y,s}$ l'ensemble  des familles d'entiers relatifs indexées par $Y_{x,y}^{r}\cup Y_{y,x}^{s}$, que nous notons $(c_{z})_{z\in   Y_{x,y}^{r}\cup Y_{y,x}^{s}}$, telles que 
\begin{gather*}\forall z\in Y_{x,y}^{r},\forall t\in   Y_{y,x}^{s}, \ \ c_{z}+c_{t}\geq d(z,t),\\ \forall z,z' \in Y_{x,y}^{r}, \ \ |c_{z}-c_{z'}|\leq d(z,z'),  \\ \forall t,t' \in   Y_{y,x}^{s}, \ \ |c_{t}-c_{t'}|\leq d(t,t'). \end{gather*}
En quelque sorte pour $z\in Y_{x,y}^{r}\cup Y_{y,x}^{s}$, $c_{z}$ est la ``distance'' d'un point virtuel à $z$ (on écrit ``distance'' car elle appartient à $\Z$ mais pas nécessairement à $\N$). 

Maintenant nous allons imposer une condition %pour que le point virtuel soit sur une géodésique entre des points réels, cette condition 
qui impliquera en particulier que le point virtuel ``appartient'' à  $3\delta\text{-}\geod(x,y)$. L'idée naïve serait de demander qu'il existe $z\in Y_{x,y}^{r}$ et $t \in Y_{y,x}^{s}$ tels que 
$c_{z}+c_{t}=d(z,t)$, c'est-à-dire que le point virtuel ``appartient'' à $\geod(z,t)$. Cependant cette condition nous empêcherait d'appliquer l'astuce des ensembles emboîtés. 
En effet on peut définir une application 
$\alpha_{r'\leftarrow r}^{s'\leftarrow s}: \Lambda_{x,r}^{y,s}\to \Lambda_{x,r'}^{y,s'}$, pour   $r,r',s,s'$   vérifiant  
\begin{gather}\label{ineg-r'rs's} 0\leq r'\leq r\leq E(\frac{d(x,y)}{2})-3\delta\text{\  et \ }0\leq s'\leq s\leq E(\frac{d(x,y)}{2})-3\delta,
\end{gather}  de la fa\c con suivante : à $c\in \Lambda_{x,r}^{y,s}$ on associe $c'=\alpha_{r'\leftarrow r}^{s'\leftarrow s}(c)\in \Lambda_{x,r'}^{y,s'}$ tel que 
\begin{gather*}\text{pour \ \ }z'\in Y_{x,y}^{r'},\  c'_{z'}=\min_{z\in Y_{x,y}^{r}}
d(z',z)+c_{z}\\ \text{et \  pour \ \ }t'\in Y_{y,x}^{s'}, \ c'_{t'}=\min_{t\in Y_{y,x}^{s}}
d(t',t)+c_{t}.\end{gather*}
En d'autres termes, pour définir $\alpha_{r'\leftarrow r}^{s'\leftarrow s}$ on fait comme si la géodésique entre le point virtuel et tout point de $Y_{x,y}^{r'}$ (resp. $Y_{y,x}^{s'}$) rencontrait $Y_{x,y}^{r}$ (resp. $Y_{y,x}^{s}$), comme dans le dessin ci-dessous. 

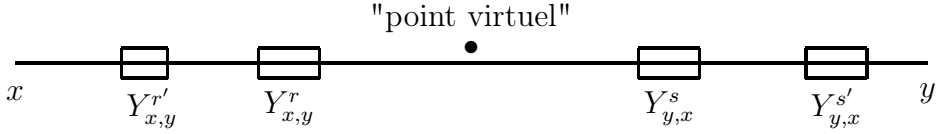
\ifx\JPicScale\undefined\def\JPicScale{1}\fi
\unitlength \JPicScale mm
\begin{picture}(140,35)(20,25)
\linethickness{0.3mm}
\put(20,40){\line(1,0){120}}
\linethickness{0.3mm}
\put(34,42){\line(1,0){6}}
\linethickness{0.3mm}
\put(40,38){\line(0,1){4}}
\linethickness{0.3mm}
\put(34,38){\line(1,0){6}}
\linethickness{0.3mm}
\put(34,38){\line(0,1){4}}
\linethickness{0.3mm}
\put(52,38){\line(0,1){4}}
\linethickness{0.3mm}
\put(52,38){\line(1,0){8}}
\linethickness{0.3mm}
\put(60,38){\line(0,1){4}}
\linethickness{0.3mm}
\put(52,42){\line(1,0){8}}
\linethickness{0.3mm}
\put(132,38){\line(0,1){4}}
\linethickness{0.3mm}
\put(124,38){\line(1,0){8}}
\linethickness{0.3mm}
\put(124,38){\line(0,1){4}}
\linethickness{0.3mm}
\put(124,42){\line(1,0){8}}
\linethickness{0.3mm}
\put(110,38){\line(0,1){4}}
\linethickness{0.3mm}
\put(102,38){\line(1,0){8}}
\linethickness{0.3mm}
\put(102,38){\line(0,1){4}}
\linethickness{0.3mm}
\put(102,42){\line(1,0){8}}
\put(20,36){\makebox(0,0)[cc]{$x$}}

\put(140,36){\makebox(0,0)[cc]{$y$}}

\put(128,34){\makebox(0,0)[cc]{$Y_{y,x}^{s'}$}}

\put(106,34){\makebox(0,0)[cc]{$Y_{y,x}^{s}$}}

\put(38,34){\makebox(0,0)[cc]{$Y_{x,y}^{r'}$}}

\put(56,34){\makebox(0,0)[cc]{$Y_{x,y}^{r}$}}

\put(80,42){\makebox(0,0)[cc]{$\bullet$}}

\put(80,46){\makebox(0,0)[cc]{$\text{"point virtuel"}$}}

\end{picture}

 Le a) du lemme~\ref{tubes} montre que ces applications se composent bien : si $r,r',r'',s,s',s''$  vérifient 
$$0\leq r''\leq r'\leq r\leq E(\frac{d(x,y)}{2})-3\delta \text{ \ et \ }
0\leq s''\leq s'\leq s\leq E(\frac{d(x,y)}{2})-3\delta ,$$ on a 
$\alpha_{r''\leftarrow r'}^{s''\leftarrow s'}\circ \alpha_{r'\leftarrow r}^{s'\leftarrow s}=\alpha_{r''\leftarrow r}^{s''\leftarrow s}$. 
Si on imposait en plus la condition évoquée ci-dessus, on ne posséderait plus de telles applications. On va imposer cette condition en d'autres entiers  $r',s'$, en utilisant précisément l'application $\alpha_{r'\leftarrow r}^{s'\leftarrow s}$, et dans l'astuce des ensembles emboîtés, $r$ et $r'$, respectivement $s$ et $s'$,  varieront en sens inverse l'un de l'autre (c'est pourquoi on ne peut pas les réunir en une seule variable).  Pour $r,r',s,s'$ vérifiant (\ref{ineg-r'rs's})
%$$0\leq  r'\leq r\leq E(\frac{d(x,y)}{2})-3\delta \text{ \ et \ }
%0\leq s'\leq s\leq E(\frac{d(x,y)}{2})-3\delta ,$$
on définit donc \label{Lambdaxysrsr}
$\Lambda_{x,r',r}^{y,s',s}$ comme l'ensemble 
des $c\in \Lambda_{x,r}^{y,s}$ dont l'image $c'=\alpha_{r'\leftarrow r}^{s'\leftarrow s}(c)\in \Lambda_{x,r'}^{y,s'}$  vérifie la condition introduite ci-dessus : il existe 
$z'\in Y_{x,y}^{r'}$ et $t' \in Y_{y,x}^{s'}$ tels que 
$c'_{z'}+c'_{t'}=d(z',t')$. 

\begin{lem}\label{inclusions-var-r2}
Pour $r,s\in \{0,\pp, E(\frac{d(x,y)}{2})-3\delta\}$ on a des inclusions 
$$\begin{matrix}\Lambda_{x,0,r}^{y,0,s} &\subset & \pp & \subset & \Lambda_{x,r-1,r}^{y,0,s} & \subset&  \Lambda_{x,r,r}^{y,0,s}\\
\cap &&\pp && \cap && \cap \\
\pp&&\pp&&\pp&&\pp \\
\cap &&\pp && \cap && \cap \\
\Lambda_{x,0,r}^{y,s-1,s} &\subset & \pp & \subset & \Lambda_{x,r-1,r}^{y,s-1,s} & \subset&  \Lambda_{x,r,r}^{y,s-1,s}\\
\cap &&\pp && \cap && \cap \\
\Lambda_{x,0,r}^{y,s,s} &\subset & \pp & \subset & \Lambda_{x,r-1,r}^{y,s,s} & \subset&  \Lambda_{x,r,r}^{y,s,s}\\
 \end{matrix}
$$ dans $\Lambda_{x,r}^{y,s}$ et toutes ces parties sont non vides. 

\end{lem}
\noindent{\bf Démonstration.} 
Soit $r'\in \{0,\pp, r-1\}$, $s'\in \{0,\pp, s\}$ et montrons l'inclusion horizontale  $\Lambda_{x,r',r}^{y,s',s} \subset \Lambda_{x,r'+1,r}^{y,s',s}$. Soit $c\in \Lambda_{x,r',r}^{y,s',s}$. On note 
$c'=\alpha_{r'\leftarrow r}^{s'\leftarrow s}(c)\in \Lambda_{x,r'}^{y,s'}$  et $c''=\alpha_{r'+1\leftarrow r}^{s'\leftarrow s}(c)\in \Lambda_{x,r'+1}^{y,s'}$, si bien que $c'=\alpha_{r'\leftarrow r'+1}^{s'\leftarrow s'}(c'')$. Par hypothèse il existe $u'\in Y_{x,y}^{r'}$ et $v'\in Y_{y,x}^{s'}$ tels que $c'_{u'}+c'_{v'}=d(u',v')$. Par définition de l'application $\alpha_{r'\leftarrow r'+1}^{s'\leftarrow s'}$ il existe 
$u''\in Y_{x,y}^{r'+1}$ tel que $c'_{u'}=c''_{u''}+d(u',u'')$. Les inégalités 
\begin{gather*}d(u',v')=c'_{u'}+c'_{v'}=c''_{u''}+d(u',u'')+c'_{v'}=d(u',u'')+c''_{u''}+c''_{v'}\\ \geq d(u',u'')+d(u'',v')\geq d(u',v')\end{gather*} sont toutes des égalités, donc $c''_{u''}+c''_{v'}=d(u'',v')$ et $c$ appartient bien à $\Lambda_{x,r'+1,r}^{y,s',s}$. Les inclusions verticales se démontrent de la même manière. 
Il reste à montrer que $\Lambda_{x,0,r}^{y,0,s}$ est non vide. Soit $t\in \geod(x,y)$ tel que $d(x,t)\in [r,d(x,y)-s]$ et soit $c\in \Lambda^{y,s}_{x,r}$ défini par $c_{z}=d(z,t)$ pour $z\in Y_{x,y}^{r}\cup Y_{y,x}^{s}$. Alors $c$ appartient à 
$\Lambda_{x,0,r}^{y,0,s}$ car $c'=\alpha_{0\leftarrow r}^{0\leftarrow s}(c)$ vérifie $c'_{x}=d(x,t)$ et $c'_{y}=d(t,y)$ (puisque $\geod(x,t)$ rencontre $Y_{x,y}^{r}$ et $\geod(t,y)$ rencontre $Y_{y,x}^{s}$) et donc $c'_{x}+c'_{y}=d(x,y)$.   
Par conséquent  $\Lambda_{x,0,r}^{y,0,s}$ est non vide. 
\cqfd

Enfin étant donné  $r_{1},r_{2},r_{3},s_{1},s_{2},s_{3}$ vérifiant 
\begin{gather}\label{rrrsss-existence}0\leq r_{1}\leq r_{2}\leq r_{3}\leq E(\frac{d(x,y)}{2})-3\delta\text{ et } 
0\leq s_{1}\leq s_{2}\leq s_{3}\leq E(\frac{d(x,y)}{2})-3\de, 
\end{gather}  \label{Lambda123}
 on définit  $\Lambda_{x,r_{1},r_{2},r_{3}}^{y,s_{1},s_{2},s_{3}}$ comme l'image de $\Lambda_{x,r_{2},r_{3}}^{y,s_{2},s_{3}}$ dans $\Lambda_{x,r_{1}}^{y,s_{1}}$ par l'application $\alpha_{r_{1}\leftarrow r_{3}}^{s_{1}\leftarrow s_{3}}$. 
Pour tout entier $r\in \{0,...,E(\frac{d(x,y)}{2})\}-3\de$
on a une application $\beta_{x,r}^{y,s}:\Lambda_{x,r}^{y,s}\to \N$ définie par  \begin{gather*}\beta_{x,r}^{y,s}(c)=\min 
_{z\in Y_{x,y}^{r},t\in Y_{y,x}^{s}}d(x,z)+c_{z}+c_{t}+d(t,y)
\\ =\min 
_{z\in Y_{x,y}^{r}}\big(d(x,z)+c_{z}\big)+\min 
_{t\in Y_{y,x}^{s}}\big(c_{t}+d(t,y)\big). \end{gather*}
On doit comprendre $\min 
_{z\in Y_{x,y}^{r}}\big(d(x,z)+c_{z}\big)$, respectivement 
$\min 
_{t\in Y_{y,x}^{s}}\big(c_{t}+d(t,y)\big)$, comme la ``distance'' de $x$, respectivement $y$,  au point virtuel. 
Si  $r,r',s,s'$ sont des entiers vérifiant (\ref{ineg-r'rs's}),
%$$0\leq r'\leq r\leq E(\frac{d(x,y)}{2})-3\delta\text{\ \ et\ \ }
%0\leq s'\leq s\leq E(\frac{d(x,y)}{2})-3\delta
%$$ 
on a $\beta_{x,r'}^{y,s'}\circ \alpha_{r'\leftarrow r}^{s'\leftarrow s}=\beta_{x,r}^{y,s}$ par le b) du lemme~\ref{tubes}. 
Cette application $\beta_{x,r}^{y,s}$ se factorise par $\tilde \Lambda_{x,r}^{y,s}$, où $\tilde \Lambda_{x,r}^{y,s}$ est le quotient de $\Lambda_{x,r}^{y,s}$ par la relation d'équivalence suivante : 
deux éléments $c$ et $c'$ sont équivalents s'il existe $k\in \Z$ avec $c'(z)=c(z)+k$ pour $z\in Y_{x,y}^{r}$ et 
$c'(t)=c(t)-k$ pour $t\in Y_{y,x}^{s}$. Toutes les constructions ci-dessus passent au quotient de cette fa\c con et on note 
en particulier $\tilde \Lambda_{x,r_{1},r_{2},r_{3}}^{y,s_{1},s_{2},s_{3}}$ l'image de 
$\tilde \Lambda_{x,r_{2},r_{3}}^{y,s_{2},s_{3}}$ dans $\tilde \Lambda_{x,r_{1}}^{y,s_{1}}$ (qui est aussi l'image de $\Lambda_{x,r_{2},r_{3}}^{y,s_{2},s_{3}}$ dans $\tilde \Lambda_{x,r_{1}}^{y,s_{1}}$). On remarque  que  si $r,r',s,s'$ vérifient (\ref{ineg-r'rs's}), 
%$$0\leq r'\leq r\leq E(\frac{d(x,y)}{2})-3\delta\text{\ \ et \ \ } 
%0\leq s'\leq s\leq E(\frac{d(x,y)}{2})-3\delta,$$ 
  $\tilde \Lambda_{x,r',r}^{y,s',s}$ est un ensemble fini, dont le cardinal est borné par une constante de la forme $C(\delta,K)$. Il en va donc de même pour $\tilde \Lambda_{x,r_{1},r_{2},r_{3}}^{y,s_{1},s_{2},s_{3}}$, lorsque  $r_{1},r_{2},r_{3},s_{1},s_{2},s_{3}$ vérifient (\ref{rrrsss-existence}). 
Grâce à $\alpha_{r_{1}-1\leftarrow r_{1}}^{s_{1}\leftarrow s_{1}},\alpha_{r_{1}\leftarrow r_{1}}^{s_{1}-1\leftarrow s_{1}},\alpha_{r_{3}-1\leftarrow r_{3}}^{s_{3}\leftarrow s_{3}}, \alpha_{r_{3}\leftarrow r_{3}}^{s_{3}-1\leftarrow s_{3}}$  et au lemme~\ref{inclusions-var-r2} on possède, pour \break $(r_{1},  r_{2},r_{3},s_{1},s_{2},s_{3})$ vérifiant (\ref{rrrsss-existence}), 
\begin{itemize} 
\item une application surjective $\tilde \Lambda_{x,r_{1},r_{2},r_{3}}^{y,s_{1},s_{2},s_{3}}\to \tilde \Lambda_{x,r_{1}-1,r_{2},r_{3}}^{y,s_{1},s_{2},s_{3}}$ si $r_{1}\geq 1$, 
\item une application surjective $\tilde \Lambda_{x,r_{1},r_{2},r_{3}}^{y,s_{1},s_{2},s_{3}}\to \tilde \Lambda_{x,r_{1},r_{2},r_{3}}^{y,s_{1}-1,s_{2},s_{3}}$ si $s_{1}\geq 1$, 
\item des applications injectives $\tilde \Lambda_{x,r_{1},r_{2},r_{3}}^{y,s_{1},s_{2},s_{3}}\to \tilde \Lambda_{x,r_{1},r_{2}+1,r_{3}}^{y,s_{1},s_{2},s_{3}}$ et $\tilde \Lambda_{x,r_{1},r_{2},r_{3}}^{y,s_{1},s_{2},s_{3}}\to \tilde \Lambda_{x,r_{1},r_{2},r_{3}-1}^{y,s_{1},s_{2},s_{3}}$ si $r_{2}< r_{3}$,  
\item des applications injectives $\tilde \Lambda_{x,r_{1},r_{2},r_{3}}^{y,s_{1},s_{2},s_{3}}\to \tilde \Lambda_{x,r_{1},r_{2},r_{3}}^{y,s_{1},s_{2}+1,s_{3}}$ et $\tilde \Lambda_{x,r_{1},r_{2},r_{3}}^{y,s_{1},s_{2},s_{3}}\to \tilde \Lambda_{x,r_{1},r_{2},r_{3}}^{y,s_{1},s_{2},s_{3}-1}$ si $s_{2}< s_{3}$. \end{itemize}

De plus ces applications sont compatibles entre elles. 

On est maintenant en mesure de construire $d^{\flat}$. 
 On pose  $d^{\flat}(x,y)=d(x,y)$ si $d(x,y)<6\delta$ et si $d(x,y)\geq 6\delta$
on pose $\Delta_{x,y}=E(d(x,y)/6)-\delta$ et on définit 
\begin{gather} \nonumber d^{\flat}(x,y)=\frac{ 1}{  (\Delta_{x,y}+1)^{6}} \sum _{%\begin{matrix}  \scriptstyle 
r_{1},s_{1}\in \{0,...,\Delta_{x,y}\},
%\\ \scriptstyle  
r_{2},s_{2}\in \{\Delta_{x,y},...,2\Delta_{x,y}\},
%\\ \scriptstyle
 r_{3},s_{3}\in \{2\Delta_{x,y},3\Delta_{x,y}\}
 %\end{matrix}
 }\\ 
\label{formule-d'}\frac{1}{\sharp (\tilde \Lambda_{x,r_{1},r_{2},r_{3}}^{y,s_{1},s_{2},s_{3}})}
\sum_{c\in \tilde \Lambda_{x,r_{1},r_{2},r_{3}}^{y,s_{1},s_{2},s_{3}}}\beta_{x,r_{1}}^{y,s_{1}}(c).\end{gather}

Pour $u_{1},u_{2},u_{3},v_{1},v_{2},v_{3}\in [0,1[$ et $x,y\in X$ on pose 
${d^{\flat}}_{u_{1},u_{2},u_{3}}^{v_{1},v_{2},v_{3}}(x,y)=d(x,y)$ si $d(x,y)<6\delta$ et si $d(x,y)\geq 6\delta$
on pose  \label{def-A-Lambda-dbemol}
\begin{gather}
\label{def-Auuuvvv}A_{x,u_{1},u_{2},u_{3}}^{y,v_{1},v_{2},v_{3}}=\tilde \Lambda_{x,E((\Delta_{x,y}+1)u_{1}),\Delta_{x,y}+E((\Delta_{x,y}+1)u_{2}),2\Delta_{x,y}+E((\Delta_{x,y}+1)u_{3})}^{y,E((\Delta_{x,y}+1)v_{1}),\Delta_{x,y}+E((\Delta_{x,y}+1)v_{2}),2\Delta_{x,y}+E((\Delta_{x,y}+1)v_{3})}, 
\\ \nonumber \text{et\ \ }{d^{\flat}}_{u_{1},u_{2},u_{3}}^{v_{1},v_{2},v_{3}}(x,y)=\frac{1}{\sharp (A_{x,u_{1},u_{2},u_{3}}^{y,v_{1},v_{2},v_{3}})}\sum_{c\in A_{x,u_{1},u_{2},u_{3}}^{y,v_{1},v_{2},v_{3}}}\beta_{x,E((\Delta_{x,y}+1)u_{1})}^{y,E((\Delta_{x,y}+1)v_{1})}(c)\end{gather} 
de sorte que d'après la formule (\ref{formule-d'}) on a toujours 
$$ d^{\flat}(x,y)=\int_{u_{1},u_{2},u_{3},v_{1},v_{2},v_{3}\in [0,1[}{d^{\flat}}_{u_{1},u_{2},u_{3}}^{v_{1},v_{2},v_{3}}(x,y)du_{1}du_{2}du_{3}dv_{1}dv_{2}dv_{3}. $$
Dans la formule précédente on intègre sur $[0,1[^{6}$ au lieu de $[0,1]^{6}$ car l'expression pourrait ne pas avoir de sens pour $u_{1}=1$ et $u_{2}=0$ par exemple.

\begin{lem}\label{123-7de} 
 a) 
Pour $r_{1},s_{1}\in \{0,...,\Delta_{x,y}\},
%\\ \scriptstyle  
r_{2},s_{2}\in \{\Delta_{x,y},...,2\Delta_{x,y}\},
%\\ \scriptstyle
 r_{3},s_{3}\in \{2\Delta_{x,y},...,3\Delta_{x,y}\},$ $ c \in \tilde \Lambda_{x,r_{1},r_{2},r_{3}}^{y,s_{1},s_{2},s_{3}}$ on a $\beta_{x,r_{1}}^{y,s_{1}}(c)\in [d(x,y),d(x,y)+7\de]$. 
 
 \noindent b) Pour  $u_{1},u_{2},u_{3},v_{1},v_{2},v_{3}\in [0,1[$ et $x,y\in X$ on a 
$$d(x,y)\leq {d^{\flat}}_{u_{1},u_{2},u_{3}}^{v_{1},v_{2},v_{3}}(x,y) \leq d(x,y)+7\de.$$
\end{lem}
\noindent{\bf Démonstration.} On montre seulement a) car b) en résulte immédiatement. 
Il existe $z\in Y_{x,y}^{r_{2}}$ et $t\in Y_{y,x}^{s_{2}}$ tels que 
$\beta_{x,r_{1}}^{y,s_{1}}(c)\leq d(x,z)+d(z,t)+d(t,y)$. On rappelle que $z,t$ appartiennent à $3\de\tg(x,y)$ et que $d(x,z)\leq \frac{d(x,y)}{2}$ et $d(t,y)\leq \frac{d(x,y)}{2}$. Par $(H_{\de}(z,x,t,y))$ on a 
\begin{gather*}d(z,t)+d(x,y)\leq \max(d(x,z)+d(t,y),d(x,t)+d(z,y))+\de
\\ =d(x,t)+d(z,y)+\de \text{\ \ \ d'où}
\\ d(x,z)+d(z,t)+d(t,y)\leq d(x,z)+d(x,t)+d(z,y)+d(t,y)+\de-d(x,y)
\\ \leq 
d(x,y)+7\de. \end{gather*}\cqfd

\begin{lem}\label{mesure-uuuvvv-dbemol}
 Pour tout $\rho\in \N$, il existe $C=C(\de,K,\rho)$ tel que pour 
 $x,x',y,y'\in X$ verifiant $d(x,x')\leq \rho$ et $d(y,y')\leq \rho$ la mesure de l'ensemble des \begin{gather*}(u_{1},u_{2},u_{3},v_{1},v_{2},v_{3})\in [0,1[^{6}\text{ tels que  }\\ {d^{\flat}}_{u_{1},u_{2},u_{3}}^{v_{1},v_{2},v_{3}}(x,y)
 -{d^{\flat}}_{u_{1},u_{2},u_{3}}^{v_{1},v_{2},v_{3}}(x',y)-{d^{\flat}}_{u_{1},u_{2},u_{3}}^{v_{1},v_{2},v_{3}}(x,y')
 + {d^{\flat}}_{u_{1},u_{2},u_{3}}^{v_{1},v_{2},v_{3}}(x',y')\neq 0\end{gather*} est 
 $\leq \frac{C}{1+d(x,y)}$.
\end{lem}
\noindent {\bf Démonstration.} 
La démonstration  consiste encore en une astuce d'ensembles emboîtés. 

\begin{souslem}\label{tubes-xyx'y'}  Soient $x,y,\tilde x,\tilde y\in X$ et $d\geq 2(d(x,\tilde x)+d(y,\tilde y)+3\delta)$. 

\noindent a) 
 Etant donnés  $r,r',r''\in \{0,...,E(\frac{d(x,y)}{2})-3\delta\}$ tels que $r+d\leq r' \leq r''-d$, 
 $z\in Y_{x,y}^{r}$ et $z''\in Y_{x,y}^{r''}$, 
$\geod(z,z'')$ rencontre $Y_{\tilde x,\tilde y}^{r'}$. 

\noindent b) Etant donnés  $r,r',r''\in \{0,...,E(\frac{d(x,y)}{2})-3\delta\}$ tels que $r\leq r' \leq r''-d$, 
 $z\in Y_{\tilde x,\tilde y}^{r}$ et $z''\in Y_{x,y}^{r''}$, 
$\geod(z,z'')$ rencontre $Y_{\tilde x,\tilde y}^{r'}$. 

\noindent c) Etant donnés  $r,r',r''\in \{0,...,E(\frac{d(x,y)}{2})-3\delta\}$ tels que $r+d\leq r' \leq r''$, 
 $z\in Y_{\tilde x,\tilde y}^{r}$ et $z''\in Y_{x,y}^{r''}$, 
$\geod(z,z'')$ rencontre $Y_{x,y}^{r'}$.

\noindent d) Etant donnés $r,r'\in \{0,...,E(\frac{d(x,y)}{2})-3\delta\}$ tels que $r+d\leq r' $, 
 $z'\in Y_{x,y}^{r'}$, 
$\geod(\tilde x,z')$ rencontre $Y_{\tilde x,\tilde y}^{r}$.
\end{souslem} 
 \noindent{\bf Démonstration.}  Pour a) on prend $z'\in \geod(z,z'')$ tel que $d(\tilde x,z')=r'$. Cela est possible car $d(\tilde x,z)\leq d(x,\tilde x)+r+3\de\leq r'-\frac{d}{2}$ et $d(\tilde x,z'')\geq r''-d(x,\tilde x)\geq r'+\frac{d}{2}$. 
  On a $z,z''\in (d-3\de)\tg(\tilde x,\tilde y)$ par le lemme~\ref{xx'yy'zz'} et $d(z,z')\geq \frac{d}{2},  d(z',z'')\geq \frac{d}{2}$ ce qui permet d'appliquer le  lemme~\ref{lemmetubes}
avec $\alpha=\beta=d-3\de$. D'où $z'\in Y_{\tilde x,\tilde y}^{r'}$. 

Pour  b) on prend $z'=z$ si $d(\tilde x,z)\geq r'$ donc on suppose $d(\tilde x,z)< r'$. On a $$d(\tilde x,z'')\geq d(x,z'')-d(x,\tilde x)\geq 
r''-\frac{d}{2}+3\de\geq r'+\frac{d}{2}+3\de$$
 donc il existe $z'\in \geod(z,z'')$ tel que $d(\tilde x,z')=r'+3\de$. On a alors $z\in 3\de\tg(\tilde x, \tilde y)$, $z''\in (d-3\de)\tg(\tilde x, \tilde y)$, $d(z,z')\geq 3\de$ et 
$d(z',z'')\geq \frac{d}{2}$, ce qui permet d'appliquer le lemme~\ref{lemmetubes} avec $\alpha=3\de$ et $\beta=d-3\de$. D'où $z'\in Y_{\tilde x,\tilde y}^{r'}$.

Pour c)  on prend $z'=z''$ si $d( x,z'')\leq r'+3\de$ donc on suppose $d( x,z'')> r'+3\de$. On a $$d(x,z)\leq d(x,\tilde x)+d(\tilde x,z)\leq r+3\de+d(x,\tilde x)\leq r+\frac{d}{2}\leq r'-\frac{d}{2}.$$ Donc il existe $z'\in \geod(z,z'')$ tel que $d(x,z')= r'$. On a $z\in (d-3\de)\tg(x,y)$, $z''\in 3\de\tg(x,y)$, $d(z,z')\geq \frac{d}{2}$, $d(z',z'')\geq 
3\de$ ce qui permet d'appliquer le lemme~\ref{lemmetubes} avec $\alpha=d-3\de$ et $\beta=3\de$. D'où $z'\in Y_{x,y}^{r'}$. 

Pour d)  on prend $z\in \geod(\tilde x,z')$ vérifiant $d(\tilde x,z)=r$. Cela est possible car $d(\tilde x,z')\geq d(x,z')-d(x,\tilde x)\geq r'-\frac{d}{2}\geq r+\frac{d}{2}$. Comme $z'\in (d-3\de)\tg(\tilde x,\tilde y)$, $z\in \geod(\tilde x,z')$ et $d(z,z')\geq \frac{d}{2}$, le lemme~\ref{contract1} appliqué à $\epsilon=d-3\de$  montre 
$z\in 3\de\tg(\tilde x, \tilde y)$ d'où $z\in Y^{r}_{\tilde x,\tilde y}$. 
 \cqfd 

\noindent{\bf Suite de la démonstration du lemme~\ref{mesure-uuuvvv-dbemol}.} 
Pour $x,y,\tilde x,\tilde y\in X$ et $d\geq 2( d(x,\tilde x)+d(y,\tilde y)+3\delta)$ comme dans le sous-lemme~\ref{tubes-xyx'y'}, et pour $r,s\in \{d,\pp, E(\frac{d(x,y)}{2})-3\delta\}$, on a une application $\gamma:\Lambda_{x,r}^{y,s}\to \Lambda_{\tilde x,r-d}^{\tilde y,s-d}$ qui envoie $c$ sur $c'$ défini par 
\begin{gather*}c'_{z'}=\min_{z\in  Y_{x,y}^{r}} c_{z}+d(z',z)\text{\  pour \ }z'\in Y_{\tilde x,\tilde y}^{r-d} \\ \text{et\ } 
c'_{t'}=\min_{t\in  Y_{y,x}^{s}} c_{t}+d(t',t)\text{\  pour \ }t'\in Y_{\tilde y,\tilde x}^{s-d}.\end{gather*} 
On note $\gamma$ sans indices pour alléger les formules car les indices  sont déterminés par l'ensemble de départ et l'ensemble d'arrivée. 
De 
 la même fa\c con, pour $r,s\in \{d,\pp, E(\frac{d(\tilde x,\tilde y)}{2})-3\delta\}$, on a une application
$\gamma: \Lambda_{\tilde x,r}^{\tilde y,s}\to \Lambda_{ x,r-d}^{y,s-d}$. Ces applications sont compatibles avec les applications $\alpha_{r'\leftarrow r}^{s'\leftarrow s}$ (associées à $(x,y)$ et $(\tilde x,\tilde y)$). Par exemple, lorsque les applications ont un sens, 
\begin{itemize}
\item 
 la composée $\Lambda_{x,r}^{y,s}\vad{\gamma} \Lambda_{\tilde x,r-d}^{\tilde y,s-d}\vad{\gamma} \Lambda_{x,r-2d}^{y,s-2d}$ coïncide avec $\alpha_{r-2d\leftarrow r}^{s-2d\leftarrow s}$ (grâce à a) du sous-lemme~\ref{tubes-xyx'y'}), 
\item 
 les deux composées $$ \Lambda_{x,r}^{y,s}\vad{\gamma} \Lambda_{\tilde x,r-d}^{\tilde y,s-d}\vad{\alpha_{r'-d\leftarrow r-d}^{s'-d\leftarrow s-d} } \Lambda_{\tilde x,r'-d} ^{\tilde y,s'-d}\text{ \ \  et \ \ } 
 \Lambda_{x,r}^{y,s}\vad{\alpha_{r'\leftarrow r}^{s'\leftarrow s}}   \Lambda_{x,r'}^{y,s'}\vad{\gamma}  \Lambda_{\tilde x,r'-d} ^{\tilde y,s'-d}$$ sont égales car, grâce à b) et c) du sous-lemme~\ref{tubes-xyx'y'}
l'image $c'\in  \Lambda_{\tilde x,r'-d} ^{\tilde y,s'-d}$ de $c\in  \Lambda_{x,r}^{y,s}$ par chacune de ces deux composées est donnée par 
\begin{gather*}c'_{z'}=\min_{z\in  Y_{x,y}^{r}} c_{z}+d(z',z)\text{\  pour \ }z'\in Y_{\tilde x,\tilde y}^{r'-d} \\ \text{et\ } 
c'_{t'}=\min_{t\in  Y_{y,x}^{s}} c_{t}+d(t',t)\text{\  pour \ }t'\in Y_{\tilde y,\tilde x}^{s'-d}.\end{gather*}
\end{itemize}

Plus généralement la compatibilité signifie que deux composées construites à partir des applications $\alpha$ et $\gamma$ sont égales lorsqu'elles possèdent le même ensemble de départ et  le même ensemble d'arrivée. En effet on démontre que si l'ensemble de départ est 
$\Lambda_{\check x,\check r}^{\check y,\check s}$ et l'ensemble d'arrivée est $\Lambda_{\hat x,\hat r}^{\hat y,\hat s}$
(avec $(\check x,\check y)$ et $(\hat x,\hat y)$ égaux à $(x,y)$ ou $(\tilde x,\tilde y)$), les deux composées sont égales à l'application qui à $\check c\in \Lambda_{\check x,\check r}^{\check y,\check s}$ associe $\hat c\in \Lambda_{\hat x,\hat r}^{\hat y,\hat s}$ défini par les formules 
\begin{gather*}\hat c_{\hat z}=\min_{\check z\in  Y_{\check x,\check y}^{\check r}} \check c_{\check z}+d(\hat z,\check z)\text{\  pour \ }\hat z\in Y_{\hat  x,\hat  y}^{\hat r} \\ \text{et\ } 
\hat c_{\hat t}=\min_{\check t\in  Y_{\check y,\check x}^{\check s}} \check c_{\check t}+d(\hat t,\check t)\text{\  pour \ }\hat t\in Y_{\hat y,\hat x}^{\hat s}.\end{gather*}

\begin{souslem}\label{Lambda-tilde-rsr's'-11j}
Pour $r,r',s,s'$ comme dans (\ref{ineg-r'rs's}) vérifiant $r'+2d\leq r$ et $s'+2d\leq s$, l'application  $\Lambda_{ x,r}^{ y,s}\vad{\gamma} \Lambda_{\tilde x,r-d}^{\tilde y,s-d}$
 envoie 
$\Lambda_{ x,r',r}^{ y,s',s}$ dans  $\Lambda_{ \tilde x,r'+d,r-d}^{\tilde y,s'+d,s-d}$. 
\end{souslem}
\noindent{\bf Démonstration.}
D'après ce qui précède, $\alpha_{r'\leftarrow r}^{s'\leftarrow s}:\Lambda_{x,r}^{y,s} \to  \Lambda_{x,r'}^{y,s'}$ est la composée 
$$\Lambda_{x,r}^{y,s}\vad{\gamma} \Lambda_{\tilde x,r-d}^{\tilde y,s-d}\vad{\alpha_{r'+d\leftarrow r-d}^{s'+d\leftarrow s-d} } \Lambda_{\tilde x,r'+d} ^{\tilde y,s'+d}
\vad{\gamma} \Lambda_{x,r'}^{y,s'}.$$
De plus pour $c\in \Lambda_{\tilde x,r'+d} ^{\tilde y,s'+d}$ d'image $c'=\gamma(c)\in \Lambda_{x,r'}^{y,s'}$, s'il existe $z\in Y_{x,y}^{r'}$ et $t\in Y_{y,x}^{s'}$ tels que $c'_{z}+c'_{t}=d(z,t)$, alors il existe $\hat z\in Y_{\tilde x,\tilde y}^{r'+d}$ et $\hat t\in Y_{\tilde y,\tilde x}^{s'+d}$ tels que $c_{\hat z}+c_{\hat t}=d(\hat z,\hat t)$ par le même argument que dans la preuve du lemme~\ref{inclusions-var-r2}. 
\cqfd

On  déduit du sous-lemme~\ref{Lambda-tilde-rsr's'-11j} que pour $(r_{1},  r_{2},r_{3},s_{1},s_{2},s_{3})$ vérifiant (\ref{rrrsss-existence}) et $(\tilde r_{1},  \tilde r_{2},\tilde r_{3},\tilde s_{1},\tilde s_{2},\tilde s_{3})$ vérifiant $$0\leq \tilde r_{1}\leq \tilde r_{2}\leq \tilde r_{3}\leq E(\frac{d(\tilde x,\tilde y)}{2})-3\delta\text{\ \ et \ \ } 
0\leq \tilde s_{1}\leq \tilde s_{2}\leq \tilde s_{3}\leq E(\frac{d(\tilde x,\tilde y)}{2})-3\de
$$
\begin{itemize}
\item si $\tilde r_{1}\leq r_{1}-d$, $\tilde s_{1}\leq s_{1}-d$, 
$\tilde r_{2}\geq r_{2}+d$, $\tilde s_{2}\geq s_{2}+d$,
$\tilde r_{3}\leq r_{3}-d$, $\tilde s_{3}\leq s_{3}-d$, on a 
une application $\Lambda_{x,r_{1},r_{2},r_{3}}^{y,s_{1},s_{2},s_{3}}\vad{\gamma} 
\Lambda_{\tilde x,\tilde r_{1},  \tilde r_{2},\tilde r_{3}}^{\tilde y,\tilde s_{1},\tilde s_{2},\tilde s_{3}}$, et en notant $\tilde c=\gamma(c)\in \Lambda_{\tilde x,\tilde r_{1},  \tilde r_{2},\tilde r_{3}}^{\tilde y,\tilde s_{1},\tilde s_{2},\tilde s_{3}}$ l'image de $c\in \Lambda_{x,r_{1},r_{2},r_{3}}^{y,s_{1},s_{2},s_{3}}$ on a, grâce au d) du sous-lemme~\ref{tubes-xyx'y'},  
\begin{gather}\label{compat-alpha-beta-22oct09}
\beta_{\tilde x,\tilde r_{1}}^{\tilde y, \tilde s_{1}}(\tilde c)=\min_{z\in Y_{x,y}^{r_{1}}}\big(d(\tilde x,z)+c_{z}\big)+\min_{t\in Y_{y,x}^{s_{1}}}\big(c_{t}+d(t,\tilde y)\big),\end{gather}
\item si $r_{1}\leq \tilde  r_{1}-d$, $s_{1}\leq \tilde s_{1}-d$, 
$r_{2}\geq \tilde r_{2}+d$, $s_{2}\geq \tilde s_{2}+d$,
$r_{3}\leq \tilde r_{3}-d$, $s_{3}\leq \tilde s_{3}-d$, on a 
une application $
\Lambda_{\tilde x,\tilde r_{1},  \tilde r_{2},\tilde r_{3}}^{\tilde y,\tilde s_{1},\tilde s_{2},\tilde s_{3}}\vad{\gamma} \Lambda_{x,r_{1},r_{2},r_{3}}^{y,s_{1},s_{2},s_{3}}$ verifiant une propriété semblable à (\ref{compat-alpha-beta-22oct09}).  
\end{itemize}

On rappelle que  
$\Delta_{x,y}=E(d(x,y)/6)-\delta$,    $\Delta_{\tilde x,\tilde y}=E(d(\tilde x,\tilde y)/6)-\delta$,   
et que pour $(u_{1},u_{2},u_{3},v_{1},v_{2},v_{3})\in [0,1[^{6}$, 
 \begin{gather}
 \label{def-Auuuvvv2}A_{x,u_{1},u_{2},u_{3}}^{y,v_{1},v_{2},v_{3}}=\tilde \Lambda_{x,E((\Delta_{x,y}+1)u_{1}),\Delta_{x,y}+E((\Delta_{x,y}+1)u_{2}),2\Delta_{x,y}+E((\Delta_{x,y}+1)u_{3})}^{y,E((\Delta_{x,y}+1)v_{1}),\Delta_{x,y}+E((\Delta_{x,y}+1)v_{2}),2\Delta_{x,y}+E((\Delta_{x,y}+1)v_{3})}, \\ \label{def-tildeAuuuvvv}A_{\tilde x,u_{1},u_{2},u_{3}}^{\tilde y,v_{1},v_{2},v_{3}}=\tilde \Lambda_{\tilde x,E((\Delta_{\tilde x,\tilde y}+1)u_{1}),\Delta_{\tilde x,\tilde y}+E((\Delta_{\tilde x,\tilde y}+1)u_{2}),2\Delta_{\tilde x,\tilde y}+E((\Delta_{\tilde x,\tilde y}+1)u_{3})}^{\tilde y,E((\Delta_{\tilde x,\tilde y}+1)v_{1}),\Delta_{\tilde x,\tilde y}+E((\Delta_{\tilde x,\tilde y}+1)v_{2}),2\Delta_{\tilde x,\tilde y}+E((\Delta_{\tilde x,\tilde y}+1)v_{3})} .\end{gather}
Soit $\eta\in ]0,\frac{1}{2}[$ tel que  $d(x,y)\geq \frac{9d}{\eta}+6\de$. On a alors pour tout $(u_{1},u_{2},u_{3},v_{1},v_{2},v_{3})\in [\eta,1-\eta[^{6}$, les six  inégalités 
\begin{gather*}E((\Delta_{\tilde x,\tilde y}+1)(u_{1}-\eta))\leq E((\Delta_{x,y}+1)u_{1})-d,\\
\Delta_{\tilde x,\tilde y}+E((\Delta_{\tilde x,\tilde y}+1)(u_{2}+\eta))\geq \Delta_{x,y}+E((\Delta_{x,y}+1)u_{2})+d,   ...,\\ 
 2\Delta_{\tilde x,\tilde y}+E((\Delta_{\tilde x,\tilde y}+1)(v_{3}-\eta))\leq 2\Delta_{x,y}+E((\Delta_{x,y}+1)v_{3})-d.\end{gather*} 
 Ces inégalités ont lieu car 
  \begin{gather*}|\Delta_{x,y}-\Delta_{\tilde x,\tilde y}|\leq 
  \frac{|d(x,y)-d(\tilde x,\tilde y)|}{6}+1\leq 
  \frac{d(x,\tilde x)+d(y,\tilde y)}{6}+1\leq \frac{d}{6} \\ \text{ 
et } (\Delta_{x,y}+1)\eta\geq \frac{d(x,y)-6\de}{6}\eta \geq \frac{3d}{2} \geq d+3|\Delta_{x,y}-\Delta_{\tilde x,\tilde y}|.\end{gather*}
%  
%$(\Delta_{\tilde x,\tilde y}+1)\eta\geq 4d$ et $(\Delta_{\tilde x,\tilde y}+1)(v_{3}-\eta)\leq (\Delta_{x,y}+1)v_{3}-3d$, et les autres inégalités sont similaires ou plus faciles. 
 On montre aussi les $6$ inégalités analogues obtenues en permutant les rôles de $x,y$ et $\tilde x,\tilde y$. 
On possède donc, pour $(u_{1},u_{2},u_{3},v_{1},v_{2},v_{3})\in [\eta,1-\eta[^{6}$, des applications 
$$ A_{x,u_{1},u_{2},u_{3}}^{y,v_{1},v_{2},v_{3}}\vad{\gamma} 
A_{\tilde  x,u_{1}-\eta,u_{2}+\eta,u_{3}-\eta}^{\tilde  y,v_{1}-\eta,v_{2}+\eta,v_{3}-\eta}\text{\ \  et \ \ 
}A_{\tilde  x,u_{1},u_{2},u_{3}}^{\tilde  y,v_{1},v_{2},v_{3}}\vad{\gamma}  
A_{x,u_{1}-\eta,u_{2}+\eta,u_{3}-\eta}^{y,v_{1}-\eta,v_{2}+\eta,v_{3}-\eta}.$$ 
Les phrases précédentes, le a) du lemme~\ref{tubes} et les parties a), b), c) du sous-lemme~\ref{tubes-xyx'y'}  garantissent que ces applications vérifient les conditions de compatibilité  supposées dans le sous-lemme suivant (en prenant 
$ A_{u_{1},u_{2},u_{3}}^{v_{1},v_{2},v_{3}}=A_{x,u_{1},u_{2},u_{3}}^{y,v_{1},v_{2},v_{3}}$ 
et
$ {\tilde A}_{u_{1},u_{2},u_{3}}^{v_{1},v_{2},v_{3}}=A_{\tilde  x,u_{1},u_{2},u_{3}}^{\tilde  y,v_{1},v_{2},v_{3}}$).

 \begin{souslem}\label{lemme-asssrrr}
Soit $C\in \N^{*}$. Alors pour tout $\eta\in ]0,\frac{1}{4}[$  et pour toutes familles d'ensembles finis non vides,  de cardinaux inférieurs ou égaux à $C$, 
$$(A_{u_{1},u_{2},u_{3}}^{v_{1},v_{2},v_{3}})_{(u_{1},u_{2},u_{3},v_{1},v_{2},v_{3})\in [0,1[^{6} }
\text{\ \  et\ \ } 
({\tilde A}_{u_{1},u_{2},u_{3}}^{v_{1},v_{2},v_{3}})_{(u_{1},u_{2},u_{3},v_{1},v_{2},v_{3})\in [0,1[^{6} }
\text{\ \  munis}$$ 
\begin{itemize}
\item d'applications surjectives $A_{u_{1},u_{2},u_{3}}^{v_{1},v_{2},v_{3}}\to A_{u'_{1},u_{2},u_{3}}^{v'_{1},v_{2},v_{3}}$ et $\tilde A_{u_{1},u_{2},u_{3}}^{v_{1},v_{2},v_{3}}\to \tilde  A_{u'_{1},u_{2},u_{3}}^{v'_{1},v_{2},v_{3}}$ pour $u'_{1}\leq u_{1}$ et $v'_{1}\leq v_{1}$, 
\item d'applications injectives $A_{u_{1},u_{2},u_{3}}^{v_{1},v_{2},v_{3}}\to A_{u_{1},u'_{2},u'_{3}}^{v_{1},v'_{2},v'_{3}}$ et $\tilde  A_{u_{1},u_{2},u_{3}}^{v_{1},v_{2},v_{3}}\to \tilde  A_{u_{1},u'_{2},u'_{3}}^{v_{1},v'_{2},v'_{3}}$
pour $u'_{2}\geq u_{2}$, $v'_{2}\geq v_{2}$, $u'_{3}\leq u_{3}$ et $v'_{3}\leq v_{3}$,
\item et pour $(u_{1},u_{2},u_{3},v_{1},v_{2},v_{3})\in [\eta,1-\eta[^{6}$, 
d'applications $$ A_{u_{1},u_{2},u_{3}}^{v_{1},v_{2},v_{3}}\to 
\tilde A_{u_{1}-\eta,u_{2}+\eta,u_{3}-\eta}^{v_{1}-\eta,v_{2}+\eta,v_{3}-\eta}\text{\ \  et \ \ 
}\tilde A_{u_{1},u_{2},u_{3}}^{v_{1},v_{2},v_{3}}\to 
A_{u_{1}-\eta,u_{2}+\eta,u_{3}-\eta}^{v_{1}-\eta,v_{2}+\eta,v_{3}-\eta}$$
\end{itemize}
 telles que toutes ces applications soient compatibles entre elles
 (c'est-à-dire que deux composées d'applications comme ci-dessus sont égales lorsqu'elles ont le même ensemble de départ et le même ensemble d'arrivée), alors la mesure de l'ensemble des  $(u_{1},u_{2},u_{3},v_{1},v_{2},v_{3})\in [\eta,1-\eta[^{6}$ tels que l'application 
 %$$ A_{u_{1},u_{2},u_{3}}^{v_{1},v_{2},v_{3}}\to 
 %A_{u_{1}-\eta,u_{2}+\eta,u_{3}-\eta}^{v_{1}-\eta,v_{2}+\eta,v_{3}-\eta}\text{\ \  et \ \ 
$$A_{u_{1},u_{2},u_{3}}^{v_{1},v_{2},v_{3}}\to 
\tilde A_{u_{1}-\eta,u_{2}+\eta,u_{3}-\eta}^{v_{1}-\eta,v_{2}+\eta,v_{3}-\eta}$$
  soit bijective  est $\geq 1-50C\eta$.
\end{souslem}
 \noindent{\bf Démonstration.}  %Nous allons montrer que la mesure de l'ensemble des  $$(u_{1},u_{2},u_{3},v_{1},v_{2},v_{3})\in [\eta,1-\eta[^{6}\text{  tels que }$$ $$\tilde A_{u_{1},u_{2},u_{3}}^{v_{1},v_{2},v_{3}}\to 
%A_{u_{1}-\eta,u_{2}+\eta,u_{3}-\eta}^{v_{1}-\eta,v_{2}+\eta,v_{3}-\eta}\text{
%  soit bijective}$$  est $\geq 1-12(2C+4)\eta$ (cela suffit car en rempla\c cant $\tilde A$ par $A$ le même énoncé est évidemment vrai, d'où l'on déduit facilement le lemme). 
Pour 
$(u_{1},u_{2},u_{3},v_{1},v_{2},v_{3})\in [2\eta,1-2\eta[^{6}$ on considère  
les applications 
$$\tilde A_{u_{1}+\eta,u_{2}-\eta,u_{3}+\eta}^{v_{1}+\eta,v_{2}-\eta,v_{3}+\eta}\vad{f } A_{u_{1},u_{2},u_{3}}^{v_{1},v_{2},v_{3}}\vad{g} 
\tilde A_{u_{1}-\eta,u_{2}+\eta,u_{3}-\eta}^{v_{1}-\eta,v_{2}+\eta,v_{3}-\eta}\vad{h} A_{u_{1}-2\eta,u_{2}+2\eta,u_{3}-2\eta}
^{v_{1}-2\eta,v_{2}+2\eta,v_{3}-2\eta}.$$ 
Si $g\circ f $ et $h\circ g$ sont bijectives, $g$ est bijective. Or $g\circ f $ est égale à la composée 
\begin{gather}\nonumber  \tilde A_{u_{1}+\eta,u_{2}-\eta,u_{3}+\eta}^{v_{1}+\eta,v_{2}-\eta,v_{3}+\eta}\to 
 \tilde A_{u_{1}-\eta,u_{2}-\eta,u_{3}+\eta}^{v_{1}+\eta,v_{2}-\eta,v_{3}+\eta}\to
 \tilde A_{u_{1}-\eta,u_{2}+\eta,u_{3}+\eta}^{v_{1}+\eta,v_{2}-\eta,v_{3}+\eta}\to
\tilde  A_{u_{1}-\eta,u_{2}+\eta,u_{3}-\eta}^{v_{1}+\eta,v_{2}-\eta,v_{3}+\eta}\to \\ \label{comp-six1}
 \tilde A_{u_{1}-\eta,u_{2}+\eta,u_{3}-\eta}^{v_{1}-\eta,v_{2}-\eta,v_{3}+\eta}\to
\tilde  A_{u_{1}-\eta,u_{2}+\eta,u_{3}-\eta}^{v_{1}-\eta,v_{2}+\eta,v_{3}+\eta}\to
\tilde  A_{u_{1}-\eta,u_{2}+\eta,u_{3}-\eta}^{v_{1}-\eta,v_{2}+\eta,v_{3}-\eta}\end{gather}
et $h\circ g$ est égale à la composée 
\begin{gather} \label{comp-six2}
 \begin{split} A_{u_{1},u_{2},u_{3}}^{v_{1},v_{2},v_{3}}\to 
 A_{u_{1}-2\eta,u_{2},u_{3}}^{v_{1},v_{2},v_{3}}\to
 A_{u_{1}-2\eta,u_{2}+2\eta,u_{3}}^{v_{1},v_{2},v_{3}}\to
 A_{u_{1}-2\eta,u_{2}+2\eta,u_{3}-2\eta}^{v_{1},v_{2},v_{3}}\to\\  A_{u_{1}-2\eta,u_{2}+2\eta,u_{3}-2\eta}^{v_{1}-2\eta,v_{2},v_{3}}\to
 A_{u_{1}-2\eta,u_{2}+2\eta,u_{3}-2\eta}^{v_{1}-2\eta,v_{2}+2\eta,v_{3}}\to A_{u_{1}-2\eta,u_{2}+2\eta,u_{3}-2\eta}^{v_{1}-2\eta,v_{2}+2\eta,v_{3}-2\eta}.\end{split}\end{gather}
Pour $u_{2},u_{3},v_{1},v_{2},v_{3}\in [2\eta,1-2\eta[$, la mesure de l'ensemble des $u_{1}\in [2\eta,1-2\eta[$ tels que 
\begin{gather}\label{applic+-++-+'}\tilde A_{u_{1}+\eta,u_{2}-\eta,u_{3}+\eta}^{v_{1}+\eta,v_{2}-\eta,v_{3}+\eta}\to 
\tilde A_{u_{1}-\eta,u_{2}-\eta,u_{3}+\eta}^{v_{1}+\eta,v_{2}-\eta,v_{3}+\eta}\end{gather} 
(qui est la première application dans (\ref{comp-six1}))
ne soit pas une bijection est inférieure ou égale à 
$4\eta(C-1)$. 
En effet   l'application  (\ref{applic+-++-+'}) est surjective donc est une bijection en cas d'égalité des cardinaux des deux parties et d'autre part 
 l'application de $[0,1[$ dans $\N$ qui à $u$ associe le cardinal de 
$\tilde A_{u,u_{2}-\eta,u_{3}+\eta}^{v_{1}+\eta,v_{2}-\eta,v_{3}+\eta}$ est décroissante et prend ses valeurs dans $\{1,...,C\}$. 
 On a $11$ autres énoncés correspondant aux autres applications de (\ref{comp-six1}) et (\ref{comp-six2}). Donc la mesure de l'ensemble des  $$(u_{1},u_{2},u_{3},v_{1},v_{2},v_{3})\in [2\eta,1-2\eta[^{6}$$ tels que $g$ soit une bijection est supérieure ou égale à 
 $$(1-4\eta)^{6}-48\eta(C-1)(1-4\eta)^{5}\geq 1-24\eta-48\eta(C-1)\geq 1-48C\eta\geq 1-50C\eta$$
 et le sous-lemme~\ref{lemme-asssrrr} est démontré. 
\cqfd

 %  
 %Les parties a) des lemmes~\ref{tubes} et~\ref{tubes-xyx'y'}  garantissent la compatibilité entre les applications qui est supposée dans le lemme suivant. 

\noindent{\bf Fin de la démonstration du lemme~\ref{mesure-uuuvvv-dbemol}.} 
Soit $C=C(\delta,K)$ tel que pour $x,y\in X$ et 
$r_{1},r_{2},r_{3},s_{1},s_{2},s_{3}$ vérifiant (\ref{rrrsss-existence}) on ait 
$\sharp(\tilde \Lambda_{x,r_{1},r_{2},r_{3}}^{y,s_{1},s_{2},s_{3}})\leq C$. 
Soit $\rho\in \N$. On pose $d=4\rho+6\delta$. Soient $x,y,x',y'\in X$ vérifiant $d(x,x')\leq \rho$ et $d(y,y')\leq \rho$. Il suffit de montrer le lemme~\ref{mesure-uuuvvv-dbemol} en supposant  $d(x,y)> 36d+6\de$. On choisit alors 
$\eta\in ]0,\frac{1}{4}[$ vérifiant $d(x,y)=\frac{9d}{\eta}+6\de$, c'est-à-dire $\eta=\frac{9d}{d(x,y)-6\de}$. 
On applique le sous-lemme~\ref{lemme-asssrrr}  
  avec $(\tilde x,\tilde y)$ égal à $(x,y)$, $(x',y)$, $(x,y')$ ou $(x',y')$ et avec les familles $(A_{u_{1},u_{2},u_{3}}^{v_{1},v_{2},v_{3}})$ et  $(\tilde A_{u_{1},u_{2},u_{3}}^{v_{1},v_{2},v_{3}})$ égales à 
  $(A_{x,u_{1},u_{2},u_{3}}^{y,v_{1},v_{2},v_{3}})$ et 
  $(A_{\tilde x,u_{1},u_{2},u_{3}}^{\tilde y,v_{1},v_{2},v_{3}})$. 
  
Il existe donc une partie $J\subset [\eta,1-\eta[^{6}$ de mesure $\geq 1-200C\eta$ telle que pour 
$(u_{1},u_{2},u_{3},v_{1},v_{2},v_{3})\in J$  les applications 
\begin{gather*} A_{x,u_{1},u_{2},u_{3}}^{y,v_{1},v_{2},v_{3}}\to 
 A_{x,u_{1}-\eta,u_{2}+\eta,u_{3}-\eta}^{y,v_{1}-\eta,v_{2}+\eta,v_{3}-\eta},\ \ \ A_{x,u_{1},u_{2},u_{3}}^{y,v_{1},v_{2},v_{3}}\to 
 A_{x',u_{1}-\eta,u_{2}+\eta,u_{3}-\eta}^{y,v_{1}-\eta,v_{2}+\eta,v_{3}-\eta},\\    A_{x,u_{1},u_{2},u_{3}}^{y,v_{1},v_{2},v_{3}}\to 
 A_{x,u_{1}-\eta,u_{2}+\eta,u_{3}-\eta}^{y',v_{1}-\eta,v_{2}+\eta,v_{3}-\eta}\text{\ \  et \ \ }A_{x,u_{1},u_{2},u_{3}}^{y,v_{1},v_{2},v_{3}}\to 
 A_{x',u_{1}-\eta,u_{2}+\eta,u_{3}-\eta}^{y',v_{1}-\eta,v_{2}+\eta,v_{3}-\eta}\end{gather*} soient bijectives.

Soient $(u_{1},u_{2},u_{3},v_{1},v_{2},v_{3})\in J$, $c\in A_{x,u_{1},u_{2},u_{3}}^{y,v_{1},v_{2},v_{3}}$ et \begin{gather*}c_{x,y}\in A_{x,u_{1}-\eta,u_{2}+\eta,u_{3}-\eta}^{y,v_{1}-\eta,v_{2}+\eta,v_{3}-\eta},\ c_{x',y}\in A_{x',u_{1}-\eta,u_{2}+\eta,u_{3}-\eta}^{y,v_{1}-\eta,v_{2}+\eta,v_{3}-\eta},\\ 
c_{x,y'}\in A_{x,u_{1}-\eta,u_{2}+\eta,u_{3}-\eta}^{y',v_{1}-\eta,v_{2}+\eta,v_{3}-\eta} \text{\ et }
c_{x',y'}\in A_{x',u_{1}-\eta,u_{2}+\eta,u_{3}-\eta}^{y',v_{1}-\eta,v_{2}+\eta,v_{3}-\eta}\end{gather*}
les  images de $c$   par les quatre bijections ci-dessus.
 Alors  $$\beta_{x,u_{1}-\eta}^{y,v_{1}-\eta}(c_{x,y})-\beta_{x',u_{1}-\eta}^{y,v_{1}-\eta}(c_{x',y})-\beta_{x,u_{1}-\eta}^{y',v_{1}-\eta}(c_{x,y'})+\beta_{x',u_{1}-\eta}^{y',v_{1}-\eta}(c_{x',y'})=0$$ car par (\ref{compat-alpha-beta-22oct09}) on a 
$$\beta_{\tilde x ,u_{1}-\eta}^{\tilde y ,v_{1}-\eta}(c_{\tilde x ,\tilde y })=\min 
_{z\in Y_{x,y}^{E((\Delta_{x,y}+1)u_{1})}}\big(d(\tilde x ,z)+c_{z}\big)+\min 
_{t\in Y_{y,x}^{E((\Delta_{x,y}+1)v_{1})}}\big(c_{t}+d(t,\tilde y )\big)$$ pour $(\tilde x,\tilde y)$ égal à $(x,y)$, $(x',y)$, $(x,y')$ ou $(x',y')$.
Pour tout $(u_{1},u_{2},u_{3},v_{1},v_{2},v_{3})\in J$ on a donc 
\begin{gather*}{d^{\flat}}_{u_{1}-\eta,u_{2}+\eta,u_{3}-\eta}^{v_{1}-\eta,v_{2}+\eta,v_{3}-\eta}(x,y) 
 -{d^{\flat}}_{u_{1}-\eta,u_{2}+\eta,u_{3}-\eta}^{v_{1}-\eta,v_{2}+\eta,v_{3}-\eta}(x',y)\\  -{d^{\flat}}_{u_{1}-\eta,u_{2}+\eta,u_{3}-\eta}^{v_{1}-\eta,v_{2}+\eta,v_{3}-\eta}(x,y')
 + {d^{\flat}}_{u_{1}-\eta,u_{2}+\eta,u_{3}-\eta}^{v_{1}-\eta,v_{2}+\eta,v_{3}-\eta}(x',y')= 0.\end{gather*}
 Comme la mesure de $J$ est supérieure ou égale à $$ 1-200C\eta=1-\frac{1800Cd}{d(x,y)-6\de}=1-\frac{1800C(4\rho+6\delta)}{d(x,y)-6\delta},$$ cela termine la démonstration du lemme~\ref{mesure-uuuvvv-dbemol}.  \cqfd

\begin{prop}\label{prop-d'}
Pour tout $x,y\in X$ on a $d(x,y)\leq d^{\flat}(x,y)\leq d(x,y)+7\delta$ et 
$d^{\flat}$ vérifie la condition (\ref{condition-d'}).
\end{prop}
Il résulte de la première assertion que l'on pourrait facilement remplacer $d^{\flat}$ par une vraie distance (comme on avait obtenu $d''$ à partir de $d'$ dans la preuve de la proposition~\ref{myu}), mais cela n'est pas nécessaire pour la suite. 

 \noindent{\bf Démonstration.}
 La première assertion résulte du lemme~\ref{123-7de}. La seconde assertion résulte des lemmes~\ref{123-7de} et~\ref{mesure-uuuvvv-dbemol}.
 \cqfd 
 
 \label{def-rhobemol}
On définit alors, pour $x\in X$,   la fonction $\rho^{\flat}_{x}:X\to \R_{+}$ 
par $$\rho^{\flat}_{x}(a)=d^{\flat}(x,a).$$ On l'étend ensuite en une fonction 
$\rho^{\flat}_{x}:\Delta\to \R_{+}$ par la formule 
$$\rho^{\flat}_{x}(S)=\frac{\sum _{a\in S}\rho^{\flat}(a)}{\sharp S}\text{ si }S\text{ est non vide et }\rho^{\flat}_{x}(\emptyset)=0.$$ 

\section{Construction des normes}\label{construction-normes}

Soit $s\in ]0,1]$. 
Nous allons donner la formule pour la norme 
pré-hilbertienne 
$\|.\|_{\H_{x,s}(\Delta_{p})}$ sur $\C^{(\Delta_{p})}$ qui servira pour la partie difficile de l'homotopie de $1$ à $\gamma$. 
Nous montrerons d'abord que cette norme est bien définie.
Nous montrerons ensuite la continuité des opérateurs, puis les propriétés d'équivariance de cette norme, et enfin l'équivariance  à compacts près des opérateurs.

\subsection{Formule pour la norme}\label{sous-para-formule-norme}

On rappelle que la constante $F$ définie dans (\ref{def-F}) est majorée par une constante de la forme $C(\de,K,N)$.  
On fixe un entier $P\in \N^{*}$ tel que 
$$(H_{P}) : P \text{ soit divisible par  }3 \text{ et assez grand en fonction de }\de,K,N,Q.$$
Dans la suite nous utiliserons un nombre {\it fini} de fois l'inégalité $P\geq C$ avec $C$ de la forme $C(\de,K,N,Q)$. 

Soient $p\in \{1,...,p_{\max}\}$ et $k,m,l_{0},...,l_{m}\in \N$. 
\begin{defi}\label{defi-Y}
On note $Y_{x}^{p,k,m,(l_{0},...,l_{m})}$ l'ensemble des $(p+m+1+\sum_{i=0}^{m}l_{i})$-uplets $$(a_{1},\dots,a_{p},S_{0},...,S_{m},(\mathcal Y_{i}^{j})_{i\in \{0,\dots,m\}, j\in \{1,\dots ,l_{i}\}})$$ tels que 
 \begin{itemize}
\item i) $a_{1},\dots,a_{p}\in X$ sont deux à deux distincts, $S_{0}=
\{a_{1},\dots ,a_{p}\} $, $S_{0}$  appartient à $\Delta_{p}$, $S_{i}\in \Delta\setminus \{\emptyset\}$  pour $i\in \{1,\dots,m\}$, et
 pour tout $i\in \{0,\dots, m-1\}$, on a
\begin{gather*}S_{i+1}\subset S_{i}
%\cup B(x,3\de+N+5\de p_{\max}) $$ $$
\cup 
\bigcup _{\tilde x\in B(x,k),a\in S_{i}} \{y\in 4\de\tg(\tilde x,a), d(y,a)\in ]N-2\de, QN]\}\\ 
\cup 
\bigcup _{\tilde x\in B(x,k),a\in S_{i}}
\{z\in F\tg(\tilde x,a), d(z,a) \geq \frac{Q}{F}
\},\end{gather*}
\item ii) pour tout $i\in  \{1,\dots,m\}$,   $d(x,S_{i})>k+P$, 
\item iii) pour $i\in \{0,\dots,m-1\}$ et $j\in \{1,\dots ,l_{i}\}$, $\mathcal Y_{i}^{j}$ est une partie non vide de $X$ de diamètre inférieur ou égal à $P$  et    $$\mathcal Y_{i}^{j}\subset 
\bigcup _{y\in S_{i},z\in S_{i+1}}P\tg(y,z).$$
\item iv) pour tout $j\in \{1,...,l_{m}\}$,  $\mathcal Y_{m}^{j}$ est une partie non vide de $X$ de diamètre inférieur ou égal à $P$ et $$\mathcal Y_{m}^{j}\subset 
\bigcup _{y\in S_{m}, \tilde x\in B(x,k)}2P\tg(\tilde x,y)\text{\ \ 
et\ \ } d(x,\mathcal Y_{m}^{j})\geq k+3P.$$
\end{itemize}
\end{defi}
\noindent {\bf Remarque.} La condition iv) implique que  $l_{m}=0$ si $ d(x,S_{m})\leq k$ (ce qui ne peut se produire que si $m=0$ à cause de ii)). En effet pour $y\in B(x,k+N)$ et $\tilde x\in B(x,k)$, on a $2P\tg(\tilde x,y)\subset B(x,k+N+2P+\de)$ par le lemme~\ref{boules-geod} et on suppose $P> N+\de$ (ce qui est permis par $(H_{P})$).

Dans la définition précédente, 
la seule raison pour laquelle on veut connaître $(a_{1},\pp,a_{p})$ en plus  de $S_{0}$ est que cela détermine $e_{a_{1}}\wedge \pp \wedge e_{a_{p}}$ alors que $S_{0}$ permet seulement de le connaître au signe près (le choix de $e_{S_{0}}=\pm 
e_{a_{1}}\wedge \pp \wedge e_{a_{p}}$, qui avait pour but de simplifier certaines formules, ne doit pas être utilisé ici, bien sûr). En contrepartie, chaque partie 
$S_{0}$ est comptée $p!$ fois, mais ce n'est pas grave car $p$ est borné par $p_{\max}$, qui est de la forme $C(\de,K,N)$. 
Dans la définition précédente la donnée de $S_{0}$ est redondante mais nous préférons la garder car elle fournit une notation commode pour 
$\{a_{1},\dots,a_{p}\}$.

On fixe  un entier $M\in \N^{*}$ tel que 
$$(H_{M}) : M \text{ soit pair et  soit assez grand en fonction de }\de,K,N,Q,P.$$
Dans la suite nous utiliserons un nombre {\it fini} de fois l'inégalité $M\geq C$ avec $C$ de la forme $C(\de,K,N,Q,P)$.

On introduit maintenant une partition de $Y_{x}^{p,k,m,(l_{0},...,l_{m})}$ pour la relation d'équivalence suivante : $$(a_{1},\dots,a_{p},S_{0},...,S_{m},(\mathcal Y_{i}^{j})_{i\in \{0,\dots,m\}, j\in \{1,\dots ,l_{i}\}})$$ et $$(\hat a_{1},\dots,\hat a_{p},\hat S_{0},...,\hat S_{m},(\hat {\mathcal Y}_{i}^{j})_{i\in \{0,\dots,m\}, j\in \{1,\dots ,l_{i}\}})$$
sont en relation 
s'il existe une isométrie de $$
\bigcup _{ i\in \{0,\dots ,m\}}
B(S_{i}, M)
\cup  \bigcup _{i\in \{0,\dots,m\}, j\in \{1,\dots ,l_{i}\}} B(\mathcal Y_{i}^{j},M)\cup B(x,k+2M)$$ vers 
$$\bigcup _{ i\in \{0,\dots ,m\}}
B(\hat S_{i}, M)
\cup  \bigcup _{i\in \{0,\dots,m\}, j\in \{1,\dots ,l_{i}\}} B(\hat {\mathcal Y}_{i}^{j},M)\cup B(x,k+2M)$$ 
 qui envoie 
$a_{i}$ sur $\hat a_{i}$ pour $i\in \{1,\dots,p\}$,  
 $S_{i}$ sur $\hat S_{i}$
 pour $i\in \{0,\dots ,m\}$, $\mathcal Y_{i}^{j}$ sur $\hat {\mathcal Y}_{i}^{j}$  pour 
 $i\in \{0,\dots,m\}, j\in \{1,\dots ,l_{i}\}$  
 et est l'identité sur $B(x,k+2M)$.
 On rappelle que la notation $B(A,r)$ pour $A$ une partie de $X$ a été introduite dans (\ref{def-B-30dec1647}). 
% \noindent {\bf Remarque.} La raison pour laquelle on a pris $B(x,k+2M)$ au lieu de $B(x,k+M)$ est que  $$
%\bigcup _{ i\in \{0,\dots ,m\}}
%\{y\in X,  d(y,S_{i})\leq M\}
%\cup  \bigcup _{i\in \{0,\dots,m\}, j\in \{1,\dots ,l_{i}\}} \{y\in X,  d(y,\mathcal Y_{i}^{j})\leq M\} $$  est inclus dans $B(x,k+2M)$ lorsque les parties  $S_{i}$ et  $\mathcal Y_{i}^{j}$ sont incluses dans $B(x,k+M)$. 

\label{overlinemathcalY}
On note $\overline Y_{x}^{p,k,m,(l_{0},...,l_{m})}$
le quotient de $Y_{x}^{p,k,m,(l_{0},...,l_{m})}$ pour cette relation d'équivalence, et $\pi_{x}^{p,k,m,(l_{0},...,l_{m})}$ l'application quotient. 

\noindent{\bf Notations.} 
Pour $Z\in \overline Y_{x}^{p,k,m,(l_{0},...,l_{m})}$  on note $r_{0}(Z),\dots,r_{m}(Z),s_{0}(Z),\dots ,  s_{m}(Z)$ les entiers tels que 
\begin{itemize}
\item $r_{i}(Z)=d(x,S_{i})$ pour  $i\in \{0,\pp,m\}$, 
 \item $s_{i}(Z)=d(S_{i},S_{i+1})+2M$ pour  $i\in \{0,\pp,m-1\}$, 
\item $s_{m}(Z)=d(x,S_{m})-k$ 
\end{itemize}
%et on note $r_{i}(Z)\in \N$ l'unique entier tel que  $r_{i}(Z)=d(x,S_{i})$ 
 pour tout $(a_{1},\dots,a_{p},S_{0},...,S_{m},(\mathcal Y_{i}^{j})_{i\in \{0,\dots,m\}, j\in \{1,\dots ,l_{i}\}})\in (\pi_{x}^{p,k,m,(l_{0},...,l_{m})})^{-1}(Z). $
  
 %De même pour $Z\in \overline Y_{x}^{p,k,m,(l_{0},...,l_{m})}$, $i\in \{0,\pp,m\}$ et $j\in \{1,\pp,l_{i}\}$ on note  $s_{i}^{j}(Z)\in \{0,\pp,d(x,S_{i})-k+N\}$ l'unique entier tel que, pour tout $$(a_{1},\dots,a_{p},S_{1},...,S_{m},(\mathcal Y_{i}^{j})_{i\in \{0,\dots,m\}, j\in \{1,\dots ,l_{i}\}})$$
 %dans 
 %$(\pi_{x}^{p,k,m,(l_{0},...,l_{m})})^{-1}(Z)$, on ait  $$\mathcal Y_{i}^{j}\subset \{y\in 
%\bigcup_{\tilde x\in B(x,k),a\in S_{i}} \geod(\tilde x,a),d(y,S_{i})=s_{i}^{j}(Z)\}$$ (l'unicité de $s_{i}^{j}(Z)$ est garantie par la condition que $\mathcal Y_{i}^{j}$ est non vide). 

On fixe $B\in \R_{+}^{*}$  et $\alpha\in ]0,1[$ tels que 
$$(H_{B}) : B \text{   soit assez grand en fonction de }\de,K,N,Q,P,M,s\text{\ \ \ \ et}$$
$$(H_{\alpha}) : \alpha \text{  soit assez petit en fonction de }\de,K,N,Q,P,M,s,B.$$

Pour $Z\in \overline Y_{x}^{p,k,m,(l_{0},...,l_{m})}$ on note $\xi_{Z}$ la forme linéaire sur $\C^{(\Delta_{p})}$  définie par 
\begin{gather}\label{defxi-31dec1554}
\xi_{Z}(f)=\sum _{(a_{1},\dots,a_{p},S_{0},...,S_{m},(\mathcal Y_{i}^{j})_{i\in \{0,\dots,m\}, j\in \{1,\dots ,l_{i}\}}) \in 
(\pi_{x}^{p,k,m,(l_{0},...,l_{m})})^{-1}(Z)} f(a_1,...,a_p). 
\end{gather}
Dans cette formule, comme dans la suite, pour $f\in \C^{(\Delta_{p})}$ on note $f(a_1,...,a_p)$ le coefficient de $e_{a_1}\wedge ...\wedge
e_{a_p}$ lorsque qu'on écrit $f$ dans la base $(\pm e_{S})_{S\in \Delta_{p}}$ (dont les vecteurs sont définis au signe près). On écrira aussi $f(S)=\pm 
f(a_1,...,a_p)$ si $S=\{a_1,...,a_p\}$.

On munit alors $\C^{(\Delta_{p})}$ de la norme pré-hilbertienne, définie par la formule suivante : 
\begin{gather}\nonumber \|f\|_{\H_{x,s}(\Delta_{p})}^{2}
=\sum _{k,m,l_{0},\dots ,l_{m}\in \N}  
%$$\|\sum_{\{a_1,a_2,...,a_p\}\in \Delta_p}
%f(a_1,...,a_p)e_{a_1}\wedge ...\wedge
%e_{a_p}\|_{\H_{x,
%}(\Delta_{p})}^{2}
%=\sum _{k,m,l_{0},\dots ,l_{m}\in \N} 
%$$
B^{-(m+\sum_{i=0}^{m}l_{i})}
\sum_{Z\in \overline Y_{x}^{p,k,m,(l_{0},...,l_{m})}} \\ \label{formule-norme}  e^{2s(r_{0}(Z)-k)}\Big(\prod_{i=0}^{m}s_{i}(Z)^{-l_{i}} \Big)
\sharp \big((\pi_{x}^{p,k,m,(l_{0},...,l_{m})})^{-1}(Z)\big)^{-\alpha }
\big|\xi_{Z}(f)\big|^{2}.\end{gather}

\noindent {\bf Remarque.} On a vu que $l_{m}=0$ lorsque $s_{m}(Z) \leq 0$ (ce qui ne peut se produire que si $m=0$) et dans ce cas on convient que $s_{m}(Z)^{l_{m}} =1$. En revanche pour $i\in \{0,\pp,m-1\}$, on a toujours 
$s_{i}(Z)\geq 1$. 

\noindent {\bf Remarque. } On verra dans la démonstration des
 propositions~\ref{continuite-del}    et~\ref{continuite-Jx} que l'on majore la partie de $\|\del(f)\|_{\H_{x,s}(\Delta_{p})}^{2}$ ou
  de $\|J_{x}(f)\|_{\H_{x,s}(\Delta_{p})}^{2}$ correspondant à une valeur donnée de $m$ par la partie de $\|f\|_{\H_{x,s}(\Delta_{p})}^{2}$ correspondant à $m+1$  
  (en fait c'est un peu plus compliqué mais  on renvoie aux 
 propositions~\ref{continuite-del} et~\ref{continuite-Jx} pour les détails). D'autre part la partie de $\|f\|_{\H_{x,s}(\Delta_{p})}^{2}$ correspondant à $m=0$ et $l_{0}=0$ est essentiellement égale au carré de  la norme introduite dans le premier paragraphe de~\cite{duke} pour montrer que les groupes hyperboliques n'ont pas la propriété (T) renforcée. Enfin les $\mathcal Y_{i}^{j}$ assurent la connaissance des points intermédiaires sur lesquels on moyenne dans la construction de $h_{x},u_{x},\rho'_{x}$ (le but de ces moyennes est d'assurer l'équivariance à compacts près des opérateurs par l'astuce des ensembles emboîtés, comme dans~\cite{ks}).  %Les $\mathcal Y_{i}^{j}$ sont des parties de $X$ de diamètre $\leq 2N+2\de$ par le lemme~\ref{abxyz-trapeze}. Comme $M$ est grand par rapport à $\de,K,N,Q$, on a l'impression qu'on pourrait remplacer ces parties par des singletons. Cependant, dans la preuve du lemme~\ref{continuite-Jx}, on construit certaines parties $\mathcal Y_{i}^{j}$ comme $\{y\in 
%\bigcup_{\tilde x\in B(x,k),a\in S_{i}} \geod(\tilde x,a),d(y,S_{i})=s\}$ pour une certaine valeur de $s$ (égale  à $s_{i}^{j}(Z)$ par définition de $s_{i}^{j}(Z)$) et   la nécessité de choisir un  point dans cette partie  empêcherait la preuve de fonctionner. En même temps on ne peut pas imposer  $$\mathcal Y_{i}^{j}=\{y\in 
%\bigcup_{\tilde x\in B(x,k),a\in S_{i}} \geod(\tilde x,a),d(y,S_{i})=s_{i}^{j}(Z)\}$$
%car cela empêcherait de montrer la proposition~\ref{equiv-normeYZ} qui donne la propriété d'équivariance de la norme
%(en effet on verra comment, pour $x'\in X$ et sous certaines conditions, on peut choisir $k'$ tel que 
%$$\{y\in 
%\bigcup_{\tilde x\in B(x,k),a\in S_{i}} \geod(\tilde x,a),d(y,S_{i})=s_{i}^{j}(Z)\}$$ $$\subset \{y\in 
%\bigcup_{\tilde x'\in B(x',k'),a\in S_{i}} \geod(\tilde x',a),d(y,S_{i})=s_{i}^{j}(Z)\}$$ mais on ne peut certainement pas avoir l'égalité). 
On ne peut pas demander que les parties $\mathcal Y_{i}^{j}$ soient des singletons car ces parties seront construites de fa\c con naturelle   dans les preuves et la nécessité de choisir un point dans chacune empêcherait les preuves de fonctionner.

\noindent {\bf Remarque.} 
Le facteur $e^{2s(r_{0}(Z)-k)}$ est utile pour la continuité de $J_{x}$, comme on le verra dans la proposition~\ref{continuite-Jx}. Le facteur 
$\prod_{i=0}^{m}s_{i}(Z)^{-l_{i}} $ est motivé par le fait qu'étant donnés $a_{1},\pp,a_{p},S_{0},\pp,S_{m}$ le nombre de possibilités pour $(\mathcal Y_{i}^{j})_{i\in \{0,\dots,m\}, j\in \{1,\dots ,l_{i}\}}$ vérifiant  
$$(a_{1},\dots,a_{p},S_{0},...,S_{m},(\mathcal Y_{i}^{j})_{i\in \{0,\dots,m\}, j\in \{1,\dots ,l_{i}\}})\in Y_{x}^{p,k,m,(l_{0},...,l_{m})} $$
est borné par $C^{\sum_{i=0}^{m}l_{i}} \prod_{i=0}^{m}s_{i}(Z)^{l_{i}}$, avec $C=C(\de,K,N,Q,P)$, comme on le verra dans le lemme~\ref{existence-C}. 
La constante $2M$ qui apparaît dans la définition de $s_{i}(Z)$ pour $i\in \{0,...,m-1\}$ trouvera son utilité 
dans la démonstration du lemme~\ref{somme-aZZ'-ter-2}. 
Le  facteur $\sharp ((\pi_{x}^{p,k,m,(l_{0},...,l_{m})})^{-1}(Z))^{-\alpha }$ servira  pour montrer que les opérateurs $\del$ et $J_{x}$ conjugués par $e^{\tau  \theta^{\flat}_{x}}$  sont équivariants à compact près. Enfin le facteur 
 $B^{-(m+\sum_{i=0}^{m}l_{i})}$ (avec $B$ assez grand) est nécessaire pour que la somme converge, comme le montre la preuve de la proposition~\ref{norme-bien-definie}.

 \subsection{Enoncé des résultats} \label{sous-para-enonce-resultats}
 
 La proposition suivante montre que la norme pré-hilbertienne définie dans le sous-paragraphe précédent a bien un sens. 
 
 \begin{prop}\label{norme-bien-definie}
 Pour tout $p\in \{1,\dots  ,p_{\max}\}$ et pour tout $f\in \C^{(\Delta_{p})}$, 
 $\|f\|_{\H_{x,s}(\Delta_{p})}$ est fini. 
 \end{prop}
 
 Cette proposition sera démontrée dans le sous-paragraphe~\ref{premieres-proprietes} et la démonstration utilisera l'hypothèse (qui est une partie de $(H_{B})$) selon laquelle  $B$ est assez grand en fonction de $\de,K,N,Q,P,M$ (sans cette hypothèse la proposition serait fausse). 
 
 On note 
$\H_{x,s}(\Delta_{p})$ le complété de $\C^{(\Delta_{p})}$  pour la norme pré-hilbertienne   
$\|.\|_{\H_{x,s}(\Delta_{p})}$. 

On note $\H_{x,s}=\bigoplus_{p=1}^{p_{\max}}\H_{x,s}(\Delta_{p})$. 

 La proposition suivante sera démontrée dans le sous-paragraphe~\ref{action-G-normes}. 
 
\begin{prop}\label{equiv-normeYZ}
L'action de $G$ sur $\bigoplus_{p=1}^{p_{\max}}\C^{(\Delta_{p})}$
s'étend en une action continue de $G$ sur $\H_{x,s}$ et il existe une constante $C$ telle que pour tout $g\in G$ on ait 
$\|\pi(g)\|_{\L(\H_{x,s})}\leq e^{2s\ell(g)+C}$. 
\end{prop}
La démonstration de cette proposition utilisera entièrement l'hypothèse $(H_{B})$ (y compris la condition  que $B$   est assez grand en fonction de $s$).

La proposition suivante est le résultat principal de ce paragraphe. 

Pour  $p\in \{1,\pp,p_{\max}\}$, on définit $\theta^{\flat}_{x}:\C^{(\Delta_{p})}\to \C^{(\Delta_{p})}$ par $\theta^{\flat}_{x}(e_{S})=\rho^{\flat}_{x}(S)e_{S}$ pour tout $S\in \Delta_{p}$. Pour tout $\tau \in \R$ on note $e^{\tau  \theta^{\flat}_{x}}: \C^{(\Delta_{p})}\to \C^{(\Delta_{p})}$ l'opérateur défini par $e^{\tau  \theta^{\flat}_{x}}(e_{S})=e^{\tau \rho^{\flat}_{x}(S)}e_{S}$. 

\begin{prop}\label{enonce-ppal-KKC01}
Pour tout $T\in \R_{+}$, l'opérateur $$(e^{\tau \theta^{\flat}_{x}}(\del + J_{x}\del J_{x}) e^{-\tau \theta^{\flat}_{x}})_{\tau \in [0,T]}$$ s'étend  en un opérateur 
continu sur le  $\C[0,T]$-module hilbertien 
$\H_{x,s}[0,T]$, et 
$$(\H_{x,s}[0,T],(e^{\tau \theta^{\flat}_{x}}(\del + J_{x}\del J_{x}) e^{-\tau \theta^{\flat}_{x}})_{\tau \in [0,T]})$$ appartient à $KK_{G,2s\ell+C}(\C,\C[0,T])$. 
\end{prop}

En particulier $(\H_{x,s},\del +J_{x}\del J_{x})$ est un élément de $KK_{G,2s\ell+C}(\C,\C)$.

\begin{prop}\label{homot-1-142}
Cet élément est égal à l'image de $1\in KK_{G}(\C,\C)$. 
\end{prop}
 \noindent{\bf Démonstration.}
La preuve est analogue à celle de la proposition  1.4.2 de~\cite{kkban}.  On reprend  les arguments car on travaille ici dans un cadre hilbertien plutôt que banachique. Soit $\hat \H_{x,s}=\C\oplus \H_{x,s}$ (avec une graduation inversée pour $\H_{x,s}$). On veut montrer que 
$(\hat  \H_{x,s},\del +J_{x}\del J_{x})$ est nul dans $KK_{G,2s\ell+C}(\C,\C)$. Il résulte de l'égalité $\del J_x  + J_x  \del=\mathrm{Id}$ que $\del J_x  $ et $ J_x  \del$ commutent
et $ J_x  \del\del J_x  =0$ car $\del^2=0$. Donc $\del J_x  $ et $ J_x  \del$ sont deux projecteurs qui
commutent, de somme 
$\mathrm{Id}$, de produit nul, pairs. On pose $T= J_x  \del J_x  $. Alors  $T^2=0$. 
On a $\del T=\del  J_x  $ et $T\del = J_x  \del $. 
D'autre part, $\del $ et $\del T$
ont même image, de même que $T$ et $T\del $, d'où 
l'identification $\hat\H_{x,s}=\del (\hat\H_{x,s})\oplus T(\hat\H_{x,s})$ en tant qu'espaces de Hilbert {\it à équivalence des normes près} (car les deux idempotents $\del J_x  $ et $ J_x  \del$ ne sont pas nécessairement auto-adjoints). Dans cette décomposition, l'action de  $g\in G$ s'écrit sous la forme 
$\begin{pmatrix} c_{1,1}(g) & c_{1,2}(g)\\ 0 & c_{2,2}(g)\end{pmatrix}$. 
Considérons alors le $\C[0,1]$-module hilbertien 
$\Z/2\Z$-gradué $\hat \H_{x,s}[0,1]$, où l'action de  $g\in G$ est donnée par 
  $\Big(\begin{pmatrix} c_{1,1}(g) & tc_{1,2}(g)\\ 0 &
c_{2,2}(g)\end{pmatrix}\Big)_{t\in [0,1]}$. 
Alors $(\hat\H_{x,s}[0,1], \del +T)$ fournit une homotopie entre $(\hat\H_{x,s},\del +T)$, en $t=1$ 
et un élément dégénéré, en $t=0$. En effet, quand l'action de
$G$ est diagonale, les opérateurs $\del $ et $T$ commutent exactement à
cette action. \cqfd

Donc  $(\H_{x,s}[0,T],(e^{\tau \theta^{\flat}_{x}}(\del + J_{x}\del J_{x}) e^{-\tau \theta^{\flat}_{x}})_{\tau \in [0,T]})$ réalise  une homotopie
entre $1$ et $(\H_{x,s},e^{T\theta^{\flat}_{x}}(\del +J_{x}\del J_{x}) e^{-T\theta^{\flat}_{x}})$ et montre  donc l'égalité entre ces deux éléments dans $KK_{G,2s\ell+C}(\C,\C)$. On fixera $T$ assez grand et cela constituera  la partie difficile de l'homotopie de $1$ à $\gamma$. La partie facile (qui fera l'objet du paragraphe~\ref{construction-fin}), sera une homotopie entre $(\H_{x,s},e^{T\theta^{\flat}_{x}}(\del +J_{x}\del J_{x}) e^{-T\theta^{\flat}_{x}})$ et  $\gamma$ (qui montrera donc l'égalité entre  ces deux éléments  dans $KK_{G,2s\ell+C}(\C,\C)$). 

\subsection{Premières propriétés de la norme}\label{premieres-proprietes}

Le but de ce sous-paragraphe est de montrer la proposition~\ref{norme-bien-definie}.  
Les lemmes~\ref{approx-arbres},  \ref{approx-arbres2}, \ref{x-a,b-y-geod}, \ref{Bxk-x-t-z-y},  \ref{Bxkmunu},  \ref{Bxkmunu-cor}, \ref{m-fini} et  \ref{existence-C}
  sont  des préliminaires 
à la preuve de la proposition~\ref{norme-bien-definie}.

Nous commen\c cons par rappeler le lemme d'approximation par les arbres. 

\begin{lem}\label{approx-arbres}
Soit $(Y,d_{Y})$ un espace métrique fini et $\delta$-hyperbolique, et $w\in Y$ un point base. Soit $l\in \N$ tel que $\sharp Y\leq 2^{l}+2
$. Alors il existe un arbre métrique   fini $(T,d_{T})$ et une application $\Psi:Y\to T$ telle que 
\begin{itemize}
\item pour  $y\in Y$,  $d_{T}(\Psi w, \Psi y)=d_{Y}(w,y)$, 
\item pour $y,z\in Y$, $d_{Y}(y,z)-l\de\leq d_{T}(\Psi y, \Psi z) \leq d_{Y}(y,z)$.
\end{itemize}
\end{lem}
{\bf Démonstration.} C'est exactement le (i) du théorème 12 du chapitre 2  de~\cite{harpehyper} car un espace métrique est $\de$-hyperbolique au sens de la définition~\ref{defi-hyperb-Hxyzt}  si et seulement s'il est $\frac{\de}{2}$-hyperbolique au sens de la définition 3 (reformulée dans 4) du chapitre 2  de~\cite{harpehyper}. \cqfd

\begin{lem}\label{approx-arbres2} Dans les notations du lemme précédent, soient $y,z,t\in Y$.

\noindent a) Si  $t\in \alpha\tg(y,z)$, alors $\Psi t\in (\alpha+l\de)\tg(\Psi y,\Psi z)$  et si $y=w$, $\Psi t\in \alpha \tg(\Psi y,\Psi z)$. 

b) Si  $\Psi t\in \alpha\tg(\Psi y,\Psi z)$, alors $t\in (\alpha+2l\de)\tg(y,z)$ et si $t=w$, $t\in \alpha\tg(y,z)$. 

c) On a $|d(t,\geod(y,z))-d(\Psi t,\geod(\Psi y,\Psi z))|\leq (l+1)\de+1$. 
\end{lem}
{\bf Démonstration.} Seul c) demande une démonstration. On a $$
d(\Psi t,\geod(\Psi y,\Psi z))=\frac{d(\Psi t,\Psi y)+d(\Psi t,\Psi z)-d(\Psi y,\Psi z)}{2}.$$  Il est évident que $\frac{d( t, y)+d( t, z)-d( y, z)}{2} \leq d( t,\geod( y, z))$. Enfin 
$$d( t,\geod( y, z))\leq \frac{d( t, y)+d( t, z)-d( y, z)}{2}+\de+1$$ car si  $v\in \geod( y, z)$ est tel que $d(y,v)=E(\frac{d( t, y)+d( y, z)-d( t, z)}{2})$ on a $d(t,v)\leq  \frac{d( t, y)+d( t, z)-d( y, z)}{2}+\de+1$ par $(H_{\de}^{0}(t,y,v,z))$. Le c) en résulte facilement. \cqfd

\begin{lem}\label{x-a,b-y-geod}
Soient $\alpha,\beta\in \N$, $\rho\in\Z$  et $x,y,a,b,\in X$ tels que 
$$a\in \alpha\tg(x,y), b\in \beta\tg(x,y)\text{\  et\ } d(x,b)\leq d(x,a)+\rho. $$ Alors 
$b\in (\max(\alpha+2\rho,\beta)+\de)\tg(x,a)$. 
\end{lem}

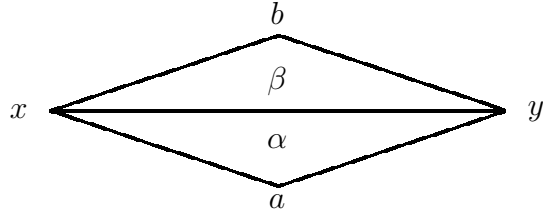
\ifx\JPicScale\undefined\def\JPicScale{1}\fi
\unitlength \JPicScale mm
\begin{picture}(90,30)(20,13)
\linethickness{0.3mm}
\put(50,30){\line(1,0){60}}
\linethickness{0.3mm}
\multiput(80,40)(0.36,-0.12){83}{\line(1,0){0.36}}
\linethickness{0.3mm}
\multiput(50,30)(0.36,0.12){84}{\line(1,0){0.36}}
\linethickness{0.3mm}
\multiput(50,30)(0.36,-0.12){84}{\line(1,0){0.36}}
\linethickness{0.3mm}
\multiput(80,20)(0.36,0.12){83}{\line(1,0){0.36}}
\put(46,30){\makebox(0,0)[cc]{$x$}}

\put(114,30){\makebox(0,0)[cc]{$y$}}

\put(80,43){\makebox(0,0)[cc]{$b$}}

\put(80,34){\makebox(0,0)[cc]{$\beta$}}

\put(80,26){\makebox(0,0)[cc]{$\alpha$}}

\put(80,18){\makebox(0,0)[cc]{$a$}}

\end{picture}
 
 \noindent{\bf Démonstration.}
Par $(H_{\de}(a,x,b,y))$ on a  
$$d(a,b)\leq \max(d(a,y)+d(x,b)-d(x,y),d(a,x)+d(b,y)-d(x,y))+\de$$ donc 
\begin{gather*}d(a,b)+d(b,x)-d(a,x)\leq \max\big(d(a,y)+2d(x,b)-d(a,x)-d(x,y),\\ d(b,x)+d(b,y)-d(x,y)\big)+\de
\leq \max\big(\alpha+2\rho,\beta\big)+\de.\end{gather*}\cqfd

\begin{lem}\label{Bxk-x-t-z-y}
Soient $k\in \N$ et $x,y,z,t$ des points de $X$  tels que $z$ appartienne à $B(x,k)$ et que $t$ soit un point de $B(x,k)$ à distance minimale de $ y$. Alors $t\in \de\tg(z,y)$. 
\end{lem}

 \noindent{\bf Démonstration.}
 L'énoncé est clair si $y\in B(x,k)$ car alors $t=y$. 
 
 \ifx\JPicScale\undefined\def\JPicScale{1}\fi
\unitlength \JPicScale mm
\begin{picture}(80,50)(0,15)
\linethickness{0.3mm}
\put(50,39.75){\line(0,1){0.5}}
\multiput(49.99,40.75)(0.01,-0.5){1}{\line(0,-1){0.5}}
\multiput(49.96,41.25)(0.02,-0.5){1}{\line(0,-1){0.5}}
\multiput(49.92,41.74)(0.04,-0.5){1}{\line(0,-1){0.5}}
\multiput(49.87,42.24)(0.05,-0.5){1}{\line(0,-1){0.5}}
\multiput(49.81,42.73)(0.06,-0.49){1}{\line(0,-1){0.49}}
\multiput(49.74,43.23)(0.07,-0.49){1}{\line(0,-1){0.49}}
\multiput(49.65,43.72)(0.09,-0.49){1}{\line(0,-1){0.49}}
\multiput(49.55,44.21)(0.1,-0.49){1}{\line(0,-1){0.49}}
\multiput(49.44,44.69)(0.11,-0.49){1}{\line(0,-1){0.49}}
\multiput(49.32,45.18)(0.12,-0.48){1}{\line(0,-1){0.48}}
\multiput(49.18,45.66)(0.14,-0.48){1}{\line(0,-1){0.48}}
\multiput(49.04,46.13)(0.15,-0.48){1}{\line(0,-1){0.48}}
\multiput(48.88,46.61)(0.16,-0.47){1}{\line(0,-1){0.47}}
\multiput(48.71,47.07)(0.17,-0.47){1}{\line(0,-1){0.47}}
\multiput(48.52,47.54)(0.09,-0.23){2}{\line(0,-1){0.23}}
\multiput(48.33,48)(0.1,-0.23){2}{\line(0,-1){0.23}}
\multiput(48.13,48.45)(0.1,-0.23){2}{\line(0,-1){0.23}}
\multiput(47.91,48.9)(0.11,-0.22){2}{\line(0,-1){0.22}}
\multiput(47.68,49.35)(0.11,-0.22){2}{\line(0,-1){0.22}}
\multiput(47.44,49.78)(0.12,-0.22){2}{\line(0,-1){0.22}}
\multiput(47.19,50.22)(0.12,-0.22){2}{\line(0,-1){0.22}}
\multiput(46.93,50.64)(0.13,-0.21){2}{\line(0,-1){0.21}}
\multiput(46.66,51.06)(0.14,-0.21){2}{\line(0,-1){0.21}}
\multiput(46.38,51.47)(0.14,-0.21){2}{\line(0,-1){0.21}}
\multiput(46.09,51.88)(0.15,-0.2){2}{\line(0,-1){0.2}}
\multiput(45.79,52.27)(0.1,-0.13){3}{\line(0,-1){0.13}}
\multiput(45.48,52.66)(0.1,-0.13){3}{\line(0,-1){0.13}}
\multiput(45.16,53.05)(0.11,-0.13){3}{\line(0,-1){0.13}}
\multiput(44.83,53.42)(0.11,-0.12){3}{\line(0,-1){0.12}}
\multiput(44.49,53.79)(0.11,-0.12){3}{\line(0,-1){0.12}}
\multiput(44.14,54.14)(0.12,-0.12){3}{\line(0,-1){0.12}}
\multiput(43.79,54.49)(0.12,-0.12){3}{\line(1,0){0.12}}
\multiput(43.42,54.83)(0.12,-0.11){3}{\line(1,0){0.12}}
\multiput(43.05,55.16)(0.12,-0.11){3}{\line(1,0){0.12}}
\multiput(42.66,55.48)(0.13,-0.11){3}{\line(1,0){0.13}}
\multiput(42.27,55.79)(0.13,-0.1){3}{\line(1,0){0.13}}
\multiput(41.88,56.09)(0.13,-0.1){3}{\line(1,0){0.13}}
\multiput(41.47,56.38)(0.2,-0.15){2}{\line(1,0){0.2}}
\multiput(41.06,56.66)(0.21,-0.14){2}{\line(1,0){0.21}}
\multiput(40.64,56.93)(0.21,-0.14){2}{\line(1,0){0.21}}
\multiput(40.22,57.19)(0.21,-0.13){2}{\line(1,0){0.21}}
\multiput(39.78,57.44)(0.22,-0.12){2}{\line(1,0){0.22}}
\multiput(39.35,57.68)(0.22,-0.12){2}{\line(1,0){0.22}}
\multiput(38.9,57.91)(0.22,-0.11){2}{\line(1,0){0.22}}
\multiput(38.45,58.13)(0.22,-0.11){2}{\line(1,0){0.22}}
\multiput(38,58.33)(0.23,-0.1){2}{\line(1,0){0.23}}
\multiput(37.54,58.52)(0.23,-0.1){2}{\line(1,0){0.23}}
\multiput(37.07,58.71)(0.23,-0.09){2}{\line(1,0){0.23}}
\multiput(36.61,58.88)(0.47,-0.17){1}{\line(1,0){0.47}}
\multiput(36.13,59.04)(0.47,-0.16){1}{\line(1,0){0.47}}
\multiput(35.66,59.18)(0.48,-0.15){1}{\line(1,0){0.48}}
\multiput(35.18,59.32)(0.48,-0.14){1}{\line(1,0){0.48}}
\multiput(34.69,59.44)(0.48,-0.12){1}{\line(1,0){0.48}}
\multiput(34.21,59.55)(0.49,-0.11){1}{\line(1,0){0.49}}
\multiput(33.72,59.65)(0.49,-0.1){1}{\line(1,0){0.49}}
\multiput(33.23,59.74)(0.49,-0.09){1}{\line(1,0){0.49}}
\multiput(32.73,59.81)(0.49,-0.07){1}{\line(1,0){0.49}}
\multiput(32.24,59.87)(0.49,-0.06){1}{\line(1,0){0.49}}
\multiput(31.74,59.92)(0.5,-0.05){1}{\line(1,0){0.5}}
\multiput(31.25,59.96)(0.5,-0.04){1}{\line(1,0){0.5}}
\multiput(30.75,59.99)(0.5,-0.02){1}{\line(1,0){0.5}}
\multiput(30.25,60)(0.5,-0.01){1}{\line(1,0){0.5}}
\put(29.75,60){\line(1,0){0.5}}
\multiput(29.25,59.99)(0.5,0.01){1}{\line(1,0){0.5}}
\multiput(28.75,59.96)(0.5,0.02){1}{\line(1,0){0.5}}
\multiput(28.26,59.92)(0.5,0.04){1}{\line(1,0){0.5}}
\multiput(27.76,59.87)(0.5,0.05){1}{\line(1,0){0.5}}
\multiput(27.27,59.81)(0.49,0.06){1}{\line(1,0){0.49}}
\multiput(26.77,59.74)(0.49,0.07){1}{\line(1,0){0.49}}
\multiput(26.28,59.65)(0.49,0.09){1}{\line(1,0){0.49}}
\multiput(25.79,59.55)(0.49,0.1){1}{\line(1,0){0.49}}
\multiput(25.31,59.44)(0.49,0.11){1}{\line(1,0){0.49}}
\multiput(24.82,59.32)(0.48,0.12){1}{\line(1,0){0.48}}
\multiput(24.34,59.18)(0.48,0.14){1}{\line(1,0){0.48}}
\multiput(23.87,59.04)(0.48,0.15){1}{\line(1,0){0.48}}
\multiput(23.39,58.88)(0.47,0.16){1}{\line(1,0){0.47}}
\multiput(22.93,58.71)(0.47,0.17){1}{\line(1,0){0.47}}
\multiput(22.46,58.52)(0.23,0.09){2}{\line(1,0){0.23}}
\multiput(22,58.33)(0.23,0.1){2}{\line(1,0){0.23}}
\multiput(21.55,58.13)(0.23,0.1){2}{\line(1,0){0.23}}
\multiput(21.1,57.91)(0.22,0.11){2}{\line(1,0){0.22}}
\multiput(20.65,57.68)(0.22,0.11){2}{\line(1,0){0.22}}
\multiput(20.22,57.44)(0.22,0.12){2}{\line(1,0){0.22}}
\multiput(19.78,57.19)(0.22,0.12){2}{\line(1,0){0.22}}
\multiput(19.36,56.93)(0.21,0.13){2}{\line(1,0){0.21}}
\multiput(18.94,56.66)(0.21,0.14){2}{\line(1,0){0.21}}
\multiput(18.53,56.38)(0.21,0.14){2}{\line(1,0){0.21}}
\multiput(18.12,56.09)(0.2,0.15){2}{\line(1,0){0.2}}
\multiput(17.73,55.79)(0.13,0.1){3}{\line(1,0){0.13}}
\multiput(17.34,55.48)(0.13,0.1){3}{\line(1,0){0.13}}
\multiput(16.95,55.16)(0.13,0.11){3}{\line(1,0){0.13}}
\multiput(16.58,54.83)(0.12,0.11){3}{\line(1,0){0.12}}
\multiput(16.21,54.49)(0.12,0.11){3}{\line(1,0){0.12}}
\multiput(15.86,54.14)(0.12,0.12){3}{\line(1,0){0.12}}
\multiput(15.51,53.79)(0.12,0.12){3}{\line(0,1){0.12}}
\multiput(15.17,53.42)(0.11,0.12){3}{\line(0,1){0.12}}
\multiput(14.84,53.05)(0.11,0.12){3}{\line(0,1){0.12}}
\multiput(14.52,52.66)(0.11,0.13){3}{\line(0,1){0.13}}
\multiput(14.21,52.27)(0.1,0.13){3}{\line(0,1){0.13}}
\multiput(13.91,51.88)(0.1,0.13){3}{\line(0,1){0.13}}
\multiput(13.62,51.47)(0.15,0.2){2}{\line(0,1){0.2}}
\multiput(13.34,51.06)(0.14,0.21){2}{\line(0,1){0.21}}
\multiput(13.07,50.64)(0.14,0.21){2}{\line(0,1){0.21}}
\multiput(12.81,50.22)(0.13,0.21){2}{\line(0,1){0.21}}
\multiput(12.56,49.78)(0.12,0.22){2}{\line(0,1){0.22}}
\multiput(12.32,49.35)(0.12,0.22){2}{\line(0,1){0.22}}
\multiput(12.09,48.9)(0.11,0.22){2}{\line(0,1){0.22}}
\multiput(11.87,48.45)(0.11,0.22){2}{\line(0,1){0.22}}
\multiput(11.67,48)(0.1,0.23){2}{\line(0,1){0.23}}
\multiput(11.48,47.54)(0.1,0.23){2}{\line(0,1){0.23}}
\multiput(11.29,47.07)(0.09,0.23){2}{\line(0,1){0.23}}
\multiput(11.12,46.61)(0.17,0.47){1}{\line(0,1){0.47}}
\multiput(10.96,46.13)(0.16,0.47){1}{\line(0,1){0.47}}
\multiput(10.82,45.66)(0.15,0.48){1}{\line(0,1){0.48}}
\multiput(10.68,45.18)(0.14,0.48){1}{\line(0,1){0.48}}
\multiput(10.56,44.69)(0.12,0.48){1}{\line(0,1){0.48}}
\multiput(10.45,44.21)(0.11,0.49){1}{\line(0,1){0.49}}
\multiput(10.35,43.72)(0.1,0.49){1}{\line(0,1){0.49}}
\multiput(10.26,43.23)(0.09,0.49){1}{\line(0,1){0.49}}
\multiput(10.19,42.73)(0.07,0.49){1}{\line(0,1){0.49}}
\multiput(10.13,42.24)(0.06,0.49){1}{\line(0,1){0.49}}
\multiput(10.08,41.74)(0.05,0.5){1}{\line(0,1){0.5}}
\multiput(10.04,41.25)(0.04,0.5){1}{\line(0,1){0.5}}
\multiput(10.01,40.75)(0.02,0.5){1}{\line(0,1){0.5}}
\multiput(10,40.25)(0.01,0.5){1}{\line(0,1){0.5}}
\put(10,39.75){\line(0,1){0.5}}
\multiput(10,39.75)(0.01,-0.5){1}{\line(0,-1){0.5}}
\multiput(10.01,39.25)(0.02,-0.5){1}{\line(0,-1){0.5}}
\multiput(10.04,38.75)(0.04,-0.5){1}{\line(0,-1){0.5}}
\multiput(10.08,38.26)(0.05,-0.5){1}{\line(0,-1){0.5}}
\multiput(10.13,37.76)(0.06,-0.49){1}{\line(0,-1){0.49}}
\multiput(10.19,37.27)(0.07,-0.49){1}{\line(0,-1){0.49}}
\multiput(10.26,36.77)(0.09,-0.49){1}{\line(0,-1){0.49}}
\multiput(10.35,36.28)(0.1,-0.49){1}{\line(0,-1){0.49}}
\multiput(10.45,35.79)(0.11,-0.49){1}{\line(0,-1){0.49}}
\multiput(10.56,35.31)(0.12,-0.48){1}{\line(0,-1){0.48}}
\multiput(10.68,34.82)(0.14,-0.48){1}{\line(0,-1){0.48}}
\multiput(10.82,34.34)(0.15,-0.48){1}{\line(0,-1){0.48}}
\multiput(10.96,33.87)(0.16,-0.47){1}{\line(0,-1){0.47}}
\multiput(11.12,33.39)(0.17,-0.47){1}{\line(0,-1){0.47}}
\multiput(11.29,32.93)(0.09,-0.23){2}{\line(0,-1){0.23}}
\multiput(11.48,32.46)(0.1,-0.23){2}{\line(0,-1){0.23}}
\multiput(11.67,32)(0.1,-0.23){2}{\line(0,-1){0.23}}
\multiput(11.87,31.55)(0.11,-0.22){2}{\line(0,-1){0.22}}
\multiput(12.09,31.1)(0.11,-0.22){2}{\line(0,-1){0.22}}
\multiput(12.32,30.65)(0.12,-0.22){2}{\line(0,-1){0.22}}
\multiput(12.56,30.22)(0.12,-0.22){2}{\line(0,-1){0.22}}
\multiput(12.81,29.78)(0.13,-0.21){2}{\line(0,-1){0.21}}
\multiput(13.07,29.36)(0.14,-0.21){2}{\line(0,-1){0.21}}
\multiput(13.34,28.94)(0.14,-0.21){2}{\line(0,-1){0.21}}
\multiput(13.62,28.53)(0.15,-0.2){2}{\line(0,-1){0.2}}
\multiput(13.91,28.12)(0.1,-0.13){3}{\line(0,-1){0.13}}
\multiput(14.21,27.73)(0.1,-0.13){3}{\line(0,-1){0.13}}
\multiput(14.52,27.34)(0.11,-0.13){3}{\line(0,-1){0.13}}
\multiput(14.84,26.95)(0.11,-0.12){3}{\line(0,-1){0.12}}
\multiput(15.17,26.58)(0.11,-0.12){3}{\line(0,-1){0.12}}
\multiput(15.51,26.21)(0.12,-0.12){3}{\line(0,-1){0.12}}
\multiput(15.86,25.86)(0.12,-0.12){3}{\line(1,0){0.12}}
\multiput(16.21,25.51)(0.12,-0.11){3}{\line(1,0){0.12}}
\multiput(16.58,25.17)(0.12,-0.11){3}{\line(1,0){0.12}}
\multiput(16.95,24.84)(0.13,-0.11){3}{\line(1,0){0.13}}
\multiput(17.34,24.52)(0.13,-0.1){3}{\line(1,0){0.13}}
\multiput(17.73,24.21)(0.13,-0.1){3}{\line(1,0){0.13}}
\multiput(18.12,23.91)(0.2,-0.15){2}{\line(1,0){0.2}}
\multiput(18.53,23.62)(0.21,-0.14){2}{\line(1,0){0.21}}
\multiput(18.94,23.34)(0.21,-0.14){2}{\line(1,0){0.21}}
\multiput(19.36,23.07)(0.21,-0.13){2}{\line(1,0){0.21}}
\multiput(19.78,22.81)(0.22,-0.12){2}{\line(1,0){0.22}}
\multiput(20.22,22.56)(0.22,-0.12){2}{\line(1,0){0.22}}
\multiput(20.65,22.32)(0.22,-0.11){2}{\line(1,0){0.22}}
\multiput(21.1,22.09)(0.22,-0.11){2}{\line(1,0){0.22}}
\multiput(21.55,21.87)(0.23,-0.1){2}{\line(1,0){0.23}}
\multiput(22,21.67)(0.23,-0.1){2}{\line(1,0){0.23}}
\multiput(22.46,21.48)(0.23,-0.09){2}{\line(1,0){0.23}}
\multiput(22.93,21.29)(0.47,-0.17){1}{\line(1,0){0.47}}
\multiput(23.39,21.12)(0.47,-0.16){1}{\line(1,0){0.47}}
\multiput(23.87,20.96)(0.48,-0.15){1}{\line(1,0){0.48}}
\multiput(24.34,20.82)(0.48,-0.14){1}{\line(1,0){0.48}}
\multiput(24.82,20.68)(0.48,-0.12){1}{\line(1,0){0.48}}
\multiput(25.31,20.56)(0.49,-0.11){1}{\line(1,0){0.49}}
\multiput(25.79,20.45)(0.49,-0.1){1}{\line(1,0){0.49}}
\multiput(26.28,20.35)(0.49,-0.09){1}{\line(1,0){0.49}}
\multiput(26.77,20.26)(0.49,-0.07){1}{\line(1,0){0.49}}
\multiput(27.27,20.19)(0.49,-0.06){1}{\line(1,0){0.49}}
\multiput(27.76,20.13)(0.5,-0.05){1}{\line(1,0){0.5}}
\multiput(28.26,20.08)(0.5,-0.04){1}{\line(1,0){0.5}}
\multiput(28.75,20.04)(0.5,-0.02){1}{\line(1,0){0.5}}
\multiput(29.25,20.01)(0.5,-0.01){1}{\line(1,0){0.5}}
\put(29.75,20){\line(1,0){0.5}}
\multiput(30.25,20)(0.5,0.01){1}{\line(1,0){0.5}}
\multiput(30.75,20.01)(0.5,0.02){1}{\line(1,0){0.5}}
\multiput(31.25,20.04)(0.5,0.04){1}{\line(1,0){0.5}}
\multiput(31.74,20.08)(0.5,0.05){1}{\line(1,0){0.5}}
\multiput(32.24,20.13)(0.49,0.06){1}{\line(1,0){0.49}}
\multiput(32.73,20.19)(0.49,0.07){1}{\line(1,0){0.49}}
\multiput(33.23,20.26)(0.49,0.09){1}{\line(1,0){0.49}}
\multiput(33.72,20.35)(0.49,0.1){1}{\line(1,0){0.49}}
\multiput(34.21,20.45)(0.49,0.11){1}{\line(1,0){0.49}}
\multiput(34.69,20.56)(0.48,0.12){1}{\line(1,0){0.48}}
\multiput(35.18,20.68)(0.48,0.14){1}{\line(1,0){0.48}}
\multiput(35.66,20.82)(0.48,0.15){1}{\line(1,0){0.48}}
\multiput(36.13,20.96)(0.47,0.16){1}{\line(1,0){0.47}}
\multiput(36.61,21.12)(0.47,0.17){1}{\line(1,0){0.47}}
\multiput(37.07,21.29)(0.23,0.09){2}{\line(1,0){0.23}}
\multiput(37.54,21.48)(0.23,0.1){2}{\line(1,0){0.23}}
\multiput(38,21.67)(0.23,0.1){2}{\line(1,0){0.23}}
\multiput(38.45,21.87)(0.22,0.11){2}{\line(1,0){0.22}}
\multiput(38.9,22.09)(0.22,0.11){2}{\line(1,0){0.22}}
\multiput(39.35,22.32)(0.22,0.12){2}{\line(1,0){0.22}}
\multiput(39.78,22.56)(0.22,0.12){2}{\line(1,0){0.22}}
\multiput(40.22,22.81)(0.21,0.13){2}{\line(1,0){0.21}}
\multiput(40.64,23.07)(0.21,0.14){2}{\line(1,0){0.21}}
\multiput(41.06,23.34)(0.21,0.14){2}{\line(1,0){0.21}}
\multiput(41.47,23.62)(0.2,0.15){2}{\line(1,0){0.2}}
\multiput(41.88,23.91)(0.13,0.1){3}{\line(1,0){0.13}}
\multiput(42.27,24.21)(0.13,0.1){3}{\line(1,0){0.13}}
\multiput(42.66,24.52)(0.13,0.11){3}{\line(1,0){0.13}}
\multiput(43.05,24.84)(0.12,0.11){3}{\line(1,0){0.12}}
\multiput(43.42,25.17)(0.12,0.11){3}{\line(1,0){0.12}}
\multiput(43.79,25.51)(0.12,0.12){3}{\line(1,0){0.12}}
\multiput(44.14,25.86)(0.12,0.12){3}{\line(0,1){0.12}}
\multiput(44.49,26.21)(0.11,0.12){3}{\line(0,1){0.12}}
\multiput(44.83,26.58)(0.11,0.12){3}{\line(0,1){0.12}}
\multiput(45.16,26.95)(0.11,0.13){3}{\line(0,1){0.13}}
\multiput(45.48,27.34)(0.1,0.13){3}{\line(0,1){0.13}}
\multiput(45.79,27.73)(0.1,0.13){3}{\line(0,1){0.13}}
\multiput(46.09,28.12)(0.15,0.2){2}{\line(0,1){0.2}}
\multiput(46.38,28.53)(0.14,0.21){2}{\line(0,1){0.21}}
\multiput(46.66,28.94)(0.14,0.21){2}{\line(0,1){0.21}}
\multiput(46.93,29.36)(0.13,0.21){2}{\line(0,1){0.21}}
\multiput(47.19,29.78)(0.12,0.22){2}{\line(0,1){0.22}}
\multiput(47.44,30.22)(0.12,0.22){2}{\line(0,1){0.22}}
\multiput(47.68,30.65)(0.11,0.22){2}{\line(0,1){0.22}}
\multiput(47.91,31.1)(0.11,0.22){2}{\line(0,1){0.22}}
\multiput(48.13,31.55)(0.1,0.23){2}{\line(0,1){0.23}}
\multiput(48.33,32)(0.1,0.23){2}{\line(0,1){0.23}}
\multiput(48.52,32.46)(0.09,0.23){2}{\line(0,1){0.23}}
\multiput(48.71,32.93)(0.17,0.47){1}{\line(0,1){0.47}}
\multiput(48.88,33.39)(0.16,0.47){1}{\line(0,1){0.47}}
\multiput(49.04,33.87)(0.15,0.48){1}{\line(0,1){0.48}}
\multiput(49.18,34.34)(0.14,0.48){1}{\line(0,1){0.48}}
\multiput(49.32,34.82)(0.12,0.48){1}{\line(0,1){0.48}}
\multiput(49.44,35.31)(0.11,0.49){1}{\line(0,1){0.49}}
\multiput(49.55,35.79)(0.1,0.49){1}{\line(0,1){0.49}}
\multiput(49.65,36.28)(0.09,0.49){1}{\line(0,1){0.49}}
\multiput(49.74,36.77)(0.07,0.49){1}{\line(0,1){0.49}}
\multiput(49.81,37.27)(0.06,0.49){1}{\line(0,1){0.49}}
\multiput(49.87,37.76)(0.05,0.5){1}{\line(0,1){0.5}}
\multiput(49.92,38.26)(0.04,0.5){1}{\line(0,1){0.5}}
\multiput(49.96,38.75)(0.02,0.5){1}{\line(0,1){0.5}}
\multiput(49.99,39.25)(0.01,0.5){1}{\line(0,1){0.5}}

\linethickness{0.3mm}
\put(30,40){\line(1,0){50}}
\linethickness{0.3mm}
\multiput(40,50)(0.48,-0.12){83}{\line(1,0){0.48}}
\linethickness{0.3mm}
\multiput(30,40)(0.12,0.12){83}{\line(1,0){0.12}}
\put(27,40){\makebox(0,0)[cc]{$x$}}

\put(37,53){\makebox(0,0)[cc]{$z$}}

\put(53,36){\makebox(0,0)[cc]{$t$}}

\put(80,35){\makebox(0,0)[cc]{$y$}}

\end{picture}

\noindent  Sinon on a $t\in \geod(x,y)$, $d(x,z)\leq k=d(x,t)$ et $d(y,z)\geq d(y,t)$ d'où par $(H_{\de}^{0}(z,x,t,y))$, $d(z,t)\leq \max(d(z,x)-d(x,t),d(z,y)-d(t,y))+\de=d(z,y)-d(t,y)+\de$ d'où $t\in \de\tg(z,y)$. 
 \cqfd

\begin{lem}\label{Bxkmunu} 
Soient $k,\mu,\nu\in \N$,  $x, y,y'\in X$ 
et $z\in B(x,k)$
vérifiant \begin{itemize}
\item $y'\not\in B(x,k+\frac{\mu+\de}{2})$, 
\item  $y'\in \mu\text{-}\geod(z,y)$ et $d(z,y')\leq d(z,y)-\nu$.
\end{itemize}
Alors $$d(x,y')\leq d(x,y)-\nu+\de\text{\ \  et \ \ }y'\in (\mu+\de)\tg(x,y).$$  
\end{lem}

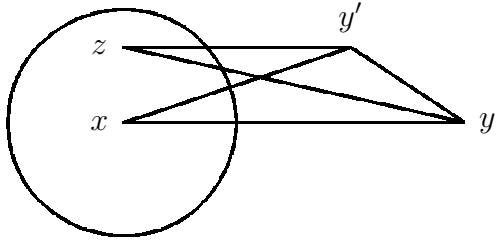
\ifx\JPicScale\undefined\def\JPicScale{1}\fi
\unitlength \JPicScale mm
\begin{picture}(95,40)(20,30)
\linethickness{0.3mm}
\put(60,49.75){\line(0,1){0.5}}
\multiput(59.98,50.75)(0.02,-0.5){1}{\line(0,-1){0.5}}
\multiput(59.95,51.25)(0.03,-0.5){1}{\line(0,-1){0.5}}
\multiput(59.9,51.75)(0.05,-0.5){1}{\line(0,-1){0.5}}
\multiput(59.83,52.25)(0.07,-0.5){1}{\line(0,-1){0.5}}
\multiput(59.75,52.74)(0.08,-0.49){1}{\line(0,-1){0.49}}
\multiput(59.65,53.23)(0.1,-0.49){1}{\line(0,-1){0.49}}
\multiput(59.53,53.72)(0.12,-0.49){1}{\line(0,-1){0.49}}
\multiput(59.4,54.2)(0.13,-0.48){1}{\line(0,-1){0.48}}
\multiput(59.25,54.68)(0.15,-0.48){1}{\line(0,-1){0.48}}
\multiput(59.09,55.16)(0.16,-0.47){1}{\line(0,-1){0.47}}
\multiput(58.91,55.62)(0.09,-0.23){2}{\line(0,-1){0.23}}
\multiput(58.71,56.09)(0.1,-0.23){2}{\line(0,-1){0.23}}
\multiput(58.5,56.54)(0.11,-0.23){2}{\line(0,-1){0.23}}
\multiput(58.27,56.99)(0.11,-0.22){2}{\line(0,-1){0.22}}
\multiput(58.03,57.43)(0.12,-0.22){2}{\line(0,-1){0.22}}
\multiput(57.78,57.86)(0.13,-0.22){2}{\line(0,-1){0.22}}
\multiput(57.51,58.28)(0.13,-0.21){2}{\line(0,-1){0.21}}
\multiput(57.22,58.69)(0.14,-0.21){2}{\line(0,-1){0.21}}
\multiput(56.93,59.1)(0.15,-0.2){2}{\line(0,-1){0.2}}
\multiput(56.62,59.49)(0.1,-0.13){3}{\line(0,-1){0.13}}
\multiput(56.29,59.87)(0.11,-0.13){3}{\line(0,-1){0.13}}
\multiput(55.96,60.25)(0.11,-0.12){3}{\line(0,-1){0.12}}
\multiput(55.61,60.61)(0.12,-0.12){3}{\line(0,-1){0.12}}
\multiput(55.25,60.96)(0.12,-0.12){3}{\line(1,0){0.12}}
\multiput(54.87,61.29)(0.12,-0.11){3}{\line(1,0){0.12}}
\multiput(54.49,61.62)(0.13,-0.11){3}{\line(1,0){0.13}}
\multiput(54.1,61.93)(0.13,-0.1){3}{\line(1,0){0.13}}
\multiput(53.69,62.22)(0.2,-0.15){2}{\line(1,0){0.2}}
\multiput(53.28,62.51)(0.21,-0.14){2}{\line(1,0){0.21}}
\multiput(52.86,62.78)(0.21,-0.13){2}{\line(1,0){0.21}}
\multiput(52.43,63.03)(0.22,-0.13){2}{\line(1,0){0.22}}
\multiput(51.99,63.27)(0.22,-0.12){2}{\line(1,0){0.22}}
\multiput(51.54,63.5)(0.22,-0.11){2}{\line(1,0){0.22}}
\multiput(51.09,63.71)(0.23,-0.11){2}{\line(1,0){0.23}}
\multiput(50.62,63.91)(0.23,-0.1){2}{\line(1,0){0.23}}
\multiput(50.16,64.09)(0.23,-0.09){2}{\line(1,0){0.23}}
\multiput(49.68,64.25)(0.47,-0.16){1}{\line(1,0){0.47}}
\multiput(49.2,64.4)(0.48,-0.15){1}{\line(1,0){0.48}}
\multiput(48.72,64.53)(0.48,-0.13){1}{\line(1,0){0.48}}
\multiput(48.23,64.65)(0.49,-0.12){1}{\line(1,0){0.49}}
\multiput(47.74,64.75)(0.49,-0.1){1}{\line(1,0){0.49}}
\multiput(47.25,64.83)(0.49,-0.08){1}{\line(1,0){0.49}}
\multiput(46.75,64.9)(0.5,-0.07){1}{\line(1,0){0.5}}
\multiput(46.25,64.95)(0.5,-0.05){1}{\line(1,0){0.5}}
\multiput(45.75,64.98)(0.5,-0.03){1}{\line(1,0){0.5}}
\multiput(45.25,65)(0.5,-0.02){1}{\line(1,0){0.5}}
\put(44.75,65){\line(1,0){0.5}}
\multiput(44.25,64.98)(0.5,0.02){1}{\line(1,0){0.5}}
\multiput(43.75,64.95)(0.5,0.03){1}{\line(1,0){0.5}}
\multiput(43.25,64.9)(0.5,0.05){1}{\line(1,0){0.5}}
\multiput(42.75,64.83)(0.5,0.07){1}{\line(1,0){0.5}}
\multiput(42.26,64.75)(0.49,0.08){1}{\line(1,0){0.49}}
\multiput(41.77,64.65)(0.49,0.1){1}{\line(1,0){0.49}}
\multiput(41.28,64.53)(0.49,0.12){1}{\line(1,0){0.49}}
\multiput(40.8,64.4)(0.48,0.13){1}{\line(1,0){0.48}}
\multiput(40.32,64.25)(0.48,0.15){1}{\line(1,0){0.48}}
\multiput(39.84,64.09)(0.47,0.16){1}{\line(1,0){0.47}}
\multiput(39.38,63.91)(0.23,0.09){2}{\line(1,0){0.23}}
\multiput(38.91,63.71)(0.23,0.1){2}{\line(1,0){0.23}}
\multiput(38.46,63.5)(0.23,0.11){2}{\line(1,0){0.23}}
\multiput(38.01,63.27)(0.22,0.11){2}{\line(1,0){0.22}}
\multiput(37.57,63.03)(0.22,0.12){2}{\line(1,0){0.22}}
\multiput(37.14,62.78)(0.22,0.13){2}{\line(1,0){0.22}}
\multiput(36.72,62.51)(0.21,0.13){2}{\line(1,0){0.21}}
\multiput(36.31,62.22)(0.21,0.14){2}{\line(1,0){0.21}}
\multiput(35.9,61.93)(0.2,0.15){2}{\line(1,0){0.2}}
\multiput(35.51,61.62)(0.13,0.1){3}{\line(1,0){0.13}}
\multiput(35.13,61.29)(0.13,0.11){3}{\line(1,0){0.13}}
\multiput(34.75,60.96)(0.12,0.11){3}{\line(1,0){0.12}}
\multiput(34.39,60.61)(0.12,0.12){3}{\line(1,0){0.12}}
\multiput(34.04,60.25)(0.12,0.12){3}{\line(0,1){0.12}}
\multiput(33.71,59.87)(0.11,0.12){3}{\line(0,1){0.12}}
\multiput(33.38,59.49)(0.11,0.13){3}{\line(0,1){0.13}}
\multiput(33.07,59.1)(0.1,0.13){3}{\line(0,1){0.13}}
\multiput(32.78,58.69)(0.15,0.2){2}{\line(0,1){0.2}}
\multiput(32.49,58.28)(0.14,0.21){2}{\line(0,1){0.21}}
\multiput(32.22,57.86)(0.13,0.21){2}{\line(0,1){0.21}}
\multiput(31.97,57.43)(0.13,0.22){2}{\line(0,1){0.22}}
\multiput(31.73,56.99)(0.12,0.22){2}{\line(0,1){0.22}}
\multiput(31.5,56.54)(0.11,0.22){2}{\line(0,1){0.22}}
\multiput(31.29,56.09)(0.11,0.23){2}{\line(0,1){0.23}}
\multiput(31.09,55.62)(0.1,0.23){2}{\line(0,1){0.23}}
\multiput(30.91,55.16)(0.09,0.23){2}{\line(0,1){0.23}}
\multiput(30.75,54.68)(0.16,0.47){1}{\line(0,1){0.47}}
\multiput(30.6,54.2)(0.15,0.48){1}{\line(0,1){0.48}}
\multiput(30.47,53.72)(0.13,0.48){1}{\line(0,1){0.48}}
\multiput(30.35,53.23)(0.12,0.49){1}{\line(0,1){0.49}}
\multiput(30.25,52.74)(0.1,0.49){1}{\line(0,1){0.49}}
\multiput(30.17,52.25)(0.08,0.49){1}{\line(0,1){0.49}}
\multiput(30.1,51.75)(0.07,0.5){1}{\line(0,1){0.5}}
\multiput(30.05,51.25)(0.05,0.5){1}{\line(0,1){0.5}}
\multiput(30.02,50.75)(0.03,0.5){1}{\line(0,1){0.5}}
\multiput(30,50.25)(0.02,0.5){1}{\line(0,1){0.5}}
\put(30,49.75){\line(0,1){0.5}}
\multiput(30,49.75)(0.02,-0.5){1}{\line(0,-1){0.5}}
\multiput(30.02,49.25)(0.03,-0.5){1}{\line(0,-1){0.5}}
\multiput(30.05,48.75)(0.05,-0.5){1}{\line(0,-1){0.5}}
\multiput(30.1,48.25)(0.07,-0.5){1}{\line(0,-1){0.5}}
\multiput(30.17,47.75)(0.08,-0.49){1}{\line(0,-1){0.49}}
\multiput(30.25,47.26)(0.1,-0.49){1}{\line(0,-1){0.49}}
\multiput(30.35,46.77)(0.12,-0.49){1}{\line(0,-1){0.49}}
\multiput(30.47,46.28)(0.13,-0.48){1}{\line(0,-1){0.48}}
\multiput(30.6,45.8)(0.15,-0.48){1}{\line(0,-1){0.48}}
\multiput(30.75,45.32)(0.16,-0.47){1}{\line(0,-1){0.47}}
\multiput(30.91,44.84)(0.09,-0.23){2}{\line(0,-1){0.23}}
\multiput(31.09,44.38)(0.1,-0.23){2}{\line(0,-1){0.23}}
\multiput(31.29,43.91)(0.11,-0.23){2}{\line(0,-1){0.23}}
\multiput(31.5,43.46)(0.11,-0.22){2}{\line(0,-1){0.22}}
\multiput(31.73,43.01)(0.12,-0.22){2}{\line(0,-1){0.22}}
\multiput(31.97,42.57)(0.13,-0.22){2}{\line(0,-1){0.22}}
\multiput(32.22,42.14)(0.13,-0.21){2}{\line(0,-1){0.21}}
\multiput(32.49,41.72)(0.14,-0.21){2}{\line(0,-1){0.21}}
\multiput(32.78,41.31)(0.15,-0.2){2}{\line(0,-1){0.2}}
\multiput(33.07,40.9)(0.1,-0.13){3}{\line(0,-1){0.13}}
\multiput(33.38,40.51)(0.11,-0.13){3}{\line(0,-1){0.13}}
\multiput(33.71,40.13)(0.11,-0.12){3}{\line(0,-1){0.12}}
\multiput(34.04,39.75)(0.12,-0.12){3}{\line(0,-1){0.12}}
\multiput(34.39,39.39)(0.12,-0.12){3}{\line(1,0){0.12}}
\multiput(34.75,39.04)(0.12,-0.11){3}{\line(1,0){0.12}}
\multiput(35.13,38.71)(0.13,-0.11){3}{\line(1,0){0.13}}
\multiput(35.51,38.38)(0.13,-0.1){3}{\line(1,0){0.13}}
\multiput(35.9,38.07)(0.2,-0.15){2}{\line(1,0){0.2}}
\multiput(36.31,37.78)(0.21,-0.14){2}{\line(1,0){0.21}}
\multiput(36.72,37.49)(0.21,-0.13){2}{\line(1,0){0.21}}
\multiput(37.14,37.22)(0.22,-0.13){2}{\line(1,0){0.22}}
\multiput(37.57,36.97)(0.22,-0.12){2}{\line(1,0){0.22}}
\multiput(38.01,36.73)(0.22,-0.11){2}{\line(1,0){0.22}}
\multiput(38.46,36.5)(0.23,-0.11){2}{\line(1,0){0.23}}
\multiput(38.91,36.29)(0.23,-0.1){2}{\line(1,0){0.23}}
\multiput(39.38,36.09)(0.23,-0.09){2}{\line(1,0){0.23}}
\multiput(39.84,35.91)(0.47,-0.16){1}{\line(1,0){0.47}}
\multiput(40.32,35.75)(0.48,-0.15){1}{\line(1,0){0.48}}
\multiput(40.8,35.6)(0.48,-0.13){1}{\line(1,0){0.48}}
\multiput(41.28,35.47)(0.49,-0.12){1}{\line(1,0){0.49}}
\multiput(41.77,35.35)(0.49,-0.1){1}{\line(1,0){0.49}}
\multiput(42.26,35.25)(0.49,-0.08){1}{\line(1,0){0.49}}
\multiput(42.75,35.17)(0.5,-0.07){1}{\line(1,0){0.5}}
\multiput(43.25,35.1)(0.5,-0.05){1}{\line(1,0){0.5}}
\multiput(43.75,35.05)(0.5,-0.03){1}{\line(1,0){0.5}}
\multiput(44.25,35.02)(0.5,-0.02){1}{\line(1,0){0.5}}
\put(44.75,35){\line(1,0){0.5}}
\multiput(45.25,35)(0.5,0.02){1}{\line(1,0){0.5}}
\multiput(45.75,35.02)(0.5,0.03){1}{\line(1,0){0.5}}
\multiput(46.25,35.05)(0.5,0.05){1}{\line(1,0){0.5}}
\multiput(46.75,35.1)(0.5,0.07){1}{\line(1,0){0.5}}
\multiput(47.25,35.17)(0.49,0.08){1}{\line(1,0){0.49}}
\multiput(47.74,35.25)(0.49,0.1){1}{\line(1,0){0.49}}
\multiput(48.23,35.35)(0.49,0.12){1}{\line(1,0){0.49}}
\multiput(48.72,35.47)(0.48,0.13){1}{\line(1,0){0.48}}
\multiput(49.2,35.6)(0.48,0.15){1}{\line(1,0){0.48}}
\multiput(49.68,35.75)(0.47,0.16){1}{\line(1,0){0.47}}
\multiput(50.16,35.91)(0.23,0.09){2}{\line(1,0){0.23}}
\multiput(50.62,36.09)(0.23,0.1){2}{\line(1,0){0.23}}
\multiput(51.09,36.29)(0.23,0.11){2}{\line(1,0){0.23}}
\multiput(51.54,36.5)(0.22,0.11){2}{\line(1,0){0.22}}
\multiput(51.99,36.73)(0.22,0.12){2}{\line(1,0){0.22}}
\multiput(52.43,36.97)(0.22,0.13){2}{\line(1,0){0.22}}
\multiput(52.86,37.22)(0.21,0.13){2}{\line(1,0){0.21}}
\multiput(53.28,37.49)(0.21,0.14){2}{\line(1,0){0.21}}
\multiput(53.69,37.78)(0.2,0.15){2}{\line(1,0){0.2}}
\multiput(54.1,38.07)(0.13,0.1){3}{\line(1,0){0.13}}
\multiput(54.49,38.38)(0.13,0.11){3}{\line(1,0){0.13}}
\multiput(54.87,38.71)(0.12,0.11){3}{\line(1,0){0.12}}
\multiput(55.25,39.04)(0.12,0.12){3}{\line(1,0){0.12}}
\multiput(55.61,39.39)(0.12,0.12){3}{\line(0,1){0.12}}
\multiput(55.96,39.75)(0.11,0.12){3}{\line(0,1){0.12}}
\multiput(56.29,40.13)(0.11,0.13){3}{\line(0,1){0.13}}
\multiput(56.62,40.51)(0.1,0.13){3}{\line(0,1){0.13}}
\multiput(56.93,40.9)(0.15,0.2){2}{\line(0,1){0.2}}
\multiput(57.22,41.31)(0.14,0.21){2}{\line(0,1){0.21}}
\multiput(57.51,41.72)(0.13,0.21){2}{\line(0,1){0.21}}
\multiput(57.78,42.14)(0.13,0.22){2}{\line(0,1){0.22}}
\multiput(58.03,42.57)(0.12,0.22){2}{\line(0,1){0.22}}
\multiput(58.27,43.01)(0.11,0.22){2}{\line(0,1){0.22}}
\multiput(58.5,43.46)(0.11,0.23){2}{\line(0,1){0.23}}
\multiput(58.71,43.91)(0.1,0.23){2}{\line(0,1){0.23}}
\multiput(58.91,44.38)(0.09,0.23){2}{\line(0,1){0.23}}
\multiput(59.09,44.84)(0.16,0.47){1}{\line(0,1){0.47}}
\multiput(59.25,45.32)(0.15,0.48){1}{\line(0,1){0.48}}
\multiput(59.4,45.8)(0.13,0.48){1}{\line(0,1){0.48}}
\multiput(59.53,46.28)(0.12,0.49){1}{\line(0,1){0.49}}
\multiput(59.65,46.77)(0.1,0.49){1}{\line(0,1){0.49}}
\multiput(59.75,47.26)(0.08,0.49){1}{\line(0,1){0.49}}
\multiput(59.83,47.75)(0.07,0.5){1}{\line(0,1){0.5}}
\multiput(59.9,48.25)(0.05,0.5){1}{\line(0,1){0.5}}
\multiput(59.95,48.75)(0.03,0.5){1}{\line(0,1){0.5}}
\multiput(59.98,49.25)(0.02,0.5){1}{\line(0,1){0.5}}

\linethickness{0.3mm}
\put(45,50){\line(1,0){45}}
\linethickness{0.3mm}
\multiput(75,60)(0.18,-0.12){83}{\line(1,0){0.18}}
\linethickness{0.3mm}
\multiput(45,50)(0.36,0.12){83}{\line(1,0){0.36}}
\linethickness{0.3mm}
\put(45,60){\line(1,0){30}}
\linethickness{0.3mm}
\multiput(45,60)(0.54,-0.12){83}{\line(1,0){0.54}}
\put(42,50){\makebox(0,0)[cc]{$x$}}

\put(42,60){\makebox(0,0)[cc]{$z$}}

\put(75,64){\makebox(0,0)[cc]{$y'$}}

\put(93,50){\makebox(0,0)[cc]{$y$}}

\end{picture}

 \noindent{\bf Démonstration.}
On applique le lemme~\ref{approx-arbres} à  $\{x,y,y',z\}$  avec $l=1$ et $y'$ comme point base.  Soit $T$ et $\Psi$ comme dans le lemme~\ref{approx-arbres}.

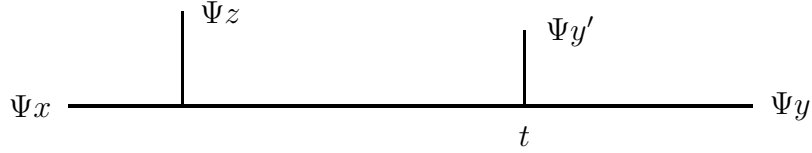
\ifx\JPicScale\undefined\def\JPicScale{1}\fi
\unitlength \JPicScale mm
\begin{picture}(120,30)(20,30)
\linethickness{0.3mm}
\put(30,40){\line(1,0){90}}
\linethickness{0.3mm}
\put(45,40){\line(0,1){12.5}}
\linethickness{0.3mm}
\put(90,40){\line(0,1){10}}
\put(25,40){\makebox(0,0)[cc]{$\Psi x$}}

\put(50,52.5){\makebox(0,0)[cc]{$\Psi z$}}

\put(96,50){\makebox(0,0)[cc]{$\Psi y'$}}

\put(90,36){\makebox(0,0)[cc]{$t$}}

\put(125,40){\makebox(0,0)[cc]{$\Psi y$}}

\end{picture}

\noindent Soit $t$ le point de $\geod(\Psi z, \Psi y)$ à distance minimale de $\Psi  y'$. 
On a donc 
\begin{gather}\label{t-geod-z1-y0-y1}t\in 
\geod(\Psi z, \Psi y) \text{\ \ et  \ \ } t\in 
\geod(\Psi z, \Psi y') \end{gather}
On a  $\Psi  y'\in (\mu+\de)\tg(\Psi z, \Psi y)$ par le a) du lemme~\ref{approx-arbres2}, donc  $d(t,\Psi  y')\leq \frac{\mu+\de}{2}$. Comme $d(\Psi  x, \Psi  y')=d(x,y')> k+\frac{\mu+\de}{2}$ par hypothèse, on a $d(\Psi  x,t)> k$ et comme $d(\Psi  x, \Psi   z)\leq k$, 
\begin{gather}\label{t-geod-x-y0-y1}t\text{\ \ appartient \ \ à \ \ }
\geod(\Psi  x, \Psi y) \text{\ \ et  \ \ à\ \ } 
\geod(\Psi  x, \Psi y').  \end{gather}
Il résulte de (\ref{t-geod-z1-y0-y1}) et  (\ref{t-geod-x-y0-y1}) que $$d(\Psi  x, \Psi y')-d(\Psi  x, \Psi y)=d(t, \Psi y')-d(t, \Psi y)=
d(\Psi  z, \Psi y')-d(\Psi  z, \Psi y)$$ et donc $d(x,y')-d(x,y)\leq d(z,y')-d(z,y)+\de
\leq 
-\nu+\de$ où la dernière inégalité a lieu par hypothèse.  
Enfin par la première partie de (\ref{t-geod-x-y0-y1}), et comme $d(t,\Psi  y')\leq \frac{\mu+\de}{2}$, on a $\Psi y'\in (\mu+\de)\tg (\Psi x,\Psi y)$,  d'où $y'\in (\mu+\de)\tg (x,y)$ par le b) du  lemme~\ref{approx-arbres2}. 
 \cqfd

Le lemme suivant est une conséquence  du précédent. 

\begin{lem}\label{Bxkmunu-cor}
Soient $\mu_{1},\nu_{1},\mu_{2},\nu_{2}\in \N$ vérifiant $$\nu_{1}> \mu_{1}+\frac{3\delta}{2}\text{\ \ \ et\ \ \ }\nu_{2}>  
\mu_{2}+\frac{3\delta}{2}.$$ Soient $k\in \N$ et  $x\in X$. 
Soit $y_{0},...,y_{j}$ une suite de points de $X$ telle que 
$y_{1},\dots ,y_{j}$ n'appartiennent pas à $B(x,k+\frac{\max(\mu_{1},\mu_{2})+\delta}{2})$ et que 
pour tout $i\in \{0,...,j-1\}$, 
\begin{itemize}
\item ou bien $y_{i+1}=y_{i}$, 
\item ou bien il existe $z\in B(x,k)$ tel que  $$y_{i+1}\in \mu_{1}\text{-}\geod(z,y_{i})\text{\ \ et \ \ }d(z,y_{i+1})\leq d(z,y_{i})-\nu_{1},$$ 
\item ou bien il existe $z\in B(x,k)$ tel que  $$y_{i+1}\in \mu_{2}\text{-}\geod(z,y_{i})\text{\ \ et \ \ }d(z,y_{i+1})\leq d(z,y_{i})-\nu_{2}. $$
\end{itemize}

\noindent a) On a $d(x,y_{j})\leq \pp \leq d(x,y_{1})\leq d(x,y_{0})$ et 
pour $i\in \{0,\pp,j-1\}$ on a  $d(x,y_{i+1})<d(x,y_{i})$ si $y_{i}\neq y_{i+1}$. 

\noindent b) 
Soit $h$ le nombre de valeurs prises par la suite $y_{0},\pp,y_{j}$. 
Alors  \begin{gather*}y_{j}\in \big(\max(\mu_{1},\mu_{2})+2\delta\big)\text{-}\geod(x,y_{0})\\ \text{et\ \ \ }d(x,y_{j})\leq d(x,y_{0})-(h-1)(\min(\nu_{1},\nu_{2})-\de) .\end{gather*} 
\end{lem}
 \noindent{\bf Démonstration.}
Montrons a). Soit $i\in \{0,\pp,j-1\}$ tel que $y_{i}\neq y_{i+1}$. On va montrer $d(x,y_{i+1})<d(x,y_{i})$. Il existe  $c\in \{1,2\}$ et  $z\in B(x,k)$ tels que  $$y_{i+1}\in \mu_{c}\text{-}\geod(z,y_{i})\text{\ \ et \ \ }d(z,y_{i+1})\leq d(z,y_{i})-\nu_{c}.$$ 
Comme $y_{i+1}\not\in B(x,k+\frac{\mu_{c}+\de}{2})$,
en appliquant le  lemme~\ref{Bxkmunu} à $(y_{i},y_{i+1},z)$ au lieu de 
$(y,y',z)$ et $(\mu_{c},\nu_{c})$ au lieu de  $(\mu,\nu)$, on obtient $d(x,y_{i+1})\leq d(x,y_{i})-\nu_{c}+\de< d(x,y_{i})$ puisque $\nu_{c}>\de$.

\noindent Pour montrer b), 
on procède par récurrence ascendante. On pose $$\lambda=\max(\mu_{1},\mu_{2})+2\delta. $$
Pour $i\in \{0,\pp,j\}$ on note $h_{i}$ le nombre de valeurs prises par la suite $y_{0},\pp,y_{i}$, de sorte que $1=h_{0}\leq h_{1}\leq \pp \leq h_{j}=h$. 
Par l'hypothèse de récurrence on a   \begin{gather}\label{hrec-Bxkmunu-cor}y_{i}\in \lambda\text{-}\geod(x,y_{0})\text{\  et \ } d(x,y_{i})\leq d(x,y_{0})-(h_{i}-1)(\min(\nu_{1},\nu_{2})-\de).\end{gather} Si $y_{i+1}=y_{i}$, on a $h_{i+1}=h_{i}$ et $y_{i+1}$ satisfait l'hypothèse de récurrence. Sinon, soit 
 $c\in\{1,2\}$ et $z\in B(x,k)$ tel que  $$y_{i+1}\in \mu_{c}\text{-}\geod(z,y_{i})\text{ \  et \ }
d(z,y_{i+1})\leq d(z,y_{i})-\nu_{c}.$$ 
En appliquant le lemme~\ref{Bxkmunu} à \begin{gather*}(y_{i},y_{i+1},z)\text{ au lieu de }(y,y',z)\text{\  et \ }  (\mu_{c},\nu_{c})\text{ au lieu de }(\mu,\nu)\end{gather*}   on obtient  
\begin{gather}\label{etape-rec-Bxkmunu-cor}y_{i+1}\in (\mu_{c}+\de)\tg(x,y_{i}) \text{\ \ et \ \ }d(x,y_{i+1})\leq d(x,y_{i})-\nu_{c}+\de. \end{gather} 
Comme $d(y_{i},y_{i+1})\geq d(z,y_{i})-d(z,y_{i+1})\geq \nu_{c}$ et grâce aux premières parties de  (\ref{hrec-Bxkmunu-cor}) et (\ref{etape-rec-Bxkmunu-cor}), le lemme~\ref{iter2} appliqué à $(x,y_{0},y_{i},y_{i+1})$ au lieu de $(x,a,b,c)$ et $(\lambda,\mu_{c}+\de)$ au lieu de $(\alpha,\beta)$ montre que 
$y_{i+1}\in \lambda \tg(x,y_{0})$ puisque $$\max(\lambda+2(\mu_{c}+\de)-2\nu_{c}, \mu_{c}+\de)+\de\leq \lambda.$$ D'autre part les deuxièmes parties de 
(\ref{hrec-Bxkmunu-cor}) et (\ref{etape-rec-Bxkmunu-cor}) impliquent immédiatement $$d(x,y_{i+1})\leq d(x,y_{0})-(h_{i+1}-1)(\min(\nu_{1},\nu_{2})-\de).$$
 \cqfd

On rappelle que  pour $x\in X$ et $S\in \Delta$, on note $d_{\max}(x,S)=\max_{y\in S}d(x,y)$. 

\begin{lem} \label{m-fini} Il existe une constante $D=C(\de,K,N,Q,P)$ telle que le résultat suivant soit vrai. 
Soit $S_{0}\in \Delta$, $x\in X$,  $k\in \N$. 
Alors pour tout $m\in \N$ et pour toute suite $S_{1},\pp,S_{m}$ de $\Delta$ vérifiant les conditions i) et ii) de la définition~\ref{defi-Y}, c'est-à-dire 
\begin{itemize}
\item i) pour tout $i\in \{0,\pp,m-1\}$, \begin{gather*}S_{i+1}\subset S_{i}\cup 
\bigcup _{\tilde x\in B(x,k),a\in S_{i}} \{y\in 4\de\tg(\tilde x,a), d(y,a)\in ]N-2\de, QN]\}\\ 
\cup 
\bigcup _{\tilde x\in B(x,k),a\in S_{i}}
\{z\in F\tg(\tilde x,a), d(z,a) \geq \frac{Q}{F}
\},\end{gather*}
\item ii) pour tout $i\in  \{1,\dots,m\}$,   $d(x,S_{i})>k+P$, 
\end{itemize}
on a 
\begin{gather}\label{S0-...Sm}S_{0}\cup \pp \cup S_{m}\subset  \bigcup_{a\in S_{0}}
  \big(F+2\de\big)\text{-}\geod(x,a),\\ \label{maxS0-...Sm}
  d_{\max}(x,S_{0})\geq d_{\max}(x,S_{1})\geq \pp \geq d_{\max}(x,S_{m}),\end{gather}
  et 
  \begin{itemize}
  \item si $d(x,S_{0})\leq k$ on a $m=0$ 
  \item si $d(x,S_{0})> k$, 
    le nombre de valeurs prises par la suite $S_{0},\pp,S_{m}$ est 
inférieur ou égal à $ D(d(x,S_{0})-k)$  et le nombre de possibilités pour  $(S_{1},\pp,S_{m})$ est fini et majoré par $e^{D(d(x,S_{0})-k+m)}$. 
\end{itemize}
\end{lem}
 \noindent{\bf Démonstration.}
On  applique le  lemme~\ref{Bxkmunu-cor}  à $$\mu_{1}=4\de, \ \nu_{1}=N-6\de, \ \mu_{2}=F \text{\  et   }  
\nu_{2}=Q/F-F.$$
Pour tout $i\in \{0,\pp,m-1\}$ et tout $y_{i+1}\in S_{i+1}$ il existe $y_{i}\in S_{i}$ tel que 
\begin{equation*} (C_{i}) \ \left\{
\begin{array}{rl}
\text{ -ou bien }& y_{i+1}=y_{i},\\  
\text{ -ou bien }&\text{il existe }\tilde x\in B(x,k)\text{  tel que  }\\ & y_{i+1}\in \mu_{1}\text{-}\geod(\tilde x,y_{i})\text{\ \ et \ \ }d(\tilde x,y_{i+1})\leq d(\tilde x,y_{i})-\nu_{1},\\
\text{ -ou bien }&\text{il existe }\tilde x\in B(x,k)\text{  tel que  }\\ & y_{i+1}\in \mu_{2}\text{-}\geod(\tilde x,y_{i})\text{\ \ et \ \ }d(\tilde x,y_{i+1})\leq d(\tilde x,y_{i})-\nu_{2}. 
\end{array} \right.
\end{equation*}
On suppose $\nu_{1}> \mu_{1}+\frac{3\delta}{2}$, $\nu_{2}>  
\mu_{2}+\frac{3\delta}{2}$ et $P\geq \frac{\max(\mu_{1},\mu_{2})+\de}{2}$, ce qui est permis par $(H_{N})$,  $(H_{Q})$ et $(H_{P})$ respectivement.  
Soit $i\in \{1,...,m\}$ et $y_{i}\in S_{i}$. Il existe $y_{i-1}\in S_{i-1},...,y_{0}\in S_{0}$ tels que les conditions $(C_{i-1}),...,(C_{0})$ soient satisfaites. 
Le a) du lemme~\ref{Bxkmunu-cor}   montre alors que $d(x,y_{i})\leq d(x,y_{i-1})$ et comme $y_{i}\in S_{i}$ est arbitraire  il résulte que $d_{\max}(x,S_{i})\leq d_{\max}(x,S_{i-1})$ et on a  montré (\ref{maxS0-...Sm}). 
Comme $\max(\mu_{1},\mu_{2})+2\delta=F+2\de$, 
 le b) du lemme~\ref{Bxkmunu-cor}   montre que 
   $$y_{i}\in   \bigcup_{a\in S_{0}}
  \big(F+2\de\big)\text{-}\geod(x,a)$$ 
  %Par (\ref{def-F}) on a $F+2\de\geq 2N$ donc $S_{0}\subset     \bigcup_{a\in S_{0}}
  %\big(F+2\de\big)\text{-}\geod(x,a)$, 
  et on a  montré (\ref{S0-...Sm}).

  Pour montrer la suite de l'énoncé on suppose d'abord $d(x,S_{0})\leq k$. 
  Alors $m=0$ par (\ref{maxS0-...Sm}) et par la condition ii), car $P\geq N$. 
  On suppose   $d(x,S_{0})> k$ dans toute la suite de la démonstration. 
    Soit $y_{m}\in S_{m}$ et soient $y_{m-1}\in S_{m-1}, ... , y_{0}\in S_{0}$ tels que les conditions $(C_{m-1}),...,(C_{0})$ soient satisfaites. 
  On note $h$ le nombre de valeurs différentes prises par la suite $y_{0},\pp,y_{m}$.
  On suppose $\min(\nu_{1},\nu_{2})=N-6\de$, ce qui est permis par $(H_{Q})$. 
   Le b) du lemme~\ref{Bxkmunu-cor} montre que $$y_{m}\in \big(F+2\de\big)\tg(x,y_{0})\text{\ \  et\ \ } 
    d(x,y_{m})\leq d(x,y_{0})-(h-1)(N-7\de).$$
   On a $d(x,y_{m})\geq k+P$ puisque $m\geq 1$. On en déduit \begin{gather}\label{ineg-h-1-N-8}(h-1)(N-7\de)\leq d(x,S_{0})+N-(k+P). \end{gather}

 On suppose $
  N-7\de \geq 1$ et $P\geq N+1$, ce qui est permis par $(H_{N})$ et $(H_{P})$. Alors  (\ref{ineg-h-1-N-8}) implique   $h\leq d(x,S_{0})-k$. 
  Il existe une constante $C_{1}=C(\de,K,N)$ telle que tout point de $X$ appartienne au plus à $C_{1}$ éléments de $\Delta$. En notant $l$ le nombre de valeurs prises par la suite $S_{0},\pp,S_{m}$, on a $h\geq\frac{l}{C_{1}}$, d'où $$l\leq C_{1}h\leq C_{1}\big(d(x,S_{0})-k\big).$$

  Il existe une constante $C_{2}=C(\de,K,N)$ telle que, pour $i\in \{0,...,m-1\}$, connaissant $S_{i}$ et $d_{\max}(x,S_{i+1}) $ le nombre de possibilités pour $S_{i+1}$ vérifiant les conditions de l'énoncé est inférieur ou égal  $C_{2}$.  En effet pour $y_{i+1}\in S_{i+1}$ il existe $y_{i}\in S_{i}$ 
  tel que $y_{i+1}\in (F+2\de)\tg(x,y_{i})$ d'après \eqref{S0-...Sm} appliqué à 
  $(S_{i},...,S_{m})$, et on déduit l'existence de $C_{2}$ du  lemme~\ref{cardinal-tranche-geod} appliqué à $(x,y_{i})$ au lieu de $(x,y)$.

  Etant donné $m$, le
   nombre de possibilités pour  les entiers $d_{\max}(x,S_{i})$ (pour $i=1,\pp,m$) qui vérifient nécessairement $$d_{\max}(x,S_{0})\geq d_{\max}(x,S_{1}) \geq \pp   \geq d_{\max}(x,S_{m})\geq k+P$$  est inférieur ou égal à 
   $$\max\Big(1,\binom{m+d_{\max}(x,S_{0})-k-P}{ m}\Big)\leq 
 2^{d(x,S_{0})-k+m}$$ car $P\geq N$. 
  Le
   nombre de possibilités pour $(S_{1},\pp,S_{m})$ est donc inférieur ou égal à $(C_{2})^{m}
   2^{d(x,S_{0})-k+m}$. 
   Ceci  termine la démonstration du lemme~\ref{m-fini}. \cqfd
  
  Les deux lemmes suivants  sont des conséquences du lemme~\ref{m-fini} et  serviront ultérieurement. 
  
  \begin{lem}\label{dmaxY-dmaxS}
 Pour $m,l_{0},...,l_{m}\in \N$ et 
 $$(a_{1},\dots,a_{p},S_{0},...,S_{m},(\mathcal Y_{i}^{j})_{i\in \{0,\dots,m\}, j\in \{1,\dots ,l_{i}\}})\in 
 Y_{x}^{p,k,m,(l_{0},...,l_{m})}$$ 
 on a 
  pour $i\in \{0,...,m\}$ et  $ j\in \{1,...,l_{i}\}$, $$  d_{\max}(x,\mathcal Y_{i}^{j})\leq d_{\max}(x,S_{i})+2P+\de.$$
 \end{lem}
 \noindent{\bf Démonstration.}
Le lemme~\ref{boules-geod} appliqué à $\alpha=2P$ donne
\begin{itemize}
\item pour $i\in \{0,...,m-1\}$, $$d_{\max}(x,\mathcal Y_{i}^{j})\leq \max(d_{\max}(x,S_{i}),d_{\max}(x,S_{i+1}))+2P+\de,$$ d'où le résultat puisque 
$d_{\max}(x,S_{i})\geq d_{\max}(x,S_{i+1})$ d'après le lemme~\ref{m-fini},
\item pour $i=m$, $d_{\max}(x,\mathcal Y_{m}^{j})\leq \max(d_{\max}(x,S_{m}),k)+2P+\de$ d'où le résultat puisque $l_{m}=0$ si $d_{\max}(x,S_{m})\leq k$, par la remarque qui suit la définition~\ref{defi-Y}. \cqfd
\end{itemize}

\begin{lem}\label{lemme-S0-...Sm}
Soient $p\in \{1,...,p_{\max}\}$, $k,m,l_{0},...,l_{m}\in \N$ et  
$$(a_{1},\dots,a_{p},S_{0},...,S_{m},(\mathcal Y_{i}^{j})_{i\in \{0,\dots,m\}, j\in \{1,\dots ,l_{i}\}})\in Y_{x}^{p,k,m,(l_{0},...,l_{m})}.$$ 
Soit $b\in S_{0}$ et $u$ un point de $B(x,k)$ à distance minimale de $b$. 
Alors

\noindent a) $S_{0}\cup \pp\cup S_{m}\subset 2F\tg(b,u)\subset P\tg(b,u)$, 

%\noindent (b) pour tout $i\geq 0$,  $$\bigcup _{j\in \{1,\dots ,l_{i}\}}\mathcal Y_{i}^{j}
%\subset \bigcup_{a\in S_{i},t\in B(x,k)\text{ à  distance minimale de }a} 
%(4N+5\de)\tg(t,a),$$
%
%\noindent (c) $S_{0}\cup \bigcup_{j\in \{1,\dots ,l_{0}\}}\{y_{0}^{j}\}\subset 
%(8N+20\de p_{\max}+20\de)\tg(b,u)$, 
\noindent b) $\bigcup _{j\in \{1,...,l_{m}\}}\mathcal Y_{m}^{j}\in \bigcup_{a\in S_{m}}(2P+\de)\tg(x,a)$,

\noindent c) pour tout $i\in \{0,\dots,m\}$, $\bigcup _{ j\in \{1,\dots ,l_{i}\}}\mathcal Y_{i}^{j}\subset 
4P\tg(b,u)$. 

\noindent d) pour tout  $i\in \{0,\dots,m-1\}$, et tout $j\in \{1,\dots ,l_{i}\}$, 
$d(x,\mathcal Y_{i}^{j})\geq 
d_{\max}(x,S_{i+1})-4P$
\end{lem}
 \noindent{\bf Démonstration.}
On commence par  traiter le cas où $b\in B(x,k)$. Alors $d(x,S_{0})\leq k$, d'où $m=0$ par le lemme~\ref{m-fini}. De plus $l_{0}=0$ par la remarque qui suit la définition~\ref{defi-Y}, et les assertions a), b), c)  et d) sont  évidentes. On suppose maintenant que $b\not\in B(x,k)$. En particulier $d(x,u)=k$ et $d(u,b)=d(x,b)-k$.

Montrons a). D'après le lemme~\ref{m-fini}, on a 
$$S_{0}\cup \pp \cup S_{m}\subset  \bigcup_{a\in S_{0}}
  \big(F+2\de\big)\text{-}\geod( x,a).$$
Comme $d(a,b)\leq N$ pour tout $a\in S_{0}$, le lemme~\ref{xx'yy'zz'} montre que 
$$S_{0}\cup \pp \cup S_{m}\subset  
  \big(F+2N+2\de\big)\text{-}\geod( x,b).$$
Pour tout $y\in S_{0}$ on a $d(y,b)\leq N$, donc $y\in 2F\tg(b,u)$ car $F\geq N$ par (\ref{def-F}). Soit donc $i\geq 1$ et $y\in S_{i}$, et montrons $y\in 
 2F\tg(b,u)$. On a $y\in \big(F+2N+2\de\big)\text{-}\geod( x,b)$ et $d(x,y)\geq d(x,u)+P$.

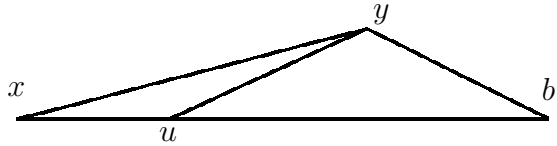
 \ifx\JPicScale\undefined\def\JPicScale{1}\fi
\unitlength \JPicScale mm
\begin{picture}(90,25)(35,25)
\linethickness{0.3mm}
\put(60,30){\line(1,0){70}}
\linethickness{0.3mm}
\multiput(106,42)(0.24,-0.12){100}{\line(1,0){0.24}}
\linethickness{0.3mm}
\multiput(60,30)(0.46,0.12){100}{\line(1,0){0.46}}
\linethickness{0.3mm}
\multiput(80,30)(0.26,0.12){100}{\line(1,0){0.26}}
\put(60,34){\makebox(0,0)[cc]{$x$}}

\put(80,28){\makebox(0,0)[cc]{$u$}}

\put(130,34){\makebox(0,0)[cc]{$b$}}

\put(108,44){\makebox(0,0)[cc]{$y$}}

\end{picture}
 
 \noindent Par $(H^{F+2N+2\de}_{\de}(u,x,y,b))$ on a $$d(u,y)\leq \max(d(u,x)-d(x,y),d(u,b)-d(b,y))+F+2N+3\de.$$ Or  $d(u,x)-d(x,y)+F+2N+3\de\leq -P+F+2N+3\de$ et on suppose $-P+F+2N+3\de<0$, ce qui est permis par $(H_{P})$. Donc  $d(u,y)\leq d(u,b)-d(b,y)+F+2N+3\de$, c'est-à-dire  $y\in 
 (F+2N+3\de)\tg(u,b)$. On en déduit $y\in  2F\tg(u,b)$ puisque $F\geq 2N+3\de$ par (\ref{def-F}). Enfin on suppose $P\geq 2F$, ce qui est permis par $(H_{P})$. 

\noindent On va montrer maintenant b), ainsi que c) dans le cas où $i=m$. 
Soit $j\in \{1,...,l_{m}\}$ et $y\in 
\mathcal Y_{m}^{j}$. On a $d(x,y)\geq k+3P$ et il existe $a\in S_{m}$  et $\tilde x\in B(x,k)$ tels  que 
$y\in 2P\tg(\tilde x,a)$. On a $a\in P\tg(u,b)$ par le a). Ensuite  $(H_{\de}^{2P}(x,\tilde x,y,a))$ implique 
$d(x,y)\leq \max(k,d(x,a)-d(a,y))+2P+\de$. Comme $d(x,y)\geq k+3P$ et $P>\de$  on en déduit $d(x,y)\leq d(x,a)-d(a,y)+2P+\de$, c'est-à-dire $y\in (2P+\de)\tg(x,a)$, ce qui montre déjà b).  
Par $(H_{\de}^{2P+\de}(u,x,y,a))$ on a $d(u,y)\leq \max (d(u,x)-d(x,y), d(u,a)-d(a,y))+2P+2\de$ et comme $d(u,x)-d(x,y)\leq -3P$ et $P>2\de$ on en déduit $d(u,y)\leq d(u,a)-d(a,y)+2P+2\de$, c'est-à-dire $y\in (2P+2\de)\tg(a,u)$. 

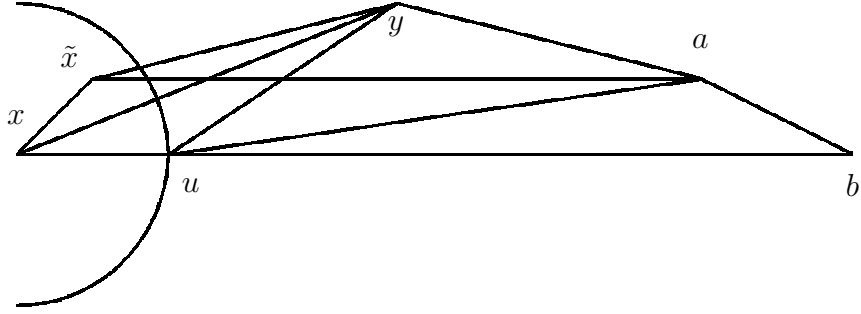
\ifx\JPicScale\undefined\def\JPicScale{1}\fi
\unitlength \JPicScale mm
\begin{picture}(120,45)(0,30)
%\begin{picture}(120,70)(0,0)
\linethickness{0.3mm}
\multiput(10,30)(0.5,0.01){1}{\line(1,0){0.5}}
\multiput(10.5,30.01)(0.5,0.02){1}{\line(1,0){0.5}}
\multiput(11,30.02)(0.5,0.03){1}{\line(1,0){0.5}}
\multiput(11.49,30.06)(0.5,0.04){1}{\line(1,0){0.5}}
\multiput(11.99,30.1)(0.5,0.06){1}{\line(1,0){0.5}}
\multiput(12.49,30.16)(0.49,0.07){1}{\line(1,0){0.49}}
\multiput(12.98,30.22)(0.49,0.08){1}{\line(1,0){0.49}}
\multiput(13.47,30.3)(0.49,0.09){1}{\line(1,0){0.49}}
\multiput(13.96,30.4)(0.49,0.1){1}{\line(1,0){0.49}}
\multiput(14.45,30.5)(0.48,0.12){1}{\line(1,0){0.48}}
\multiput(14.94,30.62)(0.48,0.13){1}{\line(1,0){0.48}}
\multiput(15.42,30.75)(0.48,0.14){1}{\line(1,0){0.48}}
\multiput(15.9,30.89)(0.47,0.15){1}{\line(1,0){0.47}}
\multiput(16.37,31.04)(0.47,0.16){1}{\line(1,0){0.47}}
\multiput(16.84,31.21)(0.47,0.18){1}{\line(1,0){0.47}}
\multiput(17.31,31.38)(0.23,0.09){2}{\line(1,0){0.23}}
\multiput(17.77,31.57)(0.23,0.1){2}{\line(1,0){0.23}}
\multiput(18.23,31.77)(0.23,0.11){2}{\line(1,0){0.23}}
\multiput(18.68,31.98)(0.22,0.11){2}{\line(1,0){0.22}}
\multiput(19.12,32.2)(0.22,0.12){2}{\line(1,0){0.22}}
\multiput(19.57,32.44)(0.22,0.12){2}{\line(1,0){0.22}}
\multiput(20,32.68)(0.21,0.13){2}{\line(1,0){0.21}}
\multiput(20.43,32.93)(0.21,0.13){2}{\line(1,0){0.21}}
\multiput(20.85,33.2)(0.21,0.14){2}{\line(1,0){0.21}}
\multiput(21.27,33.48)(0.2,0.14){2}{\line(1,0){0.2}}
\multiput(21.67,33.76)(0.2,0.15){2}{\line(1,0){0.2}}
\multiput(22.08,34.06)(0.13,0.1){3}{\line(1,0){0.13}}
\multiput(22.47,34.36)(0.13,0.11){3}{\line(1,0){0.13}}
\multiput(22.86,34.68)(0.13,0.11){3}{\line(1,0){0.13}}
\multiput(23.23,35)(0.12,0.11){3}{\line(1,0){0.12}}
\multiput(23.6,35.34)(0.12,0.11){3}{\line(1,0){0.12}}
\multiput(23.96,35.68)(0.12,0.12){3}{\line(0,1){0.12}}
\multiput(24.32,36.04)(0.11,0.12){3}{\line(0,1){0.12}}
\multiput(24.66,36.4)(0.11,0.12){3}{\line(0,1){0.12}}
\multiput(25,36.77)(0.11,0.13){3}{\line(0,1){0.13}}
\multiput(25.32,37.14)(0.11,0.13){3}{\line(0,1){0.13}}
\multiput(25.64,37.53)(0.1,0.13){3}{\line(0,1){0.13}}
\multiput(25.94,37.92)(0.15,0.2){2}{\line(0,1){0.2}}
\multiput(26.24,38.33)(0.14,0.2){2}{\line(0,1){0.2}}
\multiput(26.52,38.73)(0.14,0.21){2}{\line(0,1){0.21}}
\multiput(26.8,39.15)(0.13,0.21){2}{\line(0,1){0.21}}
\multiput(27.07,39.57)(0.13,0.21){2}{\line(0,1){0.21}}
\multiput(27.32,40)(0.12,0.22){2}{\line(0,1){0.22}}
\multiput(27.56,40.43)(0.12,0.22){2}{\line(0,1){0.22}}
\multiput(27.8,40.88)(0.11,0.22){2}{\line(0,1){0.22}}
\multiput(28.02,41.32)(0.11,0.23){2}{\line(0,1){0.23}}
\multiput(28.23,41.77)(0.1,0.23){2}{\line(0,1){0.23}}
\multiput(28.43,42.23)(0.09,0.23){2}{\line(0,1){0.23}}
\multiput(28.62,42.69)(0.18,0.47){1}{\line(0,1){0.47}}
\multiput(28.79,43.16)(0.16,0.47){1}{\line(0,1){0.47}}
\multiput(28.96,43.63)(0.15,0.47){1}{\line(0,1){0.47}}
\multiput(29.11,44.1)(0.14,0.48){1}{\line(0,1){0.48}}
\multiput(29.25,44.58)(0.13,0.48){1}{\line(0,1){0.48}}
\multiput(29.38,45.06)(0.12,0.48){1}{\line(0,1){0.48}}
\multiput(29.5,45.55)(0.1,0.49){1}{\line(0,1){0.49}}
\multiput(29.6,46.04)(0.09,0.49){1}{\line(0,1){0.49}}
\multiput(29.7,46.53)(0.08,0.49){1}{\line(0,1){0.49}}
\multiput(29.78,47.02)(0.07,0.49){1}{\line(0,1){0.49}}
\multiput(29.84,47.51)(0.06,0.5){1}{\line(0,1){0.5}}
\multiput(29.9,48.01)(0.04,0.5){1}{\line(0,1){0.5}}
\multiput(29.94,48.51)(0.03,0.5){1}{\line(0,1){0.5}}
\multiput(29.98,49)(0.02,0.5){1}{\line(0,1){0.5}}
\multiput(29.99,49.5)(0.01,0.5){1}{\line(0,1){0.5}}
\multiput(29.99,50.5)(0.01,-0.5){1}{\line(0,-1){0.5}}
\multiput(29.98,51)(0.02,-0.5){1}{\line(0,-1){0.5}}
\multiput(29.94,51.49)(0.03,-0.5){1}{\line(0,-1){0.5}}
\multiput(29.9,51.99)(0.04,-0.5){1}{\line(0,-1){0.5}}
\multiput(29.84,52.49)(0.06,-0.5){1}{\line(0,-1){0.5}}
\multiput(29.78,52.98)(0.07,-0.49){1}{\line(0,-1){0.49}}
\multiput(29.7,53.47)(0.08,-0.49){1}{\line(0,-1){0.49}}
\multiput(29.6,53.96)(0.09,-0.49){1}{\line(0,-1){0.49}}
\multiput(29.5,54.45)(0.1,-0.49){1}{\line(0,-1){0.49}}
\multiput(29.38,54.94)(0.12,-0.48){1}{\line(0,-1){0.48}}
\multiput(29.25,55.42)(0.13,-0.48){1}{\line(0,-1){0.48}}
\multiput(29.11,55.9)(0.14,-0.48){1}{\line(0,-1){0.48}}
\multiput(28.96,56.37)(0.15,-0.47){1}{\line(0,-1){0.47}}
\multiput(28.79,56.84)(0.16,-0.47){1}{\line(0,-1){0.47}}
\multiput(28.62,57.31)(0.18,-0.47){1}{\line(0,-1){0.47}}
\multiput(28.43,57.77)(0.09,-0.23){2}{\line(0,-1){0.23}}
\multiput(28.23,58.23)(0.1,-0.23){2}{\line(0,-1){0.23}}
\multiput(28.02,58.68)(0.11,-0.23){2}{\line(0,-1){0.23}}
\multiput(27.8,59.12)(0.11,-0.22){2}{\line(0,-1){0.22}}
\multiput(27.56,59.57)(0.12,-0.22){2}{\line(0,-1){0.22}}
\multiput(27.32,60)(0.12,-0.22){2}{\line(0,-1){0.22}}
\multiput(27.07,60.43)(0.13,-0.21){2}{\line(0,-1){0.21}}
\multiput(26.8,60.85)(0.13,-0.21){2}{\line(0,-1){0.21}}
\multiput(26.52,61.27)(0.14,-0.21){2}{\line(0,-1){0.21}}
\multiput(26.24,61.67)(0.14,-0.2){2}{\line(0,-1){0.2}}
\multiput(25.94,62.08)(0.15,-0.2){2}{\line(0,-1){0.2}}
\multiput(25.64,62.47)(0.1,-0.13){3}{\line(0,-1){0.13}}
\multiput(25.32,62.86)(0.11,-0.13){3}{\line(0,-1){0.13}}
\multiput(25,63.23)(0.11,-0.13){3}{\line(0,-1){0.13}}
\multiput(24.66,63.6)(0.11,-0.12){3}{\line(0,-1){0.12}}
\multiput(24.32,63.96)(0.11,-0.12){3}{\line(0,-1){0.12}}
\multiput(23.96,64.32)(0.12,-0.12){3}{\line(0,-1){0.12}}
\multiput(23.6,64.66)(0.12,-0.11){3}{\line(1,0){0.12}}
\multiput(23.23,65)(0.12,-0.11){3}{\line(1,0){0.12}}
\multiput(22.86,65.32)(0.13,-0.11){3}{\line(1,0){0.13}}
\multiput(22.47,65.64)(0.13,-0.11){3}{\line(1,0){0.13}}
\multiput(22.08,65.94)(0.13,-0.1){3}{\line(1,0){0.13}}
\multiput(21.67,66.24)(0.2,-0.15){2}{\line(1,0){0.2}}
\multiput(21.27,66.52)(0.2,-0.14){2}{\line(1,0){0.2}}
\multiput(20.85,66.8)(0.21,-0.14){2}{\line(1,0){0.21}}
\multiput(20.43,67.07)(0.21,-0.13){2}{\line(1,0){0.21}}
\multiput(20,67.32)(0.21,-0.13){2}{\line(1,0){0.21}}
\multiput(19.57,67.56)(0.22,-0.12){2}{\line(1,0){0.22}}
\multiput(19.12,67.8)(0.22,-0.12){2}{\line(1,0){0.22}}
\multiput(18.68,68.02)(0.22,-0.11){2}{\line(1,0){0.22}}
\multiput(18.23,68.23)(0.23,-0.11){2}{\line(1,0){0.23}}
\multiput(17.77,68.43)(0.23,-0.1){2}{\line(1,0){0.23}}
\multiput(17.31,68.62)(0.23,-0.09){2}{\line(1,0){0.23}}
\multiput(16.84,68.79)(0.47,-0.18){1}{\line(1,0){0.47}}
\multiput(16.37,68.96)(0.47,-0.16){1}{\line(1,0){0.47}}
\multiput(15.9,69.11)(0.47,-0.15){1}{\line(1,0){0.47}}
\multiput(15.42,69.25)(0.48,-0.14){1}{\line(1,0){0.48}}
\multiput(14.94,69.38)(0.48,-0.13){1}{\line(1,0){0.48}}
\multiput(14.45,69.5)(0.48,-0.12){1}{\line(1,0){0.48}}
\multiput(13.96,69.6)(0.49,-0.1){1}{\line(1,0){0.49}}
\multiput(13.47,69.7)(0.49,-0.09){1}{\line(1,0){0.49}}
\multiput(12.98,69.78)(0.49,-0.08){1}{\line(1,0){0.49}}
\multiput(12.49,69.84)(0.49,-0.07){1}{\line(1,0){0.49}}
\multiput(11.99,69.9)(0.5,-0.06){1}{\line(1,0){0.5}}
\multiput(11.49,69.94)(0.5,-0.04){1}{\line(1,0){0.5}}
\multiput(11,69.98)(0.5,-0.03){1}{\line(1,0){0.5}}
\multiput(10.5,69.99)(0.5,-0.02){1}{\line(1,0){0.5}}
\multiput(10,70)(0.5,-0.01){1}{\line(1,0){0.5}}

\linethickness{0.3mm}
\put(10,50){\line(1,0){110}}
\linethickness{0.3mm}
\multiput(10,50)(0.12,0.12){83}{\line(1,0){0.12}}
\linethickness{0.3mm}
\multiput(20,60)(0.48,0.12){83}{\line(1,0){0.48}}
\linethickness{0.3mm}
\put(20,60){\line(1,0){80}}
\linethickness{0.3mm}
\multiput(60,70)(0.48,-0.12){83}{\line(1,0){0.48}}
\linethickness{0.3mm}
\multiput(100,60)(0.24,-0.12){83}{\line(1,0){0.24}}
\linethickness{0.3mm}
\multiput(10,50)(0.3,0.12){167}{\line(1,0){0.3}}
\linethickness{0.3mm}
\multiput(30,50)(0.84,0.12){83}{\line(1,0){0.84}}
\linethickness{0.3mm}
\multiput(30,50)(0.18,0.12){167}{\line(1,0){0.18}}
\put(10,55){\makebox(0,0)[cc]{$x$}}
\put(33,46){\makebox(0,0)[cc]{$u$}}

\put(17,63){\makebox(0,0)[cc]{$\tilde x$}}

\put(60,67){\makebox(0,0)[cc]{$y$}}

\put(100,65){\makebox(0,0)[cc]{$a$}}

\put(120,46){\makebox(0,0)[cc]{$b$}}

\end{picture}

\noindent
Comme $a\in P\tg(u,b)$,  le a) du lemme~\ref{geod-comp-xabc}   montre alors $y\in (3P+2\de)\tg(u,b)$ d'où $y\in 4P\tg(u,b)$ car $P\geq 2\de$.

\noindent On  montre maintenant  c) dans le cas où $i\in \{0,...,m-1\}$, ainsi que d). Soit $j\in \{1,...,l_{i}\}$ et $ t\in \mathcal Y_{i}^{j}$. Il existe $y\in S_{i}$ et $z\in S_{i+1}$ tels que $t\in P\tg(y,z)$.  Les hypothèses $y,z\in P\tg(u,b)$ et $
t\in P\tg(y,z)$ suffisent à impliquer $t\in (3P+\de)\tg (u,b)$.

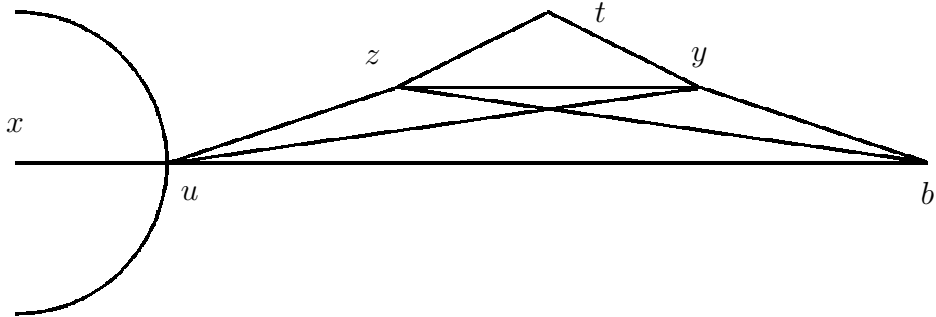
\ifx\JPicScale\undefined\def\JPicScale{1}\fi
\unitlength \JPicScale mm
\begin{picture}(120,45)(5,20)
\linethickness{0.3mm}
\multiput(10,20)(0.5,0.01){1}{\line(1,0){0.5}}
\multiput(10.5,20.01)(0.5,0.02){1}{\line(1,0){0.5}}
\multiput(11,20.02)(0.5,0.03){1}{\line(1,0){0.5}}
\multiput(11.49,20.06)(0.5,0.04){1}{\line(1,0){0.5}}
\multiput(11.99,20.1)(0.5,0.06){1}{\line(1,0){0.5}}
\multiput(12.49,20.16)(0.49,0.07){1}{\line(1,0){0.49}}
\multiput(12.98,20.22)(0.49,0.08){1}{\line(1,0){0.49}}
\multiput(13.47,20.3)(0.49,0.09){1}{\line(1,0){0.49}}
\multiput(13.96,20.4)(0.49,0.1){1}{\line(1,0){0.49}}
\multiput(14.45,20.5)(0.48,0.12){1}{\line(1,0){0.48}}
\multiput(14.94,20.62)(0.48,0.13){1}{\line(1,0){0.48}}
\multiput(15.42,20.75)(0.48,0.14){1}{\line(1,0){0.48}}
\multiput(15.9,20.89)(0.47,0.15){1}{\line(1,0){0.47}}
\multiput(16.37,21.04)(0.47,0.16){1}{\line(1,0){0.47}}
\multiput(16.84,21.21)(0.47,0.18){1}{\line(1,0){0.47}}
\multiput(17.31,21.38)(0.23,0.09){2}{\line(1,0){0.23}}
\multiput(17.77,21.57)(0.23,0.1){2}{\line(1,0){0.23}}
\multiput(18.23,21.77)(0.23,0.11){2}{\line(1,0){0.23}}
\multiput(18.68,21.98)(0.22,0.11){2}{\line(1,0){0.22}}
\multiput(19.12,22.2)(0.22,0.12){2}{\line(1,0){0.22}}
\multiput(19.57,22.44)(0.22,0.12){2}{\line(1,0){0.22}}
\multiput(20,22.68)(0.21,0.13){2}{\line(1,0){0.21}}
\multiput(20.43,22.93)(0.21,0.13){2}{\line(1,0){0.21}}
\multiput(20.85,23.2)(0.21,0.14){2}{\line(1,0){0.21}}
\multiput(21.27,23.48)(0.2,0.14){2}{\line(1,0){0.2}}
\multiput(21.67,23.76)(0.2,0.15){2}{\line(1,0){0.2}}
\multiput(22.08,24.06)(0.13,0.1){3}{\line(1,0){0.13}}
\multiput(22.47,24.36)(0.13,0.11){3}{\line(1,0){0.13}}
\multiput(22.86,24.68)(0.13,0.11){3}{\line(1,0){0.13}}
\multiput(23.23,25)(0.12,0.11){3}{\line(1,0){0.12}}
\multiput(23.6,25.34)(0.12,0.11){3}{\line(1,0){0.12}}
\multiput(23.96,25.68)(0.12,0.12){3}{\line(1,0){0.12}}
\multiput(24.32,26.04)(0.11,0.12){3}{\line(0,1){0.12}}
\multiput(24.66,26.4)(0.11,0.12){3}{\line(0,1){0.12}}
\multiput(25,26.77)(0.11,0.13){3}{\line(0,1){0.13}}
\multiput(25.32,27.14)(0.11,0.13){3}{\line(0,1){0.13}}
\multiput(25.64,27.53)(0.1,0.13){3}{\line(0,1){0.13}}
\multiput(25.94,27.92)(0.15,0.2){2}{\line(0,1){0.2}}
\multiput(26.24,28.33)(0.14,0.2){2}{\line(0,1){0.2}}
\multiput(26.52,28.73)(0.14,0.21){2}{\line(0,1){0.21}}
\multiput(26.8,29.15)(0.13,0.21){2}{\line(0,1){0.21}}
\multiput(27.07,29.57)(0.13,0.21){2}{\line(0,1){0.21}}
\multiput(27.32,30)(0.12,0.22){2}{\line(0,1){0.22}}
\multiput(27.56,30.43)(0.12,0.22){2}{\line(0,1){0.22}}
\multiput(27.8,30.88)(0.11,0.22){2}{\line(0,1){0.22}}
\multiput(28.02,31.32)(0.11,0.23){2}{\line(0,1){0.23}}
\multiput(28.23,31.77)(0.1,0.23){2}{\line(0,1){0.23}}
\multiput(28.43,32.23)(0.09,0.23){2}{\line(0,1){0.23}}
\multiput(28.62,32.69)(0.18,0.47){1}{\line(0,1){0.47}}
\multiput(28.79,33.16)(0.16,0.47){1}{\line(0,1){0.47}}
\multiput(28.96,33.63)(0.15,0.47){1}{\line(0,1){0.47}}
\multiput(29.11,34.1)(0.14,0.48){1}{\line(0,1){0.48}}
\multiput(29.25,34.58)(0.13,0.48){1}{\line(0,1){0.48}}
\multiput(29.38,35.06)(0.12,0.48){1}{\line(0,1){0.48}}
\multiput(29.5,35.55)(0.1,0.49){1}{\line(0,1){0.49}}
\multiput(29.6,36.04)(0.09,0.49){1}{\line(0,1){0.49}}
\multiput(29.7,36.53)(0.08,0.49){1}{\line(0,1){0.49}}
\multiput(29.78,37.02)(0.07,0.49){1}{\line(0,1){0.49}}
\multiput(29.84,37.51)(0.06,0.5){1}{\line(0,1){0.5}}
\multiput(29.9,38.01)(0.04,0.5){1}{\line(0,1){0.5}}
\multiput(29.94,38.51)(0.03,0.5){1}{\line(0,1){0.5}}
\multiput(29.98,39)(0.02,0.5){1}{\line(0,1){0.5}}
\multiput(29.99,39.5)(0.01,0.5){1}{\line(0,1){0.5}}
\multiput(29.99,40.5)(0.01,-0.5){1}{\line(0,-1){0.5}}
\multiput(29.98,41)(0.02,-0.5){1}{\line(0,-1){0.5}}
\multiput(29.94,41.49)(0.03,-0.5){1}{\line(0,-1){0.5}}
\multiput(29.9,41.99)(0.04,-0.5){1}{\line(0,-1){0.5}}
\multiput(29.84,42.49)(0.06,-0.5){1}{\line(0,-1){0.5}}
\multiput(29.78,42.98)(0.07,-0.49){1}{\line(0,-1){0.49}}
\multiput(29.7,43.47)(0.08,-0.49){1}{\line(0,-1){0.49}}
\multiput(29.6,43.96)(0.09,-0.49){1}{\line(0,-1){0.49}}
\multiput(29.5,44.45)(0.1,-0.49){1}{\line(0,-1){0.49}}
\multiput(29.38,44.94)(0.12,-0.48){1}{\line(0,-1){0.48}}
\multiput(29.25,45.42)(0.13,-0.48){1}{\line(0,-1){0.48}}
\multiput(29.11,45.9)(0.14,-0.48){1}{\line(0,-1){0.48}}
\multiput(28.96,46.37)(0.15,-0.47){1}{\line(0,-1){0.47}}
\multiput(28.79,46.84)(0.16,-0.47){1}{\line(0,-1){0.47}}
\multiput(28.62,47.31)(0.18,-0.47){1}{\line(0,-1){0.47}}
\multiput(28.43,47.77)(0.09,-0.23){2}{\line(0,-1){0.23}}
\multiput(28.23,48.23)(0.1,-0.23){2}{\line(0,-1){0.23}}
\multiput(28.02,48.68)(0.11,-0.23){2}{\line(0,-1){0.23}}
\multiput(27.8,49.12)(0.11,-0.22){2}{\line(0,-1){0.22}}
\multiput(27.56,49.57)(0.12,-0.22){2}{\line(0,-1){0.22}}
\multiput(27.32,50)(0.12,-0.22){2}{\line(0,-1){0.22}}
\multiput(27.07,50.43)(0.13,-0.21){2}{\line(0,-1){0.21}}
\multiput(26.8,50.85)(0.13,-0.21){2}{\line(0,-1){0.21}}
\multiput(26.52,51.27)(0.14,-0.21){2}{\line(0,-1){0.21}}
\multiput(26.24,51.67)(0.14,-0.2){2}{\line(0,-1){0.2}}
\multiput(25.94,52.08)(0.15,-0.2){2}{\line(0,-1){0.2}}
\multiput(25.64,52.47)(0.1,-0.13){3}{\line(0,-1){0.13}}
\multiput(25.32,52.86)(0.11,-0.13){3}{\line(0,-1){0.13}}
\multiput(25,53.23)(0.11,-0.13){3}{\line(0,-1){0.13}}
\multiput(24.66,53.6)(0.11,-0.12){3}{\line(0,-1){0.12}}
\multiput(24.32,53.96)(0.11,-0.12){3}{\line(0,-1){0.12}}
\multiput(23.96,54.32)(0.12,-0.12){3}{\line(0,-1){0.12}}
\multiput(23.6,54.66)(0.12,-0.11){3}{\line(1,0){0.12}}
\multiput(23.23,55)(0.12,-0.11){3}{\line(1,0){0.12}}
\multiput(22.86,55.32)(0.13,-0.11){3}{\line(1,0){0.13}}
\multiput(22.47,55.64)(0.13,-0.11){3}{\line(1,0){0.13}}
\multiput(22.08,55.94)(0.13,-0.1){3}{\line(1,0){0.13}}
\multiput(21.67,56.24)(0.2,-0.15){2}{\line(1,0){0.2}}
\multiput(21.27,56.52)(0.2,-0.14){2}{\line(1,0){0.2}}
\multiput(20.85,56.8)(0.21,-0.14){2}{\line(1,0){0.21}}
\multiput(20.43,57.07)(0.21,-0.13){2}{\line(1,0){0.21}}
\multiput(20,57.32)(0.21,-0.13){2}{\line(1,0){0.21}}
\multiput(19.57,57.56)(0.22,-0.12){2}{\line(1,0){0.22}}
\multiput(19.12,57.8)(0.22,-0.12){2}{\line(1,0){0.22}}
\multiput(18.68,58.02)(0.22,-0.11){2}{\line(1,0){0.22}}
\multiput(18.23,58.23)(0.23,-0.11){2}{\line(1,0){0.23}}
\multiput(17.77,58.43)(0.23,-0.1){2}{\line(1,0){0.23}}
\multiput(17.31,58.62)(0.23,-0.09){2}{\line(1,0){0.23}}
\multiput(16.84,58.79)(0.47,-0.18){1}{\line(1,0){0.47}}
\multiput(16.37,58.96)(0.47,-0.16){1}{\line(1,0){0.47}}
\multiput(15.9,59.11)(0.47,-0.15){1}{\line(1,0){0.47}}
\multiput(15.42,59.25)(0.48,-0.14){1}{\line(1,0){0.48}}
\multiput(14.94,59.38)(0.48,-0.13){1}{\line(1,0){0.48}}
\multiput(14.45,59.5)(0.48,-0.12){1}{\line(1,0){0.48}}
\multiput(13.96,59.6)(0.49,-0.1){1}{\line(1,0){0.49}}
\multiput(13.47,59.7)(0.49,-0.09){1}{\line(1,0){0.49}}
\multiput(12.98,59.78)(0.49,-0.08){1}{\line(1,0){0.49}}
\multiput(12.49,59.84)(0.49,-0.07){1}{\line(1,0){0.49}}
\multiput(11.99,59.9)(0.5,-0.06){1}{\line(1,0){0.5}}
\multiput(11.49,59.94)(0.5,-0.04){1}{\line(1,0){0.5}}
\multiput(11,59.98)(0.5,-0.03){1}{\line(1,0){0.5}}
\multiput(10.5,59.99)(0.5,-0.02){1}{\line(1,0){0.5}}
\multiput(10,60)(0.5,-0.01){1}{\line(1,0){0.5}}

\linethickness{0.3mm}
\put(10,40){\line(1,0){120}}
\linethickness{0.3mm}
\multiput(100,50)(0.36,-0.12){83}{\line(1,0){0.36}}
\linethickness{0.3mm}
\multiput(80,60)(0.24,-0.12){83}{\line(1,0){0.24}}
\linethickness{0.3mm}
\multiput(60,50)(0.24,0.12){83}{\line(1,0){0.24}}
\linethickness{0.3mm}
\multiput(30,40)(0.36,0.12){83}{\line(1,0){0.36}}
\linethickness{0.3mm}
\put(60,50){\line(1,0){40}}
\linethickness{0.3mm}
\multiput(30,40)(0.84,0.12){83}{\line(1,0){0.84}}
\linethickness{0.3mm}
\multiput(60,50)(0.84,-0.12){83}{\line(1,0){0.84}}
\put(10,45){\makebox(0,0)[cc]{$x$}}

\put(33,36){\makebox(0,0)[cc]{$u$}}

\put(57,54){\makebox(0,0)[cc]{$z$}}

\put(87,60){\makebox(0,0)[cc]{$t$}}

\put(100,54){\makebox(0,0)[cc]{$y$}}

\put(130,36){\makebox(0,0)[cc]{$b$}}

\end{picture}

\noindent 
 En effet quitte à permuter $y$ et $z$ on peut supposer $$d(u,z)+d(y,b)\leq d(u,y)+d(z,b).$$ Alors $(H_{\de}(z,u,y,b))$ implique $d(z,y)\leq d(z,b)+d(u,y)-d(u,b)+\de$ donc 
\begin{gather*}d(u,z)+d(z,y)+d(y,b)\leq 
(d(u,z)+d(z,b))\\ +(d(u,y)+d(y,b))-d(u,b)+\de \leq 
d(u,b)+2P+\de\end{gather*}\begin{gather}\nonumber \text{et\ \ \ }d(u,t)+d(t,b)\leq d(u,z)+d(z,t)+d(t,y)+d(y,b)\\ \label{eq2-dem-4.14}\leq d(u,z)+d(z,y)+P+d(y,b)\leq d(u,b)+3P+\de.\end{gather} On a donc $t\in (3P+\de)\tg(u,b)$ donc $t\in 4P\tg(u,b)$ puisque $P\geq\de$.

Il reste à montrer d). On garde les notations du dessin ci-dessus. Par \eqref{eq2-dem-4.14}  on a $d(x,u)+d(u,z)+d(z,t)+d(t,b)\leq d(x,b)+3P+\de$, d'où $d(x,z)+d(z,t)\leq d(x,t)+3P+\de$. Comme $d(x,z)\geq d_{\max}(x,S_{i+1})-N$, on en déduit 
$d(x,t)\geq d_{\max}(x,S_{i})-N-3P-\de\geq d_{\max}(x,S_{i+1})-4P$ car on suppose $P\geq N+\de$. 
 \cqfd

  Le lemme suivant sera utile pour la démonstration de la proposition~\ref{norme-bien-definie}.

\begin{lem}\label{existence-C}
Il existe une constante $C=C(\de,K,N,Q,P)$ telle que pour tous $x\in X$, $k,m,l_{0},\pp,l_{m}\in \N$, et $a_{1},\pp,a_{p},S_{0},\pp,S_{m}\in \Delta$ vérifiant les conditions i) et ii) de la définition~\ref{defi-Y}, 
  \begin{itemize}
 \item  pour tout $i\in \{0,\dots,m\}$ et $ j\in \{1,\dots ,l_{i}\}$, si on se donne  $d(x,\mathcal Y_{i}^{j})$, le nombre de possibilités pour $\mathcal Y_{i}^{j}$ vérifiant la condition iii) ou iv) de la définition~\ref{defi-Y} (selon que $i<m$ ou $i=m$) est inférieur ou égal à $C$, 
 \item  pour tout $i\in \{0,\dots,m\}$ et $ j\in \{1,\dots ,l_{i}\}$ le nombre de possibilités pour $\mathcal Y_{i}^{j}$ vérifiant la condition iii) ou iv) de la définition~\ref{defi-Y} est inférieur ou égal à $Cs_{i}(Z)$. 
 \item  le cardinal de l'ensemble des $(\mathcal Y_{i}^{j})_{i\in \{0,\dots,m\}, j\in \{1,\dots ,l_{i}\}})$ tels que $$(a_{1},\dots,a_{p},S_{0},...,S_{m},(\mathcal Y_{i}^{j})_{i\in \{0,\dots,m\}, j\in \{1,\dots ,l_{i}\}})$$  appartienne à 
    $Y_{x}^{p,k,m,(l_{0},...,l_{m})}$ est majoré par $C^{\sum_{i=0}^{m}l_{i}} \prod_{i=0}^{m}s_{i}(Z)^{l_{i}}$  
        \end{itemize}
\end{lem}
 \noindent{\bf Démonstration.}
La première assertion résulte immédiatement du c) du 
lemme~\ref{lemme-S0-...Sm}, qui implique que $\mathcal Y_{i}^{j}\subset 4P\tg(x,b)$, du fait que $\mathrm{diam}(\mathcal Y_{i}^{j})\leq P$, 
et du lemme~\ref{cardinal-tranche-geod} appliqué à $(x,b)$ au lieu de $(x,y)$.  
On montre  d'abord la deuxième assertion pour  $i\in \{0,...,m-1\}$. Soit $a\in S_{i}$ et $b\in S_{i+1}$. Pour $j\in \{1,...,l_{i}\}$ on a $$\mathcal Y_{i}^{j}\subset \bigcup _{y\in S_{i},z\in S_{i+1}}P\tg(y,z)\subset (P+4N)\tg(a,b)$$ par le lemme~\ref{xx'yy'zz'}. Comme  $\mathrm{diam}(\mathcal Y_{i}^{j})\leq P$,  le lemme~\ref{cardinal-tranche-geod} appliqué à $(a,b)$ au lieu de $(x,y)$ montre alors la deuxième assertion. 
On montre maintenant la deuxième assertion pour $i=m$. 
Soit  $a\in S_{m}$. On a $$\mathcal Y_{m}^{j}\subset \bigcup _{y\in S_{m}}(2P+\de)\tg(x,y)$$ par le b) du  lemme~\ref{lemme-S0-...Sm}, d'où $\mathcal Y_{m}^{j}\subset (2P+\de+2N)\tg(x,a)$ par le lemme~\ref{xx'yy'zz'}. Comme  $\mathrm{diam}(\mathcal Y_{i}^{j})\leq P$ et $d(x,\mathcal Y_{m}^{j})>k+3P$, le lemme~\ref{cardinal-tranche-geod} appliqué à $(x,a)$ au lieu de $(x,y)$ montre  la deuxième assertion. 
Enfin la troisième assertion résulte facilement de la deuxième. 
 \cqfd

\noindent{\bf Démonstration de la proposition~\ref{norme-bien-definie}.}
Soit $p\in \{1,\pp, p_{\max}\}$ et $f\in \C^{(\Delta_{p})}$. 
Soit $R=\max d_{\max }(x,S)$ où le maximum est pris sur les $S$ tels que $e_{S}$ apparaisse dans $f$ avec un coefficient non nul. 

Soit  $k\geq R$.  
Si $m,l_{0},...,l_{m}\in \N$ et $$(a_{1},\dots,a_{p},S_{0},...,S_{m},(\mathcal Y_{i}^{j})_{i\in \{0,\dots,m\}, j\in \{1,\dots ,l_{i}\}}) \in 
Y_{x}^{p,k,m,(l_{0},...,l_{m})}$$ sont tels que $f(a_{1},...,a_{p})\neq 0$, alors 
$d_{\max }(x,S_{0})\leq k$, donc $m=0$ et $l_{0}=0$ par le  lemme~\ref{m-fini} et  la remarque qui suit  la définition~\ref{defi-Y}. 
De plus  pour $ k\geq R$ la relation d'équivalence sur la partie de 
$Y_{x}^{p,k,0,(0)}$ telle que $e_{S_{0}}$ apparaisse dans $f$ avec un coefficient non nul, est triviale, donc la partie de~(\ref{formule-norme}) correspondant à $ k\geq R$ se réécrit 
$$p! \sum_{S_{0},k\geq R}e^{2s(d(x,S_{0})-k)}|f(S_{0})|^{2}$$ et elle est finie et majorée par $p! (1-e^{-2s})^{-1}\|f\|_{\ell^{2}(\Delta_{p})}^{2}$. 

On fixe  $k\in \{0,\pp,R-1\}$. Il reste donc à montrer que la partie correspondante de~(\ref{formule-norme}) est une somme convergente. 
La somme 
\begin{gather*}\sum _{m\in \N,(l_{0},\dots ,l_{m})\in \N^{m+1}} 
B^{-(m+\sum_{i=0}^{m}l_{i})}\sum_{Z\in \overline Y_{x}^{p,k,m,(l_{0},...,l_{m})}}
\\ 
e^{2s(r_{0}(Z)-k)}\Big(\prod_{i=0}^{m}s_{i}(Z)^{-l_{i}} \Big)
\sharp ((\pi_{x}^{p,k,m,(l_{0},...,l_{m})})^{-1}(Z))^{-\alpha }
\big|\xi_{Z}(f)\big|^{2}
\end{gather*} 
%$$\Big|\sum _{(a_{1},\dots,a_{p},S_{1},...,S_{m},(\mathcal Y_{i}^{j})_{i\in \{0,\dots,m\}, j\in \{1,\dots ,l_{i}\}}) \in 
%(\pi_{x}^{p,k,m,(l_{0},...,l_{m})})^{-1}(Z)} f(a_1,...,a_p)\Big|^{2}$$ 
est majorée par 
\begin{gather}\label{majoration2}e^{2sR}\sum _{m\in \N,(l_{0},\dots ,l_{m})\in \N^{m+1}} 
B^{-(m+\sum_{i=0}^{m}l_{i})}\sum_{Z\in \overline Y_{x}^{p,k,m,(l_{0},...,l_{m})}}
\Big(\prod_{i=0}^{m}s_{i}(Z)^{-l_{i}} \Big)\big|\xi_{Z}(f)\big|^{2}\end{gather} 
%$$ \Big|\sum _{(a_{1},\dots,a_{p},S_{1},...,S_{m},(\mathcal Y_{i}^{j})_{i\in \{0,\dots,m\}, j\in \{1,\dots ,l_{i}\}}) \in 
%(\pi_{x}^{p,k,m,(l_{0},...,l_{m})})^{-1}(Z)} f(a_1,...,a_p)\Big|^{2}$$
Soit $D$ comme dans le lemme~\ref{m-fini}. 
Grâce au lemme~\ref{m-fini}, pour tout $m\in \N$ et pour tout $S_{0}\in \Delta$ vérifiant $S_{0}\subset B(x,R)$,  le nombre de possibilités pour $(S_{1},\pp,S_{m})$ est inférieur ou égal à $\max(1,e^{D(R-k+m)})\leq e^{D(R+m)}$.   Soit $C$  égal à la constante $C$ du lemme~\ref{existence-C} (qui est de la forme $C(\de,K,N,Q,P)$).  
%et   à la constante $C_{3}$ du b) iv) du lemme~\ref{xaS-constante-C}.  
Pour tous $m,(l_{0},\pp,l_{m})$, et pour tout $Z\in \overline Y_{x}^{p,k,m,(l_{0},...,l_{m})}$ l'application  de $(\pi_{x}^{p,k,m,(l_{0},...,l_{m})})^{-1}(Z)$ dans $ \Delta_{p}$ qui à $$(a_{1},\dots,a_{p},S_{0},...,S_{m},(\mathcal Y_{i}^{j})_{i\in \{0,\dots,m\}, j\in \{1,\dots ,l_{i}\}})$$ associe $S_{0}$ 
a des fibres de cardinal $\leq p! e^{D(R+m)}C^{\sum_{i=0}^{m}l_{i}}$, car connaissant $S_{0}$ on a $p!$ possibilités pour $(a_{1},\pp,a_{p})$, au plus 
$e^{D(R+m)}$ possibilités pour $(S_{1},...,S_{m})$, et grâce 
au lemme~\ref{existence-C}, 
%à b) iv) du lemme~\ref{xaS-constante-C}, 
au plus $C^{\sum_{i=0}^{m}l_{i}}$ 
possibilités pour \break $(\mathcal Y_{i}^{j})_{i\in \{0,\dots,m\}, j\in \{1,\dots ,l_{i}\}}$, puisque $Z$ détermine, pour tous $i\in \{0,\dots,m\}$ et $ j\in \{1,\dots ,l_{i}\}$, l'entier $d(x,\mathcal Y_{i}^{j})$. 
%$s$ tel que $\mathcal Y_{i}^{j}\subset \{y\in X, d(y,S_{i})=s\}$. 
Donc dans  (\ref{majoration2}) on a toujours 
%$$\Big|\sum _{(a_{1},\dots,a_{p},S_{1},...,S_{m},(\mathcal Y_{i}^{j})_{i\in \{0,\dots,m\}, j\in \{1,\dots ,l_{i}\}}) \in 
%(\pi_{x}^{p,k,m,(l_{0},...,l_{m})})^{-1}(Z)} f(a_1,...,a_p)\Big|$$
$$\big|\xi_{Z}(f)\big| \leq p! e^{D(R+m)}C^{\sum_{i=0}^{m}l_{i}} \|f\|_{\ell^{1}(\Delta_{p})}. $$
Donc la somme (\ref{majoration2}) est majorée par 
\begin{gather}e^{2sR}\nonumber \sum _{m\in \N,(l_{0},\dots ,l_{m})\in \N^{m+1}} 
B^{-(m+\sum_{i=0}^{m}l_{i})}\sum_{Z\in \overline Y_{x}^{p,k,m,(l_{0},...,l_{m})}}\\ \label{majoration3}
\Big(\prod_{i=0}^{m}s_{i}(Z)^{-l_{i}} \Big)
(p!)^{2} e^{2D(R+m)}C^{2\sum_{i=0}^{m}l_{i}} \|f\|^{2}_{\ell^{1}(\Delta_{p})}
\end{gather} 
Il existe $D'=C(\de,K)$ tel que le nombre de $S_{0}\in \Delta_{p}$ inclus dans $B(x,R)$ soit inférieur ou égal à $e^{D'R}$ pour tout $R\in \N$. 
Par Cauchy-Schwarz, on en déduit $ \|f\|^{2}_{\ell^{1}(\Delta_{p})}\leq e^{D'R} \|f\|^{2}_{\ell^{2}(\Delta_{p})}$. Pour tous $m,(l_{0},\pp,l_{m})$ le cardinal de $Y_{x}^{p,k,m,(l_{0},...,l_{m})}$, et donc a fortiori celui de $\overline Y_{x}^{p,k,m,(l_{0},...,l_{m})}$, sont majorés par 
$$p! e^{D(R+m)+D'R}C^{\sum_{i=0}^{m}l_{i}}\Big(\prod_{i=0}^{m}s_{i}(Z)^{l_{i}} \Big),$$ car on a 
au plus $e^{D'R}$ possibilités pour $S_{0}$,  
$p!$ possibilités pour $(a_{1},\pp,a_{p})$, au plus 
$e^{D(R+m)}$ possibilités pour $(S_{1},...,S_{m})$, et grâce au lemme~\ref{existence-C}, au plus $C^{\sum_{i=0}^{m}l_{i}}\prod_{i=0}^{m}s_{i}(Z)^{l_{i}} $ 
possibilités pour $(\mathcal Y_{i}^{j})_{i\in \{0,\dots,m\}, j\in \{1,\dots ,l_{i}\}}$. 
Donc la somme (\ref{majoration3}) est majorée par 
$$(p!)^{3}e^{2sR}\sum _{m\in \N,(l_{0},\dots ,l_{m})\in \N^{m+1}}  B^{-(m+\sum_{i=0}^{m}l_{i})}
e^{3D(R+m)+2D'R} C^{3\sum_{i=0}^{m}l_{i}}\|f\|_{\ell^{2}(\Delta_{p})}^{2} .$$
On suppose 
$B>C^{3}+e^{3D}$, ce qui est permis par $(H_{B})$. Donc  cette somme converge et est majorée par
\begin{gather*}(p!)^{3}e^{(2s+3D+2D')R}\sum _{m\in \N} \frac{B^{-m}e^{3Dm}}{(1-B^{-1}C^{3})^{m+1}}\|f\|_{\ell^{2}(\Delta_{p})}^{2} =\\ 
(p!)^{3}e^{(2s+3D+2D')R}\frac{1}{1-B^{-1}C^{3}}\frac{1}{1-\frac{B^{-1}e^{3D}}{1-B^{-1}C^{3}}}
\|f\|_{\ell^{2}(\Delta_{p})}^{2}\\ =\frac{(p!)^{3}e^{(2s+3D+2D')R}}{1-B^{-1}(C^{3}+e^{3D})}\|f\|_{\ell^{2}(\Delta_{p})}^{2}.\end{gather*}
Cela termine la démonstration de la proposition~\ref{norme-bien-definie}.\cqfd

Au total,  pour $p\in \{1,\pp, p_{\max}\}$ et $f\in \C^{(\Delta_{p})}$ et en notant $R=\max d_{\max }(x,S)$ où le maximum est pris sur les $S$ tels que $e_{S}$ apparaisse dans $f$ avec un coefficient non nul, on a 
\begin{gather} \nonumber \|f\|_{\H_{x,s}(\Delta_{p})}^{2}\leq p!(1-e^{-2s})^{-1}\|f\|_{\ell^{2}(\Delta_{p})}^{2}+\sum_{k=0}^{R-1}\frac{(p!)^{3}e^{(2s+3D+2D')R}}{1-B^{-1}(C^{3}+e^{3D})}\|f\|_{\ell^{2}(\Delta_{p})}^{2} \\ \label{estimation-norme-Hsx} \leq \Big(p!(1-e^{-2s})^{-1}+\frac{(p!)^{3}Re^{(2s+3D+2D')R}}{1-B^{-1}(C^{3}+e^{3D})}\Big) \|f\|_{\ell^{2}(\Delta_{p})}^{2}.\end{gather}

Le lemme suivant donne au contraire une minoration de la norme
de $\H_{x,s}(\Delta_{p})$. 

\begin{lem}\label{minoration-normeHsx}
Pour tout $p\in \{1,\pp,p_{\max}\}$ et pour tout $f\in \C^{(\Delta_{p})}$ on a 
$$\|f\|_{\H_{x,s}(\Delta_{p})}^{2}\geq  p!(1-e^{-2s})^{-1} \|f\|_{\ell^{2}(\Delta_{p})}^{2}.$$
\end{lem}
 \noindent{\bf Démonstration.}
 Pour tout $S_{0}\in \Delta_{p}$, pour tout $k\geq d(x,S_{0})$ et pour toute énumération $(a_{1},\pp,a_{p})$ des points de $S_{0}$, 
$Y^{p,k,0,(0)}_{x}$ contient $(a_{1},\pp,a_{p},S_{0})$ et comme $M\geq N$, le singleton 
$\{ (a_{1},\pp,a_{p},S_{0}) \}$ est une classe d'équivalence dans $\overline Y^{p,k,0,(0)}_{x}$. 
On a donc $$\|f\|_{\H_{x,s}(\Delta_{p})}^{2}\geq p! \sum_{S_{0}\in \Delta_{p}}
\sum_{k\geq d(x,S_{0})} e^{2s(d(x,S_{0})-k)}|f(S_{0})|^{2}= p!(1-e^{-2s})^{-1} \|f\|_{\ell^{2}(\Delta_{p})}^{2}.$$ \cqfd

\subsection{Autres propriétés de la norme}

Les propriétés que nous allons établir dans ce sous-paragraphe sont des préliminaires indispensables pour les sous-paragraphes suivants. De fa\c con un peu imprécise nous allons montrer qu'il existe une constante $C$ de la forme $C(\de,K,N,Q,P,M)$ telle que pour $p\in \{1,...,p_{\max}\}$, $k,m,l_{0},...,l_{m}\in \N$ et  
$$(a_{1},\dots,a_{p},S_{0},...,S_{m},(\mathcal Y_{i}^{j})_{i\in \{0,\dots,m\}, j\in \{1,\dots ,l_{i}\}})\in Y_{x}^{p,k,m,(l_{0},...,l_{m})},$$ pour connaître les distances entre les points de 
\begin{gather*}
\bigcup _{ i\in \{0,\dots ,m\}}
B(S_{i},M)\cup   \bigcup _{i\in \{0,\dots,m\}, j\in \{1,\dots ,l_{i}\}} B(\mathcal Y_{i}^{j},  M)  \cup B(x,k+2M)\end{gather*}
il suffit de connaître certaines de ces distances, de telle sorte que pour chaque point de cet ensemble n'appartenant pas à $B(x,k+2M)$, le nombre de distances à connaître depuis de ce point soit inférieur ou égal à $C$. 
La raison est que cet ensemble est une réunion de boules de grands rayons dont les centres sont à peu près alignés le long d'une géodésique (grâce au lemme~\ref{lemme-S0-...Sm}) et qu'il suffit donc de connaître les distances entre les points de boules voisines.

\begin{lem}\label{i<j<k-distances}
Soient $c,d\in X$, $I$ un ensemble fini, et pour $i\in I$,  $\alpha_{i}\in \N$, $\rho_{i}\in \N^{*}$ et  $w_{i}\in \alpha_{i}\tg(c,d)$. On suppose que pour tout $i\in I$, $\rho_{i}\geq 4\de+\frac{\alpha_{i}}{2}$. 

\noindent a) 
Pour connaître les distances entre les points de $\bigcup_{i\in I}B(w_{i},\rho_{i})$ il suffit de connaître les distances entre les points de 
$B(w_{i},\rho_{i})$ et $B(w_{j},\rho_{j})$ pour tous les couples $(i,j)\in   \Lambda$ où $\Lambda \subset I^{2}$ est l'ensemble des couples    $(i,j)$ tels qu'il n'existe pas de $k\in I$ vérifiant \begin{gather}\label{cond-lem-i<j<k-distances}d(c,w_{j})-d(c,w_{k})\geq 2\rho_{j}+\frac{\alpha_{j}}{2}
\text{\  et\ } 
 d(c,w_{k})-d(c,w_{i})\geq 2\rho_{i}+\frac{\alpha_{i}}{2}+\frac{\alpha_{k}}{2} . \end{gather}
  
  \noindent b)  Soit $R\in \N^{*}$ tel que $\alpha_{i}\leq R$ et $\rho_{i}\leq R$ pour tout $i\in I$. Soit $i\in I$. Alors l'ensemble des $d(c,w_{j})$ pour $j\in 
  I$ tel que $(i,j)\in \Lambda$ ou $(j,i)\in \Lambda$ est inclus dans la réunion de $[d(c,w_{i})-3R, d(c,w_{i})+3R]$ et de deux intervalles  de longueur $\leq 3R$. 
  
  \noindent c) Soit $R$ comme dans b). Il existe une constante $C=C(\de,K,R)$ (indépendante de $\sharp I$ en particulier) telle que pour toute partie $J\subset  
  I$, 
  connaissant les distances entre les points de $\bigcup_{i\in I\setminus J }B(w_{i},\rho_{i})$, 
    les distances entre les points de $\bigcup_{i\in  J }B(w_{i},\rho_{i})$ et ceux de $\bigcup_{i\in I\setminus J }B(w_{i},\rho_{i})$ soient déterminées par la donnée des distances entre les points de $\bigcup_{i\in  J }B(w_{i},\rho_{i})$ et les points d'une partie  de $\bigcup_{i\in I\setminus J }B(w_{i},\rho_{i})$
    \begin{itemize}
  \item dont le cardinal est borné par $C(\sharp J)$, 
  \item qui est déterminée par la connaissance des distances entre les points de  $$\bigcup_{i\in I\setminus J }B(w_{i},\rho_{i})\cup \{c\}$$ et par les entiers $d(c,w_{i})$ pour $i\in J$. 
  \end{itemize}
\end{lem}
\noindent{\bf Remarque. } Dans toutes les situations où on appliquera le c) de ce lemme, $\sharp J$ sera majoré par une constante de la forme $C(\de,K,N,Q,P,M)$. 

 \noindent{\bf Démonstration.}
On commence par montrer a). 
On suppose $I=\{1,...,r\}$ et $w_{1},...,w_{r}$ ordonnés de telle sorte que $d(c,w_{1})\leq \pp\leq d(c,w_{r})$. Soit $i<k<j$ des entiers vérifiant (\ref{cond-lem-i<j<k-distances}). 
Soient $y\in B(w_{i},\rho_{i})$ et $z\in B(w_{j},\rho_{j})$. On a \begin{gather*}d(c,y)\leq d(c,w_{i})+\rho_{i}\leq d(c,w_{k})-\rho_{i}-\frac{\alpha_{i}}{2}-\frac{\alpha_{k}}{2}\\ \text{\  et\  }d(c,z)\geq d(c,w_{j})-\rho_{j}\geq d(c,w_{k})+\rho_{j}+\frac{\alpha_{j}}{2}.\end{gather*} Donc il existe $v\in \geod(y,z)$ tel que $d(c,v)=d(c,w_{k})-E(\frac{\alpha_{k}}{2})$. On choisit un tel $v$.

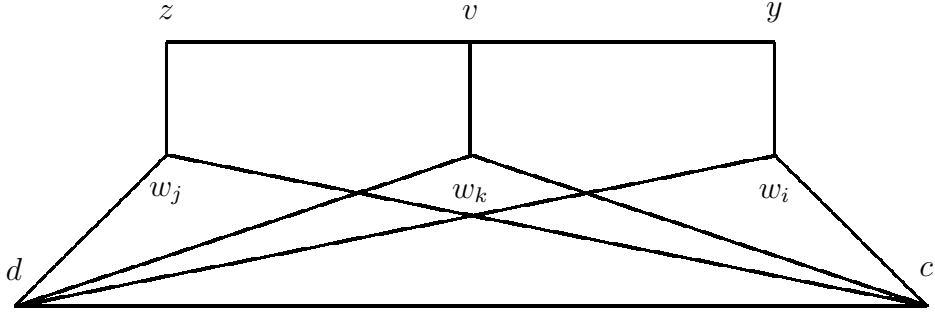
\ifx\JPicScale\undefined\def\JPicScale{1}\fi
\unitlength \JPicScale mm
\begin{picture}(130,52)(10,25)
\linethickness{0.3mm}
\put(10,30){\line(1,0){120}}
\linethickness{0.3mm}
\multiput(10,30)(0.12,0.12){167}{\line(1,0){0.12}}
\linethickness{0.3mm}
\multiput(30,50)(0.6,-0.12){167}{\line(1,0){0.6}}
\linethickness{0.3mm}
\multiput(110,50)(0.12,-0.12){167}{\line(1,0){0.12}}
\linethickness{0.3mm}
\multiput(10,30)(0.6,0.12){167}{\line(1,0){0.6}}
\linethickness{0.3mm}
\multiput(10,30)(0.36,0.12){167}{\line(1,0){0.36}}
\linethickness{0.3mm}
\multiput(70,50)(0.36,-0.12){167}{\line(1,0){0.36}}
\linethickness{0.3mm}
\put(30,50){\line(0,1){15}}
\linethickness{0.3mm}
\put(110,50){\line(0,1){15}}
\linethickness{0.3mm}
\put(30,65){\line(1,0){80}}
\linethickness{0.3mm}
\put(70,50){\line(0,1){15}}
\put(10,35){\makebox(0,0)[cc]{$d$}}

\put(130,35){\makebox(0,0)[cc]{$c$}}

\put(30,45){\makebox(0,0)[cc]{$w_j$}}

\put(70,45){\makebox(0,0)[cc]{$w_k$}}

\put(110,45){\makebox(0,0)[cc]{$w_i$}}

\put(30,69){\makebox(0,0)[cc]{$z$}}

\put(70,69){\makebox(0,0)[cc]{$v$}}

\put(110,69){\makebox(0,0)[cc]{$y$}}

\end{picture}

\noindent 
On a $$d(y,v)\geq d(c,v)-d(c,y)\geq \rho_{i}+\alpha_{i}/2\text{\  et\ } 
d(z,v)\geq d(c,z)-d(c,v)\geq \rho_{j}+\alpha_{j}/2.$$ D'après le lemme~\ref{xx'yy'zz'}, on a $y\in (\alpha_{i}+2\rho_{i})\tg(c,d)$ et $z\in (\alpha_{j}+2\rho_{j})\tg(c,d)$. Le lemme~\ref{lemmetubes} implique alors $v\in 3\de\tg(c,d)$. Comme $w_{k}\in \alpha_{k}\tg(c,d)$ et $d(c,v)=d(c,w_{k})-E(\frac{\alpha_{k}}{2})$, $(H_{\de}(v,c,w_{k},d))$ implique
\begin{gather*}d(v,w_{k})\leq
\max(d(c,v)+d(w_{k},d)-d(c,d),d(c,w_{k})+d(v,d)-d(c,d))+\de
\leq\\  \max(d(c,w_{k})-E(\frac{\alpha_{k}}{2})+d(w_{k},d)-d(c,d),d(c,v)+E(\frac{\alpha_{k}}{2})+d(v,d)-d(c,d))\\ +\de
 \leq 
\max(\alpha_{k}-E(\frac{\alpha_{k}}{2}),E(\frac{\alpha_{k}}{2})+3\de)+\de\leq 4\de+\frac{\alpha_{k}}{2}\leq 
\rho_{k},\end{gather*} donc $v\in B(w_{k},\rho_{k})\cap \geod(y,z)$. On a donc montré que pour $i<k<j$ vérifiant (\ref{cond-lem-i<j<k-distances}) toute géodésique entre un point de $B(w_{i},\rho_{i})$ et un point de $B(w_{j},\rho_{j})$ intersecte 
$ B(w_{k},\rho_{k})$, ce qui implique que la connaissance des distances entre les points de $ B(w_{k},\rho_{k})$ et ceux de $B(w_{i},\rho_{i})\cup B(w_{j},\rho_{j})$ permet de déterminer les distances entre les points de $B(w_{i},\rho_{i})$ et ceux de $B(w_{j},\rho_{j})$. Ceci termine la preuve du a). 
Le b) résulte facilement du a). Le cas particulier de c) où $\sharp J=1$ résulte du b) et du lemme~\ref{cardinal-tranche-geod}.  Pour montrer c) dans le cas général on se ramène au cas particulier déjà démontré  de la fa\c con suivante. 
 Pour tout $j\in J$ on  
 note $I'=I\setminus J \cup \{j\}$ et $J'=\{j\}$ et on 
 applique le cas particulier de c) déjà démontré avec $I'$  et $J'$ au lieu $I$ et $J$. \cqfd

\begin{lem}\label{distances-Bxk-C}
Soient $x,b\in X,l\in \N, \alpha\in \N$ et $u$ un point de $B(x,l)$ à distance minimale de $b$.  Pour tout $z\in \alpha\tg(b,u)$, les distances entre $z$ et les points de $B(x,l)$ sont déterminées par les distances entre $z$ et les points de $B(u,\alpha+4\de)\cap B(x,l)$. Plus précisément  pour tout $\tilde x\in B(x,l)$ et pour tout $z\in \alpha\tg(b,u)$, $\geod(\tilde x,z)$ intersecte $B(u,\alpha+4\de)\cap B(x,l)$.   
\end{lem}

 \noindent{\bf Démonstration.}
Le lemme est évident si $b\in B(x,l)$. On suppose donc $b\not\in B(x,l)$. 
Soit $\tilde x\in B(x,l)$ et $z\in \alpha\tg(b,u)$. On veut montrer qu'il existe $w\in \geod(\tilde x,z)\cap B(u,\alpha+4\de)\cap B(x,l)$. Si $d(u,\tilde x)\leq \alpha+4\de$ on prend $w=\tilde x$. On suppose donc $d(u,\tilde x)> \alpha+4\de$. 
D'après le lemme~\ref{Bxk-x-t-z-y} appliqué à $(x,b,u,\tilde x)$ au lieu de $(x,y,t,z)$ et $l$ au lieu de $k$,  on a $u\in \de\tg(\tilde x,b)$.  Donc  $ d(b,\tilde x)\geq d(b,u)+\alpha+3\de$  et comme de plus $d(b,z)\leq d(b,u)+\alpha$, il existe $w\in \geod(z,\tilde x)$ tel que $d(b,w)=d(b,u)+\alpha+2\de$. On choisit un tel $w$. 

\ifx\JPicScale\undefined\def\JPicScale{1}\fi
\unitlength \JPicScale mm
\begin{picture}(130,50)(10,15)
\linethickness{0.3mm}
\put(30,40){\line(1,0){100}}
\linethickness{0.3mm}
\put(50,39.75){\line(0,1){0.5}}
\multiput(49.99,40.75)(0.01,-0.5){1}{\line(0,-1){0.5}}
\multiput(49.96,41.25)(0.02,-0.5){1}{\line(0,-1){0.5}}
\multiput(49.92,41.74)(0.04,-0.5){1}{\line(0,-1){0.5}}
\multiput(49.87,42.24)(0.05,-0.5){1}{\line(0,-1){0.5}}
\multiput(49.81,42.73)(0.06,-0.49){1}{\line(0,-1){0.49}}
\multiput(49.74,43.23)(0.07,-0.49){1}{\line(0,-1){0.49}}
\multiput(49.65,43.72)(0.09,-0.49){1}{\line(0,-1){0.49}}
\multiput(49.55,44.21)(0.1,-0.49){1}{\line(0,-1){0.49}}
\multiput(49.44,44.69)(0.11,-0.49){1}{\line(0,-1){0.49}}
\multiput(49.32,45.18)(0.12,-0.48){1}{\line(0,-1){0.48}}
\multiput(49.18,45.66)(0.14,-0.48){1}{\line(0,-1){0.48}}
\multiput(49.04,46.13)(0.15,-0.48){1}{\line(0,-1){0.48}}
\multiput(48.88,46.61)(0.16,-0.47){1}{\line(0,-1){0.47}}
\multiput(48.71,47.07)(0.17,-0.47){1}{\line(0,-1){0.47}}
\multiput(48.52,47.54)(0.09,-0.23){2}{\line(0,-1){0.23}}
\multiput(48.33,48)(0.1,-0.23){2}{\line(0,-1){0.23}}
\multiput(48.13,48.45)(0.1,-0.23){2}{\line(0,-1){0.23}}
\multiput(47.91,48.9)(0.11,-0.22){2}{\line(0,-1){0.22}}
\multiput(47.68,49.35)(0.11,-0.22){2}{\line(0,-1){0.22}}
\multiput(47.44,49.78)(0.12,-0.22){2}{\line(0,-1){0.22}}
\multiput(47.19,50.22)(0.12,-0.22){2}{\line(0,-1){0.22}}
\multiput(46.93,50.64)(0.13,-0.21){2}{\line(0,-1){0.21}}
\multiput(46.66,51.06)(0.14,-0.21){2}{\line(0,-1){0.21}}
\multiput(46.38,51.47)(0.14,-0.21){2}{\line(0,-1){0.21}}
\multiput(46.09,51.88)(0.15,-0.2){2}{\line(0,-1){0.2}}
\multiput(45.79,52.27)(0.1,-0.13){3}{\line(0,-1){0.13}}
\multiput(45.48,52.66)(0.1,-0.13){3}{\line(0,-1){0.13}}
\multiput(45.16,53.05)(0.11,-0.13){3}{\line(0,-1){0.13}}
\multiput(44.83,53.42)(0.11,-0.12){3}{\line(0,-1){0.12}}
\multiput(44.49,53.79)(0.11,-0.12){3}{\line(0,-1){0.12}}
\multiput(44.14,54.14)(0.12,-0.12){3}{\line(0,-1){0.12}}
\multiput(43.79,54.49)(0.12,-0.12){3}{\line(1,0){0.12}}
\multiput(43.42,54.83)(0.12,-0.11){3}{\line(1,0){0.12}}
\multiput(43.05,55.16)(0.12,-0.11){3}{\line(1,0){0.12}}
\multiput(42.66,55.48)(0.13,-0.11){3}{\line(1,0){0.13}}
\multiput(42.27,55.79)(0.13,-0.1){3}{\line(1,0){0.13}}
\multiput(41.88,56.09)(0.13,-0.1){3}{\line(1,0){0.13}}
\multiput(41.47,56.38)(0.2,-0.15){2}{\line(1,0){0.2}}
\multiput(41.06,56.66)(0.21,-0.14){2}{\line(1,0){0.21}}
\multiput(40.64,56.93)(0.21,-0.14){2}{\line(1,0){0.21}}
\multiput(40.22,57.19)(0.21,-0.13){2}{\line(1,0){0.21}}
\multiput(39.78,57.44)(0.22,-0.12){2}{\line(1,0){0.22}}
\multiput(39.35,57.68)(0.22,-0.12){2}{\line(1,0){0.22}}
\multiput(38.9,57.91)(0.22,-0.11){2}{\line(1,0){0.22}}
\multiput(38.45,58.13)(0.22,-0.11){2}{\line(1,0){0.22}}
\multiput(38,58.33)(0.23,-0.1){2}{\line(1,0){0.23}}
\multiput(37.54,58.52)(0.23,-0.1){2}{\line(1,0){0.23}}
\multiput(37.07,58.71)(0.23,-0.09){2}{\line(1,0){0.23}}
\multiput(36.61,58.88)(0.47,-0.17){1}{\line(1,0){0.47}}
\multiput(36.13,59.04)(0.47,-0.16){1}{\line(1,0){0.47}}
\multiput(35.66,59.18)(0.48,-0.15){1}{\line(1,0){0.48}}
\multiput(35.18,59.32)(0.48,-0.14){1}{\line(1,0){0.48}}
\multiput(34.69,59.44)(0.48,-0.12){1}{\line(1,0){0.48}}
\multiput(34.21,59.55)(0.49,-0.11){1}{\line(1,0){0.49}}
\multiput(33.72,59.65)(0.49,-0.1){1}{\line(1,0){0.49}}
\multiput(33.23,59.74)(0.49,-0.09){1}{\line(1,0){0.49}}
\multiput(32.73,59.81)(0.49,-0.07){1}{\line(1,0){0.49}}
\multiput(32.24,59.87)(0.49,-0.06){1}{\line(1,0){0.49}}
\multiput(31.74,59.92)(0.5,-0.05){1}{\line(1,0){0.5}}
\multiput(31.25,59.96)(0.5,-0.04){1}{\line(1,0){0.5}}
\multiput(30.75,59.99)(0.5,-0.02){1}{\line(1,0){0.5}}
\multiput(30.25,60)(0.5,-0.01){1}{\line(1,0){0.5}}
\put(29.75,60){\line(1,0){0.5}}
\multiput(29.25,59.99)(0.5,0.01){1}{\line(1,0){0.5}}
\multiput(28.75,59.96)(0.5,0.02){1}{\line(1,0){0.5}}
\multiput(28.26,59.92)(0.5,0.04){1}{\line(1,0){0.5}}
\multiput(27.76,59.87)(0.5,0.05){1}{\line(1,0){0.5}}
\multiput(27.27,59.81)(0.49,0.06){1}{\line(1,0){0.49}}
\multiput(26.77,59.74)(0.49,0.07){1}{\line(1,0){0.49}}
\multiput(26.28,59.65)(0.49,0.09){1}{\line(1,0){0.49}}
\multiput(25.79,59.55)(0.49,0.1){1}{\line(1,0){0.49}}
\multiput(25.31,59.44)(0.49,0.11){1}{\line(1,0){0.49}}
\multiput(24.82,59.32)(0.48,0.12){1}{\line(1,0){0.48}}
\multiput(24.34,59.18)(0.48,0.14){1}{\line(1,0){0.48}}
\multiput(23.87,59.04)(0.48,0.15){1}{\line(1,0){0.48}}
\multiput(23.39,58.88)(0.47,0.16){1}{\line(1,0){0.47}}
\multiput(22.93,58.71)(0.47,0.17){1}{\line(1,0){0.47}}
\multiput(22.46,58.52)(0.23,0.09){2}{\line(1,0){0.23}}
\multiput(22,58.33)(0.23,0.1){2}{\line(1,0){0.23}}
\multiput(21.55,58.13)(0.23,0.1){2}{\line(1,0){0.23}}
\multiput(21.1,57.91)(0.22,0.11){2}{\line(1,0){0.22}}
\multiput(20.65,57.68)(0.22,0.11){2}{\line(1,0){0.22}}
\multiput(20.22,57.44)(0.22,0.12){2}{\line(1,0){0.22}}
\multiput(19.78,57.19)(0.22,0.12){2}{\line(1,0){0.22}}
\multiput(19.36,56.93)(0.21,0.13){2}{\line(1,0){0.21}}
\multiput(18.94,56.66)(0.21,0.14){2}{\line(1,0){0.21}}
\multiput(18.53,56.38)(0.21,0.14){2}{\line(1,0){0.21}}
\multiput(18.12,56.09)(0.2,0.15){2}{\line(1,0){0.2}}
\multiput(17.73,55.79)(0.13,0.1){3}{\line(1,0){0.13}}
\multiput(17.34,55.48)(0.13,0.1){3}{\line(1,0){0.13}}
\multiput(16.95,55.16)(0.13,0.11){3}{\line(1,0){0.13}}
\multiput(16.58,54.83)(0.12,0.11){3}{\line(1,0){0.12}}
\multiput(16.21,54.49)(0.12,0.11){3}{\line(1,0){0.12}}
\multiput(15.86,54.14)(0.12,0.12){3}{\line(1,0){0.12}}
\multiput(15.51,53.79)(0.12,0.12){3}{\line(0,1){0.12}}
\multiput(15.17,53.42)(0.11,0.12){3}{\line(0,1){0.12}}
\multiput(14.84,53.05)(0.11,0.12){3}{\line(0,1){0.12}}
\multiput(14.52,52.66)(0.11,0.13){3}{\line(0,1){0.13}}
\multiput(14.21,52.27)(0.1,0.13){3}{\line(0,1){0.13}}
\multiput(13.91,51.88)(0.1,0.13){3}{\line(0,1){0.13}}
\multiput(13.62,51.47)(0.15,0.2){2}{\line(0,1){0.2}}
\multiput(13.34,51.06)(0.14,0.21){2}{\line(0,1){0.21}}
\multiput(13.07,50.64)(0.14,0.21){2}{\line(0,1){0.21}}
\multiput(12.81,50.22)(0.13,0.21){2}{\line(0,1){0.21}}
\multiput(12.56,49.78)(0.12,0.22){2}{\line(0,1){0.22}}
\multiput(12.32,49.35)(0.12,0.22){2}{\line(0,1){0.22}}
\multiput(12.09,48.9)(0.11,0.22){2}{\line(0,1){0.22}}
\multiput(11.87,48.45)(0.11,0.22){2}{\line(0,1){0.22}}
\multiput(11.67,48)(0.1,0.23){2}{\line(0,1){0.23}}
\multiput(11.48,47.54)(0.1,0.23){2}{\line(0,1){0.23}}
\multiput(11.29,47.07)(0.09,0.23){2}{\line(0,1){0.23}}
\multiput(11.12,46.61)(0.17,0.47){1}{\line(0,1){0.47}}
\multiput(10.96,46.13)(0.16,0.47){1}{\line(0,1){0.47}}
\multiput(10.82,45.66)(0.15,0.48){1}{\line(0,1){0.48}}
\multiput(10.68,45.18)(0.14,0.48){1}{\line(0,1){0.48}}
\multiput(10.56,44.69)(0.12,0.48){1}{\line(0,1){0.48}}
\multiput(10.45,44.21)(0.11,0.49){1}{\line(0,1){0.49}}
\multiput(10.35,43.72)(0.1,0.49){1}{\line(0,1){0.49}}
\multiput(10.26,43.23)(0.09,0.49){1}{\line(0,1){0.49}}
\multiput(10.19,42.73)(0.07,0.49){1}{\line(0,1){0.49}}
\multiput(10.13,42.24)(0.06,0.49){1}{\line(0,1){0.49}}
\multiput(10.08,41.74)(0.05,0.5){1}{\line(0,1){0.5}}
\multiput(10.04,41.25)(0.04,0.5){1}{\line(0,1){0.5}}
\multiput(10.01,40.75)(0.02,0.5){1}{\line(0,1){0.5}}
\multiput(10,40.25)(0.01,0.5){1}{\line(0,1){0.5}}
\put(10,39.75){\line(0,1){0.5}}
\multiput(10,39.75)(0.01,-0.5){1}{\line(0,-1){0.5}}
\multiput(10.01,39.25)(0.02,-0.5){1}{\line(0,-1){0.5}}
\multiput(10.04,38.75)(0.04,-0.5){1}{\line(0,-1){0.5}}
\multiput(10.08,38.26)(0.05,-0.5){1}{\line(0,-1){0.5}}
\multiput(10.13,37.76)(0.06,-0.49){1}{\line(0,-1){0.49}}
\multiput(10.19,37.27)(0.07,-0.49){1}{\line(0,-1){0.49}}
\multiput(10.26,36.77)(0.09,-0.49){1}{\line(0,-1){0.49}}
\multiput(10.35,36.28)(0.1,-0.49){1}{\line(0,-1){0.49}}
\multiput(10.45,35.79)(0.11,-0.49){1}{\line(0,-1){0.49}}
\multiput(10.56,35.31)(0.12,-0.48){1}{\line(0,-1){0.48}}
\multiput(10.68,34.82)(0.14,-0.48){1}{\line(0,-1){0.48}}
\multiput(10.82,34.34)(0.15,-0.48){1}{\line(0,-1){0.48}}
\multiput(10.96,33.87)(0.16,-0.47){1}{\line(0,-1){0.47}}
\multiput(11.12,33.39)(0.17,-0.47){1}{\line(0,-1){0.47}}
\multiput(11.29,32.93)(0.09,-0.23){2}{\line(0,-1){0.23}}
\multiput(11.48,32.46)(0.1,-0.23){2}{\line(0,-1){0.23}}
\multiput(11.67,32)(0.1,-0.23){2}{\line(0,-1){0.23}}
\multiput(11.87,31.55)(0.11,-0.22){2}{\line(0,-1){0.22}}
\multiput(12.09,31.1)(0.11,-0.22){2}{\line(0,-1){0.22}}
\multiput(12.32,30.65)(0.12,-0.22){2}{\line(0,-1){0.22}}
\multiput(12.56,30.22)(0.12,-0.22){2}{\line(0,-1){0.22}}
\multiput(12.81,29.78)(0.13,-0.21){2}{\line(0,-1){0.21}}
\multiput(13.07,29.36)(0.14,-0.21){2}{\line(0,-1){0.21}}
\multiput(13.34,28.94)(0.14,-0.21){2}{\line(0,-1){0.21}}
\multiput(13.62,28.53)(0.15,-0.2){2}{\line(0,-1){0.2}}
\multiput(13.91,28.12)(0.1,-0.13){3}{\line(0,-1){0.13}}
\multiput(14.21,27.73)(0.1,-0.13){3}{\line(0,-1){0.13}}
\multiput(14.52,27.34)(0.11,-0.13){3}{\line(0,-1){0.13}}
\multiput(14.84,26.95)(0.11,-0.12){3}{\line(0,-1){0.12}}
\multiput(15.17,26.58)(0.11,-0.12){3}{\line(0,-1){0.12}}
\multiput(15.51,26.21)(0.12,-0.12){3}{\line(0,-1){0.12}}
\multiput(15.86,25.86)(0.12,-0.12){3}{\line(1,0){0.12}}
\multiput(16.21,25.51)(0.12,-0.11){3}{\line(1,0){0.12}}
\multiput(16.58,25.17)(0.12,-0.11){3}{\line(1,0){0.12}}
\multiput(16.95,24.84)(0.13,-0.11){3}{\line(1,0){0.13}}
\multiput(17.34,24.52)(0.13,-0.1){3}{\line(1,0){0.13}}
\multiput(17.73,24.21)(0.13,-0.1){3}{\line(1,0){0.13}}
\multiput(18.12,23.91)(0.2,-0.15){2}{\line(1,0){0.2}}
\multiput(18.53,23.62)(0.21,-0.14){2}{\line(1,0){0.21}}
\multiput(18.94,23.34)(0.21,-0.14){2}{\line(1,0){0.21}}
\multiput(19.36,23.07)(0.21,-0.13){2}{\line(1,0){0.21}}
\multiput(19.78,22.81)(0.22,-0.12){2}{\line(1,0){0.22}}
\multiput(20.22,22.56)(0.22,-0.12){2}{\line(1,0){0.22}}
\multiput(20.65,22.32)(0.22,-0.11){2}{\line(1,0){0.22}}
\multiput(21.1,22.09)(0.22,-0.11){2}{\line(1,0){0.22}}
\multiput(21.55,21.87)(0.23,-0.1){2}{\line(1,0){0.23}}
\multiput(22,21.67)(0.23,-0.1){2}{\line(1,0){0.23}}
\multiput(22.46,21.48)(0.23,-0.09){2}{\line(1,0){0.23}}
\multiput(22.93,21.29)(0.47,-0.17){1}{\line(1,0){0.47}}
\multiput(23.39,21.12)(0.47,-0.16){1}{\line(1,0){0.47}}
\multiput(23.87,20.96)(0.48,-0.15){1}{\line(1,0){0.48}}
\multiput(24.34,20.82)(0.48,-0.14){1}{\line(1,0){0.48}}
\multiput(24.82,20.68)(0.48,-0.12){1}{\line(1,0){0.48}}
\multiput(25.31,20.56)(0.49,-0.11){1}{\line(1,0){0.49}}
\multiput(25.79,20.45)(0.49,-0.1){1}{\line(1,0){0.49}}
\multiput(26.28,20.35)(0.49,-0.09){1}{\line(1,0){0.49}}
\multiput(26.77,20.26)(0.49,-0.07){1}{\line(1,0){0.49}}
\multiput(27.27,20.19)(0.49,-0.06){1}{\line(1,0){0.49}}
\multiput(27.76,20.13)(0.5,-0.05){1}{\line(1,0){0.5}}
\multiput(28.26,20.08)(0.5,-0.04){1}{\line(1,0){0.5}}
\multiput(28.75,20.04)(0.5,-0.02){1}{\line(1,0){0.5}}
\multiput(29.25,20.01)(0.5,-0.01){1}{\line(1,0){0.5}}
\put(29.75,20){\line(1,0){0.5}}
\multiput(30.25,20)(0.5,0.01){1}{\line(1,0){0.5}}
\multiput(30.75,20.01)(0.5,0.02){1}{\line(1,0){0.5}}
\multiput(31.25,20.04)(0.5,0.04){1}{\line(1,0){0.5}}
\multiput(31.74,20.08)(0.5,0.05){1}{\line(1,0){0.5}}
\multiput(32.24,20.13)(0.49,0.06){1}{\line(1,0){0.49}}
\multiput(32.73,20.19)(0.49,0.07){1}{\line(1,0){0.49}}
\multiput(33.23,20.26)(0.49,0.09){1}{\line(1,0){0.49}}
\multiput(33.72,20.35)(0.49,0.1){1}{\line(1,0){0.49}}
\multiput(34.21,20.45)(0.49,0.11){1}{\line(1,0){0.49}}
\multiput(34.69,20.56)(0.48,0.12){1}{\line(1,0){0.48}}
\multiput(35.18,20.68)(0.48,0.14){1}{\line(1,0){0.48}}
\multiput(35.66,20.82)(0.48,0.15){1}{\line(1,0){0.48}}
\multiput(36.13,20.96)(0.47,0.16){1}{\line(1,0){0.47}}
\multiput(36.61,21.12)(0.47,0.17){1}{\line(1,0){0.47}}
\multiput(37.07,21.29)(0.23,0.09){2}{\line(1,0){0.23}}
\multiput(37.54,21.48)(0.23,0.1){2}{\line(1,0){0.23}}
\multiput(38,21.67)(0.23,0.1){2}{\line(1,0){0.23}}
\multiput(38.45,21.87)(0.22,0.11){2}{\line(1,0){0.22}}
\multiput(38.9,22.09)(0.22,0.11){2}{\line(1,0){0.22}}
\multiput(39.35,22.32)(0.22,0.12){2}{\line(1,0){0.22}}
\multiput(39.78,22.56)(0.22,0.12){2}{\line(1,0){0.22}}
\multiput(40.22,22.81)(0.21,0.13){2}{\line(1,0){0.21}}
\multiput(40.64,23.07)(0.21,0.14){2}{\line(1,0){0.21}}
\multiput(41.06,23.34)(0.21,0.14){2}{\line(1,0){0.21}}
\multiput(41.47,23.62)(0.2,0.15){2}{\line(1,0){0.2}}
\multiput(41.88,23.91)(0.13,0.1){3}{\line(1,0){0.13}}
\multiput(42.27,24.21)(0.13,0.1){3}{\line(1,0){0.13}}
\multiput(42.66,24.52)(0.13,0.11){3}{\line(1,0){0.13}}
\multiput(43.05,24.84)(0.12,0.11){3}{\line(1,0){0.12}}
\multiput(43.42,25.17)(0.12,0.11){3}{\line(1,0){0.12}}
\multiput(43.79,25.51)(0.12,0.12){3}{\line(1,0){0.12}}
\multiput(44.14,25.86)(0.12,0.12){3}{\line(0,1){0.12}}
\multiput(44.49,26.21)(0.11,0.12){3}{\line(0,1){0.12}}
\multiput(44.83,26.58)(0.11,0.12){3}{\line(0,1){0.12}}
\multiput(45.16,26.95)(0.11,0.13){3}{\line(0,1){0.13}}
\multiput(45.48,27.34)(0.1,0.13){3}{\line(0,1){0.13}}
\multiput(45.79,27.73)(0.1,0.13){3}{\line(0,1){0.13}}
\multiput(46.09,28.12)(0.15,0.2){2}{\line(0,1){0.2}}
\multiput(46.38,28.53)(0.14,0.21){2}{\line(0,1){0.21}}
\multiput(46.66,28.94)(0.14,0.21){2}{\line(0,1){0.21}}
\multiput(46.93,29.36)(0.13,0.21){2}{\line(0,1){0.21}}
\multiput(47.19,29.78)(0.12,0.22){2}{\line(0,1){0.22}}
\multiput(47.44,30.22)(0.12,0.22){2}{\line(0,1){0.22}}
\multiput(47.68,30.65)(0.11,0.22){2}{\line(0,1){0.22}}
\multiput(47.91,31.1)(0.11,0.22){2}{\line(0,1){0.22}}
\multiput(48.13,31.55)(0.1,0.23){2}{\line(0,1){0.23}}
\multiput(48.33,32)(0.1,0.23){2}{\line(0,1){0.23}}
\multiput(48.52,32.46)(0.09,0.23){2}{\line(0,1){0.23}}
\multiput(48.71,32.93)(0.17,0.47){1}{\line(0,1){0.47}}
\multiput(48.88,33.39)(0.16,0.47){1}{\line(0,1){0.47}}
\multiput(49.04,33.87)(0.15,0.48){1}{\line(0,1){0.48}}
\multiput(49.18,34.34)(0.14,0.48){1}{\line(0,1){0.48}}
\multiput(49.32,34.82)(0.12,0.48){1}{\line(0,1){0.48}}
\multiput(49.44,35.31)(0.11,0.49){1}{\line(0,1){0.49}}
\multiput(49.55,35.79)(0.1,0.49){1}{\line(0,1){0.49}}
\multiput(49.65,36.28)(0.09,0.49){1}{\line(0,1){0.49}}
\multiput(49.74,36.77)(0.07,0.49){1}{\line(0,1){0.49}}
\multiput(49.81,37.27)(0.06,0.49){1}{\line(0,1){0.49}}
\multiput(49.87,37.76)(0.05,0.5){1}{\line(0,1){0.5}}
\multiput(49.92,38.26)(0.04,0.5){1}{\line(0,1){0.5}}
\multiput(49.96,38.75)(0.02,0.5){1}{\line(0,1){0.5}}
\multiput(49.99,39.25)(0.01,0.5){1}{\line(0,1){0.5}}

\linethickness{0.3mm}
\put(20,50){\line(1,0){80}}
\linethickness{0.3mm}
\multiput(100,50)(0.36,-0.12){83}{\line(1,0){0.36}}
\linethickness{0.3mm}
\multiput(50,40)(0.6,0.12){83}{\line(1,0){0.6}}
\linethickness{0.3mm}
\multiput(20,50)(0.12,-0.12){83}{\line(1,0){0.12}}
\linethickness{0.3mm}
\multiput(40,50)(1.08,-0.12){83}{\line(1,0){1.08}}
\put(30,36){\makebox(0,0)[cc]{$x$}}

\put(17,47){\makebox(0,0)[cc]{$\tilde x$}}

\put(39,54){\makebox(0,0)[cc]{$w$}}

\put(53,36){\makebox(0,0)[cc]{$u$}}

\put(100,55){\makebox(0,0)[cc]{$z$}}

\put(130,36){\makebox(0,0)[cc]{$b$}}

\end{picture}

\noindent
Comme $u\in \de\tg(\tilde x,b)$ et $z\in \alpha\tg(b,u)$
le a) du lemme~\ref{geod-comp-xabc} montre que 
$z\in (\alpha+\de)\tg(\tilde x,b)$.
On en déduit que $w\in (\alpha+\de)\tg(\tilde x,b)$.
 Alors $(H_{\de}(w,\tilde x,u,b))$ montre que 
$$d(w,u)\leq \max(d(w,\tilde x)-d(b,\tilde x)+d(b,u)+\de, d(w,b)-d(b,\tilde x)+d(\tilde x,u)+\de) .$$
Or $d(w,\tilde x)-d(b,\tilde x)+d(b,u)+\de\leq -d(b,w)+\alpha+\de+d(b,u)+\de
= 0$ et $d(w,b)-d(b,\tilde x)+d(\tilde x,u)+\de\leq d(w,b)-d(b,u)+2\de
= \alpha+4\de$. On en déduit $d(w,u)\leq \alpha+4\de$. De plus 
$(H_{\de}(x,\tilde x,w,b))$ montre que 
$$d(x,w)\leq \max(d(x,\tilde x)-d(b,\tilde x)+d(b,w)+\de, d(x,b)-d(b,\tilde x)+d(\tilde x,w)+\de).$$
Or $d(x,\tilde x)-d(b,\tilde x)+d(b,w)+\de\leq l$ car $d(x,\tilde x)\leq l$ et 
$d(b,\tilde x)\geq d(b,w)+\de$ et $d(x,b)-d(b,\tilde x)+d(\tilde x,w)+\de \leq
d(x,b)-d(b,w)+\alpha+2\de= l$ car $d(b,w)=d(b,u)+\alpha+2\de$ et $d(x,b)=l+d(u,b)$. Donc on a bien $w\in B(x,l)$. 
\cqfd

Le lemme suivant est une conséquence des deux précédents. 

\begin{lem}\label{nombre-dist-connaitre-par-point}
Il existe une constante $C=C(\de,K,N,Q,P,M)$ telle que pour $p\in \{1,...,p_{\max}\}$, $k,m,l_{0},...,l_{m}\in \N$ et 
$$(a_{1},\dots,a_{p},S_{0},...,S_{m},(\mathcal Y_{i}^{j})_{i\in \{0,\dots,m\}, j\in \{1,\dots ,l_{i}\}})\in Y_{x}^{p,k,m,(l_{0},...,l_{m})},$$  les distances entre les points de \begin{gather}\label{28dec1223}
B(S_{0}, M)
\cup   \bigcup _{ j\in \{1,\dots ,l_{0}\}} B(\mathcal Y_{0}^{j},  M) \end{gather} et ceux de 
\begin{gather}\label{28dec1225}
\bigcup _{ i\in \{1,\dots ,m\}}
B(S_{i},  M)
\cup   \bigcup _{i\in \{1,\dots,m\}, j\in \{1,\dots ,l_{i}\}} B(\mathcal Y_{i}^{j},  M)  \cup B(x,k+2M)\end{gather} sont déterminées par 
\begin{itemize}
\item a) 
les distances entre les points de (\ref{28dec1225}), 
\item b) les entiers $d(x,S_{0})$ et $d(x,\mathcal Y_{0}^{j})$, 
\item c) 
les distances
entre les points de (\ref{28dec1223}) et $C(1+l_{0})$ points de 
(\ref{28dec1225}) (qui sont eux-mêmes déterminés par a) et b))
\end{itemize}
 et de plus 
  les distances entre les points de \eqref{28dec1223}
et ceux de 
\eqref{28dec1225}
  sont déterminées à $C$ près par  a) et b). 
  \end{lem}
\noindent{\bf Remarque. } Dans toutes les situations où on appliquera ce lemme, $l_{0}$ sera majoré par une constante de la forme $C(\de,K,N,Q,P,M)$. 

 \noindent{\bf Démonstration.}
Soient $p\in \{1,...,p_{\max}\}$, $k,m,l_{0},...,l_{m}\in \N$ et  
$$(a_{1},\dots,a_{p},S_{0},...,S_{m},(\mathcal Y_{i}^{j})_{i\in \{0,\dots,m\}, j\in \{1,\dots ,l_{i}\}})\in Y_{x}^{p,k,m,(l_{0},...,l_{m})}.$$ 
On applique le lemme~\ref{i<j<k-distances} (complété par la remarque que, dans les notations de ce lemme, $4P\tg(u,b)\subset 4P\tg(x,b)$)   avec 
\begin{itemize}
\item $(x,b)$ au lieu de $(c,d)$, 
\item $\{w_{i},i\in I\}$ égal à $\bigcup_{i\in \{0,\dots,m\}}S_{i}\cup  \bigcup _{i\in \{0,\dots,m\}, j\in \{1,\dots ,l_{i}\}} \mathcal Y_{i}^{j}$,
\item $J$ la partie de $I$ telle que $\{w_{j},j\in J\}=S_{0}\cup  \bigcup _{ j\in \{1,\dots ,l_{0}\}} \mathcal Y_{0}^{j}$,
\item  et $(\alpha_{i},\rho_{i})$ égal à $(4P,M)$ pour tout $i$
(les hypothèses sont satisfaites grâce au lemme~\ref{lemme-S0-...Sm}). 
\end{itemize}
  Puis on  applique le lemme~\ref{distances-Bxk-C} à $z$ parcourant (\ref{28dec1223}),   
$\alpha=4P+2M$ et $l=k+2M$ (grâce  aux lemmes~\ref{lemme-S0-...Sm} et ~\ref{xx'yy'zz'}, 
(\ref{28dec1223}) est inclus dans $(4P+2M)\tg(b,u)$). \cqfd

\subsection{Continuité de  $\del$ et $J_{x}$}\label{par-cont-del-Jx}

On introduit  d'abord une variante $\|.\|_{\H^{\rightarrow}_{x,s}}$ de la norme $\|.\|_{\H_{x,s}}$ et on montre que ces deux normes sont équivalentes. 

Soient $p\in \{1,...,p_{\max}\}$ et $k,m,l_{0},...,l_{m}\in \N$. 
On note $Y_{x}^{\rightarrow,p,k,m,(l_{0},...,l_{m})} $ l'ensemble défini de la même fa\c con que $ Y_{x}^{p,k,m,(l_{0},...,l_{m})}$ mais en ajoutant la condition 
\begin{itemize}
\item v)  ou bien $d(x,S_{0})> k+P$, ou bien $k=0,m=0,l_{0}=0,d(x,S_{0})\leq P$. 
\end{itemize}
  \label{Yxrightarrow}

On définit le quotient $\overline Y_{x}^{\rightarrow,p,k,m,(l_{0},...,l_{m})}$ et l'application $\pi_{x}^{\rightarrow,p,k,m,(l_{0},...,l_{m})}:Y_{x}^{\rightarrow,p,k,m,(l_{0},...,l_{m})}\to \overline Y_{x}^{\rightarrow,p,k,m,(l_{0},...,l_{m})}$ 
de la même fa\c con que $\overline Y_{x}^{p,k,m,(l_{0},...,l_{m})}$ et $\pi_{x}^{p,k,m,(l_{0},...,l_{m})}$. 
Cela fournit un diagramme commutatif
$$\begin{matrix}Y_{x}^{\rightarrow,p,k,m,(l_{0},...,l_{m})}& \hookrightarrow &Y_{x}^{p,k,m,(l_{0},...,l_{m})} \\ \downarrow & & \downarrow \\
\overline Y_{x}^{\rightarrow,p,k,m,(l_{0},...,l_{m})}& \hookrightarrow &\overline Y_{x}^{p,k,m,(l_{0},...,l_{m})}\end{matrix}$$
où les flèches verticales sont des surjections. Ce diagramme est  cartésien au sens où les flèches horizontales induisent des bijections sur les fibres des flèches verticales.

On note $\|.\|_{\H^{\rightarrow}_{x,s}(\Delta_{p})}$ la norme sur  $\C^{(\Delta_{p})}$  \label{norme-arrow}
 donnée par la formule (\ref{formule-norme}) en rempla\c cant 
 $Y_{x}^{p,k,m,(l_{0},...,l_{m})}$, $\overline Y_{x}^{p,k,m,(l_{0},...,l_{m})}$ et 
 $\pi_{x}^{p,k,m,(l_{0},...,l_{m})}$ 
  par  $Y_{x}^{\rightarrow,p,k,m,(l_{0},...,l_{m})}$,   $\overline Y_{x}^{\rightarrow,p,k,m,(l_{0},...,l_{m})}$ et $\pi_{x}^{\rightarrow,p,k,m,(l_{0},...,l_{m})}$. On remarque que les sous-espaces de $\C^{(\Delta_{p})}$ engendrés par les $e_{S}$ pour $d(x,S)\leq P$, resp. $d(x,S)> P$ sont orthogonaux pour la norme pré-hilbertienne $\|.\|_{\H^{\rightarrow}_{x,s}(\Delta_{p})}$, et que sur le premier la norme est donnée par 
  $\|f\|^{2}_{\H^{\rightarrow}_{x,s}(\Delta_{p})}=p! \sum_{S}e^{2sd(x,S)}|f(S)|^{2}$,  car si $S\in \Delta_{p}$ vérifie $d(x,S)\leq P$, on a $S\subset B(x,M)$. 
  En effet on suppose $M\geq P+N$, ce qui est permis par $(H_{M})$.

 \begin{lem}\label{equiv-gauche}
 Les normes $\|.\|_{\H_{x,s}(\Delta_{p})}$ et  $\|.\|_{\H^{\rightarrow}_{x,s}(\Delta_{p})}$ sont équivalentes. 
 \end{lem}
 \noindent{\bf Démonstration.}
 D'abord il est évident que $ \|.\|_{\H^{\rightarrow}_{x,s}(\Delta_{p})}\leq \|.\|_{\H_{x,s}(\Delta_{p})}$. 
 Soit $$(a_{1},\dots,a_{p},S_{0},...,S_{m},(\mathcal Y_{i}^{j})_{i\in \{0,\dots,m\}, j\in \{1,\dots ,l_{i}\}})\in Y_{x}^{p,k,m,(l_{0},...,l_{m})}-Y_{x}^{\rightarrow,p,k,m,(l_{0},...,l_{m})}.$$ On a alors $ d(x,S_{0})\leq k+P$ et $d(x,S_{i})> k+P$ pour $i\geq 1$. 
 Le lemme~\ref{m-fini} montre que 
 $ d_{\max}(x,S_{0})\geq d_{\max}(x,S_{1})\geq \pp \geq d_{\max}(x,S_{m})$. 
Donc 
$$S_{0}\cup \pp \cup S_{m}\subset B(x,k+N+P).$$
 Le lemme~\ref{dmaxY-dmaxS} implique alors 
  $$\bigcup _{i\in \{0,\dots,m\}, j\in \{1,\dots ,l_{i}\}}\mathcal Y_{i}^{j}\subset 
B(x,k+N+3P+\de).$$

Grâce à  $(H_{M})$ on  suppose  
$M\geq N+3P+\de$,  d'où $$\bigcup _{ i\in \{0,\dots ,m\}}
B(S_{i}, M)
\cup   \bigcup _{i\in \{0,\dots,m\}, j\in \{1,\dots ,l_{i}\}} B(\mathcal Y_{i}^{j},  M) \subset  B(x,k+2M).$$
 Donc l'ensemble des points entre lesquels on veut connaître les distances, à savoir 
 $$\bigcup _{ i\in \{0,\dots ,m\}}
B(S_{i}, M)
\cup  \bigcup _{i\in \{0,\dots,m\}, j\in \{1,\dots ,l_{i}\}} B(\mathcal Y_{i}^{j},  M) \cup  B(x,k+2M)$$
 est simplement égal à 
 $B(x,k+2M)$. 
Il en résulte que le singleton $$\{(a_{1},\dots,a_{p},S_{0},...,S_{m},(\mathcal Y_{i}^{j})_{i\in \{0,\dots,m\}, j\in \{1,\dots ,l_{i}\}})\}$$  est une classe d'équivalence dans $\overline Y_{x}^{p,k,m,(l_{0},...,l_{m})}$. On a donc 
\begin{gather*}\|f\|^{2}_{\H_{x,s}(\Delta_{p})}-\|f\|^{2}_{\H^{\rightarrow}_{x,s}(\Delta_{p})}
\leq \sum _{S_{0}\in \Delta_{p}}\Bigg(
\sum_{k\geq d(x,S_{0})-P}
e^{2s(d(x,S_{0})-k)} 
\sum _{m,l_{0},\dots ,l_{m}\in \N} \\  B^{-(m+\sum_{i=0}^{m}l_{i})} 
  \sum_{\substack{(a_{1},\dots,a_{p},S_{0},...,S_{m},(\mathcal Y_{i}^{j})_{i\in \{0,\dots,m\}, j\in \{1,\dots ,l_{i}\}})\\ \in Y_{x}^{p,k,m,(l_{0},...,l_{m})}\text{\ tel que \ }S_{0}=\{a_{1},\dots,a_{p}\}}} \Big(\prod_{i=0}^{m}s_{i}(Z)^{-l_{i}} \Big)
 \Bigg)
|f(S_{0})|^{2}.\end{gather*}
Soient $D$ et $C$ comme dans les lemmes~\ref{m-fini}       
 et~\ref{existence-C} (on rappelle que ce sont des constantes de la forme $C(\de,K,N,Q,P)$). Etant donnés $S_{0},k,m,l_{0},\dots ,l_{m}$
 tels que $k\geq d(x,S_{0})-P$, 
 le nombre de possibilités pour  $(S_{1},\pp,S_{m})$ est inférieur ou égal à $ 
 e^{D(P+m)}$ d'après le lemme~\ref{m-fini} 
et  le nombre de possibilités pour $(\mathcal Y_{i}^{j})_{i\in \{0,\dots,m\}, j\in \{1,\dots ,l_{i}\}}$ est inférieur ou égal à  $ C^{\sum_{i=0}^{m}l_{i}}\prod_{i=0}^{m}s_{i}(Z)^{l_{i}} $ d'après le lemme~\ref{existence-C}. Donc 
\begin{gather}\nonumber \|f\|^{2}_{\H_{x,s}(\Delta_{p})}-\|f\|^{2}_{\H^{\rightarrow}_{x,s}(\Delta_{p})}
\leq p! \sum _{S_{0}\in \Delta_{p}}\Bigg(
\sum_{k\geq d(x,S_{0})-P}\\ \nonumber 
e^{2s(d(x,S_{0})-k)} 
\sum _{m,l_{0},\dots ,l_{m}\in \N} B^{-(m+\sum_{i=0}^{m}l_{i})} 
e^{D(P+m)}C^{\sum_{i=0}^{m}l_{i}}
 \Bigg)
|f(S_{0})|^{2}\\ \nonumber 
=p! \frac{e^{DP}}{1-B^{-1}(C+e^{D})}\sum _{S_{0}\in \Delta_{p}}
\sum_{k\geq d(x,S_{0})-P}
e^{2s(d(x,S_{0})-k)} 
|f(S_{0})|^{2}\\ \label{ineg-28dec1237}\leq p! \frac{e^{(D+2s)P}}{(1-B^{-1}(C+e^{D}))(1-e^{-2s})}\|f\|^{2}_{\ell^{2}(\Delta_{p})}.\end{gather}

 Pour $(a_{1},\pp,a_{p},S_{0})\in Y^{p,k,0,(0)}_{x}$ tel que $S_{0}$ vérifie 
 \begin{itemize}
 \item ou bien $d(x,S_{0})=k+P+1$ 
 \item ou bien $k=0$ et 
 $d(x,S_{0})\leq P$,
 \end{itemize} le singleton 
 $\{(a_{1},\pp,a_{p},S_{0})\}$ est une classe d'équivalence dans $\overline Y^{\rightarrow,p,k,0,(0)}_{x}$. 
 En limitant la somme qui définit $\|f\|^{2}_{\H^{\rightarrow}_{x,s}(\Delta_{p})}$ à ces éléments-là, 
  on voit que 
 \begin{gather}\label{ineg-28dec1238}\|f\|^{2}_{\H^{\rightarrow}_{x,s}(\Delta_{p})}\geq p!  \|f\|_{\ell^{2}(\Delta_{p})}^{2}.\end{gather}

On déduit des inégalités (\ref{ineg-28dec1237}) et (\ref{ineg-28dec1238}) qu'il existe une constante $C=C(\de,K,N,Q,P,M,s,B)$ telle que 
$ \|f\|^{2}_{\H_{x,s}(\Delta_{p})}\leq 
C \|f\|^{2}_{\H^{\rightarrow}_{x,s}(\Delta_{p})}$
pour tout $f\in \C^{(\Delta_{p})}$. \cqfd

Soit $\P$ \label{def-mathcalP} le projecteur orthogonal sur le sous-espace vectoriel de $\H^{\rightarrow}_{x,s}(\Delta_{p})$ engendré par les $e_{S}$ pour $S\in \Delta_{p}$ tel que $d(x,S)\leq P$, de sorte que $(\P f)(S)=f(S)$ si 
$d(x,S)\leq P$ et $(\P f)(S)=0$ sinon.  

Pour $f\in \C^{(\Delta_{p})}$ on a 
\begin{gather}\nonumber \|(1-\P)f\|_{\H^{\rightarrow}_{x,s}(\Delta_{p})}^{2} =\sum _{k,m,l_{0},\dots ,l_{m}\in \N} 
B^{-(m+\sum_{i=0}^{m}l_{i})}
\sum_{Z\in \overline Y_{x}^{p,k,m,(l_{0},...,l_{m})},  r_{0}(Z)> k+P} 
\\ \label{utilite-rightarrow-11j}
e^{2s(r_{0}(Z)-k)}
\Big(\prod_{i=0}^{m}s_{i}(Z)^{-l_{i}} \Big)
\sharp \big((\pi_{x}^{p,k,m,(l_{0},...,l_{m})})^{-1}(Z)\big)^{-\alpha }
\big|\xi_{Z}(f)\big|^{2} \end{gather}
Cette formule est la raison pour laquelle on a introduit la norme 
$\|.\|_{\H^{\rightarrow}_{x,s}(\Delta_{p})}$. En effet pour montrer la continuité de $\del$ ou de $J_{x}$ on cherchera à majorer $|\xi_{Z}(\del f)|^{2}$ ou 
 $|\xi_{Z}(J_{x} f)|^{2}$ par une combinaison de $|\xi_{\tilde Z}(f)|^{2}$ avec $\tilde Z$ vérifiant notamment $r_{1}(\tilde Z)=r_{0}(Z)$. Comme la condition 
 $r_{1}(\tilde Z)> k+P$ est imposée par la condition ii) de la définition~\ref{defi-Y}, 
 il est très utile d'avoir $r_{0}(Z)> k+P$.

 \begin{prop}\label{continuite-del}
 Pour tout $p\in \{1,\pp,p_{\max}\}$, $\del$ se prolonge en un opérateur continu de $\H_{x,s}(\Delta_{p})$ dans $\H_{x,s}(\Delta_{p-1})$. 
 \end{prop}
 
 \noindent{\bf Démonstration.} 
Supposons d'abord $p=1$. On va montrer qu'il existe une constante $C=C(\de,K,N,Q,P,M,s,B)$ telle que $|\sum_{a\in X}f(a)| \leq C\|f\|_{\H_{x,s}(\Delta_{1})}$. 
D'après la formule~(\ref{formule-norme}), $\|f\|^{2}_{\H_{x,s}(\Delta_{1})}$ est une somme sur $k,m,l_{0},\pp,l_{m}$ et en limitant cette somme à $k=0,m=0,l_{0}=0$ on voit que 
$$\|f\|^{2}_{\H_{x,s}(\Delta_{1})}\geq 
\sum_{Z\in \overline Y_{x}^{1,0,0,(0)}}  e^{2sr_{0}(Z)}\sharp \big((\pi_{x}^{1,0,0,(0)})^{-1}(Z)\big)^{-\alpha }\Big|\sum _{(a,\{a\}) \in 
(\pi_{x}^{1,0,0,(0)})^{-1}(Z)} f(a)\Big|^{2}.$$
De plus $\overline Y_{x}^{1,0,0,(0)}$ s'identifie au  quotient de $X$ pour la relation d'équivalence suivante : $a$ et $b$ sont équivalents s'il existe une isométrie de $B(a,M)\cup B(x,2M)$ vers $B(b,M)\cup B(x,2M)$ qui est l'identité sur $B(x,2M)$ et applique $a$ sur $b$. 
Cette relation d'équivalence détermine $d(x,a)$ et inversement il existe $C=C(\de,K,N,Q,P,M)$ telle que pour tout $r\in \N$, l'ensemble $\{a\in X,d(x,a)=r\}
$ est réunion d'au plus $C$ classes d'équivalences de $\overline  Y_{x}^{1,0,0,(0)}$. 
 Il existe $D'=C(\de,K)$ telle que pour tout $r\in \N$, le cardinal de 
l'ensemble $\{a\in X,d(x,a)=r\}$
 soit inférieur ou égal à $e^{D'r}$. On a donc, par Cauchy-Schwarz, 
 $$\|f\|^{2}_{\H_{x,s}(\Delta_{1})}\geq C^{-1}
\sum_{r\in \N}  e^{(2s-\alpha D')r}\Big|\sum _{a \in X, d(x,a)=r
} f(a)\Big|^{2}.$$
 Grâce à  $(H_{\alpha})$  on  suppose   $2s-\alpha D'\geq s$. 
 Par Cauchy-Schwarz, 
 \begin{gather*}|\sum_{a\in X}f(a)|^{2} \leq \Big(\sum  _{r\in \N}e^{-(2s-\alpha D')r} \Big)
\Big(\sum_{r\in \N}  e^{(2s-\alpha D')r}\Big|\sum _{a \in X, d(x,a)=r
} f(a)\Big|^{2}\Big)\\ \leq \frac{C}{1-e^{-(2s-\alpha D')}} \|f\|^{2}_{\H_{x,s}(\Delta_{1})}  \leq \frac{C}{1-e^{-s}} \|f\|^{2}_{\H_{x,s}(\Delta_{1})}.\end{gather*}
 
 Soit maintenant $p\in \{2,\pp,p_{\max}\}$. Grâce au lemme~\ref{equiv-gauche}, il suffit de  montrer 
 qu'il existe une constante $C=C(\de,K,N,Q,P,M,s,B)$ telle que 
 $$\|\del f\|_{\H^{\rightarrow}_{x,s}(\Delta_{p-1})} \leq C\|f\|_{\H_{x,s}(\Delta_{p})}$$ pour tout $f\in \C^{(\Delta_{p})}$. Grâce à (\ref{estimation-norme-Hsx}) et au lemme~\ref{minoration-normeHsx}, 
 il est  clair que $$\|\P(\del f)\|_{\H^{\rightarrow}_{x,s}(\Delta_{p-1})} \leq C\|f\|_{\H_{x,s}(\Delta_{p})}$$ pour une  constante $C=C(\de,K,N,Q,P,M,s,B)$. Il reste donc à montrer qu'il existe une constante $C=C(\de,K,N,Q,P,M,s,B)$ telle que 
 \begin{gather}\label{ineg-(1-P)-del-f}\|(1-\P)(\del f)\|_{\H^{\rightarrow}_{x,s}(\Delta_{p-1})} \leq C\|f\|_{\H_{x,s}(\Delta_{p})}\end{gather} pour tout $f\in \C^{(\Delta_{p})}$. 
 Par (\ref{utilite-rightarrow-11j}) on a 
\begin{gather*}\|(1-\P)(\del f)\|_{\H^{\rightarrow}_{x,s}(\Delta_{p-1})}^{2} =\sum _{k,m,l_{0},\dots ,l_{m}\in \N} 
B^{-(m+\sum_{i=0}^{m}l_{i})}
\sum_{Z\in \overline Y_{x}^{p-1,k,m,(l_{0},...,l_{m})},  r_{0}(Z)> k+P} \\ e^{2s(r_{0}(Z)-k)}
\Big(\prod_{i=0}^{m}s_{i}(Z)^{-l_{i}} \Big)
\sharp \big((\pi_{x}^{p-1,k,m,(l_{0},...,l_{m})})^{-1}(Z)\big)^{-\alpha }
\big|\xi_{Z}(\del f)\big|^{2} \end{gather*}
%\Big|\sum _{(a_{1},\dots,a_{p-1},S_{1},...,S_{m},(\mathcal Y_{i}^{j})_{i\in \{0,\dots,m\}, j\in \{1,\dots ,l_{i}\}}) \in 
%(\pi_{x}^{p-1,k,m,(l_{0},...,l_{m})})^{-1}(Z)} (\del f)(a_1,...,a_{p-1})\Big|^{2}.$$
Soient   $k,m,l_{0},\dots ,l_{m}\in \N $ et $Z\in \overline Y_{x}^{p-1,k,m,(l_{0},...,l_{m})}$ vérifiant $r_{0}(Z)> k+P$.  
On pose  $\tilde l_{0}=0, \tilde l_{i}=l_{i-1}$ pour $i\in \{1,\pp,m+1\}$. Alors  
 on a 
 \begin{gather}\nonumber \xi_{Z}(\del f)=\sum _{\substack{(a_{1},\dots,a_{p-1},S_{0},...,S_{m},(\mathcal Y_{i}^{j})_{i\in \{0,\dots,m\}, j\in \{1,\dots ,l_{i}\}})\\  \in 
(\pi_{x}^{p-1,k,m,(l_{0},...,l_{m})})^{-1}(Z)}} (\del f)(a_1,...,a_{p-1}) =\sum_{\tilde Z\in \Lambda_{Z}} \\ \label{somme-tildeZ-Y-24oct09}   
 \sum _{\substack{(\tilde a_{1},\dots,\tilde a_{p},\tilde S_{0},...,\tilde S_{m+1},(\tilde {\mathcal Y}_{i}^{j})_{i\in \{0,\dots,m+1\}, j\in \{1,\dots ,\tilde l_{i}\}}) \\ \in 
(\pi_{x}^{p,k,m+1,(\tilde l_{0},...,\tilde l_{m+1})})^{-1}(\tilde Z)}}  f(\tilde a_{1},\dots,\tilde a_{p})=\sum_{\tilde Z\in \Lambda_{Z}}\xi_{\tilde Z}(f)\end{gather}
 où $\Lambda_{Z}$ est la partie de  $ \overline Y_{x}^{p,k,m+1,(\tilde l_{0},...,\tilde l_{m+1})}$ formée des $\tilde Z$ tels que pour tout 
 $$(\tilde a_{1},\dots,\tilde a_{p},\tilde S_{0},...,\tilde S_{m+1},(\tilde {\mathcal Y}_{i}^{j})_{i\in \{0,\dots,m+1\}, j\in \{1,\dots ,\tilde l_{i}\}})\in (\pi_{x}^{p,k,m+1,(\tilde l_{0},...,\tilde l_{m+1})})^{-1}(\tilde Z)$$ on ait  
 $\tilde S_{1}=\{\tilde a_{2},\dots,\tilde a_{p}\}$ et 
 $$(\tilde a_{2},\dots,\tilde a_{p},\tilde S_{1},...,\tilde S_{m+1},(\tilde {\mathcal Y}_{i+1}^{j})_{i\in \{0,\dots,m\}, j\in \{1,\dots , l_{i}\}})\in (\pi_{x}^{p-1,k,m,(l_{0},...,l_{m})})^{-1}(Z).$$  
% On a toujours $d(x,\tilde S_{0})\leq d(x,\tilde S_{1})\leq d(x,\tilde S_{0})+N$.  
Il existe $C_{1}=C(\de,K,N,Q,P,M)$ telle que $\sharp \Lambda_{Z}\leq C_{1}$. 
 En effet, grâce au lemme~\ref{nombre-dist-connaitre-par-point}, pour connaître les distances entre les points de 
 $B(\tilde S_{0}, M)$ 
 et ceux de 
 \begin{gather}\label{ens-28dec1242}  
\bigcup _{ i\in \{1,\dots ,m+1\}}
B(\tilde S_{i}, M)
\cup  \bigcup _{i\in \{1,\dots,m+1\}, j\in \{1,\dots ,\tilde l_{i}\}}B(\tilde {\mathcal Y}_{i}^{j}, M) \cup B(x,k+2M)\end{gather} 
il suffit de connaître les distances entre les points de $B(\tilde S_{0}, M)$  
 et $C$ points de (\ref{ens-28dec1242}), avec $C=C(\de,K,N,Q,P,M)$ et 
 comme $\tilde S_{0}=\tilde S_{1}\cup \{\tilde a_{1}\}$ et $d(\tilde a_{1}, \tilde a_{2})\leq N$, 
  ces distances sont déterminées à $ N+M$ près par les distances de $\tilde a_{2}$ à ces $C$ points (qui font partie de la donnée de $Z$). 
  Comme $\sharp \Lambda_{Z}\leq C_{1}$,  grâce à  (\ref{somme-tildeZ-Y-24oct09}) et par  
 Cauchy-Schwarz, on obtient  
 \begin{gather}\label{CS-xi-xi-10fev11}\big|\xi_{Z}(\del f)\big|^{2}\leq C_{1}\sum_{\tilde Z\in \Lambda_{Z}}|\xi_{\tilde Z}(f)|^{2}.\end{gather} 
% $$ C_{1 }\sum_{\tilde Z} \Big| \sum _{\substack{(\tilde a_{1},\dots,\tilde a_{p},\tilde S_{1},...,\tilde S_{m+1},(\tilde {\mathcal Y}_{i}^{j})_{i\in \{0,\dots,m+1\}, j\in \{1,\dots ,\tilde l_{i}\}})\\ \in 
%(\pi_{x}^{p,k,m+1,(\tilde l_{0},...,\tilde l_{m+1})})^{-1}(\tilde Z)}}  f(\tilde a_{1},\dots,\tilde a_{p})\Big|^{2}$$
 %où la somme est prise sur le même ensemble de  $\tilde Z$ que dans 
 %(\ref{somme-tildeZ-Y-24oct09}). 
 De plus pour $\tilde Z\in \Lambda_{Z}$ on a \begin{gather*}
 \prod_{i=0}^{m+1}s_{i}(\tilde Z)^{-\tilde l_{i}} =\prod_{i=0}^{m}s_{i}(Z)^{-l_{i}} , 
 \\ 
 \sharp (\pi_{x}^{p,k,m+1,(\tilde l_{0},...,\tilde l_{m+1})})^{-1}(\tilde Z)\leq C\sharp (\pi_{x}^{p-1,k,m,(l_{0},...,l_{m})})^{-1}(Z)\end{gather*} avec $C=C(\de,K,N)$ et $|r_{0}(Z)-r_{0}(\tilde Z)|\leq N$. Enfin $\tilde Z$ détermine $Z$, donc quand on somme sur $Z$, chaque $\tilde Z$ ne peut apparaître qu'une fois.  L'inégalité (\ref{ineg-(1-P)-del-f}) en résulte facilement et ceci termine la démonstration 
 de la proposition~\ref{continuite-del}.   \cqfd
 
 \noindent{\bf Remarque.} Le coeur de  la démonstration ci-dessus est formé par
\begin{itemize}
 \item l'égalité \eqref{somme-tildeZ-Y-24oct09}, que l'on peut mettre sous la forme ${}^{t}\del(\xi_{Z})=\sum_{\tilde Z\in \Lambda_{Z}}\xi_{\tilde Z}$, 
 \item l'application de Cauchy-Schwarz qui fournit l'inégalité \eqref{CS-xi-xi-10fev11}, 
 \item le fait que les normes sont des sommes pondérées des $|\xi_{Z}|^{2}$.
 \end{itemize}
 Dans la suite les arguments seront plus compliqués mais ils reposeront tous sur ce principe. 
 
 Avant de montrer la continuité de $J_{x}$ on va introduire des nouvelles normes pré-hilbertiennes $\|.\|_{\H_{x,s}^{\natural,\mu_{0},\mu_{1}}}$  sur $\C^{(\Delta_{p})}$ (pour $\mu_{0},\mu_{1}\in \N$) et montrer qu'elles sont équivalentes à $\|.\|_{\H_{x,s}}$. 
 Ces normes seront obtenues en ajoutant aux parties $\mathcal Y_{i}^{j}$ qui intervenaient dans la définition~\ref{defi-Y} de nouvelles parties, notées  $\mathcal Z_{i}^{j}$ dans la définition ci-dessous. On verra dans la démonstration de la continuité de $J_{x}$ (proposition~\ref{continuite-Jx}) que la connaissance des points de $B(\mathcal Z_{i}^{j},M)$ détermine exactement les différentes moyennes intervenant dans la formule pour $J_{x}$, d'où l'intérêt de ces parties supplémentaires. 
 On va voir dans la preuve de l'équivalence des normes $\|.\|_{\H_{x,s}^{\natural,\mu_{0},\mu_{1}}}$ et $\|.\|_{\H_{x,s}}$ (lemme~\ref{eq-normes-natural}) que chaque partie $\mathcal Z_{i}^{j}$ peut être oubliée ou reconsidérée comme une partie  $\mathcal Y_{i'}^{j'}$ supplémentaire. 
 Cependant $i'$ est déterminé par $d(x,\mathcal Z_{i}^{j})$ d'une fa\c con assez compliquée. Ces nouvelles normes rendent donc la preuve de la continuité de $J_{x}$ (proposition~\ref{continuite-Jx}) beaucoup plus lisible et elles resserviront de plus pour montrer la continuité des autres opérateurs (proposition~\ref{continuite-del-J-conj}) et l'équivariance à compact près de tous les opérateurs (proposition~\ref{lemme-compacite-equiv}). 
 Inversement on n'a pas inclus ces parties $\mathcal Z_{i}^{j}$ dans la 
 définition~\ref{defi-Y}  car elles auraient rendu beaucoup plus difficile la preuve des propriétés d'équivariance de la norme $\|.\|_{\H_{x,s}}$ (proposition~\ref{equiv-normeYZ}, démontrée dans le 
  sous-paragraphe~\ref{action-G-normes}).

 Soient $p\in \{1,...,p_{\max}\}$ et $k,m,l_{0},...,l_{m},\lambda_{0},\lambda_{1}\in \N$ vérifiant $\lambda_{1}=0$ si $m=0$. 
\begin{defi}\label{defi-Y-natural}
On note $Y_{x}^{\natural,p,k,m,(l_{0},...,l_{m}),\lambda_{0},\lambda_{1}}$ l'ensemble des $(p+m+1+\sum_{i=0}^{m}l_{i}+\lambda_{0}+\lambda_{1})$-uplets $$(a_{1},\dots,a_{p},S_{0},...,S_{m},(\mathcal Y_{i}^{j})_{i\in \{0,\dots,m\}, j\in \{1,\dots ,l_{i}\}}, (\mathcal Z_{i}^{j})_{i\in \{0,1\}, j\in \{1,\dots ,\lambda_{i}\}})$$ tels que  
 \begin{itemize}   
\item 
$(a_{1},\dots,a_{p},S_{0},...,S_{m},(\mathcal Y_{i}^{j})_{i\in \{0,\dots,m\}, j\in \{1,\dots ,l_{i}\}} )$ appartient à $Y_{x}^{p,k,m,(l_{0},...,l_{m})}$, c'est-à-dire 
vérifie les conditions i), ii), iii), iv) de la définition~\ref{defi-Y}, 
\item  pour $i\in \{0,1\}$ et $j\in \{1,\dots ,\lambda_{i}\}$, $\mathcal Z_{i}^{j}$ est une partie non vide de $X$ de diamètre inférieur ou égal à $P/3$  et  $\mathcal Z_{i}^{j}\subset \bigcup_{a\in S_{i}}\geod(x,a)$. 
\end{itemize}
\end{defi}

On introduit  une partition de $Y_{x}^{\natural,p,k,m,(l_{0},...,l_{m}),\lambda_{0},\lambda_{1}}$  pour la relation d'équivalence  suivante : 
$$(a_{1},\dots,a_{p},S_{0},...,S_{m},(\mathcal Y_{i}^{j})_{i\in \{0,\dots,m\}, j\in \{1,\dots ,l_{i}\}}, (\mathcal Z_{i}^{j})_{i\in \{0,1\}, j\in \{1,\dots ,\lambda_{i}\}}
)$$ 
et $$(\hat a_{1},\dots,\hat a_{p},\hat S_{0},...,\hat S_{m},(\hat {\mathcal Y}_{i}^{j})_{i\in \{0,\dots,m\}, j\in \{1,\dots ,l_{i}\}}, (\hat {\mathcal Z}_{i}^{j})_{i\in \{0,1\}, j\in \{1,\dots ,\lambda_{i}\}}
)$$
sont en relation 
s'il existe une isométrie de $$
\bigcup _{ i\in \{0,\dots ,m\}}
B(S_{i}, M)
\cup  \bigcup _{i\in \{0,\dots,m\}, j\in \{1,\dots ,l_{i}\}} B(\mathcal Y_{i}^{j},M) $$ $$\cup \bigcup _{i\in \{0,1\}, j\in \{1,\dots ,\lambda_{i}\}} B(\mathcal Z_{i}^{j}, M)
\cup B(x,k+2M)$$ vers 
$$\bigcup _{ i\in \{0,\dots ,m\}}
B(\hat S_{i},M)
\cup  \bigcup _{i\in \{0,\dots,m\}, j\in \{1,\dots ,l_{i}\}} B(\hat {\mathcal Y}_{i}^{j}, M)$$ $$ \cup \bigcup _{i\in \{0,1\}, j\in \{1,\dots ,\lambda_{i}\}} B(\hat{\mathcal Z}_{i}^{j}, M)\cup B(x,k+2M)$$ 
 qui envoie 
$a_{i}$ sur $\hat a_{i}$ pour $i\in \{1,\dots,p\}$,  
 $S_{i}$ sur $\hat S_{i}$
 pour $i\in \{0,\dots ,m\}$, $\mathcal Y_{i}^{j}$ sur $\hat {\mathcal Y}_{i}^{j}$  pour 
 $i\in \{0,\dots,m\}, j\in \{1,\dots ,l_{i}\}$, $\mathcal Z_{i}^{j}$ sur $\hat {\mathcal Z}_{i}^{j}$  pour 
 $i\in \{0,1\}, j\in \{1,\dots ,\lambda_{i}\}$
 et est l'identité sur $B(x,k+2M)$.

On note $\overline Y_{x}^{\natural,p,k,m,(l_{0},...,l_{m}),\lambda_{0},\lambda_{1}}$
le quotient de $Y_{x}^{\natural,p,k,m,(l_{0},...,l_{m}),\lambda_{0},\lambda_{1}}$ pour cette relation d'équivalence, et $\pi_{x}^{\natural,p,k,m,(l_{0},...,l_{m}),\lambda_{0},\lambda_{1}}$ l'application quotient. 

\noindent{\bf Notations.} \label{rijmax}
Pour $Z\in \overline Y_{x}^{\natural,p,k,m,(l_{0},...,l_{m}),\lambda_{0},\lambda_{1}}$  on note $r_{0}(Z),\dots,r_{m}(Z),s_{0}(Z),\dots,$   $s_{m}(Z)$,
%$r_{0}^{\max}(Z), ...,r_{m+1}^{\max}(Z)$ 
$(r_{i,j}^{\max}(Z))_{i\in \{0,1\}, j\in \{i,...,m+1\}}$ 
et $(t_{i}^{j}(Z))_{i\in \{0,1\}, j\in \{1,\dots ,\lambda_{i}\}}$ les entiers tels que 
\begin{itemize}
\item $r_{i}(Z)=d(x,S_{i})$ pour  $i\in \{0,\pp,m\}$, 
 \item $s_{i}(Z)=d(S_{i},S_{i+1})+2M$ pour  $i\in \{0,\pp,m-1\}$, 
\item $s_{m}(Z)=d(x,S_{m})-k$,  
\item pour $i\in \{0,1\}$, si $d_{\max}(x,S_{i})\leq k+3P$, 
$r_{i,j}^{\max}(Z)=d_{\max}(x,S_{i})$ pour tout $j\in \{i,...,m+1\}$, 
\item pour $i\in \{0,1\}$, si $d_{\max}(x,S_{i})\geq  k+3P$, 
$r_{i,j}^{\max}(Z)=\max(k+3P,d_{\max}(x,S_{j}))$ pour tout $j\in \{i,...,m\}$
et $r_{i,m+1}^{\max}(Z)=k+3P$, 
%\item pour $i\in \{0,1\}$, $r_{i,m+1}^{\max}(Z)=\min(k+3P,d_{\max}(x,S_{i}))$  et $$r_{i,j}^{\max}(Z)=\max(r_{i,m+1}^{\max}(Z),d_{\max}(x,S_{j}))$$ pour $j\in \{i,...,m\}$
 \item   $t_{i}^{j}(Z)=d(x,\mathcal Z_{i}^{j})$ pour $   i\in \{0,1\}, j\in \{1,\dots ,\lambda_{i}\}$
\end{itemize}
pour tout $(a_{1},\dots,a_{p},S_{0},...,S_{m},(\mathcal Y_{i}^{j})_{i\in \{0,\dots,m\}, j\in \{1,\dots ,l_{i}\}}, (\mathcal Z_{i}^{j})_{i\in \{0,1\}, j\in \{1,\dots ,\lambda_{i}\}})\in (\pi_{x}^{\natural,p,k,m,(l_{0},...,l_{m}),\lambda_{0},\lambda_{1}})^{-1}(Z). $
  
  Pour clarifier le sens des notations ci-dessus, on rappelle 
  que d'après  le lemme~\ref{m-fini}, on a toujours 
  \begin{gather}\label{suite-egal-dmax-clarifie-24oct09}d_{\max}(x,S_{0})\geq d_{\max}(x,S_{1})\geq \pp \geq d_{\max}(x,S_{m}). \end{gather}
  D'autre part,  pour $i\in \{0,1\}$, on a toujours
  \begin{gather}\label{rem-avant-prop-r-i-max-11j}r_{i,i}^{\max}(Z)=d_{\max}(x,S_{i})\text{ \ \  et \ \ }r_{i,m+1}^{\max}(Z)=\min(k+3P,d_{\max}(x,S_{i})).\end{gather}

Le lemme suivant indique quelques propriétés de ces entiers, qui nous seront utiles ensuite. 

\begin{lem}\label{prop-r-i-max}
  Pour $k,m,l_{0},...,l_{m},\lambda_{0},\lambda_{1}\in \N$ (vérifiant $\lambda_{1}=0$ si $m=0$) et  $Z\in \overline Y_{x}^{\natural,p,k,m,(l_{0},...,l_{m}),\lambda_{0},\lambda_{1}}$ on a 

\noindent a) $r_{0,0}^{\max}(Z)\geq    r_{0,1}^{\max}(Z)\geq  ... \geq  r_{0,m+1}^{\max}(Z)$ et $r_{1,1}^{\max}(Z)\geq    r_{1,2}^{\max}(Z)\geq  ... \geq  r_{1,m+1}^{\max}(Z)$, 
  
\noindent  b) pour tout $i\in \{0,1\}$, $r_{i}(Z)\leq r_{i,i}^{\max}(Z) \leq r_{i}(Z)+N$,  
  
\noindent  c) pour $i\in\{0,1\}$ et $j\in \{1,...,\lambda_{i}\}$, 
  $t_{i}^{j}(Z)\leq r_{i,i}^{\max}(Z)$,

  \noindent  d) pour $i\in\{0,1\}$, $r_{i,m+1}^{\max}(Z)=\min(k+3P,r_{i,i}^{\max}(Z))$. 
  \end{lem}
 \noindent{\bf Démonstration.}
L'assertion  a) découle de (\ref{suite-egal-dmax-clarifie-24oct09}),  b) est évidente et pour montrer c) on 
remarque que dans les notations précédentes on a $d(x,\mathcal Z_{i}^{j})\leq d_{\max}(x,S_{i})$ par 
 la dernière condition de la définition~\ref{defi-Y-natural}. 
 Enfin d) résulte de (\ref{rem-avant-prop-r-i-max-11j}). 
 \cqfd

Pour $Z\in \overline Y_{x}^{\natural,p,k,m,(l_{0},...,l_{m}),\lambda_{0},\lambda_{1}}$ on note $\xi_{Z}$ la forme linéaire sur $\C^{(\Delta_{p})}$ définie par 
\begin{gather}\label{defxi-31dec1632}\!\!\!\!\!\!\xi_{Z}(f)=\sum _{\substack{(a_{1},\dots,a_{p},S_{0},...,S_{m},(\mathcal Y_{i}^{j})_{i\in \{0,\dots,m\}, j\in \{1,\dots ,l_{i}\}}, (\mathcal Z_{i}^{j})_{i\in \{0,1\}, j\in \{1,\dots ,\lambda_{i}\}}) \\ \in 
(\pi_{x}^{\natural,p,k,m,(l_{0},...,l_{m}),\lambda_{0},\lambda_{1}})^{-1}(Z)}} 
f(a_1,...,a_p).
\end{gather}
Pour $\mu_{0},\mu_{1}\in \N$ on munit alors $\C^{(\Delta_{p})}$ de la norme pré-hilbertienne, définie par la formule suivante : 
\begin{gather} \nonumber \|f\|_{\H^{\natural,\mu_{0},\mu_{1}}_{x,s}(\Delta_{p})}^{2}
=\sum _{k,m,l_{0},\dots ,l_{m}, \lambda_{0},\lambda_{1}}  
B^{-(m+\sum_{i=0}^{m}l_{i})}
\sum_{Z\in  \overline Y_{x}^{\natural,p,k,m,(l_{0},...,l_{m}),\lambda_{0},\lambda_{1}}}  e^{2s(r_{0}(Z)-k)}
\\ \nonumber  \Big(\prod_{i=0}^{m}s_{i}(Z)^{-l_{i}} \Big)
(r_{0}(Z)+1)^{-\lambda_{0}}(r_{1}(Z)+1)^{-\lambda_{1}}
\\ \label{formule-norme-mu} \sharp \big((\pi_{x}^{\natural,p,k,m,(l_{0},...,l_{m}),\lambda_{0},\lambda_{1}})^{-1}(Z)\big)^{-\alpha }\big|\xi_{Z}(f)\big|^{2}\end{gather}
où la première somme porte sur $k,m,l_{0},...,l_{m},\lambda_{0},\lambda_{1}\in \N$ vérifiant $\lambda_{1}=0$ si $m=0$ et satisfaisant les conditions $$\lambda_{0}\leq \mu_{0}\text{\ \ \ \ et\ \ \ \ } \lambda_{1}\leq \mu_{1}.$$ 

 En utilisant   l'hypothèse $(H_{B})$ nous allons montrer  le lemme suivant.  
 \begin{lem}\label{eq-normes-natural}
  Il existe $C=C(\de,K,N,Q,P,M,s,B)$ tel que pour $\mu_{0},\mu_{1}\in \N$ et  $f\in \C^{(\Delta_{p})}$, 
  $$\|f\|_{\H_{x,s}(\Delta_{p})}^{2}
\leq \|f\|_{\H_{x,s}^{\natural,\mu_{0},\mu_{1}}(\Delta_{p})}^{2}
\leq C^{\mu_{0}+\mu_{1}}\|f\|_{\H_{x,s}(\Delta_{p})}^{2}.$$
 \end{lem}
  \noindent{\bf Démonstration.}
L'inégalité de gauche est évidente, car la somme (\ref{formule-norme}) qui donne $\|f\|_{\H_{x,s}(\Delta_{p})}^{2}$ est une partie de la somme (\ref{formule-norme-mu}) qui donne $\|f\|_{\H_{x,s}^{\natural,\mu_{0},\mu_{1}}(\Delta_{p})}^{2}$ (c'est la partie qui correspond à $\lambda_{0}=\lambda_{1}=0$). Pour montrer l'inégalité de droite on a besoin de deux lemmes préliminaires. L'idée est simplement de reconsidérer chaque partie $\mathcal Z^{j}_{i}$ comme une partie $\mathcal Y^{j'}_{i'}$ supplémentaire, avec $i'$ déterminé par $t^{j}_{i}(Z)=d(x,\mathcal Z^{j}_{i})$, si $t^{j}_{i}(Z)>k+3P$, et d'oublier $\mathcal Z^{j}_{i}$ si $t^{j}_{i}(Z)\leq 
k+3P$. 
 
  On prend $k,m,l_{0},...,l_{m},\lambda_{0},\lambda_{1}\in \N$ (vérifiant $\lambda_{1}=0$ si $m=0$). Soit $\sigma\in\{0,1\}$ (avec $\sigma=0$ si $m=0$). On pose 
  \begin{itemize}
  \item 
  $(\tilde \lambda_{0},\tilde \lambda_{1})= (\lambda_{0}+1 ,\lambda_{1})$ si $
  \sigma=0$ 
  \item et $(\tilde \lambda_{0},\tilde \lambda_{1})= (\lambda_{0} ,\lambda_{1}+1)$ si $
  \sigma=1$.
  \end{itemize} 
 \begin{souslem}\label{remplacement-lambda0+1}
 Soit 
 $Z\in  \overline Y_{x}^{\natural,p,k,m,(l_{0},...,l_{m}),\tilde \lambda_{0},\tilde \lambda_{1}}$ et 
  \begin{gather*}(a_{1},\dots,a_{p},S_{0},...,S_{m},(\mathcal Y_{i}^{j})_{i\in \{0,\dots,m\}, j\in \{1,\dots ,l_{i}\}}, (\mathcal Z_{i}^{j})_{i\in \{0,1\}, j\in \{1,\dots ,\tilde \lambda_{i}\}}
)\\ \in (\pi_{x}^{\natural,p,k,m,(l_{0},...,l_{m}),\tilde \lambda_{0},\tilde \lambda_{1}})^{-1}(Z).
%
%Y_{x}^{\natural,p,k,m,(l_{0},...,l_{m}),\tilde \lambda_{0},\tilde \lambda_{1}}.
 \end{gather*} 
Alors 
\begin{itemize}
\item si $0\leq t_{\sigma}^{\lambda_{\sigma}+1}(Z) \leq r_{\sigma,m+1}^{\max}(Z)$, on a 
\begin{gather*}(a_{1},\dots,a_{p},S_{0},...,S_{m},(\mathcal Y_{i}^{j})_{i\in \{0,\dots,m\}, j\in \{1,\dots ,l_{i}\}}, (\mathcal Z_{i}^{j})_{i\in \{0,1\}, j\in \{1,\dots ,\lambda_{i}\}}
)\\ \in Y_{x}^{\natural,p,k,m,(l_{0},...,l_{m}),\lambda_{0},\lambda_{1}}\end{gather*}
\item si $r_{\sigma,i+1}^{\max}(Z)< t_{\sigma}^{\lambda_{\sigma}+1}(Z) \leq r_{\sigma,i}^{\max}(Z)$, 
pour $i\in \{\sigma,...,m\}$, 
  en posant $\tilde l_{j}=l_{j}$ pour $j\in \{0,...,m\}\setminus \{i\}$, $\tilde l_{i}=l_{i}+1$ et $
\mathcal Y_{i}^{l_{i}+1}=\mathcal Z_{\sigma}^{\lambda_{\sigma}+1}$ on a 
\begin{gather*}(a_{1},\dots,a_{p},S_{0},...,S_{m},(\mathcal Y_{i}^{j})_{i\in \{0,\dots,m\}, j\in \{1,\dots ,\tilde l_{i}\}}, (\mathcal Z_{i}^{j})_{i\in \{0,1\}, j\in \{1,\dots ,\lambda_{i}\}}
)\\ \in Y_{x}^{\natural,p,k,m,(\tilde l_{0},...,\tilde l_{m}),\lambda_{0},\lambda_{1}}\end{gather*}
\end{itemize}
\end{souslem} 
\noindent{\bf Remarque.} Par le c) et le d) du lemme~\ref{prop-r-i-max}, la condition $t_{\sigma}^{\lambda_{\sigma}+1}(Z) \leq r_{\sigma,m+1}^{\max}(Z)$ qui détermine le premier cas 
 est équivalente à 
$t_{\sigma}^{\lambda_{\sigma}+1}(Z) \leq k+3P$. 

 \noindent{\bf Démonstration.}
D'après le c) du lemme~\ref{prop-r-i-max}, on a $t_{\sigma}^{\lambda_{\sigma}+1}(Z) \leq r_{\sigma,\sigma}^{\max}(Z)$ donc on se trouve toujours exactement dans l'un des cas ci-dessus.  
Il n'y a rien à montrer dans le premier cas. 
Supposons maintenant 
$$t_{\sigma}^{\lambda_{\sigma}+1}(Z)\in ]r_{\sigma,m+1}^{\max}(Z),r_{\sigma,m}^{\max}(Z)].$$ Alors nécessairement $r_{\sigma,m+1}^{\max}(Z)< r_{\sigma,m}^{\max}(Z)$, ce qui implique 
$$r_{\sigma,m}^{\max}(Z)=d_{\max}(x,S_{m})\text{\ \  et\ \ }r_{\sigma,m+1}^{\max}(Z)= k+3P.$$ On a  donc $d(x,\mathcal Z_{\sigma}^{\lambda_{\sigma}+1})> k+3P$. Montrons $
\mathcal Z_{\sigma}^{\lambda_{\sigma}+1}\subset  \bigcup_{y\in S_{m}}2P\tg(x,y)$. Soit $z\in 
\mathcal Z_{\sigma}^{\lambda_{\sigma}+1}$ . 
Il existe $a\in S_{0}$ tel que $z\in 2F\tg(x,a)$.
Cela est clair si $\sigma=0$. Si $\sigma=1$ il existe 
 $b\in S_{1}$ tel que $z\in \geod(x,b)$. Soit $a\in S_{0}$. Alors 
 le a) du lemme~\ref{lemme-S0-...Sm} montre que $b\in 2F\tg(x,a)$, d'où 
$z\in 2F\tg(x,a)$. 

Soit $y\in S_{m}$ tel que $d(x,y)=d_{\max}(x,S_{m})$. On a $y\in 2F\tg(x,a)$ d'après le a) du lemme~\ref{lemme-S0-...Sm}. D'autre part $d(x,z)\leq d(x,y)+P/3$ puisque 
$t_{\sigma}^{\lambda_{\sigma}+1}(Z)\leq d_{\max}(x,S_{m})$ et $Z_{\sigma}^{\lambda_{\sigma}+1}$ est de diamètre inférieur ou égal à $P/3$. En appliquant le lemme~\ref{x-a,b-y-geod}  à $(x,a,y,z)$ au lieu de $(x,y,a,b) $ et $(2F, 2F,P/3)$ au lieu de $(\alpha,\beta,\rho)$ on trouve 
$$z\in (2F+2P/3+\de)\tg(x,y)\subset 2P\tg(x,y)$$ 
car on suppose $2F+2P/3+\de \leq 2P$, ce qui est permis par $(H_{P})$. 
On a donc montré $
\mathcal Z_{\sigma}^{\lambda_{\sigma}+1}\subset  \bigcup_{y\in S_{m}}2P\tg(x,y)$ et $\mathcal Y_{m}^{l_{m}+1}=\mathcal Z_{\sigma}^{\lambda_{\sigma}+1}$ vérifie la condition iv) de la définition~\ref{defi-Y}. 

Supposons maintenant $$t_{\sigma}^{\lambda_{\sigma}+1}(Z)\in ]r_{\sigma,i+1}^{\max}(Z), r_{\sigma,i}^{\max}(Z)]\text{\ \  avec \ \ }i\in \{\sigma,...,m-1\}.$$ 
Alors nécessairement $r_{\sigma,i+1}^{\max}(Z)< r_{\sigma,i}^{\max}(Z)$, ce qui implique 
$$r_{\sigma,i}^{\max}(Z)=d_{\max}(x,S_{i})\text{\ \  et\ \ }r_{\sigma,i+1}^{\max}(Z)=\max( k+3P,d_{\max}(x,S_{i+1})) .$$ On a donc 
$t_{\sigma}^{\lambda_{\sigma}+1}(Z)\in ]d_{\max}(x,S_{i+1}),d_{\max}(x,S_{i})]$
Soit $z\in \mathcal Z_{\sigma}^{\lambda_{\sigma}+1}$. 
On a vu qu'il existe  $a\in S_{0}$ tel que $z\in 2F\tg(x,a)$. Soit $y\in S_{i}$ tel que 
$d(x,y)=d_{\max}(x,S_{i})$ et soit $y'\in S_{i+1}$. D'après  le a) du lemme~\ref{lemme-S0-...Sm}, $y$ et $y'$ appartiennent à  $ 2F\tg(x,a)$. Il est clair que $d(x,y')\leq d(x,z)\leq d(x,y)+P/3$.

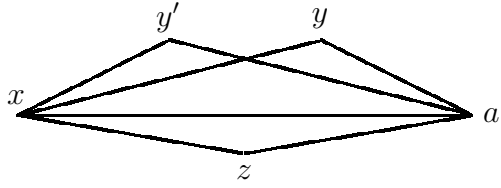
\ifx\JPicScale\undefined\def\JPicScale{1}\fi
\unitlength \JPicScale mm
\begin{picture}(80,30)(20,20)
\linethickness{0.3mm}
\put(40,30){\line(1,0){60}}
\linethickness{0.3mm}
\multiput(80,40)(0.24,-0.12){83}{\line(1,0){0.24}}
\linethickness{0.3mm}
\multiput(40,30)(0.48,0.12){83}{\line(1,0){0.48}}
\linethickness{0.3mm}
\multiput(40,30)(0.24,0.12){83}{\line(1,0){0.24}}
\linethickness{0.3mm}
\multiput(60,40)(0.48,-0.12){83}{\line(1,0){0.48}}
\linethickness{0.3mm}
\multiput(70,25)(0.71,0.12){42}{\line(1,0){0.71}}
\linethickness{0.3mm}
\multiput(40,30)(0.71,-0.12){42}{\line(1,0){0.71}}
\put(40,32.5){\makebox(0,0)[cc]{$x$}}

\put(60,43){\makebox(0,0)[cc]{$y'$}}

\put(80,43){\makebox(0,0)[cc]{$y$}}

\put(102.5,30){\makebox(0,0)[cc]{$a$}}

\put(70,22.5){\makebox(0,0)[cc]{$z$}}

\end{picture}

\noindent
Le lemme~\ref{x-a,b-y-geod} appliqué à $(x,a,z,y')$ au lieu de $(x,y,a,b)$ et $(2F,2F,0)$ au lieu de $(\alpha,\beta,\rho)$ montre 
\begin{gather}\label{geod-29dec8000}y'\in (2F+\de)\tg(x,z).\end{gather} Le lemme~\ref{x-a,b-y-geod} appliqué à $(x,a,y,z)$ au lieu de $(x,y,a,b)$ et $(2F,2F,P/3)$ au lieu de $(\alpha,\beta,\rho)$ montre 
\begin{gather}\label{geod-29dec8001}z\in (2P/3+2F+\de)\tg(x,y).\end{gather} Par le b) du lemme~\ref{geod-comp-xabc} on  déduit de (\ref{geod-29dec8000}) et (\ref{geod-29dec8001}) que 
$$z\in (2P/3+4F+2\de)\tg(y',y)\subset P\tg(y',y)$$
car on suppose $2P/3+4F+2\de \leq P$, ce qui est permis par $(H_{P})$. On a donc montré $
\mathcal Z_{\sigma}^{\lambda_{\sigma}+1}\subset  \bigcup_{y\in S_{i},y'\in S_{i+1}}P\tg(y,y')$ et $\mathcal Y_{i}^{l_{i}+1}=\mathcal Z_{\sigma}^{\lambda_{\sigma}+1}$ vérifie la condition iii) de la définition~\ref{defi-Y}. 
 \cqfd

Soit $\sigma\in\{0,1\}$ (avec $\sigma=0$ si $m=0$). On pose 
  \begin{itemize}
  \item 
  $(\tilde \lambda_{0},\tilde \lambda_{1})= (\lambda_{0}+1 ,\lambda_{1})$ si $
  \sigma=0$ 
  \item et $(\tilde \lambda_{0},\tilde \lambda_{1})= (\lambda_{0} ,\lambda_{1}+1)$ si $
  \sigma=1$.
  \end{itemize}

 On note $\kappa_{\sigma}:Y_{x}^{\natural,p,k,m,(l_{0},...,l_{m}), \tilde \lambda_{0},\tilde \lambda_{1}} \to Y_{x}^{\natural,p,k,m,(l_{0},...,l_{m}),\lambda_{0},\lambda_{1}}$ $$\cup \bigcup _{i\in \{\sigma,...,m\}} Y_{x}^{\natural,p,k,m,(l_{0},...,l_{i-1},l_{i}+1,l_{i+1},...,l_{m}),\lambda_{0},\lambda_{1}}$$  l'application définie par le sous-lemme~\ref{remplacement-lambda0+1}. On vérifie facilement que $\kappa_{\sigma}$ passe au quotient et définit \label{kappa-sigma}
 \begin{gather*}\overline \kappa_{\sigma,\infty}: \{Z\in \overline Y_{x}^{\natural,p,k,m,(l_{0},...,l_{m}), \tilde \lambda_{0},\tilde \lambda_{1}}, 0\leq t_{\sigma}^{\lambda_{\sigma}+1}(Z) \leq k+3P\}\to \overline Y_{x}^{\natural,p,k,m,(l_{0},...,l_{m}), \lambda_{0},\lambda_{1}}\\ \text{et \ }\overline \kappa_{\sigma,i}: \{Z\in \overline Y_{x}^{\natural,p,k,m,(l_{0},...,l_{m}), \tilde \lambda_{0},\tilde \lambda_{1}},  r_{\sigma,i+1}^{\max}(Z)< t_{\sigma}^{\lambda_{\sigma}+1}(Z) \leq  r_{\sigma,i}^{\max}(Z)\}\\ \to \overline Y_{x}^{\natural,p,k,m,(l_{0},...,l_{i-1},l_{i}+1,l_{i+1},...,l_{m}), \lambda_{0}, \lambda_{1}}\text{ pour } i\in \{\sigma,...,m\}. \end{gather*}

\begin{souslem}\label{prop-kappa-23oct09}
Il existe $C=C(\de,K,N,Q,P)$ tel que 

\noindent a) pour tout $i\in \{\sigma,...,m\}$, $\overline \kappa_{\sigma,i}$ est injective,
 et pour $Z\in \overline Y_{x}^{\natural,p,k,m,(l_{0},...,l_{m}), \tilde \lambda_{0},\tilde \lambda_{1}}$ tel que $r_{\sigma,i+1}^{\max}(Z)< t_{\sigma}^{\lambda_{\sigma}+1}(Z) \leq  r_{\sigma,i}^{\max}(Z)$, $\kappa_{\sigma}$ 
  induit une bijection 
  de $$(\pi_{x}^{\natural,p,k,m,(l_{0},...,l_{m}), \tilde \lambda_{0},\tilde \lambda_{1}})^{-1}(Z)$$ dans $$(\pi _{x}^{\natural,p,k,m,(l_{0},...,l_{i-1},l_{i}+1,l_{i+1},...,l_{m}), \lambda_{0}, \lambda_{1}})^{-1}(\overline \kappa_{\sigma,i}(Z)),$$

\noindent b) pour tout élément $Z_{\infty}\in  \overline Y_{x}^{\natural,p,k,m,(l_{0},...,l_{m}), \lambda_{0},\lambda_{1}}$,  $\overline  \kappa_{\sigma,\infty}^{-1}(Z_{\infty})$ est de cardinal inférieur ou égal à $C(r_{\sigma,m+1}^{\max}(Z_{\infty})+1)$
et pour tout entier $t\leq k+3P$ le nombre de $Z\in \overline Y_{x}^{\natural,p,k,m,(l_{0},...,l_{m}), \tilde \lambda_{0},\tilde \lambda_{1}}$ tels que $t_{\sigma}^{\lambda_{\sigma}+1}(Z)=t$
 et $\overline \kappa_{\sigma,\infty }(Z)=Z_{\infty}$ est inférieur ou égal à $C$. De plus pour tout 
  $Z\in \overline Y_{x}^{\natural,p,k,m,(l_{0},...,l_{m}), \tilde \lambda_{0},\tilde \lambda_{1}}$ tel que $t_{\sigma}^{\lambda_{\sigma}+1}(Z) \leq  k+3P$, 
$\kappa_{\sigma}$ induit 
  une bijection 
  de $(\pi_{x}^{\natural,p,k,m,(l_{0},...,l_{m}), \tilde \lambda_{0},\tilde \lambda_{1}})^{-1}(Z)$ dans 
$ (\pi _{x}^{\natural,p,k,m,(l_{0},...,l_{m}), \lambda_{0},\lambda_{1}})^{-1}(\overline \kappa_{\sigma,\infty }(Z))$. 
\end{souslem}
\noindent{\bf Démonstration.}
La preuve de a) est immédiate. Montrons  b). 
Soit $Z\in \overline Y_{x}^{\natural,p,k,m,(l_{0},...,l_{m}), \tilde \lambda_{0},\tilde \lambda_{1}}$ tel que $t_{\sigma}^{\lambda_{\sigma}+1}(Z) \leq  k+3P$. 
Pour tout $$(a_{1},\dots,a_{p},S_{0},...,S_{m},(\mathcal Y_{i}^{j})_{i\in \{0,\dots,m\}, j\in \{1,\dots ,l_{i}\}}, (\mathcal Z_{i}^{j})_{i\in \{0,1\}, j\in \{1,\dots ,\tilde \lambda_{i}\}})$$ dans $
(\pi_{x}^{\natural,p,k,m,(l_{0},...,l_{m}),\tilde \lambda_{0},\tilde \lambda_{1}})^{-1}(Z)$ on a  $\mathcal Z_{\sigma}^{\lambda_{\sigma}+1}\subset B(x,k+3P+P/3)$ donc 
$\mathcal Z_{\sigma}^{\lambda_{\sigma}+1}\subset B(x,k+M)$ car on suppose $3P+P/3\leq M$, ce qui est permis par $(H_{M})$,  et  donc $B(\mathcal Z_{\sigma}^{\lambda_{\sigma}+1}, M)\subset B(x,k+2M)$. Par conséquent, si on note $Z_{\infty}=\overline \kappa_{\sigma,\infty }(Z)$, la donnée de $Z$ est équivalente à celle de $Z_{\infty}$ et de $\mathcal Z_{\sigma}^{\lambda_{\sigma}+1}\subset B(x,k+2M)$. 
La condition $\mathcal Z_{\sigma}^{\lambda_{\sigma}+1}\subset \bigcup_{a\in S_{\sigma}}\geod(x,a)$ 
s'exprime uniquement en termes des distances entre les points de $S_{\sigma}$ et ceux de $B(x,k+2M)$ (qui font partie de la donnée de $Z_{\infty}$) et le nombre de possibilités pour $\mathcal Z_{\sigma}^{\lambda_{\sigma}+1}\subset B(x,k+2M)$ de diamètre $\leq P/3$ vérifiant cette condition et $d(x,\mathcal Z_{\sigma}^{\lambda_{\sigma}+1})\leq k+3P$  est borné par $C(r_{\sigma,m+1}^{\max}(Z_{\infty})+1)$ avec $C=C(\de,N,K,Q,P)$ à cause 
du  lemme~\ref{cardinal-tranche-geod}. Pour tout $t$,  l'ensemble des $\mathcal Z_{\sigma}^{\lambda_{\sigma}+1}$  vérifiant les conditions précédentes et $d(x,\mathcal Z_{\sigma}^{\lambda_{\sigma}+1})=t$ ne dépend que de $Z_{\infty}$ et de $t$  et le lemme~\ref{cardinal-tranche-geod} montre que le cardinal de cet ensemble   est inférieur ou égal à $C=C(\de,N,K,Q,P)$.     \cqfd

 \noindent{\bf Suite de la démonstration du lemme~\ref{eq-normes-natural}.} Grâce au sous-lemme~\ref{prop-kappa-23oct09},  il existe  $C=C(\de,K,N,Q,P,M)$ tel que pour $k,m,l_{0},\dots,$ $l_{m},$ $\lambda_{0},\lambda_{1}\in \N$, $\sigma\in\{0,1\}$ (vérifiant $\lambda_{1}=0$ et $\sigma=0$ si $m=0$), 
et en posant  
  \begin{itemize}
  \item 
  $(\tilde \lambda_{0},\tilde \lambda_{1})= (\lambda_{0}+1 ,\lambda_{1})$ si $
  \sigma=0$ 
  \item et $(\tilde \lambda_{0},\tilde \lambda_{1})= (\lambda_{0} ,\lambda_{1}+1)$ si $
  \sigma=1$
    \end{itemize}
    on ait l'inégalité suivante :  
 \begin{gather}\nonumber 
B^{-(m+\sum_{i=0}^{m}l_{i})}
\sum_{Z\in  \overline Y_{x}^{\natural,p,k,m,(l_{0},...,l_{m}),\tilde \lambda_{0},\tilde \lambda_{1}}}  e^{2s(r_{0}(Z)-k)}
\Big(\prod_{i=0}^{m}s_{i}(Z)^{-l_{i}} \Big)
 (r_{0}(Z)+1)^{-\tilde \lambda_{0}}
\\ \nonumber
(r_{1}(Z)+1)^{-\tilde \lambda_{1}}
\sharp \big((\pi_{x}^{\natural,p,k,m,(l_{0},...,l_{m}),\tilde \lambda_{0},\tilde  \lambda_{1}})^{-1}(Z)\big)^{-\alpha }
\big|\xi_{Z}(f)\big|^{2}
%\Big|\sum _{\substack{(a_{1},\dots,a_{p},S_{1},...,S_{m},(\mathcal Y_{i}^{j})_{i\in \{0,\dots,m\}, j\in \{1,\dots ,l_{i}\}}, (\mathcal Z_{i}^{j})_{i\in \{0,1\}, j\in \{1,\dots ,\tilde \lambda_{i}\}}) \\\in 
%(\pi_{x}^{\natural,p,k,m,(l_{0},...,l_{m}),\tilde \lambda_{0},\tilde \lambda_{1}})^{-1}(Z)}} f(a_1,...,a_p)\Big|^{2}
\\ \nonumber  
\leq  C
\Bigg(B^{-(m+\sum_{i=0}^{m}l_{i})}
\sum_{Z\in  \overline Y_{x}^{\natural,p,k,m,(l_{0},...,l_{m}),\lambda_{0},\lambda_{1}}}  \Big(\frac{r_{\sigma,m+1}^{\max}(Z)+1}{r_{\sigma,\sigma}^{\max}(Z)+1}\Big)e^{2s(r_{0}(Z)-k)}
\\ \nonumber 
 \Big(\prod_{i=0}^{m}s_{i}(Z)^{-l_{i}} \Big)
(r_{0}(Z)+1)^{- \lambda_{0}}(r_{1}(Z)+1)^{- \lambda_{1}}
\\ \nonumber 
\sharp \big((\pi_{x}^{\natural,p,k,m,(l_{0},...,l_{m}), \lambda_{0},\lambda_{1}})^{-1}(Z)\big)^{-\alpha }
\big|\xi_{Z}(f)\big|^{2}
% \Big|\sum _{\substack{(a_{1},\dots,a_{p},S_{1},...,S_{m},(\mathcal Y_{i}^{j})_{i\in \{0,\dots,m\}, j\in \{1,\dots ,l_{i}\}}, (\mathcal Z_{i}^{j})_{i\in \{0,1\}, j\in \{1,\dots ,\lambda_{i}\}}) \\ \in 
%(\pi_{x}^{\natural,p,k,m,(l_{0},...,l_{m}), \lambda_{0},\lambda_{1}})^{-1}(Z)}}f(a_1,...,a_p)\Big|^{2}
\Bigg)
 \\ \nonumber  +CB \sum_{i=\sigma}^{m}\Bigg(B^{-(m+\sum_{j=0}^{m}l_{j}+1)}
\sum_{Z\in  \overline Y_{x}^{\natural,p,k,m,(l_{0},...,l_{i-1},l_{i}+1, l_{i+1},...,l_{m}), \lambda_{0}, \lambda_{1}}} 
\\ \nonumber
 \Big(\frac{r_{\sigma,i}^{\max}(Z)-r_{\sigma,i+1}^{\max}(Z)}{r_{\sigma,\sigma}^{\max}(Z)+1}\Big)
 e^{2s(r_{0}(Z)-k)}
\Big(s_{i}(Z)^{-(l_{i}+1)} \prod_{j\in \{0,...,m\}\setminus \{i\}}s_{j}(Z)^{-l_{j}} \Big)
\\ \nonumber 
(r_{0}(Z)+1)^{- \lambda_{0}}(r_{1}(Z)+1)^{- \lambda_{1}}
\\ \label{eq-29dec8006} 
\sharp \big((\pi_{x}^{\natural,p,k,m,(l_{0},...,l_{i-1},l_{i}+1, l_{i+1},...,l_{m}), \lambda_{0},  \lambda_{1}})^{-1}(Z)\big)^{-\alpha }
\big|\xi_{Z}(f)\big|^{2}
% \Big|\sum _{\substack{(a_{1},\dots,a_{p},S_{1},...,S_{m},(\mathcal Y_{i}^{j})_{i\in \{0,\dots,m\}, j\in \{1,\dots ,\tilde l_{i}\}}, (\mathcal Z_{i}^{j})_{i\in \{0,1\}, j\in \{1,\dots , \lambda_{i}\}}) \\ \in 
%(\pi_{x}^{\natural,p,k,m,(\tilde l_{0},...,\tilde l_{m}),\lambda_{0}, \lambda_{1}})^{-1}(Z)}} f(a_1,...,a_p)\Big|^{2}
\Bigg). 
\end{gather}
En effet le  b) du lemme~\ref{prop-r-i-max}
assure que $r_{\sigma,\sigma}^{\max}(Z)+1\leq (N+1)(r_{\sigma}(Z)+1)$
d'où $(r_{\sigma}(Z)+1)^{-1}\leq (N+1)(r_{\sigma,\sigma}^{\max}(Z)+1)^{-1}$. De plus 
le $i^{\text{ème}}$ terme dans la somme qui constitue la deuxième moitié du membre de droite n'apparaît pas si $ r_{\sigma,i}^{\max}(Z)=r_{\sigma,i+1}^{\max}(Z)$ et si $ r_{\sigma,i}^{\max}(Z)>r_{\sigma,i+1}^{\max}(Z)$ on a $$s_{i}(Z)\leq C(r_{\sigma,i}^{\max}(Z)-r_{\sigma,i+1}^{\max}(Z))$$ avec $C=C(\de,K,N,Q,P)$.

%On a un énoncé similaire au lemme~\ref{remplacement-lambda0+1} avec 
% $\tilde \lambda_{0}= \lambda_{0} $ et 
% $\tilde \lambda_{1}= \lambda_{1}+1$
% (c'est $r_{1,1}^{\max}(Z)+1$ qui apparaît au dénominateur au lieu de  $r_{0,0}^{\max}(Z)+1$ et les $(r_{1,j}^{\max}(Z))_{j\in \{1,...,m+1\}}$ remplacent les $(r_{0,j}^{\max}(Z))_{j\in \{0,...,m+1\}}$). On en déduit alors  le lemme~\ref{eq-normes-natural} en appliquant le lemme suivant.
 
 \begin{souslem}\label{general-A-lambda0-lambda1}
 Soient  $C\in \R_{+}^{*}$,  $m\in \N$, $\epsilon_{0},...,\epsilon_{m},\eta_{1},...,\eta_{m}\in [0,1]$ vérifiant $\sum_{i=0}^{m}\epsilon_{i}\leq 1$ et $\sum_{i=1}^{m}\eta_{i}\leq 1$. Soient  $\mu_{0},\mu_{1}\in \N$  et  $(A_{l_{0},...,l_{m},\lambda_{0},\lambda_{1}})$ une famille d'éléments de $\R_{+}$ indexée par les $(l_{0},...,l_{m},\lambda_{0},\lambda_{1})\in \N^{m+3}$ vérifiant  $\lambda_{0}\leq \mu_{0}$ et $\lambda_{1}\leq \mu_{1}$. On suppose  \begin{itemize}
 \item pour $\lambda_{0}\in \{0,...,\mu_{0}-1\}$ et $\lambda_{1}\in \{0,...,\mu_{1}\}$, 
 $$A_{l_{0},...,l_{m},\lambda_{0}+1,\lambda_{1}}\leq C\big(A_{l_{0},...,l_{m},\lambda_{0},\lambda_{1}}+\sum_{i=0}^{m}\epsilon_{i}A_{l_{0},...,l_{i-1},l_{i}+1,l_{i+1},...,l_{m},\lambda_{0},\lambda_{1}}\big),
 $$
\item et pour $\lambda_{0}\in \{0,...,\mu_{0}\}$ et $\lambda_{1}\in \{0,...,\mu_{1}-1\}$, 
 $$A_{l_{0},...,l_{m},\lambda_{0},\lambda_{1}+1}\leq C\big(A_{l_{0},...,l_{m},\lambda_{0},\lambda_{1}}+\sum_{i=1}^{m}\eta_{i}A_{l_{0},...,l_{i-1},l_{i}+1,l_{i+1},...,l_{m},\lambda_{0},\lambda_{1}}\big).
 $$
  \end{itemize}
 Alors 
 \begin{gather}\label{concl-lemA-11j}\sum_{l_{0},...,l_{m}\in \N,\lambda_{0}\in \{0,...,\mu_{0}\},\lambda_{1}\in \{0,...,\mu_{1}\}}A_{l_{0},...,l_{m},\lambda_{0},\lambda_{1}}\leq (2C+1)^{\mu_{0}+\mu_{1}} \sum_{l_{0},...,l_{m}\in \N}A_{l_{0},...,l_{m},0,0}\end{gather}
  où l'on sous-entend que si le membre de droite converge, le membre de gauche converge aussi. 
 \end{souslem}
  \noindent{\bf Démonstration.}
 Posons pour $\lambda_{0}\in \{0,...,\mu_{0}\}$ et $\lambda_{1}\in \{0,...,\mu_{1}\}$, 
 $$\mathcal A_{\lambda_{0},\lambda_{1}}=\sum_{l_{0},...,l_{m}}A_{l_{0},...,l_{m},\lambda_{0},\lambda_{1}}.$$
 Alors pour $\lambda_{0}\in \{0,...,\mu_{0}-1\}$ et $\lambda_{1}\in \{0,...,\mu_{1}\}$, on a $$\mathcal A_{\lambda_{0}+1,\lambda_{1}} \leq 2C
 \mathcal A_{\lambda_{0},\lambda_{1}}$$ et de même pour $\lambda_{0}\in \{0,...,\mu_{0}\}$ et $\lambda_{1}\in \{0,...,\mu_{1}-1\}$, $$\mathcal A_{\lambda_{0},\lambda_{1}+1} \leq 2C\mathcal A_{\lambda_{0},\lambda_{1}}.$$ Le sous-lemme en résulte facilement car 
 $1+(2C)+...(2C)^{\mu_{c}}\leq (2C+1)^{\mu_{c}}$ pour $c=0,1$.  
      \cqfd
      
 \noindent{\bf Fin de  la démonstration du lemme~\ref{eq-normes-natural}.} 
 On applique le sous-lemme~\ref{general-A-lambda0-lambda1} de la fa\c con suivante. On fixe $m,k$ et des entiers $R_{0,0}\geq ...\geq R_{0,m+1}\geq 0$ et 
 $R_{1,1}\geq ... \geq R_{1,m+1}\geq 0$. On applique le sous-lemme~\ref{general-A-lambda0-lambda1} en prenant  $$\epsilon_{j}=\frac{R_{0,j}-R_{0,j+1}}{R_{0,0}+1}\text{\  pour\  }j=0,...,m,  \ \ \eta_{j}=\frac{R_{1,j}-R_{1,j+1}}{R_{1,1}+1}\text{\  pour \ }j=1,...,m$$ et, pour $l_{0},...,l_{m}\in \N$, $\lambda_{0}\in \{0,...,\mu_{0}\}$ et $\lambda_{1}\in \{0,...,\mu_{1}\}$, 
 \begin{gather*}A_{l_{0},...,l_{m},\lambda_{0},\lambda_{1}}=0\text{\ \  si \ \  }m=0\text{\  et \  }\lambda_{1}>0\text{,\ et\  sinon} \\ A_{l_{0},...,l_{m},\lambda_{0},\lambda_{1}}=B^{-(m+\sum_{i=0}^{m}l_{i})}
 \sum_{\substack{Z\in \overline Y_{x}^{\natural,p,k,m,(l_{0},...,l_{m}),\lambda_{0},\lambda_{1}}\text{ tel que } \\ r_{i,j}^{\max}(Z)=R_{i,j}\text{ pour } i\in \{0,1\}\text{ et }j\in \{i,...,m+1\}}}\\ e^{2s(r_{0}(Z)-k)}
 \Big(\prod_{i=0}^{m}s_{i}(Z)^{-l_{i}} \Big)
(r_{0}(Z)+1)^{- \lambda_{0}}(r_{1}(Z)+1)^{- \lambda_{1}}\\  
\sharp \big((\pi_{x}^{\natural,p,k,m,(l_{0},...,l_{m}), \lambda_{0},\lambda_{1}})^{-1}(Z)\big)^{-\alpha }
\big| \xi_{Z}(f) \big|^{2}
 %\Big|\sum _{\substack{(a_{1},\dots,a_{p},S_{1},...,S_{m},(\mathcal Y_{i}^{j})_{i\in \{0,\dots,m\}, j\in \{1,\dots ,l_{i}\}}, (\mathcal Z_{i}^{j})_{i\in \{0,1\}, j\in \{1,\dots ,\lambda_{i}\}}) \\ \in 
%\pi_{x}^{\natural,p,k,m,(l_{0},...,l_{m}), \lambda_{0},\lambda_{1}})^{-1}(Z)}}f(a_1,...,a_p)\Big|^{2}
.\end{gather*}
 Grâce à (\ref{eq-29dec8006}) et comme $\frac{r_{\sigma,m+1}^{\max}(Z)+1}{r_{\sigma,\sigma}^{\max}(Z)+1}\leq 1$,  les hypothèses du sous-lemme~\ref{general-A-lambda0-lambda1} sont satisfaites pour une constante $C=C(\de,K,N,Q,P,M,s,B)$. 
     Puis on somme l'inégalité (\ref{concl-lemA-11j}) sur $$m,k, (R_{0,0}, ...,R_{0,m+1}),(R_{1,1}, ... , R_{1,m+1}). $$
   Ceci termine la démonstration du lemme~\ref{eq-normes-natural}. \cqfd
 
 Le lemme suivant est une variante du lemme~\ref{nombre-dist-connaitre-par-point}. 
 
 \begin{lem}\label{nombre-dist-connaitre-par-point-natural}
Il existe une constante $C=C(\de,K,N,Q,P,M)$ telle que pour $p\in \{1,...,p_{\max}\}$, $k,m,l_{0},...,l_{m},\lambda_{0},\lambda_{1}\in \N$  (vérifiant $\lambda_{1}=0$ si $m=0$) et 
\begin{gather*}(a_{1},\dots,a_{p},S_{0},...,S_{m},(\mathcal Y_{i}^{j})_{i\in \{0,\dots,m\}, j\in \{1,\dots ,l_{i}\}}, (\mathcal Z_{i}^{j})_{i\in \{0,1\}, j\in \{1,\dots ,\lambda_{i}\}}) \\ \in Y_{x}^{\natural,p,k,m,(l_{0},...,l_{m}),\lambda_{0},\lambda_{1}},\end{gather*} les distances entre les points de 
\begin{gather} \label{29dec8013}
B(S_{0},M)
\cup   \bigcup _{ j\in \{1,\dots ,l_{0}\}} B(\mathcal Y_{0}^{j}, M) \cup \bigcup _{i\in \{0,1\}, j\in \{1,\dots ,\lambda_{i}\}}
B(\mathcal Z_{i}^{j},M)\end{gather} et ceux de 
 \begin{gather}\label{29dec8014} 
\bigcup _{ i\in \{1,\dots ,m\}}
B(S_{i},M)
\cup   \bigcup _{i\in \{1,\dots,m\}, j\in \{1,\dots ,l_{i}\}} B(\mathcal Y_{i}^{j},M) \cup B(x,k+2M)\end{gather}
 sont déterminées par 
 \begin{itemize}
 \item a)
  les distances entre les points de (\ref{29dec8014}), 
  \item b)  les  entiers $d(x,S_{0})$, $d(x,\mathcal Y_{0}^{j})$, $d(x,\mathcal Z_{0}^{j})$ et $d(x,\mathcal Z_{1}^{j})$, 
  \item c)  les distances
entre les points de (\ref{29dec8013}) et $C(1+l_{0}+\lambda_{0}+\lambda_{1})$ points de 
(\ref{29dec8014}) (qui sont eux-mêmes déterminés par   a) et b))
\end{itemize}
 et de plus  les distances entre les points de (\ref{29dec8013}) et ceux de 
 (\ref{29dec8014})  sont déterminées à $C$ près par a) et b). 
\end{lem}
\noindent{\bf Remarque. } Dans toutes les situations où on appliquera ce lemme, $l_{0}+\lambda_{0}+\lambda_{1}$ sera majoré par une constante de la forme $C(\de,K,N,Q,P,M)$. 

 \noindent{\bf Démonstration.} 
 Soient 
 $p\in \{1,...,p_{\max}\}$, $k,m,l_{0},...,l_{m},\lambda_{0},\lambda_{1}\in \N$  (vérifiant $\lambda_{1}=0$ si $m=0$) et 
\begin{gather*}(a_{1},\dots,a_{p},S_{0},...,S_{m},(\mathcal Y_{i}^{j})_{i\in \{0,\dots,m\}, j\in \{1,\dots ,l_{i}\}}, (\mathcal Z_{i}^{j})_{i\in \{0,1\}, j\in \{1,\dots ,\lambda_{i}\}}) \\ \in Y_{x}^{\natural,p,k,m,(l_{0},...,l_{m}),\lambda_{0},\lambda_{1}}. \end{gather*} 
  Soit $b\in S_{0}$ et $u$ un point de $B(x,k)$ à distance minimale de $b$.  
  On commence par montrer que pour $\sigma\in \{0,1\}$ et $j\in \{1,...,\lambda_{\sigma}\}$, on a 
  \begin{gather}
  \label{incl-29dec0832} 
  \mathcal Z_{\sigma}^{j}\subset 2F\tg(x,b) \end{gather} et \begin{gather}\!\!\!\!\!\!
  \label{incl-29dec8026}
 B(\mathcal Z_{\sigma}^{j},M) \subset B(x,k+2M) \text{ \ ou \ } \mathcal Z_{\sigma}^{j}\subset (2F+\de)\tg(u,b). 
 \end{gather}
 
 Soit $j\in \{1,...,\lambda_{0}\}$.  On a $$\mathcal Z_{0}^{j}\subset \bigcup_{a\in S_{0}}\geod(x,a)\subset 2N\tg(x,b)\subset 2F\tg(x,b)$$ car $F\geq N$. Si $d(x,
 \mathcal Z_{0}^{j})\leq k+3P$, $$B(\mathcal Z_{0}^{j}, M)\subset B(x,k+2M)$$ car on suppose $3P+P/3\leq M$, ce qui est permis par $(H_{M})$. 
 Si $d(x,
 \mathcal Z_{0}^{j})> k+3P$, 
 pour $z\in  \mathcal Z_{0}^{j}$ on a  
 \begin{gather}\label{cond-z-0-4.28}z\in 2N\tg(x,b)\text{\  et\  }d(x,z)>  d(x,u)+3P. \end{gather}
 
 Soit $j\in \{1,...,\lambda_{1}\}$.  On a $\mathcal Z_{1}^{j}\subset \bigcup_{a\in S_{1}}\geod(x,a)$.  Par le lemme~\ref{lemme-S0-...Sm}, $S_{1}\subset 2F\tg(u,b)$, donc   $\mathcal Z_{1}^{j} \subset  
 2F\tg(x,b)$. On a déjà prouvé (\ref{incl-29dec0832}).  Si $d(x,
 \mathcal Z_{1}^{j})\leq k+3P$, $$B(\mathcal Z_{1}^{j}, M)\subset B(x,k+2M)$$ car  $3P+P/3\leq M$. 
 Si $d(x,
 \mathcal Z_{1}^{j})> k+3P$, 
 soit  $z\in  \mathcal Z_{1}^{j}$. 
 On a donc   \begin{gather}\label{cond-z-1-4.28}z\in 
 2F\tg(x,b)\text{\  et\  }d(x,z)> d(x,u)+3P.\end{gather}
 
 Soit maintenant $z$ vérifiant (\ref{cond-z-0-4.28}) ou (\ref{cond-z-1-4.28}). 
 Alors $z$ vérifie (\ref{cond-z-1-4.28}) car $F\geq N$ par (\ref{def-F}). 
 
 \ifx\JPicScale\undefined\def\JPicScale{1}\fi
\unitlength \JPicScale mm
\begin{picture}(100,30)(25,20)
\linethickness{0.3mm}
\put(50,30){\line(1,0){80}}
\linethickness{0.3mm}
\multiput(92.5,42.5)(0.36,-0.12){104}{\line(1,0){0.36}}
\linethickness{0.3mm}
\multiput(50,30)(0.41,0.12){104}{\line(1,0){0.41}}
\linethickness{0.3mm}
\multiput(72.5,30)(0.19,0.12){104}{\line(1,0){0.19}}
\put(47,32.5){\makebox(0,0)[cc]{$x$}}

\put(92.5,45){\makebox(0,0)[cc]{$z$}}

\put(72.5,27.5){\makebox(0,0)[cc]{$u$}}

\put(132.5,32.5){\makebox(0,0)[cc]{$b$}}

\end{picture}

\noindent Par $(H_{\de}
^{2F}(u,x,z,b))$ on a 
$$d(u,z)\leq \max(d(u,x)-d(x,z),d(u,b)-d(b,z))+2F+\de.$$
Or $d(u,x)-d(x,z)+2F+\de\leq -3P+2F+\de$ et on suppose $-3P+2F+\de<0$, ce qui est permis par 
 $(H_{P})$. 
 Donc $d(u,z)\leq d(u,b)-d(b,z)+2F+\de$, c'est-à-dire 
    $z\in (2F+\de)\tg(u,b)$. Ceci termine la preuve de (\ref{incl-29dec8026}). 
    On suppose $2F+\de\leq 4P$, ce qui est permis par $(H_{P})$. 
    
      Pour montrer le lemme~\ref{nombre-dist-connaitre-par-point-natural} on répète alors les arguments de la preuve du lemme~\ref{nombre-dist-connaitre-par-point}. 
   
Plus précisément on applique le lemme~\ref{i<j<k-distances}   avec 
\begin{itemize}
\item $(x,b)$ au lieu de $(c,d)$, 
\item $\{w_{i},i\in I\}$ égal à $$\bigcup_{i\in \{0,\dots,m\}}S_{i}\cup  \bigcup _{i\in \{0,\dots,m\}, j\in \{1,\dots ,l_{i}\}} \mathcal Y_{i}^{j}\cup \bigcup _{i\in \{0,1\}, j\in \{1,\dots ,\lambda_{i}\}}
\mathcal Z_{i}^{j},$$
\item $J$ la partie de $I$ telle que $$\{w_{j},j\in J\}=S_{0}\cup  \bigcup _{ j\in \{1,\dots ,l_{0}\}} \mathcal Y_{0}^{j}\cup \bigcup _{i\in \{0,1\}, j\in \{1,\dots ,\lambda_{i}\}}
\mathcal Z_{i}^{j},$$
%\item $\{w_{1},...,w_{r}\}$ égal à la réunion disjointe de $S_{1},...,S_{m}$,  des $\mathcal Y_{i}^{j}$ pour $i\in \{1,\dots,m\}$ et $ j\in \{1,\dots ,l_{i}\}$, et des  $  \mathcal Z_{\sigma}^{j}$ pour $\sigma\in \{0,1\}$ et $j\in \{1,...,\lambda_{\sigma}\}$, 
\item  et $(\alpha_{i},\rho_{i})$ égal à $(4P,M)$ pour tout $i$
(les hypothèses sont satisfaites grâce au lemme~\ref{lemme-S0-...Sm} et à (\ref{incl-29dec0832}) et car on a supposé $2F\leq 4P$). 
\end{itemize}
  Puis on  applique le lemme~\ref{distances-Bxk-C} à $z$ parcourant 
  $$B(S_{0},M)\cup   \bigcup _{ j\in \{1,\dots ,l_{0}\}} B(\mathcal Y_{0}^{j}, M)
  \cup \bigcup_{\substack{i\in \{0,1\},j\in \{1,\dots ,\lambda_{i}\} \\
 B(\mathcal Z_{i}^{j},M)\not \subset B(x,k+2M) }}B(\mathcal Z_{i}^{j},M)
    $$ %et des 
%$B(\mathcal Z_{i}^{j},M)$ pour 
%$i\in \{0,1\}$ et $ j\in \{1,\dots ,\lambda_{i}\}$  tels que 
%$B(\mathcal Z_{i}^{j},M)\not \subset B(x,k+2M)$, 
avec 
$\alpha=4P+2M$ et $l=k+2M$ (grâce  au lemme~\ref{lemme-S0-...Sm} et à (\ref{incl-29dec8026}) et comme on a supposé $2F+\de\leq 4P$, cet ensemble  est inclus dans $(4P+2M)\tg(b,u)$). \cqfd
 
 Le lemme suivant nous sera utile ensuite. 
 \begin{lem}\label{lemme-cardinaux}
 Il existe une constante $C=C(\de,K,N,Q,P)$ telle que pour $m,l_{0},...,l_{m},\lambda_{0},\lambda_{1}\in \N$ (avec $\lambda_{1}=0$ si $m=0$) et pour 
 $Z\in \overline Y_{x}^{p,k,m,(l_{0},..., l_{m})}$ et $\tilde Z\in 
 \overline Y_{x}^{\natural, p,k,m,(l_{0},..., l_{m}), \lambda_{0},\lambda_{1}}$ tels que pour tout \begin{gather*}(a_{1},\dots,a_{p},S_{0},...,S_{m},(\mathcal Y_{i}^{j})_{i\in \{0,\dots,m\}, j\in \{1,\dots ,l_{i}\}}, (\mathcal Z_{i}^{j})_{i\in \{0,1\}, j\in \{1,\dots , \lambda_{i}\}})\\  \in 
(\pi_{x}^{\natural, p,k,m,(l_{0},..., l_{m}), \lambda_{0},\lambda_{1}})^{-1}(\tilde Z)\end{gather*} on ait 
$$(a_{1},\dots,a_{p},S_{0},...,S_{m},(\mathcal Y_{i}^{j})_{i\in \{0,\dots,m\}, j\in \{1,\dots ,l_{i}\}})\in (\pi_{x}^{p,k,m,(l_{0},..., l_{m})})^{-1}(Z)$$
 alors 
 $\sharp\big((\pi_{x}^{\natural, p,k,m,(l_{0},..., l_{m}), \lambda_{0},\lambda_{1}})^{-1}(\tilde Z)\big)\leq C^{\lambda_{0}+\lambda_{1}}\sharp\big((\pi_{x}^{p,k,m,(l_{0},..., l_{m})})^{-1}(Z)\big)$.
 \end{lem}
 \noindent{\bf Démonstration.} 
Soit $$(a_{1},\dots,a_{p},S_{0},...,S_{m},(\mathcal Y_{i}^{j})_{i\in \{0,\dots,m\}, j\in \{1,\dots ,l_{i}\}})\in (\pi_{x}^{p,k,m,(l_{0},..., l_{m})})^{-1}(Z)$$
et soient $b_{0}\in S_{0},b_{1}\in S_{1}$. Si $(\mathcal Z_{i}^{j})_{i\in \{0,1\}, j\in \{1,\dots , \lambda_{i}\}}$ sont tels que 
\begin{gather*}(a_{1},\dots,a_{p},S_{0},...,S_{m},(\mathcal Y_{i}^{j})_{i\in \{0,\dots,m\}, j\in \{1,\dots ,l_{i}\}}, (\mathcal Z_{i}^{j})_{i\in \{0,1\}, j\in \{1,\dots ,\lambda_{i}\}})\\  \in 
(\pi_{x}^{\natural, p,k,m,(l_{0},..., l_{m}), \lambda_{0},\lambda_{1}})^{-1}(\tilde Z)\end{gather*}  on a 
$\mathcal Z_{i}^{j}\subset 2N\tg(x,b_{i})$ pour $i\in \{0,1\}, j\in \{1,\dots , \lambda_{i}\}$. De plus la donnée de $\tilde Z$ détermine $t_{i}^{j}(\tilde Z)=d(x,\mathcal Z_{i}^{j})$ et le diamètre de $\mathcal Z_{i}^{j}$ doit être inférieur ou égal à $P/3$. On applique alors le lemme~\ref{cardinal-tranche-geod}. \cqfd

 Voici quelques rappels et notations pour la proposition suivante. On a 
 $$J_x=\tilde H_x+u_xK_{x}
 =\tilde H_x+\sum_{r=1}^{+\infty} u_{x,r}K_{x}
 $$ en notant  $u_{x,r}=\int_{0}^{1}u_{x,r,t}dt$ (de sorte que $u_{x}=\sum_{r=1}^{+\infty} u_{x,r}$). 
 Pour $q\in \{1,...,Q\}$ on note $\tilde H_{x,q}=h_{x}(1-\del h_{x}-h_{x}\del)^{q-1}$, de sorte que 
 $$\tilde H_{x,q}=
\int_{(t_{1},\dots ,t_{q})\in [0,1]^{q}}\tilde H_{x,q,(t_{1},\dots ,t_{q})} dt_{1}\dots dt_{q}\text{\ \ \ et\ \ \ } \tilde H_{x}=\sum_{q=1}^{Q}\tilde H_{x,q}. 
$$
 On rappelle aussi que $K_{x}=\int_{(t_{1},\dots ,t_{Q})\in [0,1]^{Q}}  K_{x,Q,(t_{1},\dots ,t_{Q})}  dt_{1}\dots dt_{Q}$.

\begin{prop}\label{continuite-Jx}
 Pour tout $p\in \{1,\pp,p_{\max}\}$, $J_{x}$ se prolonge en un opérateur continu de $\H_{x,s}(\Delta_{p-1})$ dans $\H_{x,s}(\Delta_{p})$. Plus précisément, pour $p\in \{2,\pp,p_{\max}\}$ il existe $C=C(\de,K,N,Q,P,M,s,B)$ tel que pour tout $q\in \{1,...,Q\}$, 
 $$\|\tilde H_{x,q}\|_{\L(\H_{x,s}(\Delta_{p-1}),\H_{x,s}(\Delta_{p}))}\leq C$$ et pour tout $r\in \N$, 
 $$\|u_{x,r} K_x\|_{\L(\H_{x,s}(\Delta_{p-1}),\H_{x,s}(\Delta_{p}))}\leq 
 Ce^{-\frac{s}{2}r}.$$
  \end{prop}
 On remarque que $\tilde H_{x,1}=h_{x}$ et grâce à la proposition~\ref{continuite-del} la continuité de $h_{x}$ implique celle de $ \tilde H_{x,q}$ pour tout $q\in \{1,...,Q\}$. Cependant nous préférons montrer directement la continuité de 
  $ \tilde H_{x,q}$ pour tout $q\in \{1,...,Q\}$ car cela prépare à la démonstration de la continuité de $u_{x,r} K_x$. 
  
  \noindent {\bf Démonstration. }
Le cas où $p=1$ est  trivial. 

Soit $p\in \{2,\pp,p_{\max}\}$ et $q\in \{1,...,Q\}$. Comme dans la démonstration de la proposition~\ref{continuite-del}, on note  $\P$ le projecteur orthogonal sur le sous-espace vectoriel de $\H^{\rightarrow}_{x,s}(\Delta_{p})$ engendré par les $e_{S}$ pour $S\in \Delta_{p}$ tel que $d(x,S)\leq P$, de sorte que $(\P f)(S)=f(S)$ si 
$d(x,S)\leq P$ et $(\P f)(S)=0$ sinon. La proposition~\ref{continuite-Jx} résulte donc des lemmes~\ref{sl1-1jan0923}, \ref{sl2-1jan0923}, \ref{sl3-1jan0923} et \ref{sl4-1jan0923} que nous allons montrer successivement. \cqfd

\begin{lem}\label{sl1-1jan0923}
Il existe $C=C(\de,K,N,Q,P,M,s,B)$ tel que $$
\|\P \tilde H_{x,q} \|_{\L(\H_{x,s}(\Delta_{p-1}),\H^{\rightarrow}_{x,s}(\Delta_{p}))} \leq C.$$
\end{lem}
%On commence par montrer  que $$
%\|\P \tilde H_{x,q} \|_{\L(\H_{x,s}(\Delta_{p-1}),\H^{\rightarrow}_{x,s}(\Delta_{p}))} \leq C\text{\ \ et \ } 
%\|\P u_{x,r} K_x\|_{\L(\H_{x,s}(\Delta_{p-1}),\H^{\rightarrow}_{x,s}(\Delta_{p}))} \leq Ce^{-\frac{s}{2}r}$$ pour une certaine constante $C=C(\de,K,N,Q,P,M,s,B)$. 
\noindent{\bf Démonstration.} 
Il suffit  de montrer que pour tout $U\in \Delta_{p}$ vérifiant $d(x,U)\leq P$, il existe 
$C=C(\de,K,N,Q,P,M,s,B)$ tel que pour $f\in \C^{(\Delta_{p})}$, 
\begin{gather}\label{cont-Jx-partie1-1}|(\tilde H_{x,q}  (f))(U)| \leq C\|f\|_{\H_{x,s}(\Delta_{p-1})}
.\end{gather}
Soit $U\in \Delta_{p}$ vérifiant $d(x,U)\leq P$. 
L'inégalité (\ref{cont-Jx-partie1-1})  est évidente, car, d'après le 1)a) de la proposition~\ref{recap-supp-connaiss-H-uK}, 
pour tout $S\in \Delta_{p-1}$, $\tilde H_{x,q}  (e_{S})$ est supporté par les 
$T\in \Delta_{p}$ tels que 
$T\subset   \bigcup_{a\in S} B(a,QN)$, donc $(\tilde H_{x,q}  (e_{S}))(U)$ est nul sauf si $d(x,S)\leq QN+P$ et le nombre de telles parties $S$ est majoré par une constante $C=C(\de,K,N,Q,P)$. De plus  pour une telle partie $S$, $|(\tilde H_{x,q}  (e_{S}))(U)|$ est également majoré par une telle constante par le 3) de la proposition~\ref{recap-supp-connaiss-H-uK}. On conclut en utilisant le lemme~\ref{minoration-normeHsx}.  \cqfd

\begin{lem}\label{sl2-1jan0923}
Il existe $C=C(\de,K,N,Q,P,M,s,B)$ tel que, pour tout $r\in \N$, 
 $$
\|\P u_{x,r} K_x\|_{\L(\H_{x,s}(\Delta_{p-1}),\H^{\rightarrow}_{x,s}(\Delta_{p}))} \leq Ce^{-\frac{s}{2}r}.$$
\end{lem}
\noindent{\bf Démonstration.}  
Il suffit  de montrer que pour tout $U\in \Delta_{p}$ vérifiant $d(x,U)\leq P$, il existe 
$C=C(\de,K,N,Q,P,M,s,B)$ tel que, pour $r\in \N$ et  $f\in \C^{(\Delta_{p})}$, 
\begin{gather}\label{cont-Jx-partie1-2}
|(u_{x,r} K_x (f))(U)| \leq Ce^{-\frac{s}{2}r}\|f\|_{\H_{x,s}(\Delta_{p-1})}
.\end{gather}
Soit $U\in \Delta_{p}$ vérifiant $d(x,U)\leq P$ et $r\in \N$. Nous allons montrer (\ref{cont-Jx-partie1-2}). Soit 
$t,t_{1},\pp,t_{Q}\in [0,1]$. 
D'après le 2)a) de la proposition~\ref{recap-supp-connaiss-H-uK}, pour $S\in \Delta_{p-1}$, $u_{x,r,t} K_{x,Q,(t_{1},\dots ,t_{Q})}  (e_{S}) $ est une combinaison de $e_{T}$ pour $T$ vérifiant 
$$d(x,T)\in [d(x,S)-r-QF,d(x,S)-r+N+F-\frac{Q}{F}].$$ On suppose $\frac{Q}{F}\geq F+N$, ce qui est permis par $(H_{Q})$.  Pour que 
$$(u_{x,r,t} K_{x,Q,(t_{1},\dots ,t_{Q})} (e_{S}))(U)$$ soit non nul il est donc nécessaire que 
\begin{gather}\label{condition-r-Q4-24oct09}d(x,S)\in [r,r+QF+P].\end{gather}

On note  $\Lambda_{t_{1},...,t_{Q}}$  la partie de $\overline Y_{x}^{\natural,p-1,0,0,(0),Q,0}$ formée des $Z$ tels que 
\begin{itemize}
\item $r_{0}(Z)\in [r,r+QF+P]$, 
\item pour tout $(a_{1},\dots,a_{p-1},S_{0},(\mathcal Z_{0}^{j})_{ j\in \{1,\dots ,Q\}})   \in 
(\pi_{x}^{\natural,p-1,0,0,(0),Q,0})^{-1}(Z)$, et pour tout $j\in \{1,\pp,Q\}$, on a 
\begin{gather}\label{def-29dec2021}\mathcal Z_{0}^{j}=\bigcup_{b\in S_{0}}\{z\in \geod(x,b), d(x,z)=E(t_{j}r_{0}(Z))\}\end{gather} 
(on rappelle que  $r_{0}(Z)=d(x,S_{0})$).   
 \end{itemize}
 La  condition (\ref{def-29dec2021}) implique que pour $Z\in \Lambda_{t_{1},...,t_{Q}}$
 et $j\in \{1,\pp,Q\}$ on a $t_{0}^{j}(Z)=E(t_{j}r_{0}(Z))$.  

\begin{souslem}\label{slem0-2eme-2j1317}
    Soit  $( a_{1},\dots, a_{p-1},S_{0})\in Y_{x}^{p-1,0,0,(0)}$
  tel que $d(x, S_{0})\in [r,r+QF+P]$.  
  Pour $j\in \{1,...,Q\}$ on définit ${\mathcal Z}_{0}^{j}$ par (\ref{def-29dec2021}). Alors ${\mathcal Z}_{0}^{j}$ est de diamètre inférieur ou égal à $P/3$  et  il existe $Z\in \Lambda_{t_{1},...,t_{Q}}$ tel que 
$$( a_{1},\dots, a_{p-1},S_{0},( {\mathcal Z}_{0}^{j})_{ j\in \{1,\dots ,q\}}
)\in (\pi_{x}^{\natural,p-1,0,0,(0),Q,0})^{-1}(Z).$$ \end{souslem}
\noindent{\bf Démonstration.}
Soient $j\in \{1,...,Q\}$ et $z,z'\in \mathcal Z_{0}^{j}$. Soit  $b\in  S_{0}$. On a $z,z'\in 2N\tg(x,b)$ et $d(x,z)=d(x,z')$ donc par $(H_{\de}(z,x,z',b))$, $d(z,z')\leq 2N+\de$ et on suppose $2N+\de\leq P/3$, ce qui est permis  par $(H_{P})$.  Comme les parties $\mathcal Z_{0}^{j}$ sont non vides   l'argument que nous venons de donner montre aussi que la condition (\ref{def-29dec2021}) est vérifiée par les autres éléments de la classe d'équivalence $Z$ de $( a_{1},\dots, a_{p-1},S_{0},( {\mathcal Z}_{0}^{j})_{ j\in \{1,\dots ,q\}}
)$ (car on suppose $P/3\leq M$, ce qui est permis par $(H_{M})$)
et donc que $Z\in  \Lambda_{t_{1},...,t_{Q}}$. \cqfd

\begin{souslem}\label{slem1-2eme-2j1323}
Soit 
$Z\in \Lambda_{t_{1},...,t_{Q}}$ et $$(a_{1},\dots,a_{p-1},S_{0},(\mathcal Z_{0}^{j})_{ j\in \{1,\dots ,Q\}})   \in 
(\pi_{x}^{\natural,p-1,0,0,(0),Q,0})^{-1}(Z). $$ Alors $u_{x,r,t} K_{x,Q,(t_{1},\dots ,t_{Q})}  (e_{a_{1}}\wedge ...\wedge e_{a_{p-1}})$ ne dépend que de la connaissance des points de  
\begin{gather} \label{slem-2eme-ens-3j1557}B(S_{0},M)\cup B(x,2M)\cup  \bigcup _{ j\in \{1,\dots ,Q\}}B(\mathcal Z_{0}^{j}, M) \end{gather}
et des distances entre ces points. 
\end{souslem}
\noindent{\bf Démonstration.}
Le 2)b) de la proposition~\ref{recap-supp-connaiss-H-uK} montre que 
$$u_{x,r,t} K_{x,Q,(t_{1},\dots ,t_{Q})}  (e_{a_{1}}\wedge ...\wedge e_{a_{p-1}})$$ dépend seulement 
de la connaissance des points de 
\begin{gather}\label{ens1-3j1602}
B(x,F)\cup B(S_{0},QN) \cup \bigcup_{a\in S_{0}}\{y\in F\tg(x,a), d(y,a)\in [r,r+QF]\}
\\  \label{ens2-3j1602}\cup
\bigcup_{a\in S_{0}}\{y\in F\tg (x,a), |d(x,y)-(1-t)(d(x,a)-r)|\leq QF \}
\\ \label{ens3-3j1602}
\cup \bigcup_{a\in S_{0},j\in \{1,...,Q\}}\{y\in F\text{-}\geod(x,a), |d(x,y)-t_{j}d(x,a)|\leq QF \} 
\end{gather}
et des distances entre ces points. Il suffit donc de montrer que cet ensemble est inclus dans (\ref{slem-2eme-ens-3j1557}).

On a  $B(x,F)\subset B(x,2M)$ et  
$B(S_{0},QN) \subset 
B(S_{0}, M)$
car on suppose $F\leq 2M$ et $QN\leq M$, ce qui est permis par 
 $(H_{M})$. 
 Comme $d(x,S_{0})\in [r,r+QF+P]$  on a 
\begin{gather*} 
\bigcup_{a\in S_{0}}\{y\in F\tg(x,a), d(y,a)\in [r,r+QF]\}\\ \subset  B(x, d(x,S_{0})+N+F-r)
\subset B(x, QF+P+N+F)
\subset B(x,2M)\end{gather*}
car on suppose $QF+P+N+F\leq 2M$, ce qui est permis par 
 $(H_{M})$. 
 Donc (\ref{ens1-3j1602}) est inclus dans (\ref{slem-2eme-ens-3j1557}). 
 
  Comme $d(x,S_{0})\in [r,r+QF+P]$  on a 
  \begin{gather*} \bigcup_{a\in S_{0}}\{y\in F\tg (x,a), |d(x,y)-(1-t)(d(x,a)-r)|\leq QF \} \\
\subset B(x,N+2QF+P)\subset B(x,2M)\end{gather*}
car on suppose $N+2QF+P\leq 2M$, ce qui est permis par 
 $(H_{M})$. Donc (\ref{ens2-3j1602}) est inclus dans (\ref{slem-2eme-ens-3j1557}).

%Il résulte de ce qui précède que l'ensemble figurant dans le 2)b) de la proposition~\ref{recap-supp-connaiss-H-uK} est inclus dans 
% $$B(S_{0}, M)\cup B(x,2M)
% \cup 
%\bigcup_{a\in S_{0},j\in \{1,...,Q\}}\{y\in F\text{-}\geod(x,a), |d(x,y)-t_{j}d(x,a)|\leq QF \}.$$

Enfin, soit $j\in \{1,...,Q\}$,  $a\in S_{0}$, et 
 $y\in F\text{-}\geod(x,a)$ vérifiant  $|d(x,y)-t_{j}d(x,a)|\leq QF$. 
 Soit  $z\in \geod(x,a)$ vérifiant $d(x,z)=E(t_{j}d(x,S_{0}))$, 
 si bien que $z$ appartient à $\mathcal Z_{0}^{j}$. 
 Comme 
 $y,z\in F\tg(x,a)$ et $|d(x,y)-d(x,z)|\leq QF+N+1$, 
 $(H_{\de}(y,x,z,a))$ montre que 
 $d(y,z)\leq  (QF+N+1)+F+\de$. On suppose $(QF+N+1)+F+\de\leq M$, ce qui est permis par $(H_{M})$. Donc (\ref{ens3-3j1602}) est inclus dans 
 $\bigcup_{j\in \{1,...,Q\}}B(\mathcal Z_{0}^{j}, M)$ et a fortiori dans 
 (\ref{slem-2eme-ens-3j1557}). 
 %
 %
 %  On a donc, pour tout $j\in \{1,\pp,Q\}$, 
%\begin{gather*}\bigcup_{a\in S_{0}}\{y\in F\text{-}\geod(x,a), |d(x,y)-t_{j}d(x,a)|\leq QF\} 
%\subset    B(\mathcal Z_{0}^{j}, M). 
%\end{gather*}
%Donc, sous la condition~(\ref{condition-r-Q4-24oct09}), l'ensemble figurant dans le 2)b) de la proposition~\ref{recap-supp-connaiss-H-uK} (qui a la propriété que  que $u_{x,r,t} K_{x,Q,(t_{1},\dots ,t_{Q})}  (e_{S_{0}})$ ne dépend que de la connaissance des points de de cet ensemble) est inclus dans 
%$$B(S_{0},M)\cup B(x,2M)\cup  \bigcup _{ j\in \{1,\dots ,Q\}}B(\mathcal Z_{0}^{j}, M) $$ où les $\mathcal Z_{0}^{j}$ sont définis par (\ref{def-29dec2021}). 
\cqfd

\begin{souslem}\label{slem2-2eme-2j1323}
Le cardinal de $\Lambda_{t_{1},...,t_{Q}}$ est majoré par une constante de la forme $C(\de,K,N,Q,P,M)$. 
\end{souslem}
\noindent{\bf Démonstration.} Cela résulte du lemme~\ref{nombre-dist-connaitre-par-point-natural} (ou même d'un argument plus simple car le cardinal de l'ensemble (\ref{slem-2eme-ens-3j1557}) est borné par $C=C(\de,K,N,Q,P,M)$ et les distances entre les points de (\ref{slem-2eme-ens-3j1557}) sont déterminées à $C'=C(\de,K,N,Q,P,M)$ près par la donnée de $r,t_{1},...,t_{Q}$). 
\cqfd

\noindent{\bf Fin de la démonstration du lemme~\ref{sl2-1jan0923}.}
On écrit $U=\{b_{1},...,b_{p}\}$ pour lever l'ambiguïté de signe. On rappelle, pour $Z\in \Lambda_{t_{1},...,t_{Q}} $,  la notation 
$$\xi_{Z}(f)=\sum _{(a_{1},\dots,a_{p-1},S_{0},(\mathcal Z_{0}^{j})_{ j\in \{1,\dots ,Q\}})  \in 
(\pi_{x}^{\natural,p-1,0,0,(0),Q,0})^{-1}(Z)} f(a_1,...,a_{p-1}).$$ 
On a $f=\frac{1}{(p-1)!}\sum_{(a_{1},...,a_{p-1})} f(a_{1},...,a_{p-1}) e_{a_{1}}\wedge ... \wedge e_{a_{p-1}}$ où la somme porte que les $(a_{1},...,a_{p-1})$ tel que $\{a_{1},...,a_{p-1}\}\in \Delta_{p-1}$. 
Le sous-lemme~\ref{slem0-2eme-2j1317} montre donc que 
\begin{gather}\label{eg-slem0-2eme-11j}(u_{x,r} K_{x,Q,(t_{1},\dots ,t_{Q})}  (f))(b_{1},...,b_{p})=\frac{1}{(p-1)!}\sum_{Z\in \Lambda_{t_{1},...,t_{Q}} }\alpha_{Z,(t_{1},...,t_{Q}),(b_{1},...,b_{p})}\xi_{Z}(f)\end{gather} 
 où $\alpha_{Z,(t_{1},...,t_{Q}),(b_{1},...,b_{p})}\in \C$ est défini de la fa\c con suivante : 
 \begin{gather*}\text{pour tout\ \ }(a_{1},\dots,a_{p-1},S_{0},(\mathcal Z_{0}^{j})_{ j\in \{1,\dots ,Q\}})  \in 
(\pi_{x}^{\natural,p-1,0,0,(0),Q,0})^{-1}(Z), \\ \alpha_{Z,(t_{1},...,t_{Q}),(b_{1},...,b_{p})}=\big(u_{x,r} K_{x,Q,(t_{1},\dots ,t_{Q})} (e_{a_{1}}\wedge ...\wedge e_{a_{p-1}})\big)(b_{1},...,b_{p})\end{gather*} 
(d'après le sous-lemme~\ref{slem1-2eme-2j1323} ce nombre  ne dépend que de $Z$). 
D'après le 3) de la proposition~\ref{recap-supp-connaiss-H-uK},  $|\alpha_{Z,(t_{1},...,t_{Q}),(b_{1},...,b_{p})}|$
est majoré par une constante de la forme $C(\de,K,N,Q)$. 
Par Cauchy-Schwarz et grâce au sous-lemme~\ref{slem2-2eme-2j1323}, 
on a donc $$|(u_{x,r} K_{x,Q,(t_{1},\dots ,t_{Q})}  (f))(U)|^{2}\leq C\sum_{Z\in \Lambda_{t_{1},...,t_{Q}} } |\xi_{Z}(f)|^{2}$$ pour une certaine constante $C=C(\de,K,N,Q,P,M)$.

On en déduit, pour $t_{1},\pp,t_{Q}\in [0,1]$,  
 \begin{gather*} |(u_{x,r} K_{x,Q,(t_{1},\dots ,t_{Q})}  (f))(U)|^{2} \\ \leq Ce^{-sr}
\sum_{Z\in 
%\overline Y_{x}^{\natural,p-1,0,0,(0),Q,0}, \forall j\in \{1,\pp,Q\}, t_{0}^{j}(Z)=E(t_{j}r_{0}(Z))
\Lambda_{t_{1},...,t_{Q}} }   e^{2sr_{0}(Z)}
\sharp \big((\pi_{x}^{\natural,p-1,0,0,(0),Q,0})^{-1}(Z)\big)^{-\alpha }
|\xi_{Z}(f)|^{2}
%\Big|\sum _{\substack{(a_{1},\dots,a_{p-1},(\mathcal Z_{0}^{j})_{ j\in \{1,\dots ,Q\}})\\  \in 
%(\pi_{x}^{\natural,p-1,0,0,(0),Q,0})^{-1}(Z)}} f(a_1,...,a_{p-1})\Big|^{2}
 \end{gather*}
pour une certaine constante $C=C(\de,K,N,Q,P,M)$. 
 En effet il existe une constante $D=C(\de,K,N,Q,P,M)$ telle que pour tout $Z\in \Lambda_{t_{1},...,t_{Q}}$, $$\sharp \big((\pi_{x}^{\natural,p-1,0,0,(0),Q,0})^{-1}(Z)\big)\leq e^{D(r_{0}(Z)+1)} $$ et on suppose $\alpha D\leq s$, ce qui est permis par $(H_{\alpha})$. De plus pour tout $Z\in \Lambda_{t_{1},...,t_{Q}}$ on a  
$r\leq r_{0}(Z)$, d'où $e^{-sr}\geq e^{-sr_{0}(Z)}$.  

On a vu que 
$\Lambda_{t_{1},...,t_{Q}}$  est inclus dans la partie de $\overline Y_{x}^{\natural,p-1,0,0,(0),Q,0}$ formée des $Z$ tels que 
\begin{itemize}
\item $r_{0}(Z)\in [r,r+QF+P]$, 
\item pour tout  $j\in \{1,\pp,Q\}$ on a $t_{0}^{j}(Z)=E(t_{j}r_{0}(Z))$.  
\end{itemize}
Pour $r_{0}\in \N$, quand  $(t_{1},\pp,t_{Q})$ parcourt  $[0,1]^{Q}$ muni de  la mesure de Lebesgue, 
 $(E(t_{j}r_{0}))_{j=1,\pp,Q}$ parcourt  $\big\{0,\pp,\max(0,r_{0}-1)\big\}^{Q}$ avec la probabilité uniforme $\max(1,r_{0})^{-Q}$ et comme $$  \max(1,r_{0})^{-1}\leq \frac{2}{r_{0}+1},$$ par Cauchy-Schwarz 
on obtient  l'inégalité 
\begin{gather} \nonumber  |(u_{x,r} K_x (f))(U)|^{2} \leq Ce^{-sr}\Big(
\sum_{Z\in \overline Y_{x}^{\natural,p-1,0,0,(0),Q,0}}  e^{2sr_{0}(Z)}\big(r_{0}(Z)+1\big)^{-Q} \\ \label{avant-derniere-(1-P)-del}
\sharp \big((\pi_{x}^{\natural,p-1,0,0,(0),Q,0})^{-1}(Z)\big)^{-\alpha }
|\xi_{Z}(f)|^{2}
%\Big|\sum _{\substack{(a_{1},\dots,a_{p-1},(\mathcal Z_{0}^{j})_{ j\in \{1,\dots ,Q\}})\\ \in 
%(\pi_{x}^{\natural,p-1,0,0,(0),Q,0})^{-1}(Z)}} f(a_1,...,a_{p-1})\Big|^{2}
 \Big)\end{gather}
pour une certaine constante $C=C(\de,K,N,Q,P,M)$. A fortiori on a 
$$|(u_{x,r} K_x (f))(U)|^{2} \leq Ce^{-sr}\|f\|_{\H^{\natural,Q,0}_{x,s}(\Delta_{p-1})}^{2}$$ puisque l'expression entre parenthèses dans  (\ref{avant-derniere-(1-P)-del}) est la partie  de la somme (\ref{formule-norme-mu}) donnant  
 $\|f\|_{\H^{\natural,Q,0}_{x,s}(\Delta_{p-1})}^{2}$  qui correspond à $\overline Y_{x}^{\natural,p-1,0,0,(0),Q,0}$. Grâce au lemme~\ref{eq-normes-natural} ceci termine la démonstration du lemme~\ref{sl2-1jan0923}. \cqfd

 \begin{lem}\label{sl3-1jan0923}
Il existe $C=C(\de,K,N,Q,P,M,s,B)$ tel que 
 $$\|(1-\P ) \tilde H_{x,q} \|_{\L(\H_{x,s}(\Delta_{p-1}),\H^{\rightarrow}_{x,s}(\Delta_{p}))} \leq C.$$
 \end{lem}
 \noindent{\bf Démonstration.} 
 Grâce au lemme~\ref{eq-normes-natural}, il suffit de montrer l'inégalité suivante : il existe $C=C(\de,K,N,Q,P,M,B)$ tel que pour tout $f\in \C^{(\Delta_{p-1})}$, 
\begin{gather}\label{ineg-1-P-H-natural}\|(1-\P ) (\tilde H_{x,q} f)\|^{2}_{\H^{\rightarrow}_{x,s}(\Delta_{p})}\leq 
C\|f\|^{2}_{\H_{x,s}^{\natural,q,0}(\Delta_{p-1})}.\end{gather}
On rappelle que  
\begin{gather*}\|(1-\P )(\tilde H_{x,q} f)\|_{\H^{\rightarrow}_{x,s}(\Delta_{p})}^{2} =\sum _{k,m,l_{0},\dots ,l_{m}\in \N} 
B^{-(m+\sum_{i=0}^{m}l_{i})}\sum_{\substack{Z\in \overline Y_{x}^{p,k,m,(l_{0},...,l_{m})}, \\ r_{0}(Z)> k+P}} \\ 
e^{2s(r_{0}(Z)-k)}\Big(\prod_{i=0}^{m}\big(s_{i}(Z)\big)^{-l_{i}} \Big)
\sharp \big((\pi_{x}^{p,k,m,(l_{0},...,l_{m})})^{-1}(Z)\big)^{-\alpha }\big|\xi_{Z}(\tilde H_{x,q} f)\big|^{2}.\end{gather*}
%$$
%\Big|\sum _{(a_{1},\dots,a_{p},S_{1},...,S_{m},(\mathcal Y_{i}^{j})_{i\in \{0,\dots,m\}, j\in \{1,\dots ,l_{i}\}}) \in 
%(\pi_{x}^{p,k,m,(l_{0},...,l_{m})})^{-1}(Z)} (\tilde H_{x,q} f)(a_1,...,a_p)\Big|^{2}.$$

 On va voir que l'inégalité (\ref{ineg-1-P-H-natural}) résulte de l'inégalité (\ref{ineg-1-P-H-natural-etape1}) ci-dessous. 
 
 Soient $k,m,l_{0},\dots ,l_{m}\in \N $ et $Z\in \overline Y_{x}^{p,k,m,(l_{0},...,l_{m})}$ vérifiant 
 $r_{0}(Z)> k+P$.  On pose $\tilde l_{0}=0$ et $ \tilde l_{i}=l_{i-1}$ pour $i\in \{1,\pp,m+1\}$. On va montrer 
  qu'il existe $C=C(\de,K,N,Q,P,M)$ tel que 
  \begin{gather}\label{ineg-1-P-H-natural-etape1}\Big| \xi_{Z}(\tilde H_{x,q}  f)\Big|^{2}
\leq C
\sum_{\tilde Z\in \Lambda_{Z}} \big(r_{0}(\tilde Z)+1\big)^{-q}
\Big|\xi_{\tilde Z}(f)\big|^{2}\end{gather}
 %$$ \Big| \sum _{\substack{ (\tilde a_{1},\dots,\tilde a_{p-1},\tilde S_{1},...,\tilde S_{m+1},(\tilde {\mathcal Y}_{i}^{j})_{i\in \{0,\dots,m+1\}, j\in \{1,\dots ,\tilde l_{i}\}},
%(\tilde {\mathcal Z}_{0}^{j})_{j\in \{1,\dots ,q\}})\\ \in 
%(\pi_{x}^{\natural,p-1,k,m+1,(\tilde l_{0},...,\tilde l_{m+1}),q,0})^{-1}(\tilde Z)}}  f(\tilde a_{1},\dots,\tilde a_{p-1})\Big|^{2}$$
 où $\Lambda_{Z}$ est la partie de  $\overline Y_{x}^{\natural,p-1,k,m+1,(\tilde l_{0},...,\tilde l_{m+1}),q,0}$ formée des $\tilde Z$ 
 vérifiant \begin{gather}\label{eq-29dec2048}|r_{0}(\tilde Z)-r_{1}(\tilde Z)|\leq (q+1)N\end{gather} et 
tels que pour tout  \begin{gather*}(\tilde a_{1},\dots,\tilde a_{p-1},\tilde S_{0}, ...,\tilde S_{m+1},(\tilde {\mathcal Y}_{i}^{j})_{i\in \{0,\dots,m+1\}, j\in \{1,\dots ,\tilde l_{i}\}},(\tilde {\mathcal Z}_{0}^{j})_{j\in \{1,\dots ,q\}})\\ \in (\pi_{x}^{\natural,p-1,k,m+1,(\tilde l_{0},...,\tilde l_{m+1}),q,0})^{-1}(\tilde Z)\end{gather*}
il existe une énumération $(a_{1},\pp,a_{p})$ de $\tilde S_{1}$ vérifiant 
%
% il existe $(a_{1},\pp,a_{p})$ (vérifiant 
% $\tilde S_{1}=\{a_{1},\pp,a_{p}\}$) tel que  
 \begin{gather}\label{notations-2j1240}(a_{1},\pp,a_{p},\tilde S_{1},...,\tilde S_{m+1},(\tilde {\mathcal Y}_{i+1}^{j})_{i\in \{0,\dots,m\}, j\in \{1,\dots , l_{i}\}})\in (\pi_{x}^{p,k,m,(l_{0},..., l_{m})})^{-1}(Z).\end{gather} 
 On rappelle que 
 $$ \xi_{\tilde Z}(f)= \sum _{\substack{ (\tilde a_{1},\dots,\tilde a_{p-1},\tilde S_{0},...,\tilde S_{m+1},(\tilde {\mathcal Y}_{i}^{j})_{i\in \{0,\dots,m+1\}, j\in \{1,\dots ,\tilde l_{i}\}},
(\tilde {\mathcal Z}_{0}^{j})_{j\in \{1,\dots ,q\}})\\ \in 
(\pi_{x}^{\natural,p-1,k,m+1,(\tilde l_{0},...,\tilde l_{m+1}),q,0})^{-1}(\tilde Z)}}  f(\tilde a_{1},\dots,\tilde a_{p-1}). $$

 On  justifie maintenant  le fait que (\ref{ineg-1-P-H-natural-etape1}) implique (\ref{ineg-1-P-H-natural}).
 D'abord $\tilde Z$ détermine $Z$ à permutation près de $a_{1},...,a_{p}$ donc connaissant $\tilde Z$ il y a au plus $p!$ possibilités pour $Z$. 
  Dans les notations précédentes, soit $b\in \tilde S_{0}$. D'après  (\ref{eq-29dec2048}),  pour tout $y\in \tilde S_{1}$ on a $|d(x,y)-d(x,b)|\leq (Q+2)N$ et on a $y\in 2F\tg(x,b)$ par le a) du lemme~\ref{lemme-S0-...Sm}, donc $d(y,b)\leq (Q+2)N+2F$. Connaissant $\tilde S_{1}$ on a donc au plus $C=C(\de,K,N,Q)$ possibilités pour $\tilde S_{0}$. En utilisant de plus le  lemme~\ref{lemme-cardinaux} on en déduit que pour $\tilde Z\in \Lambda_{Z}$ on a 
 $$\sharp\big((\pi_{x}^{\natural,p-1,k,m+1,(\tilde l_{0},...,\tilde l_{m+1}),q,0})^{-1}(\tilde Z)\big)\leq C\sharp\big((\pi_{x}^{p,k,m,(l_{0},..., l_{m})})^{-1}(Z)\big)$$ avec $C=C(\de,K,N,Q,P,M)$.
 Il est clair que pour $\tilde Z\in \Lambda_{Z}$ on a 
$ \prod_{i=0}^{m+1}s_{i}(\tilde Z)^{-\tilde l_{i}} =\prod_{i=0}^{m}s_{i}(Z)^{-l_{i}}$.    Donc (\ref{ineg-1-P-H-natural-etape1}) implique (\ref{ineg-1-P-H-natural}). 
 
 L'inégalité (\ref{ineg-1-P-H-natural-etape1}) résulte de 
 l'inégalité plus précise (\ref{ineg-(1-P)-tilde-l0}) ci-dessous. 
 
 Soient 
 $t_{1},\pp,t_{q}\in [0,1]$. On va montrer 
  qu'il existe $C=C(\de,K,N,Q,P,M)$ tel que 
  \begin{gather}\label{ineg-(1-P)-tilde-l0}
\Big|\xi_{Z} (\tilde H_{x,q,(t_{1},\dots ,t_{q})}  f)\Big|^{2}
\leq C
\sum_{\tilde Z\in \Lambda_{Z,(t_{1},...,t_{q})}}  \Big| \xi_{\tilde Z}( f)\Big|^{2}\end{gather}
 où $\Lambda_{Z,(t_{1},...,t_{q})}$ est l'ensemble des   $\tilde Z\in 
 \Lambda_{Z}$ 
% \overline Y_{x}^{\natural,p-1,k,m+1,(\tilde l_{0},...,\tilde l_{m+1}),q,0}$    
 %vérifiant 
 %$$|r_{0}(\tilde Z)-r_{1}(\tilde Z)|\leq (q+1)N$$ et  pour tout $j\in \{1,...,q\}$, 
%
tels que pour tout  \begin{gather*}(\tilde a_{1},\dots,\tilde a_{p-1},\tilde S_{0},...,\tilde S_{m+1},(\tilde {\mathcal Y}_{i}^{j})_{i\in \{0,\dots,m+1\}, j\in \{1,\dots ,\tilde l_{i}\}},(\tilde {\mathcal Z}_{0}^{j})_{ j\in \{1,\dots ,q\}}
)\\ \in (\pi_{x}^{\natural,p-1,k,m+1,(\tilde l_{0},...,\tilde l_{m+1}),q,0})^{-1}(\tilde Z)\end{gather*} et pour tout $j\in \{1,...,q\}$ on ait 
 \begin{gather}\label{def-mathcalZ-0j-24oct09}
\tilde{\mathcal Z}_{0}^{j}=
\bigcup_{b\in \tilde S_{0}}\{z\in \geod(x,b), d(x,z)=E(t_{j}r_{0}(\tilde Z))\}.\end{gather} 
%et   il existe $(a_{1},\pp,a_{p})$ vérifiant 
% $\tilde S_{1}=\{a_{1},\pp,a_{p}\}$ et 
% $$(a_{1},\pp,a_{p},\tilde S_{2},...,\tilde S_{m+1},(\tilde {\mathcal Y}_{i+1}^{j})_{i\in \{0,\dots,m\}, j\in \{1,\dots , l_{i}\}})\in (\pi_{x}^{p,k,m,(l_{0},..., l_{m})})^{-1}(Z).$$ 
La condition (\ref{def-mathcalZ-0j-24oct09}) implique  que pour  $\tilde Z\in \Lambda_{Z,(t_{1},...,t_{q})}$ et $j\in \{1,...,q\}$  on a 
 $t_{0}^{j}(\tilde Z) =E(t_{j}r_{0}(\tilde Z))$. 
  % On rappelle que $r_{0}(\tilde Z)=d(x,\tilde S_{0})$. 
  
 \begin{souslem}\label{slem-33-2j1245}
    Soit  $$(\tilde a_{1},\dots,\tilde a_{p-1},\tilde S_{0},...,\tilde S_{m+1},(\tilde {\mathcal Y}_{i}^{j})_{i\in \{0,\dots,m+1\}, j\in \{1,\dots ,\tilde l_{i}\}})\in Y_{x}^{p-1,k,m+1,(\tilde l_{0},...,\tilde l_{m+1})}$$
  tel que $|d(x,\tilde S_{0})-d(x,\tilde S_{1})|\leq (q+1)N$ et 
   qu'il existe $(a_{1},\pp,a_{p})$ vérifiant 
  (\ref{notations-2j1240}).
% $$(a_{1},\pp,a_{p},\tilde S_{2},...,\tilde S_{m+1},(\tilde {\mathcal Y}_{i+1}^{j})_{i\in \{0,\dots,m\}, j\in \{1,\dots , l_{i}\}})\in (\pi_{x}^{p,k,m,(l_{0},..., l_{m})})^{-1}(Z).$$ 
Pour $j\in \{1,...,q\}$ on définit $\tilde{\mathcal Z}_{0}^{j}$ par (\ref{def-mathcalZ-0j-24oct09}). Alors $\tilde{\mathcal Z}_{0}^{j}$ est de diamètre inférieur ou égal à $P/3$  et  il existe $\tilde Z\in \Lambda_{Z,(t_{1},...,t_{q})}$ tel que 
\begin{gather}\label{element-slem-33-2j1245}(\tilde a_{1},\dots,\tilde a_{p-1},\tilde S_{0},...,\tilde S_{m+1},(\tilde {\mathcal Y}_{i}^{j})_{i\in \{0,\dots,m+1\}, j\in \{1,\dots ,\tilde l_{i}\}},(\tilde {\mathcal Z}_{0}^{j})_{ j\in \{1,\dots ,q\}}
)\end{gather} appartienne à 
%$ Y_{x}^{\natural,p-1,k,m+1,(\tilde l_{0},...,\tilde l_{m+1}),q,0}$
$(\pi_{x}^{\natural,p-1,k,m+1,(\tilde l_{0},...,\tilde l_{m+1}),q,0})^{-1}(\tilde Z)$. 
%De plus la condition (\ref{def-mathcalZ-0j-24oct09}) ne dépend que de $\tilde Z$. 
\end{souslem}
\noindent{\bf Démonstration.}
Soient $j\in \{1,...,q\}$ et $z,z'\in \mathcal Z_{0}^{j}$. Soit  $b\in \tilde S_{0}$. On a $z,z'\in 2N\tg(x,b)$ et $d(x,z)=d(x,z')$ donc par $(H_{\de}(z,x,z',b))$, $d(z,z')\leq 2N+\de\leq P/3$. Comme les parties $\tilde{ \mathcal Z}_{0}^{j}$ sont  non vides et $P/3\leq M$, l'argument que nous venons de donner montre aussi que la condition (\ref{def-mathcalZ-0j-24oct09}) est vérifiée par les autres éléments de la classe d'équivalence  $\tilde Z$ de l'élément (\ref{element-slem-33-2j1245}) et donc $\tilde Z\in \Lambda_{Z,(t_{1},...,t_{q})}$. \cqfd
    
    \begin{souslem}\label{slem1-33-2j1248}
   Soit $\tilde Z\in \Lambda_{Z,(t_{1},...,t_{q})}$, et \begin{gather*} 
(\tilde a_{1},\dots,\tilde a_{p-1},\tilde S_{0},...,\tilde S_{m+1},(\tilde {\mathcal Y}_{i}^{j})_{i\in \{0,\dots,m+1\}, j\in \{1,\dots ,\tilde l_{i}\}},(\tilde {\mathcal Z}_{0}^{j})_{ j\in \{1,\dots ,q\}}
)\\ \in 
(\pi_{x}^{\natural, p-1,k,m+1,(\tilde l_{0},...,\tilde l_{m+1}),q,0})^{-1}(\tilde Z).
 \end{gather*}  
 Alors $\tilde H_{x,q,(t_{1},\dots ,t_{q})}(e_{\tilde a_{1}}\wedge ...\wedge e_{\tilde a_{p-1}})$ ne dépend que de la connaissance des points de   
\begin{gather}\label{ens-3eme-3j2013}B(\tilde S_{0},M)\cup B(x,k+2M)\cup  \bigcup _{ j\in \{1,\dots ,q\}}   B(\tilde{\mathcal Z}_{0}^{j},M) \end{gather} et des  distances entre ces points. 
\end{souslem}
\noindent{\bf Démonstration.}  D'après le 1)b) de la proposition~\ref{recap-supp-connaiss-H-uK}, $$\tilde H_{x,q,(t_{1},\dots ,t_{q})}(e_{\tilde a_{1}}\wedge ...\wedge e_{\tilde a_{p-1}})$$ ne dépend que de la connaissance des points de 
\begin{gather}\label{ens1-3eme-3j2013} B(x,7\de)\cup B(\tilde S_{0},QN) 
\\ \label{ens2-3eme-3j2013}
\cup 
\bigcup_{a\in \tilde S_{0},j\in \{1,...,q\}}\{y\in F\text{-}\geod(x,a), |d(x,y)-t_{j}d(x,a)|\leq QF\}\end{gather}
et des distances entre ces points. Il suffit donc de montrer que cet ensemble est inclus dans (\ref{ens-3eme-3j2013}). 
D'abord on suppose $7\de\leq 2M$ et $QN\leq M$, ce qui est permis par $(H_{M})$, et (\ref{ens1-3eme-3j2013})  est inclus dans (\ref{ens-3eme-3j2013}).

%
%$$B(x,2M)\cup B(\tilde S_{0}, M)
%$$ 
%
Soit  $a\in \tilde S_{0}$, $j\in \{1,...,q\}$, et $ y\in F\text{-}\geod(x,a)$  vérifiant $|d(x,y)-t_{j}d(x,a)|\leq QF$. Soit  $z\in \geod(x,a)$ vérifiant $d(x,z)=E(t_{j}r_{0}(\tilde Z))$, 
si bien que $z$ appartient à $\tilde{\mathcal Z}_{0}^{j}$. 
  On a $|t_{j}d(x,a)-E(t_{j}r_{0}(\tilde Z))|\leq 
N+1$ puisque 
$|d(x,a)-r_{0}(\tilde Z)|\leq N$,  d'où  
$|d(x,y)-d(x,z)|\leq QF+N+1$, et comme $y$ et $z$ appartiennent à $F\text{-}\geod(x,a)$, 
 $(H_{\de}(y,x,z,a))$ montre que  $$d(y,z)\leq (QF+N+1)+F+\de.$$ On suppose 
$(QF+N+1)+F+\de
\leq M$, ce qui est permis par $(H_{M})$. Donc (\ref{ens2-3eme-3j2013})  est inclus dans (\ref{ens-3eme-3j2013}).
%
% On a donc
%\begin{gather*}\bigcup_{a\in \tilde S_{0}, j\in \{1,\pp,q\}}\{y\in F\text{-}\geod(x,a), |d(x,y)-t_{j}d(x,a)|\leq QF\} 
%\subset  \bigcup_{j\in\{1,\pp,q\}}  B(\tilde{\mathcal Z}_{0}^{j}, M). \end{gather*}
%Donc   l'ensemble figurant dans  le 1)b) de la proposition~\ref{recap-supp-connaiss-H-uK} avec $\tilde S_{0}$ au lieu de $S$ (et qui a la propriété que $\tilde H_{x,q,(t_{1},\dots ,t_{q})}(e_{\tilde S_{0}})$ ne dépend que de la connaissance des distances entre ses points) est inclus dans  
%$$B(\tilde S_{0},M)\cup B(x,k+2M)\cup  \bigcup _{ j\in \{1,\dots ,q\}}   B(\tilde{\mathcal Z}_{0}^{j},M). $$
%Ceci termine la preuve du sous-lemme~\ref{slem1-33-2j1248}. 
\cqfd

 %On a toujours $|d(x,\tilde S_{0})- d(x,\tilde S_{1})|\leq (Q+1)N$.  
 
 \begin{souslem}\label{slem2-33-2j1252}
 Le cardinal de $ \Lambda_{Z,(t_{1},...,t_{q})}$ est majoré par une constante de la forme $C(\de,K,N,Q,P,M)$. 
  \end{souslem}
 \noindent{\bf Démonstration.}
 Les parties $\tilde {\mathcal Z}_{0}^{j}$ sont déterminées de manière unique par (\ref{def-mathcalZ-0j-24oct09}) et  grâce au lemme~\ref{nombre-dist-connaitre-par-point-natural}, pour connaître les distances entre les points de 
\begin{gather}\label{ens1-p102-24oct09}B(\tilde S_{0}, M)\cup   \bigcup _{ j\in \{1,\dots ,q\}}B(\tilde {\mathcal Z}_{0}^{j}, M)\end{gather} 
 et ceux de 
 \begin{gather}\label{ens2-p102-24oct09} 
\bigcup _{ i\in \{1,\dots ,m+1\}}
B(\tilde S_{i}, M)\cup  \bigcup _{i\in \{1,\dots,m+1\}, j\in \{1,\dots ,\tilde l_{i}\}}\B(\tilde {\mathcal Y}_{i}^{j}, M) \cup B(x,k+2M)\end{gather}
il suffit de connaître les distances entre les points de (\ref{ens1-p102-24oct09}) et $C$ points de (\ref{ens2-p102-24oct09}), avec $C=C(\de,K,N,Q,P,M)$ et grâce à (\ref{eq-29dec2048}) ces distances sont déterminées à $C'=C(\de,K,N,Q,P,M)$ près par les distances de $\tilde S_{1}$ à ces $C$ points (qui font partie de la donnée de $Z$) et les entiers $(t_{0}^{j}(\tilde Z))_{j\in \{1,\pp,q\}}$, qui grâce à (\ref{eq-29dec2048}) sont eux-mêmes déterminés à $C''=C(\de,K,N,Q,P,M)$ près par $r_{0}(Z),t_{1},...,t_{q}$. \cqfd 

\noindent{\bf  Fin de la démonstration du lemme~\ref{sl3-1jan0923}.}
On termine la démonstration de  (\ref{ineg-(1-P)-tilde-l0}). 
Pour 
$\tilde Z\in \Lambda_{Z,(t_{1},...,t_{q})}$ et 
$$ 
(\tilde a_{1},\dots,\tilde a_{p-1},\tilde S_{0},...,\tilde S_{m+1},(\tilde {\mathcal Y}_{i}^{j})_{i\in \{0,\dots,m+1\}, j\in \{1,\dots ,\tilde l_{i}\}},(\tilde {\mathcal Z}_{0}^{j})_{ j\in \{1,\dots ,q\}}
)$$ $$\in 
(\pi_{x}^{\natural, p-1,k,m+1,(\tilde l_{0},...,\tilde l_{m+1}),q,0})^{-1}(\tilde Z)
 $$  
on considère 
\begin{gather}
\label{somme-30dec0949}
\sum _{( b_{1},...,b_{p})}
(\tilde H_{x,q,(t_{1},\dots ,t_{q})} (e_{\tilde a_{1}}\wedge ...\wedge e_{\tilde a_{p-1}}))( b_{1},..., b_{p}), 
\end{gather} 
où la somme porte sur les énumérations $(b_{1},..., b_{p})$ de $\tilde S_{1}$ telles que 
$$( b_{1},\dots, b_{p},\tilde S_{1},...,\tilde S_{m+1},(\tilde {\mathcal Y}_{i+1}^{j})_{i\in \{0,\dots,m\}, j\in \{1,\dots ,\tilde l_{i}\}}
)\in (\pi_{x}^{p,k,m,(l_{0},..., l_{m})})^{-1}(Z). $$
Comme la somme (\ref{somme-30dec0949}) a au plus $p!$ termes, 
 le 3) de la proposition~\ref{recap-supp-connaiss-H-uK} montre qu'elle 
est majorée par une constante de la forme $C(\de,K,N,Q,P,M)$. 
D'après le sous-lemme~\ref{slem1-33-2j1248} la somme (\ref{somme-30dec0949}) ne dépend que de $\tilde Z$  et on peut donc la noter $\alpha_{Z,\tilde Z,(t_{1},...,t_{q})}$. 
D'après le sous-lemme~\ref{slem-33-2j1245} 
on a \begin{gather}\label{eg-slem-33-11j}\xi_{Z}(\tilde H_{x,q,(t_{1},\dots ,t_{q})}f)=\frac{1}{(p-1)!}\sum_{\tilde Z\in \Lambda_{Z,(t_{1},...,t_{q})}}
\alpha_{Z,\tilde Z,(t_{1},...,t_{q})}\xi_{\tilde Z}(f). \end{gather}
 Par Cauchy-Schwarz  et grâce au sous-lemme~\ref{slem2-33-2j1252} 
 on en déduit que 
 $$|\xi_{Z}(\tilde H_{x,q,(t_{1},\dots ,t_{q})}f)|^{2}\leq C\sum_{\tilde Z\in \Lambda_{Z,(t_{1},...,t_{q})}}
|\xi_{\tilde Z}(f)|^{2}$$ avec $C=C(\de,K,N,Q,P,M)$. On a montré 
  (\ref{ineg-(1-P)-tilde-l0}), donc  (\ref{ineg-1-P-H-natural-etape1}) et  (\ref{ineg-1-P-H-natural}). 
 Ceci termine la preuve du lemme~\ref{sl3-1jan0923}. \cqfd
  %tels que $$\mathcal Y_{0}^{j}\subset \{y\in 
%\bigcup_{\tilde x\in B(x,k),a\in \tilde S_{0}} \geod(\tilde x,a),d(y,\tilde S_{0})=s_{j}\}$$ qui appartiennent à $\{0,\pp,r_{0}(\tilde S_{0})-k+N\}$.  

  \begin{lem}\label{sl4-1jan0923}
Il existe $C=C(\de,K,N,Q,P,M,s,B)$ tel que, pour tout $r\in \N$, 
$$\|(1-\P )u_{x,r} K_x\|_{\L(\H_{x,s}(\Delta_{p-1}),\H^{\rightarrow}_{x,s}(\Delta_{p}))} \leq Ce^{-\frac{sr}{2}}.$$
\end{lem}
\noindent{\bf Démonstration.} 
%  Il reste à montrer   \begin{gather}\label{ineg-(1-P)tildeH}
%\|(1-\P ) \tilde H_{x,q} \|_{\L(\H_{x,s}(\Delta_{p-1}),\H^{\rightarrow}_{x,s}(\Delta_{p}))} \leq C\\ \label{ineg-(1-P)uK}\text{et \ \ \ \ }\end{gather} pour une certaine constante $C=C(\de,K,N,Q,P,M,s,B)$.
%
%  
%  Pour terminer la démonstration du lemme~\ref{continuite-Jx}   il reste à établir l'inégalité (\ref{ineg-(1-P)uK}). 
  Grâce au lemme~\ref{eq-normes-natural}, il suffit de montrer  l'inégalité (\ref{ineg-(1-P)uK-natural}) ci-dessous. 
  Soit $r\in \N$. On va montrer qu'il existe $C=C(\de,K,N,Q,P,M,B)$ tel que pour   $f\in \C^{(\Delta_{p-1})}$, 
\begin{gather}\label{ineg-(1-P)uK-natural}\|(1-\P )(u_{x,r} K_x f)\|^{2}_{\H^{\rightarrow}_{x,s}(\Delta_{p})}\leq 
Ce^{-sr}\|f\|^{2}_{\H_{x,s}^{\natural,Q,1}(\Delta_{p-1})}.\end{gather}
   On a 
\begin{gather*}\|(1-\P )( u_{x,r}K_{x} f)\|_{\H^{\rightarrow}_{x,s}(\Delta_{p})}^{2} =\sum _{k,m,l_{0},\dots ,l_{m}\in \N} 
B^{-(m+\sum_{i=0}^{m}l_{i})} 
\sum_{\substack{Z\in \overline Y_{x}^{p,k,m,(l_{0},...,l_{m})},\\  r_{0}(Z)> k+P}} \\ e^{2s(r_{0}(Z)-k)}\Big(\prod_{i=0}^{m}s_{i}(Z)^{-l_{i}} \Big)
\sharp \big((\pi_{x}^{p,k,m,(l_{0},...,l_{m})})^{-1}(Z)\big)^{-\alpha }
\big|\xi_{Z}(u_{x,r}K_{x} f)\big|^{2}.
%\\ 
%\Big|\sum _{(a_{1},\dots,a_{p},S_{1},...,S_{m},(\mathcal Y_{i}^{j})_{i\in \{0,\dots,m\}, j\in \{1,\dots ,l_{i}\}}) \in 
%(\pi_{x}^{p,k,m,(l_{0},...,l_{m})})^{-1}(Z)} (u_{x,r}K_{x} f)(a_1,...,a_p)\Big|^{2}.
 \end{gather*}
On va voir que l'inégalité (\ref{ineg-(1-P)uK-natural}) résulte de l'inégalité (\ref{ineg-(1-P)uK-natural-etape1}) ci-dessous. 

Soient $k,m,l_{0},\dots ,l_{m}\in \N $ et $Z\in \overline Y_{x}^{p,k,m,(l_{0},...,l_{m})}$ vérifiant $r_{0}(Z)> k+P$. On pose $ \tilde l_{0}=0$ et $ \tilde l_{i}=l_{i-1}$ pour $i\in \{1,\pp,m+1\}$. On va montrer qu'il existe $C=C(\de,K,N,Q,P,M)$ tel que 
   \begin{gather}\nonumber e^{2s(r_{0}(Z)-k)}
\sharp \big((\pi_{x}^{p,k,m,(l_{0},...,l_{m})})^{-1}(Z)\big)^{-\alpha }
\big|\xi_{Z}(u_{x,r}K_{x} f)\big|^{2}
%\\ \Big|\sum _{(a_{1},\dots,a_{p},S_{1},...,S_{m},(\mathcal Y_{i}^{j})_{i\in \{0,\dots,m\}, j\in \{1,\dots ,l_{i}\}}) \in 
%(\pi_{x}^{p,k,m,(l_{0},...,l_{m})})^{-1}(Z)} (u_{x,r}K_{x}  f)(a_1,...,a_p)\Big|^{2}$$
  \\ \nonumber
 \leq Ce^{-sr}
\sum_{\tilde Z\in \Lambda_{Z}} \big(r_{0}(\tilde Z)+1\big)^{-Q}\big(r_{1}(\tilde Z)+1\big)^{-1} e^{2s(r_{0}(\tilde Z)-k)}\\ \label{ineg-(1-P)uK-natural-etape1}
\sharp \big((\pi_{x}^{\natural, p-1,k,m+1,(\tilde l_{0},...,\tilde l_{m+1}),Q,1})^{-1}(\tilde Z)\big)^{-\alpha }
\big|\xi_{\tilde Z}( f)\big|^{2}
%
%
%\Big| \sum _{\substack{(\tilde a_{1},\dots,\tilde a_{p-1},\tilde S_{1},...,\tilde S_{m+1},(\tilde {\mathcal Y}_{i}^{j})_{i\in \{0,\dots,m+1\}, j\in \{1,\dots ,\tilde l_{i}\}},(\tilde {\mathcal Z}_{0}^{j})_{ j\in \{1,\dots ,Q\}},\tilde {\mathcal Z}_{1}^{1})\\ \in 
%(\pi_{x}^{\natural, p-1,k,m+1,(\tilde l_{0},...,\tilde l_{m+1}),Q,1})^{-1}(\tilde Z)}}  f(\tilde a_{1},\dots,\tilde a_{p-1})\Big|^{2}
\end{gather}
 où $\Lambda_{Z}$ est la partie de   $\overline Y_{x}^{\natural, p-1,k,m+1,(\tilde l_{0},...,\tilde l_{m+1}),Q,1}$ formée des $\tilde Z$ 
 vérifiant \begin{gather}\label{r0-r1-r} |r_{0}(\tilde Z)-r_{1}(\tilde Z)-r|\leq QF, \end{gather} et 
tels que pour tout  \begin{gather*}(\tilde a_{1},\dots,\tilde a_{p-1},\tilde S_{0},...,\tilde S_{m+1},(\tilde {\mathcal Y}_{i}^{j})_{i\in \{0,\dots,m+1\}, j\in \{1,\dots ,\tilde l_{i}\}},(\tilde {\mathcal Z}_{0}^{j})_{ j\in \{1,\dots ,Q\}},\tilde {\mathcal Z}_{1}^{1})\\ \in 
(\pi_{x}^{\natural, p-1,k,m+1,(\tilde l_{0},...,\tilde l_{m+1}),Q,1})^{-1}(\tilde Z)\end{gather*} 
il existe une énumération $(a_{1},\pp,a_{p})$ de $\tilde S_{1}$ vérifiant 
%il existe $(a_{1},\pp,a_{p})$ (vérifiant 
% $\tilde S_{1}=\{a_{1},\pp,a_{p}\}$) tel que 
 \begin{gather}\label{notation-34-2j1259}(a_{1},\pp,a_{p},\tilde S_{1},...,\tilde S_{m+1},(\tilde {\mathcal Y}_{i+1}^{j})_{i\in \{0,\dots,m\}, j\in \{1,\dots , l_{i}\}})\in (\pi_{x}^{p,k,m,(l_{0},...,l_{m})})^{-1}(Z).\end{gather}
 Pour $\tilde Z\in \Lambda_{Z}$ on a  $\prod_{i=0}^{m}s_{i}(Z)^{-l_{i}}=\prod_{i=0}^{m+1}s_{i}(\tilde Z)^{-\tilde l_{i}}$. De plus  $\tilde Z$ détermine $Z$ à permutation près de $a_{1},...,a_{p}$ donc connaissant $\tilde Z$ il y a au plus $p!$ possibilités pour $Z$ tels que $\tilde Z\in \Lambda_{Z}$. 
 Donc en sommant sur $Z$ on voit que (\ref{ineg-(1-P)uK-natural-etape1}) implique (\ref{ineg-(1-P)uK-natural}).  
   
 L'inégalité (\ref{ineg-(1-P)uK-natural-etape1}) résulte de l'inégalité plus précise (\ref{ineg-(1-P)-tilde-l0-bis}) ci-dessous (en reprenant les arguments de la fin de la démonstration du lemme~\ref{sl2-1jan0923}).  
 
 Soient $t,t_{1},\pp,t_{Q}\in [0,1]$.  
On va montrer qu'il existe   $C=C(\de,K,N,Q,P,M)$ tel que   \begin{gather}\nonumber e^{2s(r_{0}(Z)-k)}
\sharp \big((\pi_{x}^{p,k,m,(l_{0},...,l_{m})})^{-1}(Z)\big)^{-\alpha }
 \big|\xi_{Z}(u_{x,r,t}K_{x,Q,(t_{1},\dots ,t_{Q})}  f)\big|^{2}
 %\Big|\sum _{\substack{(a_{1},\dots,a_{p},S_{1},...,S_{m},(\mathcal Y_{i}^{j})_{i\in \{0,\dots,m\}, j\in \{1,\dots ,l_{i}\}}) \\\in 
%(\pi_{x}^{p,k,m,(l_{0},...,l_{m})})^{-1}(Z)}} (
%%
% u_{x,r,t}K_{x,Q,(t_{1},\dots ,t_{Q})}  f)(a_1,...,a_p)\Big|^{2}$$ $$
 \leq Ce^{-sr}\\ \label{ineg-(1-P)-tilde-l0-bis}
 \sum_{\tilde Z\in \Lambda_{Z,t,(t_{1},...,t_{Q})}}  e^{2s(r_{0}(\tilde Z)-k)}
\sharp \big((\pi_{x}^{\natural, p-1,k,m+1,(\tilde l_{0},...,\tilde l_{m+1}),Q,1})^{-1}(\tilde Z)\big)^{-\alpha }
\big|\xi_{\tilde Z}(f)\big|^{2}
%\Big| \sum _{\substack{(\tilde a_{1},\dots,\tilde a_{p-1},\tilde S_{1},...,\tilde S_{m+1},(\tilde {\mathcal Y}_{i}^{j})_{i\in \{0,\dots,m+1\}, j\in \{1,\dots ,\tilde l_{i}\}},(\tilde {\mathcal Z}_{0}^{j})_{ j\in \{1,\dots ,Q\}},\tilde {\mathcal Z}_{1}^{1}) \\ \in 
%(\pi_{x}^{\natural, p-1,k,m+1,(\tilde l_{0},...,\tilde l_{m+1}),Q,1})^{-1}(\tilde Z)}} f(\tilde a_{1},\dots,\tilde a_{p-1})\Big|^{2}
\end{gather}
 où $\Lambda_{Z,t,(t_{1},...,t_{Q})}$ est l'ensemble des  $\tilde Z\in 
 \Lambda_{Z}$ 
% \overline Y_{x}^{\natural, p-1,k,m+1,(\tilde l_{0},...,\tilde l_{m+1}),Q,1}$ vérifiant (\ref{r0-r1-r}) et 
 %$$ |r_{0}(\tilde Z)-r_{1}(\tilde Z)-r|\leq (Q+1)N+5\de p_{\max}, $$ 
tels que pour tout  \begin{gather*}(\tilde a_{1},\dots,\tilde a_{p-1},\tilde S_{0},...,\tilde S_{m+1},(\tilde {\mathcal Y}_{i}^{j})_{i\in \{0,\dots,m+1\}, j\in \{1,\dots ,\tilde l_{i}\}},(\tilde {\mathcal Z}_{0}^{j})_{ j\in \{1,\dots ,Q\}},\tilde {\mathcal Z}_{1}^{1})\\ \in (\pi_{x}^{\natural, p-1,k,m+1,(\tilde l_{0},...,\tilde l_{m+1}),Q,1})^{-1}(\tilde Z)\end{gather*} 
 on ait 
 \begin{gather}\label{defZ0j-1j1040}\tilde{\mathcal Z}_{0}^{j}=\bigcup_{b\in \tilde S_{0}}\{z\in \geod(x,b), d(x,z)=E(t_{j}r_{0}(\tilde Z))\}\text{ pour }j\in \{1,...,Q\}\\ \label{defZ1j-1j1040}\text{et\ \ \ \ }  \tilde{\mathcal Z}_{1}^{1}=\bigcup_{b\in \tilde S_{1}}\{z\in \geod(x,b), d(x,z)=E((1-t)r_{1}(\tilde Z))\}.\end{gather} 
 %et 
%il existe $(a_{1},\pp,a_{p})$ vérifiant 
% $\tilde S_{1}=\{a_{1},\pp,a_{p}\}$ et 
% $$(a_{1},\pp,a_{p},\tilde S_{2},...,\tilde S_{m+1},(\tilde {\mathcal Y}_{i+1}^{j})_{i\in \{0,\dots,m\}, j\in \{1,\dots , l_{i}\}})\in (\pi_{x}^{p,k,m,(l_{0},...,l_{m})})^{-1}(Z).$$ 
 Pour $\tilde Z\in \Lambda_{Z,t,(t_{1},...,t_{Q})}$, les conditions (\ref{defZ0j-1j1040}) et (\ref{defZ1j-1j1040}) impliquent 
 $$t_{0}^{j}(\tilde Z)=E(t_{j}r_{0}(\tilde Z))\text{\ \ pour\ \ }j\in \{1,...,Q\} \text{\ \  et \ \ }t_{1}^{1}(\tilde Z)
 =E((1-t)r_{1}(\tilde Z)).$$ 

\begin{souslem}\label{slem0-4eme-2j1302} 
Soit  $$(\tilde a_{1},\dots,\tilde a_{p-1},\tilde S_{0},...,\tilde S_{m+1},(\tilde {\mathcal Y}_{i}^{j})_{i\in \{0,\dots,m+1\}, j\in \{1,\dots ,\tilde l_{i}\}})\in Y_{x}^{p-1,k,m+1,(\tilde l_{0},...,\tilde l_{m+1})}$$
  tel que $  |d(x,\tilde S_{0})-d(x,\tilde S_{1})-r|\leq QF
  $ et 
   qu'il existe $(a_{1},\pp,a_{p})$ vérifiant  (\ref{notation-34-2j1259}). 
   On définit $(\tilde{\mathcal Z}_{0}^{j})_{j\in \{1,...,Q\}}$ et $\tilde{\mathcal Z}_{1}^{1}$ par (\ref{defZ0j-1j1040})  et (\ref{defZ1j-1j1040}). Alors les parties $\tilde{\mathcal Z}_{0}^{j}$ et $\tilde{\mathcal Z}_{1}^{1}$ sont  de diamètre inférieur ou égal à $P/3$  et il existe $\tilde Z\in \Lambda_{Z,t,(t_{1},...,t_{Q})}$ tel que 
\begin{gather}\label{element-slem0-4eme-2j1302}
(\tilde a_{1},\dots,\tilde a_{p-1},\tilde S_{0},...,\tilde S_{m+1},(\tilde {\mathcal Y}_{i}^{j})_{i\in \{0,\dots,m+1\}, j\in \{1,\dots ,\tilde l_{i}\}},(\tilde {\mathcal Z}_{0}^{j})_{ j\in \{1,\dots ,Q\}},\tilde {\mathcal Z}_{1}^{1})\end{gather}
 appartienne à 
 %$ Y_{x}^{\natural, p-1,k,m+1,(\tilde l_{0},...,\tilde l_{m+1}),Q,1}$. 
$(\pi_{x}^{\natural, p-1,k,m+1,(\tilde l_{0},...,\tilde l_{m+1}),Q,1})^{-1}(\tilde Z)$.
%De plus les conditions (\ref{defZ0j-1j1040}) et (\ref{defZ1j-1j1040})  ne dépendent que de $\tilde Z$. 
\end{souslem}
\noindent{\bf Démonstration.} Pour $\sigma\in \{0,1\}$ et $z,z'\in \tilde{\mathcal Z}_{\sigma}^{j}$, on choisit $b\in \tilde S_{\sigma}$, d'où $z,z'\in 2N\tg(x,b)$ et comme $d(x,z)=d(x,z')$, $(H_{\de}(z,x,z',b))$ donne $d(z,z')\leq 2N+\de\leq P/3$. 
 Comme les parties $\tilde{\mathcal Z}_{0}^{j}$ et $\tilde{\mathcal Z}_{1}^{1}$  sont non vides et $P/3\leq M$,  l'argument que nous venons de donner montre aussi que les conditions (\ref{defZ0j-1j1040}) et (\ref{defZ1j-1j1040}) sont vérifiées par les autres éléments de la classe d'équivalence  $\tilde Z$ de l'élément (\ref{element-slem0-4eme-2j1302}) et donc $\tilde Z\in \Lambda_{Z,t,(t_{1},...,t_{Q})}$.  
\cqfd

 Le sous-lemme suivant explique d'où vient la condition (\ref{r0-r1-r}). 
 \begin{souslem}\label{ineg73-3j2030}
 Pour $S\in \Delta_{p-1}$ et $T\in \Delta_{p}$ tels que $e_{T}$ apparaisse avec un coefficient non nul dans $u_{x,r,t}K_{x,Q,(t_{1},\dots ,t_{Q})} (e_{S})$, on a  \begin{gather*} |d(x,S)-d(x,T)-r|\leq QF\end{gather*}
 \end{souslem}
 \noindent{\bf Démonstration.} 
   D'après le 2)a) de la proposition~\ref{recap-supp-connaiss-H-uK},  on a  \begin{gather*}\nonumber  T \subset \bigcup_{a\in S}\{y\in F\tg(x,a), d(y,a)\in [r+\frac{Q}{F},r+QF]\}\end{gather*} d'où l'énoncé du sous-lemme 
 car $ \frac{Q}{F}\geq N+F$. \cqfd

\begin{souslem}\label{slem1-4eme-2j1305}
Soit $\tilde Z\in \Lambda_{Z,t,(t_{1},...,t_{Q})}$, et  \begin{gather*}(\tilde a_{1},\dots,\tilde a_{p-1},\tilde S_{0},...,\tilde S_{m+1},(\tilde {\mathcal Y}_{i}^{j})_{i\in \{0,\dots,m+1\}, j\in \{1,\dots ,\tilde l_{i}\}},(\tilde {\mathcal Z}_{0}^{j})_{ j\in \{1,\dots ,Q\}},\tilde {\mathcal Z}_{1}^{1})\\ \in (\pi_{x}^{\natural, p-1,k,m+1,(\tilde l_{0},...,\tilde l_{m+1}),Q,1})^{-1}(\tilde Z).\end{gather*} 
Alors 
 $u_{x,r,t} K_{x,Q,(t_{1},\dots ,t_{Q})}  (e_{\tilde S_{0}})$ ne dépend que de la connaissance des 
 points de 
 \begin{gather}\label{ineg-apres55-page105-22dec} B(x,k+2M)\cup B(\tilde S_{0}, M) \cup 
B( \tilde S_{1},M) 
\cup  \bigcup _{ j\in \{1,\dots ,Q\}} B(\tilde{\mathcal Z}_{0}^{j}, M)
\cup  B(\tilde{\mathcal Z}_{1}^{1}, M) \end{gather}
et des distances entre ces points. 
\end{souslem}
\noindent{\bf Démonstration.} 
 D'après  le 2)b) de la proposition~\ref{recap-supp-connaiss-H-uK},  $u_{x,r,t} K_{x,Q,(t_{1},\dots ,t_{Q})}  (e_{\tilde S_{0}})$ ne dépend que de la connaissance des points de 
 \begin{gather}
 \label{ens1-4eme-3j2034}
 B(x,F )\cup B(\tilde S_{0},QN) 
 \\ \label{ens2-4eme-3j2034}
\cup 
\bigcup_{a\in \tilde S_{0},j\in \{1,...,Q\}}\{y\in F\text{-}\geod(x,a), |d(x,y)-t_{j}d(x,a)|\leq QF\}
\\ \label{ens3-4eme-3j2034}
\cup \bigcup_{a\in \tilde S_{0}} \{z\in F\text{-}\geod(x,a), d(z,a)\in [r,r+QF]\}
\\ \label{ens4-4eme-3j2034}
 \cup \bigcup_{a\in \tilde S_{0}}\big\{z\in F\tg(x,a), |d(x,z)-(1-t)(d(x,a)-r)|\leq QF \big\}
 \end{gather}
et des distances entre ces points. Il suffit donc de montrer que cet ensemble est inclus dans (\ref{ineg-apres55-page105-22dec}). 

   On suppose $F \leq 2M$ et $QN\leq M$, ce qui est permis par $(H_{M})$. 
   Donc (\ref{ens1-4eme-3j2034}) 
   est inclus dans (\ref{ineg-apres55-page105-22dec}).

Soient  $a\in \tilde S_{0}$, $j\in \{1,...,Q\}$, $y\in F\text{-}\geod(x,a)$ vérifiant  $$|d(x,y)-t_{j}d(x,a)|\leq QF .$$ Soit  $z\in \geod(x,a)$ vérifiant $d(x,z)=E(t_{j}r_{0}(\tilde Z))$,  
si bien que $z$ appartient à $\tilde{\mathcal Z}_{0}^{j}$. 
On a 
$$|t_{j}d(x,a)-E(t_{j}r_{0}(\tilde Z))|\leq N+1,$$ d'où 
$|d(x,y)-d(x,z)|\leq QF+N+1$ et  grâce à $(H_{\de}(y,x,z,a))$, $d(y,z)\leq 
(QF+N+1)+F+\de$. On suppose $(QF+N+1)+F+\de \leq M$, ce qui est permis par $(H_{M})$.  Donc (\ref{ens2-4eme-3j2034}) 
   est inclus dans $\bigcup_{j\in \{1,...,Q\}}B(\tilde{\mathcal Z}_{0}^{j}, M)$ et a fortiori dans (\ref{ineg-apres55-page105-22dec}).

 Soit $a\in \tilde S_{0}$ et $z\in  F\tg(x,a)$ vérifiant $d(z,a)\in [r,r+QF]$. Soit  $y\in \tilde S_{1}$.
Par le a) du lemme~\ref{lemme-S0-...Sm}, on a 
 $y\in  2F\tg(a,x)$. 
 La condition  (\ref{r0-r1-r})  implique 
$$|d(x,a)-d(x,y)-r|\leq QF+N. $$ 
Comme  $d(x,a)-d(x,z)\in [r-F ,r+QF]$, on en déduit $$|d(x,y)-d(x,z)|\leq 2QF+N. $$
Comme $z\in  F\tg(x,a)$ et $y\in  2F\tg(x,a)$, $(H_{\de}(z,x,y,a))$ montre 
$$d(y,z)\leq (2QF+N)+2F+\de.$$
On suppose $(2QF+N)+2F+\de \leq M$, ce qui est permis par $(H_{M})$. 
Donc (\ref{ens3-4eme-3j2034}) 
   est inclus dans $B(\tilde S_{1}, M)$ et a fortiori dans (\ref{ineg-apres55-page105-22dec}). 
%
%d'où $$
%\bigcup_{a\in \tilde S_{0}} \{z\in F\text{-}\geod(x,a), d(z,a)\in [r,r+QF]\}$$
%$$\subset B(\tilde S_{1}, 2QF+N+2F+\de)\subset B(\tilde S_{1}, M)$$
%car 

Enfin soit $a\in \tilde S_{0}$ et $z\in F\tg(x,a)$ vérifiant $$ |d(x,z)-(1-t)(d(x,a)-r)|\leq QF.$$ Soit $b\in \tilde S_{1}$ et $y\in \geod(x,b)$ vérifiant  $d(x,y)=E((1-t)r_{1}(\tilde Z))$, si bien que $y$ appartient à $\tilde{\mathcal Z}_{1}^{1}$. Comme $y\in \geod(x,b)$ et $b\in 2F\tg(x,a)$, on a $y\in 2F\tg(x,a)$. 
Comme $d(x,a)\in [r_{0}(\tilde Z),r_{0}(\tilde Z)+N]$ et grâce à (\ref{r0-r1-r}), on a 
$$|d(x,z)-(1-t)r_{1}(\tilde Z)|\leq 2QF+N,$$
d'où 
$$|d(x,y)-d(x,z)|\leq 2QF+1+N.$$
Comme $z\in F\tg(x,a)$ et  $y\in 2F\tg(x,a)$, $(H_{\de}(y,x,z,a))$ montre 
$$d(y,z)\leq (2QF+1+N)+2F+\de.$$
On suppose $(2QF+1+N)+2F+\de\leq M$, ce qui est permis
 par $(H_{M})$. 
 Donc (\ref{ens4-4eme-3j2034}) 
   est inclus dans $B(\tilde{\mathcal Z}_{1}^{1},M)$ et a fortiori dans (\ref{ineg-apres55-page105-22dec}). 
  % On en déduit 
%\begin{gather*} \bigcup_{a\in S}\big\{z\in F\tg(x,a), |d(x,z)-(1-t)(d(x,a)-r)|\leq QF \big\} 
%\subset   B(\tilde{\mathcal Z}_{1}^{1},M)
%  .\end{gather*}
%
%
%
%Donc   l'ensemble figurant dans le 2)b) de la proposition~\ref{recap-supp-connaiss-H-uK} avec $\tilde S_{0}$ au lieu de $S$  (et qui a la propriété que 
% $u_{x,r,t} K_{x,Q,(t_{1},\dots ,t_{Q})}  (e_{\tilde S_{0}})$ ne dépend que de la connaissance des distances entre les points de cet ensemble) 
%  est inclus dans 
%\begin{gather}\label{ineg-apres55-page105-22dec} B(x,k+2M)\cup B(\tilde S_{0}, M) \cup 
%B( \tilde S_{1},M) 
%\cup  \bigcup _{ j\in \{1,\dots ,Q\}} B(\tilde{\mathcal Z}_{0}^{j}, M)
%\cup  B(\tilde{\mathcal Z}_{1}^{1}, M). \end{gather}
%Ceci termine la démonstration du sous-lemme~\ref{slem1-4eme-2j1305}. 
\cqfd

\begin{souslem}\label{slem2-4eme-2j1308}
Le cardinal de $\Lambda_{Z,t,(t_{1},...,t_{Q})}$ est majoré par une constante de la forme $C(\de,K,N,Q,P,M)$. 
  \end{souslem}
  \noindent{\bf Démonstration.} 
Grâce au lemme~\ref{nombre-dist-connaitre-par-point-natural}, pour connaître les distances entre les points de 
\begin{gather}\label{ens1-p107-24oct09} B(\tilde S_{0}, M)\cup   \bigcup _{ j\in \{1,\dots ,Q\}}B(\tilde {\mathcal Z}_{0}^{j}, M)\cup B(\tilde{\mathcal Z}_{1}^{1}, M)\end{gather} 
 et ceux de 
\begin{gather}\label{ens2-p107-24oct09}
\bigcup _{ i\in \{1,\dots ,m+1\}}
B(\tilde S_{i},M)
\cup  \bigcup _{i\in \{0,\dots,m\}, j\in \{1,\dots , l_{i}\}}B(\tilde {\mathcal Y}_{i+1}^{j}, M) \cup B(x,k+2M)\end{gather}
il suffit de connaître les distances entre les points de (\ref{ens1-p107-24oct09})  et $C$ points de (\ref{ens2-p107-24oct09}), avec $C=C(\de,K,N,Q,P,M)$ et grâce à (\ref{r0-r1-r}) ces distances sont déterminées à $C'=C(\de,K,N,Q,P,M)$ près par les distances de $\tilde S_{1}$ à ces $C$ points (qui font partie de la donnée de $Z$) et les entiers $r,(t_{0}^{j}(\tilde Z))_{j\in \{1,\pp,Q\}}$ et $t_{1}^{1}(\tilde Z)$, qui sont eux-mêmes déterminés à $C''=C(\de,K,N,Q,P,M)$ près par
$r_{0}(Z),r,t,t_{1},...,t_{Q}$. 
 %tels que $$\mathcal Y_{0}^{j}\subset \{y\in 
%\bigcup_{\tilde x\in B(x,k),a\in \tilde S_{0}} \geod(\tilde x,a),d(y,\tilde S_{0})=s_{j}\}$$ qui appartiennent à $\{0,\pp,r_{0}(\tilde S_{0})-k+N\}$.   
\cqfd
 
 \noindent{\bf Suite de la démonstration du lemme~\ref{sl4-1jan0923}.}
 On termine  maintenant la preuve de l'inégalité (\ref{ineg-(1-P)-tilde-l0-bis}). 
 Pour $\tilde Z\in \Lambda_{Z,t,(t_{1},...,t_{Q})}$, et  \begin{gather*}(\tilde a_{1},\dots,\tilde a_{p-1},\tilde S_{0},...,\tilde S_{m+1},(\tilde {\mathcal Y}_{i}^{j})_{i\in \{0,\dots,m+1\}, j\in \{1,\dots ,\tilde l_{i}\}},(\tilde {\mathcal Z}_{0}^{j})_{ j\in \{1,\dots ,Q\}},\tilde {\mathcal Z}_{1}^{1})\\ \in (\pi_{x}^{\natural, p-1,k,m+1,(\tilde l_{0},...,\tilde l_{m+1}),Q,1})^{-1}(\tilde Z)\end{gather*} 
 on considère 
\begin{gather}
\label{somme-30dec0949-2}
\sum _{(b_{1},...,b_{p})}
(u_{x,r,t}K_{x,Q,(t_{1},\dots ,t_{Q})} (e_{\tilde a_{1}}\wedge ...\wedge e_{\tilde a_{p-1}}))(b_{1},...,b_{p}), 
\end{gather} 
où la somme porte sur les énumérations de $\tilde S_{1}=\{b_{1},...,b_{p}\}$ telles que 
$$(b_{1},\dots,b_{p},\tilde S_{1},...,\tilde S_{m+1},(\tilde {\mathcal Y}_{i+1}^{j})_{i\in \{0,\dots,m\}, j\in \{1,\dots , l_{i}\}}
)\in (\pi_{x}^{p,k,m,(l_{0},..., l_{m})})^{-1}(Z). $$
Comme la somme (\ref{somme-30dec0949-2}) a au plus $p!$ termes, 
 le 3) de la proposition~\ref{recap-supp-connaiss-H-uK} montre qu'elle 
est majorée par une constante de la forme $C(\de,K,N,Q,P,M)$. 
D'après le sous-lemme~\ref{slem1-4eme-2j1305}
 la somme (\ref{somme-30dec0949-2}) ne dépend que de $\tilde Z$  et on peut donc la noter $\alpha_{Z,\tilde Z,t,(t_{1},...,t_{Q})}$. 
D'après les sous-lemmes~\ref{slem0-4eme-2j1302} et~\ref{ineg73-3j2030} on a \begin{gather}\label{eg-slem4-11j}\xi_{Z}(u_{x,r,t}K_{x,Q,(t_{1},\dots ,t_{Q})}f)=\frac{1}{(p-1)!}\sum_{\tilde Z\in \Lambda_{Z,t,(t_{1},...,t_{Q})}}
\alpha_{Z,\tilde Z,t,(t_{1},...,t_{Q})}\xi_{\tilde Z}(f). \end{gather}
 Par Cauchy-Schwarz et grâce au sous-lemme~\ref{slem2-4eme-2j1308}, 
 on en déduit que 
 \begin{gather}\label{ineg-1j1105}|\xi_{Z}(u_{x,r,t}K_{x,Q,(t_{1},\dots ,t_{Q})}f)|^{2}\leq C\sum_{\tilde Z\in \Lambda_{Z,t,(t_{1},...,t_{Q})}}
|\xi_{\tilde Z}(f)|^{2}\end{gather} avec $C=C(\de,K,N,Q,P,M)$. 

\begin{souslem}\label{slem3-b-2j1836-0}
  Il existe $D'=C(\de,K)$  et $C= C(\de,K,N,Q,P,M)$ tels que pour tout 
$$(a_{1},\dots,a_{p},S_{0},...,S_{m},(\mathcal Y_{i}^{j})_{i\in \{0,\dots,m\}, j\in \{1,\dots ,l_{i}\}}) \in 
Y_{x}^{p,k,m,(l_{0},...,l_{m})}$$
 le nombre de possibilités pour $(\tilde a_{1},\dots,\tilde a_{p-1})$ tels que 
 $$(\tilde a_{1},\dots,\tilde a_{p-1},\{\tilde a_{1},\dots,\tilde a_{p-1}\}, S_{0},...,S_{m},
( {\mathcal Y}_{i-1}^{j})_{i\in \{0,\dots,m+1\}, j\in \{1,\dots ,\tilde l_{i}\}})$$  appartienne à $
Y_{x}^{p-1,k,m+1,(\tilde l_{0},...,\tilde l_{m+1})}
$ 
 et vérifie 
\begin{gather}\label{ineg-slem-3j2118}|d(x,\{\tilde a_{1},\dots,\tilde a_{p-1}\})-d(x,S_{0})-r|\leq QF\end{gather} soit $\leq C e^{D'r}$.  
\end{souslem}
 \noindent{\bf Démonstration.} 
 Il existe  $D'=C(\de,K)$ tel que tout  $y\in X$ et pour tout $R\in \N$,   $\sharp B(y,R)\leq e^{D'R}$.
Pour tout $a\in \{\tilde a_{1},\dots,\tilde a_{p-1}\}$, grâce à (\ref{ineg-slem-3j2118}) et au fait que 
  $S_{0}\subset 2F\tg(a,x) $ on a
 $d(a,S_{0})\leq r+QF+N+2F$. \cqfd
 
 \noindent{\bf Fin de la démonstration du lemme~\ref{sl4-1jan0923}.}
 Grâce à (\ref{r0-r1-r}), on peut appliquer le sous-lemme~\ref{slem3-b-2j1836-0} et 
 grâce au lemme~\ref{lemme-cardinaux}, on en  déduit qu'il existe $D'=C(\de,K)$ et $C=C(\de,K,N,Q,P,M)$ tels que pour $\tilde Z\in \Lambda_{Z,t,(t_{1},...,t_{Q})}$, $$\sharp \big((\pi_{x}^{\natural, p-1,k,m+1,(\tilde l_{0},...,\tilde l_{m+1}),Q,1})^{-1}(\tilde Z)\big)\leq \sharp \big((\pi_{x}^{p,k,m,(l_{0},...,l_{m})})^{-1}(Z)\big)
 C e^{D'r}.$$  
  Grâce à (\ref{r0-r1-r}) on a $$ e^{2s(r_{0}(Z)-k)}
\leq e^{2s(r_{0}(\tilde Z)-k)} e^{QF}e^{-2sr}.$$ 
 On suppose  $\alpha D'\leq s$, ce qui est permis par   $(H_{\alpha})$. Donc 
il existe $C=C(\de,K,N,Q,P,M)$ tel que
\begin{gather}\nonumber \Big(e^{2s(r_{0}(Z)-k)}\sharp \big((\pi_{x}^{p,k,m,(l_{0},...,l_{m})})^{-1}(Z)\big)^{-\alpha}\Big)\\ 
\label{ineg-1j1106}
\leq Ce^{-sr}\Big(e^{2s(r_{0}(\tilde Z)-k)}\sharp \big((\pi_{x}^{p-1,k,m+1,(\tilde l_{0},...,\tilde l_{m+1})})^{-1}(\tilde Z)\big)^{-\alpha}\Big).\end{gather}  
L'inégalité (\ref{ineg-(1-P)-tilde-l0-bis})   résulte  
de (\ref{ineg-1j1105}) et (\ref{ineg-1j1106}). On a  montré (\ref{ineg-(1-P)-tilde-l0-bis}) et donc (\ref{ineg-(1-P)uK-natural-etape1}) et  (\ref{ineg-(1-P)uK-natural}). Ceci termine la preuve  du lemme~\ref{sl4-1jan0923}. \cqfd
  
  On a donc montré la proposition~\ref{continuite-Jx}.

\subsection{Continuité de $e^{\tau  \theta^{\flat}_{x}}\del  e^{-\tau  \theta^{\flat}_{x}}$ et 
$e^{\tau  \theta^{\flat}_{x}}J_{x}e^{-\tau  \theta^{\flat}_{x}}$}

Le but de ce sous-paragraphe est de montrer la proposition~\ref{continuite-del-J-conj}. 

\begin{prop}\label{continuite-del-J-conj}
Pour tout $T\in \R_{+}$ et tout $r\in \N$,   \begin{gather*}(e^{\tau  \theta^{\flat}_{x}}\del e^{-\tau  \theta^{\flat}_{x}})_{\tau  \in [0,T]}, \ (e^{\tau  \theta^{\flat}_{x}}J_{x}e^{-\tau  \theta^{\flat}_{x}})_{\tau  \in [0,T]}, \ (e^{\tau  \theta^{\flat}_{x}}h_{x}e^{-\tau  \theta^{\flat}_{x}})_{\tau  \in [0,T]}\\ \text{\  et\ } (e^{\tau  \theta^{\flat}_{x}}u_{x,r}K_{x}e^{-\tau  \theta^{\flat}_{x}})_{\tau  \in [0,T]}\end{gather*}
  s'étendent   en des   opérateurs 
continus sur le  $\C[0,T]$-module hilbertien 
$\H_{x,s}[0,T]$. 
\end{prop}
On a inclus les opérateurs $h_{x}$ et $u_{x,r}K_{x}$ dans l'énoncé de  cette proposition  pour un usage ultérieur. 

On rappelle que pour tout $t\in \R_{+}$ et pour $p\in \{1,\pp,p_{\max}\}$, on a défini $\theta^{\flat}_{x}:\C^{(\Delta_{p})}\to \C^{(\Delta_{p})}$ par $\theta^{\flat}_{x}(e_{S})=\rho^{\flat}_{x}(S)e_{S}$ pour tout $S\in \Delta_{p}$ et que  la fonction $\rho^{\flat}_{x}:X\to \R_{+}$ avait été définie 
par $\rho^{\flat}_{x}(a)=d^{\flat}(x,a)$ et étendue  en une fonction 
$\rho^{\flat}_{x}:\Delta\to \R_{+}$ par la formule 
$$\rho^{\flat}_{x}(S)=\frac{\sum _{a\in S}\rho_{x}^{\flat}(a)}{\sharp S}\text{ si }S\text{ est non vide et }\rho^{\flat}_{x}(\emptyset)=0.$$ 

De fa\c con analogue on définit  $\rho_{x}:X\to \R_{+}$ en posant  
$\rho_{x}(a)=d(x,a)$
et on étend cette fonction en 
$\rho_{x}: \Delta\to \R_{+}$ par la formule \label{def-rhox-rhox0Z}
$$\rho_{x}(S)=\frac{\sum _{a\in S}\rho_{x}(a)}{\sharp S}\text{ si }S\text{ est non vide et }\rho_{x}(\emptyset)=0.$$ 
D'après la proposition~\ref{prop-d'},  pour $x,y\in X$, on a $d(x,y)\leq d^{\flat}(x,y)\leq d(x,y)+7\de$. Il en résulte que pour tout $S\in \Delta$ on a $\rho_{x}(S)\leq \rho^{\flat}_{x}(S)\leq \rho_{x}(S)+7\de $.
On définit l'opérateur $\theta_{x}:\C^{(\Delta_{p})}\to \C^{(\Delta_{p})}$ par $\theta_{x}(e_{S})=\rho_{x}(S)e_{S}$. 

\noindent{\bf Démonstration de la proposition~\ref{continuite-del-J-conj} en admettant les lemmes~\ref{continuite-del-J-conj-d} et~\ref{continuite-conj-conj}.}
Pour $\tau \in \R_{+}$ et $p\in \{1,\pp,p_{\max}\}$, on introduit les opérateurs 
$e^{\tau (\theta^{\flat}_{x}-\theta_{x})}:\C^{(\Delta_{p})}\to \C^{(\Delta_{p})}$ et 
$e^{\tau  \theta_{x}}:\C^{(\Delta_{p})}\to \C^{(\Delta_{p})}$
en posant 
$$e^{\tau (\theta^{\flat}_{x}-\theta_{x})}(e_{S})=e^{\tau (\theta^{\flat}_{x}(S)-\theta_{x}(S))}e_{S} \text{\ \ et\ \ }e^{\tau  \theta_{x}}(e_{S})=e^{\tau  \theta_{x}(S)}e_{S}.$$
Ces opérateurs commutent entre eux et on a bien sûr $ e^{\tau  \theta^{\flat}_{x}}=e^{\tau (\theta^{\flat}_{x}-\theta_{x})}e^{\tau  \theta_{x}}$. Pour $p=0$ on définit $\theta^{\flat}_{x}$ et $\theta_{x}$ comme $0$ sur  $\C^{(\Delta_{0})}=\C$. 

La proposition~\ref{continuite-del-J-conj}
 résulte des lemmes~\ref{continuite-del-J-conj-d} et~\ref{continuite-conj-conj}.   \cqfd

\begin{lem}\label{continuite-del-J-conj-d}
Pour tout $T\in \R_{+}$ et tout $r\in \N$,  
\begin{gather*}(e^{\tau  \theta_{x}}\del e^{-\tau  \theta_{x}})_{\tau  \in [0,T]}, \ (e^{\tau  \theta_{x}}J_{x}e^{-\tau  \theta_{x}})_{\tau  \in [0,T]}, \ (e^{\tau  \theta_{x}}h_{x}e^{-\tau  \theta_{x}})_{\tau  \in [0,T]}\\ \text{\  et\ } (e^{\tau  \theta_{x}}u_{x,r}K_{x}e^{-\tau  \theta_{x}})_{\tau  \in [0,T]}\end{gather*}
%$(e^{\tau  \theta_{x}}\del e^{-\tau  \theta_{x}})_{\tau  \in [0,T]}$ 
%et $(e^{\tau  \theta_{x}}J_{x}e^{-\tau  \theta_{x}})_{\tau  \in [0,T]}$
s'étendent   en des   opérateurs 
continus sur le  $\C[0,T]$-module hilbertien 
$\H_{x,s}[0,T]$. 
\end{lem}
Pour la démonstration de ce lemme on a besoin de la notation suivante. 

\noindent{\bf Notation.} Pour 
$Z\in \overline Y_{x}^{p,k,m,(l_{0},...,l_{m})}$ et $$(a_{1},\dots,a_{p},S_{0},...,S_{m},(\mathcal Y_{i}^{j})_{i\in \{0,\dots,m\}, j\in \{1,\dots ,l_{i}\}})\in (\pi_{x}^{p,k,m,(l_{0},...,l_{m})})^{-1}(Z),$$ 
$\rho_{x}(S_{0})$ et $\rho_{x}(S_{1})$ ne dépendent  que de $Z$ et on les note  $\rho_{x}^{0}(Z)$ et $\rho_{x}^{1}(Z)$. On adopte une notation similaire pour $Z\in \overline Y_{x}^{\natural,p,k,m,(l_{0},...,l_{m}),\lambda_{0},\lambda_{1}}$. 

\noindent{\bf Remarque.} Il est évidemment faux que  $\rho_{x}^{\flat}(S_{0})$ et $\rho_{x}^{\flat}(S_{1})$ ne dépendent que de $Z$    et c'est pour cette raison que  la preuve de  la proposition~\ref{continuite-del-J-conj} n'est pas aussi simple que celle du lemme~\ref{continuite-del-J-conj-d}. 

\noindent{\bf Démonstration du lemme~\ref{continuite-del-J-conj-d}.}
On reprend la démonstration des 
propositions~\ref{continuite-del} et~\ref{continuite-Jx} du  sous-paragraphe précédent, qui affirmaient 
la continuité de  $\del, J_{x}, h_{x}$ et $u_{x,r}K_{x}$. 
Les seuls ingrédients supplémentaires sont les faits suivants : 
\begin{itemize}
\item  dans les notations ci-dessus $\rho_{x}(S_{0})$ et $\rho_{x}(S_{1})$ ne dépendent  que de $Z$
\item d'après le 1)a) et le 2)a) de la proposition~\ref{recap-supp-connaiss-H-uK}, il existe $C=C(\de,K,N,Q)$ tel que si $e_{T}$ apparaît dans $\del(e_{S})$ ou $J_{x}(e_{S})$ avec un coefficient non nul, on a $d(x,T)\leq d(x,S)+C$, d'où $\rho_{x}(T)\leq \rho_{x}(S)+C+N$. 
\end{itemize}

Voici de fa\c con plus précise les modifications à apporter  : 
\begin{itemize}
\item pour la proposition~\ref{continuite-del}, dans (\ref{somme-tildeZ-Y-24oct09}) on remplace $\xi_{\tilde Z}(f)$ par $e^{\tau (\rho_{x}^{1}(\tilde Z)-\rho_{x}^{0}(\tilde Z))}\xi_{\tilde Z}(f)$ et on remarque que $\rho_{x}^{1}(\tilde Z)-\rho_{x}^{0}(\tilde Z)\leq N$, 
\item pour le lemme~\ref{sl2-1jan0923}, dans (\ref{eg-slem0-2eme-11j}) on remplace $\xi_{Z}(f)$ par $e^{\tau (\rho_{x}(U)-\rho_{x}^{0}(Z))}\xi_{Z}(f)$ et on remarque que $\rho_{x}(U)\leq P+N$, 
\item pour  le lemme~\ref{sl3-1jan0923}, dans (\ref{eg-slem-33-11j})  on remplace $\xi_{\tilde Z}(f)$ par $e^{\tau (\rho_{x}^{1}(\tilde Z)-\rho_{x}^{0}(\tilde Z))}\xi_{\tilde Z}(f)$ et on remarque que $\rho_{x}^{1}(\tilde Z)-\rho_{x}^{0}(\tilde Z)\leq (q+2)N$ grâce à (\ref{eq-29dec2048}), 
\item pour  le lemme~\ref{sl4-1jan0923}, dans (\ref{eg-slem4-11j}) 
on remplace $\xi_{\tilde Z}(f)$ par $e^{\tau (\rho_{x}^{1}(\tilde Z)-\rho_{x}^{0}(\tilde Z))}\xi_{\tilde Z}(f)$ et on remarque que $\rho_{x}^{1}(\tilde Z)-\rho_{x}^{0}(\tilde Z)\leq QF+N$ par (\ref{r0-r1-r}) 
\end{itemize}
et en plus on remplace les opérateurs par les opérateurs conjugués à de nombreux endroits (notamment  dans les égalités (\ref{somme-tildeZ-Y-24oct09}), (\ref{eg-slem0-2eme-11j}), (\ref{eg-slem-33-11j}) et (\ref{eg-slem4-11j})). 

En utilisant le fait que les fonctions à support fini sont denses dans $\H_{x,s}$ on montre que les opérateurs du lemme~\ref{continuite-del-J-conj-d} sont continus en $\tau$ pour la topologie forte et leurs adjoints aussi. Ceci justifie le fait qu'ils s'étendent en des morphismes de $\C[0,T]$-modules hilbertiens. 
\cqfd

\begin{lem}\label{continuite-conj-conj}
Pour tout $T\in \R_{+}$, l'opérateur $(e^{\tau (\theta^{\flat}_{x}-\theta_{x})})_{\tau  \in [0,T]}$ 
s'étend en un automorphisme du   $\C[0,T]$-module hilbertien 
$\H_{x,s}[0,T]$. 
\end{lem}
\noindent{\bf Démonstration du  lemme~\ref{continuite-conj-conj} en admettant le lemme~\ref{continuite-eta}.} 
Soit $p\in \{1,\pp,p_{\max}\}$. Le lemme~\ref{continuite-conj-conj}  résulte immédiatement du lemme suivant.  \cqfd

\begin{lem}\label{continuite-eta}
L'opérateur $\theta^{\flat}_{x}-\theta_{x}$  
s'étend en un opérateur continu sur $\H_{x,s}(\Delta_{p})$, dont la norme est bornée par une constante du type $C(\de,N,K,Q,P,M,s,B)$.  
\end{lem}
Bien entendu cet énoncé n'est pas vrai pour les opérateurs $\theta^{\flat}_{x}$ et $\theta_{x}$ séparément.  
Pour $u_{1},u_{2},u_{3},v_{1},v_{2},v_{3}\in [0,1[$ on définit  
$(\rho^{\flat}_{x})_{u_{1},u_{2},u_{3}}^{v_{1},v_{2},v_{3}}:X\to \R_{+}$ en posant  $(\rho^{\flat}_{x})_{u_{1},u_{2},u_{3}}^{v_{1},v_{2},v_{3}}(a)={d^{\flat}}_{u_{1},u_{2},u_{3}}^{v_{1},v_{2},v_{3}}(x,a)$ et on l'étend en une fonction \label{rhoflatxuuuvvv}
$(\rho^{\flat}_{x})_{u_{1},u_{2},u_{3}}^{v_{1},v_{2},v_{3}}:\Delta\to \R_{+}$ par la formule 
$$(\rho^{\flat}_{x})_{u_{1},u_{2},u_{3}}^{v_{1},v_{2},v_{3}}(S)=\frac{\sum _{a\in S}(\rho^{\flat}_{x})_{u_{1},u_{2},u_{3}}^{v_{1},v_{2},v_{3}}(a)}{\sharp S}\text{ si }S\text{ est non vide et }(\rho^{\flat}_{x})_{u_{1},u_{2},u_{3}}^{v_{1},v_{2},v_{3}}(\emptyset)=0.$$ 
Pour  $u_{1},u_{2},u_{3},v_{1},v_{2},v_{3}\in [0,1[$ et $S\in \Delta$ on a 
$$\rho_{x}(S)\leq (\rho^{\flat}_{x})_{u_{1},u_{2},u_{3}}^{v_{1},v_{2},v_{3}}(S) \leq \rho_{x}(S)+7\de.$$ De plus, pour $S\in \Delta$,  
$$ (\rho^{\flat}_{x})(S)=\int_{u_{1},u_{2},u_{3},v_{1},v_{2},v_{3}\in [0,1[}(\rho^{\flat}_{x})_{u_{1},u_{2},u_{3}}^{v_{1},v_{2},v_{3}}(S)du_{1}du_{2}du_{3}dv_{1}dv_{2}dv_{3}. $$

\begin{lem}\label{dependance-rho'-123123}
Pour $S\in \Delta_{p}$, 
  $(\rho^{\flat}_{x})_{u_{1},u_{2},u_{3}}^{v_{1},v_{2},v_{3}}(S)$ ne dépend que de la connaissance des points de 
  \begin{gather}\label{eta-conn-deuxieme-partie}\{x\}\cup S\cup \bigcup_{a\in S, j\in \{1,...,6\}}
  \{y\in 3\de\tg(x,a), |d(x,y)-w_{j}|\leq N+6\de+4\}\end{gather}
  et des distances entre ces points, où l'on note 
  \begin{gather}\nonumber  
  w_{1}=E(\frac{u_{1}}{6}d(x,S)),
  w_{2}=E((\frac{1}{6}+\frac{u_{2}}{6})d(x,S)),
  w_{3}=E((\frac{2}{6}+\frac{u_{3}}{6})d(x,S)),\\ \label{def-wj-30dec1027} 
   w_{4}=E((1-\frac{v_{1}}{6})d(x,S)),
     w_{5}=E((\frac{5}{6}-\frac{v_{2}}{6})d(x,S)), 
    w_{6}=E((\frac{4}{6}-\frac{v_{3}}{6})d(x,S)). \end{gather}
        \end{lem}
        
    \noindent{\bf Démonstration.}  
    L'énoncé est évident si $d(x,S)\leq 6\de$ et on suppose donc 
   $d(x,S)>6\de$. Par la construction même de  $d^{\flat}$,  $(\rho^{\flat}_{x})_{u_{1},u_{2},u_{3}}^{v_{1},v_{2},v_{3}}(S)$ ne dépend que de la connaissance de $x$, de $S$ et de la réunion pour $a\in S$ de
    \begin{gather}\nonumber   Y_{x,a}^{E((\Delta_{x,a}+1)u_{1})}\cup Y_{x,a}^{\Delta_{x,a}+E((\Delta_{x,a}+1)u_{2})}\cup Y_{x,a}^{2\Delta_{x,a}+E((\Delta_{x,a}+1)u_{3})}
\\ \label{eta-30dec1037} \cup Y_{a,x}^{E((\Delta_{x,a}+1)v_{1})}\cup    Y_{a,x}^{ \Delta_{x,a}+E((\Delta_{x,a}+1)v_{2})}\cup     Y_{a,x}^{2\Delta_{x,a}+E((\Delta_{x,a}+1)v_{3})}. 
      \end{gather} 
      Soit $a\in S$. On 
 rappelle que $\Delta_{x,a}=E(\frac{d(x,a)}{6})-\de$ et que 
    pour $u,v\in X$ et 
    $r\in \{0,...,E(\frac{d(u,v)}{2})-3\delta\}$, on note $Y_{u,v}^{r}$ l'ensemble  des points $z\in 3\delta\text{-}\geod(u,v)$ tels que $d(u,z)\in \{r,...,r+3\delta\}$. Il suffit donc de montrer que l'ensemble (\ref{eta-30dec1037}) est inclus dans l'ensemble      \begin{gather}\label{eta-conn-deuxieme-partie2} \bigcup_{ j\in \{1,...,6\}}
  \{y\in 3\de\tg(x,a), |d(x,y)-w_{j}|\leq N+6\de+4\}.\end{gather}
   On a $\big|\Delta_{x,a}-\frac{d(x,a)}{6}\big|\leq \de+1$ et pour $i\in \{1,...,3\}$, 
    $$\big|E((\Delta_{x,a}+1)u_{i})-\frac{u_{i}d(x,a)}{6}\big|\leq \de+1\text{ et } 
\big|E((\Delta_{x,a}+1)v_{i})-\frac{v_{i}d(x,a)}{6}\big|\leq \de+1. $$
Il en résulte facilement que pour $i\in \{1,...,3\}$ et 
$y\in Y_{x,a}^{(i-1)\Delta_{x,a}+E((\Delta_{x,a}+1)u_{i})}$ on a 
\begin{gather*}d(x,y)\in [(i-1)\Delta_{x,a}+E((\Delta_{x,a}+1)u_{i}),(i-1)\Delta_{x,a}+E((\Delta_{x,a}+1)u_{i})+3\de]\\ \subset [\frac{((i-1)+u_{i})d(x,a)}{6}-(3\de+3), 
\frac{((i-1)+u_{i})d(x,a)}{6}+(6\de+3)].\end{gather*} 
D'autre part $|\frac{((i-1)+u_{i})d(x,a)}{6}-w_{i}|\leq N+1$ car $|d(x,a)-d(x,S)|\leq N$. On en déduit que $y$ appartient à (\ref{eta-conn-deuxieme-partie2}).

Soit maintenant $i\in \{1,...,3\}$ et 
$y\in Y_{a,x}^{ (i-1)\Delta_{x,a}+E((\Delta_{x,a}+1)v_{i})}$. 
Donc $d(a,y)$ appartient à \begin{gather}\nonumber  [(i-1)\Delta_{x,a}+E((\Delta_{x,a}+1)v_{i}),(i-1)\Delta_{x,a}+E((\Delta_{x,a}+1)v_{i})+3\de]\\ \label{ineg-30dec1058}\subset [\frac{((i-1)+v_{i})d(x,a)}{6}-(3\de+3), 
\frac{((i-1)+v_{i})d(x,a)}{6}+(6\de+3)].\end{gather} 
Comme $y\in 3\de\tg(x,a)$ on a $d(x,y)\in [d(x,a)-d(a,y), d(x,a)-d(a,y)+3\de]$, et comme $d(a,y)$ appartient à    (\ref{ineg-30dec1058}) on en déduit que $d(x,y)$ appartient à 
\begin{gather*}[\frac{(7-i-v_{i})d(x,a)}{6}-(6\de+3), 
\frac{(7-i-v_{i})d(x,a)}{6}+(6\de+3)].\end{gather*} 
D'autre part $|\frac{(7-i-v_{i})d(x,a)}{6}-w_{i+3}|\leq N+1$ car $|d(x,a)-d(x,S)|\leq N$. On en déduit que $y$ appartient à (\ref{eta-conn-deuxieme-partie2}).
 \cqfd

Enfin on note $(\theta_{x}^{\flat})_{u_{1},u_{2},u_{3}}^{v_{1},v_{2},v_{3}}: \C^{(\Delta_{p})}\to  \C^{(\Delta_{p})}$ l'opérateur défini par $$(\theta_{x}^{\flat})_{u_{1},u_{2},u_{3}}^{v_{1},v_{2},v_{3}}(e_{S})=(\rho^{\flat}_{x})_{u_{1},u_{2},u_{3}}^{v_{1},v_{2},v_{3}}(S)e_{S}, $$ 
de sorte que 
$$ \theta_{x}^{\flat}=\int_{u_{1},u_{2},u_{3},v_{1},v_{2},v_{3}\in [0,1[}(\theta_{x}^{\flat})_{u_{1},u_{2},u_{3}}^{v_{1},v_{2},v_{3}}du_{1}du_{2}du_{3}dv_{1}dv_{2}dv_{3}. $$

\noindent{\bf Démonstration du lemme~\ref{continuite-eta}.} 
Soit $f\in \C^{(\Delta_{p})}$. 
Par définition  
\begin{gather*}\|(\theta^{\flat}_{x}-\theta_{x})(f)\|_{\H_{x,s}(\Delta_{p})}^{2} =\sum _{k,m,l_{0},\dots ,l_{m}\in \N} 
B^{-(m+\sum_{i=0}^{m}l_{i})}
\sum_{Z\in \overline Y_{x}^{p,k,m,(l_{0},...,l_{m})}} \\  e^{2s(r_{0}(Z)-k)}\Big(\prod_{i=0}^{m}s_{i}(Z)^{-l_{i}} \Big)
\sharp \big((\pi_{x}^{p,k,m,(l_{0},...,l_{m})})^{-1}(Z)\big)^{-\alpha }
\big|\xi_{Z}\big((\theta^{\flat}_{x}-\theta_{x}) (f)\big)\big|^{2}.
%\Big|\sum _{(a_{1},\dots,a_{p},S_{1},...,S_{m},(\mathcal Y_{i}^{j})_{i\in \{0,\dots,m\}, j\in \{1,\dots ,l_{i}\}}) \in 
%(\pi_{x}^{p,k,m,(l_{0},...,l_{m})})^{-1}(Z)} ((\theta^{\flat}_{x}-\theta_{x}) (f))(a_1,...,a_p)\Big|^{2}.
\end{gather*}

On va voir que le lemme~\ref{continuite-eta} résulte de l'inégalité (\ref{cont-eta-etape}) ci-dessous. 

  Soient $k,m,l_{0},\dots ,l_{m}\in \N $ et $Z\in \overline Y_{x}^{p,k,m,(l_{0},...,l_{m})}$. On va montrer qu'il existe  $C=C(\de,K,N,Q,P,M)$ tel que    \begin{gather}\label{cont-eta-etape}\big|\xi_{Z}\big((\theta^{\flat}_{x}-\theta_{x}) (f)\big)\big|^{2}
% \Big|\sum _{(a_{1},\dots,a_{p},S_{1},...,S_{m},(\mathcal Y_{i}^{j})_{i\in \{0,\dots,m\}, j\in \{1,\dots ,l_{i}\}}) \in 
%(\pi_{x}^{p,k,m,(l_{0},...,l_{m})})^{-1}(Z)} ((\theta^{\flat}_{x}-\theta_{x}) (f))(a_1,...,a_p)\Big|^{2}$$ 
 \leq C
\sum_{\tilde Z\in \Lambda_{Z}} \big(r_{0}(\tilde Z)+1\big)^{-6}
\big|\xi_{\tilde Z}(f)\big|^{2}
%
% \Big| \sum _{\substack{(\tilde a_{1},\dots,\tilde a_{p},\tilde S_{1},...,\tilde S_{m},(\tilde {\mathcal Y}_{i}^{j})_{i\in \{0,\dots,m\}, j\in \{1,\dots ,l_{i}\}},(\tilde {\mathcal Z}_{0}^{j})_{j\in \{1,\dots ,6\}}) \\ \in 
%(\pi_{x}^{\natural, p,k,m,( l_{0},...,l_{m}),6,0})^{-1}(\tilde Z)}}  f(\tilde a_{1},\dots,\tilde a_{p})\Big|^{2}
\end{gather}
 où $\Lambda_{Z}$ est la partie de  $\overline Y_{x}^{\natural, p,k,m,( l_{0},...,l_{m}),6,0}$ formée des $\tilde Z$ 
tels que pour tout  \begin{gather*}(\tilde a_{1},\dots,\tilde a_{p},\tilde S_{0},...,\tilde S_{m},(\tilde {\mathcal Y}_{i}^{j})_{i\in \{0,\dots,m\}, j\in \{1,\dots ,l_{i}\}},(\tilde {\mathcal Z}_{0}^{j})_{j\in \{1,\dots ,6\}}) \\  \in 
(\pi_{x}^{\natural, p,k,m,( l_{0},...,l_{m}),6,0})^{-1}(\tilde Z)\end{gather*} on ait 
 $$(\tilde a_{1},\pp,\tilde a_{p},\tilde S_{0},...,\tilde S_{m},(\tilde {\mathcal Y}_{i}^{j})_{i\in \{0,\dots,m\}, j\in \{1,\dots , l_{i}\}})\in (\pi_{x}^{p,k,m,(l_{0},...,l_{m})})^{-1}(Z).$$

 On justifie maintenant le fait que l'inégalité (\ref{cont-eta-etape}) implique l'énoncé du lemme. Le lemme~\ref{lemme-cardinaux} montre 
 que pour $\tilde Z\in \Lambda_{Z}$ on a 
 $$\sharp (\pi_{x}^{\natural, p,k,m,( l_{0},...,l_{m}),6,0})^{-1}(\tilde Z)  \leq C\sharp  (\pi_{x}^{p,k,m,(l_{0},...,l_{m})})^{-1}(Z) $$ avec $C=C(\de,K,N,Q,P,M)$. En sommant sur $Z$  
 et en appliquant le  lemme~\ref{eq-normes-natural} à $\mu_{0}=6$ et $\mu_{1}=0$, on voit que l'inégalité (\ref{cont-eta-etape}) implique l'énoncé du lemme.

L'inégalité (\ref{cont-eta-etape}) découle de l'inégalité  (\ref{ineg-eta-uuuvvv}) plus précise ci-dessous. 
Soient $u_{1},u_{2},u_{3},v_{1},v_{2},v_{3}\in [0,1[$. 
On va montrer qu'il existe $C=C(\de,K,N,Q,P,M)$ tel  que  
 \begin{gather}\label{ineg-eta-uuuvvv}
 \big|\xi_{Z}\big(((\theta^{\flat}_{x})_{u_{1},u_{2},u_{3}}^{v_{1},v_{2},v_{3}}-\theta_{x})  (f)\big)\big|^{2} 
% \Big|\sum _{\substack{(a_{1},\dots,a_{p},S_{1},...,S_{m},(\mathcal Y_{i}^{j})_{i\in \{0,\dots,m\}, j\in \{1,\dots ,l_{i}\}})\\  \in 
%(\pi_{x}^{p,k,m,(l_{0},...,l_{m})})^{-1}(Z)}} (((\theta^{\flat}_{x})_{u_{1},u_{2},u_{3}}^{v_{1},v_{2},v_{3}}-\theta_{x})  f)(a_1,...,a_p)\Big|^{2}\\ 
\leq C\sum_{\tilde Z\in (\Lambda_{Z})_{u_{1},u_{2},u_{3}}^{v_{1},v_{2},v_{3}}}  
\big|\xi_{\tilde Z}(f)\big|^{2}
%\Big| \sum _{\substack{(\tilde a_{1},\dots,\tilde a_{p},\tilde S_{1},...,\tilde S_{m},(\tilde {\mathcal Y}_{i}^{j})_{i\in \{0,\dots,m\}, j\in \{1,\dots ,l_{i}\}},(\tilde {\mathcal Z}_{0}^{j})_{j\in \{1,\dots ,6\}}) \\ \in 
%(\pi_{x}^{\natural, p,k,m,( l_{0},...,l_{m}),6,0})^{-1}(\tilde Z)}}  f(\tilde a_{1},\dots,\tilde a_{p})\Big|^{2}
\end{gather}
 où $(\Lambda_{Z})_{u_{1},u_{2},u_{3}}^{v_{1},v_{2},v_{3}}$  est l'ensemble des 
 $\tilde Z\in \Lambda_{Z}
 %\overline Y_{x}^{\natural, p,k,m,( l_{0},...,l_{m}),6,0}
 $ 
tels que, en notant 
\begin{gather*}w_{1}=E(\frac{u_{1}}{6}r_{0}(\tilde Z)),
  \ w_{2}=E((\frac{1}{6}+\frac{u_{2}}{6})r_{0}(\tilde Z)),
  \ w_{3}=E((\frac{2}{6}+\frac{u_{3}}{6})r_{0}(\tilde Z)),\\ 
 w_{4}=E((1-\frac{v_{1}}{6})r_{0}(\tilde Z)),
    \ w_{5}=E((\frac{5}{6}-\frac{v_{2}}{6})r_{0}(\tilde Z)),\ 
    w_{6}=E((\frac{4}{6}-\frac{v_{3}}{6})r_{0}(\tilde Z)) \end{gather*}  
on ait, 
 pour tout  \begin{gather*}(\tilde a_{1},\dots,\tilde a_{p},\tilde S_{0},...,\tilde S_{m},(\tilde {\mathcal Y}_{i}^{j})_{i\in \{0,\dots,m\}, j\in \{1,\dots ,l_{i}\}},(\tilde {\mathcal Z}_{0}^{j})_{j\in \{1,\dots ,6\}}) \\  \in 
(\pi_{x}^{\natural, p,k,m,( l_{0},...,l_{m}),6,0})^{-1}(\tilde Z)\end{gather*}  
% on ait 
% $(\tilde a_{1},\pp,\tilde a_{p},\tilde S_{1},...,\tilde S_{m},(\tilde {\mathcal Y}_{i}^{j})_{i\in \{0,\dots,m\}, j\in \{1,\dots , l_{i}\}})\in Z$ 
% et en notant 
%       on ait, 
et pour tout $j\in \{1,...,6\}$, 
\begin{gather}\label{defZ0j-1j1124}\tilde{\mathcal Z}_{0}^{j}=\bigcup_{b\in \tilde S_{0}}\{z\in \geod(x,b), d(x,z)=w_{j}\}.\end{gather} 
La condition (\ref{defZ0j-1j1124}) implique que pour $\tilde Z\in (\Lambda_{Z})_{u_{1},u_{2},u_{3}}^{v_{1},v_{2},v_{3}}$ et 
pour $j\in \{1,...,6\}$, on a 
  $t_{0}^{j}(\tilde Z)=w_{j}$. 

\begin{souslem}\label{slem0-rho-2j1517}  Soit  $$(\tilde a_{1},\dots,\tilde a_{p},\tilde S_{0},...,\tilde S_{m},(\tilde {\mathcal Y}_{i}^{j})_{i\in \{0,\dots,m\}, j\in \{1,\dots ,l_{i}\}})  \in 
 (\pi_{x}^{p,k,m,(l_{0},..., l_{m})})^{-1}(Z).$$  On définit $(\tilde{\mathcal Z}_{0}^{j})_{j\in \{1,...,6\}}$ par  (\ref{defZ0j-1j1124}). Alors les parties $\tilde{\mathcal Z}_{0}^{j}$  sont  de diamètre inférieur ou égal à $P/3$  et il existe $\tilde Z\in (\Lambda_{Z})_{u_{1},u_{2},u_{3}}^{v_{1},v_{2},v_{3}}$ tel que  
\begin{gather}\label{elem-slem0-rho-11j}(\tilde a_{1},\dots,\tilde a_{p},\tilde S_{0},...,\tilde S_{m},(\tilde {\mathcal Y}_{i}^{j})_{i\in \{0,\dots,m\}, j\in \{1,\dots ,l_{i}\}},(\tilde {\mathcal Z}_{0}^{j})_{j\in \{1,\dots ,6\}})\end{gather} appartienne à 
$(\pi_{x}^{\natural, p,k,m,( l_{0},...,l_{m}),6,0})^{-1}(\tilde Z)$.
\end{souslem}
\noindent{\bf Démonstration.}
Soit $j\in \{1,...,6\}$  et $z,z'\in \tilde {\mathcal Z}_{0}^{j}$. Soit  $b\in \tilde S_{0}$. On a $z,z'\in 2N\tg(x,b)$ et $d(x,z)=d(x,z')$ donc par $(H_{\de}(z,x,z',b))$, $d(z,z')\leq 2N+\de \leq P/3$. 
 Comme les parties $\tilde{\mathcal Z}_{0}^{j}$  sont   non vides  et  $P/3\leq M$, l'argument que nous venons de donner montre aussi que la condition (\ref{defZ0j-1j1124}) est vérifiée par les autres éléments de la classe d'équivalence  $\tilde Z$ de l'élément (\ref{elem-slem0-rho-11j}) et donc 
  $\tilde Z\in (\Lambda_{Z})_{u_{1},u_{2},u_{3}}^{v_{1},v_{2},v_{3}}$. \cqfd

\begin{souslem}\label{slem1-rho-2j1518}
Pour $\tilde Z\in (\Lambda_{Z})_{u_{1},u_{2},u_{3}}^{v_{1},v_{2},v_{3}}$ et  \begin{gather*}(\tilde a_{1},\dots,\tilde a_{p},\tilde S_{0},...,\tilde S_{m},(\tilde {\mathcal Y}_{i}^{j})_{i\in \{0,\dots,m\}, j\in \{1,\dots ,l_{i}\}},(\tilde {\mathcal Z}_{0}^{j})_{j\in \{1,\dots ,6\}}) \\  \in 
(\pi_{x}^{\natural, p,k,m,( l_{0},...,l_{m}),6,0})^{-1}(\tilde Z),\end{gather*}  
$(\rho^{\flat}_{x})_{u_{1},u_{2},u_{3}}^{v_{1},v_{2},v_{3}}(\tilde S_{0})$ ne dépend que de la connaissance des points de 
\begin{gather}\label{ens-slem1-rho-11j}B(x,k+2M)\cup B(\tilde S_{0}, M) 
\cup  \bigcup_{j\in\{1,... ,6\}}  B(\tilde{\mathcal Z}_{0}^{j}, M)\end{gather} 
et des distances entre ces points. 
\end{souslem}
\noindent{\bf Démonstration.} Il suffit de montrer que 
  l'ensemble (\ref{eta-conn-deuxieme-partie}) figurant dans le lemme~\ref{dependance-rho'-123123} (avec $\tilde S_{0}$  au lieu de $S$) est inclus dans  (\ref{ens-slem1-rho-11j}). Soit $a\in \tilde S_{0}$, $j\in \{1,...,6\}$, et $y\in 3\de\tg (x,a)$ vérifiant 
$$|d(x,y)-w_{j}|\leq N+6\de+4.$$ Soit $z\in \geod(x,a)$ vérifiant $d(x,z)=w_{j}$, si bien que $z$ appartient à $\tilde{\mathcal Z}_{0}^{j}$. Comme 
$|d(x,y)-d(x,z)|\leq N+6\de+4$ et que $y,z\in 3\de\tg (x,a)$, $(H_{\de}(y,x,z,a))$ implique $d(y,z)\leq (N+6\de+4)+3\de+\de$. On suppose $(N+6\de+4)+3\de+\de\leq M$, ce qui est permis par $(H_{M})$. On a donc $d(y,\tilde{\mathcal Z}_{0}^{j})\leq M$. \cqfd

\begin{souslem}\label{slem2-rho-2j1519}
Le cardinal de $(\Lambda_{Z})_{u_{1},u_{2},u_{3}}^{v_{1},v_{2},v_{3}} $ est majoré par une constante de la forme $C(\de,K,N,Q,P,M)$. 
 \end{souslem}
 \noindent{\bf Démonstration.} Grâce au lemme~\ref{nombre-dist-connaitre-par-point-natural}, pour connaître les distances entre les points de 
\begin{gather}\label{ens1-p113-1j1138} \bigcup_{j\in\{1,... ,6\}}  B(\tilde{\mathcal Z}_{0}^{j}, M)\end{gather}  et ceux de 
 \begin{gather}\label{ens2-p113-24oct09}
\bigcup _{ i\in \{0,\dots ,m\}}
B(\tilde S_{i},M)
\cup  \bigcup _{i\in \{0,\dots,m\}, j\in \{1,\dots , l_{i}\}}B(\tilde {\mathcal Y}_{i}^{j},M)  \cup B(x,k+2M)\end{gather}
il suffit de connaître les distances entre les points de (\ref{ens1-p113-1j1138})  et $C$ points de (\ref{ens2-p113-24oct09}), avec $C=C(\de,K,N,Q,P,M)$ et ces distances sont déterminées à $C'=C(\de,K,N,Q,P,M)$ par $Z$ et $u_{1},u_{2},u_{3}, v_{1},v_{2},v_{3}  $
(en fait le lemme~\ref{nombre-dist-connaitre-par-point-natural} fournit  un énoncé légèrement différent où $
B(\tilde S_{0},M)$ figure dans  (\ref{ens1-p113-1j1138})  et non dans  (\ref{ens2-p113-24oct09}), mais il est clair que cet énoncé implique le nôtre). 
\cqfd

\noindent{\bf Fin de la démonstration du lemme~\ref{continuite-eta}.} 
D'après le sous-lemme~\ref{slem1-rho-2j1518}, pour $\tilde Z\in (\Lambda_{Z})_{u_{1},u_{2},u_{3}}^{v_{1},v_{2},v_{3}}$ on a 
$$\xi_{\tilde Z}\big(((\theta^{\flat}_{x})_{u_{1},u_{2},u_{3}}^{v_{1},v_{2},v_{3}}-\theta_{x})  (f)\big)=
(\alpha_{\tilde Z})_{u_{1},u_{2},u_{3}}^{v_{1},v_{2},v_{3}}
\xi_{\tilde Z}  (f)$$ avec 
$\big|(\alpha_{\tilde Z})_{u_{1},u_{2},u_{3}}^{v_{1},v_{2},v_{3}}
\big|\leq 7\de$. Le sous-lemme~\ref{slem0-rho-2j1517} montre que  pour $Z\in \overline Y_{x}^{p,k,m,(l_{0},...,l_{m})}$, 
$$\xi_{Z}\big(((\theta^{\flat}_{x})_{u_{1},u_{2},u_{3}}^{v_{1},v_{2},v_{3}}-\theta_{x})  (f)\big)=\sum_{\tilde Z\in (\Lambda_{Z})_{u_{1},u_{2},u_{3}}^{v_{1},v_{2},v_{3}}}
\xi_{\tilde Z}\big(((\theta^{\flat}_{x})_{u_{1},u_{2},u_{3}}^{v_{1},v_{2},v_{3}}-\theta_{x})  (f)\big).$$
Par Cauchy-Schwarz et grâce au sous-lemme~\ref{slem2-rho-2j1519}
 on en déduit  (\ref{ineg-eta-uuuvvv}). On a montré (\ref{ineg-eta-uuuvvv}) donc (\ref{cont-eta-etape}) et ceci termine la démonstration du lemme~\ref{continuite-eta}.   \cqfd

\subsection{Propriétés d'équivariance de la norme}\label{action-G-normes}

Ce sous-paragraphe a pour but de montrer la proposition~\ref{equiv-normeYZ}. 
On recommande au lecteur de commencer par lire le premier paragraphe   de~\cite{duke}, car la démonstration de la proposition~\ref{equiv-normeYZ}  repose sur les mêmes idées que  celle de la proposition 1.10 de~\cite{duke} mais sur des calculs  beaucoup plus compliqués. 

On fixe $p\in \{1,\pp,p_{\max}\}$ et $x,x'\in X$. 
On rappelle qu'après la définition~\ref{defi-Y}, pour $k,m\in \N$ et $  (l_{0},...,l_{m})\in \N^{m+1}$  on a introduit une relation d'équivalence sur $Y_{x}^{p,k,m,(l_{0},...,l_{m})}$,  et noté $\overline Y_{x}^{p,k,m,(l_{0},...,l_{m})}$ le quotient, et $\pi_{x}^{p,k,m,(l_{0},...,l_{m})}$ l'application quotient. On va introduire maintenant une relation d'équivalence plus fine
$$Y_{x}^{p,k,m,(l_{0},...,l_{m})}\vad{\pi_{x,x'}^{p,k,m,(l_{0},...,l_{m})}} \overline Y_{x,x'}^{p,k,m,(l_{0},...,l_{m})}$$
 et une autre encore plus fine 
$$Y_{x}^{p,k,m,(l_{0},...,l_{m})}\vad{\pi_{x,x',\star}^{p,k,m,(l_{0},...,l_{m})}} \overline Y_{x,x',\star}^{p,k,m,(l_{0},...,l_{m})}$$
telles que  tout élément de $ \overline Y_{x}^{p,k,m,(l_{0},...,l_{m})}$
 a au plus $2d(x,x')+1$
antécédents dans $ \overline Y_{x,x'}^{p,k,m,(l_{0},...,l_{m})}$ et tout 
 élément de $ \overline Y_{x,x'}^{p,k,m,(l_{0},...,l_{m})}$ a au plus $C$  antécédents dans 
$ \overline  Y_{x,x',\star}^{p,k,m,(l_{0},...,l_{m})}$, pour une certaine constante $C=C(\de,K,N,Q,P,M)$. 

On définit la relation d'équivalence 
$$Y_{x}^{p,k,m,(l_{0},...,l_{m})}\vad{\pi_{x,x'}^{p,k,m,(l_{0},...,l_{m})}} \overline Y_{x,x'}^{p,k,m,(l_{0},...,l_{m})}$$ de la fa\c con suivante : 
   $$(a_{1},\dots,a_{p},S_{0},...,S_{m},(\mathcal Y_{i}^{j})_{i\in \{0,\dots,m\}, j\in \{1,\dots ,l_{i}\}})$$ et $$(\hat a_{1},\dots,\hat a_{p},\hat S_{0},...,\hat S_{m},(\hat {\mathcal Y}_{i}^{j})_{i\in \{0,\dots,m\}, j\in \{1,\dots ,l_{i}\}})$$
ont même image dans $\overline Y_{x,x'}^{p,k,m,(l_{0},...,l_{m})}$ s'ils ont même image dans $\overline Y_{x}^{p,k,m,(l_{0},...,l_{m})}$ et si $d(x',S_{0})=d(x',\hat S_{0})$. L'application 
$$(a_{1},\dots,a_{p},S_{0},...,S_{m},(\mathcal Y_{i}^{j})_{i\in \{0,\dots,m\}, j\in \{1,\dots ,l_{i}\}})\mapsto d(x',S_{0})$$ se factorise donc en une application $r_{0}':\overline Y_{x,x'}^{p,k,m,(l_{0},...,l_{m})}\to \N$. 
On définit une application  $$k': \overline Y_{x,x'}^{p,k,m,(l_{0},...,l_{m})}\to \N$$ en posant, pour  
 $Z\in \overline Y_{x,x'}^{p,k,m,(l_{0},...,l_{m})}$ \label{def-r0'-k'-xx'}
 %et 
%$$(a_{1},\pp,a_{p},S_{1},...,S_{m},(\mathcal Y_{i}^{j})_{i\in \{0,\dots,m\}, j\in \{1,\dots ,l_{i}\}})\in (\pi_{x,x'}^{p,k,m,(l_{0},...,l_{m})} )^{-1}(Z)$$
%le plus grand entier $l$ tel que les distances de $a_{1},\pp,a_{p}$ à $B(x',l)$ déterminent les distances de $a_{1},\pp,a_{p}$  à $B(x,k+2M)$ (c'est-à-dire plus précisément pour tout $i\in \{1,\pp,p\}$ et pour tout $z\in B(x,k+2M)$, $\geod(a_{i},z)$ rencontre $B(x',l)$) vérifie $l\in [k'(Z)-M-N-2\de,k'(Z)]$ (cependant nous montrerons seulement $l\leq k'(Z)$ dans le lemme~\ref{estim-1-k'} car l'autre inégalité ne nous sera pas utile, et le cas des arbres discuté ci-dessus suffit à justifier que la formule pour  $k'(Z)$ 
%est la bonne). 
%
%
%On définit $k'(Z)$ par la relation suivante 
%$$r_{0}'(Z)-k'(Z)=\min(r_{0}(Z)-k,E(\frac{r_{0}(Z)+r_{0}'(Z)-d(x,x')}{2}))+10\de.$$
\begin{gather}\nonumber k'(Z)=\min\big(r'_{0}(Z),\max(  r'_{0}(Z)-r_{0}(Z)+k,    E(\frac{r'_{0}(Z)+d(x,x')-r_{0}(Z)}{2}))\big)\\ \label{def-k'Z-24oct09}+\frac{M}{2}.\end{gather}
Voici le dessin dans le cas où $X$ est un arbre et $p=1$ (si bien que $S_{0}=\{a_{1}\}$ et on note $b=a_{1}$ pour que le dessin serve de nouveau dans la suite) et où l'on prend $M=0$ dans la formule précédente.   Dans ce dessin  on a choisi $$(b,\{b\})\in (\pi_{x,x'}^{1,k,0,(0)} )^{-1}(Z),$$
on a représenté $B(x,k)$ par une moitié de boule et on a noté $u$ le point de $B(x,k)$ à distance minimale de $b$. 
On définit $t$ comme le point de $$ \bigcap _{\tilde x\in B(x,k)} \geod(\tilde x,b)$$ le plus proche de $x'$ (ce point est unique et appartient à $\geod(x',b)$). 
On remarque que $t$ est aussi le point central du triangle $ubx'$. 
La formule précédente pour $k'(Z)$ devient 
$$k'(Z)=\min\big(d(x',b),\max(d(x',b)-d(x,b)+k, 
\frac{d(x',b)+d(x,x')-d(x,b)}{2})\big)$$ $$=d(x',t).$$
%$$d(x',b)-k'(Z)=\min(d(x,b)-k,\frac{d(x,b)+d(x',b)-d(x,x')}{2}))$$
%$$=\min(d(u,b),\frac{d(x,b)+d(x',b)-d(x,x')}{2})).$$
%$$d(t,b)=d(x',b)-k'(Z)=\min(d(x,b)-k,\frac{d(x,b)+d(x',b)-d(x,x')}{2}))$$
%$$=\min(d(u,b),\frac{d(x,b)+d(x',b)-d(x,x')}{2})).$$
 La boule $B(x',k'(Z))$
 (dont le bord contient $t$) est la boule de centre $x'$ de plus petit rayon telle que toute géodésique entre $b$ et un point de $B(x,k)$ la rencontre. On remarque que $t$ dépend de $S_{0}=\{b\}$ mais que $d(x',t)$ ne dépend que de $Z$. Dans le dessin ci-dessous, le premier cas est le cas où $b\in B(x,k)$ et le   deuxième cas et le troisième cas sont distingués par l'appartenance ou non de $u$ à $\geod(x,x')$ (dans le premier et le deuxième cas $x'$ n'appartient pas forcément à $B(x,k)$). 
  
\ifx\JPicScale\undefined\def\JPicScale{1}\fi
\unitlength \JPicScale mm


On revient maintenant au cas général et on va distinguer trois cas comme dans le dessin pour les arbres. Soit $Z\in \overline Y_{x,x'}^{p,k,m,(l_{0},...,l_{m})}$

\noindent {\bf Premier cas.} On suppose $r_{0}(Z)\leq k$.

Alors, par (\ref{def-k'Z-24oct09}),  
 \begin{gather}\label{form-r'-cas1-4j} k'(Z)=r'_{0}(Z)+\frac{M}{2}. \end{gather}

\noindent {\bf Deuxième cas.} On suppose $r_{0}(Z)> k$ et 
$r_{0}(Z)-r_{0}'(Z)+d(x,x')\leq 2k$. 

Cette dernière condition implique 
$$r'_{0}(Z)-r_{0}(Z)+k\geq E(\frac{r'_{0}(Z)-r_{0}(Z)+d(x,x')}{2})$$
d'où, par (\ref{def-k'Z-24oct09}),  
  \begin{gather}\label{form-r'-cas2-4j}
 k'(Z)=r'_{0}(Z)-r_{0}(Z)+k+\frac{M}{2}.\end{gather}

\noindent {\bf Troisième cas.}
 On suppose $r_{0}(Z)> k$ et 
$r_{0}(Z)-r_{0}'(Z)+d(x,x')>  2k$. 

Cette dernière condition implique 
$$r'_{0}(Z)-r_{0}(Z)+k\leq E(\frac{r'_{0}(Z)-r_{0}(Z)+d(x,x')}{2})$$
et comme $d(x,x')\leq r_{0}(Z)+r'_{0}(Z)+N$ on en déduit, par (\ref{def-k'Z-24oct09}),  
 \begin{gather}\label{form-r'-cas3-4j}\Big| k'(Z)-\Big(E\big(\frac{r'_{0}(Z)-r_{0}(Z)+d(x,x')}{2}\big)+\frac{M}{2}\Big)\Big|\leq \frac{N}{2}.
 \end{gather}

On introduit maintenant une partition de $Y_{x}^{p,k,m,(l_{0},...,l_{m})}$ pour la relation suivante : $$(a_{1},\dots,a_{p},S_{0},...,S_{m},(\mathcal Y_{i}^{j})_{i\in \{0,\dots,m\}, j\in \{1,\dots ,l_{i}\}})$$ et $$(\hat a_{1},\dots,\hat a_{p},\hat S_{0},...,\hat S_{m},(\hat {\mathcal Y}_{i}^{j})_{i\in \{0,\dots,m\}, j\in \{1,\dots ,l_{i}\}})$$
sont en relation s'ils ont même image $Z$ dans $\overline Y_{x,x'}^{p,k,m,(l_{0},...,l_{m})}$ par $\pi_{x,x'}^{p,k,m,(l_{0},...,l_{m})}$ et 
s'il existe une isométrie de $$
\bigcup _{ i\in \{0,\dots ,m\}}
B(S_{i}, M)
\cup  \bigcup _{\substack{i\in \{0,\dots,m\},\\  j\in \{1,\dots ,l_{i}\}}} B(\mathcal Y_{i}^{j},M)\cup B(x,k+2M)\cup B(x',k'(Z)+2M)$$ vers 
$$\bigcup _{ i\in \{0,\dots ,m\}}
B(\hat S_{i}, M)
\cup  \bigcup _{\substack{i\in \{0,\dots,m\},\\  j\in \{1,\dots ,l_{i}\}}} B(\hat {\mathcal Y}_{i}^{j}, M)\cup B(x,k+2M)\cup B(x',k'(Z)+2M)$$ 
 qui envoie 
$a_{i}$ sur $\hat a_{i}$ pour $i\in \{1,\dots,p\}$,  
 $S_{i}$ sur $\hat S_{i}$
 pour $i\in \{0,\dots ,m\}$, $\mathcal Y_{i}^{j}$ sur $\hat {\mathcal Y}_{i}^{j}$  pour 
 $i\in \{0,\dots,m\}, j\in \{1,\dots ,l_{i}\}$  
 et est l'identité sur $B(x,k+2M)$ et sur $B(x',k'(Z)+2M)$.

On note $\overline Y_{x,x',\star}^{p,k,m,(l_{0},...,l_{m})}$
le quotient de $Y_{x}^{p,k,m,(l_{0},...,l_{m})}$ pour cette relation d'équivalence, et $\pi_{x,x',\star}^{p,k,m,(l_{0},...,l_{m})}$ l'application quotient.

\begin{lem}\label{est-card-fibres-star}
a) Les fibres de l'application surjective $$\overline Y_{x,x'}^{p,k,m,(l_{0},...,l_{m})}\to\overline Y_{x}^{p,k,m,(l_{0},...,l_{m})} $$ sont toutes de cardinal $\leq 2d(x,x')+1$. 

\noindent b) Il existe une constante $C=C(\de,K,N,Q,P,M)$ telle que les fibres de l'application surjective $$\overline Y_{x,x',\star}^{p,k,m,(l_{0},...,l_{m})}\to \overline Y_{x,x'}^{p,k,m,(l_{0},...,l_{m})} $$ soient toutes de cardinal $\leq C$. 
\end{lem}
Avant de montrer ce lemme, expliquons la stratégie de la 
preuve  de la proposition~\ref{equiv-normeYZ}. On veut majorer $\|f\|_{\H_{x,s}(\Delta_{p})}^{2}$ en fonction de $\|f\|_{\H_{x',s}(\Delta_{p})}^{2}$. 
D'abord $\|f\|_{\H_{x,s}(\Delta_{p})}^{2}$ est une somme pondérée indexée par certaines classes d'équivalence $Z\in \overline Y_{x}^{p,k,m,(l_{0},...,l_{m})}$ des carrés des formes linéaires $$\xi_{Z}:f\mapsto \sum _{(a_{1},\dots,a_{p},S_{0},...,S_{m},(\mathcal Y_{i}^{j})_{i\in \{0,\dots,m\}, j\in \{1,\dots ,l_{i}\}}) \in 
(\pi_{x}^{p,k,m,(l_{0},...,l_{m})})^{-1}(Z)} f(a_1,...,a_p) .$$ Grâce au lemme précédent on peut couper chaque classe d'équivalence $Z$ en  morceaux $\tilde Z$ (paramétrés par $\overline Y_{x,x',\star}^{p,k,m,(l_{0},...,l_{m})}$) dont le nombre est inférieur ou égal à $C(2d(x,x')+1)$. Par Cauchy-Schwarz on a donc pour tout $f$, $|\xi_{Z}(f)|^{2}\leq 
C(2d(x,x')+1) \sum_{\tilde Z }|\xi_{\tilde Z}(f)|^{2}$
où la somme porte sur les $\tilde Z\in \overline Y_{x,x',\star}^{p,k,m,(l_{0},...,l_{m})}$ dont l'image dans $\overline Y_{x}^{p,k,m,(l_{0},...,l_{m})}$ est $Z$. On montrera ensuite (dans le lemme~\ref{oubli-x-xx'-x'-Z}) que chaque forme linéaire $\xi_{\tilde Z}$ est proportionnelle à une forme linéaire $\xi_{Z'}$ apparaissant dans la formule pour $\|f\|_{\H_{x',s}(\Delta_{p})}^{2}$. 
Par définition  $\tilde Z$ est la donnée des distances entre les points de $B(x,k+2M)$,  ceux de $B(x',k'+2M)$ et les points à distance $\leq M$ de la réunion de toutes les  parties $S_{i}$ et $\mathcal Y_{i}^{j}$. Le dessin ci-dessous correspond au troisième cas du dessin pour les arbres (qui est le cas le plus intéressant) et  les deux boules y représentent $B(x,k)$ et  $B(x',k')$. Les parties $S_{i}$ sont représentées par des cercles plus grands que les parties $\mathcal Y_{i}^{j}$ dans un souci de clarté, bien que leur diamètre maximal soit plus petit ($N$ au lieu de $P$).

\ifx\JPicScale\undefined\def\JPicScale{1}\fi
\unitlength \JPicScale mm


\noindent 
La classe d'équivalence $Z'$  sera la donnée des distances entre les points de $B(x',k'+2M)$ et les points à distance $\leq M$ de la réunion des parties $S_{i}$ et $\mathcal Y_{i}^{j}$ situées à droite de la boule $B(x',k')$, qui sont indiquées sur le dessin comme ``parties conservées''. On montrera qu'il existe $m'$ tel que les parties conservées soient les $S_{i}$ pour $i\in \{0,...,m'\}$, tous les $\mathcal Y_{i}^{j}$ pour $i\in \{0,...,m'-1\}$ et certains des $\mathcal Y_{m'}^{j}$. 
La raison pour laquelle $\xi_{\tilde Z}$ est proportionnelle à $\xi_{Z'}$
est que la connaissance de $\tilde Z$ et des parties $S_{i}$ et $\mathcal Y_{i}^{j}$ éliminées (qui se trouvent dans $B(x,k)$ ou $B(x',k')$ ou  entre ces deux boules) détermine $Z$, car 
si $\Gamma$ est une partie conservée et $\Delta $ une partie éliminée
 toute géodésique entre $B(\Gamma,M)$ et   $B(\Delta,M)\cup B(x,k+2M)$ traverse $ B(x',k'+2M)$, comme on le verra dans le lemme~\ref{oubli-x-xx'-x'}. 
Par conséquent, pour connaître les distances entre $B(\Gamma,M)$ et   $B(\Delta,M)\cup B(x,k+2M)$ , il suffit de connaître les distances entre 
$B(\Gamma,M)$ et  $ B(x',k'+2M)$, qui  font partie de la donnée de $\tilde Z$. 
La difficulté sera ensuite que chaque $\xi_{Z'}$ est proportionnelle à  $\xi_{\tilde Z}$ pour une infinité de   $\tilde Z\in \bigcup_{k,m,(l_{0},...,l_{m})}\overline Y_{x,x',\star}^{p,k,m,(l_{0},...,l_{m})} $ et on devra vérifier que, compte tenu des pondérations et grâce à $(H_{B})$,  la somme sur $\tilde Z$ n'introduit pas de divergence. 

\noindent {\bf Démonstration du lemme~\ref{est-card-fibres-star}.} Le a) vient simplement du fait que pour 
$$(a_{1},\dots,a_{p},S_{0},...,S_{m},(\mathcal Y_{i}^{j})_{i\in \{0,\dots,m\}, j\in \{1,\dots ,l_{i}\}})\in Y_{x}^{p,k,m,(l_{0},...,l_{m})}$$ on a $d(x',S_{0})\in [d(x,S_{0})-d(x,x'),d(x,S_{0})+d(x,x')]$. 
Montrons b). Soit $Z\in \overline Y_{x,x'}^{p,k,m,(l_{0},...,l_{m})}$. On doit montrer que le nombre d'antécédents de $Z$ dans $\overline Y_{x,x',\star}^{p,k,m,(l_{0},...,l_{m})}$ est majoré par une constante $C=C(\de,K,N,Q,P,M)$. On note $u$ un point de $B(x,k)$ tel que pour tout $$(a_{1},\dots,a_{p},S_{0},...,S_{m},(\mathcal Y_{i}^{j})_{i\in \{0,\dots,m\}, j\in \{1,\dots ,l_{i}\}})\in (\pi_{x,x'}^{p,k,m,(l_{0},...,l_{m})})^{-1}(Z)$$  $u$ soit à distance minimale de $S_{0}$ (cela ne dépend que de $Z$  car les distances entre les points de  $S_{0}$ et ceux de $B(x,k)$ font partie de la donnée  de $Z$). On a alors  $d(u,S_{0})=\max(0,r_{0}(Z)-k)$. On rappelle que d'après le lemme~\ref{lemme-S0-...Sm}, pour tout 
$$(a_{1},\dots,a_{p},S_{0},...,S_{m},(\mathcal Y_{i}^{j})_{i\in \{0,\dots,m\}, j\in \{1,\dots ,l_{i}\}})\in (\pi_{x,x'}^{p,k,m,(l_{0},...,l_{m})})^{-1}(Z)$$ 
et pour  $b\in S_{0}$  à distance minimale de $u$ (si bien que $u$ est un point de $B(x,k)$ à distance minimale de $b$) on a 
 $$S_{0}\cup \pp\cup S_{m}\cup \bigcup _{i\in \{0,\dots,m\}, j\in \{1,\dots ,l_{i}\}}\mathcal Y_{i}^{j}\subset 
4P\tg(b,u) $$ d'où par le lemme~\ref{xx'yy'zz'}, 
$$\bigcup _{ i\in \{0,\dots ,m\}}
B(S_{i}, M)
\cup  \bigcup _{i\in \{0,\dots,m\}, j\in \{1,\dots ,l_{i}\}} B(\mathcal Y_{i}^{j},M)\subset (2M+4P)\tg(b,u). $$ 
On va distinguer trois cas qui correspondent à peu près aux trois cas envisagés dans le dessin pour les arbres (auquel le lecteur peut se reporter 
pour lire la démonstration). 

\noindent {\bf Premier cas.} On suppose $r_{0}(Z)\leq k$. 

Soit 
$$(a_{1},\dots,a_{p},S_{0},...,S_{m},(\mathcal Y_{i}^{j})_{i\in \{0,\dots,m\}, j\in \{1,\dots ,l_{i}\}})\in (\pi_{x,x'}^{p,k,m,(l_{0},...,l_{m})})^{-1}(Z)$$ 
et  $b\in S_{0}$  à distance minimale de $u$. 
On a alors $b=u$ d'où 
\begin{gather*}\bigcup _{ i\in \{0,\dots ,m\}}
B(S_{i}, M)
\cup  \bigcup _{i\in \{0,\dots,m\}, j\in \{1,\dots ,l_{i}\}} B(\mathcal Y_{i}^{j}, M)\\ \subset (2M+4P)\tg(u,u)=  B(u, M+2P)\subset B(x,k+2M)\end{gather*}
où la dernière inclusion vient de l'inégalité $2P\leq M$ que l'on suppose  grâce à $(H_{M})$. Donc $(\pi_{x,x'}^{p,k,m,(l_{0},...,l_{m})})^{-1}(Z)$ est un singleton et a fortiori l'image inverse de $Z$ par l'application 
$\overline Y_{x,x',\star}^{p,k,m,(l_{0},...,l_{m})}\to \overline Y_{x,x'}^{p,k,m,(l_{0},...,l_{m})} $
 est un singleton. 

\noindent {\bf Deuxième cas.} On suppose $r_{0}(Z)> k$ et 
$r_{0}(Z)-r_{0}'(Z)+d(x,x')\leq 2k$. 
%Cette dernière condition implique 
%$$r'_{0}(Z)-r_{0}(Z)+k\geq E(\frac{r'_{0}(Z)-r_{0}(Z)+d(x,x')}{2})$$
%d'où, par (\ref{def-k'Z-24oct09}), 

Par (\ref{form-r'-cas2-4j}) on a 
 $$k'(Z)=r'_{0}(Z)-r_{0}(Z)+k+\frac{M}{2}.$$
Soit 
$$(a_{1},\dots,a_{p},S_{0},...,S_{m},(\mathcal Y_{i}^{j})_{i\in \{0,\dots,m\}, j\in \{1,\dots ,l_{i}\}})\in (\pi_{x,x'}^{p,k,m,(l_{0},...,l_{m})})^{-1}(Z)$$ 
et   $b\in S_{0}$  à distance minimale de $u$. On a donc
\begin{gather}\label{ineg-xbx'-12dec}d(x,b)-d(x',b)+d(x,x')\leq 2k+N, \ \  d(x,b)> k, \\ \nonumber 
\big|      k'(Z)-(d(x',b)-d(x,b)+k+\frac{M}{2})  \big|\leq N.\end{gather}
  Soit $t\in B(x',k'(Z)+2M)$ à distance minimale de $b$. 
    On a $u\in \geod(x,b)$, et $d(x,u)=k$, $d(u,b)=d(x,b)-k$. 
Par $(H_{\de}^{0}(x',x,u,b))$ et (\ref{ineg-xbx'-12dec}) on a \begin{gather*}d(x',u)\leq \max(d(x,x')-k,d(x',b)-d(x,b)+k)+\de\\ \leq d(x',b)-d(x,b)+k+N+\de=d(x',b)-d(u,b)+N+\de\end{gather*}
donc $u\in (N+\de)\tg(x',b)$. 
On a $t\in \geod(x',b)$, $$d(x',t)=\min(d(x',b),k'(Z)+2M)$$
et $|k'(Z)-(d(x',b)-d(u,b)+\frac{M}{2})|\leq N$ donc 
\begin{gather*}d(x',t)\leq k'(Z)+2M \leq  d(x',u)+\frac{5M}{2}+N,\\ 
%$$d(x',t)\leq \min(d(x',b),d(x',u)+\frac{5M}{2}+2N+\de),$$
d(x',t)\geq \min(d(x',b),d(x',u)+\frac{5M}{2}-2N-\de),\end{gather*}
et comme $d(x',u)\leq d(x',b)+N+\de$,   on en déduit 
\begin{gather}\label{ineg-x'tu-14dec}|d(x',t)-d(x',u)|\leq \frac{5M}{2}+N\end{gather}
car on suppose $N+\de\leq \frac{5M}{2}$, ce qui est permis grâce à $(H_{M})$. 

Comme $u\in (N+\de)\tg(x',b)$ et $t\in \geod(x',b)$, et grâce à (\ref{ineg-x'tu-14dec}),  $(H_{\de}(u,x',t,b))$ implique alors  
$$d(u,t)\leq (\frac{5M}{2}+N)+(N+\de)+\de=\frac{5M}{2}+2N+2\de.$$
Cette inégalité reflète l'égalité $u=t$ dans le  dessin pour les arbres qui correspond au deuxième cas. 

On en déduit 
$$(2M+4P)\tg(u,b)\subset \alpha\tg(t,b)$$ avec $\alpha=2(\frac{5M}{2}+2N+2\de)+(2M+4P)$. 
Le lemme~\ref{distances-Bxk-C} montre que  les distances entre un point de $\alpha\tg(t,b)$  et les points de 
$B(x',k'(Z)+2M)$ sont déterminées par les distances entre ce point et les points de 
\begin{gather*}B(x',k'(Z)+2M)\cap B(t,\alpha+4\de) \subset B(t,\alpha+4\de)\subset 
B(u,\beta) \\ \text{avec \ \ } \beta=(\frac{5M}{2}+2N+2\de)+\alpha+4\de=
3 (\frac{5M}{2}+2N+2\de)+2M+4P+4\de.\end{gather*}
Donc  les distances entre 
les points de l'ensemble  \begin{gather}\label{ens-prem-cas-12dec}\bigcup _{ i\in \{0,\dots ,m\}}
B(S_{i}, M)
\cup  \bigcup _{i\in \{0,\dots,m\}, j\in \{1,\dots ,l_{i}\}} B(\mathcal Y_{i}^{j}, M)\end{gather}
(qui est inclus dans $(4P+2M)\tg(u,b) \subset \alpha\tg(t,b)$) et les points de $B(x',k'(Z)+2M)$ sont déterminées par les distances entre les  points de (\ref{ens-prem-cas-12dec}) et les points  de $B(u,\beta)$. 
D'autre part le  c) du lemme~\ref{i<j<k-distances}, 
appliqué à 
\begin{itemize}
\item $(u,b)$ au lieu de $(c,d)$, 
\item $\{w_{i}, i\in I\}$ égal à $\{u\}\cup \bigcup _{ i\in \{0,\dots ,m\}}
S_{i}\cup  \bigcup _{i\in \{0,\dots,m\}, j\in \{1,\dots ,l_{i}\}} \mathcal Y_{i}^{j}$,
\item $J\subset I$ le singleton tel que $\{w_{j}, j\in J\}=\{u\}$, 
 \item $(\alpha_{i},\rho_{i})$ égal à $(0,\beta)$ pour $i\in J$ et à 
 %\item $w_{1}=\{u\}$ et 
 %$\{w_{2},...,w_{r}\}$ égal à la réunion disjointe $$\bigcup _{ i\in \{0,\dots ,m\}}
%S_{i}\cup  \bigcup _{i\in \{0,\dots,m\}, j\in \{1,\dots ,l_{i}\}} \mathcal Y_{i}^{j},$$ \item $(\alpha_{1},\rho_{1})=(0,\beta)$ et $(\alpha_{j},\rho_{j})$ égal à 
$(4P,M)$ pour $i\in I\setminus J$
%$j\neq 1$, et $I=\{1\}$, 
\end{itemize}
 montre  que les distances entre 
les points de l'ensemble (\ref{ens-prem-cas-12dec}) et ceux de $B(u,\beta)$ sont déterminées par les distances entre $C=C(\de,K,N,Q,P,M)$ points de (\ref{ens-prem-cas-12dec}) (déterminés par la donnée de $Z$) et les points  de $B(u,\beta)$. 
 De plus  le cardinal de $B(u,\beta)$ est majoré par une  constante de la forme $C(\de,N,K,Q,P,M)$ et les distances entre ces $C$ points de (\ref{ens-prem-cas-12dec}) et les points  de $B(u,\beta)$ sont déterminées 
à $\beta$ près par la donnée de $Z$. 
Cela termine l'étude du  deuxième cas. 

\noindent {\bf Troisième cas.}
 On suppose $r_{0}(Z)> k$ et 
$r_{0}(Z)-r_{0}'(Z)+d(x,x')>  2k$. 
%Cette dernière condition implique 
%$$r'_{0}(Z)-r_{0}(Z)+k\leq E(\frac{r'_{0}(Z)-r_{0}(Z)+d(x,x')}{2})$$
%et comme $d(x,x')\leq r_{0}(Z)+r'_{0}(Z)+N$ on en déduit, par (\ref{def-k'Z-24oct09}),   

Par (\ref{form-r'-cas3-4j}) on a $$\Big| k'(Z)-\Big(E\big(\frac{r'_{0}(Z)-r_{0}(Z)+d(x,x')}{2}\big)+\frac{M}{2}\Big)\Big|\leq \frac{N}{2}.$$
Soit 
$$(a_{1},\dots,a_{p},S_{0},...,S_{m},(\mathcal Y_{i}^{j})_{i\in \{0,\dots,m\}, j\in \{1,\dots ,l_{i}\}})\in (\pi_{x,x'}^{p,k,m,(l_{0},...,l_{m})})^{-1}(Z)$$ 
et   $b\in S_{0}$  à distance minimale de $u$. On a donc
\begin{gather*}d(x,b)-d(x',b)+d(x,x')\geq 2k-N, \ \  d(x,b)> k, \\ \big|      k'(Z)-(E(\frac{d(x',b)-d(x,b)+d(x,x')}{2})+\frac{M}{2})  \big|\leq N.\end{gather*}
  Comme 
    $d(x,b)-d(x',b)+d(x,x')\geq 2k-N$, $(H_{\de}^{0}(x',x,u,b))$ donne 
    \begin{gather*}d(x',u)\leq \max(d(x',x)-k, d(x',b)-d(x,b)+k)+\de\leq d(x,x')-k+N+\de\\ =d(x,x')-d(x,u)+N+\de\text{\ \  donc \ \ }u\in (N+\de)\tg(x',x).\end{gather*} 
Donc 
 \begin{gather}\label{ineg-k'-xx'bu-14dec}\big|      k'(Z)-(E(\frac{d(x',b)+d(x',u)-d(u,b)}{2})+\frac{M}{2})  \big|\leq 2N+\de.\end{gather}

  Soit 
  $u'\in B(x',k'(Z)+2M)$ à distance minimale de $u$ et 
    $b'\in B(x',k'(Z)+2M)$ à distance minimale de $b$. 
  Soit $w\in 4P\tg(u,b)$. On applique le lemme~\ref{approx-arbres} à $u,b,x',u',b',w$ avec $l=2$ et $x'$ comme point base. On note $t$ le point central du triangle de sommets $\Psi u$, $\Psi b$, $\Psi x'$. 
 Dans le dessin ci-dessous le point de $\geod(\Psi u, \Psi b)$ le plus proche de  $\Psi w$ est arbitraire.   
  
  \ifx\JPicScale\undefined\def\JPicScale{1}\fi
\unitlength \JPicScale mm
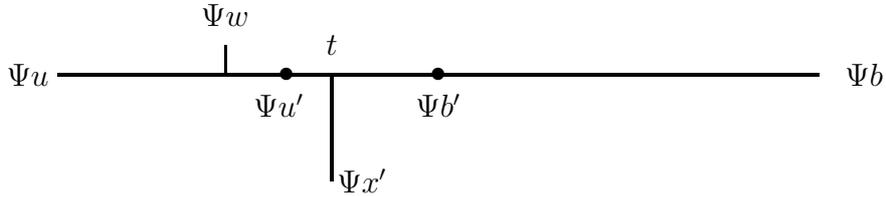
\begin{picture}(130,30)(30,13)
\linethickness{0.3mm}
\put(40,30){\line(1,0){100}}
\linethickness{0.3mm}
\put(76,16){\line(0,1){14}}
\linethickness{0.3mm}
\put(62,30){\line(0,1){4}}
\linethickness{0.3mm}
\put(36,30){\makebox(0,0)[cc]{$\Psi u$}}

\put(146,30){\makebox(0,0)[cc]{$\Psi b$}}

\put(80,16){\makebox(0,0)[cc]{$\Psi x'$}}

\put(76,34){\makebox(0,0)[cc]{$t$}}

\put(62,38){\makebox(0,0)[cc]{$\Psi w$}}

\put(70,30){\makebox(0,0)[cc]{$\bullet$}}

\put(90,30){\makebox(0,0)[cc]{$\bullet$}}

\put(69,26){\makebox(0,0)[cc]{$\Psi u'$}}

\put(90,26){\makebox(0,0)[cc]{$\Psi b'$}}

\end{picture}

  \noindent 
 On a \begin{gather}\nonumber d(\Psi x',t)= \frac{d(\Psi x',\Psi u)+d(\Psi x',\Psi b)-d(\Psi u,\Psi b)}{2}\\ \label{ineg-phix't-15dec}\in [\frac{d(x', u)+d( x', b)-d( u, b)}{2}, \frac{d(x', u)+d( x', b)-d( u, b)}{2}+\de].\end{gather} Grâce au  a) du lemme~\ref{approx-arbres2} on a de plus \begin{gather*}\Psi u'\in \geod(\Psi x', \Psi u), \ \ d(\Psi x',\Psi u')=\min
  (k'(Z)+2M, d(\Psi x',\Psi u))\\ \text{et\ \ \ }\Psi b'\in \geod(\Psi x', \Psi b),\ \ \ d(\Psi x',\Psi b')=\min
  (k'(Z)+2M, d(\Psi x',\Psi b)).\end{gather*}
  Par (\ref{ineg-k'-xx'bu-14dec}) on a $k'(Z)+2M\geq \frac{d(x', u)+d( x', b)-d( u, b)}{2}+\de$ car on suppose $\frac{5M}{2}\geq 2N+2\de+1$, ce qui est permis par $(H_{M})$. On  déduit alors de (\ref{ineg-phix't-15dec}) que 
  \begin{gather*}\Psi u'\in \geod(t, \Psi u), \ \ d(t,\Psi u')\leq \frac{5M}{2}+2N+\de \\ \text{et \ }
  \Psi b'\in \geod(t, \Psi b),  \ \ d(t,\Psi b')\leq \frac{5M}{2}+2N+\de.\end{gather*}
 Grâce au  a) du lemme~\ref{approx-arbres2}  on a $\Psi w\in (4P+2\de)\tg(\Psi b, \Psi u) $ donc $\Psi w$ est à distance $\leq 2P+\de$ de $\geod(\Psi b, \Psi u) $. Donc l'une au moins des assertions suivante est vraie (suivant que le point de $\geod(\Psi b, \Psi u) $
le plus proche de $\Psi w$ appartient à $\geod(\Psi b, t) $ ou $\geod(t, \Psi u) $) : 

 \begin{itemize}
 \item $\Psi w$ est à distance $\leq \frac{5M}{2}+2N+2P+2\de $ de $\geod(\Psi b', \Psi b) $ 
 \item ou $\Psi w$ est à distance $\leq \frac{5M}{2}+2N+2P+2\de $ de $\geod(\Psi u', \Psi u) $
 \end{itemize}
  Il résulte facilement de ce qui précède et du b) du   lemme~\ref{approx-arbres2} que l'une au moins des assertions suivante est vraie 
 \begin{itemize}
 \item $w\in (5M+4N+4P+8\de)\tg(b',b) $,  
 \item ou $w\in (5M+4N+4P+8\de)\tg(u',u) $. 
 \end{itemize}
 On suppose 
 $4N+4P+8\de\leq M$, ce qui est possible par $(H_{M})$. 
 Il résulte de ce qui précède que $w\in 6M\tg(b',b) \cup 6M\tg(u',u) $.   
 On applique ceci à $w\in \bigcup _{ i\in \{0,\dots ,m\}}
S_{i} \cup 
 \bigcup _{i\in \{0,\dots,m\}, j\in \{1,\dots ,l_{i}\}}\mathcal Y_{i}^{j}\subset 4P\tg(u,b)$. 
 Donc l'ensemble   \begin{gather}\label{reunion-Mboules-14dec}\bigcup _{ i\in \{0,\dots ,m\}}
B(S_{i}, M) 
\cup  \bigcup _{i\in \{0,\dots,m\}, j\in \{1,\dots ,l_{i}\}} B(\mathcal Y_{i}^{j}, M) \end{gather}
est inclus dans $8M\tg(b',b) \cup 8M\tg(u',u) $. 
 
  En appliquant le lemme~\ref{distances-Bxk-C} à $\alpha=8M$, $(u,u')$  (resp.  $(b,b')$) au lieu de $(b,u)$  et $(x',k'(Z)+2M)$ au lieu de $(x,l)$ on voit que les distances entre n'importe quel point   $y\in 8M\tg(u',u)$ (resp. $y\in 8M\tg(b',b)$)  et les points de $B(x',k'(Z)+2M)$ sont déterminées par les distances entre $y$ et les points de \begin{gather}\label{incl1-15dec}B(x',k'(Z)+2M)\cap B(u',8M+4\de)\subset B(u',8M+4\de)\\ \label{incl2-15dec} \text{et }B(x',k'(Z)+2M)\cap B(b',8M+4\de)\subset B(b',8M+4\de)\end{gather}
 respectivement. 
 De plus $$d(u',b')\leq d(\Psi u',\Psi b')+2\de\leq d(\Psi u',t)+d(t,\Psi b')
 +2\de\leq5M+4N+4\de$$ donc 
 $B(b',8M+4\de)\subset B(u',13M+4N+8\de)$. On a $\Psi u'\in \geod(\Psi u,\Psi b)$ donc $ u'\in 4\de\tg(u,b)$ par le b)   du lemme~\ref{approx-arbres2}. 
 En appliquant le   c) du lemme~\ref{i<j<k-distances}  à 
\begin{itemize}
\item $(u,b)$ au lieu de $(c,d)$, 

\item $\{w_{i}, i\in I\}$ égal à $\{u'\}\cup \bigcup _{ i\in \{0,\dots ,m\}}
S_{i}\cup  \bigcup _{i\in \{0,\dots,m\}, j\in \{1,\dots ,l_{i}\}} \mathcal Y_{i}^{j}$,
\item $J\subset I$ le singleton tel que $\{w_{j}, j\in J\}=\{u'\}$, 
 \item $(\alpha_{i},\rho_{i})$ égal à $(4\de,13M+4N+8\de)$ pour $i\in J$ et à 
 $(4P,M)$ pour $i\in I\setminus J$,
%
%
%\item  
% $\{w_{1},...,w_{r}\}$ égal à la réunion disjointe de $\{u'\}$ et de $$\bigcup _{ i\in \{0,\dots ,m\}}
%S_{i}\cup  \bigcup _{i\in \{0,\dots,m\}, j\in \{1,\dots ,l_{i}\}} \mathcal Y_{i}^{j},$$ \item $I=\{i\}$ où $i$ est tel que $w_{i}=u'$,  
%\item $(\alpha_{i},\rho_{i})=(4\de,13M+4N+8\de)$ et $(\alpha_{j},\rho_{j})$ égal à 
%$(4P,M)$ pour $j\neq i$, 
\end{itemize}
on obtient 
   que les distances entre 
les points de l'ensemble (\ref{reunion-Mboules-14dec}) et ceux de $B(u',13M+4N+8\de)$ sont déterminées par les distances entre $C=C(\de,K,N,Q,P,M)$ points de (\ref{reunion-Mboules-14dec}) (déterminés par la donnée de $Z$) et les points  de $B(u',13M+4N+8\de)$. 
 De plus ces distances sont déterminées à $C'=C(\de,N,K,Q,P,M)$ près par la donnée de $Z\in \overline Y_{x,x'}^{p,k,m,(l_{0},...,l_{m})}$. 

Cela termine l'étude du  troisième cas et achève donc la démonstration du lemme~\ref{est-card-fibres-star}. 
\cqfd

On définit la norme pré-hilbertienne $\|.\|_{\H_{x,x',\star,s}(\Delta_{p})}$ sur $\C^{(\Delta_{p})}$ de la même fa\c con que $\|.\|_{\H_{x,s}(\Delta_{p})}$, mais en rempla\c cant 
$\overline Y_{x}^{p,k,m,(l_{0},...,l_{m})}$ par $\overline Y_{x,x',\star}^{p,k,m,(l_{0},...,l_{m})}$. Plus précisément, 
pour $Z\in \overline Y_{x,x',\star}^{p,k,m,(l_{0},...,l_{m})}$ on note $
\xi_{Z}$ la forme linéaire sur $\C^{(\Delta_{p})}$ définie par 
$$\xi_{Z}(f)=\sum _{(a_{1},\dots,a_{p},S_{0},...,S_{m},(\mathcal Y_{i}^{j})_{i\in \{0,\dots,m\}, j\in \{1,\dots ,l_{i}\}}) \in 
(\pi_{x,x',\star}^{p,k,m,(l_{0},...,l_{m})})^{-1}(Z)} f(a_1,...,a_p).$$
Puis pour  $f\in \C^{(\Delta_{p})}$ on pose 
\begin{gather} \nonumber \|f\|_{\H_{x,x',\star,s}(\Delta_{p})}^{2}
=\sum _{k,m,l_{0},\dots ,l_{m}\in \N}  
%$$\|\sum_{\{a_1,a_2,...,a_p\}\in \Delta_p}
%f(a_1,...,a_p)e_{a_1}\wedge ...\wedge
%e_{a_p}\|_{\H_{x,s}(\Delta_{p})}^{2}
%=\sum _{k,m,l_{0},\dots ,l_{m}\in \N} 
%$$
B^{-(m+\sum_{i=0}^{m}l_{i})}
\sum_{Z\in \overline Y_{x,x',\star}^{p,k,m,(l_{0},...,l_{m})}} \\ \label{formule-norme-x'-star} e^{2s(r_{0}(Z)-k)}\Big(\prod_{i=0}^{m}s_{i}(Z)^{-l_{i}} \Big)
\sharp \big((\pi_{x,x',\star}^{p,k,m,(l_{0},...,l_{m})})^{-1}(Z)\big)^{-\alpha }
\big|\xi_{Z}(f)\big|^{2}.\end{gather}

\begin{lem}\label{est-x-x-x'-star}
En notant $C$ la constante qui apparaît dans le b) du lemme~\ref{est-card-fibres-star}, on a pour tout $f\in \C^{(\Delta_{p})}$, 
$$\|f\|_{\H_{x,s}(\Delta_{p})}^{2}\leq C(2d(x,x')+1)\|f\|_{\H_{x,x',\star,s}(\Delta_{p})}^{2}.$$
\end{lem}
\noindent{\bf Démonstration.} Cela résulte immédiatement du lemme~\ref{est-card-fibres-star}, de l'inégalité de Cauchy-Schwarz et du fait que 
pour $Z\in \overline Y_{x}^{p,k,m,(l_{0},...,l_{m})}$ et  $\tilde Z\in \overline Y_{x,x',\star}^{p,k,m,(l_{0},...,l_{m})}$ antécédent de $Z$, c'est-à-dire que 
$$(\pi_{x,x',\star}^{p,k,m,(l_{0},...,l_{m})})^{-1}(\tilde Z)\subset 
(\pi_{x}^{p,k,m,(l_{0},...,l_{m})})^{-1}( Z)$$ on a évidemment 
$\sharp \big( (\pi_{x,x',\star}^{p,k,m,(l_{0},...,l_{m})})^{-1}(\tilde Z) \big)\leq \sharp \big( (\pi_{x}^{p,k,m,(l_{0},...,l_{m})})^{-1}( Z) \big)$. \cqfd

\noindent{\bf Démonstration de la proposition~\ref{equiv-normeYZ} en admettant le lemme~\ref{deuxieme-ineg-equiv-norme}. } Grâce au lemme~\ref{est-x-x-x'-star}, pour montrer la proposition~\ref{equiv-normeYZ} on est ramené à montrer le lemme suivant.  \cqfd

\begin{lem}\label{deuxieme-ineg-equiv-norme}
Il existe une constante $C=C(\de,K,N,Q,P,M,s,B)$, tel que pour tout $f\in \C^{(\Delta_{p})}$, $$\|f\|_{\H_{x,x',\star,s}(\Delta_{p})}^{2}\leq Ce^{3sd(x,x')}\|f\|_{\H_{x',s}(\Delta_{p})}^{2}.$$
\end{lem}

La preuve du lemme~\ref{deuxieme-ineg-equiv-norme} occupe toute la suite de ce sous-paragraphe. 

Afin de mieux comprendre le lemme suivant, on peut se référer au dessin pour les arbres, en remarquant que $\max(k,\frac{d(x,x')+r_{0}(Z)-r'_{0}(Z)}{2})$ vaut $k$ dans les deux premiers dessins et $d(x,t)$ dans le troisième. 

\begin{lem}\label{oubli-x-xx'-x'}
Soit $$(a_{1},\dots,a_{p},S_{0},...,S_{m},(\mathcal Y_{i}^{j})_{i\in \{0,\dots,m\}, j\in \{1,\dots ,l_{i}\}})\in Y_{x}^{p,k,m,(l_{0},...,l_{m})}.$$

\noindent a)  Il existe un unique $m'\in \{0,...,m\}$ tel que \begin{gather*}d_{\max}(x,S_{i})\leq \max(k,\frac{d(x,x')+r_{0}(Z)-r'_{0}(Z)}{2})+M\text{\ \ pour\ \ } i>m' \\ \text{et\ \ \ } d_{\max}(x,S_{i})>\max(k,\frac{d(x,x')+r_{0}(Z)-r'_{0}(Z)}{2})+M\text{\ \ pour\ \ } 1\leq i\leq m'.\end{gather*}

\noindent b) Pour $i>m'$ et $ j\in \{1,...,l_{i}\}$ on a $$d_{\max}(x,\mathcal Y_{i}^{j})\leq \max(k,\frac{d(x,x')+r_{0}(Z)-r'_{0}(Z)}{2})+M+2P+\de.$$

\noindent c) Soit $J\subset \{1,...,l_{m'}\}$ l'ensemble des $j$ tels que $$d_{\max}(x,\mathcal Y_{m'}^{j})> \max(k,\frac{d(x,x')+r_{0}(Z)-r'_{0}(Z)}{2})+M+2P
+\de.$$ On écrit $J=\{j_{1},...,j_{l'_{m'}}\}$ avec $l'_{m'}=\sharp J\in \{0,..., l_{m'}\}$
et $j_{1}< ... < j_{l'_{m'}}$ et on note $k'=k'(Z)$
où $Z\in \overline Y_{x,x'}^{p,k,m,(l_{0},...,l_{m})}$ est tel que $$(a_{1},\dots,a_{p},S_{0},...,S_{m},(\mathcal Y_{i}^{j})_{i\in \{0,\dots,m\}, j\in \{1,\dots ,l_{i}\}})\in (\pi_{x,x'}^{p,k,m,(l_{0},...,l_{m})})^{-1}(Z).$$
 Alors 
$$(a_{1},\dots,a_{p},S_{0},...,S_{m'},(\mathcal Y_{i}^{j})_{i\in \{0,\dots,m'-1\}, j\in \{1,\dots ,l_{i}\}}, (\mathcal Y_{m'}^{j_{\lambda}})_{ \lambda\in \{1,\dots ,l'_{m'}\}})$$ appartient  à 
$Y_{x'}^{p,k',m',(l_{0},...,l_{m'-1},l'_{m'})}$. 

\noindent d) Pour tout $z$ dans \begin{gather*} B(x,k+2M)\cup 
\bigcup _{ i\in \{m'+1,\dots ,m\}}
B(S_{i}, M)\\  
\cup  \bigcup _{i\in \{m'+1,\dots,m\}, j\in \{1,\dots ,l_{i}\}} B(\mathcal Y_{i}^{j},M) \cup \bigcup _{j\not\in J} B(\mathcal Y_{m'}^{j}, M)\end{gather*} et pour tout $z'$ dans $$ \bigcup _{ i\in \{0,\dots ,m'\}}
B(S_{i}, M)
\cup \bigcup _{i\in \{0,\dots,m'-1\}, j\in \{1,...,l_{i}\}} B(\mathcal Y_{i}^{j},M)\cup \bigcup _{j\in J} B(\mathcal Y_{m'}^{j}, M),$$ 
$\geod(z,z')$ rencontre $B(x',k'+2M-2\de)$. 

\noindent e) Pour tout 
$$(\tilde a_{1},\dots,\tilde a_{p},\tilde S_{0},...,\tilde S_{m'},(\tilde {\mathcal Y}_{i}^{j})_{i\in \{0,\dots,m'-1\}, j\in \{1,\dots ,l_{i}\}}, (\tilde {\mathcal Y}_{m'}^{j_{\lambda}})_{ \lambda\in \{1,\dots ,l'_{m'}\}})$$  
appartenant à la même classe que 
$$(a_{1},\dots,a_{p},S_{0},...,S_{m'},(\mathcal Y_{i}^{j})_{i\in \{0,\dots,m'-1\}, j\in \{1,\dots ,l_{i}\}}, (\mathcal Y_{m'}^{j_{\lambda}})_{ \lambda\in \{1,\dots ,l'_{m'}\}})$$  
dans $\overline Y_{x'}^{p,k',m',(l_{0},...,l_{m'-1},l'_{m'})}$, pour tout $z$ dans \begin{gather*} B(x,k+2M)\cup 
\bigcup _{ i\in \{m'+1,\dots ,m\}}
B(S_{i}, M)\\  
\cup  \bigcup _{i\in \{m'+1,\dots,m\}, j\in \{1,\dots ,l_{i}\}} B(\mathcal Y_{i}^{j},M) \cup \bigcup _{j\not\in J} B(\mathcal Y_{m'}^{j}, M)\end{gather*} et pour tout $\tilde z'$ dans $$ \bigcup _{ i\in \{0,\dots ,m'\}}
B(\tilde S_{i}, M)
\cup \bigcup _{i\in \{0,\dots,m'-1\}, j\in \{1,...,l_{i}\}} B(\tilde {\mathcal Y}_{i}^{j},M)\cup \bigcup _{\lambda\in \{1,...,l'_{m'}\}} B(\tilde {\mathcal Y}_{m'}^{j_{\lambda}}, M),$$ 
$\geod(z,\tilde z')$ rencontre $B(x',k'+2M)$. 
\end{lem}
\noindent{\bf Remarque.} L'intérêt de e) est que la classe de 
$$(a_{1},\dots,a_{p},S_{0},...,S_{m},(\mathcal Y_{i}^{j})_{i\in \{0,\dots,m\}, j\in \{1,\dots ,l_{i}\}})$$  dans $\overline Y_{x,x',\star}^{p,k,m,(l_{0},...,l_{m})}$
est déterminée par celle de 
$$(a_{1},\dots,a_{p},S_{0},...,S_{m'},(\mathcal Y_{i}^{j})_{i\in \{0,\dots,m'-1\}, j\in \{1,\dots ,l_{i}\}}, (\mathcal Y_{m'}^{j_{\lambda}})_{ \lambda\in \{1,\dots ,l'_{m'}\}})$$  
dans $\overline Y_{x'}^{p,k',m',(l_{0},...,l_{m'-1},l'_{m'})}$ et par la connaissance des parties ``éliminées'' $$(S_{i})_{ i\in \{m'+1,\dots,m\}}, (\mathcal Y_{i}^{j})_{ i\in \{m'+1,\dots,m\}, j\in \{1,\dots ,l_{i}\}}, (\mathcal Y_{m'}^{j})_{ j\not\in J},$$ qui se trouvent en gros entre les boules $B(x,k+2M)$ et $B(x',k'+2M)$.

\noindent{\bf Démonstration.}
Le a) et le b) résultent des lemmes~\ref{m-fini} et~\ref{dmaxY-dmaxS}. 
Comme la preuve de c) et d) est longue on commence par montrer que d) implique e). Cela résulte du sous-lemme suivant, appliqué à $r=k'+2M$, et $x'$ au lieu de $x$. 

\begin{souslem}\label{souslem-rayon-2de}
Soit $r\geq 2\de$ et $z,z',\tilde z'\in X$ tels que 
\begin{itemize}
\item $\geod(z,z')$ rencontre $B(x,r-2\de)$, 
\item pour tout $y\in B(x,r)$, $d(y,z')=d(y,\tilde z')$. 
\end{itemize}
Alors $\geod(z,\tilde z')$ rencontre $B(x,r)$. 
\end{souslem}

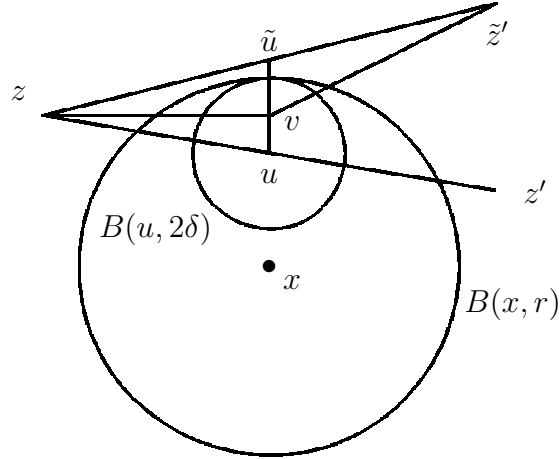
\ifx\JPicScale\undefined\def\JPicScale{1}\fi
\unitlength \JPicScale mm


\noindent{\bf Démonstration. }  
Soit $u\in \geod(z,z')\cap B(x,r-2\de)$. 
On veut montrer $u\in \geod(z,\tilde z')$. Supposons par l'absurde que cela ne soit pas vrai. Soit  $\tilde u \in B(z,d(z,u))$ à distance minimale de 
$\tilde z'$. On a alors $d(\tilde u,\tilde z')< d(u, \tilde z')=d(u,z')\leq d(\tilde u,z')$  donc $\tilde u \not\in B(x,r)$ et en particulier $d(u,\tilde u)\geq 2\de+1$. Soit $v\in \geod(u,\tilde u)$ vérifiant $d(u,v)=\de+1$. 
Par $(H_{\de}^{0}(z,u,v,\tilde u)) $ on a $d(z,v)\leq d(z,u)$. Par $
(H_{\de}^{0}(\tilde z',u,v,\tilde u)) $ on a
$d(\tilde z',v)< d(\tilde z',u)$. Comme $v\in B(x,r)$ on a $d(z',v)=d(\tilde z',v)$. On en déduit que $d(z,v)+d(z',v)< d(z,u)+d(z',u)=d(z,z')$ ce qui est impossible. \cqfd

\noindent{\bf Suite de la démonstration du lemme~\ref{oubli-x-xx'-x'}.} Il reste à montrer  c) et d). Pour le faire  on distingue trois cas comme dans le dessin pour les arbres et dans la démonstration du lemme~\ref{est-card-fibres-star}.

\noindent {\bf Premier cas.} 
On suppose $r_{0}(Z)\leq k$. 

Alors $k'=r'_{0}(Z)+\frac{M}{2}$ par (\ref{form-r'-cas1-4j}). 
%(\ref{def-k'Z-24oct09}). 
D'après le lemme~\ref{m-fini}  
et la remarque qui suit la définition~\ref{defi-Y},  on a $m=0$ et $l_{0}=0$, donc 
 c) est évident. Comme $S_{0}\subset B(x',r_{0}'(Z)+N)$  on a $B(S_{0}, M)\subset B(x',k'+2M-2\de)$ car $M\geq N+2\de$ et 
  d) en résulte, puisque pour $z'$ comme dans d) on a $z'\in B(x',k'+2M-2\de)$.

\noindent {\bf Deuxième  cas.} 
On suppose $r_{0}(Z)> k$ et $\frac{r_{0}(Z)-r_{0}'(Z)+d(x,x')}{2}\leq k$. 

Par (\ref{form-r'-cas2-4j})
%(\ref{def-k'Z-24oct09}) 
on a alors 
$$k'=r_{0}'(Z)-r_{0}(Z)+k+\frac{M}{2}.$$ 
Soit $u\in B(x,k)$ à distance minimale de $S_{0}$ et 
$b\in S_{0}$ à distance minimale de $u$. 
Pour la suite on note le fait suivant. Comme $d(x,b)-d(x',b)+d(x,x')\leq 2k+N$ et $u\in \geod(x,b)$, $(H_{\de}^{0}(x',x,u,b))$ implique $$d(x',u)\leq \max(d(x,x')-k,d(x',b)-d(x,b)+k)+\de\leq d(x',b)-d(b,u)+N+\de$$
et \begin{gather}\label{uNdetg-12dec}u\in (N+\de)\tg(x',b).\end{gather} 
On a $k'\geq d(x',b)-d(x,b)+k+\frac{M}{2}-N=d(x',b)-d(u,b)+\frac{M}{2}-N$ d'où, grâce à (\ref{uNdetg-12dec}), 
\begin{gather}\label{ineg-k'-deuxiemecas-12dec}k'\geq d(x',u)+\frac{M}{2}-2N-\de\end{gather}

Comme  $\max(k,\frac{d(x,x')+r_{0}(Z)-r'_{0}(Z)}{2})=k$,  
pour passer de \begin{gather*}(a_{1},\dots,a_{p},S_{0},...,S_{m},(\mathcal Y_{i}^{j})_{i\in \{0,\dots,m\}, j\in \{1,\dots ,l_{i}\}})\\ \text{à \ \ \ }(a_{1},\dots,a_{p},S_{0},...,S_{m'},(\mathcal Y_{i}^{j})_{i\in \{0,\dots,m'-1\}, j\in \{1,\dots ,l_{i}\}}, (\mathcal Y_{m'}^{j_{\lambda}})_{ \lambda\in \{1,\dots ,l'_{m'}\}})\end{gather*} on enlève les $S_{i}$ et les $\mathcal Y_{i}^{j}$ pour $i\geq1$ vérifiant  $d_{\max}(x,S_{i})\leq k+M$ et les $\mathcal Y_{m'}^{j}$ vérifiant $d_{\max}(x,\mathcal Y_{m'}^{j})\leq k+M+2P+\de$. 
Pour vérifier c) on va montrer que les conditions i), ii), iii) et iv)  de la définition~\ref{defi-Y} sont satisfaites par $$(a_{1},\dots,a_{p},S_{0},...,S_{m'},(\mathcal Y_{i}^{j})_{i\in \{0,\dots,m'-1\}, j\in \{1,\dots ,l_{i}\}}, (\mathcal Y_{m'}^{j_{\lambda}})_{ \lambda\in \{1,\dots ,l'_{m'}\}}).$$  On vérifie d'abord i). Soit $i\in \{0,...,m'-1\}$, $y\in S_{i+1}$, $a\in S_{i}$ et $ \tilde x\in B(x,k)$    tels que 
\begin{itemize}
\item $y\in 4\de\tg(\tilde x,a)$ et $d(y,a)\in ]N-2\de,QN]$, 
\item ou $y\in F\tg(\tilde x,a)$ et $d(y,a)\geq \frac{Q}{F}$. 
\end{itemize}
On note $\alpha=4\de$ ou $F$ suivant le cas, de sorte que $y\in \alpha\tg(\tilde x,a)$. 
Comme $d(x,S_{i+1})\geq d_{\max}(x,S_{i+1})-N$ et $i+1\in \{1,...,m'\}$, on a 
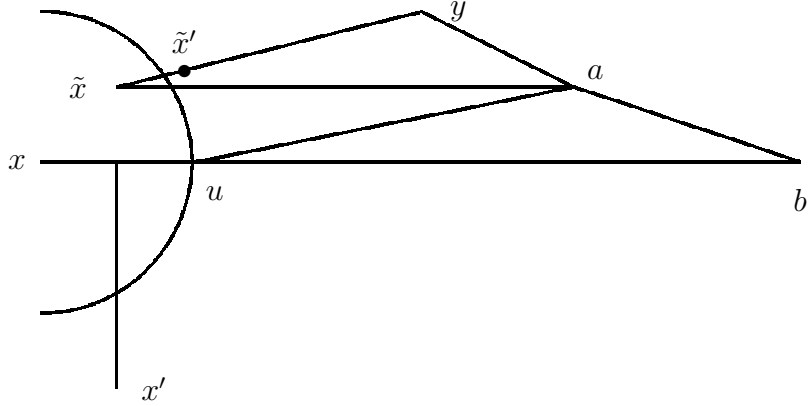
\begin{gather}\label{cond-dxy-p122-24oct09}d(x,y)> k+M-N.\end{gather}

\ifx\JPicScale\undefined\def\JPicScale{1}\fi
\unitlength \JPicScale mm
\begin{picture}(90,60)(0,25)
\linethickness{0.3mm}
\multiput(10,40)(0.5,0.01){1}{\line(1,0){0.5}}
\multiput(10.5,40.01)(0.5,0.02){1}{\line(1,0){0.5}}
\multiput(11,40.02)(0.5,0.03){1}{\line(1,0){0.5}}
\multiput(11.49,40.06)(0.5,0.04){1}{\line(1,0){0.5}}
\multiput(11.99,40.1)(0.5,0.06){1}{\line(1,0){0.5}}
\multiput(12.49,40.16)(0.49,0.07){1}{\line(1,0){0.49}}
\multiput(12.98,40.22)(0.49,0.08){1}{\line(1,0){0.49}}
\multiput(13.47,40.3)(0.49,0.09){1}{\line(1,0){0.49}}
\multiput(13.96,40.4)(0.49,0.1){1}{\line(1,0){0.49}}
\multiput(14.45,40.5)(0.48,0.12){1}{\line(1,0){0.48}}
\multiput(14.94,40.62)(0.48,0.13){1}{\line(1,0){0.48}}
\multiput(15.42,40.75)(0.48,0.14){1}{\line(1,0){0.48}}
\multiput(15.9,40.89)(0.47,0.15){1}{\line(1,0){0.47}}
\multiput(16.37,41.04)(0.47,0.16){1}{\line(1,0){0.47}}
\multiput(16.84,41.21)(0.47,0.18){1}{\line(1,0){0.47}}
\multiput(17.31,41.38)(0.23,0.09){2}{\line(1,0){0.23}}
\multiput(17.77,41.57)(0.23,0.1){2}{\line(1,0){0.23}}
\multiput(18.23,41.77)(0.23,0.11){2}{\line(1,0){0.23}}
\multiput(18.68,41.98)(0.22,0.11){2}{\line(1,0){0.22}}
\multiput(19.12,42.2)(0.22,0.12){2}{\line(1,0){0.22}}
\multiput(19.57,42.44)(0.22,0.12){2}{\line(1,0){0.22}}
\multiput(20,42.68)(0.21,0.13){2}{\line(1,0){0.21}}
\multiput(20.43,42.93)(0.21,0.13){2}{\line(1,0){0.21}}
\multiput(20.85,43.2)(0.21,0.14){2}{\line(1,0){0.21}}
\multiput(21.27,43.48)(0.2,0.14){2}{\line(1,0){0.2}}
\multiput(21.67,43.76)(0.2,0.15){2}{\line(1,0){0.2}}
\multiput(22.08,44.06)(0.13,0.1){3}{\line(1,0){0.13}}
\multiput(22.47,44.36)(0.13,0.11){3}{\line(1,0){0.13}}
\multiput(22.86,44.68)(0.13,0.11){3}{\line(1,0){0.13}}
\multiput(23.23,45)(0.12,0.11){3}{\line(1,0){0.12}}
\multiput(23.6,45.34)(0.12,0.11){3}{\line(1,0){0.12}}
\multiput(23.96,45.68)(0.12,0.12){3}{\line(0,1){0.12}}
\multiput(24.32,46.04)(0.11,0.12){3}{\line(0,1){0.12}}
\multiput(24.66,46.4)(0.11,0.12){3}{\line(0,1){0.12}}
\multiput(25,46.77)(0.11,0.13){3}{\line(0,1){0.13}}
\multiput(25.32,47.14)(0.11,0.13){3}{\line(0,1){0.13}}
\multiput(25.64,47.53)(0.1,0.13){3}{\line(0,1){0.13}}
\multiput(25.94,47.92)(0.15,0.2){2}{\line(0,1){0.2}}
\multiput(26.24,48.33)(0.14,0.2){2}{\line(0,1){0.2}}
\multiput(26.52,48.73)(0.14,0.21){2}{\line(0,1){0.21}}
\multiput(26.8,49.15)(0.13,0.21){2}{\line(0,1){0.21}}
\multiput(27.07,49.57)(0.13,0.21){2}{\line(0,1){0.21}}
\multiput(27.32,50)(0.12,0.22){2}{\line(0,1){0.22}}
\multiput(27.56,50.43)(0.12,0.22){2}{\line(0,1){0.22}}
\multiput(27.8,50.88)(0.11,0.22){2}{\line(0,1){0.22}}
\multiput(28.02,51.32)(0.11,0.23){2}{\line(0,1){0.23}}
\multiput(28.23,51.77)(0.1,0.23){2}{\line(0,1){0.23}}
\multiput(28.43,52.23)(0.09,0.23){2}{\line(0,1){0.23}}
\multiput(28.62,52.69)(0.18,0.47){1}{\line(0,1){0.47}}
\multiput(28.79,53.16)(0.16,0.47){1}{\line(0,1){0.47}}
\multiput(28.96,53.63)(0.15,0.47){1}{\line(0,1){0.47}}
\multiput(29.11,54.1)(0.14,0.48){1}{\line(0,1){0.48}}
\multiput(29.25,54.58)(0.13,0.48){1}{\line(0,1){0.48}}
\multiput(29.38,55.06)(0.12,0.48){1}{\line(0,1){0.48}}
\multiput(29.5,55.55)(0.1,0.49){1}{\line(0,1){0.49}}
\multiput(29.6,56.04)(0.09,0.49){1}{\line(0,1){0.49}}
\multiput(29.7,56.53)(0.08,0.49){1}{\line(0,1){0.49}}
\multiput(29.78,57.02)(0.07,0.49){1}{\line(0,1){0.49}}
\multiput(29.84,57.51)(0.06,0.5){1}{\line(0,1){0.5}}
\multiput(29.9,58.01)(0.04,0.5){1}{\line(0,1){0.5}}
\multiput(29.94,58.51)(0.03,0.5){1}{\line(0,1){0.5}}
\multiput(29.98,59)(0.02,0.5){1}{\line(0,1){0.5}}
\multiput(29.99,59.5)(0.01,0.5){1}{\line(0,1){0.5}}
\multiput(29.99,60.5)(0.01,-0.5){1}{\line(0,-1){0.5}}
\multiput(29.98,61)(0.02,-0.5){1}{\line(0,-1){0.5}}
\multiput(29.94,61.49)(0.03,-0.5){1}{\line(0,-1){0.5}}
\multiput(29.9,61.99)(0.04,-0.5){1}{\line(0,-1){0.5}}
\multiput(29.84,62.49)(0.06,-0.5){1}{\line(0,-1){0.5}}
\multiput(29.78,62.98)(0.07,-0.49){1}{\line(0,-1){0.49}}
\multiput(29.7,63.47)(0.08,-0.49){1}{\line(0,-1){0.49}}
\multiput(29.6,63.96)(0.09,-0.49){1}{\line(0,-1){0.49}}
\multiput(29.5,64.45)(0.1,-0.49){1}{\line(0,-1){0.49}}
\multiput(29.38,64.94)(0.12,-0.48){1}{\line(0,-1){0.48}}
\multiput(29.25,65.42)(0.13,-0.48){1}{\line(0,-1){0.48}}
\multiput(29.11,65.9)(0.14,-0.48){1}{\line(0,-1){0.48}}
\multiput(28.96,66.37)(0.15,-0.47){1}{\line(0,-1){0.47}}
\multiput(28.79,66.84)(0.16,-0.47){1}{\line(0,-1){0.47}}
\multiput(28.62,67.31)(0.18,-0.47){1}{\line(0,-1){0.47}}
\multiput(28.43,67.77)(0.09,-0.23){2}{\line(0,-1){0.23}}
\multiput(28.23,68.23)(0.1,-0.23){2}{\line(0,-1){0.23}}
\multiput(28.02,68.68)(0.11,-0.23){2}{\line(0,-1){0.23}}
\multiput(27.8,69.12)(0.11,-0.22){2}{\line(0,-1){0.22}}
\multiput(27.56,69.57)(0.12,-0.22){2}{\line(0,-1){0.22}}
\multiput(27.32,70)(0.12,-0.22){2}{\line(0,-1){0.22}}
\multiput(27.07,70.43)(0.13,-0.21){2}{\line(0,-1){0.21}}
\multiput(26.8,70.85)(0.13,-0.21){2}{\line(0,-1){0.21}}
\multiput(26.52,71.27)(0.14,-0.21){2}{\line(0,-1){0.21}}
\multiput(26.24,71.67)(0.14,-0.2){2}{\line(0,-1){0.2}}
\multiput(25.94,72.08)(0.15,-0.2){2}{\line(0,-1){0.2}}
\multiput(25.64,72.47)(0.1,-0.13){3}{\line(0,-1){0.13}}
\multiput(25.32,72.86)(0.11,-0.13){3}{\line(0,-1){0.13}}
\multiput(25,73.23)(0.11,-0.13){3}{\line(0,-1){0.13}}
\multiput(24.66,73.6)(0.11,-0.12){3}{\line(0,-1){0.12}}
\multiput(24.32,73.96)(0.11,-0.12){3}{\line(0,-1){0.12}}
\multiput(23.96,74.32)(0.12,-0.12){3}{\line(0,-1){0.12}}
\multiput(23.6,74.66)(0.12,-0.11){3}{\line(1,0){0.12}}
\multiput(23.23,75)(0.12,-0.11){3}{\line(1,0){0.12}}
\multiput(22.86,75.32)(0.13,-0.11){3}{\line(1,0){0.13}}
\multiput(22.47,75.64)(0.13,-0.11){3}{\line(1,0){0.13}}
\multiput(22.08,75.94)(0.13,-0.1){3}{\line(1,0){0.13}}
\multiput(21.67,76.24)(0.2,-0.15){2}{\line(1,0){0.2}}
\multiput(21.27,76.52)(0.2,-0.14){2}{\line(1,0){0.2}}
\multiput(20.85,76.8)(0.21,-0.14){2}{\line(1,0){0.21}}
\multiput(20.43,77.07)(0.21,-0.13){2}{\line(1,0){0.21}}
\multiput(20,77.32)(0.21,-0.13){2}{\line(1,0){0.21}}
\multiput(19.57,77.56)(0.22,-0.12){2}{\line(1,0){0.22}}
\multiput(19.12,77.8)(0.22,-0.12){2}{\line(1,0){0.22}}
\multiput(18.68,78.02)(0.22,-0.11){2}{\line(1,0){0.22}}
\multiput(18.23,78.23)(0.23,-0.11){2}{\line(1,0){0.23}}
\multiput(17.77,78.43)(0.23,-0.1){2}{\line(1,0){0.23}}
\multiput(17.31,78.62)(0.23,-0.09){2}{\line(1,0){0.23}}
\multiput(16.84,78.79)(0.47,-0.18){1}{\line(1,0){0.47}}
\multiput(16.37,78.96)(0.47,-0.16){1}{\line(1,0){0.47}}
\multiput(15.9,79.11)(0.47,-0.15){1}{\line(1,0){0.47}}
\multiput(15.42,79.25)(0.48,-0.14){1}{\line(1,0){0.48}}
\multiput(14.94,79.38)(0.48,-0.13){1}{\line(1,0){0.48}}
\multiput(14.45,79.5)(0.48,-0.12){1}{\line(1,0){0.48}}
\multiput(13.96,79.6)(0.49,-0.1){1}{\line(1,0){0.49}}
\multiput(13.47,79.7)(0.49,-0.09){1}{\line(1,0){0.49}}
\multiput(12.98,79.78)(0.49,-0.08){1}{\line(1,0){0.49}}
\multiput(12.49,79.84)(0.49,-0.07){1}{\line(1,0){0.49}}
\multiput(11.99,79.9)(0.5,-0.06){1}{\line(1,0){0.5}}
\multiput(11.49,79.94)(0.5,-0.04){1}{\line(1,0){0.5}}
\multiput(11,79.98)(0.5,-0.03){1}{\line(1,0){0.5}}
\multiput(10.5,79.99)(0.5,-0.02){1}{\line(1,0){0.5}}
\multiput(10,80)(0.5,-0.01){1}{\line(1,0){0.5}}

\linethickness{0.3mm}
\put(10,60){\line(1,0){100}}
\linethickness{0.3mm}
\multiput(80,70)(0.36,-0.12){83}{\line(1,0){0.36}}
\linethickness{0.3mm}
\multiput(30,60)(0.6,0.12){83}{\line(1,0){0.6}}
\linethickness{0.3mm}
\put(20,70){\line(1,0){60}}
\linethickness{0.3mm}
\multiput(20,70)(0.48,0.12){83}{\line(1,0){0.48}}
\linethickness{0.3mm}
\multiput(60,80)(0.24,-0.12){83}{\line(1,0){0.24}}
\linethickness{0.3mm}
\put(20,30){\line(0,1){30}}

\put(15,70){\makebox(0,0)[cc]{$\tilde x$}}

\put(29,76){\makebox(0,0)[cc]{$\tilde x'$}}

\put(29,72){\makebox(0,0)[cc]{$\bullet$}}

\put(07,60){\makebox(0,0)[cc]{$x$}}

\put(25,30){\makebox(0,0)[cc]{$x'$}}

\put(33,56){\makebox(0,0)[cc]{$u$}}

\put(65,80){\makebox(0,0)[cc]{$y$}}

\put(83,72){\makebox(0,0)[cc]{$a$}}

\put(110,55){\makebox(0,0)[cc]{$b$}}

\end{picture}
 
  \noindent Pour vérifier i) il suffit de montrer qu'il existe $\tilde x'\in B(x', k')$ tel que $y\in \alpha\tg(\tilde x',a)$. Cela résulte du sous-lemme suivant appliqué à $\alpha=4\de$ ou $\alpha=F$ et à $\beta=M-N$. On a $F\geq 4\de$ et on suppose  $F+\frac{3\de+F}{2}<  M-N$, ce qui est permis par $(H_{M})$. Les hypothèses du sous-lemme suivant sont donc satisfaites et cela termine la preuve de i). 

\begin{souslem}\label{souslem1}
Soit $\alpha,\beta \in \N$ vérifiant $F+\frac{3\de+\alpha}{2}< \beta$. 
Soit $a\in 2F\tg(u,b)$, $\tilde x\in B(x,k)$ et  $y \in \alpha\tg(\tilde x,a)$ tel que $d(x,y)\geq k+\beta$. Alors il existe $\tilde x'\in B(x',k')$ tel que $y\in \alpha\tg(\tilde x',a)$.
\end{souslem}

\noindent{\bf Démonstration du sous-lemme~\ref{souslem1}.}
Comme $a\in 2F\tg(b,u)$ et $u\in \de\tg(\tilde x,b)$ par le lemme~\ref{Bxk-x-t-z-y}, on a $a\in (2F+\de)\tg(\tilde x,b)$  par le lemme~\ref{geod-comp-xabc}, 
et comme $y\in \alpha\tg(\tilde x,a)$ on en déduit  $$y\in (2F+\de+\alpha)\tg(\tilde x,b).$$ 
On applique le lemme~\ref{approx-arbres} à $\{x,u,b,\tilde x,y\}$ avec $l=2$ et $x$ comme point base. 
Soit $t$ le point de $\geod(\Psi \tilde x, \Psi b)$ le plus proche de $\Psi y$. 

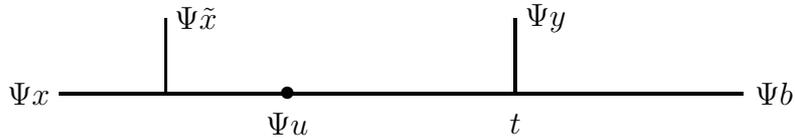
\ifx\JPicScale\undefined\def\JPicScale{1}\fi
\unitlength \JPicScale mm
\begin{picture}(100,30)(20,20)
\linethickness{0.3mm}
\put(40,30){\line(1,0){90}}
\linethickness{0.3mm}
\put(100,30){\line(0,1){10}}
\linethickness{0.3mm}
\put(54,30){\line(0,1){10}}
\put(36,30){\makebox(0,0)[cc]{$\Psi x$}}

\put(58,40){\makebox(0,0)[cc]{$\Psi \tilde x$}}

\put(70,30){\makebox(0,0)[cc]{$\bullet$}}

\put(70,26){\makebox(0,0)[cc]{$\Psi u$}}

\put(100,26){\makebox(0,0)[cc]{$t$}}

\put(104,40){\makebox(0,0)[cc]{$\Psi y$}}

\put(134,30){\makebox(0,0)[cc]{$\Psi b$}}

\end{picture}

\noindent
Grâce au a) du  lemme~\ref{approx-arbres2} on a $\Psi y\in (2F+3\de+\alpha)\tg(\Psi \tilde x,\Psi b)$, d'où 
$d(\Psi y,t)\leq  F+\frac{3\de+\alpha}{2}$. De plus $d(\Psi x,\Psi y)=d(x,y)\geq  
k+\beta$. Comme $F+\frac{3\de+\alpha}{2}< \beta$ on a 
$d(\Psi x,t)> k=d(\Psi x,\Psi u)$ et comme $d(\Psi x,\Psi \tilde x)\leq k$ on 
en déduit $\Psi u\in \geod(\Psi\tilde x, \Psi y)$ et $u\in 4\de\tg(\tilde x, y)$ par le b) du  lemme~\ref{approx-arbres2} . 

Soit $\tilde x'\in \geod(\tilde x,y)$ vérifiant $d(y,\tilde x')=\min(d(y,u), d(y,\tilde x))$. On a alors, par $(H_{\de}^{0}(u,y,\tilde x',\tilde x))$, 
$d(u,\tilde x')\leq 5\de$. Par (\ref{ineg-k'-deuxiemecas-12dec}) 
on a $
k'\geq d(x',u)+\frac{M}{2}-2N-\de$ et on suppose $2N+6\de\leq \frac{M}{2}$, ce qui est permis par $(H_{M})$. Alors $d(x',\tilde x')\leq k'$ et ceci termine la démonstration du sous-lemme~\ref{souslem1}. \cqfd
% Par le lemme~\ref{geod-comp-xabc},  $y\in \alpha\tg(\tilde x',a)$. 
% Comme $u$ appartient à $\geod(x,b)$  et $ (N+\de)\tg(x',b)$ on a $$k'\geq d(x',b)-d(x,b)+k+\frac{M}{2}-N$$ $$=d(x',b)-d(u,b)+\frac{M}{2}-N\geq d(x',u)+\frac{M}{2}-2N-\de$$ d'où $\tilde x'\in B(x',k')$ car $ \frac{M}{2}\geq (2N+\de)+6\de$ grâce à $(H_{M})$. Ceci termine la vérification de i). 

\noindent{\bf Fin de l'étude du deuxième cas.} 
On vérifie maintenant ii). Soit $i\in \{1,...,m'\}$. On a $$S_{i}\subset 2F\tg(b,u)$$ par le a) du lemme~\ref{lemme-S0-...Sm}. On a $d_{\max}(x,S_{i})>k+M$ par hypothèse. 
 Soit $a\in S_{i}$. On a alors $d(x,a)> k+M-N$ d'où $d(u,a)> M-N$.  
   Par (\ref{uNdetg-12dec}) et  le b) du lemme~\ref{geod-comp-xabc}, on a $u\in (2F+N+\de)\tg(x',a)$. D'où \begin{gather*}d(x',a)
 \geq  d(x',u)+d(u,a)-(2F+N+\de)\\ 
 > d(x',u)+M-(2F+2N+\de).\end{gather*} Or \begin{gather*}k'\leq d(x',b)-d(x,b)+k+\frac{M}{2}+N\\ =d(x',b)-d(b,u)+\frac{M}{2}+N\leq d(x',u)+\frac{M}{2}+N.\end{gather*} 
 D'où $$d(x',a)> k'+\frac{M}{2}-(2F+3N+\de)\geq k'+P$$ car on suppose $(2F+3N+\de)+P\leq \frac{M}{2}$, ce qui est permis par  $(H_{M})$. Cela achève la preuve de ii).
 
 La propriété  iii) est immédiate. 
 
 On vérifie maintenant iv).   Soit $j\in J$. On commence par montrer $d(x',
 \mathcal Y_{m'}^{j})\geq k'+3P$. On rappelle que $\mathcal Y_{m'}^{j}\subset 4P\tg(x,b)$ et $d(x,\mathcal Y_{m'}^{j})\geq k+M+P+\de$. Soit $y\in 
 \mathcal Y_{m'}^{j}$. On a alors $$d(b,y)\leq d(x,b)-d(x,y)+4P\leq d(b,x)-k-M-\de+3P$$ et $$d(x',y)\geq d(x',b)-d(b,y)\geq d(x',b)-d(x,b)+k+M+\de-3P$$ $$\geq k'+\frac{M}{2}-N+\de-3P\geq k'+3P.$$ 
car on suppose $\frac{M}{2}-N+\de-3P\geq 3P$, ce qui est permis par  $(H_{M})$. 

 Il reste donc à montrer que pour tout  $y\in 
 \mathcal Y_{m'}^{j}$ il existe $a\in S_{m'}$ et  
 $\tilde x'\in B(x',k')$ tels que  $y\in 2P\tg(\tilde x',a)$. Soit $y\in 
 \mathcal Y_{m'}^{j}$. 
 On va distinguer deux cas. 
 
  On suppose d'abord $m'<m$. Alors il existe $a\in S_{m'}$ et $w\in S_{m'+1}$ tels que $y\in P\tg(w,a)$. Par le a) du lemme~\ref{lemme-S0-...Sm}, $w$ et $a$ appartiennent à  $ 2F\tg(u,b)$. Comme $u
  \in \geod(x,b)$ le b) du  lemme~\ref{geod-comp-xabc} montre $u\in 2F\tg(x,w)$. Donc $d(u,w)\leq d(x,w)-k+2F$.

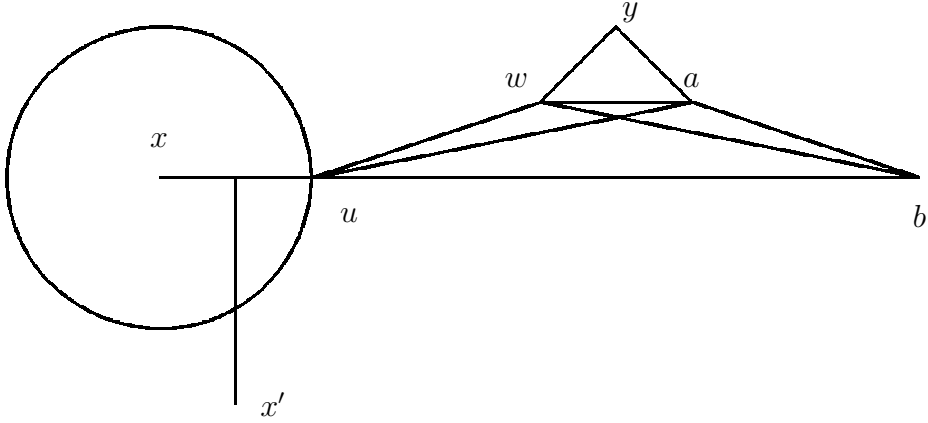
  \ifx\JPicScale\undefined\def\JPicScale{1}\fi
\unitlength \JPicScale mm
\begin{picture}(110,63)(10,15)
\linethickness{0.3mm}
\put(50,49.75){\line(0,1){0.5}}
\multiput(49.99,50.75)(0.01,-0.5){1}{\line(0,-1){0.5}}
\multiput(49.96,51.25)(0.02,-0.5){1}{\line(0,-1){0.5}}
\multiput(49.92,51.74)(0.04,-0.5){1}{\line(0,-1){0.5}}
\multiput(49.87,52.24)(0.05,-0.5){1}{\line(0,-1){0.5}}
\multiput(49.81,52.73)(0.06,-0.49){1}{\line(0,-1){0.49}}
\multiput(49.74,53.23)(0.07,-0.49){1}{\line(0,-1){0.49}}
\multiput(49.65,53.72)(0.09,-0.49){1}{\line(0,-1){0.49}}
\multiput(49.55,54.21)(0.1,-0.49){1}{\line(0,-1){0.49}}
\multiput(49.44,54.69)(0.11,-0.49){1}{\line(0,-1){0.49}}
\multiput(49.32,55.18)(0.12,-0.48){1}{\line(0,-1){0.48}}
\multiput(49.18,55.66)(0.14,-0.48){1}{\line(0,-1){0.48}}
\multiput(49.04,56.13)(0.15,-0.48){1}{\line(0,-1){0.48}}
\multiput(48.88,56.61)(0.16,-0.47){1}{\line(0,-1){0.47}}
\multiput(48.71,57.07)(0.17,-0.47){1}{\line(0,-1){0.47}}
\multiput(48.52,57.54)(0.09,-0.23){2}{\line(0,-1){0.23}}
\multiput(48.33,58)(0.1,-0.23){2}{\line(0,-1){0.23}}
\multiput(48.13,58.45)(0.1,-0.23){2}{\line(0,-1){0.23}}
\multiput(47.91,58.9)(0.11,-0.22){2}{\line(0,-1){0.22}}
\multiput(47.68,59.35)(0.11,-0.22){2}{\line(0,-1){0.22}}
\multiput(47.44,59.78)(0.12,-0.22){2}{\line(0,-1){0.22}}
\multiput(47.19,60.22)(0.12,-0.22){2}{\line(0,-1){0.22}}
\multiput(46.93,60.64)(0.13,-0.21){2}{\line(0,-1){0.21}}
\multiput(46.66,61.06)(0.14,-0.21){2}{\line(0,-1){0.21}}
\multiput(46.38,61.47)(0.14,-0.21){2}{\line(0,-1){0.21}}
\multiput(46.09,61.88)(0.15,-0.2){2}{\line(0,-1){0.2}}
\multiput(45.79,62.27)(0.1,-0.13){3}{\line(0,-1){0.13}}
\multiput(45.48,62.66)(0.1,-0.13){3}{\line(0,-1){0.13}}
\multiput(45.16,63.05)(0.11,-0.13){3}{\line(0,-1){0.13}}
\multiput(44.83,63.42)(0.11,-0.12){3}{\line(0,-1){0.12}}
\multiput(44.49,63.79)(0.11,-0.12){3}{\line(0,-1){0.12}}
\multiput(44.14,64.14)(0.12,-0.12){3}{\line(0,-1){0.12}}
\multiput(43.79,64.49)(0.12,-0.12){3}{\line(1,0){0.12}}
\multiput(43.42,64.83)(0.12,-0.11){3}{\line(1,0){0.12}}
\multiput(43.05,65.16)(0.12,-0.11){3}{\line(1,0){0.12}}
\multiput(42.66,65.48)(0.13,-0.11){3}{\line(1,0){0.13}}
\multiput(42.27,65.79)(0.13,-0.1){3}{\line(1,0){0.13}}
\multiput(41.88,66.09)(0.13,-0.1){3}{\line(1,0){0.13}}
\multiput(41.47,66.38)(0.2,-0.15){2}{\line(1,0){0.2}}
\multiput(41.06,66.66)(0.21,-0.14){2}{\line(1,0){0.21}}
\multiput(40.64,66.93)(0.21,-0.14){2}{\line(1,0){0.21}}
\multiput(40.22,67.19)(0.21,-0.13){2}{\line(1,0){0.21}}
\multiput(39.78,67.44)(0.22,-0.12){2}{\line(1,0){0.22}}
\multiput(39.35,67.68)(0.22,-0.12){2}{\line(1,0){0.22}}
\multiput(38.9,67.91)(0.22,-0.11){2}{\line(1,0){0.22}}
\multiput(38.45,68.13)(0.22,-0.11){2}{\line(1,0){0.22}}
\multiput(38,68.33)(0.23,-0.1){2}{\line(1,0){0.23}}
\multiput(37.54,68.52)(0.23,-0.1){2}{\line(1,0){0.23}}
\multiput(37.07,68.71)(0.23,-0.09){2}{\line(1,0){0.23}}
\multiput(36.61,68.88)(0.47,-0.17){1}{\line(1,0){0.47}}
\multiput(36.13,69.04)(0.47,-0.16){1}{\line(1,0){0.47}}
\multiput(35.66,69.18)(0.48,-0.15){1}{\line(1,0){0.48}}
\multiput(35.18,69.32)(0.48,-0.14){1}{\line(1,0){0.48}}
\multiput(34.69,69.44)(0.48,-0.12){1}{\line(1,0){0.48}}
\multiput(34.21,69.55)(0.49,-0.11){1}{\line(1,0){0.49}}
\multiput(33.72,69.65)(0.49,-0.1){1}{\line(1,0){0.49}}
\multiput(33.23,69.74)(0.49,-0.09){1}{\line(1,0){0.49}}
\multiput(32.73,69.81)(0.49,-0.07){1}{\line(1,0){0.49}}
\multiput(32.24,69.87)(0.49,-0.06){1}{\line(1,0){0.49}}
\multiput(31.74,69.92)(0.5,-0.05){1}{\line(1,0){0.5}}
\multiput(31.25,69.96)(0.5,-0.04){1}{\line(1,0){0.5}}
\multiput(30.75,69.99)(0.5,-0.02){1}{\line(1,0){0.5}}
\multiput(30.25,70)(0.5,-0.01){1}{\line(1,0){0.5}}
\put(29.75,70){\line(1,0){0.5}}
\multiput(29.25,69.99)(0.5,0.01){1}{\line(1,0){0.5}}
\multiput(28.75,69.96)(0.5,0.02){1}{\line(1,0){0.5}}
\multiput(28.26,69.92)(0.5,0.04){1}{\line(1,0){0.5}}
\multiput(27.76,69.87)(0.5,0.05){1}{\line(1,0){0.5}}
\multiput(27.27,69.81)(0.49,0.06){1}{\line(1,0){0.49}}
\multiput(26.77,69.74)(0.49,0.07){1}{\line(1,0){0.49}}
\multiput(26.28,69.65)(0.49,0.09){1}{\line(1,0){0.49}}
\multiput(25.79,69.55)(0.49,0.1){1}{\line(1,0){0.49}}
\multiput(25.31,69.44)(0.49,0.11){1}{\line(1,0){0.49}}
\multiput(24.82,69.32)(0.48,0.12){1}{\line(1,0){0.48}}
\multiput(24.34,69.18)(0.48,0.14){1}{\line(1,0){0.48}}
\multiput(23.87,69.04)(0.48,0.15){1}{\line(1,0){0.48}}
\multiput(23.39,68.88)(0.47,0.16){1}{\line(1,0){0.47}}
\multiput(22.93,68.71)(0.47,0.17){1}{\line(1,0){0.47}}
\multiput(22.46,68.52)(0.23,0.09){2}{\line(1,0){0.23}}
\multiput(22,68.33)(0.23,0.1){2}{\line(1,0){0.23}}
\multiput(21.55,68.13)(0.23,0.1){2}{\line(1,0){0.23}}
\multiput(21.1,67.91)(0.22,0.11){2}{\line(1,0){0.22}}
\multiput(20.65,67.68)(0.22,0.11){2}{\line(1,0){0.22}}
\multiput(20.22,67.44)(0.22,0.12){2}{\line(1,0){0.22}}
\multiput(19.78,67.19)(0.22,0.12){2}{\line(1,0){0.22}}
\multiput(19.36,66.93)(0.21,0.13){2}{\line(1,0){0.21}}
\multiput(18.94,66.66)(0.21,0.14){2}{\line(1,0){0.21}}
\multiput(18.53,66.38)(0.21,0.14){2}{\line(1,0){0.21}}
\multiput(18.12,66.09)(0.2,0.15){2}{\line(1,0){0.2}}
\multiput(17.73,65.79)(0.13,0.1){3}{\line(1,0){0.13}}
\multiput(17.34,65.48)(0.13,0.1){3}{\line(1,0){0.13}}
\multiput(16.95,65.16)(0.13,0.11){3}{\line(1,0){0.13}}
\multiput(16.58,64.83)(0.12,0.11){3}{\line(1,0){0.12}}
\multiput(16.21,64.49)(0.12,0.11){3}{\line(1,0){0.12}}
\multiput(15.86,64.14)(0.12,0.12){3}{\line(1,0){0.12}}
\multiput(15.51,63.79)(0.12,0.12){3}{\line(0,1){0.12}}
\multiput(15.17,63.42)(0.11,0.12){3}{\line(0,1){0.12}}
\multiput(14.84,63.05)(0.11,0.12){3}{\line(0,1){0.12}}
\multiput(14.52,62.66)(0.11,0.13){3}{\line(0,1){0.13}}
\multiput(14.21,62.27)(0.1,0.13){3}{\line(0,1){0.13}}
\multiput(13.91,61.88)(0.1,0.13){3}{\line(0,1){0.13}}
\multiput(13.62,61.47)(0.15,0.2){2}{\line(0,1){0.2}}
\multiput(13.34,61.06)(0.14,0.21){2}{\line(0,1){0.21}}
\multiput(13.07,60.64)(0.14,0.21){2}{\line(0,1){0.21}}
\multiput(12.81,60.22)(0.13,0.21){2}{\line(0,1){0.21}}
\multiput(12.56,59.78)(0.12,0.22){2}{\line(0,1){0.22}}
\multiput(12.32,59.35)(0.12,0.22){2}{\line(0,1){0.22}}
\multiput(12.09,58.9)(0.11,0.22){2}{\line(0,1){0.22}}
\multiput(11.87,58.45)(0.11,0.22){2}{\line(0,1){0.22}}
\multiput(11.67,58)(0.1,0.23){2}{\line(0,1){0.23}}
\multiput(11.48,57.54)(0.1,0.23){2}{\line(0,1){0.23}}
\multiput(11.29,57.07)(0.09,0.23){2}{\line(0,1){0.23}}
\multiput(11.12,56.61)(0.17,0.47){1}{\line(0,1){0.47}}
\multiput(10.96,56.13)(0.16,0.47){1}{\line(0,1){0.47}}
\multiput(10.82,55.66)(0.15,0.48){1}{\line(0,1){0.48}}
\multiput(10.68,55.18)(0.14,0.48){1}{\line(0,1){0.48}}
\multiput(10.56,54.69)(0.12,0.48){1}{\line(0,1){0.48}}
\multiput(10.45,54.21)(0.11,0.49){1}{\line(0,1){0.49}}
\multiput(10.35,53.72)(0.1,0.49){1}{\line(0,1){0.49}}
\multiput(10.26,53.23)(0.09,0.49){1}{\line(0,1){0.49}}
\multiput(10.19,52.73)(0.07,0.49){1}{\line(0,1){0.49}}
\multiput(10.13,52.24)(0.06,0.49){1}{\line(0,1){0.49}}
\multiput(10.08,51.74)(0.05,0.5){1}{\line(0,1){0.5}}
\multiput(10.04,51.25)(0.04,0.5){1}{\line(0,1){0.5}}
\multiput(10.01,50.75)(0.02,0.5){1}{\line(0,1){0.5}}
\multiput(10,50.25)(0.01,0.5){1}{\line(0,1){0.5}}
\put(10,49.75){\line(0,1){0.5}}
\multiput(10,49.75)(0.01,-0.5){1}{\line(0,-1){0.5}}
\multiput(10.01,49.25)(0.02,-0.5){1}{\line(0,-1){0.5}}
\multiput(10.04,48.75)(0.04,-0.5){1}{\line(0,-1){0.5}}
\multiput(10.08,48.26)(0.05,-0.5){1}{\line(0,-1){0.5}}
\multiput(10.13,47.76)(0.06,-0.49){1}{\line(0,-1){0.49}}
\multiput(10.19,47.27)(0.07,-0.49){1}{\line(0,-1){0.49}}
\multiput(10.26,46.77)(0.09,-0.49){1}{\line(0,-1){0.49}}
\multiput(10.35,46.28)(0.1,-0.49){1}{\line(0,-1){0.49}}
\multiput(10.45,45.79)(0.11,-0.49){1}{\line(0,-1){0.49}}
\multiput(10.56,45.31)(0.12,-0.48){1}{\line(0,-1){0.48}}
\multiput(10.68,44.82)(0.14,-0.48){1}{\line(0,-1){0.48}}
\multiput(10.82,44.34)(0.15,-0.48){1}{\line(0,-1){0.48}}
\multiput(10.96,43.87)(0.16,-0.47){1}{\line(0,-1){0.47}}
\multiput(11.12,43.39)(0.17,-0.47){1}{\line(0,-1){0.47}}
\multiput(11.29,42.93)(0.09,-0.23){2}{\line(0,-1){0.23}}
\multiput(11.48,42.46)(0.1,-0.23){2}{\line(0,-1){0.23}}
\multiput(11.67,42)(0.1,-0.23){2}{\line(0,-1){0.23}}
\multiput(11.87,41.55)(0.11,-0.22){2}{\line(0,-1){0.22}}
\multiput(12.09,41.1)(0.11,-0.22){2}{\line(0,-1){0.22}}
\multiput(12.32,40.65)(0.12,-0.22){2}{\line(0,-1){0.22}}
\multiput(12.56,40.22)(0.12,-0.22){2}{\line(0,-1){0.22}}
\multiput(12.81,39.78)(0.13,-0.21){2}{\line(0,-1){0.21}}
\multiput(13.07,39.36)(0.14,-0.21){2}{\line(0,-1){0.21}}
\multiput(13.34,38.94)(0.14,-0.21){2}{\line(0,-1){0.21}}
\multiput(13.62,38.53)(0.15,-0.2){2}{\line(0,-1){0.2}}
\multiput(13.91,38.12)(0.1,-0.13){3}{\line(0,-1){0.13}}
\multiput(14.21,37.73)(0.1,-0.13){3}{\line(0,-1){0.13}}
\multiput(14.52,37.34)(0.11,-0.13){3}{\line(0,-1){0.13}}
\multiput(14.84,36.95)(0.11,-0.12){3}{\line(0,-1){0.12}}
\multiput(15.17,36.58)(0.11,-0.12){3}{\line(0,-1){0.12}}
\multiput(15.51,36.21)(0.12,-0.12){3}{\line(0,-1){0.12}}
\multiput(15.86,35.86)(0.12,-0.12){3}{\line(1,0){0.12}}
\multiput(16.21,35.51)(0.12,-0.11){3}{\line(1,0){0.12}}
\multiput(16.58,35.17)(0.12,-0.11){3}{\line(1,0){0.12}}
\multiput(16.95,34.84)(0.13,-0.11){3}{\line(1,0){0.13}}
\multiput(17.34,34.52)(0.13,-0.1){3}{\line(1,0){0.13}}
\multiput(17.73,34.21)(0.13,-0.1){3}{\line(1,0){0.13}}
\multiput(18.12,33.91)(0.2,-0.15){2}{\line(1,0){0.2}}
\multiput(18.53,33.62)(0.21,-0.14){2}{\line(1,0){0.21}}
\multiput(18.94,33.34)(0.21,-0.14){2}{\line(1,0){0.21}}
\multiput(19.36,33.07)(0.21,-0.13){2}{\line(1,0){0.21}}
\multiput(19.78,32.81)(0.22,-0.12){2}{\line(1,0){0.22}}
\multiput(20.22,32.56)(0.22,-0.12){2}{\line(1,0){0.22}}
\multiput(20.65,32.32)(0.22,-0.11){2}{\line(1,0){0.22}}
\multiput(21.1,32.09)(0.22,-0.11){2}{\line(1,0){0.22}}
\multiput(21.55,31.87)(0.23,-0.1){2}{\line(1,0){0.23}}
\multiput(22,31.67)(0.23,-0.1){2}{\line(1,0){0.23}}
\multiput(22.46,31.48)(0.23,-0.09){2}{\line(1,0){0.23}}
\multiput(22.93,31.29)(0.47,-0.17){1}{\line(1,0){0.47}}
\multiput(23.39,31.12)(0.47,-0.16){1}{\line(1,0){0.47}}
\multiput(23.87,30.96)(0.48,-0.15){1}{\line(1,0){0.48}}
\multiput(24.34,30.82)(0.48,-0.14){1}{\line(1,0){0.48}}
\multiput(24.82,30.68)(0.48,-0.12){1}{\line(1,0){0.48}}
\multiput(25.31,30.56)(0.49,-0.11){1}{\line(1,0){0.49}}
\multiput(25.79,30.45)(0.49,-0.1){1}{\line(1,0){0.49}}
\multiput(26.28,30.35)(0.49,-0.09){1}{\line(1,0){0.49}}
\multiput(26.77,30.26)(0.49,-0.07){1}{\line(1,0){0.49}}
\multiput(27.27,30.19)(0.49,-0.06){1}{\line(1,0){0.49}}
\multiput(27.76,30.13)(0.5,-0.05){1}{\line(1,0){0.5}}
\multiput(28.26,30.08)(0.5,-0.04){1}{\line(1,0){0.5}}
\multiput(28.75,30.04)(0.5,-0.02){1}{\line(1,0){0.5}}
\multiput(29.25,30.01)(0.5,-0.01){1}{\line(1,0){0.5}}
\put(29.75,30){\line(1,0){0.5}}
\multiput(30.25,30)(0.5,0.01){1}{\line(1,0){0.5}}
\multiput(30.75,30.01)(0.5,0.02){1}{\line(1,0){0.5}}
\multiput(31.25,30.04)(0.5,0.04){1}{\line(1,0){0.5}}
\multiput(31.74,30.08)(0.5,0.05){1}{\line(1,0){0.5}}
\multiput(32.24,30.13)(0.49,0.06){1}{\line(1,0){0.49}}
\multiput(32.73,30.19)(0.49,0.07){1}{\line(1,0){0.49}}
\multiput(33.23,30.26)(0.49,0.09){1}{\line(1,0){0.49}}
\multiput(33.72,30.35)(0.49,0.1){1}{\line(1,0){0.49}}
\multiput(34.21,30.45)(0.49,0.11){1}{\line(1,0){0.49}}
\multiput(34.69,30.56)(0.48,0.12){1}{\line(1,0){0.48}}
\multiput(35.18,30.68)(0.48,0.14){1}{\line(1,0){0.48}}
\multiput(35.66,30.82)(0.48,0.15){1}{\line(1,0){0.48}}
\multiput(36.13,30.96)(0.47,0.16){1}{\line(1,0){0.47}}
\multiput(36.61,31.12)(0.47,0.17){1}{\line(1,0){0.47}}
\multiput(37.07,31.29)(0.23,0.09){2}{\line(1,0){0.23}}
\multiput(37.54,31.48)(0.23,0.1){2}{\line(1,0){0.23}}
\multiput(38,31.67)(0.23,0.1){2}{\line(1,0){0.23}}
\multiput(38.45,31.87)(0.22,0.11){2}{\line(1,0){0.22}}
\multiput(38.9,32.09)(0.22,0.11){2}{\line(1,0){0.22}}
\multiput(39.35,32.32)(0.22,0.12){2}{\line(1,0){0.22}}
\multiput(39.78,32.56)(0.22,0.12){2}{\line(1,0){0.22}}
\multiput(40.22,32.81)(0.21,0.13){2}{\line(1,0){0.21}}
\multiput(40.64,33.07)(0.21,0.14){2}{\line(1,0){0.21}}
\multiput(41.06,33.34)(0.21,0.14){2}{\line(1,0){0.21}}
\multiput(41.47,33.62)(0.2,0.15){2}{\line(1,0){0.2}}
\multiput(41.88,33.91)(0.13,0.1){3}{\line(1,0){0.13}}
\multiput(42.27,34.21)(0.13,0.1){3}{\line(1,0){0.13}}
\multiput(42.66,34.52)(0.13,0.11){3}{\line(1,0){0.13}}
\multiput(43.05,34.84)(0.12,0.11){3}{\line(1,0){0.12}}
\multiput(43.42,35.17)(0.12,0.11){3}{\line(1,0){0.12}}
\multiput(43.79,35.51)(0.12,0.12){3}{\line(1,0){0.12}}
\multiput(44.14,35.86)(0.12,0.12){3}{\line(0,1){0.12}}
\multiput(44.49,36.21)(0.11,0.12){3}{\line(0,1){0.12}}
\multiput(44.83,36.58)(0.11,0.12){3}{\line(0,1){0.12}}
\multiput(45.16,36.95)(0.11,0.13){3}{\line(0,1){0.13}}
\multiput(45.48,37.34)(0.1,0.13){3}{\line(0,1){0.13}}
\multiput(45.79,37.73)(0.1,0.13){3}{\line(0,1){0.13}}
\multiput(46.09,38.12)(0.15,0.2){2}{\line(0,1){0.2}}
\multiput(46.38,38.53)(0.14,0.21){2}{\line(0,1){0.21}}
\multiput(46.66,38.94)(0.14,0.21){2}{\line(0,1){0.21}}
\multiput(46.93,39.36)(0.13,0.21){2}{\line(0,1){0.21}}
\multiput(47.19,39.78)(0.12,0.22){2}{\line(0,1){0.22}}
\multiput(47.44,40.22)(0.12,0.22){2}{\line(0,1){0.22}}
\multiput(47.68,40.65)(0.11,0.22){2}{\line(0,1){0.22}}
\multiput(47.91,41.1)(0.11,0.22){2}{\line(0,1){0.22}}
\multiput(48.13,41.55)(0.1,0.23){2}{\line(0,1){0.23}}
\multiput(48.33,42)(0.1,0.23){2}{\line(0,1){0.23}}
\multiput(48.52,42.46)(0.09,0.23){2}{\line(0,1){0.23}}
\multiput(48.71,42.93)(0.17,0.47){1}{\line(0,1){0.47}}
\multiput(48.88,43.39)(0.16,0.47){1}{\line(0,1){0.47}}
\multiput(49.04,43.87)(0.15,0.48){1}{\line(0,1){0.48}}
\multiput(49.18,44.34)(0.14,0.48){1}{\line(0,1){0.48}}
\multiput(49.32,44.82)(0.12,0.48){1}{\line(0,1){0.48}}
\multiput(49.44,45.31)(0.11,0.49){1}{\line(0,1){0.49}}
\multiput(49.55,45.79)(0.1,0.49){1}{\line(0,1){0.49}}
\multiput(49.65,46.28)(0.09,0.49){1}{\line(0,1){0.49}}
\multiput(49.74,46.77)(0.07,0.49){1}{\line(0,1){0.49}}
\multiput(49.81,47.27)(0.06,0.49){1}{\line(0,1){0.49}}
\multiput(49.87,47.76)(0.05,0.5){1}{\line(0,1){0.5}}
\multiput(49.92,48.26)(0.04,0.5){1}{\line(0,1){0.5}}
\multiput(49.96,48.75)(0.02,0.5){1}{\line(0,1){0.5}}
\multiput(49.99,49.25)(0.01,0.5){1}{\line(0,1){0.5}}

\linethickness{0.3mm}
\put(30,50){\line(1,0){100}}
\linethickness{0.3mm}
\put(40,20){\line(0,1){30}}
\linethickness{0.3mm}
\multiput(50,50)(0.36,0.12){83}{\line(1,0){0.36}}
\linethickness{0.3mm}
\put(80,60){\line(1,0){20}}
\linethickness{0.3mm}
\multiput(100,60)(0.36,-0.12){83}{\line(1,0){0.36}}
\linethickness{0.3mm}
\multiput(80,60)(0.12,0.12){83}{\line(1,0){0.12}}
\linethickness{0.3mm}
\multiput(90,70)(0.12,-0.12){83}{\line(1,0){0.12}}
\linethickness{0.3mm}
\multiput(80,60)(0.6,-0.12){83}{\line(1,0){0.6}}
\linethickness{0.3mm}
\multiput(50,50)(0.6,0.12){83}{\line(1,0){0.6}}
\put(30,55){\makebox(0,0)[cc]{$x$}}

\put(45,20){\makebox(0,0)[cc]{$x'$}}

\put(55,45){\makebox(0,0)[cc]{$u$}}

\put(77,63){\makebox(0,0)[cc]{$w$}}

\put(92,72){\makebox(0,0)[cc]{$y$}}

\put(100,63){\makebox(0,0)[cc]{$a$}}

\put(130,45){\makebox(0,0)[cc]{$b$}}

\end{picture}

  D'autre part $d(x,w)\leq d_{\max}(x,S_{m'+1})\leq d_{\max}(x,S_{m'})\leq d(x,a)+N$ et $d(x,a)-k\leq d(u,a)$. On en déduit 
  $d(u,w)\leq d(u,a)+(2F+N)$. Le lemme~\ref{x-a,b-y-geod} 
  appliqué à $(u,b,a,w)$ au lieu de $(x,y,a,b)$ et à $(2F,2F,2F+N)$ au lieu de $(\alpha,\beta,\rho)$ 
  montre alors $w\in (6F+2N+\de)\tg(u,a)$. On en déduit $y\in (6F+2N+\de+P)\tg(u,a)$ par le lemme~\ref{geod-comp-xabc}
  et donc $y\in 2P\tg(u,a)$ car on suppose $P\geq 6F+2N+\de$ grâce à  
  $(H_{P})$.  On  suppose $\frac{M}{2}\geq 2N+\de$ grâce à  $(H_{M})$. 
  Par (\ref{ineg-k'-deuxiemecas-12dec}) on a alors $u\in B(x',k')$. En 
     prenant $\tilde x'=u$ on a fini. 
 
 On suppose maintenant $m'=m$. Il existe  $a\in S_{m}$ et $\tilde x\in B(x,k)$ tels que $y\in 2P\tg(\tilde x,a)$. On veut montrer qu'il existe 
 $\tilde x'\in B(x',k')$ tel que $y\in 2P\tg(\tilde x',a)$.  
 On a $d(x,y)\geq k+M+P+\de$. On applique le sous-lemme~\ref{souslem1} avec $\alpha=2P$ et $\beta=M+P+\de$. On suppose $F+\frac{3\de+\alpha}{2}< \beta$, ce qui est permis par $(H_{M})$.

 Il reste à montrer d). 
 Le premier ensemble de d) est  inclus dans $B(x,k+2M+2P+\de)$ et 
 le second ensemble de d) est inclus dans $(2M+4P)\tg(u,b)$. 
On a  $B(u,\frac{5M}{2}-2N-3\de)
 \subset B(x',k'+2M-2\de)$ grâce à (\ref{ineg-k'-deuxiemecas-12dec}). 
 Il suffit donc de montrer que pour $z\in B(x,k+2M+2P+\de)$ et $z'\in (2M+4P)\tg(u,b)$, $\geod(z,z')$ rencontre $B(u,\frac{5M}{2}-2N-3\de)$. 
 Par $(H_{\de}^{0}(z,x,u,b))$, on a  $u\in \de\tg(x,z)$ ou $ u\in \de\tg(z,b)$. 
 Si $u\in \de\tg(x,z)$, $d(u,z)\leq 2M+2P+2\de$, donc $z\in \geod(z,z')$ appartient à $B(u,\frac{5M}{2}-2N-3\de)$
 car on  suppose $\frac{M}{2}\geq 2P+2N+5\de$ 
  grâce à $(H_{M})$. Si $ u\in \de\tg(z,b)$, comme $z'\in (2M+4P)\tg(u,b)$, on a $u\in (2M+4P+\de)\tg(z,z')$ par le b) du lemme~\ref{geod-comp-xabc}, donc $t\in \geod(z,z')$ vérifiant $d(z,t)=\min(d(z,u),d(z,z'))$ satisfait, grâce à $(H_{\de}^{0}(u,z,t,z'))$,  
 $$d(u,t)\leq 2M+4P+2\de\leq \frac{5M}{2}-2N-3\de$$ car on suppose $2M+4P+2\de\leq \frac{5M}{2}-2N-3\de$ grâce à $(H_{M})$. Cela termine la preuve de d) et donc l'étude du deuxième cas.

\noindent {\bf Troisième  cas.} 
On suppose $r_{0}(Z)> k$ et $\frac{r_{0}(Z)-r_{0}'(Z)+d(x,x')}{2}> k$. 

Par 
(\ref{form-r'-cas3-4j})
%(\ref{def-k'Z-24oct09}) 
on a alors 
\begin{gather}\label{ineg-k'-troisiemecas-19dec}\big|k'-\big(E(\frac{r'_{0}(Z)-r_{0}(Z)+d(x,x')}{2})+\frac{M}{2}\big)\big|\leq \frac{N}{2}.\end{gather}
Comme $\max(k,\frac{d(x,x')+r_{0}(Z)-r'_{0}(Z)}{2})=\frac{d(x,x')+r_{0}(Z)-r'_{0}(Z)}{2}$, 
pour passer de \begin{gather*}(a_{1},\dots,a_{p},S_{0},...,S_{m},(\mathcal Y_{i}^{j})_{i\in \{0,\dots,m\}, j\in \{1,\dots ,l_{i}\}})\\ \text{à  \ \ \ }(a_{1},\dots,a_{p},S_{0},...,S_{m'},(\mathcal Y_{i}^{j})_{i\in \{0,\dots,m'-1\}, j\in \{1,\dots ,l_{i}\}}, (\mathcal Y_{m'}^{j_{\lambda}})_{ \lambda\in \{1,\dots ,l'_{m'}\}})\end{gather*} on enlève les $S_{i}$ et les $\mathcal Y_{i}^{j}$ pour $i\geq1$ vérifiant  $$d_{\max}(x,S_{i})\leq \frac{r_{0}(Z)-r'_{0}(Z)+d(x,x')}{2}+M$$ ainsi que  les $\mathcal Y_{m'}^{j}$ vérifiant $$d_{\max}(x,\mathcal Y_{m'}^{j})\leq \frac{r_{0}(Z)-r'_{0}(Z)+d(x,x')}{2}+M+2P+\de. $$
Soit $u\in B(x,k)$ à distance minimale de $S_{0}$ et 
$b\in S_{0}$ à distance minimale de $u$. Il résulte de (\ref{ineg-k'-troisiemecas-19dec}) que 
\begin{gather}\label{ineg1-k'-troisiemecas-19dec} k'\geq \frac{d(x',b)-d(x,b)+d(x,x')}{2}+\frac{M}{2}-N \\ \label{ineg2-k'-troisiemecas-19dec} \text{et\ \ \ } k'\leq \frac{d(x',b)-d(x,b)+d(x,x')}{2}+\frac{M}{2}+N\end{gather}
Pour vérifier c) on va montrer que les conditions i), ii), iii) et iv)  de la définition~\ref{defi-Y} sont satisfaites par $$(a_{1},\dots,a_{p},S_{0},...,S_{m'},(\mathcal Y_{i}^{j})_{i\in \{0,\dots,m'-1\}, j\in \{1,\dots ,l_{i}\}}, (\mathcal Y_{m'}^{j_{\lambda}})_{ \lambda\in \{1,\dots ,l'_{m'}\}}).$$  On vérifie d'abord i). 
Soit $i\in \{0,...,m'-1\}$, $y\in S_{i+1}$, $a\in S_{i}$ et $ \tilde x\in B(x,k)$    tels que 
\begin{itemize}
\item $y\in 4\de\tg(\tilde x,a)$ et $d(y,a)\in ]N-2\de,QN]$, 
\item ou $y\in F\tg(\tilde x,a)$ et $d(y,a)\geq \frac{Q}{F}$. 
\end{itemize}
On note $\alpha=4\de$ ou $F$ suivant le cas, de sorte que $y\in \alpha\tg(\tilde x,a)$. 
Comme $d(x,S_{i+1})\geq d_{\max}(x,S_{i+1})-N$ et $i+1\in \{1,...,m'\}$,  on a 
\begin{gather}\nonumber d(x,y)\geq \frac{r_{0}(Z)-r'_{0}(Z)+d(x,x')}{2}+M-N\\ \label{cond-dxy-p126-24oct09} \geq \frac{d(x,b)-d(x',b)+d(x,x')}{2}+M-\frac{3N}{2} .\end{gather} 
%Soit $\alpha=4\de$ ou $\alpha=5\de$. Soit $\tilde x\in B(x,k)$ et $y\in \alpha\tg(\tilde x,a)$, tel que $$d(x,y)\geq \frac{r_{0}(Z)-r'_{0}(Z)+d(x,x')}{2}+M-2N-5\de p_{\max}$$ 
%(comme $d(x,S_{i+1})\geq d_{\max}(x,S_{i+1})-N$, cette condition est vérifiée si $\alpha=4\de$ et $y\in S_{i+1}$ ou si $\alpha=5\de$ et 
%$B(y,N+5\de p_{\max})\cap S_{i+1}\neq \emptyset$). 
   \noindent Pour vérifier i) il suffit de montrer qu'il existe $\tilde x'\in B(x', k')$ tel que $y\in \alpha\tg(\tilde x',a)$.  
  
\ifx\JPicScale\undefined\def\JPicScale{1}\fi
\unitlength \JPicScale mm
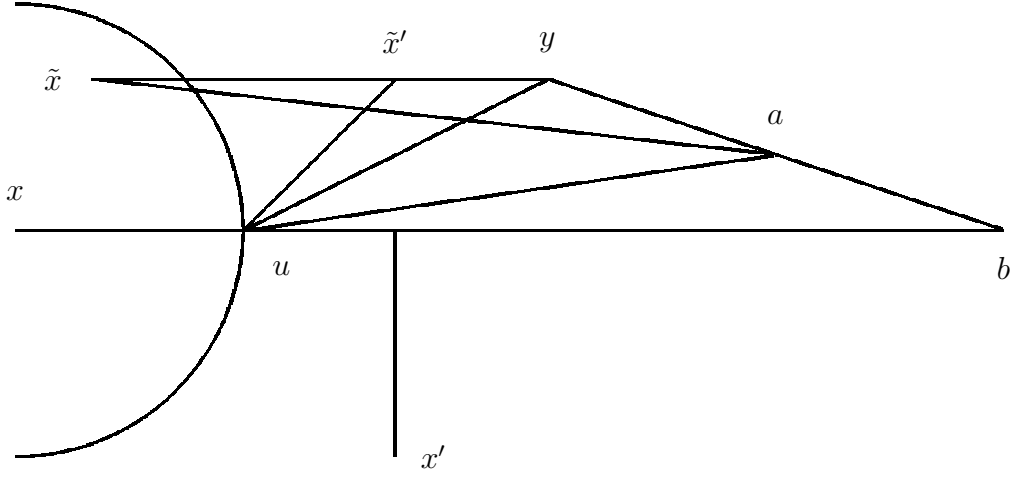
\begin{picture}(120,80)(10,10)
\linethickness{0.3mm}
\multiput(10,20)(0.5,0){1}{\line(1,0){0.5}}
\multiput(10.5,20)(0.5,0.01){1}{\line(1,0){0.5}}
\multiput(11,20.02)(0.5,0.02){1}{\line(1,0){0.5}}
\multiput(11.5,20.04)(0.5,0.03){1}{\line(1,0){0.5}}
\multiput(12,20.07)(0.5,0.04){1}{\line(1,0){0.5}}
\multiput(12.5,20.1)(0.5,0.05){1}{\line(1,0){0.5}}
\multiput(13,20.15)(0.5,0.05){1}{\line(1,0){0.5}}
\multiput(13.5,20.21)(0.5,0.06){1}{\line(1,0){0.5}}
\multiput(14,20.27)(0.5,0.07){1}{\line(1,0){0.5}}
\multiput(14.49,20.34)(0.5,0.08){1}{\line(1,0){0.5}}
\multiput(14.99,20.42)(0.49,0.09){1}{\line(1,0){0.49}}
\multiput(15.48,20.51)(0.49,0.1){1}{\line(1,0){0.49}}
\multiput(15.98,20.6)(0.49,0.1){1}{\line(1,0){0.49}}
\multiput(16.47,20.71)(0.49,0.11){1}{\line(1,0){0.49}}
\multiput(16.95,20.82)(0.49,0.12){1}{\line(1,0){0.49}}
\multiput(17.44,20.94)(0.48,0.13){1}{\line(1,0){0.48}}
\multiput(17.93,21.07)(0.48,0.14){1}{\line(1,0){0.48}}
\multiput(18.41,21.2)(0.48,0.14){1}{\line(1,0){0.48}}
\multiput(18.89,21.35)(0.48,0.15){1}{\line(1,0){0.48}}
\multiput(19.37,21.5)(0.47,0.16){1}{\line(1,0){0.47}}
\multiput(19.84,21.66)(0.47,0.17){1}{\line(1,0){0.47}}
\multiput(20.31,21.83)(0.47,0.18){1}{\line(1,0){0.47}}
\multiput(20.78,22)(0.23,0.09){2}{\line(1,0){0.23}}
\multiput(21.25,22.19)(0.23,0.1){2}{\line(1,0){0.23}}
\multiput(21.71,22.38)(0.23,0.1){2}{\line(1,0){0.23}}
\multiput(22.17,22.58)(0.23,0.1){2}{\line(1,0){0.23}}
\multiput(22.63,22.79)(0.23,0.11){2}{\line(1,0){0.23}}
\multiput(23.08,23)(0.22,0.11){2}{\line(1,0){0.22}}
\multiput(23.53,23.22)(0.22,0.11){2}{\line(1,0){0.22}}
\multiput(23.98,23.45)(0.22,0.12){2}{\line(1,0){0.22}}
\multiput(24.42,23.69)(0.22,0.12){2}{\line(1,0){0.22}}
\multiput(24.86,23.94)(0.22,0.13){2}{\line(1,0){0.22}}
\multiput(25.29,24.19)(0.21,0.13){2}{\line(1,0){0.21}}
\multiput(25.72,24.45)(0.21,0.13){2}{\line(1,0){0.21}}
\multiput(26.14,24.71)(0.21,0.14){2}{\line(1,0){0.21}}
\multiput(26.56,24.99)(0.21,0.14){2}{\line(1,0){0.21}}
\multiput(26.98,25.27)(0.21,0.14){2}{\line(1,0){0.21}}
\multiput(27.39,25.55)(0.2,0.15){2}{\line(1,0){0.2}}
\multiput(27.8,25.85)(0.13,0.1){3}{\line(1,0){0.13}}
\multiput(28.2,26.15)(0.13,0.1){3}{\line(1,0){0.13}}
\multiput(28.59,26.46)(0.13,0.1){3}{\line(1,0){0.13}}
\multiput(28.98,26.77)(0.13,0.11){3}{\line(1,0){0.13}}
\multiput(29.37,27.09)(0.13,0.11){3}{\line(1,0){0.13}}
\multiput(29.75,27.42)(0.12,0.11){3}{\line(1,0){0.12}}
\multiput(30.12,27.75)(0.12,0.11){3}{\line(1,0){0.12}}
\multiput(30.49,28.09)(0.12,0.12){3}{\line(1,0){0.12}}
\multiput(30.86,28.44)(0.12,0.12){3}{\line(1,0){0.12}}
\multiput(31.21,28.79)(0.12,0.12){3}{\line(0,1){0.12}}
\multiput(31.56,29.14)(0.12,0.12){3}{\line(0,1){0.12}}
\multiput(31.91,29.51)(0.11,0.12){3}{\line(0,1){0.12}}
\multiput(32.25,29.88)(0.11,0.12){3}{\line(0,1){0.12}}
\multiput(32.58,30.25)(0.11,0.13){3}{\line(0,1){0.13}}
\multiput(32.91,30.63)(0.11,0.13){3}{\line(0,1){0.13}}
\multiput(33.23,31.02)(0.1,0.13){3}{\line(0,1){0.13}}
\multiput(33.54,31.41)(0.1,0.13){3}{\line(0,1){0.13}}
\multiput(33.85,31.8)(0.1,0.13){3}{\line(0,1){0.13}}
\multiput(34.15,32.2)(0.15,0.2){2}{\line(0,1){0.2}}
\multiput(34.45,32.61)(0.14,0.21){2}{\line(0,1){0.21}}
\multiput(34.73,33.02)(0.14,0.21){2}{\line(0,1){0.21}}
\multiput(35.01,33.44)(0.14,0.21){2}{\line(0,1){0.21}}
\multiput(35.29,33.86)(0.13,0.21){2}{\line(0,1){0.21}}
\multiput(35.55,34.28)(0.13,0.21){2}{\line(0,1){0.21}}
\multiput(35.81,34.71)(0.13,0.22){2}{\line(0,1){0.22}}
\multiput(36.06,35.14)(0.12,0.22){2}{\line(0,1){0.22}}
\multiput(36.31,35.58)(0.12,0.22){2}{\line(0,1){0.22}}
\multiput(36.55,36.02)(0.11,0.22){2}{\line(0,1){0.22}}
\multiput(36.78,36.47)(0.11,0.22){2}{\line(0,1){0.22}}
\multiput(37,36.92)(0.11,0.23){2}{\line(0,1){0.23}}
\multiput(37.21,37.37)(0.1,0.23){2}{\line(0,1){0.23}}
\multiput(37.42,37.83)(0.1,0.23){2}{\line(0,1){0.23}}
\multiput(37.62,38.29)(0.1,0.23){2}{\line(0,1){0.23}}
\multiput(37.81,38.75)(0.09,0.23){2}{\line(0,1){0.23}}
\multiput(38,39.22)(0.18,0.47){1}{\line(0,1){0.47}}
\multiput(38.17,39.69)(0.17,0.47){1}{\line(0,1){0.47}}
\multiput(38.34,40.16)(0.16,0.47){1}{\line(0,1){0.47}}
\multiput(38.5,40.63)(0.15,0.48){1}{\line(0,1){0.48}}
\multiput(38.65,41.11)(0.14,0.48){1}{\line(0,1){0.48}}
\multiput(38.8,41.59)(0.14,0.48){1}{\line(0,1){0.48}}
\multiput(38.93,42.07)(0.13,0.48){1}{\line(0,1){0.48}}
\multiput(39.06,42.56)(0.12,0.49){1}{\line(0,1){0.49}}
\multiput(39.18,43.05)(0.11,0.49){1}{\line(0,1){0.49}}
\multiput(39.29,43.53)(0.1,0.49){1}{\line(0,1){0.49}}
\multiput(39.4,44.02)(0.1,0.49){1}{\line(0,1){0.49}}
\multiput(39.49,44.52)(0.09,0.49){1}{\line(0,1){0.49}}
\multiput(39.58,45.01)(0.08,0.5){1}{\line(0,1){0.5}}
\multiput(39.66,45.51)(0.07,0.5){1}{\line(0,1){0.5}}
\multiput(39.73,46)(0.06,0.5){1}{\line(0,1){0.5}}
\multiput(39.79,46.5)(0.05,0.5){1}{\line(0,1){0.5}}
\multiput(39.85,47)(0.05,0.5){1}{\line(0,1){0.5}}
\multiput(39.9,47.5)(0.04,0.5){1}{\line(0,1){0.5}}
\multiput(39.93,48)(0.03,0.5){1}{\line(0,1){0.5}}
\multiput(39.96,48.5)(0.02,0.5){1}{\line(0,1){0.5}}
\multiput(39.98,49)(0.01,0.5){1}{\line(0,1){0.5}}
\multiput(40,49.5)(0,0.5){1}{\line(0,1){0.5}}
\multiput(40,50.5)(0,-0.5){1}{\line(0,-1){0.5}}
\multiput(39.98,51)(0.01,-0.5){1}{\line(0,-1){0.5}}
\multiput(39.96,51.5)(0.02,-0.5){1}{\line(0,-1){0.5}}
\multiput(39.93,52)(0.03,-0.5){1}{\line(0,-1){0.5}}
\multiput(39.9,52.5)(0.04,-0.5){1}{\line(0,-1){0.5}}
\multiput(39.85,53)(0.05,-0.5){1}{\line(0,-1){0.5}}
\multiput(39.79,53.5)(0.05,-0.5){1}{\line(0,-1){0.5}}
\multiput(39.73,54)(0.06,-0.5){1}{\line(0,-1){0.5}}
\multiput(39.66,54.49)(0.07,-0.5){1}{\line(0,-1){0.5}}
\multiput(39.58,54.99)(0.08,-0.5){1}{\line(0,-1){0.5}}
\multiput(39.49,55.48)(0.09,-0.49){1}{\line(0,-1){0.49}}
\multiput(39.4,55.98)(0.1,-0.49){1}{\line(0,-1){0.49}}
\multiput(39.29,56.47)(0.1,-0.49){1}{\line(0,-1){0.49}}
\multiput(39.18,56.95)(0.11,-0.49){1}{\line(0,-1){0.49}}
\multiput(39.06,57.44)(0.12,-0.49){1}{\line(0,-1){0.49}}
\multiput(38.93,57.93)(0.13,-0.48){1}{\line(0,-1){0.48}}
\multiput(38.8,58.41)(0.14,-0.48){1}{\line(0,-1){0.48}}
\multiput(38.65,58.89)(0.14,-0.48){1}{\line(0,-1){0.48}}
\multiput(38.5,59.37)(0.15,-0.48){1}{\line(0,-1){0.48}}
\multiput(38.34,59.84)(0.16,-0.47){1}{\line(0,-1){0.47}}
\multiput(38.17,60.31)(0.17,-0.47){1}{\line(0,-1){0.47}}
\multiput(38,60.78)(0.18,-0.47){1}{\line(0,-1){0.47}}
\multiput(37.81,61.25)(0.09,-0.23){2}{\line(0,-1){0.23}}
\multiput(37.62,61.71)(0.1,-0.23){2}{\line(0,-1){0.23}}
\multiput(37.42,62.17)(0.1,-0.23){2}{\line(0,-1){0.23}}
\multiput(37.21,62.63)(0.1,-0.23){2}{\line(0,-1){0.23}}
\multiput(37,63.08)(0.11,-0.23){2}{\line(0,-1){0.23}}
\multiput(36.78,63.53)(0.11,-0.22){2}{\line(0,-1){0.22}}
\multiput(36.55,63.98)(0.11,-0.22){2}{\line(0,-1){0.22}}
\multiput(36.31,64.42)(0.12,-0.22){2}{\line(0,-1){0.22}}
\multiput(36.06,64.86)(0.12,-0.22){2}{\line(0,-1){0.22}}
\multiput(35.81,65.29)(0.13,-0.22){2}{\line(0,-1){0.22}}
\multiput(35.55,65.72)(0.13,-0.21){2}{\line(0,-1){0.21}}
\multiput(35.29,66.14)(0.13,-0.21){2}{\line(0,-1){0.21}}
\multiput(35.01,66.56)(0.14,-0.21){2}{\line(0,-1){0.21}}
\multiput(34.73,66.98)(0.14,-0.21){2}{\line(0,-1){0.21}}
\multiput(34.45,67.39)(0.14,-0.21){2}{\line(0,-1){0.21}}
\multiput(34.15,67.8)(0.15,-0.2){2}{\line(0,-1){0.2}}
\multiput(33.85,68.2)(0.1,-0.13){3}{\line(0,-1){0.13}}
\multiput(33.54,68.59)(0.1,-0.13){3}{\line(0,-1){0.13}}
\multiput(33.23,68.98)(0.1,-0.13){3}{\line(0,-1){0.13}}
\multiput(32.91,69.37)(0.11,-0.13){3}{\line(0,-1){0.13}}
\multiput(32.58,69.75)(0.11,-0.13){3}{\line(0,-1){0.13}}
\multiput(32.25,70.12)(0.11,-0.12){3}{\line(0,-1){0.12}}
\multiput(31.91,70.49)(0.11,-0.12){3}{\line(0,-1){0.12}}
\multiput(31.56,70.86)(0.12,-0.12){3}{\line(0,-1){0.12}}
\multiput(31.21,71.21)(0.12,-0.12){3}{\line(0,-1){0.12}}
\multiput(30.86,71.56)(0.12,-0.12){3}{\line(1,0){0.12}}
\multiput(30.49,71.91)(0.12,-0.12){3}{\line(1,0){0.12}}
\multiput(30.12,72.25)(0.12,-0.11){3}{\line(1,0){0.12}}
\multiput(29.75,72.58)(0.12,-0.11){3}{\line(1,0){0.12}}
\multiput(29.37,72.91)(0.13,-0.11){3}{\line(1,0){0.13}}
\multiput(28.98,73.23)(0.13,-0.11){3}{\line(1,0){0.13}}
\multiput(28.59,73.54)(0.13,-0.1){3}{\line(1,0){0.13}}
\multiput(28.2,73.85)(0.13,-0.1){3}{\line(1,0){0.13}}
\multiput(27.8,74.15)(0.13,-0.1){3}{\line(1,0){0.13}}
\multiput(27.39,74.45)(0.2,-0.15){2}{\line(1,0){0.2}}
\multiput(26.98,74.73)(0.21,-0.14){2}{\line(1,0){0.21}}
\multiput(26.56,75.01)(0.21,-0.14){2}{\line(1,0){0.21}}
\multiput(26.14,75.29)(0.21,-0.14){2}{\line(1,0){0.21}}
\multiput(25.72,75.55)(0.21,-0.13){2}{\line(1,0){0.21}}
\multiput(25.29,75.81)(0.21,-0.13){2}{\line(1,0){0.21}}
\multiput(24.86,76.06)(0.22,-0.13){2}{\line(1,0){0.22}}
\multiput(24.42,76.31)(0.22,-0.12){2}{\line(1,0){0.22}}
\multiput(23.98,76.55)(0.22,-0.12){2}{\line(1,0){0.22}}
\multiput(23.53,76.78)(0.22,-0.11){2}{\line(1,0){0.22}}
\multiput(23.08,77)(0.22,-0.11){2}{\line(1,0){0.22}}
\multiput(22.63,77.21)(0.23,-0.11){2}{\line(1,0){0.23}}
\multiput(22.17,77.42)(0.23,-0.1){2}{\line(1,0){0.23}}
\multiput(21.71,77.62)(0.23,-0.1){2}{\line(1,0){0.23}}
\multiput(21.25,77.81)(0.23,-0.1){2}{\line(1,0){0.23}}
\multiput(20.78,78)(0.23,-0.09){2}{\line(1,0){0.23}}
\multiput(20.31,78.17)(0.47,-0.18){1}{\line(1,0){0.47}}
\multiput(19.84,78.34)(0.47,-0.17){1}{\line(1,0){0.47}}
\multiput(19.37,78.5)(0.47,-0.16){1}{\line(1,0){0.47}}
\multiput(18.89,78.65)(0.48,-0.15){1}{\line(1,0){0.48}}
\multiput(18.41,78.8)(0.48,-0.14){1}{\line(1,0){0.48}}
\multiput(17.93,78.93)(0.48,-0.14){1}{\line(1,0){0.48}}
\multiput(17.44,79.06)(0.48,-0.13){1}{\line(1,0){0.48}}
\multiput(16.95,79.18)(0.49,-0.12){1}{\line(1,0){0.49}}
\multiput(16.47,79.29)(0.49,-0.11){1}{\line(1,0){0.49}}
\multiput(15.98,79.4)(0.49,-0.1){1}{\line(1,0){0.49}}
\multiput(15.48,79.49)(0.49,-0.1){1}{\line(1,0){0.49}}
\multiput(14.99,79.58)(0.49,-0.09){1}{\line(1,0){0.49}}
\multiput(14.49,79.66)(0.5,-0.08){1}{\line(1,0){0.5}}
\multiput(14,79.73)(0.5,-0.07){1}{\line(1,0){0.5}}
\multiput(13.5,79.79)(0.5,-0.06){1}{\line(1,0){0.5}}
\multiput(13,79.85)(0.5,-0.05){1}{\line(1,0){0.5}}
\multiput(12.5,79.9)(0.5,-0.05){1}{\line(1,0){0.5}}
\multiput(12,79.93)(0.5,-0.04){1}{\line(1,0){0.5}}
\multiput(11.5,79.96)(0.5,-0.03){1}{\line(1,0){0.5}}
\multiput(11,79.98)(0.5,-0.02){1}{\line(1,0){0.5}}
\multiput(10.5,80)(0.5,-0.01){1}{\line(1,0){0.5}}
\multiput(10,80)(0.5,-0){1}{\line(1,0){0.5}}

\linethickness{0.3mm}
\put(10,50){\line(1,0){130}}
\linethickness{0.3mm}
\multiput(40,50)(0.84,0.12){83}{\line(1,0){0.84}}
\linethickness{0.3mm}
\multiput(110,60)(0.36,-0.12){83}{\line(1,0){0.36}}
\linethickness{0.3mm}
\put(20,70){\line(1,0){60}}
\linethickness{0.3mm}
\multiput(80,70)(0.36,-0.12){83}{\line(1,0){0.36}}
\linethickness{0.3mm}
\multiput(40,50)(0.24,0.12){167}{\line(1,0){0.24}}
\linethickness{0.3mm}
\multiput(40,50)(0.12,0.12){167}{\line(1,0){0.12}}
\linethickness{0.3mm}
\multiput(20,70)(1.08,-0.12){83}{\line(1,0){1.08}}
\linethickness{0.3mm}
\put(60,20){\line(0,1){30}}
\linethickness{0.3mm}
\put(10,55){\makebox(0,0)[cc]{$x$}}

\put(45,45){\makebox(0,0)[cc]{$u$}}

\put(15,70){\makebox(0,0)[cc]{$\tilde x$}}

\put(60,75){\makebox(0,0)[cc]{$\tilde x'$}}

\put(80,75){\makebox(0,0)[cc]{$y$}}

\put(110,65){\makebox(0,0)[cc]{$a$}}

\put(140,45){\makebox(0,0)[cc]{$b$}}

\put(65,20){\makebox(0,0)[cc]{$x'$}}

\end{picture}

Cela résulte du sous-lemme suivant appliqué à $\alpha=4\de$ ou $\alpha=F$ et à $\beta=M-\frac{3N}{2}$. On a $F\geq 4\de$ et on suppose  $F+\frac{3\de+F}{2}+\de+N\leq  M-\frac{3N}{2}$, ce qui est permis par $(H_{M})$. Les hypothèses du sous-lemme suivant sont donc satisfaites et cela termine la preuve de i). 

\begin{souslem}\label{souslem2}
Soit $\alpha,\beta \in \N$ vérifiant $F+\frac{3\de+\alpha}{2}+\de+N\leq \beta$. 
Soit $a\in 2F\tg(u,b)$, $\tilde x\in B(x,k)$ et  $y \in \alpha\tg(\tilde x,a)$ tel que $$d(x,y)\geq \frac{d(x,b)-d(x',b)+d(x,x')}{2}+\beta.$$ Alors il existe $\tilde x'\in B(x',k')$ tel que $y\in \alpha\tg(\tilde x',a)$.
\end{souslem}

\noindent{\bf Démonstration du sous-lemme~\ref{souslem2}.}
Comme $a\in 2F\tg(b,u)$ et $u\in \de\tg(\tilde x,b)$ par le lemme~\ref{Bxk-x-t-z-y}, on a $a\in (2F+\de)\tg(\tilde x,b)$  par le lemme~\ref{geod-comp-xabc}, 
et comme $y\in \alpha\tg(\tilde x,a)$ on en déduit  $y\in (2F+\de+\alpha)\tg(\tilde x,b)$. On applique le lemme~\ref{approx-arbres}
 à $\{x,b,x',\tilde x, y\}$ avec $l=2$ et $x$ comme point base.  
Soit $t\in \geod(\Psi \tilde x, \Psi b)$ à distance minimale de 
$\Psi y$. Comme $y\in (2F+\de+\alpha)\tg(\tilde x,b)$, on a $\Psi y\in 
 (2F+3\de+\alpha)\tg(\Psi \tilde x,\Psi b)$ par le a) du lemme~\ref{approx-arbres2}, donc $d(\Psi y,t)\leq F+\frac{3\de+\alpha}{2}$. Or \begin{gather*}d(\Psi x, \Psi y)=d(x,y)\geq \frac{d(x,b)-d(x',b)+d(x,x')}{2}+\beta\\ \geq k-\frac{N}{2}+\beta> k+(F+\frac{3\de+\alpha}{2})\end{gather*} par hypothèse. Donc $d(\Psi x, t)> k$ et comme $d(\Psi x, \Psi \tilde x)\leq k$, $t\in \geod (\Psi x, \Psi b)$. 
 
 \ifx\JPicScale\undefined\def\JPicScale{1}\fi
\unitlength \JPicScale mm
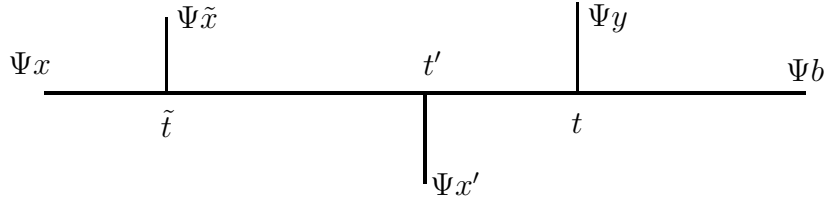
\begin{picture}(100,35)(22,15)
\linethickness{0.3mm}
\put(30,30){\line(1,0){100}}
\linethickness{0.3mm}
\put(46,30){\line(0,1){10}}
\linethickness{0.3mm}
\put(80,18){\line(0,1){12}}
\linethickness{0.3mm}
\put(100,30){\line(0,1){12}}
\put(28,34){\makebox(0,0)[cc]{$\Psi x$}}

\put(50,40){\makebox(0,0)[cc]{$\Psi \tilde x$}}

\put(46,26){\makebox(0,0)[cc]{$\tilde  t$}}

\put(81,34){\makebox(0,0)[cc]{$t'$}}

\put(84,18){\makebox(0,0)[cc]{$\Psi  x'$}}

\put(104,40){\makebox(0,0)[cc]{$\Psi  y$}}

\put(130,33){\makebox(0,0)[cc]{$\Psi  b$}}

\put(100,26){\makebox(0,0)[cc]{$t$}}

\end{picture}

\noindent De ce qui précède on retient aussi que 
\begin{gather}\label{ineg-phixt-19dec}d(\Psi x, t)\geq \frac{d(x,b)-d(x',b)+d(x,x')}{2}+\beta -(F+\frac{3\de+\alpha}{2}).\end{gather} 
Soit $t'$ le point de $\geod(\Psi x,\Psi b)$ à distance minimale de $\Psi x'$. 
On a donc \begin{gather}\nonumber d(\Psi x, t')=\frac{d(\Psi x,\Psi b)-d(\Psi x',\Psi b)+d(\Psi x,\Psi x')}{2} \\ \label{ineg-11j1831}
\in [\frac{d(x,b)-d(x',b)+d(x,x')}{2}, \frac{d(x,b)-d(x',b)+d(x,x')}{2}+\de]\end{gather}
par le lemme~\ref{approx-arbres}.  
Comme $F+\frac{3\de+\alpha}{2}+\de+N\leq \beta$ par hypothèse, il résulte  de (\ref{ineg-phixt-19dec}) et (\ref{ineg-11j1831}) que 
$d(\Psi x, t')\leq d(\Psi x, t)$. Par ailleurs $d(\Psi x, t')\geq 
\frac{d(x,b)-d(x',b)+d(x,x')}{2}> k-N$. Si $\tilde t $ désigne le point de 
$\geod(\Psi x,\Psi b)$ à distance minimale de $\Psi \tilde x$, comme 
$d(\Psi x, \Psi \tilde x)\leq k$, on a $d(\Psi x, \tilde t)\leq k$, donc 
 $d(\Psi x, t')\geq d(\Psi x, \tilde t)-N$ (dans le dessin ci-dessus  $t'$ pourrait être entre $\Psi x$ et $\tilde t$, mais dans ce cas à distance $\leq N$ de $\tilde t$). En tous cas on a $d(t', \geod(\tilde t, t))\leq N$, d'où $d(t', \geod(\Psi \tilde x,\Psi y))\leq N$.  On en déduit 
 \begin{gather*}
d(\Psi x', \geod(\Psi \tilde x,\Psi y))\leq d(\Psi x',t')+d(t', \geod(\Psi \tilde x,\Psi y))\leq d(\Psi x',t')+N\\ =\frac{d(\Psi x,\Psi x')+d(\Psi x',\Psi b)-d(\Psi x,\Psi b)}{2}+N\leq \frac{d(x,x')+d(x',b)-d(x,b)}{2}+N.\end{gather*}
Par le c) du lemme~\ref{approx-arbres2}, 
\begin{gather*}d(x', \geod(\tilde x,y))\leq d(\Psi x', \geod(\Psi \tilde x,\Psi y))+3\de+1\\ \leq \frac{d(x,x')+d(x',b)-d(x,b)}{2}+N+3\de+1\leq k'\end{gather*} grâce à (\ref{ineg1-k'-troisiemecas-19dec}) et car on suppose $2N+3\de+1\leq \frac{M}{2}$, ce qui est permis par $(H_{M})$. Il existe donc $\tilde x'\in \geod(\tilde x,y)\cap B(x',k')$. On a alors $y\in \alpha\tg(\tilde x', a)$. \cqfd

\noindent{\bf Fin de l'étude du troisième cas.} 
On vérifie maintenant ii). Soit $i\in \{1,...,m'\}$. On a $$S_{i}\subset 2F\tg(b,u)\subset 2F\tg(x,b)$$ par le a) du lemme~\ref{lemme-S0-...Sm}. On a $d_{\max}(x,S_{i})> \frac{r_{0}(Z)-r'_{0}(Z)+d(x,x')}{2}+M$ par hypothèse. 
 Soit $a\in S_{i}$. On a donc \begin{gather}\nonumber d(x,a)>  \frac{r_{0}(Z)-r'_{0}(Z)+d(x,x')}{2}+M-N\\ \label{ineg-dxaa-19dec}\geq \frac{d(x,b)-d(x',b)+d(x,x')}{2}+M-2N. \end{gather}
Donc $d(b,a)\leq d(x,b)-d(x,a)+2F< \frac{d(x,b)+d(x',b)-d(x,x')}{2}-M+(2N+2F)$ et  
$$d(x',a)\geq 
d(x',b)-d(b,a)> 
\frac{d(x',b)+d(x,x')-d(x,b)}{2}+M-(2N+2F).$$ 
Grâce à (\ref{ineg2-k'-troisiemecas-19dec}) on en déduit 
$$d(x',a)> k'+\frac{M}{2}-(3N+2F)
\geq k'+P$$ car on suppose $ \frac{M}{2}\geq (3N+2F)+P$, ce qui est permis par  $(H_{M})$. Ceci achève la preuve de ii). 

La propriété iii) est immédiate. 

On vérifie maintenant iv).
On note $b$ et $u$ comme dans la preuve de i).  Soit $j\in J$. On commence par montrer $d(x',
 \mathcal Y_{m'}^{j})\geq k'+3P$. On rappelle que $\mathcal Y_{m'}^{j}\subset 4P\tg(x,b)$ et \begin{gather}\nonumber d(x,\mathcal Y_{m'}^{j})\geq 
 \frac{r_{0}(Z)-r_{0}'(Z)+d(x,x')}{2}
  +M+P+\de\\ \label{ineg-dxYm'j-19dec}\geq \frac{d(x,b)-d(x',b)+d(x,x')}{2}+M+P+\de-\frac{N}{2}.\end{gather}
   Soit $y\in 
 \mathcal Y_{m'}^{j}$. On a alors, grâce à (\ref{ineg-dxYm'j-19dec}),    
 $$d(b,y)\leq d(x,b)-d(x,y)+4P\leq \frac{d(x,b)+d(x',b)-d(x,x')}{2}-M+3P-\de+\frac{N}{2}$$ et $$d(x',y)\geq d(x',b)-d(b,y)\geq \frac{d(x',b)-d(x,b)+d(x,x')}{2} +M+\de-3P-\frac{N}{2} $$ $$  \geq k'+\frac{M}{2}+\de-3P-\frac{3N}{2}\geq k'+3P$$  où l'avant-dernière inégalité a lieu par (\ref{ineg2-k'-troisiemecas-19dec}) et où, pour  la dernière, on a supposé 
 $ \frac{M}{2}+\de-3P-\frac{3N}{2}\geq 3P$, ce qui est permis par 
$(H_{M})$.

 Il reste donc à montrer que pour tout  $y\in 
 \mathcal Y_{m'}^{j}$ il existe $a\in S_{m'}$ et  
 $\tilde x'\in B(x',k')$ tels que  $y\in 2P\tg(\tilde x',a)$. Soit $y\in 
 \mathcal Y_{m'}^{j}$. 
 On va distinguer deux cas.

 On suppose d'abord $m'<m$. Alors il existe $a\in S_{m'}$ et $w\in S_{m'+1}$ tels que $y\in P\tg(w,a)$. On a $w,a\in 2F\tg(u,b)$. D'autre part $$d(x,w)\leq d_{\max}(x,S_{m'+1})\leq d_{\max}(x,S_{m'})\leq d(x,a)+N.$$
   Le lemme~\ref{x-a,b-y-geod} 
   appliqué à $(x,b,a,w)$ au lieu de $(x,y,a,b)$ et $(2F,2F,N)$ au lieu de $(\alpha,\beta,\rho)$ 
      montre alors $w\in (2F+2N+\de)\tg(x,a)$. 
      
      \ifx\JPicScale\undefined\def\JPicScale{1}\fi
\unitlength \JPicScale mm
\begin{picture}(140,55)(16,22)
\linethickness{0.3mm}
\put(50,49.75){\line(0,1){0.5}}
\multiput(49.99,50.75)(0.01,-0.5){1}{\line(0,-1){0.5}}
\multiput(49.96,51.25)(0.02,-0.5){1}{\line(0,-1){0.5}}
\multiput(49.92,51.74)(0.04,-0.5){1}{\line(0,-1){0.5}}
\multiput(49.87,52.24)(0.05,-0.5){1}{\line(0,-1){0.5}}
\multiput(49.81,52.73)(0.06,-0.49){1}{\line(0,-1){0.49}}
\multiput(49.74,53.23)(0.07,-0.49){1}{\line(0,-1){0.49}}
\multiput(49.65,53.72)(0.09,-0.49){1}{\line(0,-1){0.49}}
\multiput(49.55,54.21)(0.1,-0.49){1}{\line(0,-1){0.49}}
\multiput(49.44,54.69)(0.11,-0.49){1}{\line(0,-1){0.49}}
\multiput(49.32,55.18)(0.12,-0.48){1}{\line(0,-1){0.48}}
\multiput(49.18,55.66)(0.14,-0.48){1}{\line(0,-1){0.48}}
\multiput(49.04,56.13)(0.15,-0.48){1}{\line(0,-1){0.48}}
\multiput(48.88,56.61)(0.16,-0.47){1}{\line(0,-1){0.47}}
\multiput(48.71,57.07)(0.17,-0.47){1}{\line(0,-1){0.47}}
\multiput(48.52,57.54)(0.09,-0.23){2}{\line(0,-1){0.23}}
\multiput(48.33,58)(0.1,-0.23){2}{\line(0,-1){0.23}}
\multiput(48.13,58.45)(0.1,-0.23){2}{\line(0,-1){0.23}}
\multiput(47.91,58.9)(0.11,-0.22){2}{\line(0,-1){0.22}}
\multiput(47.68,59.35)(0.11,-0.22){2}{\line(0,-1){0.22}}
\multiput(47.44,59.78)(0.12,-0.22){2}{\line(0,-1){0.22}}
\multiput(47.19,60.22)(0.12,-0.22){2}{\line(0,-1){0.22}}
\multiput(46.93,60.64)(0.13,-0.21){2}{\line(0,-1){0.21}}
\multiput(46.66,61.06)(0.14,-0.21){2}{\line(0,-1){0.21}}
\multiput(46.38,61.47)(0.14,-0.21){2}{\line(0,-1){0.21}}
\multiput(46.09,61.88)(0.15,-0.2){2}{\line(0,-1){0.2}}
\multiput(45.79,62.27)(0.1,-0.13){3}{\line(0,-1){0.13}}
\multiput(45.48,62.66)(0.1,-0.13){3}{\line(0,-1){0.13}}
\multiput(45.16,63.05)(0.11,-0.13){3}{\line(0,-1){0.13}}
\multiput(44.83,63.42)(0.11,-0.12){3}{\line(0,-1){0.12}}
\multiput(44.49,63.79)(0.11,-0.12){3}{\line(0,-1){0.12}}
\multiput(44.14,64.14)(0.12,-0.12){3}{\line(0,-1){0.12}}
\multiput(43.79,64.49)(0.12,-0.12){3}{\line(1,0){0.12}}
\multiput(43.42,64.83)(0.12,-0.11){3}{\line(1,0){0.12}}
\multiput(43.05,65.16)(0.12,-0.11){3}{\line(1,0){0.12}}
\multiput(42.66,65.48)(0.13,-0.11){3}{\line(1,0){0.13}}
\multiput(42.27,65.79)(0.13,-0.1){3}{\line(1,0){0.13}}
\multiput(41.88,66.09)(0.13,-0.1){3}{\line(1,0){0.13}}
\multiput(41.47,66.38)(0.2,-0.15){2}{\line(1,0){0.2}}
\multiput(41.06,66.66)(0.21,-0.14){2}{\line(1,0){0.21}}
\multiput(40.64,66.93)(0.21,-0.14){2}{\line(1,0){0.21}}
\multiput(40.22,67.19)(0.21,-0.13){2}{\line(1,0){0.21}}
\multiput(39.78,67.44)(0.22,-0.12){2}{\line(1,0){0.22}}
\multiput(39.35,67.68)(0.22,-0.12){2}{\line(1,0){0.22}}
\multiput(38.9,67.91)(0.22,-0.11){2}{\line(1,0){0.22}}
\multiput(38.45,68.13)(0.22,-0.11){2}{\line(1,0){0.22}}
\multiput(38,68.33)(0.23,-0.1){2}{\line(1,0){0.23}}
\multiput(37.54,68.52)(0.23,-0.1){2}{\line(1,0){0.23}}
\multiput(37.07,68.71)(0.23,-0.09){2}{\line(1,0){0.23}}
\multiput(36.61,68.88)(0.47,-0.17){1}{\line(1,0){0.47}}
\multiput(36.13,69.04)(0.47,-0.16){1}{\line(1,0){0.47}}
\multiput(35.66,69.18)(0.48,-0.15){1}{\line(1,0){0.48}}
\multiput(35.18,69.32)(0.48,-0.14){1}{\line(1,0){0.48}}
\multiput(34.69,69.44)(0.48,-0.12){1}{\line(1,0){0.48}}
\multiput(34.21,69.55)(0.49,-0.11){1}{\line(1,0){0.49}}
\multiput(33.72,69.65)(0.49,-0.1){1}{\line(1,0){0.49}}
\multiput(33.23,69.74)(0.49,-0.09){1}{\line(1,0){0.49}}
\multiput(32.73,69.81)(0.49,-0.07){1}{\line(1,0){0.49}}
\multiput(32.24,69.87)(0.49,-0.06){1}{\line(1,0){0.49}}
\multiput(31.74,69.92)(0.5,-0.05){1}{\line(1,0){0.5}}
\multiput(31.25,69.96)(0.5,-0.04){1}{\line(1,0){0.5}}
\multiput(30.75,69.99)(0.5,-0.02){1}{\line(1,0){0.5}}
\multiput(30.25,70)(0.5,-0.01){1}{\line(1,0){0.5}}
\put(29.75,70){\line(1,0){0.5}}
\multiput(29.25,69.99)(0.5,0.01){1}{\line(1,0){0.5}}
\multiput(28.75,69.96)(0.5,0.02){1}{\line(1,0){0.5}}
\multiput(28.26,69.92)(0.5,0.04){1}{\line(1,0){0.5}}
\multiput(27.76,69.87)(0.5,0.05){1}{\line(1,0){0.5}}
\multiput(27.27,69.81)(0.49,0.06){1}{\line(1,0){0.49}}
\multiput(26.77,69.74)(0.49,0.07){1}{\line(1,0){0.49}}
\multiput(26.28,69.65)(0.49,0.09){1}{\line(1,0){0.49}}
\multiput(25.79,69.55)(0.49,0.1){1}{\line(1,0){0.49}}
\multiput(25.31,69.44)(0.49,0.11){1}{\line(1,0){0.49}}
\multiput(24.82,69.32)(0.48,0.12){1}{\line(1,0){0.48}}
\multiput(24.34,69.18)(0.48,0.14){1}{\line(1,0){0.48}}
\multiput(23.87,69.04)(0.48,0.15){1}{\line(1,0){0.48}}
\multiput(23.39,68.88)(0.47,0.16){1}{\line(1,0){0.47}}
\multiput(22.93,68.71)(0.47,0.17){1}{\line(1,0){0.47}}
\multiput(22.46,68.52)(0.23,0.09){2}{\line(1,0){0.23}}
\multiput(22,68.33)(0.23,0.1){2}{\line(1,0){0.23}}
\multiput(21.55,68.13)(0.23,0.1){2}{\line(1,0){0.23}}
\multiput(21.1,67.91)(0.22,0.11){2}{\line(1,0){0.22}}
\multiput(20.65,67.68)(0.22,0.11){2}{\line(1,0){0.22}}
\multiput(20.22,67.44)(0.22,0.12){2}{\line(1,0){0.22}}
\multiput(19.78,67.19)(0.22,0.12){2}{\line(1,0){0.22}}
\multiput(19.36,66.93)(0.21,0.13){2}{\line(1,0){0.21}}
\multiput(18.94,66.66)(0.21,0.14){2}{\line(1,0){0.21}}
\multiput(18.53,66.38)(0.21,0.14){2}{\line(1,0){0.21}}
\multiput(18.12,66.09)(0.2,0.15){2}{\line(1,0){0.2}}
\multiput(17.73,65.79)(0.13,0.1){3}{\line(1,0){0.13}}
\multiput(17.34,65.48)(0.13,0.1){3}{\line(1,0){0.13}}
\multiput(16.95,65.16)(0.13,0.11){3}{\line(1,0){0.13}}
\multiput(16.58,64.83)(0.12,0.11){3}{\line(1,0){0.12}}
\multiput(16.21,64.49)(0.12,0.11){3}{\line(1,0){0.12}}
\multiput(15.86,64.14)(0.12,0.12){3}{\line(1,0){0.12}}
\multiput(15.51,63.79)(0.12,0.12){3}{\line(0,1){0.12}}
\multiput(15.17,63.42)(0.11,0.12){3}{\line(0,1){0.12}}
\multiput(14.84,63.05)(0.11,0.12){3}{\line(0,1){0.12}}
\multiput(14.52,62.66)(0.11,0.13){3}{\line(0,1){0.13}}
\multiput(14.21,62.27)(0.1,0.13){3}{\line(0,1){0.13}}
\multiput(13.91,61.88)(0.1,0.13){3}{\line(0,1){0.13}}
\multiput(13.62,61.47)(0.15,0.2){2}{\line(0,1){0.2}}
\multiput(13.34,61.06)(0.14,0.21){2}{\line(0,1){0.21}}
\multiput(13.07,60.64)(0.14,0.21){2}{\line(0,1){0.21}}
\multiput(12.81,60.22)(0.13,0.21){2}{\line(0,1){0.21}}
\multiput(12.56,59.78)(0.12,0.22){2}{\line(0,1){0.22}}
\multiput(12.32,59.35)(0.12,0.22){2}{\line(0,1){0.22}}
\multiput(12.09,58.9)(0.11,0.22){2}{\line(0,1){0.22}}
\multiput(11.87,58.45)(0.11,0.22){2}{\line(0,1){0.22}}
\multiput(11.67,58)(0.1,0.23){2}{\line(0,1){0.23}}
\multiput(11.48,57.54)(0.1,0.23){2}{\line(0,1){0.23}}
\multiput(11.29,57.07)(0.09,0.23){2}{\line(0,1){0.23}}
\multiput(11.12,56.61)(0.17,0.47){1}{\line(0,1){0.47}}
\multiput(10.96,56.13)(0.16,0.47){1}{\line(0,1){0.47}}
\multiput(10.82,55.66)(0.15,0.48){1}{\line(0,1){0.48}}
\multiput(10.68,55.18)(0.14,0.48){1}{\line(0,1){0.48}}
\multiput(10.56,54.69)(0.12,0.48){1}{\line(0,1){0.48}}
\multiput(10.45,54.21)(0.11,0.49){1}{\line(0,1){0.49}}
\multiput(10.35,53.72)(0.1,0.49){1}{\line(0,1){0.49}}
\multiput(10.26,53.23)(0.09,0.49){1}{\line(0,1){0.49}}
\multiput(10.19,52.73)(0.07,0.49){1}{\line(0,1){0.49}}
\multiput(10.13,52.24)(0.06,0.49){1}{\line(0,1){0.49}}
\multiput(10.08,51.74)(0.05,0.5){1}{\line(0,1){0.5}}
\multiput(10.04,51.25)(0.04,0.5){1}{\line(0,1){0.5}}
\multiput(10.01,50.75)(0.02,0.5){1}{\line(0,1){0.5}}
\multiput(10,50.25)(0.01,0.5){1}{\line(0,1){0.5}}
\put(10,49.75){\line(0,1){0.5}}
\multiput(10,49.75)(0.01,-0.5){1}{\line(0,-1){0.5}}
\multiput(10.01,49.25)(0.02,-0.5){1}{\line(0,-1){0.5}}
\multiput(10.04,48.75)(0.04,-0.5){1}{\line(0,-1){0.5}}
\multiput(10.08,48.26)(0.05,-0.5){1}{\line(0,-1){0.5}}
\multiput(10.13,47.76)(0.06,-0.49){1}{\line(0,-1){0.49}}
\multiput(10.19,47.27)(0.07,-0.49){1}{\line(0,-1){0.49}}
\multiput(10.26,46.77)(0.09,-0.49){1}{\line(0,-1){0.49}}
\multiput(10.35,46.28)(0.1,-0.49){1}{\line(0,-1){0.49}}
\multiput(10.45,45.79)(0.11,-0.49){1}{\line(0,-1){0.49}}
\multiput(10.56,45.31)(0.12,-0.48){1}{\line(0,-1){0.48}}
\multiput(10.68,44.82)(0.14,-0.48){1}{\line(0,-1){0.48}}
\multiput(10.82,44.34)(0.15,-0.48){1}{\line(0,-1){0.48}}
\multiput(10.96,43.87)(0.16,-0.47){1}{\line(0,-1){0.47}}
\multiput(11.12,43.39)(0.17,-0.47){1}{\line(0,-1){0.47}}
\multiput(11.29,42.93)(0.09,-0.23){2}{\line(0,-1){0.23}}
\multiput(11.48,42.46)(0.1,-0.23){2}{\line(0,-1){0.23}}
\multiput(11.67,42)(0.1,-0.23){2}{\line(0,-1){0.23}}
\multiput(11.87,41.55)(0.11,-0.22){2}{\line(0,-1){0.22}}
\multiput(12.09,41.1)(0.11,-0.22){2}{\line(0,-1){0.22}}
\multiput(12.32,40.65)(0.12,-0.22){2}{\line(0,-1){0.22}}
\multiput(12.56,40.22)(0.12,-0.22){2}{\line(0,-1){0.22}}
\multiput(12.81,39.78)(0.13,-0.21){2}{\line(0,-1){0.21}}
\multiput(13.07,39.36)(0.14,-0.21){2}{\line(0,-1){0.21}}
\multiput(13.34,38.94)(0.14,-0.21){2}{\line(0,-1){0.21}}
\multiput(13.62,38.53)(0.15,-0.2){2}{\line(0,-1){0.2}}
\multiput(13.91,38.12)(0.1,-0.13){3}{\line(0,-1){0.13}}
\multiput(14.21,37.73)(0.1,-0.13){3}{\line(0,-1){0.13}}
\multiput(14.52,37.34)(0.11,-0.13){3}{\line(0,-1){0.13}}
\multiput(14.84,36.95)(0.11,-0.12){3}{\line(0,-1){0.12}}
\multiput(15.17,36.58)(0.11,-0.12){3}{\line(0,-1){0.12}}
\multiput(15.51,36.21)(0.12,-0.12){3}{\line(0,-1){0.12}}
\multiput(15.86,35.86)(0.12,-0.12){3}{\line(1,0){0.12}}
\multiput(16.21,35.51)(0.12,-0.11){3}{\line(1,0){0.12}}
\multiput(16.58,35.17)(0.12,-0.11){3}{\line(1,0){0.12}}
\multiput(16.95,34.84)(0.13,-0.11){3}{\line(1,0){0.13}}
\multiput(17.34,34.52)(0.13,-0.1){3}{\line(1,0){0.13}}
\multiput(17.73,34.21)(0.13,-0.1){3}{\line(1,0){0.13}}
\multiput(18.12,33.91)(0.2,-0.15){2}{\line(1,0){0.2}}
\multiput(18.53,33.62)(0.21,-0.14){2}{\line(1,0){0.21}}
\multiput(18.94,33.34)(0.21,-0.14){2}{\line(1,0){0.21}}
\multiput(19.36,33.07)(0.21,-0.13){2}{\line(1,0){0.21}}
\multiput(19.78,32.81)(0.22,-0.12){2}{\line(1,0){0.22}}
\multiput(20.22,32.56)(0.22,-0.12){2}{\line(1,0){0.22}}
\multiput(20.65,32.32)(0.22,-0.11){2}{\line(1,0){0.22}}
\multiput(21.1,32.09)(0.22,-0.11){2}{\line(1,0){0.22}}
\multiput(21.55,31.87)(0.23,-0.1){2}{\line(1,0){0.23}}
\multiput(22,31.67)(0.23,-0.1){2}{\line(1,0){0.23}}
\multiput(22.46,31.48)(0.23,-0.09){2}{\line(1,0){0.23}}
\multiput(22.93,31.29)(0.47,-0.17){1}{\line(1,0){0.47}}
\multiput(23.39,31.12)(0.47,-0.16){1}{\line(1,0){0.47}}
\multiput(23.87,30.96)(0.48,-0.15){1}{\line(1,0){0.48}}
\multiput(24.34,30.82)(0.48,-0.14){1}{\line(1,0){0.48}}
\multiput(24.82,30.68)(0.48,-0.12){1}{\line(1,0){0.48}}
\multiput(25.31,30.56)(0.49,-0.11){1}{\line(1,0){0.49}}
\multiput(25.79,30.45)(0.49,-0.1){1}{\line(1,0){0.49}}
\multiput(26.28,30.35)(0.49,-0.09){1}{\line(1,0){0.49}}
\multiput(26.77,30.26)(0.49,-0.07){1}{\line(1,0){0.49}}
\multiput(27.27,30.19)(0.49,-0.06){1}{\line(1,0){0.49}}
\multiput(27.76,30.13)(0.5,-0.05){1}{\line(1,0){0.5}}
\multiput(28.26,30.08)(0.5,-0.04){1}{\line(1,0){0.5}}
\multiput(28.75,30.04)(0.5,-0.02){1}{\line(1,0){0.5}}
\multiput(29.25,30.01)(0.5,-0.01){1}{\line(1,0){0.5}}
\put(29.75,30){\line(1,0){0.5}}
\multiput(30.25,30)(0.5,0.01){1}{\line(1,0){0.5}}
\multiput(30.75,30.01)(0.5,0.02){1}{\line(1,0){0.5}}
\multiput(31.25,30.04)(0.5,0.04){1}{\line(1,0){0.5}}
\multiput(31.74,30.08)(0.5,0.05){1}{\line(1,0){0.5}}
\multiput(32.24,30.13)(0.49,0.06){1}{\line(1,0){0.49}}
\multiput(32.73,30.19)(0.49,0.07){1}{\line(1,0){0.49}}
\multiput(33.23,30.26)(0.49,0.09){1}{\line(1,0){0.49}}
\multiput(33.72,30.35)(0.49,0.1){1}{\line(1,0){0.49}}
\multiput(34.21,30.45)(0.49,0.11){1}{\line(1,0){0.49}}
\multiput(34.69,30.56)(0.48,0.12){1}{\line(1,0){0.48}}
\multiput(35.18,30.68)(0.48,0.14){1}{\line(1,0){0.48}}
\multiput(35.66,30.82)(0.48,0.15){1}{\line(1,0){0.48}}
\multiput(36.13,30.96)(0.47,0.16){1}{\line(1,0){0.47}}
\multiput(36.61,31.12)(0.47,0.17){1}{\line(1,0){0.47}}
\multiput(37.07,31.29)(0.23,0.09){2}{\line(1,0){0.23}}
\multiput(37.54,31.48)(0.23,0.1){2}{\line(1,0){0.23}}
\multiput(38,31.67)(0.23,0.1){2}{\line(1,0){0.23}}
\multiput(38.45,31.87)(0.22,0.11){2}{\line(1,0){0.22}}
\multiput(38.9,32.09)(0.22,0.11){2}{\line(1,0){0.22}}
\multiput(39.35,32.32)(0.22,0.12){2}{\line(1,0){0.22}}
\multiput(39.78,32.56)(0.22,0.12){2}{\line(1,0){0.22}}
\multiput(40.22,32.81)(0.21,0.13){2}{\line(1,0){0.21}}
\multiput(40.64,33.07)(0.21,0.14){2}{\line(1,0){0.21}}
\multiput(41.06,33.34)(0.21,0.14){2}{\line(1,0){0.21}}
\multiput(41.47,33.62)(0.2,0.15){2}{\line(1,0){0.2}}
\multiput(41.88,33.91)(0.13,0.1){3}{\line(1,0){0.13}}
\multiput(42.27,34.21)(0.13,0.1){3}{\line(1,0){0.13}}
\multiput(42.66,34.52)(0.13,0.11){3}{\line(1,0){0.13}}
\multiput(43.05,34.84)(0.12,0.11){3}{\line(1,0){0.12}}
\multiput(43.42,35.17)(0.12,0.11){3}{\line(1,0){0.12}}
\multiput(43.79,35.51)(0.12,0.12){3}{\line(1,0){0.12}}
\multiput(44.14,35.86)(0.12,0.12){3}{\line(0,1){0.12}}
\multiput(44.49,36.21)(0.11,0.12){3}{\line(0,1){0.12}}
\multiput(44.83,36.58)(0.11,0.12){3}{\line(0,1){0.12}}
\multiput(45.16,36.95)(0.11,0.13){3}{\line(0,1){0.13}}
\multiput(45.48,37.34)(0.1,0.13){3}{\line(0,1){0.13}}
\multiput(45.79,37.73)(0.1,0.13){3}{\line(0,1){0.13}}
\multiput(46.09,38.12)(0.15,0.2){2}{\line(0,1){0.2}}
\multiput(46.38,38.53)(0.14,0.21){2}{\line(0,1){0.21}}
\multiput(46.66,38.94)(0.14,0.21){2}{\line(0,1){0.21}}
\multiput(46.93,39.36)(0.13,0.21){2}{\line(0,1){0.21}}
\multiput(47.19,39.78)(0.12,0.22){2}{\line(0,1){0.22}}
\multiput(47.44,40.22)(0.12,0.22){2}{\line(0,1){0.22}}
\multiput(47.68,40.65)(0.11,0.22){2}{\line(0,1){0.22}}
\multiput(47.91,41.1)(0.11,0.22){2}{\line(0,1){0.22}}
\multiput(48.13,41.55)(0.1,0.23){2}{\line(0,1){0.23}}
\multiput(48.33,42)(0.1,0.23){2}{\line(0,1){0.23}}
\multiput(48.52,42.46)(0.09,0.23){2}{\line(0,1){0.23}}
\multiput(48.71,42.93)(0.17,0.47){1}{\line(0,1){0.47}}
\multiput(48.88,43.39)(0.16,0.47){1}{\line(0,1){0.47}}
\multiput(49.04,43.87)(0.15,0.48){1}{\line(0,1){0.48}}
\multiput(49.18,44.34)(0.14,0.48){1}{\line(0,1){0.48}}
\multiput(49.32,44.82)(0.12,0.48){1}{\line(0,1){0.48}}
\multiput(49.44,45.31)(0.11,0.49){1}{\line(0,1){0.49}}
\multiput(49.55,45.79)(0.1,0.49){1}{\line(0,1){0.49}}
\multiput(49.65,46.28)(0.09,0.49){1}{\line(0,1){0.49}}
\multiput(49.74,46.77)(0.07,0.49){1}{\line(0,1){0.49}}
\multiput(49.81,47.27)(0.06,0.49){1}{\line(0,1){0.49}}
\multiput(49.87,47.76)(0.05,0.5){1}{\line(0,1){0.5}}
\multiput(49.92,48.26)(0.04,0.5){1}{\line(0,1){0.5}}
\multiput(49.96,48.75)(0.02,0.5){1}{\line(0,1){0.5}}
\multiput(49.99,49.25)(0.01,0.5){1}{\line(0,1){0.5}}

\linethickness{0.3mm}
\put(30,50){\line(1,0){120}}
\linethickness{0.3mm}
\multiput(130,60)(0.24,-0.12){83}{\line(1,0){0.24}}
\linethickness{0.3mm}
\multiput(100,70)(0.36,-0.12){83}{\line(1,0){0.36}}
\linethickness{0.3mm}
\multiput(60,60)(0.48,0.12){83}{\line(1,0){0.48}}
\linethickness{0.3mm}
\multiput(50,50)(0.12,0.12){83}{\line(1,0){0.12}}
\linethickness{0.3mm}
\put(60,60){\line(1,0){70}}
\linethickness{0.3mm}
\multiput(50,50)(0.96,0.12){83}{\line(1,0){0.96}}
\linethickness{0.3mm}
\multiput(60,60)(1.08,-0.12){83}{\line(1,0){1.08}}
\linethickness{0.3mm}
\put(80,30){\line(0,1){20}}
\put(30,55){\makebox(0,0)[cc]{$x$}}

\put(53,46){\makebox(0,0)[cc]{$u$}}

\put(85,30){\makebox(0,0)[cc]{$x'$}}

\put(57,63){\makebox(0,0)[cc]{$w$}}

\put(100,73){\makebox(0,0)[cc]{$y$}}

\put(131,63){\makebox(0,0)[cc]{$a$}}

\put(150,45){\makebox(0,0)[cc]{$b$}}

\put(84,46){\makebox(0,0)[cc]{$v$}}

\end{picture}

      Soit $v\in \geod(x,b)$ vérifiant $d(x,v)=E(\frac{d(x,b)+d(x,x')-d(x',b)}{2})$.  Par $(H_{\de}^{0}(x',x,v,b))$ on a  $$v\in \geod(x,b)\cap (\de+1)\tg(x,x')\cap (\de+1)\tg(x',b)$$ et 
      $$d(x',v)\leq \frac{d(x',b)-d(x,b)+d(x,x')}{2}+\de+1\leq k'$$ où la dernière égalité a lieu grâce à (\ref{ineg1-k'-troisiemecas-19dec})  
      car on  suppose   $N+\de+1\leq \frac{M}{2}$, ce qui est permis  par $(H_{M})$. On a  donc $v\in B(x',k')$.

      Comme $y\in P\tg(w,a)$ et $w\in (2F+2N+\de)\tg(x,a)$ on a $y\in (2F+2N+\de+P)\tg(x,a)$.       
         Par    (\ref{ineg-dxaa-19dec}) on a $d(x,a)\geq d(x,v)+M-2N$, et comme  $a\in 2F\tg(x,b)$ on en déduit  $d(b,a)-d(b,v)\leq -M+2N+2F<-\de$ car on suppose $-M+2N+2F<-\de$ grâce à  $(H_{M})$. Alors $(H_{\de}^{0}(a,x,v,b))$ donne       $v\in \de\tg(x,a)$. Grâce à (\ref{ineg-dxYm'j-19dec}) on a 
   \begin{gather}\nonumber d(x,y)\geq \frac{d(x,b)-d(x',b)+d(x,x')}{2}+M+P+\de-\frac{N}{2}\\ \label{ineg-11j1835}\geq d(x,v)-\frac{N}{2}+M+P+\de. \end{gather}
   Par $(H_{\de}^{2F+2N+\de+P}(v,x,y,a))$ on a 
   $$d(v,y)\leq \max(d(v,x)-d(x,y),d(v,a)-d(a,y))+2F+2N+2\de+P.$$
   Or $d(v,x)-d(x,y)+2F+2N+2\de+P<0$ par (\ref{ineg-11j1835}) et car on suppose 
   $2F+2N+2\de+P< -\frac{N}{2}+M+P+\de$, ce qui est permis par $(H_{M})$. 
   D'où 
   %
%   Comme $y\in (2F+2N+\de+P)\tg(x,a)$  on en déduit 
%     $$d(a,y)\leq d(a,v)+2F+\frac{5N}{2}-M$$ et le 
%      lemme~\ref{x-a,b-y-geod} 
%   appliqué à $(a,x,v,y)$ au lieu de $(x,y,a,b)$ et $(\de,2F+2N+\de+P,2F+\frac{5N}{2}-M)$ au lieu de $(\alpha,\beta,\rho)$   montre
       $y\in (2F+2N+2\de+P)\tg(v,a)$.
%        car on suppose $\de+2(2F+\frac{5N}{2}-M)\leq 2F+2N+\de+P$, ce qui est permis par  $(H_{M})$. 
        Grâce à $(H_{P})$ on  suppose $P\geq 2F+2N+2\de$, d'où $y\in 2P\tg(v,a)$. 
   En prenant $\tilde x'=v$ on a fini.

On suppose maintenant $m'=m$. Il existe  $a\in S_{m'}$ et $\tilde x\in B(x,k)$ tels que $y\in 2P\tg(\tilde x,a)$. On veut montrer qu'il existe 
 $\tilde x'\in B(x',k')$ tel que $y\in 2P\tg(\tilde x',a)$. 
 Grâce à   (\ref{ineg-dxYm'j-19dec})  on peut appliquer le sous-lemme~\ref{souslem2}   avec $\alpha
 =2P$ et $\beta=M+P+\de-\frac{N}{2}$. On suppose $F+\frac{3\de+\alpha}{2}+\de+N\leq \beta$, ce qui est permis par $(H_{M})$. Donc les hypothèses du sous-lemme~\ref{souslem2} sont satisfaites et l'existence  de $\tilde x'$ est démontrée.

 Il reste à montrer d). 
 Soit $v\in \geod(x,b)$ vérifiant $$d(x,v)=E(\frac{d(x,b)+d(x,x')-d(x',b)}{2}).$$ Par $(H_{\de}^{0}(x',x,v,b))$ on a $v\in \geod(x,b)\cap (\de+1)\tg(x,x')\cap (\de+1)\tg(x',b)$ et $$d(x',v)\leq \frac{d(x',b)-d(x,b)+d(x,x')}{2}+\de+1.$$
 Il résulte alors de (\ref{ineg1-k'-troisiemecas-19dec}) que $k'\geq d(x',v)+\frac{M}{2}-(N+\de+1)$,  donc   $$B(v,\frac{5M}{2}-(N+3\de+1))
 \subset B(x',k'+2M-2\de). $$ 
 Le premier ensemble de d) est inclus dans $$B(x,\frac{d(x,x')+r_{0}(Z)-r_{0}'(Z)}{2}+2M+2P+\de)\subset B(x,d(x,v)+2M+2P+N+\de).$$ 
 
 Si 
 \begin{gather}\label{cas-r0Z-petit-deux-ens}
 r_{0}(Z)\leq \frac{d(x,x')+r_{0}(Z)-r_{0}'(Z)}{2}+M\end{gather}
 on a $k'\geq \frac{r_{0}'(Z)-r_{0}(Z)+d(x,x')}{2}+\frac{M}{2}-N\geq 
 r_{0}'(Z) -\frac{M}{2}-N$, d'où $B(S_{0},M)\subset B(x',k'+2M)$. 
 
 Grâce au d) du lemme~\ref{lemme-S0-...Sm} le deuxième 
 ensemble de d)  (privé dans le cas où \eqref{cas-r0Z-petit-deux-ens} a lieu, de $B(S_{0},M)$ que l'on a déjà traité dans ce cas), 
   est inclus dans $$\{z'\in (2M+4P)\tg(x,b), d(x,z')\geq  
  \frac{d(x,x')+r_{0}(Z)-r_{0}'(Z)}{2}-N-4P\}.$$ Comme 
  $ \frac{d(x,x')+r_{0}(Z)-r_{0}'(Z)}{2}-N-4P\geq d(x,v)-2N-4P$, il suffit de montrer que pour \begin{gather*}z\in B(x,d(x,v)+2M+2P+N+\de)\\ \text{et \ \ }z'\in (2M+4P)\tg(x,b)\text{\ \  vérifiant \ \ }
  d(x,z')\geq d(x,v)-2N-4P,\end{gather*} $\geod(z,z')$ rencontre $B(v,\frac{5M}{2}-(N+3\de+1))$. 
 Par $(H_{\de}^{0}(z,x,v,b))$, on a  $$v\in \de\tg(x,z)\text{\ \ \  ou\ \ \  } v\in \de\tg(z,b).$$ 
 Si $v\in \de\tg(x,z)$,  $d(v,z)\leq 2M+2P+N+2\de$, donc $z\in \geod(z,z')$ appartient à $B(v,\frac{5M}{2}-(N+3\de+1))$ car on  suppose 
 $\frac{M}{2}\geq (2P+N+2\de)+(N+3\de+1)$ 
 grâce à $(H_{M})$. 
 Supposons donc  $ v\in \de\tg(z,b)$. 
 
\ifx\JPicScale\undefined\def\JPicScale{1}\fi
\unitlength \JPicScale mm
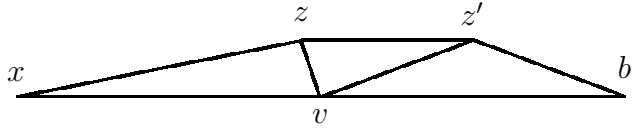
\begin{picture}(100,23)(20,24)
\linethickness{0.3mm}
\put(40,30){\line(1,0){80}}
\linethickness{0.3mm}
\multiput(80,30)(0.32,0.12){63}{\line(1,0){0.32}}
\linethickness{0.3mm}
\multiput(100,37.5)(0.32,-0.12){63}{\line(1,0){0.32}}
\linethickness{0.3mm}
\multiput(77.5,37.5)(0.12,-0.36){21}{\line(0,-1){0.36}}
\linethickness{0.3mm}
\multiput(40,30)(0.6,0.12){63}{\line(1,0){0.6}}
\linethickness{0.3mm}
\put(77.5,37.5){\line(1,0){22.5}}
\put(40,33){\makebox(0,0)[cc]{$x$}}

\put(77.5,41){\makebox(0,0)[cc]{$z$}}

\put(100,41){\makebox(0,0)[cc]{$z'$}}

\put(80,27.5){\makebox(0,0)[cc]{$v$}}

\put(120,34){\makebox(0,0)[cc]{$b$}}

\end{picture}
 
 Comme $$z'\in (2M+4P)\tg(x,b)\text{\ \  et \ \ }d(x,z')\geq d(x,v)-4P-2N,$$ on a $$d(b,z')\leq d(b,v)+2M+8P+2N. $$ Le lemme~\ref{x-a,b-y-geod} 
   appliqué à $(b,x,v,z')$ au lieu de $(x,y,a,b)$ et $(0,2M+4P, 2M+8P+2N)$ au lieu de $(\alpha,\beta,\rho)$ donne $z'\in (4M+16P+4N+\de)\tg(v,b)$. Comme $v\in \de\tg(z,b)$ le b) du lemme~\ref{geod-comp-xabc} montre $$v\in  
(4M+16P+4N+2\de)\tg(z,z').$$ Soit $t\in \geod(z,z')$ vérifiant   $d(z,t)=E(\frac{d(z,z')+d(z,v)-d(z',v)}{2})$.  Alors $(H_{\de}^{0}(v,z,t,z'))$ montre $d(v,t)\leq 2M+8P+2N+3\de+1\leq \frac{5M}{2}-(N+3\de+1)$ car on suppose $2M+8P+2N+3\de+1\leq \frac{5M}{2}-(N+3\de+1)$ grâce à $(H_{M})$. Cela termine la preuve de d) et donc l'étude du troisième cas. On a donc montré le lemme~\ref{oubli-x-xx'-x'}.  \cqfd

\begin{lem}\label{oubli-x-xx'-x'-Z}
Soit $$(a_{1},\dots,a_{p},S_{0},...,S_{m},(\mathcal Y_{i}^{j})_{i\in \{0,\dots,m\}, j\in \{1,\dots ,l_{i}\}})\in Y_{x}^{p,k,m,(l_{0},...,l_{m})}.$$
Dans les notations du lemme~\ref{oubli-x-xx'-x'}, les entiers $k',m',l'_{m'}$, la partie $J\subset \{1,...,l_{m'}\}$ et l'image $Z'$ de 
$$(a_{1},\dots,a_{p},S_{0},...,S_{m'},(\mathcal Y_{i}^{j})_{i\in \{0,\dots,m'-1\}, j\in \{1,\dots ,l_{i}\}}, (\mathcal Y_{m'}^{j_{\lambda}})_{ \lambda\in \{1,\dots ,l'_{m'}\}})$$ dans  
$\overline Y_{x'}^{p,k',m',(l_{0},...,l_{m'-1},l'_{m'})}$ par $\pi_{x'}^{p,k',m',(l_{0},...,l_{m'-1},l'_{m'})}$ ne dépendent que 
 de l'image $Z$ de $$(a_{1},\dots,a_{p},S_{0},...,S_{m},(\mathcal Y_{i}^{j})_{i\in \{0,\dots,m\}, j\in \{1,\dots ,l_{i}\}})$$ dans $\overline Y_{x,x',\star}^{p,k,m,(l_{0},...,l_{m})}$ par
$\pi_{x,x',\star}^{p,k,m,(l_{0},...,l_{m})}$. 

De plus il existe $L$ ne dépendant que de $Z$ tel que tout élément de $Z'$ a $L$ antécédents par l'application \begin{gather*}\text{de \ \ }(\pi_{x,x',\star}^{p,k,m,(l_{0},...,l_{m})})^{-1}(Z)\text{\ \ dans }(\pi_{x'}^{p,k',m',(l_{0},...,l_{m'-1},l'_{m'})})^{-1}(Z')\\  \text{qui \ à \ \ }
(a_{1},\dots,a_{p},S_{0},...,S_{m},(\mathcal Y_{i}^{j})_{i\in \{0,\dots,m\}, j\in \{1,\dots ,l_{i}\}})\text{\ \  associe}\\ 
(a_{1},\dots,a_{p},S_{0},...,S_{m'},(\mathcal Y_{i}^{j})_{i\in \{0,\dots,m'-1\}, j\in \{1,\dots ,l_{i}\}}, (\mathcal Y_{m'}^{j_{\lambda}})_{ \lambda\in \{1,\dots ,l'_{m'}\}})\end{gather*} et il existe $C=C(\de,N,K,Q,P,M)$ tel que l'on ait toujours 
\begin{gather}\label{ineg-L-24oct09}1\leq L\leq C^{(m+\sum_{i=m'}^{m}l_{i})-(m'+l'_{m'})}.\end{gather} 
\end{lem}
\noindent{\bf Démonstration.}
La première partie du lemme est évidente. L'existence de $L$ vient 
du e) du lemme~\ref{oubli-x-xx'-x'}, qui est une propriété très importante car elle ``découple'' les parties éliminées et les parties conservées en rendant superflue la connaissance des distances entre les points à distance $\leq M$ d'une partie éliminée et ceux à distance $\leq M$ d'une partie conservée.
De plus $L$ est égal au nombre de possibilités pour les parties éliminées (dont la liste est rappelée dans~\eqref{dem-460-liste-eliminees} ci-dessous) telles que les distances entre les points de $B(x,k+2M)$, $B(x',k'+2M)$ et ceux à distance $\leq M$ des parties éliminées prennent les valeurs prescrites par la donnée de $Z$. 
Il reste à montrer (\ref{ineg-L-24oct09}). 
 Les cardinaux des parties $S_{i}$ et $\mathcal Y_i^j$ sont bornés par une constante de la forme $C(\de,K,N,Q,P)$. Comme le nombre des parties éliminées est $(m+\sum_{i=m'}^{m}l_{i})-(m'+l'_{m'})$, 
il suffit  de montrer que  tout point $y$ 
d'une partie $S_{i}$ ou $\mathcal Y_i^j$ éliminée 
ne peut prendre que $C$ positions avec $C=C(\de,K,N,Q,P,M)$. Soit donc 
\begin{gather}
\label{dem-460-liste-eliminees}
y\in \bigcup _{ i\in \{m'+1,\dots ,m\}}
S_{i}
\cup  \bigcup _{i\in \{m'+1,\dots,m\}, j\in \{1,\dots ,l_{i}\}} \mathcal Y_{i}^{j}\cup \bigcup _{j\not\in J} \mathcal Y_{m'}^{j}.\end{gather}
Soit $b\in S_{0}$ et $u\in B(x,k)$ à distance minimale de $b$. 
On a alors $y\in 4P\tg(x,b)$ par le a) et le c) du lemme~\ref{lemme-S0-...Sm}, et $$d(x,y)\leq \max(k,\frac{d(x,x')+r_{0}(Z)-r_{0}'(Z)}{2})+M+2P+\de.$$ Si $d(x,y)\leq k+M+2P+\de$ on a $y\in B(x,k+2M)$ car on peut supposer 
$M\geq 2P+\de$ 
par $(H_{M})$ donc $y$ était déjà déterminé par la donnée de $Z$. Supposons donc $$d(x,y)\leq\frac{d(x,x')+r_{0}(Z)-r_{0}'(Z)}{2}+M+2P+\de.$$ Cela implique immédiatement \begin{gather}
\label{ineg-dxy-xx'xbx'b-20dec}d(x,y)\leq \frac{d(x,x')+d(x,b)-d(x',b)}{2}+M+2P+\de+N.\end{gather}    On a alors 
\begin{gather*}d(x',y)\leq \max (d(x',x)-d(b,x)+d(y,b), d(x',b)-d(b,x)+d(x,y))+\de\\ \leq \max (d(x',x)-d(x,y)+4P, d(x',b)-d(b,x)+d(x,y))+\de\\ \leq d(x',x)-d(x,y)+2M+4P+2N+3\de\end{gather*} 
où la première inégalité vient de $(H_{\de}(x',x,y,b))$, la deuxième utilise le fait que $y\in 4P\tg(x,b)$ et la dernière résulte de (\ref{ineg-dxy-xx'xbx'b-20dec}). 
On en déduit 
 $y\in (2M+4P+2N+3\de)\tg(x,x')$. 
 Comme $d(x,y)$ fait partie de  la donnée de $Z$, le lemme~\ref{cardinal-tranche-geod} montre que le nombre de possibilités pour $y$ est borné par une constante $C=C(\de,N,K,Q,P,M)$.  \cqfd

Dans les notations du lemme~\ref{oubli-x-xx'-x'-Z} on désigne par  
$$\theta : \bigcup_{k,m,l_{0},...,l_{m}}Y_{x,x',\star}^{p,k,m,(l_{0},...,l_{m})}\to  \bigcup_{k',m',l'_{0},...,l'_{m'}} Y_{x'}^{p,k',m',(l'_{0},...,l'_{m'})}$$ l'application qui \begin{gather*}\text{  à \ \ }
(a_{1},\dots,a_{p},S_{0},...,S_{m},(\mathcal Y_{i}^{j})_{i\in \{0,\dots,m\}, j\in \{1,\dots ,l_{i}\}})\ \ \text{ associe}\\ 
(a_{1},\dots,a_{p},S_{0},...,S_{m'},(\mathcal Y_{i}^{j})_{i\in \{0,\dots,m'-1\}, j\in \{1,\dots ,l_{i}\}}, (\mathcal Y_{m'}^{j_{\lambda}})_{ \lambda\in \{1,\dots ,l'_{m'}\}})\end{gather*} et on note 
 $$\overline  \theta : \bigcup_{k,m,l_{0},...,l_{m}}\overline Y_{x,x',\star}^{p,k,m,(l_{0},...,l_{m})}\to  \bigcup_{k',m',l'_{0},...,l'_{m'}}\overline Y_{x'}^{p,k',m',(l'_{0},...,l'_{m'})}$$ l'application induite par $\theta$, qui à $Z$ associe $Z'$ et dont l'existence résulte du lemme~\ref{oubli-x-xx'-x'-Z}. 
  
Pour $k',m',l_{0}',...,l'_{m'}\in \N$ et $Z'\in \overline Y_{x'}^{p,k',m',(l_{0}',...,l_{m'}')}$ on note $\xi_{Z'}$ la forme linéaire sur $\C^{(\Delta_{p})}$ définie par 
$$\xi_{Z'}(f)=\sum _{(a'_{1},\dots,a'_{p},S'_{0},...,S'_{m'},(\mathcal Y_{i}^{'j})_{i\in \{0,\dots,m'\}, j\in \{1,\dots ,l_{i}'\}}) \in 
(\pi_{x'}^{p,k',m',(l'_{0},...,l'_{m'})})^{-1}(Z')} f(a'_1,...,a'_p).$$
%
%
%$$W(Z')=\sharp \big((\pi_{x'}^{p,k',m',(l_{0}',...,l_{m'}')})^{-1}(Z')\big)^{-\alpha }
%$$ $$\bigg| \sum _{(a'_{1},\dots,a'_{p},S'_{1},...,S'_{m'},(\mathcal Y_{i}^{',j})_{i\in \{0,\dots,m'\}, j\in \{1,\dots ,l_{i}'\}}) \in 
%(\pi_{x'}^{p,k',m',(l'_{0},...,l'_{m})})^{-1}(Z')} f(a'_1,...,a'_p)\bigg|^{2}.$$
Grâce au lemme~\ref{oubli-x-xx'-x'-Z}, on a pour $Z'=\theta(Z)$ et $f\in \C^{(\Delta_{p})}$, 
\begin{gather}\nonumber \sharp \big((\pi_{x,x',\star}^{p,k,m,(l_{0},...,l_{m})})^{-1}(Z)\big)^{-\alpha }\big|\xi_{Z}(f)
%
%\sum _{(a_{1},\dots,a_{p},S_{1},...,S_{m},(\mathcal Y_{i}^{j})_{i\in \{0,\dots,m\}, j\in \{1,\dots ,l_{i}\}}) \in 
%(\pi_{x,x',\star}^{p,k,m,(l_{0},...,l_{m})})^{-1}(Z)} f(a_1,...,a_p)
\big|^{2} \\ \label{sommeZ=C.sommeZ'}
=\big(L^{2-\alpha} \big)\ \sharp \big((\pi_{x'}^{p,k',m',(l_{0}',...,l_{m'}')})^{-1}(Z')\big)^{-\alpha }\big|\xi_{Z'}(f)\big|^{2}
 \end{gather}
%\sharp \big((\pi_{x'}^{p,k',m',(l_{0}',...,l_{m'}')})^{-1}(Z')\big)^{-\alpha } $$\begin{gather}\label{sommeZ=C.sommeZ'}\sum _{(a'_{1},\dots,a'_{p},S'_{1},...,S'_{m'},(\mathcal Y_{i}^{',j})_{i\in \{0,\dots,m'\}, j\in \{1,\dots ,l_{i}'\}}) \in 
%(\pi_{x'}^{p,k',m',(l'_{0},...,l'_{m})})^{-1}(Z')} f(a'_1,...,a'_p).\end{gather}
avec $L\leq C^{(m+\sum_{i=m'}^{m}l_{i})-(m'+l'_{m'})}$, où $C$ est la constante du lemme~\ref{oubli-x-xx'-x'-Z}. Il en résulte que pour tout $f\in \C^{(\Delta_{p})}$, on a 
\begin{gather*}\|f\|_{\H_{x,x',\star,s}(\Delta_{p})}^{2}\leq \sum _{k',m',l'_{0},\dots ,l'_{m'}\in \N,Z'\in \overline Y_{x'}^{p,k',m',(l'_{0},...,l'_{m'})}}  
\\ \Bigg(\sum_{k,m,l_{0},...,l_{m}\in \N,Z\in \overline Y_{x,x',\star}^{p,k,m,(l_{0},...,l_{m})}\text{ tels que }\overline \theta(Z)=Z'}B^{-(m+\sum_{i=0}^{m}l_{i})}\\ (C^{2-\alpha})^{(m+\sum_{i=m'}^{m}l_{i})-(m'+l'_{m'})}
 e^{2s(r_{0}(Z)-k)}\Big(\prod_{i=0}^{m}s_{i}(Z)^{-l_{i}} \Big)\Bigg )\\ \sharp \big((\pi_{x'}^{p,k',m',(l_{0}',...,l_{m'}')})^{-1}(Z')\big)^{-\alpha }
 \big|\xi_{Z'}(f)\big|^{2}. \end{gather*}

D'autre part  
 \begin{gather*}\|f\|_{\H_{x',s}(\Delta_{p})}^{2}=\sum _{k',m',l'_{0},\dots ,l'_{m'}\in \N,Z'\in \overline Y_{x'}^{p,k',m',(l'_{0},...,l'_{m'})}} B^{-(m'+\sum_{i=0}^{m'}l'_{i})} e^{2s(r_{0}(Z')-k')}\\ \Big(\prod_{i=0}^{m'}s_{i}(Z')^{-l'_{i}} \Big)
 \sharp \big((\pi_{x'}^{p,k',m',(l_{0}',...,l_{m'}')})^{-1}(Z')\big)^{-\alpha }
\big|\xi_{Z'}(f)\big|^{2}\end{gather*}
où les entiers  $r_{0}(Z'), s_{0}(Z'),...,s_{m'}(Z')$ sont tels que pour tout 
$$(a'_{1},\dots,a'_{p},S'_{0},...,S'_{m'},(\mathcal Y_{i}^{'j})_{i\in \{0,\dots,m'\}, j\in \{1,\dots ,l'_{i}\}})\in (\pi_{x'}^{p,k',m',(l'_{0},...,l'_{m})})^{-1}(Z')$$
on a $r_{0}(Z')=d(x',S_{0}')$, $s_{i}(Z')=d(S'_{i},S'_{i+1})+2M$ pour $i\in \{0,...,m'-1\}$, et $s_{m'}(Z')=d(x',S'_{m'})-k'$. 
On remarque que pour   $Z'=\theta(Z)$ on a  $r_{0}(Z')=r_{0}'(Z)$ et pour $i\in \{0,...,m'-1\}$, $l'_{i}=l_{i}$ et $s_{i}(Z')=s_{i}(Z)$. 

\noindent{\bf Démonstration du lemme~\ref{deuxieme-ineg-equiv-norme} en admettant le lemme~\ref{somme-aZZ'}.} 
On pose $C_{0}=C^{2-\alpha}$ où $C$ est comme dans le lemme~\ref{oubli-x-xx'-x'-Z}. Donc $C_{0}\leq C(\de,N,K,Q,P,M)$.  
Pour montrer le lemme~\ref{deuxieme-ineg-equiv-norme}   on est ramené  à montrer le lemme suivant. \cqfd

\begin{lem}\label{somme-aZZ'}
 Il existe une constante $C=C(\de,K,N,P,Q,M,s,B)$ telle  que pour 
 tous $k',m',l_{0}',...,l'_{m'}$ et $Z'\in  \overline  Y_{x'}^{p,k',m',(l'_{0},...,l'_{m'})}$
on ait $$\sum_{k,m,l_{0},...,l_{m}\in \N,Z\in \overline Y_{x,x',\star}^{p,k,m,(l_{0},...,l_{m})}\text{ tels que }\overline \theta(Z)=Z'}(C_{0}B^{-1})^{m+\sum_{i=m'}^{m}l_{i}}
 e^{2s(r_{0}(Z)-k)}$$ $$\Big(\prod_{i=m'}^{m}s_{i}(Z)^{-l_{i}} \Big)
\leq Ce^{3sd(x,x')}(C_{0}B^{-1})^{m'+l'_{m'}}
 e^{2s(r_{0}(Z')-k')}s_{m'}(Z')^{-l'_{m'}}.$$
\end{lem}
\noindent{\bf Démonstration du lemme~\ref{somme-aZZ'} en admettant lemme~\ref{somme-aZZ'-bis}.} 
On rappelle que pour $Z'=\theta(Z)$, $r_{0}(Z')=r_{0}'(Z)$. 
Par (\ref{def-k'Z-24oct09}) on a $$k'-r_{0}'(Z)=\min(0,\max(k-r_{0}(Z), E(\frac{d(x,x')-r_{0}(Z)-r_{0}'(Z)}{2})))+\frac{M}{2}$$ d'où $$k'-r_{0}'(Z)\leq \max(k-r_{0}(Z), E(\frac{d(x,x')-r_{0}(Z)-r_{0}'(Z)}{2}))+\frac{M}{2}.$$ Comme $|r_{0}(Z)
-r_{0}'(Z)|\leq d(x,x')$ on a $$\frac{d(x,x')-r_{0}(Z)-r_{0}'(Z)}{2}\leq d(x,x')-r_{0}(Z)\leq d(x,x')+(k-r_{0}(Z))$$
d'où $(r_{0}(Z)-k)\leq (r_{0}'(Z)-k')+d(x,x')+\frac{M}{2}$. 
Les valeurs de $k$ pour lesquelles il existe $m,l_{0},...,l_{m}\in \N$ et $Z\in 
\overline Y_{x,x',\star}^{p,k,m,(l_{0},...,l_{m})}$ tels que  $\theta(Z)=Z'$ sont donc incluses dans un intervalle $[k_{0},+\infty[$ avec $k_{0}\in \N$ vérifiant 
$$(r_{0}(Z)-k_{0})\leq (r_{0}'(Z)-k')+d(x,x')+\frac{M}{2}. $$
Comme la série $1+e^{-2s}+e^{-4s}+...$ converge, 
pour montrer le lemme~\ref{somme-aZZ'} il suffit  d'établir le lemme suivant. \cqfd

\begin{lem}\label{somme-aZZ'-bis}
%Il existe une constante $C=C(\de,K,N,Q,P,M,s,B)$ telle que pour 
Pour tous $k,k',m',l_{0}',...,l'_{m'}$ et $Z'\in  \overline  Y_{x'}^{p,k',m',(l'_{0},...,l'_{m'})}$
on a \begin{gather*}\sum_{m,l_{0},...,l_{m}\in \N,Z\in \overline Y_{x,x',\star}^{p,k,m,(l_{0},...,l_{m})}\text{ tels que }\overline\theta(Z)=Z'}(C_{0}B^{-1})^{m+\sum_{i=m'}^{m}l_{i}}
 \Big(\prod_{i=m'}^{m}s_{i}(Z)^{-l_{i}} \Big)\\ 
\leq 4e^{sd(x,x')}(C_{0}B^{-1})^{m'+l'_{m'}}
 s_{m'}(Z')^{-l'_{m'}} 
.\end{gather*}
\end{lem}
 \noindent{\bf Démonstration du lemme~\ref{somme-aZZ'-bis} en admettant le lemme~\ref{somme-aZZ'-ter}.} 
D'après le lemme~\ref{oubli-x-xx'-x'-Z}, pour $Z'=\overline \theta(Z)$, $\theta$ induit une surjection $$\text{de \ \ }(\pi_{x,x',\star}^{p,k,m,(l_{0},...,l_{m})})^{-1}(Z)\text{\ \ 
dans \ \ }(\pi_{x'}^{p,k',m',(l'_{0},...,l'_{m'})})^{-1}(Z').$$ 
Pour montrer  le lemme~\ref{somme-aZZ'-bis} il suffit  donc d'établir le lemme suivant. \cqfd

\begin{lem}\label{somme-aZZ'-ter}
%Il existe une constante $C=C(\de,K,N,Q,P,M,s,B)$, telle  que pour 
Pour tous $k,k',m',l_{0}',...,l'_{m'}$ et 
$$(a'_{1},\dots,a'_{p},S'_{0},...,S'_{m'},(\mathcal Y_{i}^{'j})_{i\in \{0,\dots,m'\}, j\in \{1,\dots ,l'_{i}\}})\in Y_{x'}^{p,k',m',(l'_{0},...,l'_{m'})}$$
on a
 \begin{gather} \nonumber 
\sum_{m,l_{0},...,l_{m},(a_{1},\dots,a_{p},S_{0},...,S_{m},(\mathcal Y_{i}^{j})_{i\in \{0,\dots,m\}, j\in \{1,\dots ,l_{i}\}})}
(C_{0}B^{-1})^{(m-m')+(\sum_{i=m'}^{m}l_{i})-l'_{m'}}\\  \nonumber\Big(\prod_{i=m'}^{m-1}\big(d(S_{i},S_{i+1})+2M\big)^{-l_{i}} \Big)
 \big(d(x,S_{m})-k\big)^{-l_{m}} \big(d(x',S'_{m'})-k'\big)^{l'_{m'}}\\ \label{ineg-somme-antecedents-theta}
\leq 4e^{sd(x,x')}\end{gather}
où la somme porte sur les $m,l_{0},...,l_{m}\in \N$  et les 
$$(a_{1},\dots,a_{p},S_{0},...,S_{m},(\mathcal Y_{i}^{j})_{i\in \{0,\dots,m\}, j\in \{1,\dots ,l_{i}\}})\in Y_{x}^{p,k,m,(l_{0},...,l_{m})}$$
 tels que  $k',m',l_{0}',...,l'_{m'}$ et 
$$(a'_{1},\dots,a'_{p},S'_{0},...,S'_{m'},(\mathcal Y_{i}^{'j})_{i\in \{0,\dots,m'\}, j\in \{1,\dots ,l'_{i}\}})$$ 
soient associés à $$(a_{1},\dots,a_{p},S_{0},...,S_{m},(\mathcal Y_{i}^{j})_{i\in \{0,\dots,m\}, j\in \{1,\dots ,l_{i}\}})$$ comme dans le lemme~\ref{oubli-x-xx'-x'} (donc en particulier  $a'_{1}=a_{1}, ...,a'_{p}=a_{p},S'_{0}=S_{0},...,S'_{m'}=S_{m'}$, $k'=k'(Z)$, 
$l'_{i}=l_{i}$ pour $i\in \{0,...,m'-1\}$, 
 $\mathcal Y_{i}^{'j}=\mathcal Y_{i}^{j}$ pour $i\in \{0,...,m'-1\}$ et $j\in \{1,...,l_{i}\}$ et $\mathcal Y_{m'}^{'\lambda}=\mathcal Y_{m'}^{j_{\lambda}}$ pour $\lambda\in \{1,...,l'_{m'}\})$. 
\end{lem}
\noindent{\bf Démonstration du lemme~\ref{somme-aZZ'-ter} en admettant les lemmes~\ref{somme-aZZ'-ter-1}   et ~\ref{somme-aZZ'-ter-2}.} 
D'abord on suppose $l'_{0}=...=l'_{m'-1}=0$ car en supprimant les $\mathcal Y_{i}^{'j}$ pour $i\in \{0,...,m'-1\}$ et $j\in \{1,...,l'_{i}\}$ le nouvel énoncé est strictement équivalent à l'ancien. 
D'après le lemme~\ref{existence-C}
 il existe une constante $C_{1}=C(\de,K,N,Q,P)$ telle que, 
connaissant $S_{0},...,S_{m}$, pour $i\in \{m'+1,...,m\}$ et $j\in \{1,...,l_{i}\}$ le nombre de possibilités pour $\mathcal Y_{i}^{j}$ est borné par $C_{1}(d(S_{i},S_{i+1})+2M)$ si $i<m$ et  par $C_{1}(d(x,S_{m})-k)$ si $i=m$. Comme 
$$\sum_{l_{m'+1},...,l_{m}\in \N}(C_{1}C_{0}B^{-1})^{\sum_{i=m'+1}^{m}l_{i}}=(1-C_{1}C_{0}B^{-1})^{-(m-m')},$$ pour montrer le lemme~\ref{somme-aZZ'-ter} il suffit de  montrer les deux lemmes suivants, qui distinguent les cas $m'=m$ et $m'<m$. Dans les deux lemmes on note $l$ et $l'$ au lieu de $l_{m'}$ et $l'_{m'}$, puisque que les  $l_{i}$ et $l'_{i}$ pour $i\neq m'$ ont disparu. Dans les deux lemmes suivants on a supposé $J=\{1,...,l'\}$ et remplacé la somme sur les parties $J\subset \{1,...,l\}$ de cardinal $l'$  par la multiplication par $\binom{l}{l'}$. \cqfd

\begin{lem}\label{somme-aZZ'-ter-1}
%Il existe $C=C(\de,K,N,Q,P,M,s,B)$ telle que pour 
Pour tous $k,k',m',l'\in \N$ et $$(a_{1},\dots,a_{p},S_{0},...,S_{m'},(\mathcal Y_{m'}^{j})_{ j\in \{1,\dots ,l'\}})\in Y_{x'}^{p,k',m',(0,...,0,l')}$$
on a 
\begin{gather}\label{ineg-somme-aZZ'-ter-1}\sum_{l\geq l'}\binom{l}{l'}
\sum_{\mathcal Y_{m'}^{l'+1}, ...,\mathcal Y_{m'}^{l}  } (C_{0}B^{-1})^{l-l'}
 \big(d(x,S_{m'})-k\big)^{-l} 
 \big(d(x',S_{m'})-k'\big)^{l'}
\leq 2\end{gather}
où la somme porte sur les $\mathcal Y_{m'}^{l'+1}, ...,\mathcal Y_{m'}^{l} $ tels que 
$$(a_{1},\dots,a_{p},S_{0},...,S_{m'},(\mathcal Y_{m'}^{j})_{ j\in \{1,\dots ,l\}})\in Y_{x}^{p,k,m',(0,...,0,l)}$$
et, en notant $T=\frac{d(x,x')+d(x,S_{0})-d(x',S_{0})}{2}$ on ait 
\begin{gather*}d_{\max}(x,S_{m'})> \max(k,T)+M\text{\ \ \ si \ \ \ }m'>0,\\ d_{\max}(x,\mathcal Y_{m'}^{j})> \max(k,T)+M+2P+\de\text{\ \ \ pour\ \ \ }j\in \{1,...,l'\}\text{\ \ \ et} \\ d_{\max}(x,\mathcal Y_{m'}^{j})\leq \max(k,T)+M+2P+\de\text{\ \ \ pour\ \ \ }j\in \{l'+1,...,l\}.\end{gather*}
\end{lem}

On rappelle que les constantes $C_{0}$ et  $C_{1}$ 
(qui apparaissent avant les énoncés des lemmes~\ref{somme-aZZ'} et \ref{somme-aZZ'-ter-1})
sont majorées par  $C(\de,K,N,Q,P,M)$. 

\begin{lem}\label{somme-aZZ'-ter-2}
%Il existe $C=C(\de,K,N,Q,P,M,s,B)$ telle que pour
Pour  tous $k,k',m',l'\in \N$ et $$(a_{1},\dots,a_{p},S_{0},...,S_{m'},(\mathcal Y_{m'}^{j})_{ j\in \{1,\dots ,l'\}})\in Y_{x'}^{p,k',m',(0,...,0,l')}$$
on a 
\begin{gather}\nonumber \sum_{m>m', l\geq l'}     \binom{l}{l'} \Big(
\frac{C_{0}B^{-1}}{
1-C_{0}C_{1}B^{-1}}\Big)^{m-m'}
\sum_{S_{m'+1},...,S_{m},\mathcal Y_{m'}^{l'+1}, ...,\mathcal Y_{m'}^{l}  } (C_{0}B^{-1})^{l-l'}\\ \label{ineg-somme-aZZ'-ter-2}
 \big(d(S_{m'},S_{m'+1})+2M\big)^{-l} 
 \big(d(x',S_{m'})-k'\big)^{l'}
\leq 4e^{sd(x,x')}\end{gather}
où la deuxième somme porte sur les $S_{m'+1},...,S_{m},\mathcal Y_{m'}^{l'+1}, ...,\mathcal Y_{m'}^{l} $ tels que 
$$(a_{1},\dots,a_{p},S_{0},...,S_{m},(\mathcal Y_{m'}^{j})_{ j\in \{1,\dots ,l\}})\in Y_{x}^{p,k,m,(0,...,0,l,0,...,0)}$$
et, en notant $T=\frac{d(x,x')+d(x,S_{0})-d(x',S_{0})}{2}$ on ait 
\begin{gather*}d_{\max}(x,S_{m'})> \max(k,T)+M\text{\ \ \ si \ \ \ }m'>0,\\ d_{\max}(x,S_{m'+1})\leq \max(k,T)+M,
\\ d_{\max}(x,\mathcal Y_{m'}^{j})> \max(k,T)+M+2P+\de\text{\ \ \ pour\ \ \ }j\in \{1,...,l'\}\text{\ \ \ et} \\ d_{\max}(x,\mathcal Y_{m'}^{j})\leq \max(k,T)+M+2P+\de\text{\ \ \ pour\ \ \ }j\in \{l'+1,...,l\}.\end{gather*}
\end{lem}
Dans la démonstration des lemmes~\ref{somme-aZZ'-ter-1}   et~\ref{somme-aZZ'-ter-2} le petit calcul suivant servira plusieurs fois. 

\begin{lem}\label{calcul-A-B}
Soit $\mathcal A, \mathcal B\in \R_{+}$, $2\mathcal A+\mathcal B\leq 1$. 
Alors pour tout $l'\in \N$, 
$$\sum_{l\geq l'} \binom{l}{l'} \mathcal A^{l-l'}\mathcal B^{l'}\leq 2.$$
\end{lem}
\noindent{\bf Démonstration.} On a la formule générale, pour $x\in \C$ avec $|x|<1$ et $k\in \N$,  
\begin{gather}\label{form-binom-24oct09}\sum_{n\in \N} \binom{n+k}{k}x^{n}=(1-x)^{-(k+1)}.\end{gather} Donc 
$\sum_{l\geq l'} \binom{l}{l'} \mathcal A^{l-l'}\mathcal B^{l'}
=\frac{\mathcal B^{l'}}{(1-\mathcal A)^{l'+1} }\leq \frac{1}{1-\mathcal A} \leq 2$ puisque $\mathcal B\leq 1-\mathcal A$ et $1-\mathcal A \geq 1/2$. \cqfd

\noindent{\bf Démonstration du lemme~\ref{somme-aZZ'-ter-1}.} On va distinguer trois cas comme dans la démonstration du lemme~\ref{oubli-x-xx'-x'}, mais une partie de la démonstration est commune aux trois cas. On note $T=\frac{d(x,S_{0})-d(x',S_{0})+d(x,x')}{2}$ comme dans l'énoncé du lemme. 
On fixe $b\in S_{0}$. Soit $y\in S_{m'}$. 
Alors \begin{gather}\label{ineg-dxy-ter1-cas2}d(x,y)\geq \frac{d(x,b)+d(x,x')-d(x',b)}{2}-N.\end{gather} 
En effet cela est vrai si $m'=0$ car $d(x,y)\geq d(x,S_{0})\geq d(x,b)-N$ et $d(x,x')\leq d(x,b)+d(x',b)$, et cela est vrai si $m'>0$
 car  la condition $d_{\max}(x,S_{m'})> \max(k,T)+M$  implique $$d(x,y)\geq d_{\max}(x,S_{m'})-N\geq T+M-N\geq \frac{d(x,b)+d(x,x')-d(x',b)}{2}+M-2N$$ et on suppose $M\geq N$, ce qui est permis par 
 $(H_{M})$. 
On en déduit  \begin{gather}\label{ineg-commun-46-20dec}d(x',y)-d(x',b)+d(x,b)\leq d(x,y)+P+2N+\de.\end{gather} 
En effet  \begin{gather*}
 d(x',y)\leq \max(d(x',b)-d(x,b)+d(x,y), d(x',x)-d(x,b)+d(b,y))+\de \\
 \leq \max(d(x',b)-d(x,b)+d(x,y), d(x',x)-d(x,y)+P)+\de \\
 \leq d(x',b)-d(x,b)+d(x,y)+P+2N+\de\end{gather*}
 où la première inégalité a lieu par  $(H_{\de}(x',x,y,b))$, la deuxième inégalité utilise $ y\in P\tg(x,b)$ et la dernière inégalité résulte de (\ref{ineg-dxy-ter1-cas2}).

\noindent {\bf Premier cas.} On suppose $d(x,S_{0})\leq k$. 

Alors $m'=0$ et $l=0$ d'après le lemme~\ref{m-fini}  et la remarque qui suit la définition~\ref{defi-Y}. Donc la somme est vide ou réduite à un élément et l'inégalité est triviale. 

\noindent {\bf Deuxième cas.} On suppose $d(x,S_{0})> k$ et $T\leq k$. 

% et $u\in B(x,k)$ à distance minimale de $b$. 
Alors pour $j\in \{l'+1,...,l\}$ on a 
$$k+3P\leq d_{\max}(x,\mathcal Y_{m'}^{j})\leq k+M+2P+\de$$ 
où l'inégalité de gauche vient de la condition iv) de la définition~\ref{defi-Y}. 
Comme $
\mathcal Y_{m'}^{j}\subset 4P\tg(x,b)$ le lemme~\ref{cardinal-tranche-geod}
montre que le nombre de possibilités pour $\mathcal Y_{m'}^{j}$ est borné par $C_{2}=C(\de,K,N,Q,P,M)$. On a $$k'=d(x',S_{0})-d(x,S_{0})+k+\frac{M}{2}$$ par  
(\ref{form-r'-cas2-4j})
%(\ref{def-k'Z-24oct09}) 
donc \begin{gather}\label{ineg-k'-deux-cas-20dec}k'\geq d(x',b)-d(x,b)+k+\frac{M}{2}-N.\end{gather} 
Soit $y\in S_{m'}$. 
 On a  
 $$d(x',y)-k'\leq  d(x',y)-d(x',b)+d(x,b)+N-k-\frac{M}{2}\leq d(x,y)-k-\frac{M}{2}
+P+3N+\de $$
où la première inégalité vient de (\ref{ineg-k'-deux-cas-20dec}) et la deuxième de (\ref{ineg-commun-46-20dec}). 
Comme cela est vrai pour tout $y\in S_{m'}$ on  en déduit 
\begin{gather}\label{ineg-bis-cas2-x'Sm'k'-xSm'k}\big(d(x',S_{m'})-k'\big)\leq \big(d(x,S_{m'})-k\big)-\frac{M}{2}+P+3N+\de
\end{gather}
d'où \begin{gather}\label{ineg-diff1-deux-cas-20dec}\big(d(x',S_{m'})-k'\big) \leq \big(d(x,S_{m'})-k\big)-1\end{gather}
car on suppose $\frac{M}{2}\geq P+3N+\de+1$, ce qui est 
 permis 
par  $(H_{M})$. Le membre de gauche de (\ref{ineg-somme-aZZ'-ter-1}) est donc majoré par 
$$\sum_{l\geq l'}\binom{l}{l'}
 \Big(\frac{C_{0}C_{2}B^{-1}}{d(x,S_{m'})-k}\Big)^{(l-l')}
 \Big(\frac{d(x',S_{m'})-k'}{d(x,S_{m'})-k}\Big)^{l'}\leq 2
 $$
par le lemme~\ref{calcul-A-B}. Les hypothèses  du lemme~\ref{calcul-A-B}
 sont satisfaites car on suppose $2C_{0}C_{2}B^{-1}\leq1$ grâce à  $(H_{B})$ ce qui implique 
 $$2C_{2}C_{0}B^{-1}+(d(x',S_{m'})-k')\leq  \big(d(x',S_{m'})-k\big)+1 \leq 
(d(x,S_{m'})-k). $$

\noindent {\bf Troisième cas.} On suppose $d(x,S_{0})> k$ et $T> k$.

Pour $j\in \{l'+1,...,l\}$ on a 
$$k+3P\leq d_{\max}(x,\mathcal Y_{m'}^{j})\leq T+M+2P+\de. $$
Comme $
\mathcal Y_{m'}^{j}\subset 4P\tg(x,b)$ le lemme~\ref{cardinal-tranche-geod}
montre que  le nombre de possibilités pour $\mathcal Y_{m'}^{j}$ est borné par $C_{2}(T-k)$ avec $C_{2}=C(\de,K,N,Q,P,M)$.
Par 
(\ref{form-r'-cas3-4j})
%(\ref{def-k'Z-24oct09})  
on a \begin{gather}\label{ineg-k'-trois-cas-20dec}|k'-(\frac{d(x',S_{0})+d(x,x')-d(x,S_{0})}{2}+\frac{M}{2})|\leq \frac{N}{2}+1.\end{gather} 
 Soit $y\in S_{m'}$.
 Il résulte de (\ref{ineg-commun-46-20dec}) que 
\begin{gather}\label{ineg-x'yx'bxbxy-troiscas-20dec}
d(x',y)-d(x',S_{0})+d(x,S_{0})\leq d(x,y)+P+3N+\de.\end{gather}
Donc 
\begin{gather*}d(x',y)-k'\leq d(x',y)-\frac{d(x',S_{0})+d(x,x')-d(x,S_{0})}{2}-\frac{M}{2}+\frac{N}{2}+1 \\ \leq d(x,y)-T-\frac{M}{2}+(P+7N/2+\de+1)\end{gather*} 
où la première inégalité a lieu par (\ref{ineg-k'-trois-cas-20dec}) et la deuxième par (\ref{ineg-x'yx'bxbxy-troiscas-20dec}). 
Comme cela est vrai pour tout $y\in S_{m'}$ on  en déduit 
\begin{gather}\nonumber (d(x',S_{m'})-k') \\ \label{ineg-x'Sk'-xSk-Tk-bis-cas3} \leq (d(x,S_{m'})-k) -(T-k)-\frac{M}{2}+(P+7N/2+\de+1).\end{gather}
 Le membre de gauche de (\ref{ineg-somme-aZZ'-ter-1}) est alors  majoré par 
$$\sum_{l\geq l'}\binom{l}{l'}
 \Big(\frac{C_{0}C_{2}B^{-1}(T-k)}{d(x,S_{m'})-k}\Big)^{(l-l')}
 \Big(\frac{d(x',S_{m'})-k'}{d(x,S_{m'})-k}\Big)^{l'}\leq 2
 $$
par le lemme~\ref{calcul-A-B}. Les hypothèses  du lemme~\ref{calcul-A-B}
 sont satisfaites car on suppose  $2C_{0}C_{2}B^{-1}\leq 1 $ 
 grâce à  $(H_{B})$ et
 $\frac{M}{2}\geq P+\frac{7N}{2}+\de+1$ grâce à $(H_{M})$
d'où 
  $(T-k)+(d(x',S_{m'})-k')\leq 
 (d(x,S_{m'})-k)$ grâce  à (\ref{ineg-x'Sk'-xSk-Tk-bis-cas3}). \cqfd

\noindent{\bf Démonstration du lemme~\ref{somme-aZZ'-ter-2}.} On distingue trois cas comme dans les démonstrations des  lemmes~\ref{oubli-x-xx'-x'} et~\ref{somme-aZZ'-ter-1}. On note $T=\frac{d(x,S_{0})-d(x',S_{0})+d(x,x')}{2}$ comme dans l'énoncé du lemme~\ref{somme-aZZ'-ter-2}. On fixe $b\in S_{0}$. En reprenant mot pour mot la preuve de (\ref{ineg-commun-46-20dec}), c'est-à-dire le début de la démonstration du  lemme~\ref{somme-aZZ'-ter-1} jusqu'à la distinction des trois cas, on obtient (\ref{ineg-commun-46-20dec}), c'est-à-dire que pour tout  $y\in S_{m'}$,   \begin{gather}\label{ineg-commun-47-20dec}d(x',y)-d(x',b)+d(x,b)\leq d(x,y)+P+2N+\de.\end{gather}

\noindent {\bf Premier cas.} On suppose $d(x,S_{0})\leq k$. 

Alors $m=0$ d'après le lemme~\ref{m-fini},  et cela  est impossible, puisque $m> m'$. 

\noindent {\bf Deuxième cas.} On suppose $d(x,S_{0})> k$ et $T\leq k$.

Pour $j\in \{l'+1,...,l\}$ on a 
$$k+3P\leq d_{\max}(x,\mathcal Y_{m'}^{j})\leq k+M+2P+\de.$$ Comme $
\mathcal Y_{m'}^{j}\subset 4P\tg(x,b)$ le lemme~\ref{cardinal-tranche-geod}
montre que le nombre de possibilités pour $\mathcal Y_{m'}^{j}$ est borné par $C_{2}=C(\de,K,N,Q,P,M)$. 
Pour $i\in \{m'+1,...,m\}$ on a $$k+P \leq d(x,S_{i})\leq k+M.$$ En effet l'inégalité de gauche a lieu par la condition ii) de la définition~\ref{defi-Y}. Comme $
S_{i}\subset P\tg(x,b)$ le lemme~\ref{cardinal-tranche-geod}
montre que le nombre de possibilités pour $S_{i}$ est borné par $C_{3}=C(\de,K,N,Q,P,M)$. On a $$k'=d(x',S_{0})-d(x,S_{0})+k+\frac{M}{2}$$ d'après 
(\ref{form-r'-cas2-4j})
%(\ref{def-k'Z-24oct09}) 
donc \begin{gather}\label{ineg-k'-deux-cas-20dec-bbis}k'\geq d(x',b)-d(x,b)+k+\frac{M}{2}-N.\end{gather} 

Soit $y\in S_{m'}$. 
 On a  
 $$d(x',y)-k'\leq  d(x',y)-d(x',b)+d(x,b)+N-k-\frac{M}{2}\leq d(x,y)-k-\frac{M}{2}
+P+3N+\de $$
où la première inégalité vient de (\ref{ineg-k'-deux-cas-20dec-bbis}) et la deuxième de (\ref{ineg-commun-47-20dec}). 
Comme cela est vrai pour tout $y\in S_{m'}$ on  en déduit 
\begin{gather}\label{ineg-bis-cas2-x'Sm'k'-xSm'k-b}\big(d(x',S_{m'})-k'\big)\leq \big(d(x,S_{m'})-k\big)-\frac{M}{2}+P+3N+\de
\end{gather}
On a $d(x,S_{m'+1})\leq k+M$ donc $d(S_{m'},S_{m'+1})\geq d(x,S_{m'})-k-M-N$ et \begin{gather} \nonumber d(S_{m'},S_{m'+1})+2M\geq d(x,S_{m'})-k+M-N \\ \label{ineg-Sm'-ter-2-cas2} \geq (d(x',S_{m'})-k')+3M/2-(P+4N+\de)\geq (d(x',S_{m'})-k')+1\end{gather} où 
l'avant-dernière inégalité a lieu par (\ref{ineg-bis-cas2-x'Sm'k'-xSm'k-b}) et 
où la dernière inégalité a lieu car on suppose $3M/2\geq P+4N+\de+1$, ce qui est permis par 
$(H_{M})$.

 Le membre de gauche de (\ref{ineg-somme-aZZ'-ter-2}) est  majoré par 
\begin{gather*}\sum_{m>m', l\geq l'} \Big(\frac{C_{0}C_{3}B^{-1}}{1-C_{0}C_{1}B^{-1}}\Big)^{m-m'}\binom{l}{l'}
\Big(\frac{C_{0}C_{2}B^{-1}}{d(S_{m'},S_{m'+1})+2M}\Big)^{l-l'}\\ \Big(\frac{d(x',S_{m'})-k'}{d(S_{m'},S_{m'+1})+2M}\Big)^{l'}\leq 2 \sum_{m>m'} \Big(\frac{C_{0}C_{3}B^{-1}}{1-C_{0}C_{1}B^{-1}}\Big)^{m-m'}\end{gather*} 
grâce au  
   lemme~\ref{calcul-A-B}. En effet les hypothèses  du lemme~\ref{calcul-A-B}
 sont satisfaites  car  $2C_{0}C_{2}B^{-1}\leq1$ par $(H_{B})$ donc 
 $$2C_{0}C_{2}B^{-1}+(d(x',S_{m'})-k')\leq  d(S_{m'},S_{m'+1})+2M$$ grâce à (\ref{ineg-Sm'-ter-2-cas2}). Ensuite on calcule 
$$\sum_{m>m'} \Big(\frac{C_{0}C_{3}B^{-1}}{1-C_{0}C_{1}B^{-1}}\Big)^{m-m'}=\frac{C_{0}C_{3}B^{-1}  }{1-C_{0}(C_{1}+C_{3})B^{-1}}\leq 1$$
 où la dernière inégalité a lieu car on suppose $B\geq 2C_{0}(C_{1}+C_{3})$ grâce à  $(H_{B})$.

\noindent {\bf Troisième cas.} On suppose $d(x,S_{0})> k$ et $T> k$. 

Soit $j\in \{l'+1,...,l\}$ et $w\in \mathcal Y_{m'}^{j}$. Il existe $y\in S_{m'}$ et $z\in S_{m'+1}$ tels que $w\in P\tg(y,z)$. On a  
$d(x,w)\geq d(x,y)-d(y,w)$ et $d(y,w)\leq d(y,z)+P$ donc $d(x,w) \geq d(x,y)-d(y,z)-P$. 
On en déduit 
$$d_{\max}(x,\mathcal Y_{m'}^{j})\geq d(x,S_{m'})-d(S_{m'},S_{m'+1})-2N-P.$$
D'autre part comme $j\in \{l'+1,...,l\}$ on a  $d_{\max}(x,\mathcal Y_{m'}^{j})\leq T+M+2P+\de$. 
Donc $d_{\max}(x,\mathcal Y_{m'}^{j})$ appartient à   l'intervalle 
$$[d(x,S_{m'})-d(S_{m'},S_{m'+1})-2N-P, T+M+2P+\de]$$ qui est de longueur 
$$(M+3P+2N+\de)+d(S_{m'},S_{m'+1})+T-d(x,S_{m'}).$$
On suppose $3P+2N+\de+1\leq M$,  ce qui est permis par $(H_{M})$. 
Comme $
\mathcal Y_{m'}^{j}\subset 4P\tg(x,b)$ le lemme~\ref{cardinal-tranche-geod} montre que  le nombre de possibilités pour $\mathcal Y_{m'}^{j}$ est borné par $$C_{2}\big(2M+d(S_{m'},S_{m'+1})+T-d(x,S_{m'})\big)$$ avec $C_{2}=C(\de,K,N,Q,P)$. Pour éviter des absurdités vérifions que  
$$2M+d(S_{m'},S_{m'+1})+T-d(x,S_{m'})\geq 1.$$ 
Cela résulte  des  inégalités    
$d(S_{m'},S_{m'+1})\geq d(x,S_{m'})-d(x,S_{m'+1})-N$  
et    $d(x,S_{m'+1})\leq d_{\max}(x,S_{m'+1})\leq T+M$  car on suppose $M\geq N+1$ grâce à $(H_{M})$. 

 Par (\ref{form-r'-cas3-4j})
  %(\ref{def-k'Z-24oct09})  
  on a \begin{gather}\label{ineg-k'-trois-cas-20dec-bis}|k'-(\frac{d(x',S_{0})+d(x,x')-d(x,S_{0})}{2}+\frac{M}{2})|\leq \frac{N}{2}+1.\end{gather} 
 Soit $y\in S_{m'}$.
 Il résulte de (\ref{ineg-commun-47-20dec}) que 
\begin{gather}\label{ineg-x'yx'bxbxy-troiscas-20dec-bis}
d(x',y)-d(x',S_{0})+d(x,S_{0})\leq d(x,y)+P+3N+\de.\end{gather}
Donc 
\begin{gather*}d(x',y)-k'\leq d(x',y)-\frac{d(x',S_{0})+d(x,x')-d(x,S_{0})}{2}-\frac{M}{2}+\frac{N}{2}+1\\ \leq d(x,y)-T-\frac{M}{2}+(P+7N/2+\de+1)\end{gather*}
où la première inégalité a lieu par (\ref{ineg-k'-trois-cas-20dec-bis}) et la deuxième par (\ref{ineg-x'yx'bxbxy-troiscas-20dec-bis}). 
Comme cela est vrai pour tout $y\in S_{m'}$ on en déduit 
\begin{gather}\label{ineg-x'Sk'-xSk-Tk-bis-cas3-bis} d(x',S_{m'})-k'  \leq d(x,S_{m'}) -T-\frac{M}{2}+(P+7N/2+\de+1).\end{gather} 

D'autre part on a 
$$k+P \leq d_{\max}(x,S_{m})\leq ...\leq 
d_{\max}(x,S_{m'+1})\leq T+M$$
donc le nombre de possibilités pour ces entiers est borné par $\binom{(T-k)+M+(m-m')}{m-m'}$. Comme $S_{m'+1}, ...,S_{m}$ sont inclus dans $P\tg(x,b)$, le lemme~\ref{cardinal-tranche-geod}
montre que, connaissant les entiers 
$d_{\max}(x,S_{m}), ...,  
d_{\max}(x,S_{m'+1})$, 
  le nombre de possibilités pour $S_{m'+1}, ...,  
S_{m}$ est borné par $C_{3}^{m-m'}$ avec $C_{3}=C(\de,K,N,Q,P)$.

Le membre de gauche de (\ref{ineg-somme-aZZ'-ter-2}) est alors  majoré par
\begin{gather*}\sum_{m>m'} \Big(\frac{C_{0}C_{3}B^{-1}}{1-C_{0}C_{1}B^{-1}}\Big)^{m-m'}\binom{(T-k)+M+(m-m')}{m-m'}\\ \sum_{l\geq l'}
\binom{l}{l'}
\Big(\frac{C_{0}C_{2}B^{-1}\big(2M+d(S_{m'},S_{m'+1})+T-d(x,S_{m'})\big)}{d(S_{m'},S_{m'+1})+2M}\Big)^{l-l'}\\ \Big(\frac{d(x',S_{m'})-k'}{d(S_{m'},S_{m'+1})+2M}\Big)^{l'}\\ 
%$$\Big(\frac{(d(x,S_{m'})-k) -(T-k)-\frac{M}{2}+(2P+6N+\de+1)}{d(S_{m'},S_{m'+1})+2M}\Big)^{l'}$$
\leq 2\sum_{m>m'} \Big(\frac{C_{0}C_{3}B^{-1}}{1-C_{0}C_{1}B^{-1}}\Big)^{m-m'}\binom{(T-k)+M+(m-m')}{m-m'}\end{gather*}
 où la dernière inégalité résulte  du 
   lemme~\ref{calcul-A-B}. Les hypothèses  du lemme~\ref{calcul-A-B}
 sont satisfaites  car  on suppose $2C_{0}C_{2}B^{-1}\leq1$ grâce à  $(H_{B})$ et car 
 \begin{gather*}\big(2M+d(S_{m'},S_{m'+1})+T-d(x,S_{m'})\big)+(d(x',S_{m'})-k')\\ \leq d(S_{m'},S_{m'+1}) +\frac{3M}{2}+(2P+7N/2+\de+1)\leq d(S_{m'},S_{m'+1})+2M\end{gather*} où la première inégalité a lieu par (\ref{ineg-x'Sk'-xSk-Tk-bis-cas3-bis}) et la deuxième a lieu car on suppose 
$ \frac{M}{2}\geq 2P+7N/2+\de+1$, ce qui est permis par 
 $(H_{M})$. 
Or grâce à 
(\ref{form-binom-24oct09}), 
\begin{gather*}\sum_{m\geq m'} \Big(\frac{C_{0}C_{3}B^{-1}}{1-C_{0}C_{1}B^{-1}}\Big)^{m-m'}\binom{(T-k)+M+(m-m')}{m-m'}\\ =\Big(1-\Big(\frac{C_{0}C_{3}B^{-1}}{1-C_{0}C_{1}B^{-1}}\Big)\Big)^{-((T-k)+M+1)}\leq 2e^{sd(x,x')}\end{gather*}
  car $(T-k)\leq T\leq d(x,x')$ et car on suppose $$1-\Big(\frac{C_{0}C_{3}B^{-1}}{1-C_{0}C_{1}B^{-1}}\Big)\geq \max(e^{-s},2^{-(M+1)^{-1}})$$ grâce à    $(H_{B})$. Ceci termine l'étude du troisième cas. %, et donc la démonstration du lemme~\ref{somme-aZZ'-ter-2} 
  \cqfd

On a démontré les lemmes~\ref{somme-aZZ'-ter-1} et~\ref{somme-aZZ'-ter-2} et donc aussi les lemmes~\ref{somme-aZZ'-ter}, 
\ref{somme-aZZ'-bis}, 
\ref{somme-aZZ'}, 
\ref{deuxieme-ineg-equiv-norme}, 
et  la proposition~\ref{equiv-normeYZ}. 

\subsection{Equivariance des opérateurs à compacts près. }

Pour terminer la démonstration de la proposition~\ref{enonce-ppal-KKC01}, il suffit de montrer la proposition suivante. 

\begin{prop}\label{lemme-compacite-equiv} Soit $T\in \R_{+}$ et $x,x'\in X$ vérifiant $d(x,x')=1$. 
Alors

\noindent a) $(e^{\tau  \theta^{\flat}_{x}}\del e^{-\tau  \theta^{\flat}_{x}}-
e^{\tau  \theta^{\flat}_{x'}}\del e^{-\tau  \theta^{\flat}_{x'}}
)_{\tau  \in [0,T]}$ appartient à $\K(\H_{x,s}[0,T])$, 

\noindent b) $(e^{\tau  \theta^{\flat}_{x}}h_{x}e^{-\tau  \theta^{\flat}_{x}}-
e^{\tau  \theta^{\flat}_{x'}}h_{x'} e^{-\tau  \theta^{\flat}_{x'}}
)_{\tau  \in [0,T]}$ appartient à $\K(\H_{x,s}[0,T])$, 

\noindent c) pour tout $r\in \N$, $(e^{\tau  \theta^{\flat}_{x}}u_{x,r}K_{x}e^{-\tau  \theta^{\flat}_{x}}-
e^{\tau  \theta^{\flat}_{x'}}u_{x',r}K_{x'} e^{-\tau  \theta^{\flat}_{x'}})_{\tau  \in [0,T]}$ appartient à $\K(\H_{x,s}[0,T])$. 
\end{prop}
La proposition~\ref{equiv-normeYZ}, qui a été établie au sous-paragraphe précédent,  assure l'équivalence des normes de $\H_{x,s}$ et $\H_{x',s}$. On en déduit, grâce à  la proposition~\ref{continuite-del-J-conj}, que   tous les opérateurs apparaissant dans la proposition~\ref{lemme-compacite-equiv} sont continus. La proposition~\ref{lemme-compacite-equiv} est donc seulement un énoncé de compacité. 

\noindent{\bf Démonstration de la 
proposition~\ref{enonce-ppal-KKC01} en admettant 
 la proposition~\ref{lemme-compacite-equiv}.}  Comme $X$
  est géodésique, 
  l'énoncé de la proposition~\ref{lemme-compacite-equiv} implique évidemment le même énoncé pour $x,x'$ quelconques (c'est-à-dire sans l'hypothèse $d(x,x')=1$). 
  Comme $J_{x}=\sum_{q=1}^{Q}h_{x}(1-\del h_{x}-h_{x}\del)^{q-1}+\sum_{r=1}^{\infty}u_{x,r}K_{x}$,    on 
    déduit des propositions~\ref{continuite-del-J-conj} et \ref{lemme-compacite-equiv} que 
  pour $T\in \R_{+}$ et $x,x'\in X$,  
$$(e^{\tau  \theta^{\flat}_{x}}(\del + J_{x}\del J_{x}) e^{-\tau  \theta^{\flat}_{x}}-
e^{\tau  \theta^{\flat}_{x'}}(\del + J_{x'}\del J_{x'}) e^{-\tau  \theta^{\flat}_{x'}}
)_{\tau  \in [0,T]}$$ est   un opérateur 
compact sur le  $\C[0,T]$-module hilbertien 
$\H_{x,s}[0,T]$. \cqfd

On va voir que la proposition~\ref{lemme-compacite-equiv} résulte du lemme suivant. On commence par remarquer que $\theta_{x'}-\theta_{x}$ est un opérateur borné (et même de norme $\leq 1$) sur $\H_{x,s}$ : en effet pour tous
$p,k,m,l_{0},...,l_{m}$, pour $Z\in \overline Y_{x}^{p,k,m,(l_{0},...,l_{m})}$  et pour $$(a_{1},\dots,a_{p},S_{0},...,S_{m},(\mathcal Y_{i}^{j})_{i\in \{0,\dots,m\}, j\in \{1,\dots ,l_{i}\}})\in (\pi_{x}^{p,k,m,(l_{0},...,l_{m})})^{-1}(Z)$$ 
$d(x,S_{0})-d(x',S_{0})$ appartient à $\{-1,0,1\}$ et est déterminé par $Z$
 (on rappelle que $S_{0}=\{a_{1},...,a_{p}\}$). D'autre part $\theta^{\flat}_{x}-\theta_{x}$ est un opérateur borné sur $\H_{x,s}$ d'après le lemme~\ref{continuite-eta} et de même $\theta^{\flat}_{x'}-\theta_{x'}$ est un opérateur borné sur $\H_{x',s}$.
  Grâce à l'équivalence des normes de $\H_{x,s}$ et $\H_{x',s}$ on en déduit  que l'opérateur $\theta^{\flat}_{x}-\theta^{\flat}_{x'}$ est borné sur 
  $\H_{x,s}$. 
   On note $[u,v]=uv-vu$. 

  \begin{lem}\label{lemme-compacite-equiv2} Soit $T\in \R_{+}$ et $x,x'\in X$ vérifiant $d(x,x')=1$. 
Alors

\noindent a) $(e^{\tau  \theta_{x}}(h_{x}-h_{x'}) e^{-\tau  \theta_{x}}
)_{\tau  \in [0,T]}$ appartient à $\K(\H_{x,s}[0,T])$,

\noindent b) pour tout $r\in \N$,  $(e^{\tau  \theta_{x}}(u_{x,r}K_{x}-u_{x',r}K_{x'}) e^{-\tau  \theta_{x}}
)_{\tau  \in [0,T]}$ appartient à $\K(\H_{x,s}[0,T])$,

\noindent c) $([(\theta^{\flat}_{x}-\theta^{\flat}_{x'}), e^{\tau  \theta_{x}}\del e^{-\tau  \theta_{x}}])_{\tau  \in [0,T]}$ appartient à $\K(\H_{x,s}[0,T])$, 

\noindent d)  $([(\theta^{\flat}_{x}-\theta^{\flat}_{x'}), e^{\tau  \theta_{x}}h_{x}e^{-\tau  \theta_{x}}])_{\tau  \in [0,T]}$ appartient à $\K(\H_{x,s}[0,T])$, 

\noindent e) pour tout $r\in \N$, 
$([(\theta^{\flat}_{x}-\theta^{\flat}_{x'}),e^{\tau  \theta_{x}}u_{x,r}K_{x}e^{-\tau  \theta_{x}}])_{\tau  \in [0,T]}$ appartient à $\K(\H_{x,s}[0,T])$. 
\end{lem}
\noindent{\bf Démonstration de la proposition~\ref{lemme-compacite-equiv} en admettant le lemme~\ref{lemme-compacite-equiv2}.} 
  Montrons d'abord que le c) du lemme~\ref{lemme-compacite-equiv2} implique le a) de la proposition~\ref{lemme-compacite-equiv}. Pour tout $i\in \N$, $([(\theta^{\flat}_{x}-\theta^{\flat}_{x'})^{i}, e^{\tau  \theta_{x}}\del e^{-\tau  \theta_{x}}])_{\tau  \in [0,T]}$ appartient à $\K(\H_{x,s}[0,T])$ et donc 
  $$([e^{\tau(\theta^{\flat}_{x}-\theta^{\flat}_{x'})}, e^{\tau  \theta_{x}}\del e^{-\tau  \theta_{x}}])_{\tau  \in [0,T]}\in \K(\H_{x,s}[0,T]).$$
  Or on a 
  $$e^{\tau  \theta^{\flat}_{x}}\del e^{-\tau  \theta^{\flat}_{x}}-
e^{\tau  \theta^{\flat}_{x'}}\del e^{-\tau  \theta^{\flat}_{x'}}=e^{\tau(\theta^{\flat}_{x'}-\theta_{x})}[e^{\tau(\theta^{\flat}_{x}-\theta^{\flat}_{x'})}, e^{\tau  \theta_{x}}\del e^{-\tau  \theta_{x}}]e^{-\tau(\theta^{\flat}_{x}-\theta_{x})}$$
  ce qui montre le a) de la proposition~\ref{lemme-compacite-equiv}.
  Ensuite le a) et le d)  du lemme~\ref{lemme-compacite-equiv2} impliquent le b) de la proposition~\ref{lemme-compacite-equiv}. En effet, 
  on a 
  $$e^{\tau  \theta^{\flat}_{x}}h_{x} e^{-\tau  \theta^{\flat}_{x}}-
e^{\tau  \theta^{\flat}_{x'}}h_{x'} e^{-\tau  \theta^{\flat}_{x'}}=
(e^{\tau  \theta^{\flat}_{x'}}(h_{x}-h_{x'})e^{-\tau  \theta^{\flat}_{x'}})
+
(e^{\tau  \theta^{\flat}_{x}}h_{x} e^{-\tau  \theta^{\flat}_{x}}-
e^{\tau  \theta^{\flat}_{x'}}h_{x} e^{-\tau  \theta^{\flat}_{x'}})$$
  donc la compacité de $(e^{\tau  \theta^{\flat}_{x}}h_{x} e^{-\tau  \theta^{\flat}_{x}}-
e^{\tau  \theta^{\flat}_{x'}}h_{x'} e^{-\tau  \theta^{\flat}_{x'}})_{\tau  \in [0,T]}$ que l'on cherche à établir résulte de la compacité des deux termes du membre de droite. D'abord   en conjugant  par 
  $e^{\tau(\theta^{\flat}_{x'}-\theta_{x})}$  on voit que le a)   du lemme~\ref{lemme-compacite-equiv2}  implique que l'opérateur 
  $(e^{\tau  \theta^{\flat}_{x'}}(h_{x}-h_{x'})e^{-\tau  \theta^{\flat}_{x'}})_{\tau  \in [0,T]}$ est compact. Ensuite on a l'égalité  
  $$e^{\tau  \theta^{\flat}_{x}}h_{x} e^{-\tau  \theta^{\flat}_{x}}-
e^{\tau  \theta^{\flat}_{x'}}h_{x} e^{-\tau  \theta^{\flat}_{x'}}=e^{\tau(\theta^{\flat}_{x'}-\theta_{x})}[e^{\tau(\theta^{\flat}_{x}-\theta^{\flat}_{x'})}, e^{\tau  \theta_{x}}h_{x} e^{-\tau  \theta_{x}}]e^{-\tau(\theta^{\flat}_{x}-\theta_{x})}$$
  donc le d)  du lemme~\ref{lemme-compacite-equiv2} montre la compacité de   $$(e^{\tau  \theta^{\flat}_{x}}h_{x} e^{-\tau  \theta^{\flat}_{x}}-
e^{\tau  \theta^{\flat}_{x'}}h_{x} e^{-\tau  \theta^{\flat}_{x'}})_{\tau  \in [0,T]}. $$
%  donc $(e^{\tau  \theta^{\flat}_{x}}h_{x}e^{-\tau  \theta^{\flat}_{x}}-
%e^{\tau  \theta^{\flat}_{x'}}h_{x'} e^{-\tau  \theta^{\flat}_{x'}}
%)_{\tau  \in [0,T]}$ appartient à $\K(\H_{x,s}[0,T])$ et 
Ceci termine la preuve du b) de la proposition~\ref{lemme-compacite-equiv}. 
Enfin par un argument similaire (en rempla\c cant $h_{x}$  par 
$u_{x,r}K_{x}$ et $h_{x'}$ par  $u_{x',r}K_{x'}$)  
 le b) et le e)  du lemme~\ref{lemme-compacite-equiv2} impliquent le c) de la proposition~\ref{lemme-compacite-equiv}. \cqfd

\noindent{\bf Démonstration du lemme~\ref{lemme-compacite-equiv2} en admettant le lemme~\ref{lemme-compacite-equiv3}.}   Pour tout $p\in \{1,...,p_{\max}\}$ et $n\in \N$ on note \label{def-Pn}
   $\P_{n}$ le projecteur orthogonal sur le sous-espace vectoriel de $\H^{\rightarrow}_{x,s}(\Delta_{p})$ engendré par les $e_{S}$ pour $S\in \Delta_{p}$ tel que $d(x,S)\leq n$, de sorte que $$(\P_{n} f)(S)=f(S) \text{ si }
d(x,S)\leq n \text{ et }(\P_{n} f)(S)=0 \text{ si } d(x,S)> n. $$ Dans les notations adoptées jusqu'ici on a donc $\P=\P_P$. 
On supposera toujours $n\geq P$, de sorte que $(1-\P_{n})(1-\P)=1-\P_{n}$. 
 Il est évident que le lemme~\ref{lemme-compacite-equiv2} résulte du lemme suivant. \cqfd
 
 \begin{lem}\label{lemme-compacite-equiv3} Soit $T\in \R_{+}$ et $x,x'\in X$ vérifiant $d(x,x')=1$ et $p\in \{2,...,p_{\max}\}$. 
Alors

\noindent a) $\sup_{\tau\in [0,T]}\|(1-\P_{n})e^{\tau\theta_{x}}(h_{x}-h_{x'}) e^{-\tau\theta_{x}}
\|_{\L(\H_{x,s}(\Delta_{p-1}),\H^{\rightarrow}_{x,s}(\Delta_{p}))}$ tend vers $0$ quand $n\to \infty$,

\noindent b) pour tout $r\in \N$,  $$\sup_{\tau\in [0,T]}\|(1-\P_{n})e^{\tau\theta_{x}}(u_{x,r}K_{x}-u_{x',r}K_{x'}) e^{-\tau\theta_{x}}
\|_{\L(\H_{x,s}(\Delta_{p-1}),\H^{\rightarrow}_{x,s}(\Delta_{p}))}$$ tend vers $0$ quand $n\to \infty$,

\noindent c) $\sup_{\tau\in [0,T]}\|(1-\P_{n})[(\theta^{\flat}_{x}-\theta^{\flat}_{x'}), e^{\tau\theta_{x}}\del e^{-\tau\theta_{x}}]\|_{\L(\H_{x,s}(\Delta_{p}),\H_{x,s}^{\rightarrow}(\Delta_{p-1}))}$ tend vers $0$ quand $n\to \infty$, 

\noindent d)  $\sup_{\tau\in [0,T]}\|(1-\P_{n})[(\theta^{\flat}_{x}-\theta^{\flat}_{x'}), e^{\tau\theta_{x}}h_{x}e^{-\tau\theta_{x}}]\|_{\L(\H_{x,s}(\Delta_{p-1}),\H_{x,s}^{\rightarrow}(\Delta_{p}))}$ tend vers $0$ quand $n\to \infty$, 

\noindent e) pour tout $r\in \N$, 
$$\sup_{\tau\in [0,T]}\|(1-\P_{n})[(\theta^{\flat}_{x}-\theta^{\flat}_{x'}),e^{\tau\theta_{x}}u_{x,r}K_{x}e^{-\tau\theta_{x}}]\|_{\L(\H_{x,s}(\Delta_{p-1}),\H_{x,s}^{\rightarrow}(\Delta_{p}))}$$ tend vers $0$ quand $n\to \infty$. 
\end{lem}

  %Pour montrer le lemme~\ref{lemme-compacite-equiv3} on aura besoin du lemme suivant. On rappelle que $\alpha\in ]0,1[$ (satisfaisant $(H_{\alpha})$) a déjà été fixé. 
  %
  %\begin{lem}\label{somme-ri-alpha} 
  %Pour tout $n\in \N^{*}$ et $C\in \R_{+}$ on a 
  %$$\sup_{r_{1},...,r_{n}\in \R_{+}\text{\  vérifiant \ }\sum_{i=1}^{n}r_{i}\leq C}
 %\Big( \frac{1}{n} \sum_{i=1}^{n}r_{i}^{\alpha}\Big)=C^{\alpha}n^{-\alpha}.$$
 % \end{lem}
 % \noindent {\bf Démonstration.}
 % L'application $x\mapsto x^{\alpha}$ de $\R_{+}$ dans $\R_{+}$ est concave. Donc le supremum est atteint lorsque $r_{i}=\frac{C}{n}$ pour tout $i\in \{1,...,n\}$. \cqfd
 % 
 %   \begin{lem}\label{somme-ri-alpha-bis} 
 % Pour toute fonction mesurable $f:[0,1]\to \R_{+}$ telle que $\int_{0}^{1}f(t)dt$ $< +\infty $  on a 
 % $$\int_{0}^{1}f(t)^{\alpha}dt\leq \Big(\int_{0}^{1}f(t)dt\Big)^{\alpha}.$$
 % \end{lem}
 % \noindent {\bf Démonstration.}
 % L'application $x\mapsto x^{\alpha}$ de $\R_{+}$ dans $\R_{+}$ est concave.  \cqfd
 % 

\noindent {\bf Démonstration du a) du  lemme~\ref{lemme-compacite-equiv3}.} La preuve qui suit est assez voisine de la démonstration 
du lemme~\ref{sl3-1jan0923}.
%de l'inégalité (\ref{ineg-(1-P)tildeH}), qui était une partie de la démonstration du lemme~\ref{continuite-Jx}. 
 %Pour simplifier les notations on va montrer seulement \begin{gather}\label{variante-a)-lemme-compacite-equiv3}
%\|(1-\P_{n})(h_{x}-h_{x'}) 
%\|_{\L(\H_{x,s}(\Delta_{p-1}),\H^{\rightarrow}_{x,s}(\Delta_{p}))}\rightarrow 0\text{\ \ quand\ \ } n\to \infty.\end{gather} La preuve du lemme~\ref{continuite-del-J-conj-d} indique le tout petit argument supplémentaire nécessaire pour montrer a). 

On déduira  a)  de (\ref{ineg-(1-P)-tilde-l0-xx'}) qui est une variante de l'inégalité (\ref{ineg-(1-P)-tilde-l0}). 

\begin{souslem}\label{mesure-t-psiSxt-xx'}
Il existe $C_{1}=C(\de,K,N)$ tel que pour tout $S\in \Delta_{p-1}$, la mesure de l'ensemble des $t\in [0,1]$ tels que $(h_{x,t}-h_{x',t})(e_S)\neq  0$ est $\leq \frac{C_{1}}{1+d(x,S)} $. 
\end{souslem}
\noindent{\bf Démonstration.} C'est une conséquence immédiate du lemme~\ref{lem-hxt-hx't-8j1630}.  \cqfd

Le sous-lemme suivant est une conséquence évidente du lemme~\ref{hxt-eS-connaissance}.

\begin{souslem}\label{hxt-eS-connaissance-r-r'} Pour tout $S\in \Delta_{p-1}$, 
$h_{x,t}(e_S)$ et $ h_{x',t}(e_S)$   ne dépendent  que de la connaissance des points de 
\begin{gather}\label{ens-hxt-eS-connaissance-r-r'} B(x,2\de+1)\cup B(S, N) \\ \nonumber  \cup  \bigcup_{a\in S}\{y\in (5\delta+2)\text{-}\geod(x,a), 
d(x,y)\in [td(x,S)-2N-3,td(x,S)+2]\}\end{gather} et des distances entre ces points. 
\end{souslem}
\noindent{\bf Démonstration.} On applique le 
lemme~\ref{hxt-eS-connaissance} à $x$ et $x'$ et on remarque que $5\de\tg(x',a)\subset (5\delta+2)\text{-}\geod(x,a)$ et $|d(x,S)-d(x',S)|\leq 1$. 
\cqfd

\noindent {\bf Suite de la démonstration du a).} 
Soient 
$k,m,l_{0},\dots ,l_{m}\in \N $ et \begin{gather*}Z\in \overline Y_{x}^{p,k,m,(l_{0},...,l_{m})}\text{\ \  vérifiant \ \ }r_{0}(Z)> k+P.\end{gather*} 
On pose $ \tilde l_{0} =0$ et $ \tilde l_{i}=l_{i-1}$ pour $i\in \{1,\pp,m+1\}$ et on   note $\Lambda_{Z}$ la partie de 
 $\overline Y_{x}^{\natural,p-1,k,m+1,(\tilde l_{0},...,\tilde l_{m+1}),1,0}$ formée des $\tilde Z$ vérifiant \begin{gather}\label{ineg-a-11j1936} r_{0}( Z) \leq r_{0}(\tilde Z)\leq  r_{0}(Z)+N\end{gather}  et 
tels que pour tout  \begin{gather*}(\tilde a_{1},\dots,\tilde a_{p-1},\tilde S_{0},...,\tilde S_{m+1},(\tilde {\mathcal Y}_{i}^{j})_{i\in \{0,\dots,m+1\}, j\in \{1,\dots ,\tilde l_{i}\}},\tilde {\mathcal Z}_{0}^{1}
)\\ \in (\pi_{x}^{\natural,p-1,k,m+1,(\tilde l_{0},...,\tilde l_{m+1}),1,0})^{-1}(\tilde Z)\end{gather*} 
il existe une énumération $(a_{1},\pp,a_{p})$ de $\tilde S_{1}$ vérifiant 
%   il existe $(a_{1},\pp,a_{p})$ (vérifiant 
% $\tilde S_{1}=\{a_{1},\pp,a_{p}\}$) tel que   
 \begin{gather}\label{cond-1j1432}\!\!\!\!\!\!\!(a_{1},\pp,a_{p},\tilde S_{1},...,\tilde S_{m+1},(\tilde {\mathcal Y}_{i+1}^{j})_{i\in \{0,\dots,m\}, j\in \{1,\dots , l_{i}\}})\in (\pi_{x}^{p,k,m,(l_{0},..., l_{m})})^{-1}(Z).\end{gather} 
 
Soit $t\in [0,1]$. On note $\Lambda_{Z,t}$  l'ensemble des  $\tilde Z\in \Lambda_{Z}$ 
tels que pour tout  \begin{gather*}(\tilde a_{1},\dots,\tilde a_{p-1},\tilde S_{0},...,\tilde S_{m+1},(\tilde {\mathcal Y}_{i}^{j})_{i\in \{0,\dots,m+1\}, j\in \{1,\dots ,\tilde l_{i}\}},\tilde {\mathcal Z}_{0}^{1}
)\\ \in (\pi_{x}^{\natural,p-1,k,m+1,(\tilde l_{0},...,\tilde l_{m+1}),1,0})^{-1}(\tilde Z)\end{gather*} on ait 
\begin{gather}\label{cond-Z01-p150-25oct09}\tilde{\mathcal Z}_{0}^{1}=
\bigcup_{b\in \tilde S_{0}}\{z\in \geod(x,b), d(x,z)=E(tr_{0}(\tilde Z))\}.\end{gather}

 La condition (\ref{cond-Z01-p150-25oct09}) implique que 
 pour $\tilde Z \in \Lambda_{Z,t}$ on a  
  $t_{0}^{1}(\tilde Z) =E(tr_{0}(\tilde Z))$. 

 \begin{souslem}\label{slem0-a-2j1537} Soit  
 $$(\tilde a_{1},\dots,\tilde a_{p-1},\tilde S_{0},...,\tilde S_{m+1},(\tilde {\mathcal Y}_{i}^{j})_{i\in \{0,\dots,m+1\}, j\in \{1,\dots ,\tilde l_{i}\}}
)\\ \in Y_{x}^{p-1,k,m+1,(\tilde l_{0},...,\tilde l_{m+1})}$$
  tel que $d(x,\tilde S_{1})\leq d(x,\tilde S_{0})\leq d(x,\tilde S_{1})+N$ et 
   qu'il existe $(a_{1},\pp,a_{p})$ vérifiant 
 (\ref{cond-1j1432}). 
 %$$(a_{1},\pp,a_{p},\tilde S_{2},...,\tilde S_{m+1},(\tilde {\mathcal Y}_{i+1}^{j})_{i\in \{0,\dots,m\}, j\in \{1,\dots , l_{i}\}})\in (\pi_{x}^{p,k,m,(l_{0},..., l_{m})})^{-1}(Z).$$ 
 On définit $\tilde{\mathcal Z}_{0}^{1}$ par (\ref{cond-Z01-p150-25oct09}). Alors $\tilde{\mathcal Z}_{0}^{1}$ est de diamètre inférieur ou égal à $P/3$  et  il existe $\tilde Z\in \Lambda_{Z,t}$ tel que 
\begin{gather}\label{elem-slem0-a-2j1537}(\tilde a_{1},\dots,\tilde a_{p-1},\tilde S_{0},...,\tilde S_{m+1},(\tilde {\mathcal Y}_{i}^{j})_{i\in \{0,\dots,m+1\}, j\in \{1,\dots ,\tilde l_{i}\}},\tilde {\mathcal Z}_{0}^{1}
)\end{gather} appartienne  à 
$(\pi_{x}^{\natural,p-1,k,m+1,(\tilde l_{0},...,\tilde l_{m+1}),1,0})^{-1}(\tilde Z)$.
%$ Y_{x}^{\natural,p-1,k,m+1,(\tilde l_{0},...,\tilde l_{m+1}),1,0}$. 
\end{souslem}
\noindent{\bf Démonstration.}
Soient  $z,z'\in \tilde{\mathcal Z}_{0}^{1}$. Soit  $b\in S$. On a $z,z'\in 2N\tg(x,b)$ et $d(x,z)=d(x,z')$ donc par $(H_{\de}(z,x,z',b))$, $d(z,z')\leq 2N+\de\leq P/3$. Comme $\tilde{\mathcal Z}_{0}^{1}$ est non vide et $P/3\leq M$,  l'argument que nous venons de donner montre aussi que la condition (\ref{cond-Z01-p150-25oct09}) est vérifiée par les autres éléments de la classe d'équivalence  $\tilde Z$ de l'élément (\ref{elem-slem0-a-2j1537}) et donc  $\tilde Z\in \Lambda_{Z,t}$. \cqfd

\noindent {\bf Suite de la démonstration du a).} 
Notre but est maintenant de montrer l'inégalité suivante, qui est une variante de (\ref{ineg-(1-P)-tilde-l0}) : il existe une constante $C_{2}=C(\de,K,N,Q,P,M,T)$ telle que 
 \begin{gather}\label{ineg-(1-P)-tilde-l0-xx'}\sup_{\tau\in [0,T]}|\xi_{Z}(e^{\tau\theta_{x}}(h_{x,t}-h_{x',t})e^{-\tau\theta_{x}}f)|^{2}\leq C_{2}\sum_{\tilde Z\in \Lambda_{Z,t}^{\neq}}
|\xi_{\tilde Z}(f)|^{2}\end{gather}
  où $\Lambda_{Z,t}^{\neq}$  est l'ensemble des  $\tilde Z\in \Lambda_{Z,t}$ 
tels que pour tout  \begin{gather*}(\tilde a_{1},\dots,\tilde a_{p-1},\tilde S_{0},...,\tilde S_{m+1},(\tilde {\mathcal Y}_{i}^{j})_{i\in \{0,\dots,m+1\}, j\in \{1,\dots ,\tilde l_{i}\}},\tilde {\mathcal Z}_{0}^{1}
)\\ \in (\pi_{x}^{\natural,p-1,k,m+1,(\tilde l_{0},...,\tilde l_{m+1}),1,0})^{-1}(\tilde Z)\end{gather*} on ait 
\begin{gather}\label{condition-r-3-delta}(h_{x,t}- h_{x',t})(e_{\tilde S_{0}})\neq 0.\end{gather}

 Nous allons montrer (\ref{ineg-(1-P)-tilde-l0-xx'})
 et en même temps justifier que la condition (\ref{condition-r-3-delta}) ne dépend  que de $\tilde Z$ (c'est-à-dire que pour $\tilde Z\in \Lambda_{Z,t}$ elle est vérifiée ou non simultanément pour tous les éléments de $(\pi_{x}^{\natural,p-1,k,m+1,(\tilde l_{0},...,\tilde l_{m+1}),1,0})^{-1}(\tilde Z)$).
 
  \begin{souslem}\label{slem1-a-2j1541}
  Soit 
 $\tilde Z\in \Lambda_{Z,t}$ et  \begin{gather*} 
(\tilde a_{1},\dots,\tilde a_{p-1},\tilde S_{0},...,\tilde S_{m+1},(\tilde {\mathcal Y}_{i}^{j})_{i\in \{0,\dots,m+1\}, j\in \{1,\dots ,\tilde l_{i}\}},\tilde {\mathcal Z}_{0}^{1}
)\\ \in 
(\pi_{x}^{\natural, p-1,k,m+1,(\tilde l_{0},...,\tilde l_{m+1}),1,0})^{-1}(\tilde Z).
 \end{gather*}  
Alors  $h_{x,t}(e_{\tilde S_{0}})$ et $h_{x',t} (e_{\tilde S_{0}})$ ne dépendent   que de la connaissance des  points de   
 \begin{gather}\label{ens-a-connaissance-22dec}B(\tilde S_{0}, M)\cup B(x,k+2M)\cup B(\tilde{\mathcal Z}_{0}^{1}, M). \end{gather}
 et des distances entre ces points. 
 \end{souslem}
  \noindent{\bf Remarque.} On devrait plutôt noter 
 $h_{x,t}(
 e_{\tilde a_{1}}\wedge ... \wedge e_{\tilde a_{p-1}})$ au lieu de  $h_{x,t}(e_{\tilde S_{0}})$ mais à partir de maintenant nous commettrons cet abus. 
 
 \noindent{\bf Démonstration.}
 Grâce au sous-lemme~\ref{hxt-eS-connaissance-r-r'},  $h_{x,t}(e_{\tilde S_{0}})$ et $h_{x',t} (e_{\tilde S_{0}})$ ne dépendent  que de la connaissance des  points de  
 \begin{gather}\nonumber B(x,2\de+1)\cup B(\tilde S_{0}, N)\cup   \\ \label{ens-4.53-tildeS0-22dec} \bigcup_{a\in \tilde S_{0}}\{y\in (5\delta+2)\text{-}\geod(x,a), 
d(x,y)\in [tr_{0}(\tilde Z)-2N-3,tr_{0}(\tilde Z)+2]\}\end{gather}
et des distances entre ces points. 
Il suffit donc de montrer que (\ref{ens-4.53-tildeS0-22dec})  est inclus dans (\ref{ens-a-connaissance-22dec}). 
  Il est évident que $B(x,2\de+1)\cup B(\tilde S_{0},N)$ est inclus dans (\ref{ens-a-connaissance-22dec}). 
  Soient $a\in \tilde S_{0}$, $y\in (5\delta+2)\text{-}\geod(x,a)$ vérifiant  
$$d(x,y)\in [tr_{0}(\tilde Z)-2N-3,tr_{0}(\tilde Z)+2].$$ Soit  $z\in \geod(x,a)$ vérifiant $d(x,z)=E(tr_{0}(\tilde Z))$, si bien que $z$ appartient à $\tilde {\mathcal Z}_{0}^{1}$. Alors $|d(x,y)-d(x,z)|\leq 2N+3$ donc $(H_{\de}(y,x,z,a))$ montre $d(y,z)\leq (2N+3)+(5\de+2)+\de$ et on suppose $(2N+3)+(5\de+2)+\de\leq M$ grâce à  $(H_{M})$. Donc  
l'ensemble (\ref{ens-4.53-tildeS0-22dec}) 
est inclus dans (\ref{ens-a-connaissance-22dec}). 
\cqfd
 
 Le sous-lemme~\ref{slem1-a-2j1541} implique immédiatement 
 que  la 
 condition (\ref{condition-r-3-delta}) ne dépend que de $\tilde Z$.

  \begin{souslem}\label{slem2-a-2j1548}
  Le cardinal de $ \Lambda_{Z,t}$ est majoré par une constante de la forme $C(\de,K,N,Q,P,M)$. 
  \end{souslem}
 \noindent{\bf Démonstration.}
Grâce au lemme~\ref{nombre-dist-connaitre-par-point-natural}, pour connaître les distances entre les points de 
\begin{gather}\label{ens1-p151-25oct09}B(\tilde S_{0}, M)\cup   B(\tilde {\mathcal Z}_{0}^{1}, M)\end{gather}  
 et ceux de 
\begin{gather}\label{ens2-p151-25oct09} 
\bigcup _{ i\in \{1,\dots ,m+1\}}
B(\tilde S_{i},M)
\cup  \bigcup _{i\in \{1,\dots,m+1\}, j\in \{1,\dots ,\tilde l_{i}\}}B(\tilde {\mathcal Y}_{i}^{j},M)\cup B(x,k+2M)\end{gather}
il suffit de connaître les distances entre les points de (\ref{ens1-p151-25oct09}) 
 et $C'$ points de (\ref{ens2-p151-25oct09}), avec $C'=C(\de,K,N,Q,P,M)$ et ces distances sont déterminées à $C''=C(\de,K,N,Q,P,M)$ près par les distances de $\tilde S_{1}$ à ces $C'$ points (qui font partie de la donnée de $Z$) et l'entier $t_{0}^{1}(\tilde Z)$, qui est lui-même déterminé à $C'''= C(\de,K,N,Q,P,M)$ près par $r_{0}(Z)$ et $t$. \cqfd
 
 \noindent {\bf Suite de la démonstration du a).}
 On termine maintenant la preuve de (\ref{ineg-(1-P)-tilde-l0-xx'}). 
 Pour \begin{gather*} 
(\tilde a_{1},\dots,\tilde a_{p-1},\tilde S_{0},...,\tilde S_{m+1},(\tilde {\mathcal Y}_{i}^{j})_{i\in \{0,\dots,m+1\}, j\in \{1,\dots ,\tilde l_{i}\}},\tilde {\mathcal Z}_{0}^{1}
)\\ \in 
(\pi_{x}^{\natural, p-1,k,m+1,(\tilde l_{0},...,\tilde l_{m+1}),1,0})^{-1}(\tilde Z).
 \end{gather*} 
on considère 
\begin{gather}
\label{somme-1j1621}
\sum _{(b_{1},...,b_{p})}
\big((h_{x,t}-h_{x',t})  (e_{\tilde a_{1}}\wedge ...\wedge e_{\tilde a_{p-1}})\big)(b_{1},...,b_{p}), 
\end{gather} 
où la somme porte sur les énumérations $(b_{1},...,b_{p})$ de $\tilde S_{1}$  telles que 
$$(b_{1},\dots,b_{p},\tilde S_{1},...,\tilde S_{m+1},(\tilde {\mathcal Y}_{i+1}^{j})_{i\in \{0,\dots,m\}, j\in \{1,\dots ,\tilde l_{i}\}}
)\in (\pi_{x}^{p,k,m,(l_{0},..., l_{m})})^{-1}(Z).$$ 
Comme la somme (\ref{somme-1j1621}) a au plus $p!$ termes, 
 le 3) de la proposition~\ref{recap-supp-connaiss-H-uK} montre qu'elle 
est majorée par une constante de la forme $C(\de,K,N,Q,P,M)$. 
D'après le sous-lemme~\ref{slem1-a-2j1541} la somme (\ref{somme-1j1621}) ne dépend que de $\tilde Z$  et on peut donc la noter $\alpha_{Z,\tilde Z,t}$. 
D'après le sous-lemme~\ref{slem0-a-2j1537} on a  $$\xi_{Z}(e^{\tau\theta_{x}}(h_{x,t}-h_{x',t})e^{-\tau\theta_{x}}  f)=\frac{1}{(p-1)!}\sum_{\tilde Z\in \Lambda_{Z,t}^{\neq}} e^{\tau(\rho_{x}^{1}(\tilde Z)-\rho_{x}^{0}(\tilde Z))}
\alpha_{Z,\tilde Z,t}\xi_{\tilde Z}(f). $$
Grâce à (\ref{ineg-a-11j1936})  on a $\rho_{x}^{1}(\tilde Z)-\rho_{x}^{0}(\tilde Z)\leq 2N$. 
 Par Cauchy-Schwarz et grâce au sous-lemme~\ref{slem2-a-2j1548} 
 on en déduit (\ref{ineg-(1-P)-tilde-l0-xx'}). 
% $$|\xi_{Z}((h_{x,t}-h_{x',t})f)|^{2}\leq C\sum_{\tilde Z\in \Lambda_{Z,t}^{\neq}}
%|\xi_{\tilde Z}(f)|^{2}$$ avec $C=C(\de,K,N,Q,P,M)$.
% Ceci termine la preuve de (\ref{ineg-(1-P)-tilde-l0-xx'}). 

Montrons maintenant a) à l'aide de 
(\ref{ineg-(1-P)-tilde-l0-xx'}).  
Soit $C_{2}$ comme dans 
 (\ref{ineg-(1-P)-tilde-l0-xx'}).  Soient $k,m,l_{0},\dots ,l_{m}\in \N $. On pose  $ \tilde l_{0}=0$ et $ \tilde l_{i}=l_{i-1}$ pour $i\in \{1,\pp,m+1\}$ comme précédemment. 
 \begin{souslem}\label{slem-a-3-2j1554}
Il existe $C_{3}= C(\de,K,N,Q,P,M)$ tel que pour tout 
$$(a_{1},\dots,a_{p},S_{0},...,S_{m},(\mathcal Y_{i}^{j})_{i\in \{0,\dots,m\}, j\in \{1,\dots ,l_{i}\}}) \in 
Y_{x}^{p,k,m,(l_{0},...,l_{m})}$$
 le nombre de possibilités pour $(\tilde a_{1},\dots,\tilde a_{p-1})$ tels que 
 $$(\tilde a_{1},\dots,\tilde a_{p-1},\{\tilde a_{1},\dots,\tilde a_{p-1}\},  S_{0},...,S_{m},
( {\mathcal Y}_{i-1}^{j})_{i\in \{0,\dots,m+1\}, j\in \{1,\dots ,\tilde l_{i}\}})$$  appartienne à $
Y_{x}^{p-1,k,m+1,(\tilde l_{0},...,\tilde l_{m+1})}
$ 
 et vérifie 
$
d(x,S_{0}) \leq d(x,\{\tilde a_{1},\dots,\tilde a_{p-1}\}) \leq  d(x,S_{0}) +N$ soit $\leq C_{3}$.
\end{souslem}
\noindent{\bf Démonstration.} 
 Pour $a\in S_{0}$ et $i\in \{1,...,p-1\}$ on a $a\in 2F\tg(x,\tilde a_{i})$ et $|d(x,a)-d(x,\tilde a_{i})|\leq 2N$, d'où $d(a,\tilde a_{i})\leq 2F+2N$. \cqfd 
   
   \noindent {\bf Fin de la démonstration du a).}
Soit  $$Z\in \overline Y_{x}^{p,k,m,(l_{0},...,l_{m})}\text{\ \  vérifiant \ \ }r_{0}(Z)> k+P.$$
  Grâce au sous-lemmes~\ref{mesure-t-psiSxt-xx'} et~\ref{slem-a-3-2j1554},  et comme   $\tilde{\mathcal Z}_{0}^{1}$ est déterminé par (\ref{cond-Z01-p150-25oct09}) et $r_{0}(Z)\leq r_{0}(\tilde Z)$, 
 on  a alors 
 \begin{gather*}\int_{0}^{1}\Big(\sum_{\tilde Z\in \Lambda_{Z,t}^{\neq}}\sharp\big((\pi_{x}^{\natural,p-1,k,m+1,(\tilde l_{0},...,\tilde l_{m+1}),1,0})^{-1}(\tilde Z)\big)\Big)dt\\ \leq \frac{C_{1}C_{3}}{r_{0}(Z)+1}\sharp\big((\pi_{x}^{p,k,m,(l_{0},..., l_{m})})^{-1}(Z)\big).\end{gather*} 
  Notons $I_{Z}$ l'ensemble des $t\in [0,1]$ tels qu'il existe $\tilde Z\in \Lambda_{Z,t}^{\neq}$ vérifiant   $$\sharp\big((\pi_{x}^{\natural,p-1,k,m+1,(\tilde l_{0},...,\tilde l_{m+1}),1,0})^{-1}(\tilde Z)\big)\geq (r_{0}(Z)+1)^{-\frac{1}{2}}\sharp\big((\pi_{x}^{p,k,m,(l_{0},..., l_{m})})^{-1}(Z)\big). $$ La mesure de $I_{Z}$ est donc $\leq C_{1}C_{3}(r_{0}(Z)+1)^{-\frac{1}{2}}$. Grâce à Cauchy-Schwarz et à  (\ref{ineg-(1-P)-tilde-l0-xx'}) on obtient 
que 
\begin{gather}\nonumber \sup_{\tau\in [0,T]}\Big|\xi_{Z}\Big(e^{\tau\theta_{x}}\Big(\int_{t\in I_{Z}}(h_{x,t}-h_{x',t})dt\Big) e^{-\tau\theta_{x}} f\Big)\Big|^{2}
\\
\nonumber
\leq C_{1}C_{3}(r_{0}(Z)+1)^{-\frac{1}{2}}\sup_{\tau\in [0,T]}\int_{t\in [0,1]} \big|\xi_{Z}(e^{\tau\theta_{x}}(h_{x,t}-h_{x',t}) e^{-\tau\theta_{x}} f)\big|^{2}dt
%\\ \nonumber 
%\Big|\sum _{\substack{(a_{1},\dots,a_{p},S_{1},...,S_{m},(\mathcal Y_{i}^{j})_{i\in \{0,\dots,m\}, j\in \{1,\dots ,l_{i}\}}) \\ \in 
%(\pi_{x}^{p,k,m,(l_{0},...,l_{m})})^{-1}(Z)}} ((h_{x,t}-h_{x',t}) f)(a_1,...,a_p)\Big|^{2}dt 
\\ \label{ineg-IZ-xx'} \leq C_{1}C_{2}C_{3}(r_{0}(Z)+1)^{-\frac{1}{2}}
 \sum_{\tilde Z\in \Lambda_{Z}}(r_{0}(\tilde Z)+1)^{-1}
|\xi_{\tilde Z}(f)|^{2}.\end{gather}
%$$ 
%\Big| \sum _{\substack{(\tilde a_{1},\dots,\tilde a_{p-1},\tilde S_{1},...,\tilde S_{m+1},
%(\tilde {\mathcal Y}_{i}^{j})_{i\in \{0,\dots,m+1\}, j\in \{1,\dots ,\tilde l_{i}\}}
%,\tilde {\mathcal Z}_{0}^{1}
%)\\ \in 
%(\pi_{x}^{\natural, p-1,k,m+1,(\tilde l_{0},...,\tilde l_{m+1}),1,0})^{-1}(\tilde Z)}}  f(\tilde a_{1},\dots,\tilde a_{p-1})\Big|^{2}$$
  
D'après le sous-lemme~\ref{slem-a-3-2j1554} et le lemme~\ref{lemme-cardinaux}, 
 il existe une constante $C_{4}=C(\de,K,N,Q,P,M)$ telle   que pour $\tilde Z\in \Lambda_{Z}$  on ait 
 \begin{gather}\label{ineg-a-C4-21dec}\sharp\big((\pi_{x}^{\natural,p-1,k,m+1,(\tilde l_{0},...,\tilde l_{m+1}),1,0})^{-1}(\tilde Z)\big)\leq C_{4}\sharp\big((\pi_{x}^{p,k,m,(l_{0},..., l_{m})})^{-1}(Z)\big).\end{gather} 
 D'autre part pour $t\not\in I_{Z}$ et $\tilde Z\in \Lambda_{Z,t}^{\neq}$ on a $$\sharp\big((\pi_{x}^{\natural,p-1,k,m+1,(\tilde l_{0},...,\tilde l_{m+1}),1,0})^{-1}(\tilde Z)\big)\leq (r_{0}(Z)+1)^{-\frac{1}{2}}\sharp\big((\pi_{x}^{p,k,m,(l_{0},..., l_{m})})^{-1}(Z)\big). $$
  Par Cauchy-Schwarz on déduit alors de (\ref{ineg-(1-P)-tilde-l0-xx'}) que 
  \begin{gather}\nonumber \sup_{\tau\in [0,T]}\big(\sharp(\pi_{x}^{p,k,m,(l_{0},..., l_{m})})^{-1}(Z)\big)^{-\alpha}\Big|\xi_{Z}\Big(e^{\tau\theta_{x}}\Big(\int_{t\not \in I_{Z}}(h_{x,t}-h_{x',t})dt\Big)e^{-\tau\theta_{x}}  f\Big)\Big|^{2}  
  %$$\Big|\sum _{\substack{(a_{1},\dots,a_{p},S_{1},...,S_{m},(\mathcal Y_{i}^{j})_{i\in \{0,\dots,m\}, j\in \{1,\dots ,l_{i}\}}) \\ \in 
%(\pi_{x}^{p,k,m,(l_{0},...,l_{m})})^{-1}(Z)}} ((\int_{t\not \in I_{Z}}(h_{x,t}-h_{x',t}) dt) f)(a_1,...,a_p)\Big|^{2}$$ 
\\ \nonumber 
\leq C_{2}(r_{0}(Z)+1)^{-\frac{\alpha}{2}}\\  \label{ineg-hors-IZ-xx'} 
\sum_{\tilde Z\in \Lambda_{Z}} (r_{0}(\tilde Z)+1)^{-1}\big(\sharp(\pi_{x}^{\natural,p-1,k,m+1,(\tilde l_{0},...,\tilde l_{m+1}),1,0})^{-1}(\tilde Z)\big)^{-\alpha}|\xi_{\tilde Z}(f)|^{2}.\end{gather}
% \Big| \sum _{\substack{(\tilde a_{1},\dots,\tilde a_{p-1},\tilde S_{1},...,\tilde S_{m+1},
%(\tilde {\mathcal Y}_{i}^{j})_{i\in \{0,\dots,m+1\}, j\in \{1,\dots ,\tilde l_{i}\}}
%,\tilde {\mathcal Z}_{0}^{1}
%)\\ \in 
%(\pi_{x}^{\natural, p-1,k,m+1,(\tilde l_{0},...,\tilde l_{m+1}),1,0})^{-1}(\tilde Z)}}  f(\tilde a_{1},\dots,\tilde a_{p-1})\Big|^{2}\end{gather}

Comme     $$h_{x}-h_{x'}
=\int_{t \in I_{Z}}(h_{x,t}-h_{x',t}) dt+
\int_{t\not \in I_{Z}}(h_{x,t}-h_{x',t}) dt, $$ en combinant les inégalités (\ref{ineg-IZ-xx'}), (\ref{ineg-a-C4-21dec}) et 
  (\ref{ineg-hors-IZ-xx'}) et par Cauchy-Schwarz on obtient que 
\begin{gather}\nonumber \sup_{\tau\in [0,T]}\big(\sharp(\pi_{x}^{p,k,m,(l_{0},..., l_{m})})^{-1}(Z)\big)^{-\alpha}|\xi_{Z}(e^{\tau\theta_{x}}(h_{x}-h_{x'})e^{-\tau\theta_{x}}f)|^{2}\\ \nonumber
%$$\Big|\sum _{(a_{1},\dots,a_{p},S_{1},...,S_{m},(\mathcal Y_{i}^{j})_{i\in \{0,\dots,m\}, j\in \{1,\dots ,l_{i}\}}) \in 
%(\pi_{x}^{p,k,m,(l_{0},...,l_{m})})^{-1}(Z)} ((h_{x}-h_{x'})  f)(a_1,...,a_p)\Big|^{2}$$ 
\leq 2\Big(C_{1}C_{2}C_{3}C_{4}^{\alpha}(r_{0}(Z)+1)^{-\frac{1}{2}}+
 C_{2}(r_{0}(Z)+1)^{-\frac{\alpha}{2}}\Big)
\sum_{\tilde Z\in \Lambda_{Z}} 
(r_{0}(\tilde  Z)+1)^{-1}\\ \label{C1234-p153-25oct09}
\big(\sharp(\pi_{x}^{\natural,p-1,k,m+1,(\tilde l_{0},...,\tilde l_{m+1}),1,0})^{-1}(\tilde Z)\big)^{-\alpha}|\xi_{\tilde Z}(f)|^{2}.\end{gather} 
    De plus pour $\tilde Z\in \Lambda_{Z}$ on a $\prod_{i=0}^{m}s_{i}(Z)^{-l_{i}}=\prod_{i=0}^{m+1}s_{i}(\tilde Z)^{-\tilde l_{i}}$. 
   Pour calculer la norme de $(1-\P_{n})e^{\tau\theta_{x}}(h_{x}-h_{x'})e^{-\tau\theta_{x}}$ on peut se limiter aux $Z$ tels que $r_{0}(Z)\geq n$ et on déduit donc de (\ref{C1234-p153-25oct09}) que 
\begin{gather*}\sup_{\tau\in [0,T]}\|(1-\P_{n})e^{\tau\theta_{x}}(h_{x}-h_{x'}) e^{-\tau\theta_{x}}
\|_{\L(\H_{x,s}^{\natural,1,0}(\Delta_{p-1}),\H^{\rightarrow}_{x,s}(\Delta_{p}))}^{2}\\ \leq  2p! B
\Big(C_{1}C_{2}C_{3}C_{4}^{\alpha}(n+1)^{-\frac{1}{2}}+
 C_{2}(n+1)^{-\frac{\alpha}{2}}\Big) \end{gather*}
où le facteur $p!$ est dû au fait que $\tilde Z$ détermine $Z$ à permutation près de $a_{1},...,a_{p}$. 
Ceci termine la preuve de a).

\noindent {\bf Démonstration du b) du  lemme~\ref{lemme-compacite-equiv3}.} La preuve de b) n'introduit aucune idée nouvelle par rapport à celle de a), donc on sera bref. Soit $r\in \N$. 
On déduira  b)  de (\ref{ineg-(1-P)-tilde-l0-bis-xx'}) qui est une variante de l'inégalité (\ref{ineg-(1-P)-tilde-l0-bis}), de la même fa\c con que l'on avait montré le lemme~\ref{sl4-1jan0923} 
 à l'aide de (\ref{ineg-(1-P)-tilde-l0-bis}).

\begin{souslem}\label{mesure-t-murtxa-xx'}
Il existe $C=C(\de,K,N,r)$ tel que pour tout $a\in X$, la mesure de l'ensemble des $t\in [0,1]$ tels que $\mu_{r,t}(x,a)\neq \mu_{r,t}(x',a)$ est $\leq \frac{C}{1+d(x,a)} $. 
\end{souslem}
\noindent{\bf Démonstration.} Cela résulte du lemme~\ref{lem-murtxx'-8j1848}.  \cqfd

\begin{souslem}\label{mesure-t-uxrtK-xx'}
Il existe $C_{1}=C(\de,K,N,Q,r)$ tel que pour tout $S\in \Delta_{p-1}$, la mesure de l'ensemble des $(t,t_{1},...,t_{Q})\in [0,1]^{Q+1}$ tels que $$u_{x,r,t}K_{x,Q,(t_{1},\dots ,t_{Q})} (e_{S}) \neq u_{x',r,t}K_{x',Q,(t_{1},\dots ,t_{Q})} (e_{S})$$ est $\leq \frac{C_{1}}{1+d(x,S)} $. 
\end{souslem}
\noindent{\bf Démonstration.} Il existe $C=C(\de,K,N,Q,r)$ tel que pour $S\in \Delta_{p-1}$ la connaissance de $u_{x,r,t}K_{x,Q,(t_{1},\dots ,t_{Q})} (e_{S}) $ et celle de $u_{x',r,t}K_{x',Q,(t_{1},\dots ,t_{Q})} (e_{S}) $ ne dépendent que de la connaissance des $h_{x,t_{i}}(e_{S'})$ et $h_{x',t_{i}}(e_{S'})$ pour $i=1,...,Q$ et $S'\in \Delta$ vérifiant  $d(S,S')\leq C$ et des $\mu_{r,t}(x,a)$ et $\mu_{r,t}(x',a)$ pour $a\in X$ 
 vérifiant $d(a,S)\leq C$. On applique alors le sous-lemme~\ref{mesure-t-psiSxt-xx'}  à ces parties $S'$  et le sous-lemme~\ref{mesure-t-murtxa-xx'} à ces points $a$.  \cqfd 
 
   \noindent {\bf Suite de la démonstration du b).}
 Soient 
$k,m,l_{0},\dots ,l_{m}\in \N $ et $Z\in \overline Y_{x}^{p,k,m,(l_{0},...,l_{m})}$ vérifiant $r_{0}(Z)> k+P$.  
On pose $ \tilde l_{0} =0$ et $ \tilde l_{i}=l_{i-1}$ pour $i\in \{1,\pp,m+1\}$ et on   note $\Lambda_{Z}$ la partie de  $\overline Y_{x}^{\natural, p-1,k,m+1,(\tilde l_{0},...,\tilde l_{m+1}),Q,1}$ formée des $\tilde Z$  vérifiant 
  \begin{gather}\label{r0-r1-r-xx'}|r_{0}(\tilde Z)-r_{1}(\tilde Z)-r|\leq QF \end{gather} 
 et 
tels que pour tout  \begin{gather*}(\tilde a_{1},\dots,\tilde a_{p-1},\tilde S_{0},...,\tilde S_{m+1},(\tilde {\mathcal Y}_{i}^{j})_{i\in \{0,\dots,m+1\}, j\in \{1,\dots ,\tilde l_{i}\}},(\tilde {\mathcal Z}_{0}^{j})_{ j\in \{1,\dots ,Q\}},\tilde {\mathcal Z}_{1}^{1})\\ \in (\pi_{x}^{\natural, p-1,k,m+1,(\tilde l_{0},...,\tilde l_{m+1}),Q,1})^{-1}(\tilde Z)\end{gather*} 
il existe une énumération $(a_{1},\pp,a_{p})$ de $\tilde S_{1}$ vérifiant 
%il existe $(a_{1},\pp,a_{p})$ (vérifiant 
% $\tilde S_{1}=\{a_{1},\pp,a_{p}\}$) tel que  
 \begin{gather}\label{cond-Z-b-2j1818}(a_{1},\pp,a_{p},\tilde S_{1},...,\tilde S_{m+1},(\tilde {\mathcal Y}_{i+1}^{j})_{i\in \{0,\dots,m\}, j\in \{1,\dots , l_{i}\}})\in (\pi_{x}^{p,k,m,(l_{0},..., l_{m})})^{-1}(Z).\end{gather}

 Soient  $t,t_{1},\pp,t_{Q}\in [0,1]$. On note $\Lambda_{Z,t,(t_{1},...,t_{Q})}$  l'ensemble des  $\tilde Z\in \Lambda_{Z}$ 
 tels que pour tout  \begin{gather*}(\tilde a_{1},\dots,\tilde a_{p-1},\tilde S_{0},...,\tilde S_{m+1},(\tilde {\mathcal Y}_{i}^{j})_{i\in \{0,\dots,m+1\}, j\in \{1,\dots ,\tilde l_{i}\}},(\tilde {\mathcal Z}_{0}^{j})_{ j\in \{1,\dots ,Q\}},\tilde {\mathcal Z}_{1}^{1})\\ \in (\pi_{x}^{\natural, p-1,k,m+1,(\tilde l_{0},...,\tilde l_{m+1}),Q,1})^{-1}(\tilde Z)\end{gather*} 
 on ait 
 \begin{gather} \label{cond-b-1j1756}\tilde{\mathcal Z}_{0}^{j}=\bigcup_{b\in \tilde S_{0}}\{z\in \geod(x,b), d(x,z)=E(t_{j}r_{0}(\tilde Z))\}\text{ pour }j\in \{1,...,Q\}
 \\ \label{cond-b-1j1757}\text{et\ \ \ }\tilde{\mathcal Z}_{1}^{1}=\bigcup_{b\in \tilde S_{1}}\{z\in \geod(x,b), d(x,z)=E((1-t)r_{1}(\tilde Z))\}.\end{gather} 
 Pour $\tilde Z\in \Lambda_{Z,t,(t_{1},...,t_{Q})}$, les conditions (\ref{cond-b-1j1756}) et (\ref{cond-b-1j1757}) impliquent 
 $$t_{0}^{j}(\tilde Z)=E(t_{j}r_{0}(\tilde Z))\text{\ \ pour\ \ }j\in \{1,...,Q\}, \ \ \ \text{et}\ \ \ t_{1}^{1}(\tilde Z)
 =E((1-t)r_{1}(\tilde Z)).$$ 

\noindent{\bf Remarque.} Les notations $\Lambda_{Z}$ et $\Lambda_{Z,t,(t_{1},...,t_{Q})}$ que nous venons d'introduire sont les mêmes que dans la preuve du lemme~\ref{sl4-1jan0923} 
et les conditions (\ref{r0-r1-r-xx'}), (\ref{cond-Z-b-2j1818}), 
(\ref{cond-b-1j1756}) et (\ref{cond-b-1j1757}) coïncident avec les conditions (\ref{r0-r1-r}), (\ref{notation-34-2j1259}), (\ref{defZ0j-1j1040}) et (\ref{defZ1j-1j1040}). 

\begin{souslem} \label{slem0-b-2j1820} Soit  $$(\tilde a_{1},\dots,\tilde a_{p-1},\tilde S_{0},...,\tilde S_{m+1},(\tilde {\mathcal Y}_{i}^{j})_{i\in \{0,\dots,m+1\}, j\in \{1,\dots ,\tilde l_{i}\}})\in Y_{x}^{p-1,k,m+1,(\tilde l_{0},...,\tilde l_{m+1})}$$
  tel que $  |d(x,\tilde S_{0})-d(x,\tilde S_{1})-r|\leq QF
  $ et 
   qu'il existe $(a_{1},\pp,a_{p})$ vérifiant 
 (\ref{cond-Z-b-2j1818}). 
On définit $(\tilde{\mathcal Z}_{0}^{j})_{j\in \{1,...,Q\}}$ et $\tilde{\mathcal Z}_{1}^{1}$ par (\ref{cond-b-1j1756})  et (\ref{cond-b-1j1757}). Alors les parties $\tilde{\mathcal Z}_{0}^{j}$ et $\tilde{\mathcal Z}_{1}^{1}$ sont  de diamètre inférieur ou égal à $P/3$  et il existe $\tilde Z\in \Lambda_{Z,t,(t_{1},...,t_{Q})}$ tel que 
$$
(\tilde a_{1},\dots,\tilde a_{p-1},\tilde S_{0},...,\tilde S_{m+1},(\tilde {\mathcal Y}_{i}^{j})_{i\in \{0,\dots,m+1\}, j\in \{1,\dots ,\tilde l_{i}\}},(\tilde {\mathcal Z}_{0}^{j})_{ j\in \{1,\dots ,Q\}},\tilde {\mathcal Z}_{1}^{1})$$
 appartienne  à 
 $(\pi_{x}^{\natural, p-1,k,m+1,(\tilde l_{0},...,\tilde l_{m+1}),Q,1})^{-1}(\tilde Z)$. 
 %
%$ Y_{x}^{\natural, p-1,k,m+1,(\tilde l_{0},...,\tilde l_{m+1}),Q,1}$. 
\end{souslem}
\noindent{\bf Démonstration.} 
C'est exactement le sous-lemme~\ref{slem0-4eme-2j1302}. 
%Pour $\sigma\in \{0,1\}$ et $z,z'\in \tilde{\mathcal Z}_{\sigma}^{j}$, on choisit $b\in \tilde S_{\sigma}$, d'où $z,z'\in 2N\tg(x,b)$ et comme $d(x,z)=d(x,z')$, $(H_{\de}(z,x,z',b))$ donne $d(z,z')\leq 2N+\de\leq P/3$. 
% Comme les parties $\mathcal Z_{0}^{j}$ et $\tilde{\mathcal Z}_{1}^{1}$  sont  non vides  l'argument que nous venons de donner montre aussi que la condition (\ref{def-mathcalZ-0j-24oct09}) ne dépend que de $\tilde Z$ (car $P/3\leq M$). 
 \cqfd

 \noindent {\bf Suite de la démonstration du b).}
Notre but est maintenant de montrer l'inégalité suivante, qui est une variante de (\ref{ineg-(1-P)-tilde-l0-bis}) : il existe une constante $C_{2}=C(\de,K,N,Q,P,M,r,T)$ telle  que 
 % $$\Big|\sum _{(a_{1},\dots,a_{p},S_{1},...,S_{m},(\mathcal Y_{i}^{j})_{i\in \{0,\dots,m\}, j\in \{1,\dots ,l_{i}\}}) \in 
%(\pi_{x}^{p,k,m,(l_{0},...,l_{m})})^{-1}(Z)} $$ $$(
%(u_{x,r,t}K_{x,Q,(t_{1},\dots ,t_{Q})}-
%  u_{x',r,t}K_{x',Q,(t_{1},\dots ,t_{Q})})  f)(a_1,...,a_p)\Big|^{2}$$  
 \begin{gather}\nonumber \sup_{\tau\in [0,T]}
 \big|\xi_{Z}\big(e^{\tau\theta_{x}}
(u_{x,r,t}K_{x,Q,(t_{1},\dots ,t_{Q})}-
 u_{x',r,t}K_{x',Q,(t_{1},\dots ,t_{Q})}) e^{-\tau\theta_{x}} f\big)\big|^{2} \\ \label{ineg-(1-P)-tilde-l0-bis-xx'}
 \leq C_{2}\sum_{\tilde Z\in \Lambda_{Z,t,(t_{1},...,t_{Q})}^{\neq}}  
 \big|\xi_{\tilde Z}(f)\big|^{2}.\end{gather}
%$$\Big| \sum _{\substack{(\tilde a_{1},\dots,\tilde a_{p-1},\tilde S_{1},...,\tilde S_{m+1},(\tilde {\mathcal Y}_{i}^{j})_{i\in \{0,\dots,m+1\}, j\in \{1,\dots ,\tilde l_{i}\}},\\ (\tilde {\mathcal Z}_{0}^{j})_{ j\in \{1,\dots ,Q\}},\tilde {\mathcal Z}_{1}^{1})  \in 
%(\pi_{x}^{\natural, p-1,k,m+1,(\tilde l_{0},...,\tilde l_{m+1}),Q,1})^{-1}(\tilde Z)}} f(\tilde a_{1},\dots,\tilde a_{p-1})\Big|^{2}$$
 où $\Lambda_{Z,t,(t_{1},...,t_{Q})}^{\neq}$ est l'ensemble des  $\tilde Z\in \Lambda_{Z,t,(t_{1},...,t_{Q})}$   
tels que pour tout  \begin{gather*}(\tilde a_{1},\dots,\tilde a_{p-1},\tilde S_{0},...,\tilde S_{m+1},(\tilde {\mathcal Y}_{i}^{j})_{i\in \{0,\dots,m+1\}, j\in \{1,\dots ,\tilde l_{i}\}},(\tilde {\mathcal Z}_{0}^{j})_{ j\in \{1,\dots ,Q\}},\tilde {\mathcal Z}_{1}^{1})\\ \in (\pi_{x}^{\natural, p-1,k,m+1,(\tilde l_{0},...,\tilde l_{m+1}),Q,1})^{-1}(\tilde Z)\end{gather*} 
 on ait 
 \begin{gather} \label{condition-uxrtKxqttt-xx'}\big( u_{x,r,t}K_{x,Q,(t_{1},\dots ,t_{Q})} - u_{x',r,t}K_{x',Q,(t_{1},\dots ,t_{Q})}\big) (e_{\tilde S_{0}})\neq 0
 .\end{gather}
%et 
%il existe $(a_{1},\pp,a_{p})$ vérifiant 
% $\tilde S_{1}=\{a_{1},\pp,a_{p}\}$ et 
% $$(a_{1},\pp,a_{p},\tilde S_{2},...,\tilde S_{m+1},(\tilde {\mathcal Y}_{i+1}^{j})_{i\in \{0,\dots,m\}, j\in \{1,\dots , l_{i}\}}) \in (\pi_{x}^{p,k,m,(l_{0},..., l_{m})})^{-1}(Z).$$ 

   Le sous-lemme suivant indique d'où vient  la condition  (\ref{r0-r1-r-xx'}). 
   
   \begin{souslem}\label{just-r0-r1-r-xx'-3j2157}
   Pour $S\in \Delta_{p-1}$ et $T\in \Delta_{p}$ tels que $e_{T}$ apparaisse avec un coefficient non nul dans  $u_{x,r,t}K_{x,Q,(t_{1},\dots ,t_{Q})} (e_{S})$ ou dans $u_{x',r,t}K_{x',Q,(t_{1},\dots ,t_{Q})} (e_{S})$
   on a 
   \begin{gather*} |d(x,S)-d(x,T)-r|\leq QF.\end{gather*}
   \end{souslem}
  \noindent{\bf Démonstration.} 
  D'après le 2)a) de la proposition~\ref{recap-supp-connaiss-H-uK} on a  \begin{gather*} T \subset \bigcup_{a\in S}\{y\in (F+2)\tg(x,a), d(y,a)\in [r+\frac{Q}{F},r+QF]\} \end{gather*}
 et  on suppose $ \frac{Q}{F}+QF\geq N+(F+2)$, ce qui est permis par $(H_{Q})$. \cqfd

\begin{souslem}\label{slem1-b-2j1827}
 Soit  $\tilde Z\in \Lambda_{Z,t,(t_{1},...,t_{Q})}^{\neq}$ et \begin{gather*}(\tilde a_{1},\dots,\tilde a_{p-1},\tilde S_{0},...,\tilde S_{m+1},(\tilde {\mathcal Y}_{i}^{j})_{i\in \{0,\dots,m+1\}, j\in \{1,\dots ,\tilde l_{i}\}},(\tilde {\mathcal Z}_{0}^{j})_{ j\in \{1,\dots ,Q\}},\tilde {\mathcal Z}_{1}^{1})\\ 
\in (\pi_{x}^{\natural, p-1,k,m+1,(\tilde l_{0},...,\tilde l_{m+1}),Q,1})^{-1}(\tilde Z).\end{gather*} Alors 
$u_{x,r,t} K_{x,Q,(t_{1},\dots ,t_{Q})}  (e_{\tilde S_{0}})$ et 
 $u_{x',r,t} K_{x',Q,(t_{1},\dots ,t_{Q})}  (e_{\tilde S_{0}})$ ne dépendent  que de la connaissance des points de 
 \begin{gather}\nonumber B(x,k+2M)\cup B(\tilde S_{0}, M)\cup B(\tilde S_{1},  M) \\ 
 \label{ens-tot-b-22dec}\cup  \bigcup _{ j\in \{1,\dots ,Q\}} B(\tilde{\mathcal Z}_{0}^{j}, M)
\cup  B(\tilde{\mathcal Z}_{1}^{1},M) \end{gather}
et des distances entre ces points. 
\end{souslem}
\noindent{\bf Démonstration.} 
 La réunion des ensembles figurant dans le 2)b) de la proposition~\ref{recap-supp-connaiss-H-uK} avec $x$ et $x'$ au lieu de $x$  et $\tilde S_{0}$ au lieu de $S$ (et qui a la propriété que $u_{x,r,t} K_{x,Q,(t_{1},\dots ,t_{Q})}  (e_{\tilde S_{0}})$ et 
 $u_{x',r,t} K_{x',Q,(t_{1},\dots ,t_{Q})}  (e_{\tilde S_{0}})$ ne dépendent  que de la connaissance des distances entre les points de cette réunion) 
  est inclus dans 
\begin{gather}\label{ens3.33-22dec-1}B(x,F+1)\cup B(\tilde S_{0}, QN) 
\\ \label{ens3.33-22dec-2}
\cup 
\bigcup_{a\in \tilde S_{0},j\in \{1,...,Q\}}\{y\in (F+2)\text{-}\geod(x,a), |d(x,y)-t_{j}d(x,a)|\leq QF+2\} 
\\ \label{ens3.33-22dec-3} \cup \bigcup_{a\in \tilde S_{0}} \{z\in (F+2)\text{-}\geod(x,a), d(z,a)\in [r,r+QF]\}
\\ \nonumber 
\cup \bigcup_{a\in \tilde S_{0}}\big\{z\in (F+2)\tg(x,a),\\  \label{ens3.33-22dec-4}  |d(x,z)-(1-t)(d(x,a)-r)|\leq QF+2\big\}.\end{gather}
 Il suffit donc de  montrer que cet ensemble est inclus dans (\ref{ens-tot-b-22dec}).

   On suppose $F+1\leq 2M$ et $QN\leq M$, ce qui est permis par $(H_{M})$. Donc (\ref{ens3.33-22dec-1}) est inclus dans (\ref{ens-tot-b-22dec}). 
   
Soit  $a\in \tilde S_{0}$, $j\in \{1,...,Q\}$ et  $y\in (F+2)\text{-}\geod(x,a)$ vérifiant  $$|d(x,y)-t_{j}d(x,a)|\leq QF+2.$$ Soit  $z\in \geod(x,a)$ vérifiant $d(x,z)=E(t_{j}r_{0}(\tilde Z))$, si bien que $z$ appartient à 
$\tilde{\mathcal Z}_{0}^{j}$. 
 On a 
$$|t_{j}d(x,a)-E(t_{j}r_{0}(\tilde Z))|\leq N+1,$$ d'où 
$|d(x,y)-d(x,z)|\leq QF+N+3$ et  grâce à $(H_{\de}(y,x,z,a))$, $$d(y,z)\leq 
(QF+N+3)+(F+2)+\de.$$ On suppose $(QF+N+3)+(F+2)+\de \leq M$, ce qui est permis par $(H_{M})$. Par conséquent  (\ref{ens3.33-22dec-2}) est inclus dans $\bigcup_{j\in\{1,\pp,Q\}}  B(\tilde{\mathcal Z}_{0}^{j},M) $ 
 et donc il est inclus dans (\ref{ens-tot-b-22dec}).

Par le a) du lemme~\ref{lemme-S0-...Sm}, on a pour tout $a\in \tilde S_{0}$, 
 $$\tilde S_{1}\subset 2F\tg(a,x). $$
Soit $a\in \tilde S_{0}$ et $y\in \tilde S_{1}$. L'inégalité (\ref{r0-r1-r-xx'}) implique 
$$|d(x,a)-d(x,y)-r|\leq QF+N. $$ 
Soit $z\in  (F+2)\tg(x,a)$ vérifiant $d(z,a)\in [r,r+QF]$. 
Cela implique $d(x,a)-d(x,z)\in [r-F-2,r+QF]$. 
On en déduit $$|d(x,y)-d(x,z)|\leq 2QF+N. $$
Comme $z\in  (F+2)\tg(x,a)$, $y\in  2F\tg(x,a)$ et $2F\geq F+2$, $(H_{\de}(z,x,y,a))$ montre 
$$d(y,z)\leq (2QF+N)+2F+\de$$
d'où $$
\bigcup_{a\in S} \{z\in (F+2)\text{-}\geod(x,a), d(z,a)\in [r,r+QF]\}$$
$$\subset B(\tilde S_{1}, 2QF+N+2F+\de)\subset B(\tilde S_{1}, M)$$
car on suppose $2QF+N+2F+\de \leq M$, ce qui est permis par $(H_{M})$. Donc (\ref{ens3.33-22dec-3}) est inclus dans (\ref{ens-tot-b-22dec}).

Enfin soit $a\in \tilde S_{0}$ et $z\in (F+2)\tg(x,a)$ vérifiant $$ |d(x,z)-(1-t)(d(x,a)-r)|\leq QF+2.$$ Soit $b\in \tilde S_{1}$ et $y\in \geod(x,b)$ vérifiant  $d(x,y)=E((1-t)r_{1}(\tilde Z))$, si bien que $y$ appartient à $\tilde{\mathcal Z}_{1}^{1}$. Comme $y\in \geod(x,b)$ et $b\in 2F\tg(x,a)$, on a $y\in 2F\tg(x,a)$. 
Comme $d(x,a)\in [r_{0}(\tilde Z),r_{0}(\tilde Z)+N]$ et grâce à (\ref{r0-r1-r}), on a 
$$|d(x,z)-(1-t)r_{1}(\tilde Z)|\leq 2QF+2+N,$$
d'où 
$$|d(x,y)-d(x,z)|\leq 2QF+3+N.$$
Comme $z\in (F+2)\tg(x,a)$, $y\in 2F\tg(x,a)$ et $F+2\leq 2F$, $(H_{\de}(y,x,z,a))$ montre 
$$d(y,z)\leq (2QF+3+N)+2F+\de.$$
On suppose $(2QF+3+N)+2F+\de\leq M$, ce qui est permis
 par $(H_{M})$.  Par conséquent   (\ref{ens3.33-22dec-4}) est inclus dans 
 $B(\tilde{\mathcal Z}_{1}^{1}, M
)$ et donc  dans 
 (\ref{ens-tot-b-22dec}). 
Ceci termine la preuve du sous-lemme~\ref{slem1-b-2j1827}. 
 %
%
%Donc  la réunion des ensembles figurant dans le 2)b) de la proposition~\ref{recap-supp-connaiss-H-uK} avec $x$ et $x'$ au lieu de $x$ et $\tilde S_{0}$ au lieu de $S$  (et qui a la propriété que $u_{x,r,t} K_{x,Q,(t_{1},\dots ,t_{Q})}  (e_{\tilde S_{0}})$ et 
% $u_{x',r,t} K_{x',Q,(t_{1},\dots ,t_{Q})}  (e_{\tilde S_{0}})$ ne dépendent  que de la connaissance des distances entre les points de cet ensemble) 
%  est inclus dans 
%(\ref{ens-tot-b-22dec}).
 \cqfd

Le sous-lemme~\ref{slem1-b-2j1827} justifie le fait que la condition (\ref{condition-uxrtKxqttt-xx'}) ne dépend que de $\tilde Z$ (c'est-à-dire que pour $\tilde Z\in \Lambda_{Z,t,(t_{1},...,t_{Q})}$ elle est vérifiée ou non simultanément pour tous les éléments de $(\pi_{x}^{\natural, p-1,k,m+1,(\tilde l_{0},...,\tilde l_{m+1}),Q,1})^{-1}(\tilde Z)$). 

 \begin{souslem}\label{slem2-b-2j1831}
 Le cardinal de $\Lambda_{Z,t,(t_{1},...,t_{Q})}$ est majoré par une constante de la forme $C(\de,K,N,Q,P,M)$. 
\end{souslem}
 \noindent{\bf Démonstration.} Grâce au lemme~\ref{nombre-dist-connaitre-par-point-natural}, pour connaître les distances entre les points de 
\begin{gather}\label{ens1-p107-24oct09bis}B(\tilde S_{0}, M)\cup   \bigcup _{ j\in \{1,\dots ,Q\}}B(\tilde {\mathcal Z}_{0}^{j}, M)\cup B(\tilde{\mathcal Z}_{1}^{1}, M)\end{gather} 
 et ceux de 
\begin{gather}\label{ens2-p107-24oct09bis}
\bigcup _{ i\in \{1,\dots ,m+1\}}
B(\tilde S_{i},M)
\cup  \bigcup _{i\in \{0,\dots,m\}, j\in \{1,\dots , l_{i}\}}B(\tilde {\mathcal Y}_{i+1}^{j},M) \cup B(x,k+2M)\end{gather}
il suffit de connaître les distances entre les points de (\ref{ens1-p107-24oct09bis})  et $C$ points de (\ref{ens2-p107-24oct09bis}), avec $C=C(\de,K,N,Q,P,M)$ et grâce à (\ref{r0-r1-r-xx'}) ces distances sont déterminées à $C'=C(\de,K,N,Q,P,M)$ près par les distances de $\tilde S_{1}$ à ces $C$ points (qui font partie de la donnée de $Z$) et les entiers $(t_{0}^{j}(\tilde Z))_{j\in \{1,\pp,Q\}}$ et $t_{1}^{1}(\tilde Z)$. \cqfd

\noindent {\bf Suite de la démonstration du b).}
   On termine  maintenant la preuve de l'inégalité (\ref{ineg-(1-P)-tilde-l0-bis-xx'}). 
 Pour \begin{gather*}(\tilde a_{1},\dots,\tilde a_{p-1},\tilde S_{0},...,\tilde S_{m+1},(\tilde {\mathcal Y}_{i}^{j})_{i\in \{0,\dots,m+1\}, j\in \{1,\dots ,\tilde l_{i}\}},(\tilde {\mathcal Z}_{0}^{j})_{ j\in \{1,\dots ,Q\}},\tilde {\mathcal Z}_{1}^{1})\\ 
\in (\pi_{x}^{\natural, p-1,k,m+1,(\tilde l_{0},...,\tilde l_{m+1}),Q,1})^{-1}(\tilde Z)\end{gather*} 
on considère 
\begin{gather}
\label{somme-1j1809}
\sum _{( b_{1},..., b_{p})}
\big((u_{x,r,t}K_{x,Q,(t_{1},\dots ,t_{Q})}-
 u_{x',r,t}K_{x',Q,(t_{1},\dots ,t_{Q})})   (e_{\tilde a_{1}}\wedge ...\wedge e_{\tilde a_{p-1}})\big)( b_{1},..., b_{p}), 
\end{gather} 
où la somme porte sur les énumérations $\tilde S_{1}=\{ b_{1},..., b_{p}\}$ telles que 
$$( b_{1},\dots, b_{p},\tilde S_{1},...,\tilde S_{m+1},(\tilde {\mathcal Y}_{i+1}^{j})_{i\in \{0,\dots,m\}, j\in \{1,\dots , l_{i}\}}
)\in (\pi_{x}^{p,k,m,(l_{0},..., l_{m})})^{-1}(Z).$$ 
Comme la somme (\ref{somme-1j1809}) a au plus $p!$ termes, 
 le 3) de la proposition~\ref{recap-supp-connaiss-H-uK} montre qu'elle 
est majorée par une constante de la forme $C(\de,K,N,Q,P,M)$. 
D'après le sous-lemme~\ref{slem1-b-2j1827} la somme (\ref{somme-1j1809}) ne dépend que de $\tilde Z$  et on peut donc la noter $\alpha_{Z,\tilde Z,t,(t_{1},...,t_{Q})}$. 
D'après les sous-lemmes~\ref{slem0-b-2j1820} et~\ref{just-r0-r1-r-xx'-3j2157} on a \begin{gather*}\xi_{Z}(e^{\tau\theta_{x}}(u_{x,r,t}K_{x,Q,(t_{1},\dots ,t_{Q})}-
 u_{x',r,t}K_{x',Q,(t_{1},\dots ,t_{Q})}) e^{-\tau\theta_{x}}  f)\\ =\frac{1}{(p-1)!}\sum_{\tilde Z\in 
 \Lambda_{Z,t,(t_{1},...,t_{Q})}^{\neq}}
\alpha_{Z,\tilde Z,t,(t_{1},...,t_{Q})}
e^{\tau(\rho_{x}^{1}(\tilde Z)-\rho_{x}^{0}(\tilde Z))}
\xi_{\tilde Z}(f). \end{gather*}
Grâce à (\ref{r0-r1-r-xx'}) on a $\rho_{x}^{1}(\tilde Z)-\rho_{x}^{0}(\tilde Z)\leq QF+N$. 
 Par Cauchy-Schwarz et grâce au sous-lemme~\ref{slem2-b-2j1831} 
 on en déduit (\ref{ineg-(1-P)-tilde-l0-bis-xx'}).

Montrons maintenant 
b)   à l'aide de (\ref{ineg-(1-P)-tilde-l0-bis-xx'}).  
Soit $C_{2}$ comme dans 
 (\ref{ineg-(1-P)-tilde-l0-bis-xx'}).  Soient $k,m,l_{0},\dots ,l_{m}\in \N $. On pose   $ \tilde l_{0}=0$ et $ \tilde l_{i}=l_{i-1}$ pour $i\in \{1,\pp,m+1\}$ comme précédemment. 
 
 \begin{souslem}\label{slem3-b-2j1836}
  Il existe $C_{3}= C(\de,K,N,Q,P,M,r)$ tel que pour tout 
$$(a_{1},\dots,a_{p},S_{0},...,S_{m},(\mathcal Y_{i}^{j})_{i\in \{0,\dots,m\}, j\in \{1,\dots ,l_{i}\}}) \in 
Y_{x}^{p,k,m,(l_{0},...,l_{m})}$$
 le nombre de possibilités pour $(\tilde a_{1},\dots,\tilde a_{p-1})$ tels que 
 $$(\tilde a_{1},\dots,\tilde a_{p-1},\{\tilde a_{1},\dots,\tilde a_{p-1}\}, S_{0},...,S_{m},
( {\mathcal Y}_{i-1}^{j})_{i\in \{0,\dots,m+1\}, j\in \{1,\dots ,\tilde l_{i}\}})$$  appartienne à $
Y_{x}^{p-1,k,m+1,(\tilde l_{0},...,\tilde l_{m+1})}
$ 
 et vérifie 
$$|d(x,\{\tilde a_{1},\dots,\tilde a_{p-1}\})-d(x,S_{0})-r|\leq QF$$ soit $\leq C_{3}$.  
\end{souslem}
 \noindent{\bf Démonstration.} C'est une conséquence immédiate du sous-lemme~\ref{slem3-b-2j1836-0}. \cqfd
 
 \noindent {\bf Fin de la démonstration du b).}
Soit  $Z\in \overline Y_{x}^{p,k,m,(l_{0},...,l_{m})}$  vérifiant $r_{0}(Z)> k+P$.
  Grâce aux sous-lemmes~\ref{mesure-t-uxrtK-xx'} et~\ref{slem3-b-2j1836}, on  a  
 \begin{gather*}\int_{(t,t_{1},...,t_{Q})\in [0,1]^{Q+1}}\Big(\sum_{\tilde Z\in  \Lambda_{Z,t,(t_{1},...,t_{Q})}^{\neq}}\sharp\big((\pi_{x}^{\natural,p-1,k,m+1,(\tilde l_{0},...,\tilde l_{m+1}),Q,1})^{-1}(\tilde Z)\big)\Big)dtdt_{1}...dt_{Q}\\ \leq \frac{C_{1}C_{3}}{r_{0}(Z)+1}\sharp\big((\pi_{x}^{p,k,m,(l_{0},..., l_{m})})^{-1}(Z)\big).\end{gather*}  
Notons $I_{Z}$ l'ensemble des $(t,t_{1},...,t_{Q})\in [0,1]^{Q+1}$ tels qu'il existe $\tilde Z\in  \Lambda_{Z,t,(t_{1},...,t_{Q})}^{\neq}$ vérifiant  
$$\sharp\big((\pi_{x}^{\natural,p-1,k,m+1,(\tilde l_{0},...,\tilde l_{m+1}),Q,1})^{-1}(\tilde Z)\big)\geq (r_{0}(Z)+1)^{-\frac{1}{2}}\sharp\big((\pi_{x}^{p,k,m,(l_{0},..., l_{m})})^{-1}(Z)\big). $$ La mesure de $I_{Z}$ est donc $\leq C_{1}C_{3}(r_{0}(Z)+1)^{-\frac{1}{2}}$. Grâce à Cauchy-Schwarz et à  (\ref{ineg-(1-P)-tilde-l0-bis-xx'}) on obtient 
que 
\begin{gather}\nonumber \sup_{\tau\in [0,T]}\Big|\xi_{Z}\Big(e^{\tau\theta_{x}}\Big(\int_{(t,t_{1},...,t_{Q}) \in I_{Z}}
(u_{x,r,t}K_{x,Q,(t_{1},\dots ,t_{Q})}
\\ \nonumber
- u_{x',r,t}K_{x',Q,(t_{1},\dots ,t_{Q})})dtdt_{1}...dt_{Q}\Big) e^{-\tau\theta_{x}} f\Big)\Big|^{2}
 \\ \nonumber 
 \leq C_{1}C_{3}(r_{0}(Z)+1)^{-\frac{1}{2}}\sup_{\tau\in [0,T]} \int_{(t,t_{1},...,t_{Q}) \in [0,1]^{Q+1}}\\ \nonumber 
\big|\xi_{Z}(e^{\tau\theta_{x}}
(u_{x,r,t}K_{x,Q,(t_{1},\dots ,t_{Q})}-
 u_{x',r,t}K_{x',Q,(t_{1},\dots ,t_{Q})}) e^{-\tau\theta_{x}} f)\big|^{2}dtdt_{1}...dt_{Q}\\ \label{ineg-IZ-b-xx'}  \leq 2^{Q+1}C_{1}C_{2}C_{3}(r_{0}(Z)+1)^{-\frac{1}{2}}
\sum_{\tilde Z\in \Lambda_{Z}} (r_{0}(\tilde Z)+1)^{-Q}(r_{1}(\tilde Z)+1)^{-1}
\big| \xi_{\tilde Z}(f)\big|^{2}.\end{gather}
 D'après le sous-lemme~\ref{slem3-b-2j1836} et le lemme~\ref{lemme-cardinaux}, 
 il existe $C_{4}=C(\de,K,N,Q,P,M,r)$ tel  que pour $\tilde Z\in \Lambda_{Z}$  on ait 
 \begin{gather}\label{ineg-C4-b-21dec}\sharp\big((\pi_{x}^{\natural,p-1,k,m+1,(\tilde l_{0},...,\tilde l_{m+1}),Q,1})^{-1}(\tilde Z)\big)\leq C_{4}\sharp\big((\pi_{x}^{p,k,m,(l_{0},..., l_{m})})^{-1}(Z)\big).\end{gather} 
 D'autre part pour $(t,t_{1},...,t_{Q})\not\in I_{Z}$ et $\tilde Z\in  \Lambda_{Z,t,(t_{1},...,t_{Q})}^{\neq}$  on a $$\sharp\big((\pi_{x}^{\natural,p-1,k,m+1,(\tilde l_{0},...,\tilde l_{m+1}),Q,1})^{-1}(\tilde Z)\big)\leq (r_{0}(Z)+1)^{-\frac{1}{2}}\sharp\big((\pi_{x}^{p,k,m,(l_{0},..., l_{m})})^{-1}(Z)\big). $$
  Par Cauchy-Schwarz on déduit alors de (\ref{ineg-(1-P)-tilde-l0-bis-xx'}) que 
  \begin{gather}\nonumber \sup_{\tau\in [0,T]}\big(\sharp(\pi_{x}^{p,k,m,(l_{0},..., l_{m})})^{-1}(Z)\big)^{-\alpha} \\ \nonumber
   \Big|\xi_{Z}\Big(e^{\tau\theta_{x}}\Big(\int_{(t,t_{1},...,t_{Q}) \not \in I_{Z}}
(u_{x,r,t}K_{x,Q,(t_{1},\dots ,t_{Q})}
\\ \nonumber
-
 u_{x',r,t}K_{x',Q,(t_{1},\dots ,t_{Q})}) dtdt_{1}...dt_{Q}\Big) e^{-\tau\theta_{x}}f\Big)\Big|^{2}
 \\ \nonumber 
 \leq 2^{Q+1} C_{2}(r_{0}(Z)+1)^{-\frac{\alpha}{2}}
\sum_{\tilde Z\in \Lambda_{Z}} (r_{0}(\tilde Z)+1)^{-Q}(r_{1}(\tilde Z)+1)^{-1}
\\ \label{ineg-hors-IZ-b-xx'}
 \big(\sharp(\pi_{x}^{\natural,p-1,k,m+1,(\tilde l_{0},...,\tilde l_{m+1}),Q,1})^{-1}(\tilde Z)\big)^{-\alpha}\big| \xi_{\tilde Z}(f)\big|^{2}.\end{gather}

  Comme     \begin{gather*}u_{x,r}K_{x}-u_{x',r}K_{x'}
=\int_{(t,t_{1},...,t_{Q})  \in I_{Z}}
(u_{x,r,t}K_{x,Q,(t_{1},\dots ,t_{Q})}-
 u_{x',r,t}K_{x',Q,(t_{1},\dots ,t_{Q})})\\ dtdt_{1}...dt_{Q}+
\int_{(t,t_{1},...,t_{Q}) \not \in I_{Z}}
(u_{x,r,t}K_{x,Q,(t_{1},\dots ,t_{Q})}-
 u_{x',r,t}K_{x',Q,(t_{1},\dots ,t_{Q})})dtdt_{1}...dt_{Q},\end{gather*}  en combinant les inégalités (\ref{ineg-IZ-b-xx'}), (\ref{ineg-C4-b-21dec}) et 
  (\ref{ineg-hors-IZ-b-xx'}) et par Cauchy-Schwarz on obtient que 
\begin{gather}\nonumber \sup_{\tau\in [0,T]}\big(\sharp(\pi_{x}^{p,k,m,(l_{0},..., l_{m})})^{-1}(Z)\big)^{-\alpha}|\xi_{Z}(e^{\tau\theta_{x}}(u_{x,r}K_{x}-u_{x',r}K_{x'})e^{-\tau\theta_{x}}f)|^{2}\\ \nonumber
%$$\Big|\sum _{(a_{1},\dots,a_{p},S_{1},...,S_{m},(\mathcal Y_{i}^{j})_{i\in \{0,\dots,m\}, j\in \{1,\dots ,l_{i}\}}) \in 
%(\pi_{x}^{p,k,m,(l_{0},...,l_{m})})^{-1}(Z)} ((h_{x}-h_{x'})  f)(a_1,...,a_p)\Big|^{2}$$ 
\leq 2^{Q+2}\Big(C_{1}C_{2}C_{3}C_{4}^{\alpha}(r_{0}(Z)+1)^{-\frac{1}{2}}+
 C_{2}(r_{0}(Z)+1)^{-\frac{\alpha}{2}}\Big)
\sum_{\tilde Z\in \Lambda_{Z}} (r_{0}(\tilde Z)+1)^{-Q}
\\ \label{C1234-3j2209}
(r_{1}(\tilde Z)+1)^{-1}
 \big(\sharp(\pi_{x}^{\natural,p-1,k,m+1,(\tilde l_{0},...,\tilde l_{m+1}),Q,1})^{-1}(\tilde Z)\big)^{-\alpha}\big| \xi_{\tilde Z}(f)\big|^{2}.
 \end{gather} 
    De plus pour $\tilde Z\in \Lambda_{Z}$ on a $\prod_{i=0}^{m}s_{i}(Z)^{-l_{i}}=\prod_{i=0}^{m+1}s_{i}(\tilde Z)^{-\tilde l_{i}}$. 
   Pour calculer la norme de $(1-\P_{n})e^{\tau\theta_{x}}(u_{x,r}K_{x}-u_{x',r}K_{x'})e^{-\tau\theta_{x}}$ on peut se limiter aux $Z$ tels que $r_{0}(Z)\geq n$ et on déduit donc de (\ref{C1234-3j2209}) que 
\begin{gather*}\sup_{\tau\in [0,T]}\|(1-\P_{n})e^{\tau\theta_{x}}(u_{x,r}K_{x}-u_{x',r}K_{x'}) e^{-\tau\theta_{x}}
\|_{\L(\H_{x,s}^{\natural,Q,1}(\Delta_{p-1}),\H^{\rightarrow}_{x,s}(\Delta_{p}))}^{2}\\ \leq  2^{Q+2}p! B e^{2(QF-r)s}
\Big(C_{1}C_{2}C_{3}C_{4}^{\alpha}(n+1)^{-\frac{1}{2}}+
 C_{2}(n+1)^{-\frac{\alpha}{2}}\Big) \end{gather*}
où le facteur $p!$ est dû au fait que $\tilde Z$ détermine $Z$ à permutation près de $a_{1},...,a_{p}$. 
Ceci termine la preuve de b).

\noindent {\bf Démonstration du c) du  lemme~\ref{lemme-compacite-equiv3}.} La preuve de c) est quasiment identique à celle de d) et légèrement plus simple (du fait que  $\del$ ne  fait  pas intervenir une  moyenne, contrairement à  $h_{x}$). Nous choisissons donc de montrer d) seulement. 

\noindent {\bf Démonstration du d) du  lemme~\ref{lemme-compacite-equiv3}.}
Soient $k,m,l_{0},\dots ,l_{m}\in \N $ et $$Z\in \overline Y_{x}^{p,k,m,(l_{0},...,l_{m})}\text{\ \  vérifiant \ \ }r_{0}(Z)> k+P.$$ 
On pose $ \tilde l_{0}=0$ et $ \tilde l_{i}=l_{i-1}$ pour $i\in \{1,\pp,m+1\}$.  
 On note $\Lambda_{Z}$ la partie de 
  $\overline Y_{x}^{\natural,p-1,k,m+1,(\tilde l_{0},...,\tilde l_{m+1}),7,0}$ formée des $\tilde Z$ vérifiant \begin{gather}\label{cond-d-3j2213} r_{0}( Z) \leq r_{0}(\tilde Z)\leq  r_{0}(Z)+N\end{gather}  et 
tels que pour tout  \begin{gather*}(\tilde a_{1},\dots,\tilde a_{p-1},\tilde S_{0},...,\tilde S_{m+1},(\tilde {\mathcal Y}_{i}^{j})_{i\in \{0,\dots,m+1\}, j\in \{1,\dots ,\tilde l_{i}\}},(\tilde {\mathcal Z}_{0}^{j})_{j\in \{1,\dots ,7\}}
)\\ \in (\pi_{x}^{\natural,p-1,k,m+1,(\tilde l_{0},...,\tilde l_{m+1}),7,0})^{-1}(\tilde Z)\end{gather*} 
il existe une énumération $(a_{1},\pp,a_{p})$ de $\tilde S_{1}$ vérifiant 
%   il existe $(a_{1},\pp,a_{p})$ (vérifiant 
% $\tilde S_{1}=\{a_{1},\pp,a_{p}\}$) tel que  
 \begin{gather}\label{cond-Z-d-1j2017}(a_{1},\pp,a_{p},\tilde S_{1},...,\tilde S_{m+1},(\tilde {\mathcal Y}_{i+1}^{j})_{i\in \{0,\dots,m\}, j\in \{1,\dots , l_{i}\}})\in (\pi_{x}^{p,k,m,(l_{0},..., l_{m})})^{-1}(Z).\end{gather} 

Soient $t\in [0,1]$ et  $u_{1},u_{2},u_{3},v_{1},v_{2},v_{3}\in [0,1[$. 
On note $(\Lambda_{Z,t})_{u_{1},u_{2},u_{3}}^{v_{1},v_{2},v_{3}}$  l'ensemble des 
$\tilde Z\in \Lambda_{Z}$  tels que, en notant \begin{gather*}w_{1}=E(\frac{u_{1}}{6}r_{0}(\tilde Z)),
 \ \ w_{2}=E((\frac{1}{6}+\frac{u_{2}}{6})r_{0}(\tilde Z)),
 \ \ w_{3}=E((\frac{2}{6}+\frac{u_{3}}{6})r_{0}(\tilde Z)),\\ w_{4}=E((1-\frac{v_{1}}{6})r_{0}(\tilde Z)),
   \ \ w_{5}=E((\frac{5}{6}-\frac{v_{2}}{6})r_{0}(\tilde Z)),\ \
    w_{6}=E((\frac{4}{6}-\frac{v_{3}}{6})r_{0}(\tilde Z)) \end{gather*}  
  on ait,   pour tout  \begin{gather*}(\tilde a_{1},\dots,\tilde a_{p-1},\tilde S_{0},...,\tilde S_{m+1},(\tilde {\mathcal Y}_{i}^{j})_{i\in \{0,\dots,m+1\}, j\in \{1,\dots ,\tilde l_{i}\}},(\tilde {\mathcal Z}_{0}^{j})_{j\in \{1,\dots ,7\}}
) \\ \in (\pi_{x}^{\natural,p-1,k,m+1,(\tilde l_{0},...,\tilde l_{m+1}),7,0})^{-1}(\tilde Z),\end{gather*} les égalités 
 \begin{gather} \label{cond-d-1j1951} \tilde{\mathcal Z}_{0}^{j}=\bigcup_{b\in \tilde S_{0}}\{z\in \geod(x,b), d(x,z)=w_{j}\}\text{ pour }j\in \{1,...,6\}
 \\ \label{cond-d-1j1952} \text{\ \ \ et\ \ \ }
 \tilde{\mathcal Z}_{0}^{7}=
\bigcup_{b\in \tilde S_{0}}\{z\in \geod(x,b), d(x,z)=E(tr_{0}(\tilde Z))\}.\end{gather}

Pour $\tilde Z\in (\Lambda_{Z,t})_{u_{1},u_{2},u_{3}}^{v_{1},v_{2},v_{3}}$, les conditions (\ref{cond-d-1j1951}) et (\ref{cond-d-1j1952}) impliquent 
$$t_{0}^{j}(\tilde Z)=w_{j}\text{\ \ pour\ \ }j\in \{1,...,6\}, \text{\ \ et \ \ } t_{0}^{7}(\tilde Z) =E(tr_{0}(\tilde Z)).$$ 

 \begin{souslem}\label{slem1-d-2j1839} Soit  $$(\tilde a_{1},\dots,\tilde a_{p-1},\tilde S_{0},...,\tilde S_{m+1},(\tilde {\mathcal Y}_{i}^{j})_{i\in \{0,\dots,m+1\}, j\in \{1,\dots ,\tilde l_{i}\}})\in Y_{x}^{p-1,k,m+1,(\tilde l_{0},...,\tilde l_{m+1})}$$
  tel que  $d(x,\tilde S_{1})\leq d(x,\tilde S_{0})\leq d(x,\tilde S_{1})+N$
   et 
   qu'il existe $(a_{1},\pp,a_{p})$ vérifiant  (\ref{cond-Z-d-1j2017}). 
 On définit $(\tilde{\mathcal Z}_{0}^{j})_{j\in \{1,...,7\}}$  par (\ref{cond-d-1j1951}) et (\ref{cond-d-1j1952}). Alors les parties $\tilde{\mathcal Z}_{0}^{j}$  sont  de diamètre inférieur ou égal à $P/3$  et il existe 
 $\tilde Z\in 
  (\Lambda_{Z,t})_{u_{1},u_{2},u_{3}}^{v_{1},v_{2},v_{3}}$ tel que 
\begin{gather}\label{elem-slem1-d-11j}
(\tilde a_{1},\dots,\tilde a_{p-1},\tilde S_{0},...,\tilde S_{m+1},(\tilde {\mathcal Y}_{i}^{j})_{i\in \{0,\dots,m+1\}, j\in \{1,\dots ,\tilde l_{i}\}},(\tilde {\mathcal Z}_{0}^{j})_{ j\in \{1,\dots ,7\}})\end{gather}
 appartienne  à 
 $(\pi_{x}^{\natural,p-1,k,m+1,(\tilde l_{0},...,\tilde l_{m+1}),7,0})^{-1}(\tilde Z)$.
%$ Y_{x}^{\natural, p-1,k,m+1,(\tilde l_{0},...,\tilde l_{m+1}),7,0}$. 
\end{souslem}
\noindent{\bf Démonstration.} 
En effet soit $b\in \tilde S_{0}$. Pour  $z,z'\in \tilde{\mathcal Z}_{0}^{j}$, on a  $z,z'\in 2N\tg(x,b)$ et comme $d(x,z)=d(x,z')$, $(H_{\de}(z,x,z',b))$ donne $d(z,z')\leq 2N+\de\leq P/3$. 
 Comme les parties $\tilde{\mathcal Z}_{0}^{j}$   sont  non vides et $P/3\leq M$,   l'argument que nous venons de donner montre aussi que les conditions (\ref{cond-d-1j1951}) et (\ref{cond-d-1j1952}) 
 sont vérifiées par les autres éléments de la classe d'équivalence $\tilde Z$ de l'élément 
 (\ref{elem-slem1-d-11j})
 et donc  $\tilde Z\in 
  (\Lambda_{Z,t})_{u_{1},u_{2},u_{3}}^{v_{1},v_{2},v_{3}}$.  \cqfd

\noindent {\bf Suite de la démonstration du d).}
On déduira d)  de  l'inégalité suivante : il existe $C_{2}=C(\de,K,N,Q,P,M,T)$ tel que 
   \begin{gather}\nonumber \sup_{\tau\in [0,T]}\big|\xi_{Z}([(\theta^{\flat}_{x})_{u_{1},u_{2},u_{3}}^{v_{1},v_{2},v_{3}}-
(\theta^{\flat}_{x'})_{u_{1},u_{2},u_{3}}^{v_{1},v_{2},v_{3}}, e^{\tau\theta_{x}}h_{x,t}e^{-\tau\theta_{x}}]  f)\big|^{2}
\\ 
\label{ineg-(1-P)-tilde-l0-theta}
\leq C_{2}\sum_{\tilde Z\in (\Lambda_{Z,t}^{\neq})_{u_{1},u_{2},u_{3}}^{v_{1},v_{2},v_{3}}}   
 \big| \xi_{\tilde Z}(  f)\big|^{2}\end{gather}
 où $(\Lambda_{Z,t}^{\neq})_{u_{1},u_{2},u_{3}}^{v_{1},v_{2},v_{3}}$ est l'ensemble des  $\tilde Z\in (\Lambda_{Z,t})_{u_{1},u_{2},u_{3}}^{v_{1},v_{2},v_{3}}$  tels que
   pour tout  \begin{gather*}(\tilde a_{1},\dots,\tilde a_{p-1},\tilde S_{0},...,\tilde S_{m+1},(\tilde {\mathcal Y}_{i}^{j})_{i\in \{0,\dots,m+1\}, j\in \{1,\dots ,\tilde l_{i}\}},(\tilde {\mathcal Z}_{0}^{j})_{j\in \{1,\dots ,7\}}
) \\ \in (\pi_{x}^{\natural,p-1,k,m+1,(\tilde l_{0},...,\tilde l_{m+1}),7,0})^{-1}(\tilde Z)\end{gather*} on ait  
 \begin{gather} \nonumber 
(\rho^{\flat}_{x})_{u_{1},u_{2},u_{3}}^{v_{1},v_{2},v_{3}}(\tilde S_{0})-
(\rho^{\flat}_{x'})_{u_{1},u_{2},u_{3}}^{v_{1},v_{2},v_{3}}(\tilde S_{0}) \\ \label{condition-rho-uuuvvvxx'01}-(\rho^{\flat}_{x})_{u_{1},u_{2},u_{3}}^{v_{1},v_{2},v_{3}}(\tilde S_{1})+(\rho^{\flat}_{x'})_{u_{1},u_{2},u_{3}}^{v_{1},v_{2},v_{3}}(\tilde S_{1})\neq 0.
\end{gather}

Nous allons maintenant  montrer  (\ref{ineg-(1-P)-tilde-l0-theta})
 et justifier que la condition (\ref{condition-rho-uuuvvvxx'01}) ne dépend  que de $\tilde Z$ (c'est-à-dire que pour $\tilde Z\in 
  (\Lambda_{Z,t})_{u_{1},u_{2},u_{3}}^{v_{1},v_{2},v_{3}}$ elle est vérifiée ou non simultanément pour tous les éléments de $(\pi_{x}^{\natural,p-1,k,m+1,(\tilde l_{0},...,\tilde l_{m+1}),7,0})^{-1}(\tilde Z)$). 
 
 \begin{souslem}\label{sslem-d-1j2043} 
 Pour 
 $\tilde Z\in (\Lambda_{Z,t})_{u_{1},u_{2},u_{3}}^{v_{1},v_{2},v_{3}}$ et  
 \begin{gather*}(\tilde a_{1},\dots,\tilde a_{p-1},\tilde S_{0},...,\tilde S_{m+1},(\tilde {\mathcal Y}_{i}^{j})_{i\in \{0,\dots,m+1\}, j\in \{1,\dots ,\tilde l_{i}\}},(\tilde {\mathcal Z}_{0}^{j})_{j\in \{1,\dots ,7\}}
) \\ \in (\pi_{x}^{\natural,p-1,k,m+1,(\tilde l_{0},...,\tilde l_{m+1}),7,0})^{-1}(\tilde Z)\end{gather*} le coefficient de $e_{\tilde S_{1}}$ dans 
$[(\theta^{\flat}_{x})_{u_{1},u_{2},u_{3}}^{v_{1},v_{2},v_{3}}-
(\theta^{\flat}_{x'})_{u_{1},u_{2},u_{3}}^{v_{1},v_{2},v_{3}}, e^{\tau\theta_{x}}h_{x,t}e^{-\tau\theta_{x}}] (e_{\tilde S_{0}})$ et le membre de gauche de (\ref{condition-rho-uuuvvvxx'01}) ne dépendent  que de 
la connaissance des points de 
\begin{gather}\label{ens-total-d-21dec} B( \tilde S_{0}, M)\cup  B( \tilde S_{1}, M)\cup B(x,k+2M) \cup \bigcup_{j\in \{1,...,7\}}B(\tilde{\mathcal Z}_{0}^{j},M) \end{gather} et des distances entre ces points. 
\end{souslem}
\noindent{\bf Démonstration.} 
Le sous-lemme~\ref{slem1-a-2j1541} 
(que l'on applique avec $\tilde{\mathcal Z}_{0}^{7}$ au lieu de $\tilde{\mathcal Z}_{0}^{1}$ et en oubliant $\tilde{\mathcal Z}_{0}^{j}$ pour $j\in \{1,...,6\}$)
montre que $h_{x,t}(e_{\tilde S_{0}})$ ne dépend que de la connaissance des points de  $$B(\tilde S_{0}, M)\cup B(x,k+2M)\cup B(\tilde{\mathcal Z}_{0}^{7}, M).$$  
Pour montrer le sous-lemme,  il  suffit donc de  montrer que pour $\tilde x\in B(x,1)$, 
 $(\rho^{\flat}_{\tilde x})_{u_{1},u_{2},u_{3}}^{v_{1},v_{2},v_{3}}(\tilde S_{0})$ et $(\rho^{\flat}_{\tilde x})_{u_{1},u_{2},u_{3}}^{v_{1},v_{2},v_{3}}(\tilde S_{1})$ ne dépendent que de la connaissance des points de (\ref{ens-total-d-21dec})  et des distances entre ces points (en effet on applique ceci à $\tilde x=x$ et $\tilde x=x'$). En vertu du lemme~\ref{dependance-rho'-123123}, 
il suffit de montrer que pour $\tilde x\in B(x,1)$, $\sigma\in \{0,1\}$, $a\in \tilde S_{\sigma}$ et $j\in \{1,...,6\}$,  l'ensemble 
  \begin{gather}\label{ens-30dec1211}  \{y\in 3\de\tg(\tilde x,a), |d(\tilde x,y)-w_{j}(\sigma,\tilde x)|\leq N+6\de+4\}\end{gather}
  est inclus dans (\ref{ens-total-d-21dec}), 
    où l'on note 
  \begin{gather}\nonumber  
  w_{1}(\sigma,\tilde x)=E(\frac{u_{1}}{6}d(\tilde x,\tilde S_{\sigma})),
  w_{2}(\sigma,\tilde x)=E((\frac{1}{6}+\frac{u_{2}}{6})d(\tilde x,
  \tilde S_{\sigma})), \\ \nonumber
  w_{3}(\sigma,\tilde x)=E((\frac{2}{6}+\frac{u_{3}}{6})d(\tilde x,\tilde S_{\sigma})),   w_{4}(\sigma,\tilde x)=E((1-\frac{v_{1}}{6})d(\tilde x,\tilde S_{\sigma})), \\ \label{def-wj-30dec1028} 
     w_{5}(\sigma,\tilde x)=E((\frac{5}{6}-\frac{v_{2}}{6})d(\tilde x,\tilde S_{\sigma})), 
    w_{6}(\sigma,\tilde x)=E((\frac{4}{6}-\frac{v_{3}}{6})d(\tilde x,\tilde S_{\sigma})). \end{gather}
       Soit $\tilde x\in B(x,1)$, $\sigma\in \{0,1\}$, $a\in \tilde S_{\sigma}$ et $j\in \{1,...,6\}$.

  Soit $b\in \tilde S_{0}$. On a \begin{gather}\label{ineg-dab2N2F-21dec}d(a,b)\leq 2N+2F.\end{gather}
  En effet c'est évident si $\sigma=0$ et si $\sigma=1 $ cela résulte du fait  que 
  \begin{gather}\label{ineg-30dec1249}d(x,\tilde S_{1})=r_{0}(Z)\in [d(x,\tilde S_{0})-N,d(x,\tilde S_{0})] \\ \nonumber 
  \subset [d(x,b)-2N,d(x,b)]
 \end{gather} et $\tilde S_{1}\subset 2F\tg(x,b)$. 
 
 Soit $y$ dans l'ensemble (\ref{ens-30dec1211}). 
 Comme $y\in 3\de\tg(\tilde x,a)$ et $\tilde x\in B(x,1)$ et grâce à (\ref{ineg-dab2N2F-21dec}), le lemme~\ref{xx'yy'zz'} 
 montre 
  que 
 \begin{gather}\label{incl-30dec1250}y\in (3\de+4N+4F+2)\tg(x,b).\end{gather} D'autre part il résulte de (\ref{ineg-30dec1249}) que 
 $|w_{j}(\sigma,\tilde x)-w_{j}|\leq N+1$. On a donc 
 \begin{gather}\label{ineg-30dec1251}
 |d(x,y)-w_{j}|\leq 2N+6\de+6.\end{gather}
 Soit $z\in \geod(x,b)$ vérifiant $d(x,z)=w_{j}$, si bien que $z$ appartient à 
 $\tilde{\mathcal Z}_{0}^{j}$. 
Comme  $y$ et $z$ appartiennent à $(3\de+4N+4F+2)\tg(x,b)$ et 
  que $$|d(x,y)-d(x,z)|\leq 2N+6\de+6,$$ $(H_{\de}(y,x,z,b))$ implique 
  $$d(y,z)\leq (2N+6\de+6)+(3\de+4N+4F+2)+\de.$$ On suppose 
  $(2N+6\de+6)+(3\de+4N+4F+2)+\de\leq M$, ce qui est permis par $(H_{M})$. Donc $d(y,z)\leq M$ et 
    $y$ appartient à (\ref{ens-total-d-21dec}). Ceci termine la preuve du sous-lemme~\ref{sslem-d-1j2043}. \cqfd
  
  Le sous-lemme~\ref{sslem-d-1j2043} implique immédiatement que la condition (\ref{condition-rho-uuuvvvxx'01}) ne dépend  que de $\tilde Z$. 
  
  \begin{souslem}\label{slem-2-d-2j1843}
   Le cardinal de $ (\Lambda_{Z,t}^{\neq})_{u_{1},u_{2},u_{3}}^{v_{1},v_{2},v_{3}}$ est majoré par une constante de la forme $C(\de,K,N,Q,P,M)$. 
     \end{souslem}
     \noindent{\bf Démonstration.} 
 %Il résulte de (\ref{ineg-dab2N2F-21dec}) que pour $\tilde S_{0},\tilde S_{1}$ comme ci-dessus on a 
%$d(\tilde S_{0},\tilde S_{1})\leq 2N+2F$. 
Grâce au lemme~\ref{nombre-dist-connaitre-par-point-natural}, pour connaître les distances entre les points de  
\begin{gather}\label{23dec2129}B(\tilde S_{0}, M)\cup  \bigcup_{j\in \{1,...,7\}}B(\tilde{\mathcal Z}_{0}^{j}, M)\end{gather} 
 et ceux de 
\begin{gather}\label{23dec2130} 
\bigcup _{ i\in \{1,\dots ,m+1\}}
B(\tilde S_{i},M)
\cup  \bigcup _{i\in \{0,\dots,m\}, j\in \{1,\dots , l_{i}\}}B(\tilde {\mathcal Y}_{i+1}^{j},M) \cup B(x,k+2M)\end{gather}
il suffit de connaître les distances entre les points de (\ref{23dec2129})  et $C$ points de (\ref{23dec2130}), avec $C=C(\de,K,N,Q,P,M)$ et grâce à (\ref{cond-d-3j2213}) ces distances sont déterminées à $C'=C(\de,K,N,Q,P,M)$ près par les distances de $\tilde S_{1}$ à ces $C$ points (qui font partie de la donnée de $Z$) et les entiers $(t_{0}^{j}(\tilde Z))_{j\in \{1,\pp,7\}}$. \cqfd

\noindent {\bf Fin de la démonstration du d).}
On termine maintenant   la preuve de (\ref{ineg-(1-P)-tilde-l0-theta}).
Pour  \begin{gather*}(\tilde a_{1},\dots,\tilde a_{p-1},\tilde S_{0},...,\tilde S_{m+1},(\tilde {\mathcal Y}_{i}^{j})_{i\in \{0,\dots,m+1\}, j\in \{1,\dots ,\tilde l_{i}\}},(\tilde {\mathcal Z}_{0}^{j})_{j\in \{1,\dots ,7\}}
) \\ \in (\pi_{x}^{\natural,p-1,k,m+1,(\tilde l_{0},...,\tilde l_{m+1}),7,0})^{-1}(\tilde Z)\end{gather*} 
on considère 
\begin{gather}
\label{somme-1j2059}
\sum _{( b_{1},..., b_{p})}
\big([(\theta^{\flat}_{x})_{u_{1},u_{2},u_{3}}^{v_{1},v_{2},v_{3}}-
(\theta^{\flat}_{x'})_{u_{1},u_{2},u_{3}}^{v_{1},v_{2},v_{3}}, e^{\tau\theta_{x}}h_{x,t}e^{-\tau\theta_{x}}]    (e_{\tilde a_{1}}\wedge ...\wedge e_{\tilde a_{p-1}})\big)( b_{1},..., b_{p}), 
\end{gather} 
où la somme porte sur les énumérations $\tilde S_{1}=\{ b_{1},..., b_{p}\}$ telles que 
$$( b_{1},\dots, b_{p},\tilde S_{1},...,\tilde S_{m+1},(\tilde {\mathcal Y}_{i+1}^{j})_{i\in \{0,\dots,m\}, j\in \{1,\dots ,\tilde l_{i}\}}
)\in (\pi_{x}^{p,k,m,(l_{0},..., l_{m})})^{-1}(Z).$$ 
Comme la somme (\ref{somme-1j2059}) a au plus $p!$ termes, 
 le 3) de la proposition~\ref{recap-supp-connaiss-H-uK} montre qu'elle 
est majorée par une constante de la forme $C(\de,K,N,Q,P,M,T)$. 
D'après le sous-lemme~\ref{sslem-d-1j2043} la somme (\ref{somme-1j2059}) ne dépend que de $\tilde Z$  et on peut donc la noter $(\alpha_{Z,\tilde Z,t,\tau})_{u_{1},u_{2},u_{3}}^{v_{1},v_{2},v_{3}}$. 
D'après le sous-lemme~\ref{slem1-d-2j1839} on a \begin{gather*}\xi_{Z}([(\theta^{\flat}_{x})_{u_{1},u_{2},u_{3}}^{v_{1},v_{2},v_{3}}-
(\theta^{\flat}_{x'})_{u_{1},u_{2},u_{3}}^{v_{1},v_{2},v_{3}}, e^{\tau\theta_{x}}h_{x,t}e^{-\tau\theta_{x}}]    f)\\ =\frac{1}{(p-1)!}\sum_{\tilde Z\in 
 (\Lambda_{Z,t}^{\neq})_{u_{1},u_{2},u_{3}}^{v_{1},v_{2},v_{3}}}
 (\alpha_{Z,\tilde Z,t,\tau})_{u_{1},u_{2},u_{3}}^{v_{1},v_{2},v_{3}}\ \xi_{\tilde Z}(f). \end{gather*}
 Par Cauchy-Schwarz et grâce au sous-lemme~\ref{slem-2-d-2j1843} 
 on en déduit  (\ref{ineg-(1-P)-tilde-l0-theta}).

Montrons maintenant d)  à l'aide de 
(\ref{ineg-(1-P)-tilde-l0-theta}).  Il résulte de (\ref{ineg-dab2N2F-21dec}) que pour $\tilde S_{0},\tilde S_{1}$ comme ci-dessus on a 
$d(\tilde S_{0},\tilde S_{1})\leq 2N+2F$. Donc 
le lemme~\ref{mesure-uuuvvv-dbemol} implique facilement qu'il existe $C_{1}=C(\de,K,N)$ tel que pour $\tilde S_{0},\tilde S_{1}$ comme ci-dessus  la mesure de l'ensemble des $(u_{1},u_{2}, u_{3},v_{1},v_{2},v_{3})\in [0,1[^{6}$ vérifiant (\ref{condition-rho-uuuvvvxx'01}) est $\leq \frac{C_{1}}{1+r_{0}(Z)}$. 
Soit $C_{2}$ comme dans 
 (\ref{ineg-(1-P)-tilde-l0-xx'}).  Soient $k,m$, $l_{0},\dots ,l_{m}\in \N $. On pose  $ \tilde l_{0}=0$ et $ \tilde l_{i}=l_{i-1}$ pour $i\in \{1,\pp,m+1\}$ comme précédemment. Soit $C_{3}= C(\de,K,N,Q,P,M)$ comme dans le sous-lemme~\ref{slem-a-3-2j1554}. %Il existe $C_{3}= C(\de,K,N,Q,P,M)$ tel que pour tout 
%$$(a_{1},\dots,a_{p},S_{1},...,S_{m},(\mathcal Y_{i}^{j})_{i\in \{0,\dots,m\}, j\in \{1,\dots ,l_{i}\}}) \in 
%Y_{x}^{p,k,m,(l_{0},...,l_{m})}$$
% le nombre de possibilités pour $(\tilde a_{1},\dots,\tilde a_{p-1})$ tels que 
% $$(\tilde a_{1},\dots,\tilde a_{p-1}, S_{0},...,S_{m},
%( {\mathcal Y}_{i-1}^{j})_{i\in \{0,\dots,m+1\}, j\in \{1,\dots ,\tilde l_{i}\}})$$  appartienne à $
%Y_{x}^{p-1,k,m+1,(\tilde l_{0},...,\tilde l_{m+1})}
%$ 
% et vérifie 
%$
%d(x,S_{0}) \leq d(x,\{\tilde a_{1},\dots,\tilde a_{p-1}\}) \leq  d(x,S_{0}) +N$ %soit $\leq C_{3}$.  
Soit  $$Z\in \overline Y_{x}^{p,k,m,(l_{0},...,l_{m})}\text{\ \  vérifiant \ \ }r_{0}(Z)> k+P.$$
  On  a alors 
 \begin{gather*}\int_{(u_{1},u_{2},u_{3},v_{1},v_{2},v_{3})\in [0,1[^{6}}\Big(\sum_{\tilde Z\in (\Lambda_{Z,t}^{\neq})_{u_{1},u_{2},u_{3}}^{v_{1},v_{2},v_{3}}}\sharp\big((\pi_{x}^{\natural,p-1,k,m+1,(\tilde l_{0},...,\tilde l_{m+1}),7,0})^{-1}(\tilde Z)\big)\Big)du_{1}...dv_{3}\\ \leq \frac{C_{1}C_{3}}{r_{0}(Z)+1}\sharp\big((\pi_{x}^{p,k,m,(l_{0},..., l_{m})})^{-1}(Z)\big).\end{gather*} 
 Notons $I_{Z}$ l'ensemble des $$(t,u_{1},u_{2}, u_{3},v_{1},v_{2},v_{3})\in [0,1]\times [0,1[^{6}$$ tels qu'il existe $\tilde Z\in (\Lambda_{Z,t}^{\neq})_{u_{1},u_{2},u_{3}}^{v_{1},v_{2},v_{3}}$ vérifiant 
 $$\sharp\big((\pi_{x}^{\natural,p-1,k,m+1,(\tilde l_{0},...,\tilde l_{m+1}),7,0})^{-1}(\tilde Z)\big)\geq (r_{0}(Z)+1)^{-\frac{1}{2}}\sharp\big((\pi_{x}^{p,k,m,(l_{0},..., l_{m})})^{-1}(Z)\big). $$ La mesure de $I_{Z}$ est donc $\leq C_{1}C_{3}(r_{0}(Z)+1)^{-\frac{1}{2}}$. Grâce à Cauchy-Schwarz on déduit de (\ref{ineg-(1-P)-tilde-l0-theta})
que 
\begin{gather}\nonumber
 \sup_{\tau\in [0,T]}\Big|\xi_{Z}\Big(\Big(\int_{(t,u_{1},u_{2}, u_{3},v_{1},v_{2},v_{3})\in  I_{Z}}[(\theta^{\flat}_{x})_{u_{1},u_{2},u_{3}}^{v_{1},v_{2},v_{3}}-
(\theta^{\flat}_{x'})_{u_{1},u_{2},u_{3}}^{v_{1},v_{2},v_{3}}, e^{\tau\theta_{x}}h_{x,t}e^{-\tau\theta_{x}}]
\\  \label{ineg-IZ-theta-xx'} dt...dv_{3}\Big)  f\Big)\Big|^{2}
\leq 2^{7}6^{6} C_{1}C_{2}C_{3}(r_{0}(Z)+1)^{-\frac{1}{2}}
\sum_{\tilde Z\in \Lambda_{Z}} (r_{0}(\tilde Z)+1)^{-7}
  \big| \xi_{\tilde Z} ( f)
\big|^{2}.\end{gather}
Grâce au   sous-lemme~\ref{slem-a-3-2j1554} et au lemme~\ref{lemme-cardinaux}, il existe $C_{4}=C(\de,K,N,Q,P,M)$ tel  que pour $\tilde Z\in \Lambda_{Z}$  on ait 
 \begin{gather}\label{ineg-C4-d-21dec}
 \sharp\big((\pi_{x}^{\natural,p-1,k,m+1,(\tilde l_{0},...,\tilde l_{m+1}),7,0})^{-1}(\tilde Z)\big)\leq C_{4}\sharp\big((\pi_{x}^{p,k,m,(l_{0},..., l_{m})})^{-1}(Z)\big).\end{gather} 
 D'autre part pour $(t,u_{1},u_{2}, u_{3},v_{1},v_{2},v_{3})\not\in I_{Z}$ et $\tilde Z\in (\Lambda_{Z,t}^{\neq})_{u_{1},u_{2},u_{3}}^{v_{1},v_{2},v_{3}}$ on a $$\sharp\big((\pi_{x}^{\natural,p-1,k,m+1,(\tilde l_{0},...,\tilde l_{m+1}),7,0})^{-1}(\tilde Z)\big)\leq (r_{0}(Z)+1)^{-\frac{1}{2}}\sharp\big((\pi_{x}^{p,k,m,(l_{0},..., l_{m})})^{-1}(Z)\big). $$
  Par Cauchy-Schwarz on déduit alors de (\ref{ineg-(1-P)-tilde-l0-theta}) que 
   \begin{gather}
   \nonumber 
   \big(\sharp(\pi_{x}^{p,k,m,(l_{0},..., l_{m})})^{-1}(Z)\big)^{-\alpha}\sup_{\tau\in [0,T]} 
   \\
   \nonumber
    \Big|\xi_{Z}\Big(\Big(\int_{(t,u_{1},u_{2}, u_{3},v_{1},v_{2},v_{3})\not \in  I_{Z}}[(\theta^{\flat}_{x})_{u_{1},u_{2},u_{3}}^{v_{1},v_{2},v_{3}}-
(\theta^{\flat}_{x'})_{u_{1},u_{2},u_{3}}^{v_{1},v_{2},v_{3}}, e^{\tau\theta_{x}}h_{x,t}e^{-\tau\theta_{x}}]dt...dv_{3}\Big)  f\Big)\Big|^{2}
\\
\nonumber
\leq 2^{7}6^{6}C_{2}(r_{0}(Z)+1)^{-\frac{\alpha}{2}}
\sum_{\tilde Z\in \Lambda_{Z}} (r_{0}(\tilde Z)+1)^{-7}
\\
\label{ineg-hors-IZ-theta-xx'}  \big(\sharp(\pi_{x}^{\natural,p-1,k,m+1,(\tilde l_{0},...,\tilde l_{m+1}),7,0})^{-1}(\tilde Z)\big)^{-\alpha} \big| \xi_{\tilde Z} ( f)
\big|^{2}.\end{gather}
 En combinant les inégalités (\ref{ineg-IZ-theta-xx'}), (\ref{ineg-C4-d-21dec})  et 
  (\ref{ineg-hors-IZ-theta-xx'}) et par Cauchy-Schwarz on obtient que 
\begin{gather}
\nonumber 
 \big(\sharp(\pi_{x}^{p,k,m,(l_{0},..., l_{m})})^{-1}(Z)\big)^{-\alpha}\sup_{\tau\in [0,T]} 
|\xi_{Z}(([(\theta^{\flat}_{x}-\theta^{\flat}_{x'}), e^{\tau\theta_{x}}h_{x}e^{-\tau\theta_{x}}])f)|^{2}\\ \nonumber
%$$\Big|\sum _{(a_{1},\dots,a_{p},S_{1},...,S_{m},(\mathcal Y_{i}^{j})_{i\in \{0,\dots,m\}, j\in \{1,\dots ,l_{i}\}}) \in 
%(\pi_{x}^{p,k,m,(l_{0},...,l_{m})})^{-1}(Z)} ((h_{x}-h_{x'})  f)(a_1,...,a_p)\Big|^{2}$$ 
\leq 2^{8}6^{6}\Big(C_{1}C_{2}C_{3}C_{4}^{\alpha}(r_{0}(Z)+1)^{-\frac{1}{2}}+
 C_{2}(r_{0}(Z)+1)^{-\frac{\alpha}{2}}\Big)
\sum_{\tilde Z\in \Lambda_{Z}} 
(r_{0}(\tilde  Z)+1)^{-7}\\ \label{C1234-p153-25oct09d}
\big(\sharp(\pi_{x}^{\natural,p-1,k,m+1,(\tilde l_{0},...,\tilde l_{m+1}),7,0})^{-1}(\tilde Z)\big)^{-\alpha}|\xi_{\tilde Z}(f)|^{2}.\end{gather} 
    De plus pour $\tilde Z\in \Lambda_{Z}$ on a $\prod_{i=0}^{m}s_{i}(Z)^{-l_{i}}=\prod_{i=0}^{m+1}s_{i}(\tilde Z)^{-\tilde l_{i}}$. 
   Pour calculer la norme de $(1-\P_{n})([(\theta^{\flat}_{x}-\theta^{\flat}_{x'}), e^{\tau\theta_{x}}h_{x}e^{-\tau\theta_{x}}])$ on peut se limiter aux $Z$ tels que $r_{0}(Z)\geq n$ et on déduit donc de (\ref{C1234-p153-25oct09d}) que 
\begin{gather*}\sup_{\tau\in [0,T]}\|(1-\P_{n})([(\theta^{\flat}_{x}-\theta^{\flat}_{x'}), e^{\tau\theta_{x}}h_{x}e^{-\tau\theta_{x}}]) 
\|_{\L(\H_{x,s}^{\natural,7,0}(\Delta_{p-1}),\H^{\rightarrow}_{x,s}(\Delta_{p}))}^{2}\\ \leq  2^{8}6^{6}p!
\Big(C_{1}C_{2}C_{3}C_{4}^{\alpha}(n+1)^{-\frac{1}{2}}+
 C_{2}(n+1)^{-\frac{\alpha}{2}}\Big) \end{gather*}
où le facteur $p!$ est dû au fait que $\tilde Z$ détermine $Z$ à permutation près de $a_{1},...,a_{p}$. Ceci termine la preuve de d). 
  
\noindent {\bf Démonstration du e) du  lemme~\ref{lemme-compacite-equiv3}.} La preuve de e) est plus subtile que celle de d) pour la raison suivante. L'entier $r\in \N$ est fixé mais peut être beaucoup plus grand que $M$. On ne peut donc pas espérer que, dans les notations de la preuve de d), 
$(\rho^{\flat}_{x})_{u_{1},u_{2},u_{3}}^{v_{1},v_{2},v_{3}}(e_{\tilde S_{0}})$ et 
$(\rho^{\flat}_{x})_{u_{1},u_{2},u_{3}}^{v_{1},v_{2},v_{3}}(e_{\tilde S_{1}})$ soit déterminés par la connaissance des points de 
$$B(\tilde S_{0}, M)\cup B(x,k+2M)\cup \bigcup_{j\in \{1,...,6\}}B(\tilde{\mathcal Z}_{0}^{j}, M) $$ et des distances entre ces points, pour certaines parties $\tilde{\mathcal Z}_{0}^{j}$. Au contraire  si des parties $(\tilde{\mathcal Z}_{1}^{j})_{j\in \{1,...,6\}}$ sont choisies de telle sorte qu'elles déterminent $(\rho^{\flat}_{x})_{u_{1},u_{2},u_{3}}^{v_{1},v_{2},v_{3}}(e_{\tilde S_{1}})$, on peut affirmer (en utilisant de nouveau  $(H_{M})$) qu'elles déterminent 
$(\rho^{\flat}_{x})_{\hat u_{1},\hat u_{2},\hat u_{3}}^{\hat v_{1},\hat v_{2},\hat v_{3}}(e_{\tilde S_{0}})$ pour certains $\hat u_{1},\hat u_{2},\hat u_{3}, \hat v_{1},\hat v_{2},\hat v_{3}$ que nous allons calculer. D'abord   $\tilde S_{1}$ est situé en gros (c'est-à-dire modulo des constantes de la forme 
$C(\de,K,N,Q,P)$) sur une géodésique entre $x$ et $\tilde S_{0}$, à distance $r$ de $\tilde S_{0}$. Pour que   les parties $(\tilde{\mathcal Z}_{1}^{j})_{j\in \{1,...,6\}}$ déterminent $(\rho^{\flat}_{x})_{u_{1},u_{2},u_{3}}^{v_{1},v_{2},v_{3}}(e_{\tilde S_{1}})$, elles doivent être situées en gros sur une géodésique entre $x$ et $\tilde S_{1}$, à des distances de $x$ égales à \begin{gather*}\frac{u_{1}}{6}d(x,\tilde S_{1}), \frac{1+u_{2}}{6}d(x,\tilde S_{1}),\frac{2+u_{3}}{6}d(x,\tilde S_{1}),\\ 
\frac{6-v_{1}}{6}d(x,\tilde S_{1}),\frac{5-v_{2}}{6}d(x,\tilde S_{1}),\frac{4-v_{3}}{6}d(x,\tilde S_{1}).\end{gather*}
Mais pour déterminer $(\rho^{\flat}_{x})_{\hat u_{1},\hat u_{2},\hat u_{3}}^{\hat v_{1},\hat v_{2},\hat v_{3}}(e_{\tilde S_{0}})$  les parties $(\tilde{\mathcal Z}_{1}^{j})_{j\in \{1,...,6\}}$ doivent également être situées en gros sur une géodésique entre $x$ et $\tilde S_{0}$, à des distances de $x$ égales à \begin{gather*}\frac{\hat u_{1}}{6}d(x,\tilde S_{0}), \frac{1+\hat u_{2}}{6}d(x,\tilde S_{0}),\frac{2+\hat u_{3}}{6}d(x,\tilde S_{0}),\\ 
\frac{6-\hat v_{1}}{6}d(x,\tilde S_{0}),\frac{5-\hat v_{2}}{6}d(x,\tilde S_{0}),\frac{4-\hat v_{3}}{6}d(x,\tilde S_{0}).\end{gather*}
On pose   $\kappa=\frac{r}{1+d(x,\tilde S_{1})}$ de sorte que 
$d(x,\tilde S_{0})$ est en gros égal à $(1+\kappa)d(x,\tilde S_{1})$. On obtient donc les relations 
\begin{gather}\nonumber \hat u_{1}=\frac{u_{1}}{1+\kappa},\ 
\hat u_{2}=\frac{u_{2}-\kappa}{1+\kappa}, 
\ \hat u_{3}=\frac{u_{3}-2\kappa}{1+\kappa},\\
\label{hatuuuvvv-uuuvvv}
\hat  v_{1}=\frac{v_{1}+6\kappa}{1+\kappa},
\ \hat  v_{2}=\frac{v_{2}+5\kappa}{1+\kappa}, \
\hat  v_{3}=\frac{v_{3}+4\kappa}{1+\kappa}.\end{gather}
Pour adapter la preuve de d) on utilisera  la variante suivante du lemme~\ref{mesure-uuuvvv-dbemol} (que l'on appliquera avec $\rho$ en gros égal à $r$, $\kappa$ comme ci-dessus, $y\in \tilde S_{1}$ et $y'\in \tilde S_{0}$).  

\begin{souslem}\label{mesure-uuuvvv-dbemol-e)}
 Pour tout $\rho\in \N$, il existe $C=C(\de,K,\rho)$ tel que pour  $x,x',y,y'\in X$ verifiant $d(x,x')\leq \rho$ et $d(y,y')\leq \rho$ et 
 pour $\kappa\in [0,\frac{1}{10}]$ vérifiant $\kappa d(x,y)\leq \rho$,   
  la mesure de l'ensemble des $(u_{1},u_{2},u_{3},v_{1},v_{2},v_{3})\in [5\kappa,1-5\kappa[^{6}$  tels que, en définissant  $\hat u_{1},\hat u_{2},\hat u_{3}, \hat v_{1},\hat v_{2},\hat v_{3}$
  comme dans (\ref{hatuuuvvv-uuuvvv}), 
 $${d^{\flat}}_{u_{1},u_{2},u_{3}}^{v_{1},v_{2},v_{3}}(x,y)
 -{d^{\flat}}_{u_{1},u_{2},u_{3}}^{v_{1},v_{2},v_{3}}(x',y)-{d^{\flat}}_{\hat u_{1},\hat u_{2},\hat u_{3}}^{\hat v_{1},\hat v_{2},\hat v_{3}}(x,y')
 + {d^{\flat}}_{\hat u_{1},\hat u_{2},\hat u_{3}}^{\hat v_{1},\hat v_{2},\hat v_{3}}(x',y')\neq 0$$ est 
 $\leq \frac{C}{1+d(x,y)}$.
\end{souslem}
\noindent {\bf Démonstration.} 
La démonstration  est une adaptation de celle du lemme~\ref{mesure-uuuvvv-dbemol}. En particulier on applique le sous-lemme~\ref{lemme-asssrrr} aux familles  \begin{gather*} {\tilde A}_{a_{1},a_{2},a_{3}}^{b_{1},b_{2},b_{3}}=A_{\tilde  x,u_{1},u_{2},u_{3}}^{\tilde  y,v_{1},v_{2},v_{3}}\text{\ \ si \ }\tilde y=y\\ \text{et\ \ 
}{\tilde A}_{a_{1},a_{2},a_{3}}^{b_{1},b_{2},b_{3}}=
A_{\tilde  x,\hat u_{1},\hat u_{2},\hat u_{3}}^{\tilde  y,\hat v_{1},\hat v_{2},\hat v_{3}}\text{\  \ 
si \ }\tilde y=y',\end{gather*} pour $a_{1},a_{2},a_{3}, b_{1},b_{2},b_{3}\in [0,1[$, avec 
$$u_{i}=5\kappa+(1-10\kappa)a_{i}\text{\  et \ }v_{i}=5\kappa+(1-10\kappa)b_{i}\text{\  pour \ }i=1,2,3$$
et  $\hat u_{1},\hat u_{2},\hat u_{3}, \hat v_{1},\hat v_{2},\hat v_{3}$
  comme dans (\ref{hatuuuvvv-uuuvvv}). 
\cqfd

\noindent {\bf Suite de la démonstration du e).}
Soient $k,m$, $l_{0},\dots ,l_{m}\in \N  $ et  $$Z\in \overline Y_{x}^{p,k,m,(l_{0},...,l_{m})}\text{\ \  vérifiant \ \ }r_{0}(Z)> k+P.$$ 
On pose $ \tilde l_{0}=0$ et $ \tilde l_{i}=l_{i-1}$ pour $i\in \{1,\pp,m+1\}$.  
 On note $\Lambda_{Z}$ 
la partie de  $ \overline Y_{x}^{\natural, p-1,k,m+1,(\tilde l_{0},...,\tilde l_{m+1}),Q,7}$ formée des $\tilde Z$ vérifiant 
  \begin{gather}\label{cond-r0r1-23dec}|r_{0}(\tilde Z)-r_{1}(\tilde Z)-r|\leq QF,\end{gather}
 et 
tels que pour tout  \begin{gather*}(\tilde a_{1},\dots,\tilde a_{p-1},\tilde S_{0},...,\tilde S_{m+1},(\tilde {\mathcal Y}_{i}^{j})_{i\in \{0,\dots,m+1\}, j\in \{1,\dots ,\tilde l_{i}\}},(\tilde {\mathcal Z}_{0}^{j})_{ j\in \{1,\dots ,Q\}},(\tilde {\mathcal Z}_{1}^{j})_{j\in \{1,\dots ,7\}})\\ \in (\pi_{x}^{\natural, p-1,k,m+1,(\tilde l_{0},...,\tilde l_{m+1}),Q,7})^{-1}(\tilde Z)\end{gather*} 
il existe une énumération $(a_{1},\pp,a_{p})$ de $\tilde S_{1}$ vérifiant 
 \begin{gather}\label{cond-Z-e-1j2114}\!\!\!\!\!\!\!(a_{1},\pp,a_{p},\tilde S_{1},...,\tilde S_{m+1},(\tilde {\mathcal Y}_{i+1}^{j})_{i\in \{0,\dots,m\}, j\in \{1,\dots , l_{i}\}})\in (\pi_{x}^{p,k,m,(l_{0},..., l_{m})})^{-1}(Z).\end{gather} 
 
On pose 
  $\kappa=\frac{r}{1+r_{0}(Z)}$.  Soient 
$t,t_{1},\pp,t_{Q}\in [0,1]$   et   $u_{1},u_{2}$, $u_{3},v_{1},v_{2},v_{3}\in [5\kappa,1-5\kappa[$. On note 
  $(\Lambda_{Z,t,(t_{1},...,t_{Q})})_{u_{1},u_{2},u_{3}}^{v_{1},v_{2},v_{3}}$  l'ensemble des 
$\tilde Z\in \Lambda_{Z}$  tels que, en notant 
\begin{gather*}w_{1}=E(\frac{u_{1}}{6}r_{1}(\tilde Z)),
 \ \ w_{2}=E((\frac{1}{6}+\frac{u_{2}}{6})r_{1}(\tilde Z)),
 \ \ w_{3}=E((\frac{2}{6}+\frac{u_{3}}{6})r_{1}(\tilde Z)),\\ w_{4}=E((1-\frac{v_{1}}{6})r_{1}(\tilde Z)),
   \ \ w_{5}=E((\frac{5}{6}-\frac{v_{2}}{6})r_{1}(\tilde Z)),\ \
    w_{6}=E((\frac{4}{6}-\frac{v_{3}}{6})r_{1}(\tilde Z)) \end{gather*}  
on ait, 
pour tout  \begin{gather*}(\tilde a_{1},\dots,\tilde a_{p-1},\tilde S_{0},...,\tilde S_{m+1},(\tilde {\mathcal Y}_{i}^{j})_{i\in \{0,\dots,m+1\}, j\in \{1,\dots ,\tilde l_{i}\}},(\tilde {\mathcal Z}_{0}^{j})_{ j\in \{1,\dots ,Q\}},(\tilde {\mathcal Z}_{1}^{j})_{j\in \{1,\dots ,7\}})\\ \in (\pi_{x}^{\natural, p-1,k,m+1,(\tilde l_{0},...,\tilde l_{m+1}),Q,7})^{-1}(\tilde Z)\end{gather*} 
les égalités 
  \begin{gather}
  \label{cond-Z0-e-1j2118}
   \tilde{\mathcal Z}_{0}^{j}=\bigcup_{b\in \tilde S_{0}}\{z\in \geod(x,b), d(x,z)=E(t_{j}r_{0}(\tilde Z))\}\text{ pour }j\in \{1,...,Q\},
  \\ 
  \label{cond-Z1-e-1j2118}
  \tilde{\mathcal Z}_{1}^{j}=\bigcup_{b\in \tilde S_{1}}\{z\in \geod(x,b), d(x,z)=w_{j}\}\text{ pour }j\in \{1,...,6\},
  \\ 
  \label{cond-Z1bis-e-1j2118}
 \text{et\ \ \ \ } \tilde{\mathcal Z}_{1}^{7}=\bigcup_{b\in \tilde S_{1}}\{z\in \geod(x,b), d(x,z)=E((1-t)r_{1}(\tilde Z))\}.
\end{gather}

Pour $\tilde Z\in (\Lambda_{Z,t,(t_{1},...,t_{Q})})_{u_{1},u_{2},u_{3}}^{v_{1},v_{2},v_{3}}$, les conditions (\ref{cond-Z0-e-1j2118}), (\ref{cond-Z1-e-1j2118}) et (\ref{cond-Z1bis-e-1j2118}) impliquent 
 \begin{gather*}t_{0}^{j}(\tilde Z)=E(t_{j}r_{0}(\tilde Z))\text{\ \ pour\ \ }j\in \{1,...,Q\}, \\ t_{1}^{j}(\tilde Z)=w_{j}\text{\ \ pour\ \ }j\in \{1,...,6\} 
 \text{\ \ \ \ et \ \ \ } 
   t_{1}^{7}(\tilde Z)
 =E((1-t)r_{1}(\tilde Z)).
 \end{gather*}
 
 \begin{souslem}\label{slem0-e-2j1938}  Soit  $$(\tilde a_{1},\dots,\tilde a_{p-1},\tilde S_{0},...,\tilde S_{m+1},(\tilde {\mathcal Y}_{i}^{j})_{i\in \{0,\dots,m+1\}, j\in \{1,\dots ,\tilde l_{i}\}})\in Y_{x}^{p-1,k,m+1,(\tilde l_{0},...,\tilde l_{m+1})}$$
  tel que  $|d(x,\tilde S_{0})-d(x,\tilde S_{1})-r|\leq QF
  $
   et 
   qu'il existe $(a_{1},\pp,a_{p})$ vérifiant 
 (\ref{cond-Z-e-1j2114}). 
 On définit $(\tilde{\mathcal Z}_{0}^{j})_{j\in \{1,...,Q\}}$ et $(\tilde{\mathcal Z}_{1}^{j})_{j\in \{1,...,7\}}$ par (\ref{cond-Z0-e-1j2118}), (\ref{cond-Z1-e-1j2118}) et (\ref{cond-Z1bis-e-1j2118}). Alors les parties $\tilde{\mathcal Z}_{0}^{j}$ et $\tilde{\mathcal Z}_{1}^{j}$ sont  de diamètre inférieur ou égal à $P/3$  et  il existe $\tilde Z\in  (\Lambda_{Z,t,(t_{1},...,t_{Q})})_{u_{1},u_{2},u_{3}}^{v_{1},v_{2},v_{3}}$ tel que 
\begin{gather}\label{elem-slem0-e-11j}
(\tilde a_{1},\dots,\tilde a_{p-1},\tilde S_{0},...,\tilde S_{m+1},(\tilde {\mathcal Y}_{i}^{j})_{i\in \{0,\dots,m+1\}, j\in \{1,\dots ,\tilde l_{i}\}},(\tilde{\mathcal Z}_{0}^{j})_{j\in \{1,...,Q\}}, (\tilde{\mathcal Z}_{1}^{j})_{j\in \{1,...,7\}})\end{gather}
 appartienne  à
  $(\pi_{x}^{\natural, p-1,k,m+1,(\tilde l_{0},...,\tilde l_{m+1}),Q,7})^{-1}(\tilde Z)$.
%$ Y_{x}^{\natural, p-1,k,m+1,(\tilde l_{0},...,\tilde l_{m+1}),Q,7}$. 
\end{souslem}
\noindent{\bf Démonstration. }
Pour $\sigma\in \{0,1\}$ et $z,z'\in \tilde{\mathcal Z}_{\sigma}^{j}$, on choisit $b\in \tilde S_{\sigma}$, d'où $z,z'\in 2N\tg(x,b)$ et comme $d(x,z)=d(x,z')$, $(H_{\de}(z,x,z',b))$ donne $d(z,z')\leq 2N+\de\leq P/3$. 
 Comme les parties $\tilde{\mathcal Z}_{0}^{j}$ et $\tilde{\mathcal Z}_{1}^{j}$   sont  non vides et $P/3\leq M$, l'argument que nous venons de donner montre aussi que les conditions (\ref{cond-Z0-e-1j2118}), (\ref{cond-Z1-e-1j2118}) et (\ref{cond-Z1bis-e-1j2118}) 
 sont vérifiées par les autres éléments de la classe d'équivalence 
   $\tilde Z$ de l'élément (\ref{elem-slem0-e-11j}) et donc $\tilde Z\in  (\Lambda_{Z,t,(t_{1},...,t_{Q})})_{u_{1},u_{2},u_{3}}^{v_{1},v_{2},v_{3}}$.  \cqfd

\noindent {\bf Suite de la démonstration du e).}
On déduira e)  de  l'inégalité suivante : il existe $C_{2}=C(\de,K,N,Q,P,M,r,T)$ tel que en notant $\hat u_{1},\hat u_{2},\hat u_{3}, \hat v_{1},\hat v_{2},\hat v_{3}$
  comme dans (\ref{hatuuuvvv-uuuvvv}), on ait 
   \begin{gather}
\nonumber 
\sup_{\tau\in [0,T]}\Big|\xi_{Z}\Big(\Big(\big((\theta^{\flat}_{x})_{u_{1},u_{2},u_{3}}^{v_{1},v_{2},v_{3}}-
(\theta^{\flat}_{x'})_{u_{1},u_{2},u_{3}}^{v_{1},v_{2},v_{3}}\big) e^{\tau\theta_{x}}u_{x,r,t}K_{x,Q,(t_{1},\dots ,t_{Q})}e^{-\tau\theta_{x}}-
\\
\nonumber 
e^{\tau\theta_{x}}u_{x,r,t}K_{x,Q,(t_{1},\dots ,t_{Q})}e^{-\tau\theta_{x}}
\big((\theta^{\flat}_{x})_{\hat u_{1},\hat u_{2},\hat u_{3}}^{\hat v_{1},\hat v_{2},\hat v_{3}}-
(\theta^{\flat}_{x'})_{\hat u_{1},\hat u_{2},\hat u_{3}}^{\hat v_{1},\hat v_{2},\hat v_{3}}\big)\Big)
  f\Big)\Big|^{2} 
\\
 \label{ineg-(1-P)-tilde-l0-bis-xx'-e)}
 \leq C_{2} \sum_{\tilde Z\in 
 (\Lambda_{Z,t,(t_{1},...,t_{Q})}^{\neq})_{u_{1},u_{2},u_{3}}^{v_{1},v_{2},v_{3}}}  
\big| \xi_{\tilde Z} (f)\big|^{2}\end{gather}
 où 
 $(\Lambda_{Z,t,(t_{1},...,t_{Q})}^{\neq})_{u_{1},u_{2},u_{3}}^{v_{1},v_{2},v_{3}}$
 est l'ensemble des  $\tilde Z\in (\Lambda_{Z,t,(t_{1},...,t_{Q})})_{u_{1},u_{2},u_{3}}^{v_{1},v_{2},v_{3}}$ 
 tels que pour tout  \begin{gather*}(\tilde a_{1},\dots,\tilde a_{p-1},\tilde S_{0},...,\tilde S_{m+1},(\tilde {\mathcal Y}_{i}^{j})_{i\in \{0,\dots,m+1\}, j\in \{1,\dots ,\tilde l_{i}\}},(\tilde {\mathcal Z}_{0}^{j})_{ j\in \{1,\dots ,Q\}},(\tilde {\mathcal Z}_{1}^{j})_{j\in \{1,\dots ,7\}})\\ \in (\pi_{x}^{\natural, p-1,k,m+1,(\tilde l_{0},...,\tilde l_{m+1}),Q,7})^{-1}(\tilde Z)\end{gather*} 
 on ait 
  \begin{gather}\nonumber 
(\rho^{\flat}_{x})_{\hat u_{1},\hat u_{2},\hat u_{3}}^{\hat v_{1},\hat v_{2},\hat v_{3}}(\tilde S_{0})-
(\rho^{\flat}_{x'})_{\hat u_{1},\hat u_{2},\hat u_{3}}^{\hat v_{1},\hat v_{2},\hat v_{3}}(\tilde S_{0})\\ \label{condition-rho-uuuvvvxx'01-e)}
-(\rho^{\flat}_{x})_{u_{1},u_{2},u_{3}}^{v_{1},v_{2},v_{3}}(\tilde S_{1})+(\rho^{\flat}_{x'})_{u_{1},u_{2},u_{3}}^{v_{1},v_{2},v_{3}}(\tilde S_{1})\neq 0. 
\end{gather}

Nous allons maintenant  montrer  (\ref{ineg-(1-P)-tilde-l0-bis-xx'-e)}). Nous allons 
  justifier aussi que la condition (\ref{condition-rho-uuuvvvxx'01-e)}) ne dépend  que de $\tilde Z$  (c'est-à-dire que pour $\tilde Z\in  (\Lambda_{Z,t,(t_{1},...,t_{Q})})_{u_{1},u_{2},u_{3}}^{v_{1},v_{2},v_{3}}$ elle est vérifiée ou non simultanément pour tous les éléments de \break $(\pi_{x}^{\natural, p-1,k,m+1,(\tilde l_{0},...,\tilde l_{m+1}),Q,7})^{-1}(\tilde Z)$).

 \begin{souslem}\label{sslem-e-1j2139} 
 Pour 
 $\tilde Z\in
 (\Lambda_{Z,t,(t_{1},...,t_{Q})})_{u_{1},u_{2},u_{3}}^{v_{1},v_{2},v_{3}}$
 et
 \begin{gather*}(\tilde a_{1},\dots,\tilde a_{p-1},\tilde S_{0},...,\tilde S_{m+1},(\tilde {\mathcal Y}_{i}^{j})_{i\in \{0,\dots,m+1\}, j\in \{1,\dots ,\tilde l_{i}\}},(\tilde {\mathcal Z}_{0}^{j})_{ j\in \{1,\dots ,Q\}},(\tilde {\mathcal Z}_{1}^{j})_{j\in \{1,\dots ,7\}}) \\ \in (\pi_{x}^{\natural, p-1,k,m+1,(\tilde l_{0},...,\tilde l_{m+1}),Q,7})^{-1}(\tilde Z) \end{gather*} le coefficient de $e_{\tilde S_{1}}$ dans 
\begin{gather*}\Big(\big((\theta^{\flat}_{x})_{u_{1},u_{2},u_{3}}^{v_{1},v_{2},v_{3}}-
(\theta^{\flat}_{x'})_{u_{1},u_{2},u_{3}}^{v_{1},v_{2},v_{3}}\big) e^{\tau\theta_{x}}u_{x,r,t}K_{x,Q,(t_{1},\dots ,t_{Q})}e^{-\tau\theta_{x}}-
\\ 
e^{\tau\theta_{x}}u_{x,r,t}K_{x,Q,(t_{1},\dots ,t_{Q})}e^{-\tau\theta_{x}}
\big((\theta^{\flat}_{x})_{\hat u_{1},\hat u_{2},\hat u_{3}}^{\hat v_{1},\hat v_{2},\hat v_{3}}-
(\theta^{\flat}_{x'})_{\hat u_{1},\hat u_{2},\hat u_{3}}^{\hat v_{1},\hat v_{2},\hat v_{3}}\big)\Big)   (e_{\tilde S_{0}})\end{gather*} et le membre de gauche de (\ref{condition-rho-uuuvvvxx'01-e)}) ne dépendent  que de 
la connaissance des  points de   
\begin{gather}\nonumber B(x,k+2M)\cup B( \tilde S_{0}, M)\cup B(\tilde S_{1}, M)\\ \label{ens-conn-e)-23dec}\cup  \bigcup _{ j\in \{1,\dots ,Q\}}
B(\tilde{\mathcal Z}_{0}^{j}, M)
\cup  \bigcup_{j\in \{1,...,7\}}
B(\tilde{\mathcal Z}_{1}^{j},M)  \end{gather}
 et des distances entre ces points. 
 \end{souslem}
 \noindent{\bf Démonstration.} 
 Le coefficient de $e_{\tilde S_{1}}$ dans 
$ e^{\tau\theta_{x}}u_{x,r,t}K_{x,Q,(t_{1},\dots ,t_{Q})}e^{-\tau\theta_{x}} (e_{\tilde S_{0}})$ ne dépend que de la connaissance des points de (\ref{ens-conn-e)-23dec}) et de leurs distances mutuelles. Cela résulte 
du sous-lemme~\ref{slem1-4eme-2j1305} ou du sous-lemme~\ref{slem1-b-2j1827}
(que l'on applique avec $\tilde{\mathcal Z}_{0}^{7}$ au lieu de $\tilde{\mathcal Z}_{0}^{1}$ et en oubliant $\tilde{\mathcal Z}_{0}^{j}$ pour $j\in \{1,...,6\}$). 
%
%des arguments donnés dans la preuve de (\ref{ineg-(1-P)-tilde-l0-bis}), qui va essentiellement de la fin de l'inégalité (\ref{ineg-(1-P)-tilde-l0-bis}) à (\ref{ineg-apres55-page105-22dec}).
% Cela résulte aussi bien sûr de la démonstration de l'inclusion de la réunion des ensembles (\ref{ens3.33-22dec-1}), (\ref{ens3.33-22dec-2}), 
%(\ref{ens3.33-22dec-3})  et  (\ref{ens3.33-22dec-4}) 
%dans l'ensemble (\ref{ens-tot-b-22dec}) dans la preuve de b) de ce lemme. 

 Pour montrer    le sous-lemme il  suffit donc de  montrer que pour $\tilde x\in B(x,1)$, $(\rho^{\flat}_{\tilde x})_{u_{1},u_{2},u_{3}}^{v_{1},v_{2},v_{3}}(\tilde S_{1})$ et 
 $(\rho^{\flat}_{\tilde  x})_{\hat u_{1},\hat u_{2},\hat u_{3}}^{\hat v_{1},\hat v_{2},\hat v_{3}}(\tilde S_{0})$  ne dépendent que de la connaissance des points de (\ref{ens-conn-e)-23dec})  et des distances entre ces points (en effet on applique ceci à $\tilde x=x$ et $\tilde x=x'$). En vertu du lemme~\ref{dependance-rho'-123123} 
 il suffit de montrer que pour $\tilde x\in B(x,1)$, $\sigma\in \{0,1\}$, $a\in \tilde S_{\sigma}$ et $j\in \{1,...,6\}$,  l'ensemble 
  \begin{gather}\label{ens-30dec1355}  \{y\in 3\de\tg(\tilde x,a), |d(\tilde x,y)-w_{j}(\sigma,\tilde x)|\leq N+6\de+4\}\end{gather}
  est inclus dans (\ref{ens-conn-e)-23dec}), 
    où l'on note 
                  \begin{gather}    \nonumber  \lambda_{1}=\frac{u_{1}}{6}  , \lambda_{2}= \frac{1+u_{2}}{6} ,\lambda_{3}=\frac{2+u_{3}}{6}  ,\lambda_{4}=\frac{6-v_{1}}{6}  ,\lambda_{5}=\frac{5-v_{2}}{6}  ,\lambda_{6}=  \frac{4-v_{3}}{6},
                  \\   \nonumber  
    \hat   \lambda_{1}=\frac{\hat u_{1}}{6}  ,\hat  \lambda_{2}= \frac{1+\hat u_{2}}{6} ,\hat \lambda_{3}=\frac{2+\hat u_{3}}{6}  ,\hat \lambda_{4}=\frac{6-\hat v_{1}}{6}  ,\hat \lambda_{5}=\frac{5-\hat v_{2}}{6}  ,\hat \lambda_{6}=  \frac{4-\hat v_{3}}{6},\\ 
   \nonumber  
  w_{j}(0,\tilde x)=E(\hat \lambda_{j}d(\tilde x,\tilde S_{0}))
    \text{\ \ et\ \ \ }w_{j}(1,\tilde x)=E(\lambda_{j}d(\tilde x,\tilde S_{1})). 
\end{gather}
 % \begin{gather}\nonumber  
 % w_{1}(1,\tilde x)=E(\frac{u_{1}}{6}d(\tilde x,\tilde S_{1})),
 % w_{2}(1,\tilde x)=E((\frac{1}{6}+\frac{u_{2}}{6})d(\tilde x,
 % \tilde S_{1})), \\ \nonumber
 % w_{3}(1,\tilde x)=E((\frac{2}{6}+\frac{u_{3}}{6})d(\tilde x,\tilde S_{1})),   w_{4}(1,\tilde x)=E((1-\frac{v_{1}}{6})d(\tilde x,\tilde S_{1})), \\ \label{def-w1-30dec1357} 
 %    w_{5}(1,\tilde x)=E((\frac{5}{6}-\frac{v_{2}}{6})d(\tilde x,\tilde S_{1})), 
 %   w_{6}(1,\tilde x)=E((\frac{4}{6}-\frac{v_{3}}{6})d(\tilde x,\tilde S_{1})) \end{gather}
 %   et 
 %    \begin{gather}\nonumber  
 % w_{1}(0,\tilde x)=E(\frac{\hat u_{1}}{6}d(\tilde x,\tilde S_{0})),
 % w_{2}(0,\tilde x)=E((\frac{1}{6}+\frac{\hat u_{2}}{6})d(\tilde x,
 % \tilde S_{0})), \\ \nonumber
 % w_{3}(0,\tilde x)=E((\frac{2}{6}+\frac{\hat u_{3}}{6})d(\tilde x,\tilde S_{0})),   w_{4}(0,\tilde x)=E((1-\frac{\hat v_{1}}{6})d(\tilde x,\tilde S_{0})), \\ \label{def-w0-30dec1357} 
 %    w_{5}(0,\tilde x)=E((\frac{5}{6}-\frac{\hat v_{2}}{6})d(\tilde x,\tilde S_{0})), 
 %   w_{6}(0,\tilde x)=E((\frac{4}{6}-\frac{\hat v_{3}}{6})d(\tilde x,\tilde S_{0})). \end{gather}

 On rappelle que $w_{j}=E(\lambda_{j}d(x,\tilde S_{1}))$. 
    
    On commence par le cas où $\sigma=1$. 
       Soit $\tilde x\in B(x,1)$,  $a\in \tilde S_{1}$, $j\in \{1,...,6\}$ et $y$ dans l'ensemble (\ref{ens-30dec1355}). Soit $z\in \geod(x,a)$ vérifiant $d(x,z)=w_{j}$, si bien  que $z$ appartient à $\tilde {\mathcal Z}_{1}^{j}$. 
Comme $|d(x,\tilde S_{1})-d(\tilde x,\tilde S_{1})|\leq 1$ on a 
$|w_{j}(1,\tilde x)-w_{j}|\leq 1$ et donc $|d(x,y)-d(x,z)|\leq N+6\de+6$. Comme $y$ et $z$ appartiennent à $(3\de+2)\tg(x,a)$, $(H_{\de}(y,x,z,a))$ implique 
$d(y,z)\leq (N+6\de+6)+(3\de+2)+\de$. On suppose $(N+6\de+6)+(3\de+2)+\de\leq M$, ce qui est permis par $(H_{M})$. On a donc $d(y,z)\leq M$ et $y$ appartient à (\ref{ens-conn-e)-23dec}). 
 
 On considère maintenant le cas où $\sigma=0$. Soit $\tilde x\in B(x,1)$,  $a\in \tilde S_{0}$, $j\in \{1,...,6\}$ et $y$ dans l'ensemble (\ref{ens-30dec1355}).  On commence par montrer
 \begin{gather}\label{ineg-30dec1558}|w_{j}(0,\tilde x)-w_{j}|\leq QF+2.\end{gather}
   %  On pose $(\lambda, \hat \lambda)$ égal à 
   % \begin{gather*}(\frac{u_{1}}{6}, \frac{\hat u_{1}}{6}), (\frac{1+u_{2}}{6}, \frac{1+\hat u_{2}}{6}), (\frac{2+u_{3}}{6}, \frac{2+\hat u_{3}}{6}), \\  (\frac{6-v_{1}}{6}, \frac{6-\hat v_{1}}{6}), (\frac{5-v_{2}}{6}, \frac{5-\hat v_{2}}{6}), (\frac{4-v_{3}}{6}, \frac{4-\hat v_{3}}{6})\end{gather*}
   %       suivant que $j=1,2,3,4,5,6$. 
  On a    $\hat \lambda_{j}=\frac{\lambda_{j}}{1+\kappa}$ et $\frac{r_{1}(\tilde Z)+1+r}{1+\kappa}=r_{1}(\tilde Z)+1$, d'où 
  $$\hat\lambda_{j}(r_{1}(\tilde Z)+1+r)=\lambda_{j}(r_{1}(\tilde Z)+1).$$ 
  D'après (\ref{cond-r0r1-23dec}) on a $|(r_{0}(\tilde Z)+1)-(r_{1}(\tilde Z)+1+r)|\leq QF$. 
   Donc $$|\hat \lambda_{j} (r_{0}(\tilde Z)+1)-\lambda_{j}(r_{1}(\tilde Z)+1)|\leq QF,$$
    et  on en déduit immédiatement  
    $$|\hat \lambda_{j} r_{0}(\tilde Z)-\lambda_{j}r_{1}(\tilde Z)|\leq QF+1. $$
     Comme  $r_{0}(\tilde Z)=d(x,\tilde S_{0})$ et $r_{1}(\tilde Z)=d(x,\tilde S_{1})$, (\ref{ineg-30dec1558}) en résulte aisément. 
     
        Soit $b\in \tilde S_{1}$ et $z\in \geod(x,b)$ vérifiant $d(x,z)=w_{j}$
        si bien  que $z$ appartient à $\tilde {\mathcal Z}_{1}^{j}$. 
         Il résulte facilement de (\ref{ineg-30dec1558}) que 
        \begin{gather}\label{ineg-30dec1630}
        |d(x,y)-d(x,z)|\leq QF+N+6\de+7.
        \end{gather}
        Comme $b\in 2F\tg(x,a)$ on a $z\in 2F\tg(x,a)$. On a  $y\in (3\de+2)\tg(x,a)$ et comme $2F\geq 3\de+2$, $y$ et $z$ appartiennent à 
        $2F\tg(x,a)$. Grâce à (\ref{ineg-30dec1630}), $(H_{\de}(y,x,z,a))$ implique 
$d(y,z)\leq (QF+N+6\de+7)+2F+\de$. On suppose $(QF+N+6\de+7)+2F+\de\leq M$, ce qui est permis par $(H_{M})$. On a donc $d(y,z)\leq M$ et $y$ appartient à (\ref{ens-conn-e)-23dec}). Ceci termine la démonstration du sous-lemme~\ref{sslem-e-1j2139}. \cqfd

  Une conséquence immédiate du sous-lemme~\ref{sslem-e-1j2139} est que       
       la condition (\ref{condition-rho-uuuvvvxx'01-e)}) ne dépend  que de $\tilde Z$.

       \begin{souslem}\label{slem2-e-2j1943} 
       Le cardinal de $(\Lambda_{Z,t,(t_{1},...,t_{Q})})_{u_{1},u_{2},u_{3}}^{v_{1},v_{2},v_{3}} $ est majoré par une constante de la forme $C(\de,K,N,Q,P,M)$. 
              \end{souslem}
       \noindent{\bf Démonstration.} 
 En effet, grâce au lemme~\ref{nombre-dist-connaitre-par-point-natural}, pour connaître les distances entre les points de   
\begin{gather} \label{23dec2135} 
 B(\tilde S_{0},M)\cup  \bigcup _{ j\in \{1,\dots ,Q\}} B(\tilde{\mathcal Z}_{0}^{j}, M) \cup  \bigcup_{j\in \{1,...,7\}}B(\tilde{\mathcal Z}_{1}^{j},M)
\end{gather} 
 et ceux de 
\begin{gather}\label{23dec2137} 
\bigcup _{ i\in \{1,\dots ,m+1\}}
B(\tilde S_{i}, M)
\cup  \bigcup _{i\in \{0,\dots,m\}, j\in \{1,\dots , l_{i}\}}B(\tilde {\mathcal Y}_{i+1}^{j}, M)\cup B(x,k+2M)\end{gather}
il suffit de connaître les distances entre les points de (\ref{23dec2135})  et $C$ points de (\ref{23dec2137}), avec $C=C(\de,K,N,Q,P,M)$ et grâce à (\ref{cond-r0r1-23dec}) ces distances sont déterminées à $C'=C(\de,K,N,Q,P,M)$ près par les distances de $\tilde S_{1}$ à ces $C$ points (qui font partie de la donnée de $Z$) et par les entiers $(t_{0}^{j}(\tilde Z))_{j\in \{1,\pp,Q\}}$ et $(t_{1}^{j}(\tilde Z))_{j\in \{1,\pp,7\}}$, 
qui sont eux-mêmes déterminés à $C''=C(\de,K,N,Q,P,M)$ près par
$r_{0}(Z),r,t,t_{1},...,t_{Q},u_{1},u_{2},u_{3},v_{1},v_{2},v_{3}$. 
 \cqfd
 
 \noindent {\bf Suite de la démonstration du e).}
   On termine maintenant   la preuve de (\ref{ineg-(1-P)-tilde-l0-bis-xx'-e)}). 
Pour  \begin{gather*}(\tilde a_{1},\dots,\tilde a_{p-1},\tilde S_{0},...,\tilde S_{m+1},(\tilde {\mathcal Y}_{i}^{j})_{i\in \{0,\dots,m+1\}, j\in \{1,\dots ,\tilde l_{i}\}},(\tilde {\mathcal Z}_{0}^{j})_{ j\in \{1,\dots ,Q\}},(\tilde {\mathcal Z}_{1}^{j})_{j\in \{1,\dots ,7\}}) \\ \in (\pi_{x}^{\natural, p-1,k,m+1,(\tilde l_{0},...,\tilde l_{m+1}),Q,7})^{-1}(\tilde Z) \end{gather*} 
on considère 
\begin{gather}
\nonumber 
\sum _{(b_{1},...,b_{p})}
\Big(\Big(\big((\theta^{\flat}_{x})_{u_{1},u_{2},u_{3}}^{v_{1},v_{2},v_{3}}-
(\theta^{\flat}_{x'})_{u_{1},u_{2},u_{3}}^{v_{1},v_{2},v_{3}}\big) e^{\tau\theta_{x}}u_{x,r,t}K_{x,Q,(t_{1},\dots ,t_{Q})}e^{-\tau\theta_{x}}-
\\ \nonumber  
e^{\tau\theta_{x}}u_{x,r,t}K_{x,Q,(t_{1},\dots ,t_{Q})}e^{-\tau\theta_{x}}
\big((\theta^{\flat}_{x})_{\hat u_{1},\hat u_{2},\hat u_{3}}^{\hat v_{1},\hat v_{2},\hat v_{3}}-
(\theta^{\flat}_{x'})_{\hat u_{1},\hat u_{2},\hat u_{3}}^{\hat v_{1},\hat v_{2},\hat v_{3}}\big)\Big)  \\ \label{somme-1j2148}
  (e_{\tilde a_{1}}\wedge ...\wedge e_{\tilde a_{p-1}})\Big)(b_{1},...,b_{p}), 
\end{gather} 
où la somme porte sur les énumérations $(b_{1},...,b_{p})$ de $\tilde S_{1}$ telles que 
$$(b_{1},\dots,b_{p},\tilde S_{1},...,\tilde S_{m+1},(\tilde {\mathcal Y}_{i+1}^{j})_{i\in \{0,\dots,m\}, j\in \{1,\dots ,\tilde l_{i}\}}
)\in (\pi_{x}^{p,k,m,(l_{0},..., l_{m})})^{-1}(Z).$$ 
Comme la somme (\ref{somme-1j2148}) a au plus $p!$ termes, 
 le 3) de la proposition~\ref{recap-supp-connaiss-H-uK} montre qu'elle 
est majorée par une constante de la forme $C(\de,K,N,Q,P,M,T)$. 
D'après le sous-lemme~\ref{sslem-e-1j2139} la somme (\ref{somme-1j2148}) ne dépend que de $\tilde Z$  et on peut donc la noter $(\alpha_{Z,\tilde Z,t,\tau, (t_{1},...,t_{Q})})_{u_{1},u_{2},u_{3}}^{v_{1},v_{2},v_{3}}$. 
D'après le sous-lemme~\ref{slem0-e-2j1938} et le sous-lemme~\ref{ineg73-3j2030} (ou~\ref{just-r0-r1-r-xx'-3j2157}), on a \begin{gather*}\xi_{Z}\Big(\Big(\big((\theta^{\flat}_{x})_{u_{1},u_{2},u_{3}}^{v_{1},v_{2},v_{3}}-
(\theta^{\flat}_{x'})_{u_{1},u_{2},u_{3}}^{v_{1},v_{2},v_{3}}\big) e^{\tau\theta_{x}}u_{x,r,t}K_{x,Q,(t_{1},\dots ,t_{Q})}e^{-\tau\theta_{x}}-
\\ 
e^{\tau\theta_{x}}u_{x,r,t}K_{x,Q,(t_{1},\dots ,t_{Q})}e^{-\tau\theta_{x}}
\big((\theta^{\flat}_{x})_{\hat u_{1},\hat u_{2},\hat u_{3}}^{\hat v_{1},\hat v_{2},\hat v_{3}}-
(\theta^{\flat}_{x'})_{\hat u_{1},\hat u_{2},\hat u_{3}}^{\hat v_{1},\hat v_{2},\hat v_{3}}\big)\Big)  f\Big)
\\ 
   =\frac{1}{(p-1)!}\sum_{\tilde Z\in 
 (\Lambda_{Z,t,(t_{1},...,t_{Q})}^{\neq})_{u_{1},u_{2},u_{3}}^{v_{1},v_{2},v_{3}}}
 (\alpha_{Z,\tilde Z,t,\tau,(t_{1},...,t_{Q})})_{u_{1},u_{2},u_{3}}^{v_{1},v_{2},v_{3}}\ \xi_{\tilde Z}(f). \end{gather*}
 Par Cauchy-Schwarz et grâce au sous-lemme~\ref{slem2-e-2j1943}  
 on en déduit   (\ref{ineg-(1-P)-tilde-l0-bis-xx'-e)}). 
  %
  %
%
%
%
%
% Les arguments pour montrer  (\ref{ineg-(1-P)-tilde-l0-bis-xx'-e)})
% et justifier que la condition (\ref{condition-rho-uuuvvvxx'01-e)}) ne dépend  que de $\tilde Z$ 
% sont analogues à ceux donnés dans les preuves de b)  de d). En plus on applique le lemme~\ref{dependance-rho'-123123} à $S=\tilde S_{0}$ et $S=\tilde S_{1}$, puis on montre  qu'il existe une constante $C=C(\de,K,N,Q)$ telle que $\tilde S_{0},\tilde S_{1}$ comme ci-dessus on ait 
%$d(\tilde S_{0},\tilde S_{1})\leq C+r$ et en utilisant $(H_{M})$ (donc en particulier le fait que $M$ est grand par rapport à $C$).  
% on en déduit que 
  %$$\Big(\big((\theta^{\flat}_{x})_{u_{1},u_{2},u_{3}}^{v_{1},v_{2},v_{3}}-
%(\theta^{\flat}_{x'})_{u_{1},u_{2},u_{3}}^{v_{1},v_{2},v_{3}}\big) e^{\tau\theta_{x}}u_{x,r,t}K_{x,Q,(t_{1},\dots ,t_{Q})}e^{-\tau\theta_{x}}-$$ $$
%e^{\tau\theta_{x}}u_{x,r,t}K_{x,Q,(t_{1},\dots ,t_{Q})}e^{-\tau\theta_{x}}
%\big((\theta^{\flat}_{x})_{\hat u_{1},\hat u_{2},\hat u_{3}}^{\hat v_{1},\hat v_{2},\hat v_{3}}-
%(\theta^{\flat}_{x'})_{\hat u_{1},\hat u_{2},\hat u_{3}}^{\hat v_{1},\hat v_{2},\hat v_{3}}\big)\Big)   (e_{\tilde S_{0}})$$ ne dépend que de

 Montrons maintenant e)  à l'aide de 
(\ref{ineg-(1-P)-tilde-l0-bis-xx'-e)}).  
Le sous-lemme~\ref{mesure-uuuvvv-dbemol-e)} appliqué à $\rho=r+C$ avec $C=C(\de,K,N,Q)$  implique facilement qu'il existe $C_{1}=C(\de,K,N,Q,r)$ tel que pour $\tilde S_{0},\tilde S_{1}$ comme ci-dessus  la mesure de l'ensemble des $(u_{1},u_{2}, u_{3}$, $v_{1},v_{2},v_{3})\in [5\kappa,1-5\kappa[^{6}$ vérifiant (\ref{condition-rho-uuuvvvxx'01-e)}) est $\leq \frac{C_{1}}{1+r_{0}(Z)}$. 
Soit $C_{2}$ comme dans 
 (\ref{ineg-(1-P)-tilde-l0-bis-xx'-e)}).  Soient $k,m$, $l_{0},\dots ,l_{m}\in \N $ et posons $ \tilde l_{0}=0$ et $ \tilde l_{i}=l_{i-1}$ pour $i\in \{1,\pp,m+1\}$ comme précédemment. Soit $C_{3}= C(\de,K,N,Q,P,M,r)$ comme dans le sous-lemme~\ref{slem3-b-2j1836}.  
 Soit  $$Z\in \overline Y_{x}^{p,k,m,(l_{0},...,l_{m})}\text{\ \  vérifiant \ \ }r_{0}(Z)> k+P.$$
  On  a alors 
   \begin{gather*}\int_{(u_{1},u_{2},u_{3},v_{1},v_{2},v_{3})\in [5\kappa,1-5\kappa[^{6}}\Big(\sum_{\tilde Z\in (\Lambda_{Z,t,(t_{1},...,t_{Q})}^{\neq})_{u_{1},u_{2},u_{3}}^{v_{1},v_{2},v_{3}}}\sharp\big((\pi_{x}^{\natural,p-1,k,m+1,(\tilde l_{0},...,\tilde l_{m+1}),Q,7})^{-1}(\tilde Z)\big)\Big)\\ du_{1}...dv_{3} \leq \frac{C_{1}C_{3}}{r_{0}(Z)+1}\sharp\big((\pi_{x}^{p,k,m,(l_{0},..., l_{m})})^{-1}(Z)\big),\end{gather*} 
 Notons $I_{Z}$ l'ensemble des $$(t,t_{1},...,t_{Q},u_{1},u_{2}, u_{3},v_{1},v_{2},v_{3})\in [0,1]^{Q+1}\times [5\kappa,1-5\kappa[^{6}$$ tels qu'il existe $\tilde Z\in (\Lambda_{Z,t,(t_{1},...,t_{Q})}^{\neq})_{u_{1},u_{2},u_{3}}^{v_{1},v_{2},v_{3}}$   vérifiant   $$\sharp\big((\pi_{x}^{\natural,p-1,k,m+1,(\tilde l_{0},...,\tilde l_{m+1}),Q,7})^{-1}(\tilde Z)\big)\geq (r_{0}(Z)+1)^{-\frac{1}{2}}\sharp\big((\pi_{x}^{p,k,m,(l_{0},..., l_{m})})^{-1}(Z)\big). $$ La mesure de $I_{Z}$ est donc $\leq C_{1}C_{3}(r_{0}(Z)+1)^{-\frac{1}{2}}$. Grâce à Cauchy-Schwarz on déduit de (\ref{ineg-(1-P)-tilde-l0-bis-xx'-e)})
que 
\begin{gather}\nonumber 
\sup_{\tau\in [0,T]}\Big|\xi_{Z}\Big(\Big(\int_{(t,t_{1},...,t_{Q},u_{1},u_{2}, u_{3},v_{1},v_{2},v_{3})\in  I_{Z}} 
\\ 
\nonumber \Big(\big((\theta^{\flat}_{x})_{u_{1},u_{2},u_{3}}^{v_{1},v_{2},v_{3}}-
(\theta^{\flat}_{x'})_{u_{1},u_{2},u_{3}}^{v_{1},v_{2},v_{3}}\big) e^{\tau\theta_{x}}u_{x,r,t}K_{x,Q,(t_{1},\dots ,t_{Q})}e^{-\tau\theta_{x}}
\\
\nonumber 
-e^{\tau\theta_{x}}u_{x,r,t}K_{x,Q,(t_{1},\dots ,t_{Q})}e^{-\tau\theta_{x}}
\big((\theta^{\flat}_{x})_{\hat u_{1},\hat u_{2},\hat u_{3}}^{\hat v_{1},\hat v_{2},\hat v_{3}}-
(\theta^{\flat}_{x'})_{\hat u_{1},\hat u_{2},\hat u_{3}}^{\hat v_{1},\hat v_{2},\hat v_{3}}\big)\Big) dt...dv_{3} \Big)  f\Big)\Big|^{2} 
\\ 
\label{ineg-IZ-theta-xx'-e)} 
\leq 2^{Q+7}6^{6}C_{1}C_{2}C_{3}(r_{0}(Z)+1)^{-\frac{1}{2}}
\sum_{\tilde Z\in \Lambda_{Z}} (r_{0}(\tilde Z)+1)^{-Q}(r_{1}(\tilde Z)+1)^{-7} 
\big|\xi_{\tilde Z}(f)\big|^{2}. 
\end{gather}
D'après le sous-lemme~\ref{slem3-b-2j1836}  et le lemme~\ref{lemme-cardinaux}, 
 il existe $C_{4}=C(\de,K,N,Q,P,M,r)$ tel  que pour $\tilde Z\in \Lambda_{Z}$  on ait 
 \begin{gather}\label{ineg-C4-e-22dec}\sharp\big((\pi_{x}^{\natural,p-1,k,m+1,(\tilde l_{0},...,\tilde l_{m+1}),Q,7})^{-1}(\tilde Z)\big)\leq C_{4}\sharp\big((\pi_{x}^{p,k,m,(l_{0},..., l_{m})})^{-1}(Z)\big).\end{gather} 
  D'autre part pour $(t,t_{1},...,t_{Q},u_{1},u_{2}, u_{3},v_{1},v_{2},v_{3})\in \big(  [0,1]^{Q+1}\times [5\kappa,1-5\kappa[^{6}\setminus I_{Z}\big)$ et $\tilde Z\in (\Lambda_{Z,t,(t_{1},...,t_{Q})}^{\neq})_{u_{1},u_{2},u_{3}}^{v_{1},v_{2},v_{3}}$  on a $$\sharp\big((\pi_{x}^{\natural,p-1,k,m+1,(\tilde l_{0},...,\tilde l_{m+1}),Q,7})^{-1}(\tilde Z)\big)\leq (r_{0}(Z)+1)^{-\frac{1}{2}}\sharp\big((\pi_{x}^{p,k,m,(l_{0},..., l_{m})})^{-1}(Z)\big). $$
   Par Cauchy-Schwarz on déduit alors de (\ref{ineg-(1-P)-tilde-l0-bis-xx'-e)}) que 
   \begin{gather}
   \nonumber \big(\sharp(\pi_{x}^{p,k,m,(l_{0},..., l_{m})})^{-1}(Z)\big)^{-\alpha}
   \\
   \nonumber 
   \sup_{\tau\in [0,T]}\Big|\xi_{Z}\Big(\Big(\int_{(t,t_{1},...,t_{Q},u_{1},u_{2}, u_{3},v_{1},v_{2},v_{3}) \in  \big(  [0,1]^{Q+1}\times [5\kappa,1-5\kappa[^{6}\setminus I_{Z}\big)}
    \\
   \nonumber 
\Big(\big((\theta^{\flat}_{x})_{u_{1},u_{2},u_{3}}^{v_{1},v_{2},v_{3}}-
(\theta^{\flat}_{x'})_{u_{1},u_{2},u_{3}}^{v_{1},v_{2},v_{3}}\big) e^{\tau\theta_{x}}u_{x,r,t}K_{x,Q,(t_{1},\dots ,t_{Q})}e^{-\tau\theta_{x}}
 \\
   \nonumber 
-e^{\tau\theta_{x}}u_{x,r,t}K_{x,Q,(t_{1},\dots ,t_{Q})}e^{-\tau\theta_{x}}
\big((\theta^{\flat}_{x})_{\hat u_{1},\hat u_{2},\hat u_{3}}^{\hat v_{1},\hat v_{2},\hat v_{3}}-
(\theta^{\flat}_{x'})_{\hat u_{1},\hat u_{2},\hat u_{3}}^{\hat v_{1},\hat v_{2},\hat v_{3}}\big)\Big) dt...dv_{3} \Big)  f\Big)\Big|^{2}
 \\
   \nonumber 
\leq 2^{Q+7}6^{6}C_{2}(r_{0}(Z)+1)^{-\frac{\alpha}{2}}
\sum_{\tilde Z\in \Lambda_{Z}} (r_{0}(\tilde Z)+1)^{-Q}(r_{1}(\tilde Z)+1)^{-7}
 \\
  \label{ineg-hors-IZ-theta-xx'-e)}
  \big(\sharp(\pi_{x}^{\natural,p-1,k,m+1,(\tilde l_{0},...,\tilde l_{m+1}),Q,7})^{-1}(\tilde Z)\big)^{-\alpha} 
  \big|\xi_{\tilde Z}(f)\big|^{2}.\end{gather}
 En combinant les inégalités (\ref{ineg-IZ-theta-xx'-e)}), (\ref{ineg-C4-e-22dec})  et 
  (\ref{ineg-hors-IZ-theta-xx'-e)}) 
  et par Cauchy-Schwarz on obtient que 
\begin{gather}
\nonumber 
 \big(\sharp(\pi_{x}^{p,k,m,(l_{0},..., l_{m})})^{-1}(Z)\big)^{-\alpha}\sup_{\tau\in [0,T]}   \\
   \nonumber 
 \Big|\xi_{Z}\Big(\Big(\int_{(u_{1},u_{2}, u_{3},v_{1},v_{2},v_{3}) \in   [5\kappa,1-5\kappa[^{6}}
  \Big(\big((\theta^{\flat}_{x})_{u_{1},u_{2},u_{3}}^{v_{1},v_{2},v_{3}}-
(\theta^{\flat}_{x'})_{u_{1},u_{2},u_{3}}^{v_{1},v_{2},v_{3}}\big) e^{\tau\theta_{x}}u_{x,r}K_{x}e^{-\tau\theta_{x}}
 \\
   \nonumber 
-e^{\tau\theta_{x}}u_{x,r}K_{x}e^{-\tau\theta_{x}}
\big((\theta^{\flat}_{x})_{\hat u_{1},\hat u_{2},\hat u_{3}}^{\hat v_{1},\hat v_{2},\hat v_{3}}-
(\theta^{\flat}_{x'})_{\hat u_{1},\hat u_{2},\hat u_{3}}^{\hat v_{1},\hat v_{2},\hat v_{3}}\big)\Big) du_{1}...dv_{3} \Big)  f\Big)\Big|^{2}
\\ \nonumber 
%$$\Big|\sum _{(a_{1},\dots,a_{p},S_{1},...,S_{m},(\mathcal Y_{i}^{j})_{i\in \{0,\dots,m\}, j\in \{1,\dots ,l_{i}\}}) \in 
%(\pi_{x}^{p,k,m,(l_{0},...,l_{m})})^{-1}(Z)} ((h_{x}-h_{x'})  f)(a_1,...,a_p)\Big|^{2}$$ 
\leq 2^{Q+8}6^{6}\Big(C_{1}C_{2}C_{3}C_{4}^{\alpha}(r_{0}(Z)+1)^{-\frac{1}{2}}+
 C_{2}(r_{0}(Z)+1)^{-\frac{\alpha}{2}}\Big)\sum_{\tilde Z\in \Lambda_{Z}} (r_{0}(\tilde Z)+1)^{-Q}
 \\
  \label{ineg-theta-xx'-e-4j1117)}
(r_{1}(\tilde Z)+1)^{-7}
   \big(\sharp(\pi_{x}^{\natural,p-1,k,m+1,(\tilde l_{0},...,\tilde l_{m+1}),Q,7})^{-1}(\tilde Z)\big)^{-\alpha} 
  \big|\xi_{\tilde Z}(f)\big|^{2}.\end{gather} 
    De plus pour $\tilde Z\in \Lambda_{Z}$ on a $\prod_{i=0}^{m}s_{i}(Z)^{-l_{i}}=\prod_{i=0}^{m+1}s_{i}(\tilde Z)^{-\tilde l_{i}}$. 
  En sommant sur les  $Z$ tels que $r_{0}(Z)=n$, et en posant $\kappa=\frac{r}{1+n}$,  on déduit donc de (\ref{ineg-theta-xx'-e-4j1117)}) que pour tout $n> P$, 
  \begin{gather}
    \nonumber 
 \sup_{\tau\in [0,T]} \Big\|(\P_{n}-\P_{n-1})\Big(
\Big(\int_{  [5\kappa,1-5\kappa[^{6}}
\big((\theta^{\flat}_{x})_{u_{1},u_{2},u_{3}}^{v_{1},v_{2},v_{3}}-
(\theta^{\flat}_{x'})_{u_{1},u_{2},u_{3}}^{v_{1},v_{2},v_{3}}\big) du_{1}...dv_{3}\Big) 
\\
   \nonumber 
e^{\tau\theta_{x}} u_{x,r}K_{x}e^{-\tau\theta_{x}}
 -e^{\tau\theta_{x}}u_{x,r}K_{x}e^{-\tau\theta_{x}}
 \\
   \nonumber 
\Big(
\int_{  [5\kappa,1-5\kappa[^{6}}
\big((\theta^{\flat}_{x})_{\hat u_{1},\hat u_{2},\hat u_{3}}^{\hat v_{1},\hat v_{2},\hat v_{3}}-
(\theta^{\flat}_{x'})_{\hat u_{1},\hat u_{2},\hat u_{3}}^{\hat v_{1},\hat v_{2},\hat v_{3}}\big) du_{1}...dv_{3} \Big) 
\Big) 
 \Big\|_{\L(\H_{x,s}^{\natural,Q,7}(\Delta_{p-1}),\H^{\rightarrow}_{x,s}(\Delta_{p}))}^{2}
 \\ \nonumber
  \leq  2^{Q+8}6^{6}p! B e^{2(QF-r)s}
\Big(C_{1}C_{2}C_{3}C_{4}^{\alpha}(n+1)^{-\frac{1}{2}}+
 C_{2}(n+1)^{-\frac{\alpha}{2}}\Big) \end{gather}
  où 
 $(\hat u_{1},\hat u_{2},\hat u_{3}, \hat v_{1},\hat v_{2},\hat v_{3})$ est défini par (\ref{hatuuuvvv-uuuvvv}). 
 Le facteur $p!$ est dû au fait que $\tilde Z$ détermine $Z$ à permutation près de $a_{1},...,a_{p}$. On rappelle que $\P_{n}-\P_{n-1}$ est le projecteur orthogonal défini par $$(\P_{n}-\P_{n-1})(e_{S})=e_{S} \text{\  si \ } d(x,S)=n \text{ \  et \ }  
  (\P_{n}-\P_{n-1})(e_{S})=0 \text{ \  sinon.}$$ 
  
D'après les lemmes~\ref{eq-normes-natural}
  et~\ref{equiv-gauche}, les normes de $\H_{x,s}^{\natural,Q,7}(\Delta_{p-1})$ et $\H^{\rightarrow}_{x,s}(\Delta_{p})$ sont équivalentes à celles de 
  $\H_{x,s}(\Delta_{p-1})$ et $\H_{x,s}(\Delta_{p})$. Donc  il existe $C=C(\de,K,N,Q,P,M,s,B,r,T)$ tel  que
  pour tout $n> P$, 
  \begin{gather}
    \nonumber 
 \sup_{\tau\in [0,T]} \Big\|(\P_{n}-\P_{n-1})\Big(\Big(
\int_{  [5\kappa,1-5\kappa[^{6}}
\big((\theta^{\flat}_{x})_{u_{1},u_{2},u_{3}}^{v_{1},v_{2},v_{3}}-
(\theta^{\flat}_{x'})_{u_{1},u_{2},u_{3}}^{v_{1},v_{2},v_{3}}\big) du_{1}...dv_{3}\Big) 
\\
   \nonumber 
e^{\tau\theta_{x}} u_{x,r}K_{x}e^{-\tau\theta_{x}}
 -e^{\tau\theta_{x}}u_{x,r}K_{x}e^{-\tau\theta_{x}}
 \\
   \nonumber 
\Big(
\int_{  [5\kappa,1-5\kappa[^{6}}
\big((\theta^{\flat}_{x})_{\hat u_{1},\hat u_{2},\hat u_{3}}^{\hat v_{1},\hat v_{2},\hat v_{3}}-
(\theta^{\flat}_{x'})_{\hat u_{1},\hat u_{2},\hat u_{3}}^{\hat v_{1},\hat v_{2},\hat v_{3}}\big) du_{1}...dv_{3} \Big) 
\Big) 
\Big\|_{\L(\H_{x,s}(\Delta_{p-1}),\H_{x,s}(\Delta_{p}))}^{2}
 \\ \label{ineg-penult-e-4j1149}
\leq  C(n+1)^{-\frac{\alpha}{2}}\end{gather}
 où   $\kappa=\frac{r}{1+n}$ et où 
 $(\hat u_{1},\hat u_{2},\hat u_{3}, \hat v_{1},\hat v_{2},\hat v_{3})$ est défini par (\ref{hatuuuvvv-uuuvvv}). 

%  On conclut alors comme dans les preuves de a), b) et d), avec un  petit argument supplémentaire dû 
%  au domaine d'intégration   $[5\kappa,1-5\kappa[^{6}$ (au lieu de $[0,1[^{6}$) et 
%    au changement de variable $$( u_{1},u_{2}, u_{3},v_{1},v_{2},v_{3})\mapsto (\hat u_{1},\hat u_{2},\hat u_{3}, \hat v_{1},\hat v_{2},\hat v_{3}).$$  
    
\begin{souslem}\label{arg-suppl-e-4j1148}    
     Il existe $C=C(\de,K,N,Q,P,M,s,B)$ tel que   pour $p\in \{1,...,p_{\max}\}$ et $\kappa\in ]0,\frac{1}{10}[$, 
\begin{gather*}\Big\| \int_{ [5\kappa,1-5\kappa[^{6}}
\big((\theta^{\flat}_{x})_{\hat u_{1},\hat u_{2},\hat u_{3}}^{\hat v_{1},\hat v_{2},\hat v_{3}}-
(\theta^{\flat}_{x'})_{\hat u_{1},\hat u_{2},\hat u_{3}}^{\hat v_{1},\hat v_{2},\hat v_{3}}\big)du_{1}...dv_{3}- \\ \int_{ [0,1[^{6}}\big((\theta^{\flat}_{x})_{u_{1},u_{2},u_{3}}^{v_{1},v_{2},v_{3}}-
(\theta^{\flat}_{x'})_{u_{1},u_{2},u_{3}}^{v_{1},v_{2},v_{3}}\big)du_{1}...dv_{3}\Big\|_{\L(\H_{x,s}(\Delta_{p}),\H_{x,s}(\Delta_{p}))}^{2}\leq C\kappa\\ \text{et\ \ }\Big\| \int_{ [5\kappa,1-5\kappa[^{6}}
\big((\theta^{\flat}_{x})_{ u_{1}, u_{2}, u_{3}}^{ v_{1}, v_{2}, v_{3}}-
(\theta^{\flat}_{x'})_{ u_{1}, u_{2}, u_{3}}^{ v_{1}, v_{2}, v_{3}}\big)du_{1}...dv_{3}- 
\\
\int_{[0,1[^{6}}\big((\theta^{\flat}_{x})_{u_{1},u_{2},u_{3}}^{v_{1},v_{2},v_{3}}-
(\theta^{\flat}_{x'})_{u_{1},u_{2},u_{3}}^{v_{1},v_{2},v_{3}}\big)du_{1}...dv_{3}\Big\|_{\L(\H_{x,s}(\Delta_{p}),\H_{x,s}(\Delta_{p}))}^{2}\leq C\kappa.\end{gather*}
\end{souslem}
\noindent{\bf Démonstration.}
Comme $\|\theta_{x}-\theta_{x'} \|_{\L(\H_{x,s}(\Delta_{p}),\H_{x,s}(\Delta_{p}))}\leq 1$, il suffit de montrer qu'il existe $C=C(\de,K,N,Q,P,M,s,B)$ tel que  pour $\tilde x\in \{x,x'\}$ on ait 
\begin{gather}\nonumber \Big\| \int_{ [5\kappa,1-5\kappa[^{6}}
\big((\theta^{\flat}_{\tilde x})_{\hat u_{1},\hat u_{2},\hat u_{3}}^{\hat v_{1},\hat v_{2},\hat v_{3}}-
\theta_{\tilde x}\big)du_{1}...dv_{3}- 
\\ 
\label{slem-fin1-e-4j1249}
\int_{ [0,1[^{6}}\big((\theta^{\flat}_{\tilde x})_{u_{1},u_{2},u_{3}}^{v_{1},v_{2},v_{3}}-
\theta_{\tilde x}\big)du_{1}...dv_{3}\Big\|_{\L(\H_{x,s}(\Delta_{p}),\H_{x,s}(\Delta_{p}))}^{2}\leq C\kappa 
 \\ \label{slem-fin2-e-4j1249}
 \text{et\ \ }\Big\| \int_{ [0,1[^{6}\setminus [5\kappa,1-5\kappa[^{6}}
\big((\theta^{\flat}_{\tilde x})_{ u_{1}, u_{2}, u_{3}}^{ v_{1}, v_{2}, v_{3}}-
\theta_{\tilde x}\big)du_{1}...dv_{3}\Big\|_{\L(\H_{x,s}(\Delta_{p}),\H_{x,s}(\Delta_{p}))}^{2}\leq C\kappa.\end{gather} 
Grâce à l'équivalence des normes de $\H_{x,s}(\Delta_{p})$ et $\H_{x',s}(\Delta_{p})$, il suffit de montrer (\ref{slem-fin1-e-4j1249}) et (\ref{slem-fin2-e-4j1249}) pour $\tilde x=x$. 

Soient $h_{\kappa,1}$ et $h_{\kappa,2}$ les fonctions mesurables sur $[0,1[^{6}$ définies de la manière suivante: $h_{\kappa,1}$ est telle que pour toute fonction continue 
$f$ sur $ [0,1]^{6}$ on a 
\begin{gather*}\int_{ [5\kappa,1-5\kappa[^{6}}
f(\hat u_{1},\hat u_{2},\hat u_{3},\hat v_{1},\hat v_{2},\hat v_{3})
du_{1}...dv_{3}- 
\int_{ [0,1[^{6}}f(u_{1},u_{2},u_{3},v_{1},v_{2},v_{3})du_{1}...dv_{3}
\\ =\int _{[0,1[^{6}}fh_{\kappa,1}du_{1}...dv_{3}\end{gather*} et 
$h_{\kappa,2}$ est  la fonction caractéristique de $[0,1[^{6}\setminus [5\kappa,1-5\kappa[^{6}$. 

Il suffit donc de montrer qu'il existe $C=C(\de,K,N,Q,P,M,s,B)$ tel que  pour $ i \in \{1,2\}$ on ait 
\begin{gather} \nonumber 
\Big\|
\int_{ [0,1[^{6}}\big((\theta^{\flat}_{x})_{u_{1},u_{2},u_{3}}^{v_{1},v_{2},v_{3}}-
\theta_{x}\big)h_{\kappa,i}(u_{1},u_{2},u_{3},v_{1},v_{2},v_{3})du_{1}...dv_{3}\Big\|_{\L(\H_{x,s}(\Delta_{p}),\H_{x,s}(\Delta_{p}))}^{2}
\\ \label{fin-e-4j1540}
\leq C\kappa .\end{gather} 

A partir de maintenant on reprend les notations 
$\Lambda_{Z}$ et 
$(\Lambda_{Z})_{u_{1},u_{2},u_{3}}^{v_{1},v_{2},v_{3}}$
de la démonstration du lemme~\ref{continuite-eta}. 
Ces notations sont définies entre les formules (\ref{cont-eta-etape}) et (\ref{defZ0j-1j1124}). 

Pour montrer (\ref{fin-e-4j1540}) il suffit de montrer qu'il existe $C=C(\de,K,N,Q,P,M,s,B)$ tel que pour $i\in \{1,2\}$, 
$k,m,l_{0},\dots ,l_{m}\in \N $, $Z\in \overline Y_{x}^{p,k,m,(l_{0},...,l_{m})}$  et $f\in \C^{(\Delta_{p})}$ on ait   \begin{gather}\nonumber  \Big|\xi_{Z}\Big(\Big(\int_{ [0,1[^{6}}\big((\theta^{\flat}_{ x})_{u_{1},u_{2},u_{3}}^{v_{1},v_{2},v_{3}}-
\theta_{ x}\big)h_{\kappa,i}(u_{1},u_{2},u_{3},v_{1},v_{2},v_{3})du_{1}...dv_{3} \Big)(f)\Big)\Big|^{2}
\\
\label{ineg-e-4j1552}
% \Big|\sum _{(a_{1},\dots,a_{p},S_{1},...,S_{m},(\mathcal Y_{i}^{j})_{i\in \{0,\dots,m\}, j\in \{1,\dots ,l_{i}\}}) \in 
%(\pi_{x}^{p,k,m,(l_{0},...,l_{m})})^{-1}(Z)} ((\theta^{\flat}_{x}-\theta_{x}) (f))(a_1,...,a_p)\Big|^{2}$$ 
 \leq C \kappa 
\sum_{\tilde Z\in \Lambda_{Z}} \big(r_{0}(\tilde Z)+1\big)^{-6}
\big|\xi_{\tilde Z}(f)\big|^{2}.\end{gather}

Or  (\ref{ineg-eta-uuuvvv}) montre qu'il existe $C=C(\de,K,N,Q,P,M,s,B)$ tel que 
 \begin{gather}\label{ineg-e-4j1556}  \int_{ [0,1[^{6}}\Big|\xi_{Z}\Big(\big((\theta^{\flat}_{ x})_{u_{1},u_{2},u_{3}}^{v_{1},v_{2},v_{3}}-
\theta_{ x}\big)(f)\Big)\Big|^{2}du_{1}...dv_{3} 
\leq 
% \Big|\sum _{(a_{1},\dots,a_{p},S_{1},...,S_{m},(\mathcal Y_{i}^{j})_{i\in \{0,\dots,m\}, j\in \{1,\dots ,l_{i}\}}) \in 
%(\pi_{x}^{p,k,m,(l_{0},...,l_{m})})^{-1}(Z)} ((\theta^{\flat}_{x}-\theta_{x}) (f))(a_1,...,a_p)\Big|^{2}$$ 
 C 
\sum_{\tilde Z\in \Lambda_{Z}} \big(r_{0}(\tilde Z)+1\big)^{-6}
\big|\xi_{\tilde Z}(f)\big|^{2}.\end{gather}

D'autre part on vérifie facilement qu'il existe une constante $C$ telle que pour $i\in \{1,2\}$, on ait 
 \begin{gather}\label{ineg-e-4j1601}  \int_{ [0,1[^{6}}
 |h_{\kappa,i}(u_{1},u_{2},u_{3},v_{1},v_{2},v_{3})|^{2}du_{1}...dv_{3}\leq C\kappa.
 \end{gather}

Par Cauchy-Schwarz et grâce à (\ref{ineg-e-4j1556}) et (\ref{ineg-e-4j1601}),  on obtient (\ref{ineg-e-4j1552}). 
Ceci termine la démonstration du sous-lemme~\ref{arg-suppl-e-4j1148}. 
 \cqfd

\noindent {\bf Fin de la démonstration du e).}
D'après le lemme~\ref{continuite-del-J-conj-d}, il existe $C=C(\de,K,N,Q,P,M,s,B,T)$ tel que 
\begin{gather}\label{ineg-uK-e-4j1210}\sup_{\tau\in [0,T]}\|e^{\tau\theta_{x}} u_{x,r}K_{x}e^{-\tau\theta_{x}}\|_{\L(\H_{x,s}(\Delta_{p-1}),\H_{x,s}(\Delta_{p}))}\leq C. \end{gather}

En combinant (\ref{ineg-penult-e-4j1149}), le sous-lemme~\ref{arg-suppl-e-4j1148} et  (\ref{ineg-uK-e-4j1210}) on voit qu'il existe $C=C(\de,K,N,Q,P,M,s,B,r,T)$ tel  que
  pour tout $n> P$, 
  \begin{gather}\nonumber \sup_{\tau\in [0,T]}
    \big\|(\P_{n}-\P_{n-1})[(\theta^{\flat}_{x}-\theta^{\flat}_{x'}),e^{\tau\theta_{x}}u_{x,r}K_{x}e^{-\tau\theta_{x}}]
 \big\|_{\L(\H_{x,s}(\Delta_{p-1}),\H_{x,s}(\Delta_{p}))}^{2} \\ \label{ineg-fin-e-4j1215} \leq  C(n+1)^{-\frac{\alpha}{2}}.\end{gather}

Pour tout $i\in \Z$ soit $T_{i}\in \L(\H_{x,s}(\Delta_{p-1}),\H_{x,s}(\Delta_{p}))$ défini par 
\begin{gather*}(\P_{n}-\P_{n-1})T_{i}(\P_{n'}-\P_{n'-1})\\ = 
(\P_{n}-\P_{n-1})[(\theta^{\flat}_{x}-\theta^{\flat}_{x'}),e^{\tau\theta_{x}}u_{x,r}K_{x}e^{-\tau\theta_{x}}](\P_{n'}-\P_{n'-1})\text{ \ \ si \ \ }n'-n=i \\ \text{et \ \ \ }(\P_{n}-\P_{n-1})T_{i}(\P_{n'}-\P_{n'-1}) = 
0\text{\ \ \ sinon.}\end{gather*}  
D'après le sous-lemme~\ref{ineg73-3j2030} (ou~\ref{just-r0-r1-r-xx'-3j2157}) on a $T_{i}=0$ sauf si $|i-r|\leq QF$. Il est clair que $[(\theta^{\flat}_{x}-\theta^{\flat}_{x'}),e^{\tau\theta_{x}}u_{x,r}K_{x}e^{-\tau\theta_{x}}]=\sum_{i}T_{i}$. 

D'après (\ref{ineg-fin-e-4j1215}) 
il existe donc  $C=C(\de,K,N,Q,P,M,s,B,r,T)$ tel  que
  pour tout $i\in \Z$, 
  \begin{gather}
    \big\|(1-\P_{n})T_{i}
 \big\|_{\L(\H_{x,s}(\Delta_{p-1}),\H_{x,s}(\Delta_{p}))}^{2} \nonumber  \leq  C(n+1)^{-\frac{\alpha}{2}}.\end{gather}

En sommant sur $i\in \{r-QF, ...,r+QF\}$ on en déduit qu'il existe $C=C(\de,K,N,Q,P,M,s,B,r,T)$ tel  que
  pour tout $n> P$, 
  \begin{gather}
    \big\|(1-\P_{n})[(\theta^{\flat}_{x}-\theta^{\flat}_{x'}),e^{\tau\theta_{x}}u_{x,r}K_{x}e^{-\tau\theta_{x}}]
 \big\|_{\L(\H_{x,s}(\Delta_{p-1}),\H_{x,s}(\Delta_{p}))}^{2} \nonumber  \leq  C(n+1)^{-\frac{\alpha}{2}}.\end{gather}
 Ceci termine la preuve de e) et donc celle du lemme~\ref{lemme-compacite-equiv3}. \cqfd

On a donc montré les lemmes~\ref{lemme-compacite-equiv2} 
et~\ref{lemme-compacite-equiv}, et la 
proposition~\ref{enonce-ppal-KKC01} qui était l'énoncé principal de ce paragraphe.

\section{Fin de l'homotopie}\label{construction-fin}

Ce paragraphe n'offre guère d'intérêt car il recopie quasiment
la fin du  paragraphe 2.3.4 et le paragraphe 2.3.5 de~\cite{kkban}. 

\begin{lem}\label{Tgrand-l2-12j}
Pour $T$ assez grand  
\begin{itemize}
\item a) on a  $\|e^{T\theta_{x}^{\flat}}K_{x} e^{-T\theta_{x}^{\flat}}\|_{\L(\bigoplus _{p=1}^{p_{max}}\ell^{2}(\Delta_{p}))}\leq \frac{1}{2}$, 
\item b) il existe $C$ tel que 
$
\|e^{T\theta_{x}^{\flat}}\tilde H_{x} e^{-T\theta_{x}^{\flat}}\|_{\L(\bigoplus _{p=1}^{p_{max}}\ell^{2}(\Delta_{p}))}\leq C
$ et $$
\|e^{T\theta_{x}^{\flat}}u_{x,r} e^{-T\theta_{x}^{\flat}}\|_{\L(\bigoplus _{p=1}^{p_{max}}\ell^{2}(\Delta_{p}))}\leq C2^{-r}. $$
\end{itemize}
\end{lem}
\noindent{\bf Démonstration.} Grâce à 2) a) du lemme~\ref{support-connaissance-tildeH}, 
$K_{x}(e_{S})$ 
est une combinaison de $e_{T}$ où $T$ vérifie  
$d(x,T)\leq  d(x,S)-(Q\frac{N-6\delta}{p_{\max}}-2N-4\de)$, d'où 
$\rho_{x}^{\flat}(T)\leq \rho_{x}^{\flat}(S)-(Q\frac{N-6\delta}{p_{\max}}-3N-11\de)$. On suppose $(Q\frac{N-6\delta}{p_{\max}}-3N-11\de)\geq 1$, ce qui est permis par $(H_{Q})$. Ceci permet de montrer que la condition a) est réalisée pour $T$ assez grand et pour b) on utilise 
 1)a) et  3) de la proposition~\ref{recap-supp-connaiss-H-uK},  
la proposition~\ref{supp-uxrt} et le lemme~\ref{28dec1029}. \cqfd

On prend $T$ assez grand pour que  les conditions du lemme~\ref{Tgrand-l2-12j} soient satisfaites. En particulier 
 $e^{T\theta_{x}^{\flat}}(\del + J_{x}\del J_{x}) e^{-T\theta_{x}^{\flat}}$ est  continu sur $\bigoplus _{p=1}^{p_{max}}\ell^{2}(\Delta_{p})$. 
On peut relier les espaces de Hilbert $\H_{x,s}$ et 
$\bigoplus _{p=1}^{p_{max}}\ell^{2}(\Delta_{p})$ par un champ continu d'espaces de Hilbert $(\H_{x,s,\alpha})_{\alpha\in [0,1]}$, défini par  $\|.\|_{\H_{x,s,\alpha}}^{2}=\alpha \|.\|_{\ell^{2}}+(1-\alpha)\|.\|_{\H_{x,s}}$ sur $\bigoplus _{p=1}^{p_{max}}\C^{(\Delta_{p})}$. 

\begin{lem}\label{homotopie1-12j} 
$$((\H_{x,s,\alpha})_{\alpha\in [0,1]},e^{T\theta_{x}^{\flat}}(\del + J_{x}\del J_{x}) e^{-T\theta_{x}^{\flat}})$$ appartient à 
$KK_{G,2s\ell+C}(\C,\C[0,1])$ 
et réalise donc  une homotopie entre 
\begin{gather}\label{kk1-11j}(\H_{x,s},e^{T\theta_{x}^{\flat}}(\del + J_{x}\del J_{x})e^{-T\theta_{x}^{\flat}})\\ \label{kk2-11j}\text{\ \  \ \  et\ \ \ \  }(\bigoplus _{p=1}^{p_{max}}\ell^{2}(\Delta_{p}),e^{T\theta_{x}^{\flat}}(\del +J_{x}\del J_{x}) e^{-T\theta_{x}^{\flat}})\end{gather} 
\end{lem}
\noindent{\bf Démonstration.} La continuité de  $e^{T\theta_{x}^{\flat}}(\del + J_{x}\del J_{x})e^{-T\theta_{x}^{\flat}}$ résulte de la proposition~\ref{continuite-del-J-conj}, du lemme~\ref{Tgrand-l2-12j}
et du fait que pour tout opérateur $U$ on a 
\begin{gather}\label{fait-T-12j}\|U\|_{\L(\H_{x,s,\alpha})}\leq \max (\|U\|_{\L(\H_{x,s})},\|U\|_{\L(\ell^{2})}). \end{gather}
 Pour montrer que $e^{T\theta_{x}^{\flat}}(\del + J_{x}\del J_{x})e^{-T\theta_{x}^{\flat}}$ est équivariant à compact près, on voit, en reprenant l'argument de la démonstration de la proposition~\ref{lemme-compacite-equiv} et du lemme~\ref{lemme-compacite-equiv2}, qu'il suffit de montrer l'énoncé analogue à celui du lemme~\ref{lemme-compacite-equiv3} obtenu en prenant $\tau=T$ mais en rempla\c cant les normes $\|.\|_{\H_{x,s}}$ par les normes $\|.\|_{\H_{x,s,\alpha}}$ et en prenant le supremum sur $\alpha$. Par exemple l'énoncé analogue à a) du lemme~\ref{lemme-compacite-equiv3} est que 
$$\sup_{\alpha\in [0,1]}\|(1-\P_{n})e^{T\theta_{x}}(h_{x}-h_{x'}) e^{-T \theta_{x}}
\|_{\L(\H_{x,s,\alpha}(\Delta_{p-1}),\H_{x,s,\alpha}(\Delta_{p}))}$$ tend vers $0$ quand $n\to \infty$. Par (\ref{fait-T-12j}) 
 il suffit  de montrer que 
\begin{gather}\label{eq-compact-12j1439}\|(1-\P_{n})e^{T\theta_{x}}(h_{x}-h_{x'}) e^{-T\theta_{x}}
\|_{\L(\ell^{2}(\Delta_{p-1}),\ell^{2}(\Delta_{p}))}\to 0 \text{\ \ quand \ \ } n\to \infty. \end{gather}
Cela résulte du lemme~\ref{lem-psix-psix'-8j1630} (et même simplement du fait, mentionné avant le lemme~\ref{lem-psix-psix'-8j1630}, que 
$\|\psi_{S,x}-\psi_{S,x'}\|_{1}$ tend vers $0$ en dehors des parties finies de $\Delta$). Il reste à montrer les énoncés analogues à (\ref{eq-compact-12j1439}) correspondant aux opérateurs de b), c), d) et e) du lemme~\ref{lemme-compacite-equiv3}.  
Pour cela on utilise de plus la propriété (\ref{fortement}) et le lemme~\ref{mesure-uuuvvv-dbemol} (ou même l'énoncé plus faible qui est la propriété (\ref{condition-d'}) avec $d^{\flat}$ au lieu de $d''$). 
\cqfd

En fait (\ref{kk2-11j}) appartient à $KK_{G,0}(\C,\C)$ (c'est-à-dire qu'il vérifie 
les mêmes conditions qu'un élément de $KK_{G}(\C,\C)$ sauf celle qui assure que l'opérateur est auto-adjoint à compact près). De plus 
l'image de (\ref{kk2-11j})  dans $KK_{G,2s\ell+C}(\C,\C)$ est égale à $1$ d'après le lemme~\ref{homotopie1-12j} et les propositions~\ref{enonce-ppal-KKC01} et~\ref{homot-1-142}. Pour  terminer la preuve du théorème~\ref{hyperb-bon}, il suffit donc de montrer le lemme suivant, dont la preuve suit de très près  le paragraphe 2.3.5 de~\cite{kkban}.

  \begin{lem}\label{homotopie2-12j} 
L'élément   (\ref{kk2-11j}) 
%$$(\bigoplus _{p=1}^{p_{max}}\ell^{2}(\Delta_{p}),
%e^{T\theta_{x}^{\flat}}(\del + J_{x}\del J_{x}) e^{-T\theta_{x}^{\flat}})$$ 
 est égal à l'image de $\gamma\in KK_{G}(\C,\C)$ dans $KK_{G,0}(\C,\C)$.
\end{lem}
\noindent{\bf Démonstration.}
D'après le  a) du lemme~\ref{Tgrand-l2-12j} et comme $$K_{x}=(1-\del h_{x}-h_{x}\del )^{Q},$$ on a 
 \begin{gather}\label{rayon-spec-12j}\rho(e^{T\theta_{x}^{\flat}}(1-\del h_{x}-h_{x}\del )e^{-T\theta_{x}^{\flat}})
 \leq 2^{-Q^{-1}}\end{gather} 
 où $\rho$ désigne 
 le rayon spectral dans $\L(\ell^{2}(\Delta_{p}))$. 
  On rappelle que $$H_{x}=h_{x}(\del h_{x}+h_{x}\del)^{-1}.$$ Donc $e^{T\theta_{x}^{\flat}}H_{x}e^{-T\theta_{x}^{\flat}}$ agit continûment sur $\bigoplus _{p=1}^{p_{max}}\ell^{2}(\Delta_{p})$. 

\begin{souslem} L'élément 
\begin{gather}\label{kk3-11j}(\bigoplus _{p=1}^{p_{max}}\ell^{2}(\Delta_{p}),
e^{T\theta_{x}^{\flat}}(\del +H_{x}) e^{-T\theta_{x}^{\flat}})\end{gather}
appartient à $KK_{G,0}(\C,\C)$ et (\ref{kk2-11j}) est homotope à
(\ref{kk3-11j}). 
\end{souslem}
\noindent{\bf Démonstration.} La continuité de $e^{T\theta_{x}^{\flat}}(\del +H_{x}) e^{-T\theta_{x}^{\flat}}$ a été justifiée avant le lemme et l'équivariance à compact près résulte de (\ref{eq-compact-12j1439}). On a $H_{x}^{2}=0$ car on a vu dans le lemme~\ref{lemmeh2=0} que $h_{x}^{2}=0$. 
L'homotopie entre (\ref{kk2-11j}) et 
(\ref{kk3-11j}) résulte du 
lemme 1.4.1 de~\cite{kkban}. \cqfd 
   
   On pose $D_{x}=\del(\del h_{x}+h_{x}\del)^{-1}$. 
   
   \begin{souslem} L'élément 
\begin{gather}\label{kk4-11j}(\bigoplus _{p=1}^{p_{max}}\ell^{2}(\Delta_{p}),e^{T\theta_{x}^{\flat}}(h_{x}+D_{x}) e^{-T\theta_{x}^{\flat}})\end{gather}
appartient à $KK_{G,0}(\C,\C)$ et (\ref{kk3-11j}) est homotope à
(\ref{kk4-11j}). 
\end{souslem}
\noindent{\bf Démonstration.} La continuité de $e^{T\theta_{x}^{\flat}}(\del +D_{x}) e^{-T\theta_{x}^{\flat}}$ vient de nouveau de (\ref{rayon-spec-12j})   et  son équivariance à compact près résulte de (\ref{eq-compact-12j1439}). 
Pour réaliser l'homotopie entre (\ref{kk3-11j}) et 
(\ref{kk4-11j}) on procède 
   comme dans le lemme 2.3.11 de~\cite{kkban}. Pour $\alpha\in [0,1]$ on définit grâce à (\ref{rayon-spec-12j}) et en utilisant la détermination principale du logarithme, 
   $$H_{x}^{\alpha}=h_{x}(\del h_{x}+h_{x}\del)^{-\alpha}
   \text{\  et \ } D_{x}^{\alpha}=\del(\del h_{x}+h_{x}\del)^{-\alpha}.$$
   Alors 
   \begin{gather}\label{kk34-11j}\big(\big(\bigoplus _{p=1}^{p_{max}}\ell^{2}(\Delta_{p})\big)[0,1],(e^{T\theta_{x}^{\flat}}(H_{x}^{1-\alpha}+D_{x}^{\alpha}) e^{-T\theta_{x}^{\flat}})_{\alpha\in [0,1]}\big)\end{gather}
appartient à $KK_{G,0}(\C,\C[0,1])$ et réalise une homotopie entre (\ref{kk3-11j}) et 
(\ref{kk4-11j}). 
   \cqfd

 On rappelle maintenant quelques notations de~\cite{kkban}.  On note $\phi_{S,x}=\sqrt{\psi_{S,x}}$. On note $f,h'$  les opérateurs de $\ell^{2}(\Delta_{p})$ dans $\ell^{2}(\Delta_{p+1})$ donnés par 
 $$f(e_{S})=\phi_{S,x}\wedge e_{S}\text{\ \ \ et\ \ \ } h'(e_{S})=\frac{1}{\|\phi_{S,x}^{3}\|_{1}}\psi_{S,x}\wedge e_{S}.$$
 On note $g,g'$  les opérateurs de $\ell^{2}(\Delta_{p})$ dans $\ell^{2}(\Delta_{p-1})$ donnés par 
 $$g(e_{S})=\phi_{S,x}\lrcorner e_{S}\text{\ \ \ et\ \ \ } g'(e_{S})=\frac{1}{\|\phi_{S,x}^{3}\|_{1}}\phi_{S,x}\lrcorner e_{S}.$$ 
 On note aussi $h=h_{x}$ pour être cohérent avec~\cite{kkban}.
 
   Alors (\ref{kk4-11j}) est homotope à \begin{gather}\label{kk5-11j}(\bigoplus _{p=1}^{p_{max}}\ell^{2}(\Delta_{p}),e^{T\theta_{x}^{\flat}}(h+g') e^{-T\theta_{x}^{\flat}})\end{gather}
   grâce au 
lemme 1.4.1 de~\cite{kkban}. Puis (\ref{kk5-11j}) est homotope à \begin{gather}\label{kk6-11j}(\bigoplus _{p=1}^{p_{max}}\ell^{2}(\Delta_{p}),h+g')\end{gather}
   par l'homotopie évidente 
$$\big(\big(\bigoplus _{p=1}^{p_{max}}\ell^{2}(\Delta_{p})\big)[0,T],\big(e^{\tau \theta_{x}^{\flat}}(h+g') e^{-\tau \theta_{x}^{\flat}}\big)_{\tau\in [0,T]}\big).$$

On note que l'homotopie entre (\ref{kk4-11j}) et (\ref{kk6-11j}) correspond au  
 lemme 2.3.12 de~\cite{kkban}. 
 
 Enfin (\ref{kk6-11j}) est homotope  à
$(\bigoplus _{p=1}^{p_{max}}\ell^{2}(\Delta_{p}),h'+g)$
puis à \begin{gather}\label{kk8-11j}(\bigoplus _{p=1}^{p_{max}}\ell^{2}(\Delta_{p}),f+g)\end{gather}
(voir les lemmes 2.3.13 et 2.3.14 de~\cite{kkban}). 

Or (\ref{kk8-11j})  est égal à l'image de $\gamma\in KK_{G}(\C,\C)$ dans $KK_{G,0}(\C,\C)$. On renvoie à~\cite{ks} pour la construction complète de $\gamma$, qui est  rappelée au début de la section 2 de~\cite{kkban}. Il y a une toute petite subtilité due au fait que dans~\cite{ks} la construction utilise des mesures $   \psi_{S,x}^{KS}$ qui sont légèrement différentes des mesures $\psi_{S,x}$ que nous avons définies. Cependant on peut les relier par une homotopie
$\alpha\mapsto (1-\alpha)\psi_{S,x}+\alpha \psi_{S,x}^{KS}$. En effet grâce au lemme~\ref{73} et au 
 lemme 6.3 de~\cite{ks}, pour tout $T$ tel que $e_{T}$ apparaisse dans 
 $\psi_{S,x}\wedge e_{S}, \psi_{S,x}\lrcorner e_{S}, \psi_{S,x}^{KS}\wedge e_{S}$ ou $\psi_{S,x}^{KS}\lrcorner e_{S}$, on a 
 $\psi_{S,x}=\psi_{T,x}$ et $\psi_{S,x}^{KS}=\psi_{T,x}^{KS}$. 
 Donc les opérateurs analogues à $f,g,g',h,h'$ construits à l'aide de 
 $(1-\alpha)\psi_{S,x}+\alpha \psi_{S,x}^{KS}$ sont de carré nul. 
 Ceci termine la preuve du lemme~\ref{homotopie2-12j}.  \cqfd

On a donc terminé la preuve de 
théorème~\ref{hyperb-bon}.

\vskip 1cm

\noindent
{\Large \bf Index des notations.} 

\vskip 5mm

\noindent Page \pageref{H{delta}(x,y,z,t)}   : $H_{\delta}(x,y,z,t)$ 

\noindent Page \pageref{def-Bxr} : 
$B(x,r)$

\noindent Page \pageref{ddmax} : $d(A,B)$, $d_{\max}(A,B)$, 

\noindent Page \pageref{defHn} : 
$N$, $(H_{N})$, $\Delta$, $\Delta_{p}$, $\C^{(\Delta_{p})}$, $p_{\max}$

\noindent Page \pageref{def-del-page} : $\del$

\noindent Page \pageref{Hdelta-beta} : $(H_{\delta}^{\beta}(x,a,b,c))$, $(H_{\delta}^{0}(x,a,b,c))$, $U_{S}$

\noindent Page \pageref{def-Asx} : $A_{S,x}$, $Y_{S,x,r}$, $A_{\emptyset,x}$, $Y_{\emptyset,x,r}$

\noindent Page \pageref{psiSx} : $\chi_{A}$, $\nu_{A}$, $\psi_{S,x}$, $\psi_{S,x,t}$, $E(.)$

\noindent Page \pageref{def-hx-hxs} : $h_{x}$, $h_{x,t}$

\noindent Page \pageref{zetaxS} : $H_{x}$, $\zeta_{x}(S)$

\noindent Page \pageref{Phip-page} : $\Phi_{p}$

\noindent Page \pageref{murtxa} : $A_{x,a,r,k}$, $\mu_r(x,a)$, $\mu_{r,t}(x,a)$

\noindent Page \pageref{def-u-v} : $u_{x}$, $u_{x,r,t}$, $u^{p}_{x,r,t}$, $v^{p}_{x,r,t}$

\noindent Page \pageref{hyp-HQ} : $(H_{Q})$, $\tilde H_{x}$, $K_{x}$

\noindent Page \pageref{Hxqttt} : $K_{x,q,(t_{1},\dots ,t_{q})}$, $\tilde H_{x,q,(t_{1},\dots ,t_{q})}$, $B(A,r)$

\noindent Page \pageref{def-Jx} : $J_{x}$ 

\noindent Page \pageref{def-F} : $F$

\noindent Page \pageref{Yrxy} : $Y_{x,y}^{r}$, $\Lambda_{x,r}^{y,s}$, $\alpha_{r'\leftarrow r}^{s'\leftarrow s}$ 

\noindent Page \pageref{Lambdaxysrsr} : $\Lambda_{x,r',r}^{y,s',s}$

\noindent Page \pageref{Lambda123} : $\Lambda_{x,r_{1},r_{2},r_{3}}^{y,s_{1},s_{2},s_{3}}$, $\beta_{x,r}^{y,s}$, $\tilde \Lambda_{x,r_{1},r_{2},r_{3}}^{y,s_{1},s_{2},s_{3}}$

\noindent Page \pageref{def-Auuuvvv} : $A_{x,u_{1},u_{2},u_{3}}^{y,v_{1},v_{2},v_{3}}$, $ d^{\flat}(x,y)$, ${d^{\flat}}_{u_{1},u_{2},u_{3}}^{v_{1},v_{2},v_{3}}(x,y)$ 

\noindent Page \pageref{def-rhobemol} : $\rho^{\flat}_{x}(a)$, $\rho^{\flat}_{x}(S)$

\noindent Page \pageref{defi-Y} : $(H_{P})$, $Y_{x}^{p,k,m,(l_{0},...,l_{m})}$, 
$\mathcal Y_{i}^{j}$

\noindent Page \pageref{overlinemathcalY} : $(H_{M})$, $\overline Y_{x}^{p,k,m,(l_{0},...,l_{m})}$,  $\pi_{x}^{p,k,m,(l_{0},...,l_{m})}$, $r_{0}(Z),\dots,r_{m}(Z), s_{0}(Z),\dots ,  s_{m}(Z)$

\noindent Page \pageref{defxi-31dec1554} : $(H_{B})$, $(H_{\alpha})$, $\xi_{Z}$, $\H_{x,s}(\Delta_{p})$

\noindent Page \pageref{enonce-ppal-KKC01} : $\theta^{\flat}_{x}$

\noindent Page \pageref{approx-arbres} : $\Psi$

\noindent Page \pageref{Yxrightarrow} : 
$Y_{x}^{\rightarrow,p,k,m,(l_{0},...,l_{m})} $,  $\overline Y_{x}^{\rightarrow,p,k,m,(l_{0},...,l_{m})}$,  $\pi_{x}^{\rightarrow,p,k,m,(l_{0},...,l_{m})}$

\noindent Page \pageref{norme-arrow} : 
$\|.\|_{\H^{\rightarrow}_{x,s}(\Delta_{p})}$

\noindent Page \pageref{def-mathcalP} : 
$\P$

\noindent Page \pageref{defi-Y-natural} : 
 $Y_{x}^{\natural,p,k,m,(l_{0},...,l_{m}),\lambda_{0},\lambda_{1}}$, $\mathcal Z_{i}^{j}$

\noindent Page \pageref{rijmax} : 
$\overline Y_{x}^{\natural,p,k,m,(l_{0},...,l_{m}),\lambda_{0},\lambda_{1}}$,
 $\pi_{x}^{\natural,p,k,m,(l_{0},...,l_{m}),\lambda_{0},\lambda_{1}}$, $r_{i,j}^{\max}(Z)$,  $t_{i}^{j}(Z)$ 
 
\noindent Page \pageref{formule-norme-mu} : 
$ \H^{\natural,\mu_{0},\mu_{1}}_{x,s}(\Delta_{p})$

\noindent Page \pageref{kappa-sigma} : 
 $\kappa_{\sigma}$,  $\overline \kappa_{\sigma,\infty}$, $\overline \kappa_{\sigma,i}$
 
\noindent Page \pageref{def-rhox-rhox0Z} : $\rho_{x}$,  $\theta_{x}$, 
$\rho_{x}^{0}(Z)$, $\rho_{x}^{1}(Z)$

\noindent Page \pageref{rhoflatxuuuvvv} : 
$(\rho^{\flat}_{x})_{u_{1},u_{2},u_{3}}^{v_{1},v_{2},v_{3}}$

\noindent Page \pageref{def-r0'-k'-xx'} : $\pi_{x,x'}^{p,k,m,(l_{0},...,l_{m})}$, $\overline Y_{x,x'}^{p,k,m,(l_{0},...,l_{m})}$, 
$r_{0}'(Z)$, 
$k'(Z)$

\noindent Page \pageref{est-card-fibres-star} : 
$\overline Y_{x,x',\star}^{p,k,m,(l_{0},...,l_{m})}$, $\pi_{x,x',\star}^{p,k,m,(l_{0},...,l_{m})}$

\noindent Page \pageref{def-Pn} : 
   $\P_{n}$

\end{document}